\documentclass[11pt,a4paper]{article}

\usepackage[hyperref,english]{package_RFb}
\usepackage{numbered_constants}
\numberwithin{equation}{section}

\usepackage[top=3cm, bottom=3cm, left=3cm, right=3cm]{geometry}
\usepackage[style=alphabetic,backend=bibtex,natbib=true,url=false,doi=false]{biblatex}

\newcommand*\samethanks[1][\value{footnote}]{\footnotemark[#1]}

\title{Fluctuations for fully pushed stochastic fronts}
\author{Alison Etheridge\thanks{Department of Statistics, University of Oxford, UK}
 \and Rapha\"el Forien\thanks{BioSP, INRAE Avignon, France}
 \and Thomas Hughes\thanks{Department of Mathematical Sciences, University of Bath, UK}
  \and Sarah Penington\samethanks[3] }

\NewDocumentCommand{\Hnorm}{m m}{\| #1 \|_{\tilde{H}^{#2,\alpha}}\xspace}
\NewDocumentCommand{\enorm}{m}{\| #1 \|_{e}\xspace}

\newcommand{\cW}{\mathcal{W}}
\newcommand{\bP}{\mathbb{P}}
\newcommand{\indc}{\mathds{1}}
\newcommand{\cF}{\mathcal{F}}

\DeclareConstant{Beta}{\beta}
\DeclareConstant{Epsilon}{\varepsilon}
\DeclareConstant{cexp}{c}
\DeclareConstant{Cst}{C}

\DeclareConstant{cinit}{a}
\DeclareConstant{Cinit}{\tilde{C}}
\DeclareConstant{ggamma}{\gamma}
\NewDocumentCommand{\dispSigmaN}{m}{\sigma_{#1,N}}
\DeclareCustomConstant{sigmaN}{\dispSigmaN}
\DeclareConstant{stime}{t}
\DeclareConstant{smalla}{b}
\DeclareConstant{bigC}{A}

\newcommand{\IR}{\mathbb R}

\NewDocumentCommand{\dist}{}{\mathrm{dist}\xspace}
\NewDocumentCommand{\drift}{}{\mu_{(\alpha)}\xspace}
\NewDocumentCommand{\variance}{}{\sigma_{(\alpha)}\xspace}
\NewDocumentCommand{\specgap}{}{\lambda_0\xspace}
\NewDocumentCommand{\genericc}{}{c\xspace}
\NewDocumentCommand{\Tr}{m}{\mathrm{Tr}[#1]\xspace}

\NewDocumentCommand{\Vphi}{m}{\varphi_{#1}\xspace}
\NewDocumentCommand{\Vpsi}{m}{\psi_{#1}\xspace}

\definecolor{darkgreen}{rgb}{0,0.5,0}

\begin{document}
	\date{}
	\maketitle
	
	\begin{abstract}
		We study the asymptotic behaviour, in the small noise limit, of stochastic travelling wave solutions to reaction-diffusion equations perturbed by Wright-Fisher noise.
		Such equations are predicted to display three distinct responses to noise in three parametric regimes: fully pushed, semi-pushed, and pulled.
		We prove, for the entire fully pushed regime, that solutions are asymptotically close to a stochastic shift of the deterministic travelling wave, and characterize the limiting shift process as a Brownian motion with drift.
		This gives the first full fluctuation theorem demonstrating fully pushed phenomenology for a non-linear stochastic reaction-diffusion equation and verifies a physical conjecture of Birzu, Hallatschek and Korolev \cite{birzu_fluctuations_2018}. 
		
		The proof uses an infinite-dimensional version of a method introduced by Katzenberger \cite{katzenberger_solutions_1991}, as pioneered by Funaki \cite{funaki_scaling_1995}. 
		This approach views the dynamics as a stochastic perturbation of a dynamical system (the PDE) with strong drift towards an invariant manifold, in our case the set of shifts of the travelling wave profile,
		and gives an expression for the stochastic motion ``along'' this manifold. 
		Implementing this method in our setting requires many ingredients, including a close analysis of the dynamics of the corresponding PDE, integrability and regularity properties of solutions to the SPDE, 
		and sharp control of the position of the right endpoint of the solution's support.
	\end{abstract} 

	\tableofcontents

\section{Introduction}

\subsection{Background}

Reaction-diffusion equations are ubiquitous as mathematical models of 
spreading phenomena. They arise in 
domains ranging from population genetics to combustion. Under rather 
mild assumptions they exhibit travelling wave solutions.
There may be many such solutions; in what follows the one with the smallest
possible velocity, $\alpha$, will be written  
as $u(t,x)=m(x-\alpha t)$ for some
wave profile $m$ (which is unique up to a shift).
Such waves are
classically classified as `pulled' or `pushed' \cite{stokes_two_1976}. 
For a pulled wave, 
the minimal wave speed is determined by linearisation at the leading edge, 
whereas for 
pushed waves, the wave speed is strictly larger than the linear prediction,
being influenced by nonlinear growth in the bulk. The most famous
example of a pulled wave arises from the Fisher-KPP equation, introduced
by Fisher \cite{fisher_wave_1937} to model the spread of an advantageous genetic type
through a linear habitat, and so the linear prediction for the wave speed is
often called the Fisher speed.  

The biological and physical systems being modelled are almost always 
subject to noise, which makes it important to understand stochastic 
perturbations of the deterministic equation; for a broad survey of stochastic travelling waves, see \cite{kuehn_travelling_2020}.
 It has long been known that pulled waves are very sensitive to noisy perturbations,  
and a large literature is devoted to understanding the effect on the 
wavefront of stochastic perturbations of equations of Fisher-KPP type; see
e.g.~\cite{mueller_random_1995}, \cite{brunet_shift_1997}, 
\cite{conlon_travelling_2005}, \cite{brunet_phenomenological_2006}, 
\cite{mueller_effect_2011}.
The differing levels of 
sensitivity to stochastic perturbations of pushed waves led 
\cite{birzu_fluctuations_2018} 
to further divide this class into 
semi-pushed and fully pushed, with 
the transitions between all
three regimes (pulled, semi-pushed, and fully pushed)
governed by the ratio of the actual wave speed to the Fisher wave speed. 

%
%

Here we shall be interested in the fully pushed regime.
We consider a stochastic partial differential equation (SPDE) obtained by
perturbing the underlying reaction-diffusion equation by a 
small Wright-Fisher noise (see Definition~\ref{def:u}). The form of the noise
arises naturally in population genetics, where it captures the 
randomness due to reproduction in a finite population; see, for
example, \cite{shiga_stepping_1988}. 
The variance of the noise is proportional to $1/N$, where the parameter $N$
can be interpreted as population density,
and we are interested in the behaviour of the solution to the SPDE in
the limit as $N\to\infty$.

We shall say that a function $u:\IR\rightarrow [0,1]$ has the compact 
interface property if $\inf\{x\in\IR: u(x)<1\}>-\infty$, and 
$\sup\{x\in\IR: u(x)>0\}<\infty$. 
Our main result, Theorem~\ref{thm:main_result}, shows that 
(under some technical conditions)
in the fully pushed setting, 
if we choose an initial condition $u_0^N$ 
with the compact interface property, and
which is sufficiently close in shape to 
the travelling wavefront $m$,
then the 
solution retains the compact support property, 
and remains close to a stochastic shift of $m$.
Moreover, for any $T>0$, as $N\to\infty$,
the stochastic process describing the shift at times $\{Nt\}_{t\in [0,T]}$ 
converges to a Brownian motion with drift, and we provide explicit expressions for both its drift and its quadratic variation.
In particular, the noise primarily affects the position of 
the wave at large times, not its shape.

\subsection{Summary of the main result} \label{sec:statement}

We now turn to a more precise statement of our main result. 
We first specify our assumptions on the form of
the deterministic equation; for detailed accounts of convergence to 
travelling wave solutions for such equations
we refer to \cite{volpert_traveling_1994}, \cite{gilding_travelling_2004}, and the 
references therein.
	
	Let $ f : [0,1] \to \R $ be $ C^2 $ with $ f'' $ Lipschitz continuous, such that
	\begin{align*}
		f(0) = f(1) = 0, && f'(1) < 0,
	\end{align*}
and there exists $ (\alpha, m) \in \R_+ \times C(\R,[0,1]) $ solving
	\begin{equation} \label{def:m}
		\left\lbrace
		\begin{aligned}
			& \partial_{xx} m + \alpha \partial_x m + f(m) = 0 \\
			& \lim_{x\to -\infty} m(x) = 1, \qquad \lim_{x\to \infty} m(x) = 0
		\end{aligned}
		\right.
	\end{equation}
	with $ m(x) \sim k e^{-\lambda_+ x} $ as $ x \to + \infty $, where
	\begin{equation} \label{eq:lambda+def}
		\lambda_+ = \frac{\alpha}{2} + \left( \frac{\alpha^2}{4} - f'(0) \right)^{1/2}.
	\end{equation}
(To understand where this comes from, substitute a solution of the form $e^{-\lambda x}$
into~\eqref{def:m}.)
	We further assume that
	\begin{equation*}
		f'(0) < \frac{\alpha^2}{4}.
	\end{equation*}
	This implies that $ m $ is a \emph{pushed} travelling wave and that, when $ f'(0) \geq 0 $, it is the front of minimal velocity (when $ f'(0) < 0 $ then the travelling wave is unique).
	We also assume that
	\begin{equation*}
		\int_{0}^{1} f(u) du > 0,
	\end{equation*}
	which implies that $ \alpha > 0 $.
	Without loss of generality, we assume that $ \lambda_+ = 1 $ (by a scaling of time and space all values of $ \lambda_+ $ can be obtained).
	This implies that $ \alpha < 2 $ and
	\begin{align*}
		f'(0) = \alpha-1.
	\end{align*}
	Any $C^2$ reaction term $ f $ for which the initial value problem
	\begin{equation*}
		\left\lbrace
		\begin{aligned}
			\partial_t u &= \partial_{xx} u + f(u), \\
			u(0,\cdot) &= u_0,
		\end{aligned}
		\right.
	\end{equation*}
	admits pushed travelling waves with positive speeds as solutions is thus covered by these assumptions, including bistable, ignition and monostable reaction terms $ f $.
	Within this parametrization, the transition from fully pushed to semi-pushed waves predicted in \cite{birzu_fluctuations_2018} occurs at $\alpha = 3/2$.

The specific form of $f$ that motivated our work, and which it is helpful to
keep in mind, is
\begin{equation}
\label{f}
f(u)=u(1-u)(2u-1+\alpha),\qquad \alpha\in (0,2).
\end{equation}
With this choice,
\begin{equation*}
		m(x) = \frac{1}{1+e^x}.
	\end{equation*}
The interpretation from genetics is that we are 
modelling a gene occurring in two
forms (alleles), $a$ and $A$ say, and for which the relative 
fitness of individuals
carrying genetic types $aa$, $aA$ and $AA$ is $1+2\alpha s:1+(\alpha-1)s:1$
for some small parameter $s$ (which determines the relevant time-scaling of 
the model). The function $u$ records 
the proportion of $a$-alleles at each space-time point.
\begin{remark}
When $\alpha=0$, the reaction term~\eqref{f}
models a population in which homozygotes ($aa$
and $AA$)
are equally fit, but there is selection 
against heterozygotes ($aA$). 
In that case, $m(x)$ provides a stationary solution.

For $\alpha >1$, if we set $B=2/(\alpha-1)$ 
and $r_0=(\alpha -1)$
the reaction term~\eqref{f} becomes 
$r_0p(1-p)(1+Bp)$, which is the form 
considered by~\cite{birzu_fluctuations_2018}. The fully pushed regime then
corresponds to $B\geq 4$,
or equivalently $\alpha\leq 3/2$.

A more usual parametrization of 
relative fitnesses of $aa:aA:AA$ individuals would be 
$1+2\Sigma s:1+ h\Sigma s:1$, for $s$ small and $\Sigma>0$. 
The reaction term becomes $u(1-u)(2 (1-h)u+h)$ and
the transition from pushed to 
pulled waves is when $h=1/3$.
If $h=1$ we recover Fisher's equation, whereas
$h=2$ (corresponding to dominance of the $a$-allele) gives the 
equation considered by \cite{kolmogorov_study_1937}. 
\end{remark}

In order to take into account the randomness due to reproduction in a finite
population, we shall model the proportion of the advantageous allele by a 
random function $ u^N_t : \R \to [0,1] $ given as the solution to the 
following SPDE, which is an
extension of that introduced in \cite{shiga_stepping_1988}.

	\begin{definition} \label{def:u}
		Let $ (W(t), t \geq 0) $ be a cylindrical Wiener process on $ L^{2}(\R) $ (see for example Section~4.1.2 in \cite{prato_stochastic_2014}), and for $N\ge 1$, let $ u^N_0 : \R \to [0,1] $ be continuous.
		For $ N \geq 1 $, let $ (u^N_t, t \geq 0) $ be a mild solution to the stochastic partial differential equation
		\begin{equation} \label{spde_u_N}
		d u^N_t = \left( \partial_{xx} u^N_t + f(u^N_t) \right) dt + \sqrt{\frac{1}{N} u^N_t(1-u^N_t)} d W(t),
		\end{equation}
		with initial condition $ u^N_0 $.
	\end{definition}
	
	Recall (see e.g.~Section~7.1 in \cite{prato_stochastic_2014}) that an $ L^2(\R)$-valued process $ (u^N_t, t \geq 0) $ adapted to the filtration associated to $ (W(t), t \geq 0) $ is said to be a mild solution to~\eqref{spde_u_N} if, for any $ t \geq 0 $,
	\begin{equation*}
		\P{ \int_{0}^{t} \| u_s^N \|_{L^2}^2 ds < + \infty } = 1,
	\end{equation*}
	and
	\begin{equation} \label{mild_formulation}
		u^N_t = P(t) u^N_0 + \int_{0}^{t} P(t-s) f(u^N_s) ds + \int_{0}^{t} P(t-s) \left( \sqrt{\frac{1}{N} u^N_s (1-u^N_s)} dW(s) \right),
	\end{equation}
	where $ (P(t), t \geq 0) $ is the semigroup of linear operators generated by $ \partial_{xx} $ acting on $ L^2(\R) $, i.e.
	\begin{equation*}
		P(t) \phi(x) = \int_{\R} \frac{1}{\sqrt{4\pi t}} \exp\left( - \frac{(x-y)^2}{4 t} \right) \phi(y) dy, \quad \phi \in L^2(\R),
	\end{equation*}
	and the second integral on the right-hand side of \eqref{mild_formulation} is a stochastic integral with respect to the cylindrical Wiener process $ (W(t), t \geq 0) $ in the sense of \citep[Section~4.2.1]{prato_stochastic_2014}.
Existence and uniqueness in distribution for mild solutions of~\eqref{spde_u_N} was proved in~\cite{mueller_speed_2021}, along with a compact interface property that we state here. 

\begin{lemma}[Theorem 1.1 and Lemma 2.2 in~\cite{mueller_speed_2021}] \label{lem:compactinterface}
Suppose $u^N_0:\R\to [0,1]$ is measurable with 
\[
\inf\{x\in \R:u^N_0(x)<1\}>-\infty \quad \text{and}\quad \sup\{x\in \R:u^N_0(x)>0\}<\infty.
\]
Then~\eqref{spde_u_N} with initial condition $u_0^N$ has a mild solution $(u^N_t,t\ge 0)$, which is unique in law.
Moreover, $ x \mapsto u^N_t(x) $ is a continuous real-valued function for each $t> 0$,\
and, almost surely, for all $t\ge 0$,
\[
\inf_{s\in [0,t]}\inf\{x\in \R:u^N_s(x)<1\}>-\infty \quad \text{and}\quad \sup_{s\in [0,t]}\sup\{x\in \R:u^N_s(x)>0\}<\infty.
\]
\end{lemma}

	Even though $ (u^N_t, t \geq 0) $ depends on the value of $ \alpha $, 
we shall consider $ \alpha $ to be fixed to some arbitrary value in $(0,3/2)$ 
and suppress dependence on $\alpha$ in our notation.
We will however sometimes state results starting with ``for any 
$ \alpha \in \ldots $'' to highlight which intermediate results rely on 
$ \alpha $ being smaller than some specific value, and which ones are valid 
for all $\alpha \in (0,2)$.
	
	When $ N \to \infty $, assuming that $ u^N_0 $ converges to some measurable $ u^0 $, for any $T>0$, $ (u^N_t, t \in [0,T]) $ converges in probability to the solution $ (u_t, t \in [0,T]) $ of the deterministic partial differential equation
	\begin{equation} \label{pde_u}
		\left\lbrace
		\begin{aligned}
			&\partial_t u_t = \partial_{xx} u_t + f(u_t), \\
			&u_0 = u^0.
		\end{aligned}
		\right.
	\end{equation}
For $ \alpha \in (0,1) $, which corresponds to $ f $ being a bistable reaction term, Fife and McLeod \citep[Theorem~3.1]{fife_approach_1977} 
proved that, subject to some mild conditions on $ u^0 $ (in particular, they 
require that $0\le u^0\le 1$ and also $\limsup_{x\to \infty}u^0(x)$ and 
$\limsup_{x\to -\infty}(1-u^0(x))$ are both sufficiently small), there 
exists $ \zeta \in \R $ (depending on $ u^0 $) such that
	\begin{equation*}
		\sup_{x \in \R} | u_t(x) - m(x - \alpha t - \zeta) | \to 0 \quad \text{ as } t \to \infty.
	\end{equation*}
	This result was extended by Rothe \cite{rothe_convergence_1981} to $\alpha \in (1,2)$, assuming that  $u^0(x) \leq k e^{-\lambda x}$ for some $\lambda > \alpha-1$. 
When $\alpha=1$ the equation is sometimes called the Zeldovich equation. 
The wave is still pushed, but because of the degeneracy at $u=0$, 
the equation is not amenable to the phase plane methods exploited in parts of
the analyses above; 
an alternative approach which captures this 
case can be found in \cite{gilding_travelling_2004}.

If we set $ \Phi_t(u^0,x) := u_t(x + \alpha t) $, then $ v \mapsto \Phi_t(v,\cdot) $ is the flow of the partial differential equation \eqref{pde_u} in the moving frame (at speed $ \alpha $), i.e., for any $ v \in C(\R, [0,1]) $,
	\begin{equation} \label{eq:Phidefn}
		\left\lbrace
		\begin{aligned}
			& \partial_t \Phi_t(v,x) = \partial_{xx} \Phi_t(v,x) + \alpha \partial_x \Phi_t(v,x) + f(\Phi_t(v,x)) \\
			& \Phi_0(v,x) = v(x).
		\end{aligned}
		\right.
	\end{equation}
The results above
can then be reformulated as stating the existence of $ \zeta \in \R $ 
such that $ \Phi_t(v,\cdot) \to m(\cdot - \zeta) $ uniformly as 
$ t \to \infty $ (see also Theorem~\ref{thm:zeta} and 
Corollary~\ref{cor:uniform_cvg} below for a more precise statement under 
stronger assumptions on the initial data).
	In the remainder of the article, for $ \eta \in \R $, we set
	\begin{equation} \label{def_meta}
	m_\eta(x) := m(x-\eta)\quad \text{for }x\in \R,
	\end{equation}
	and we define
	\begin{equation} \label{def:manifoldM}
		M := \lbrace m_\eta : \eta \in \R \rbrace.
	\end{equation}
	Then $ M $ is a stable manifold of fixed points for the deterministic flow $ v \mapsto \Phi_t(v,\cdot) $.
	
The goal of this paper is to study the asymptotic dynamics of 
$ (u^N_t, t \geq 0) $ in a neighbourhood of the manifold $ M $ as 
$ N \to \infty $. In particular, we make precise the idea that 
	as $ N \to \infty $, $ u^N_t $ should remain close to the
manifold $M$, travelling along it  
at a speed which is a small perturbation of $ \alpha $ (the speed 
of the solution to the deterministic equation confined to the manifold), 
and experiencing small random fluctuations away from $M$.
	
The following family of norms will allow us to quantify the deviations 
of $ u^N_t $ away from the manifold $ M $.
	For $ p \geq 1 $ and $ q \in \R $, for $g:\R\to \R$, define
	\begin{equation} \label{eq:normpq}
		\| g \|_{p,q} := \left( \int_\R \abs{g(x)}^p e^{q x} dx \right)^{1/p},
	\end{equation}
	and let $ L^{p,q} $ denote the space of measurable functions $ g : \R \to \R $ such that $ \| g \|_{p,q} < \infty $, where we identify functions that differ on sets of Lebesgue measure zero.
	Most of our analysis will use $L^{p,q}$ with $ p = 2 $ and $ q = \alpha $.

	We now set out the assumptions on the initial condition $u^N_0$ 
under which we prove our main result. They make precise the idea that we start
from an initial condition which has the compact support property and which is
close to $M$.  
	These assumptions depend on positive constants which are later assumed to belong to a certain range of values.
	We take some positive $K_0$ whose value will be specified later,
and let $ \newcinit{big}, \newcinit{small} \in (0,1)$ and $\newCst{init} \geq 1$.

\begin{assumption} \label{assumpt:v0}
		Suppose that for all $N\ge 1$, $u^N_0:\R\to [0,1]$ is continuous and satisfies the following conditions:
		\begin{enumerate}[i)]
			\item \label{v0:dist} there exists $\eta \in [-K_0,K_0] $ such that
			\begin{equation*}
				\| u^N_0 - m_\eta \|_{2,\alpha} \leq N^{-\cinit{big}},
			\end{equation*}
			\item \label{v0:compact_support} $ u^N_0(x) = 0 $ for all $ x \geq (1+\cinit{small}) \log N $, $\, u^N_0(x) = 1 $ for all $ x \leq -N $,
			\item \label{v0:holder} for all $ x, y \in \R $ with $|x-y| \leq 1$,
			\begin{equation*}
				| u^N_0(x)-u^N_0(y) | \leq \left( \Cst{init} \wedge \left( N^{\cinit{small}} e^{-((1-\cinit{small})/2) x} \right) \right) |x-y|^{\cinit{big}}.
			\end{equation*}
			\item \label{v0:tail} for all $ x \in \R $,	
				\begin{equation*} 
					u_0^N(x) \leq N^{\cinit{small}} e^{-(1-\cinit{small})x} \quad \text{ and }  \quad u^N_0(x) \vee (1 - u^N_0(x)) \leq \Cst{init}(e^{-\cinit{big}x} + N^{-\cinit{big}}).
				\end{equation*}
		\end{enumerate}
	\end{assumption}
	
	To state our main result, we need one more piece of notation.
	For $ v : \R \to [0,1] $, $ v \in L^{2,\alpha} $, let $ p_{s,t}(v,x,\cdot) $ denote the fundamental solution of $ \partial_t - \partial_{xx} - \alpha \partial_x - f'(\Phi_t(v,x)) $, i.e.~for any $ u \in L^{2,\alpha} $, we let $ u_t(x) := \langle p_{s,t}(v,x,\cdot), u \rangle_\alpha $, for $ t > s $, denote the solution of
		\begin{equation} \label{def_pst}
			\left\lbrace
			\begin{aligned}
			& \partial_t u_t = \partial_{xx} u_t + \alpha \partial_x u_t + f'(\Phi_t(v,\cdot))u_t, \quad t > s, \\
			& \lim_{t \downarrow s} u_t = u.
			\end{aligned}
			\right.
		\end{equation}
		For $T>0$, let $ D([0,T],\R) $ denote the space of \cadlag real-valued functions endowed with the Skorokhod topology and the usual metric which makes it complete (see Theorem~12.2 in \cite{billingsley_convergence_1999}).

	Our main result is the following. We remind the reader that the range $\alpha \in (0,3/2)$ corresponds to the entire fully pushed regime as predicted by \cite{birzu_fluctuations_2018}.
	
	\begin{theorem} \label{thm:main_result}
		Assume that $ \alpha \in (0, 3/2) $, and let $(u^N_t,t\ge 0)$ be as in Definition~\ref{def:u}.
		For $N \geq 1$, set
		\begin{equation} \label{def_vN}
			v^N_t(x) := u^N_{Nt}(x + \alpha N t).
		\end{equation}
		Then for any $ K_0 > 0 $, $\cinit{big} \in (0,1)$ and $\Cst{init} \geq 1$, for sufficiently small $\cinit{small}\in (0,1)$ the following holds.
		Under Assumption~\ref{assumpt:v0}, for any fixed $ T > 0 $, there exists a sequence of real-valued processes $ (\xi^N_t, t \in [0,T]) $ such that
		\begin{align*}
			\sup_{t \in [0,T]} \| v^N_t - m_{\xi^N_t} \|_{2,\alpha} \cvgas{N} 0,
		\end{align*}
		in probability.
		Moreover, as $N\to \infty$, the sequence of processes $ (\xi^N_t - \xi^N_0, t \in [0,T]) $ converges in distribution in $ D([0,T],\R) $ to $ (\xi_t, t \geq 0) $ where
		\begin{align} \label{eds_xi}
			\xi_t =  - \drift t + \variance B_t,
		\end{align}
		where $ (B_t, t \geq 0) $ is standard Brownian motion, 
		and $ \drift $ and $ \variance $ are two constants given by
		\begin{align} \label{parameters_xi}
		\drift := \frac{\alpha}{4} A_1 + \frac{1}{2} A_2, && \variance^2 := A_1,
		\end{align}
		with
		\begin{equation} \label{eq:A1A2defn}
		\begin{aligned}
		& A_1 := \frac{1}{\| \partial_{x} m \|_{{2,\alpha}}^4} \int_\R \partial_{x} m(x)^2 m(x) (1-m(x)) e^{2\alpha x} dx, \\
		& A_2 := \frac{1}{\| \partial_{x} m \|_{{2,\alpha}}^2} \int_{0}^{\infty} \int_{\R^2} p_{0,t}(m,x,y)^2 f''(m(x)) \partial_{x} m(x) m(y)(1-m(y)) e^{2\alpha y + \alpha x} dy dx dt. 
		\end{aligned}
		\end{equation}
	\end{theorem}
	
	\begin{figure}
		\centering
		\includegraphics[width=0.8\textwidth]{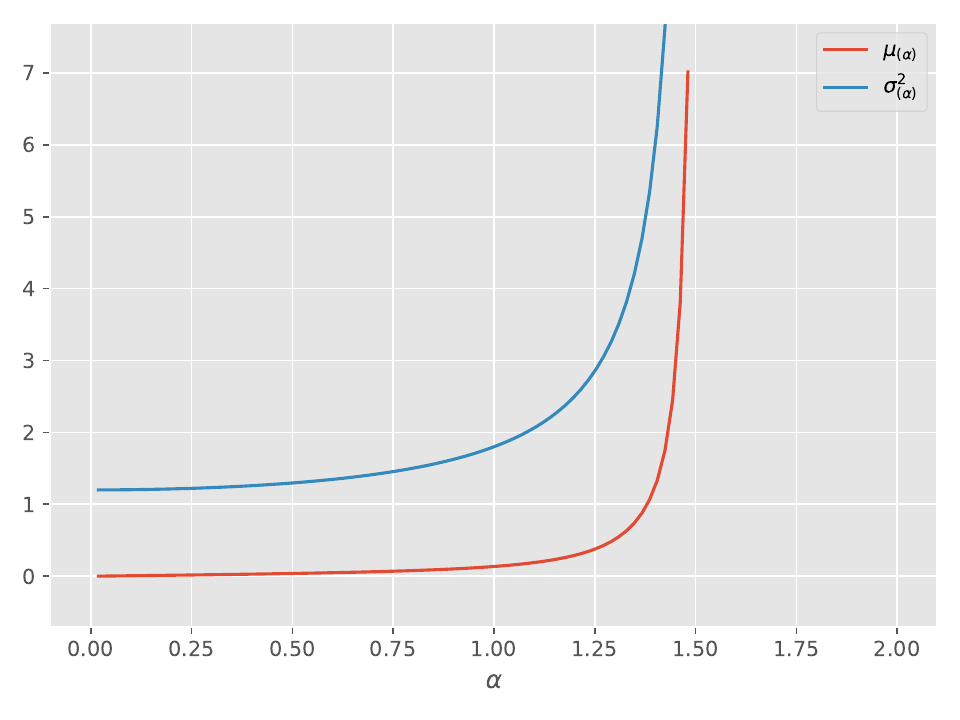}
		\caption{Values of $ \drift $ and $ \variance^2 $ as a function of $ \alpha $, for $ \alpha \in (0,3/2) $, when $f$ is given by \eqref{f}. These values were computed numerically by approximating the deterministic equation with a finite difference method on a discretised domain and solving the equation for $ p_{0,t}(m,\cdot,\cdot) $ on the resulting domain.
		When $ \alpha \downarrow 0 $, we have $ \drift \to 0 $ (as expected by symmetry, since in this case $ 1-u^N_t(-x) $ has the same dynamics as $ u^N_t $) and $ \variance \to \sigma_{(0)} > 0 $, where $ \sigma_{(0)} $ can be obtained by adapting Funaki's result in~\cite{funaki_scaling_1995}.
		The plot suggests that $ \drift > 0 $, as we conjecture. We also see that $ \drift $ and $ \variance^2 $ both appear to be increasing as functions of $ \alpha $.} \label{fig:mu_sigma2}
	\end{figure}
	
	
	Note that, since we have assumed that $\lambda_+ = 1$, $ m(x) \sim k e^{-x} $ as $ x \to \infty $, and $ |\partial_{x} m(x)| \sim k' e^{-x} $ as $ x \to \infty $ for some $k,k'>0$.
	As a result, $ \partial_{x} m \in L^{2,\alpha} $ for any $ \alpha \in (0,2) $, and $ A_1 < \infty $ for any $ \alpha \in (0, 3/2) $.
	The fact that $ |A_2| < \infty $ is less obvious, but we will show this for $ \alpha \in (0,3/2) $ in Corollary~\ref{cor:A2} below.
	
	Theorem~\ref{thm:main_result} proves a conjecture formulated in \cite{birzu_fluctuations_2018}, and covers the entire fully pushed regime.
	Figure~\ref{fig:mu_sigma2} shows numerical computations of $ \drift $ and $\variance$ for the range $\alpha \in (0,3/2)$ in the case where $f$ is given by \eqref{f}.
	These numerical computations suggest that, for this specific reaction term, $ \drift > 0 $ for any $\alpha\in (0,3/2)$, although we have not been able to prove this.
	The range $ \alpha \in (0,3/2) $ covered by Theorem~\ref{thm:main_result} is optimal, as both $\drift$ and $ \variance $ diverge as $ \alpha $ approaches $ 3/2 $.
	The behaviour of the position of the front in the range $ \alpha \in [3/2,2) $ remains fully open.
	However, it is conjectured in \cite{birzu_fluctuations_2018} that for $\alpha \in [3/2,1)$, the relevant time scale is $ N^{\frac{2-\alpha}{\alpha-1}} $. 



\subsection{Relation to existing literature}

Our proof owes a lot to~\cite{funaki_scaling_1995}, in which Funaki 
studies a stochastic partial differential equation of Ginzburg-Landau type 
for a symmetric double well potential as the temperature of the system tends
to zero, resulting in phase separation. He finds the equation of motion 
of the phase separation point. The system is tuned in such a way that the 
noise becomes small as the reaction term becomes large.
Although~\cite{funaki_scaling_1995} was really interested in the motion
of the sharp interface between phases, the proof goes via scaling and
his result can be recast as follows. Consider the equation
	\begin{equation} \label{spde_funaki}
		\partial_t u^\varepsilon = \partial_{xx} u^\varepsilon + f(u^\varepsilon) + \varepsilon^\gamma a_\varepsilon(\cdot) \dot{W},
	\end{equation}
	where $ \gamma > 5 $ and $ a_\varepsilon(x) = a(\sqrt{\varepsilon} x) $ for some smooth and compactly supported function $ a : \R \to \R $, and where $ f $ is assumed to be bistable and symmetric about $ 1/2 $, i.e.
	\begin{align*}
		f(0) = f(1/2) = f(1) = 0, && f(1-u) = -f(u).
	\end{align*}
	In particular, this implies $ \alpha = 0 $.
	Funaki showed that, under these assumptions, as $ \varepsilon \to 0 $, $ u^\varepsilon(\varepsilon^{-2\gamma-1} t, \cdot) $ approaches the standing wave profile $ m $ shifted by a random quantity $ \varepsilon^{-1/2} \xi^\varepsilon_t $, and that $ (\xi^\varepsilon_t, t \geq 0) $ converges in distribution to a limiting process $ (\xi_t, t \geq 0) $ solving
	\begin{equation*}
		d\xi_t = c_1 a(\xi_t) d B_t + c_2 a(\xi_t) a'(\xi_t) dt,
	\end{equation*}
	where $ (B_t, t \geq 0) $ is standard Brownian motion and $ c_1 $ and $ c_2 $ are two constants.
In a more recent paper, \cite{xu_interface_2025} showed that the same result holds for any $ \gamma > 0 $.

Proving results of this form for equations with
double well potentials and additive noise has been a topic of particular focus.
The paper~\cite{brassesco_brownian_1995} 
attacks a similar problem to Funaki. They specialise to $f(u)=u-u^3$
and $a_{\varepsilon}(\cdot)\equiv 1$; and instead of having the noise
vanish outside a compact set, they work on the expanding region 
$[-\varepsilon^{-1}, \varepsilon^{-1}]$ with 
Neumann boundary conditions. 
In~\cite{brassesco_interface_1998} 
this is generalised to 
the case of a non-symmetric potential, but still 
requiring that $ \int_{0}^{1} f(u) du = 0 $. Interestingly the drift of the 
limiting Brownian motion is not necessarily zero in this case, 
which raises questions of the sign of $ \mu_{(\alpha)} $ in our result
for general $ f $.
With a completely different perspective,
\cite{kruger_multiscale-analysis_2017} uses a multiscale analysis to 
decompose a stochastic bistable reaction-diffusion equation with an additive
noise into the orthogonal sum of a travelling wave moving with random 
speed and Gaussian fluctuations, extending their result 
in~\cite{kruger_front_2014} for stochastic neural field equations. 

Our equation \eqref{spde_u_N} differs from 
those considered in the works described above
in several important ways: the reaction term can correspond to any fully pushed
wave, so that the 
solutions to the underlying deterministic reaction-diffusion equation
converge to travelling waves moving with a positive speed, not to a 
stationary profile; similarly, the reaction term need not be bistable, for example $f$ as in \eqref{f} with $\alpha \in [1,3/2)$; and the 
noise is not additive, but instead depends on the solution itself. 
Concerning the first point, the effect of small noise on the asymptotic speed of travelling waves 
was considered in the more noise-sensitive pulled regime in \cite{mueller_effect_2011}, but fluctuations of the position of the front were not studied.

Of course, there are many other classes of equation that model 
the kinematics of 
phase segregation for different physical systems, and for which one
can ask about the impact of a noisy perturbation. For example,
\cite{bertini_front_2014} considers a stochastic Cahn-Hilliard equation,
and shows that under suitable conditions, 
the motion of the interface is a continuous, self-similar process
and non-Markovian. This setting is not amenable to our approach. 

In population genetics, the main interest goes beyond the speed and
shape of the
stochastic wavefront to its implications for patterns of genetic diversity. 
These will be encoded by the genealogical
relationships between individuals sampled 
from the population close to the wavefront. It has long been conjectured that
for the stochastic Fisher-KPP equation these genealogies will converge, 
after suitable
scaling, to the Bolthausen-Sznitman coalescent. This reflects the large
fluctuations in the stochastic wavefront, caused by individuals getting 
unusually far out into the tip of the wave, where their family of 
descendants will grow very quickly; when the rest of the wave `catches up'
that family will form a significant portion of the individuals in the front.
On suitable timescales, the creation of the large family will appear 
instantaneous and, as we trace backwards in time, all ancestral lineages
that are in the family will merge in a single event. What is crucial
is that the fluctuations in the front 
are much bigger than $1/N$, so we see these
events on a faster timescale than the coalescences of ancestral lineages
in the bulk which take place at rate $1/N$.

For semi-pushed waves the fluctuations in the wavefront are less dramatic, 
but still on a scale that dwarfs $1/N$. This led
\cite{birzu_fluctuations_2018} to conjecture that for 
semi-pushed waves we would see similar effects, again leading 
to a multiple merger coalescent on suitable timescales,
whereas 
for a fully pushed wave, fluctuations are order $1/N$, and we should
expect that 
genealogies converge on timescales of order $N$ to a Kingman coalescent.

Although these results are out of reach for~\eqref{spde_u_N},
there is a substantial body of literature that addresses the question
of genealogies for closely related models. 
A popular class of models considers a population of fixed size.
For example, \cite{cortines_genealogy_2018} finds an exactly solvable
model that encapsulates all three classes. 
\cite{schweinsberg_rigorous_2017-1} considers a model for a population of fixed size, and 
shows convergence of genealogies to the Bolthausen-Sznitman coalescent. His
model is closely related to so-called N-BBM, for which the hydrodynamic
limit has recently been shown by
\cite{berestycki_convergence_2025} to sit in the Fisher-KPP universality
class. They also apply their results to a more general free boundary 
problem for which they see a transition from pulled to pushed behaviours.
\cite{tourniaire_branching_2024} produced an inhomogeneous branching 
particle system that can be seen as an analytically tractable model of 
pushed fronts and in \cite{schertzer_spectral_2025} they establish the genealogy
in the fully pushed case. 

Finally we mention that
\cite{etheridge_genealogies_2022} consider a discrete analogue of the 
SPDE~\eqref{spde_u_N} for a reaction term of the form~\eqref{f} with
$\alpha<1$ and show convergence of the genealogies to a Kingman coalescent.
Note that this does not capture the whole of the fully pushed regime, but
only that part in which the growth rate in the front is negative. This
tames the fluctuations in the front. In contrast, our results in this 
paper capture the \emph{whole} of the semi-pushed regime.


	\paragraph*{Overview of the proof}
	
	The proof is largely based on the proof of Funaki's 1995 result \cite{funaki_scaling_1995}.
	The first step of the proof of both Theorem~\ref{thm:main_result} and Funaki's result is to show that, as $ N \to \infty $,
	\begin{equation*}
		\dist(v^N_t, M) \to 0,
	\end{equation*}
	locally uniformly in $ t $, in probability, where $ M $ is the stable manifold of shifts of the travelling wave profile defined in $\eqref{def:manifoldM}$
	and $\dist(v,M) = \inf_{\eta \in \R}\| v - m_\eta\|_{2,\alpha}$. This is done by controlling exponential moments of the following energy functional
	\begin{equation*}
		\mathcal{H}(u) := \int_\R \left( \frac{1}{2} |\partial_{x} u(x)|^2 + F(u(x)) \right) e^{\alpha x} dx,
	\end{equation*}
	where $ F(u) = - \int_{0}^{u} f(v) dv $.
	Because of the exponential weight $ e^{\alpha x} $ appearing in the above definition (which was absent in \cite{funaki_scaling_1995} since $ \alpha = 0 $ in his setting), we need a precise control on the exponential decay of the solution $ v^N_t(x) $ as $ x \to \infty $, which turns out to require precise estimates on the right endpoint of the support of $ v^N_t $.
	More precisely, we show that, for any $ c > 0 $ and $ T > 0 $, with high probability, $ v^N_t(x) = 0 $ for all $ x \geq (1+c) \log(N) $ and $ t \in [0,T] $. 	
	While the non-quantitative compact interface property (see Lemma~\ref{lem:compactinterface}) is well-known for stochastic reaction-diffusion equations, the precise quantitative bound above is novel in this setting.
	That the right endpoint of the support should be near $\log N$ can be guessed from the asymptotics of $m$, to which $v^N_t$ is close, and extinction probabilities for solutions to \eqref{def:u} started from small initial data. 
	We establish the uniform control of the right endpoint by showing that it has self-regulating behaviour, which we achieve using a modification of Krylov's method for interface control \cite{krylov_result_1997} and new estimates related to the integrability of solutions near their right endpoint.

	The limiting dynamics of the position of the front is obtained by 
adapting a technique introduced in the finite dimensional
setting by \cite{katzenberger_solutions_1991}. To gain some intuition, we 
briefly recall his approach. 
The idea is that we are trying to solve an equation of the form
\[
X_n(t)=X_n(0)+\int_0^t\sigma_n(X_n(s-))dZ_n(s)+\int_0^tF(X_n(s-))dA_n(s)
\]
where $Z_N$ is a well behaved semi-martingale, $F$ is a deterministic
flow with stable manifold $M$, and $A_n$ is a non-decreasing process which 
asymptotically puts infinite mass on every interval.
In our case, since our result is concerned with the solution $u^N$ at times
$Nt$, $dA_n(s)=Nds$. 
If we solve the deterministic system
\[
\phi_0(x,t)=x+\int_0^tF(\phi_0(x,s))ds,
\]
then started from a point in the domain 
of attraction of $M$, $\phi_0(x,t)$ 
converges as $t\to\infty$ to a point $\Phi(x)\in M$.
This defines a map 
$\Phi:\IR^d\to M$ which is such that $\Phi(y)=y$ for $y\in M$,
and, from an application of the chain rule, we see that
$\partial\Phi(x)F(x)=0$ for $x\in\IR^d$. Applying It\^o's
formula to $\Phi(X_n)$ gives
\[
\Phi(X_n(t))=\Phi(X_n(0))+\int_0^t\partial\Phi\sigma_ndZ_n
+\frac{1}{2}\sum_{ijkl}\int_0^t\partial_{ij}\Phi\sigma_n^{ik}
\sigma_n^{jl}d[Z_n^k,Z_n^l] +\eta_n(t),
\]
where $\eta_n$ captures corrections for the jumps (which will be 
absent in our setting). Under suitable
conditions $\eta_n\to 0$ and $d(X_n,M)\Rightarrow 0$,
so $X_n-\Phi(X_n)\Rightarrow 0$ and $X_n$ solves
\begin{equation*}
X_n(t)=X_n(0)+\int_0^t\partial\Phi\sigma_ndZ_n
+\frac{1}{2}\sum_{ijkl}\int_0^t\partial_{ij}\Phi
\sigma_n^{ik}\sigma_n^{jl}d[Z_n^k,Z_n^l]
+\epsilon_n(t)
\end{equation*}
and $\epsilon_n\Rightarrow 0$ as $n\to\infty$. 

Our proof will also involve
applying It\^o's formula to a suitable projection on the stable manifold 
$M$ and letting $ N \to \infty $.
	The infinite-dimensional setting causes some difficulties here, but we are nonetheless able to find a suitable function $ \zeta(v^N_t) \in \R $ such that (up to some well-behaved stopping time)
	\begin{equation*}
		\| v^N_t - m_{\zeta(v^N_t)} \|_{2,\alpha} \leq N^{-c},
	\end{equation*}
	for some $ c > 0 $ and
	\begin{multline*}
		\zeta(v^N_t) = \zeta(u^N_0) + \int_{[0,t] \times \R} D\zeta(v^N_s,y) \sqrt{v^N_s(y)(1-v^N_s(y))} e^{\alpha y} W(dy,ds) \\ + \frac{1}{2} \int_{0}^{t} \int_\R D^2 \zeta(v^N_s,y,y) v^N_s(y) (1-v^N_s(y)) e^{2\alpha y} dy ds.
	\end{multline*}
	The additional potentially divergent term appearing in the It\^o formula vanishes thanks to a nice property of $ \zeta $.
	To obtain \eqref{eds_xi}, we then wish to let $ N \to \infty $ in the above equation, which requires both that $ v \mapsto D\zeta(v,\cdot) $ and $ v \mapsto D^2 \zeta(v,\cdot,\cdot) $ are Lipschitz in some suitable norm ($ L^p $ norms with some exponential weight) and precise estimates on the exponential decay of $ v^N_s(y) $ as $ y \to \infty $.
	Deriving these estimates proves to be a technical challenge, in particular when dealing with the terms involving $ p_{0,t}(m,x,y) $ in $ D^2\zeta(v,\cdot,\cdot) $.
	The optimality of the range $ \alpha \in (0,3/2) $ for Theorem~\ref{thm:main_result} appears here,  
	as several quantities diverge as $\alpha$ approaches $3/2$.

	We prove Theorem~\ref{thm:main_result} in Section~\ref{sec:proof} assuming several intermediate results. 
	These intermediate results are then proved in the rest of the paper; the organization of these proofs is outlined at the end of Section~\ref{sec:proof}.
	The proof of Theorem~\ref{thm:main_result} also requires a few results on the large time behaviour of the deterministic flow $ v \mapsto \Phi_t(v,\cdot) $ in order to define the sequence of processes $ (\xi^N_t, t \geq 0) $.
	These are stated in Section~\ref{sec:coordinates}.

	\section{Coordinates along the stable manifold} \label{sec:coordinates}

	For $ v \in L^{2,\alpha} $, we define the distance between $ v $ and the manifold $M$ as
	\begin{equation*}
		\dist(v, M) = \inf_{\eta \in \R} \| v - m_\eta \|_{{2,\alpha}}
	\end{equation*}
	($ m_\eta \in L^{2,\alpha} $ for all $ \eta \in \R $).
	Note that, if the infinimum is attained for some $ \eta \in \R $, then
	\begin{equation} \label{diff_norm}
		\deriv{}{\eta} \| v - m_{\eta} \|_{2,\alpha}^2 = 2 \langle v - m_\eta, \partial_{x} m_\eta \rangle_\alpha = 0,
	\end{equation}
	where $ \langle \cdot, \cdot \rangle_{\alpha} $ denotes the scalar product in $ L^{2,\alpha} $,
	\begin{equation} \label{scalar_product}
		\langle \phi, \psi \rangle_{\alpha} := \int_{\R} \phi(x) \psi(x) e^{\alpha x} dx \quad \forall \phi, \psi \in L^{2,\alpha}.
	\end{equation}
	This motivates the following.
	
	\begin{lemma} \label{lemma:eta}
		For any $ K > 0 $, there exists $ \newBeta{eta} = \Beta{eta}(K) > 0 $ such that, for all $ v \in L^{2,\alpha} $, if there exists $ \eta_0 \in [-K,+\infty) $ with $ \| v - m_{\eta_0} \|_{2,\alpha} \leq \Beta{eta} $ then there exists a unique $ \eta = \eta(v) \in \R $ such that
		\begin{equation} \label{dist_realised}
			\| v - m_{\eta(v)} \|_{2,\alpha} = \dist(v,M).
		\end{equation}
		For such $ v $,
		\begin{equation} \label{s_orthogonal}
			\langle v - m_{\eta(v)}, \partial_{x} m_{\eta(v)} \rangle_\alpha = 0.
		\end{equation}
		Furthermore, there exists $ \newCst{eta} > 0 $ (depending on $ K $) such that, for any such $ v $ and $ \eta_0 $,
		\begin{equation} \label{bound_eta_eta_0}
			| \eta(v) - \eta_0 | \leq \Cst{eta} \| v - m_{\eta_0} \|_{2,\alpha},
		\end{equation}
		and for any $ v_1, v_2 \in L^{2,\alpha} $ such that there exists $ \eta_1, \eta_2 \in [-K,+\infty) $ with $ \| v_i - m_{\eta_i} \|_{2,\alpha} \leq \Beta{eta} $,
		\begin{equation} \label{eta_Lipschitz}
			| \eta(v_1) - \eta(v_2) | \leq \Cst{eta} \| v_1 - v_2 \|_{2,\alpha}.
		\end{equation}
	\end{lemma}
	
	We prove Lemma~\ref{lemma:eta} in Appendix~\ref{sec:fermi}.
	Let us then define a set $ \mathcal{V}_{K} \subset L^{2,\alpha} $ as
	\begin{equation*}
		\mathcal{V}_{K} := \lbrace v \in L^{2,\alpha} : \exists \eta \in [-K,K] : \| v - m_{\eta} \|_{2,\alpha} < \Beta{eta}(K) \rbrace.
	\end{equation*}
	By Lemma~\ref{lemma:eta}, for all $ v \in \mathcal{V}_{K} $, there exists exactly one $ \eta = \eta(v) \in \R $ satisfying \eqref{dist_realised} and $ | \eta(v) | < K + \Cst{eta} \Beta{eta} $.
	It is clear that $ \mathcal{V}_{K}$ is an open subset of $ L^{2,\alpha} $, as it is a union of open balls.
	Also set
	\begin{equation*}
		\mathcal{V} := \bigcup_{K > 1} \mathcal{V}_{K}.
	\end{equation*}
	For any $ v \in \mathcal{V} $, we set 
	\begin{equation} \label{eq:Fermidefn}	
		s(v) := v - m_{\eta(v)}.
	\end{equation}
	Then $ (\eta(v), s(v)) $ are called the \textit{Fermi coordinates} of $ v $.
	
	For $\lambda \geq 0$, we introduce a weighted $L^\infty$-norm defined by
	\begin{equation} \label{eq:weightedsupdefn}
	\| u \|_{\infty,\lambda} := \sup_{x \in \R} \frac{|u(x)|}{1 \wedge e^{-\lambda x}}.
	\end{equation}
	For $ K > 0 $, $ \beta \in (0,\Beta{eta}(K)) $ and $ \epsilon > 0 $, define
	\begin{align} \label{eq:VbetaKep}
		\mathcal{V}_{\beta,K,\epsilon}^{(\lambda)} := \lbrace v \in L^{2,\alpha} : 0 \leq v \leq 1, v \in \mathcal{V}_{K}, \dist(v,M) < \beta  \text{ and } \| v - m_{\eta(v)} \|_{\infty,\lambda} < \epsilon \rbrace.
	\end{align}
	Note that for $v\in \mathcal{V}_{\beta,K,\epsilon}^{(\lambda)}$, by Lemma~\ref{lemma:eta} we have $\dist(v,M)=\| v - m_{\eta(v)} \|_{2,\alpha}$.
	For $ n \in \N $, also let $ H^{n,\alpha} $ be the weighted Sobolev space equipped with the norm
	\begin{align} \label{def:Hnorm_n}
		\| g \|_{H^{n,\alpha}}^2 = \sum_{k=0}^{n} \| \partial_{x}^k g \|_{2,\alpha}^2.
	\end{align}
	
	As stated in Section~\ref{sec:statement}, Fife and McLeod proved in \cite{fife_approach_1977} that, for $\alpha \in (0,1)$, for suitable initial conditions $v$, $ \Phi_t(v,\cdot) \to m_{\zeta} $ as $ t \to \infty $ for some $ \zeta \in \R $, with convergence in the uniform topology (moreover, they prove that this convergence is exponential).
	Rothe \cite{rothe_convergence_1981} extended this to $\alpha \in (1,2)$, assuming that $v$ decays sufficiently fast at $+\infty$, and proved that the convergence takes place in the topology induced by $\| \cdot \|_{\infty,\lambda}$ for suitable values of $\lambda$.
	In our analysis, we shall need to establish that the convergence also takes place in $ L^{2,\alpha} $, and in $ H^{1,\alpha} $, at exponential speed, and uniformly for $v$ in some neighbourhood of the manifold $ M $.
	This is the content of the next theorem, which is proved in Section~\ref{sec:det_flow}.
	
	Note that, since $ f $ is defined on $ [0,1] $, $ \Phi_t(v,\cdot) $ is a priori defined only for $ v : \R \to [0,1] $.
	We can nonetheless extend $ \Phi_t $ to $ L^{2,\alpha} $ by considering a smooth extension of $ f $ defined on $ \R $ and such that $ f $, $ f' $ and $ f'' $ are all bounded and continuous.
	
	\begin{theorem} \label{thm:zeta} 
		For any $ \alpha \in (0,2) $, any $\lambda \in (\alpha-1, 1)$ with $\lambda \geq 0$ and any $ K > 0 $, there exist $ \newBeta{zeta} \in (0,\Beta{eta}) $, $ \newEpsilon{zeta} > 0 $, $ \newCst{zeta} > 0 $ and $ \genericc > 0 $ such that, for any $ v \in \mathcal{V}_{\Beta{zeta}, K, \Epsilon{zeta}}^{(\lambda)} $, $ \Phi_t(v,\cdot) \in \mathcal{V}_{K+1} $ for all $ t \geq 0 $ and there exists $ \zeta = \zeta(v) \in \R $ such that, for all $ t \geq 0 $,
		\begin{align*}
		\| \Phi_t(v,\cdot) - m_{\zeta(v)} \|_{{2,\alpha}} \leq \Cst{zeta} \, \dist(v,M) e^{-\genericc t}.
		\end{align*}
		Furthermore, the limit also holds in $ H^{1,\alpha} $; more precisely, for all $ v \in \mathcal{V}_{\Beta{zeta}, K, \Epsilon{zeta}}^{(\lambda)} $ and $ t \geq 1 $,
		\begin{align*}
		\| \Phi_t(v,\cdot) - m_{\zeta(v)} \|_{H^{1,\alpha}} \leq \Cst{zeta} \, \dist(v,M) e^{-\genericc t}.
		\end{align*}
		In addition, for all $ t \geq 0 $,
		\begin{align} \label{diff_eta_zeta}
		\abs{\zeta(v) - \eta(\Phi_t(v,\cdot))} \leq \Cst{zeta} \, \dist(v,M) e^{-\genericc t}. 
		\end{align}
		In particular, there exists $ \newCst{sup_zeta} > 0 $ (depending on $ K $) such that $ \abs{\zeta(v)} \leq \Cst{sup_zeta} $ for all $ v \in \mathcal{V}_{\Beta{zeta}, K, \Epsilon{zeta}}^{(\lambda)} $.
	\end{theorem}

	We call $ \zeta(v) $ the \textit{Katzenberger coordinate} of $ v $, after \cite{katzenberger_solutions_1991}, where it was introduced in a finite-dimensional setting (see also \cite{funaki_scaling_1995} where it was generalised to an infinite-dimensional process).
	Here and in most of the remainder of the article, when there is no ambiguity we omit the dependence of the various quantities $ \beta_i $, $ \varepsilon_i $, $ C_i $, etc.~on the constant $ K $ for notational convenience.
	In all our statements, the constant $ K $ is kept fixed but arbitrary, only to be chosen at the beginning of the proof of Theorem~\ref{thm:main_result}.
	
	The idea of the proof of Theorem~\ref{thm:main_result} (which we will implement in Section~\ref{sec:proof} below) is to define a process $ \xi^N_t $ which is equal to $ \zeta(v^N_t) $, at least for $ t\le \tau_N $, where $\tau_N$ is a stopping time such that $ v^N_t $ remains in a suitable neighbourhood of $ M $ until time $\tau_N$.
	Then, using It\^o's formula for infinite-dimensional stochastic processes, we are able to identify the infinitesimal drift and variance of this process, and show that they converge to $ -\drift $ and $ \variance^2 $ respectively.
	In order to achieve this, we need to introduce the Fr\'echet derivatives of first and second order of $ \zeta $.
	Let us first recall the definition of Fréchet derivatives.
	
	\begin{definition} \label{def:Frechet_derivative}
		Let $ V $, $ W $ be normed vector spaces and $ U \subset V $ be an open subset of $ V $.
		A function $ g : U \to W $ is Fréchet differentiable at $ v \in U $ if there exists a bounded linear operator $ A : V \to W $ such that for $ h \in V $ with $ v + h \in U $,
		\begin{equation*}
			\| g(v + h) - g(v) - A h \|_{W} = \littleO{\| h \|_{V}} \quad \text{ as } \| h \|_V \to 0.
		\end{equation*}
		We then set $ Dg(v) = A $.
		Moreover, $ g $ is twice Fréchet differentiable at $ v \in U $ if it is Fréchet differentiable at $ v $ and if there exists a continuous bilinear map $ B: V \times V \to W $ such that for $ h \in V $ with $ v + h \in U $,
		\begin{equation*}
			\left\| g(v + h) - g(v) - Dg(v) h - \frac{1}{2} B(h,h) \right\|_W = \littleO{ \|h\|_{V}^2 } \quad \text{ as } \| h \|_V \to 0.
		\end{equation*}
		In that case we set $ D^2g(v) = B $.
	\end{definition}
	
	Note that, when $ V $ is a Hilbert space endowed with the inner product $ \langle \cdot, \cdot \rangle_V $ and $ W = \R $, by the Riesz representation theorem, $ Dg(v) $ can be identified with an element of $ V $, which we also denote by $ Dg(v) $ (or $ y \mapsto Dg(v,y) $), such that
	\begin{equation} \label{eq:Dgvhalpha}
		Dg(v) h = \langle Dg(v), h \rangle_V.
	\end{equation}
	Also by the Riesz representation theorem, there exists a linear operator $ L : V \to V $ such that $ D^2g(v)(h_1, h_2) = \langle L h_1, h_2 \rangle_{V} $.
	When $ V = L^2(\R, \mu) $ for some measure $ \mu $ and $ L $ is Hilbert-Schmidt, there exists a function in $ L^{2}(\R^2, \mu \otimes \mu) $, which in that case we also denote by $ (x,y) \mapsto D^2g(v,x,y) $, such that
	\begin{equation} \label{HS_derivative}
		D^2g(v)(h_1, h_2) = \int_{\R^2} D^2g(v,x,y) h_1(x) h_2(y) \mu(dx) \mu(dy),
	\end{equation}
	see \citep[Theorem~6.12]{brezis_functional_2011}.
	
	One difficulty here is that $ \mathcal{V}_{\Beta{zeta},K,\Epsilon{zeta}} $ is not an open subset of some suitable Hilbert space.
	To circumvent this, we instead apply Ito's formula to the function $ v \mapsto \eta(\Phi_t(v,\cdot)) $, which by Proposition~\ref{prop:Deta_t} below is defined and Fréchet differentiable on an open subset of a suitable Hilbert space.
	By letting $ t \to \infty $, we will be able to recover what we would have obtained by directly applying Ito's formula to $ v \mapsto \zeta(v) $.
	In the following, we use the shorthand
	\begin{equation} \label{def:eta_t}
		\eta_t(v) := \eta(\Phi_t(v,\cdot)).
	\end{equation}
	Set
	\begin{equation} \label{def:w}
		w(x) := \1{x \geq 0} e^{\alpha x} + \1{x < 0} e^{\frac{\alpha}{3} x},
	\end{equation}
	and let $ H $ be the Hilbert space defined as the weighted $ L^2 $ space on $ \R $ with weight function $ w $, equipped with the scalar product
	\begin{equation*}
		\langle u, v \rangle_H := \int_{\R} u(x) v(x) w(x) dx.
	\end{equation*}
	Note that $ H $ is trivially embedded in $ L^{2,\alpha} $ and hence $ \Phi_t $ is defined on $ H $.
	Also let $ \| \cdot \|_H $ denote the associated norm on $ H $.

	\begin{lemma} \label{lemma:Phi_t_Lipschitz}
		For any $ T \geq 0 $, $ \Phi_t : L^{2,\alpha} \to L^{2,\alpha} $ is uniformly Lipschitz continuous for all $ t \in [0,T] $, i.e.~there exists $ C > 0 $ such that, for all $ v, v' \in L^{2,\alpha} $ and all $ t \in [0,T] $,
		\begin{equation*}
			\| \Phi_t(v,\cdot) - \Phi_t(v',\cdot) \|_{2,\alpha} \leq C \| v - v' \|_{2,\alpha}.
		\end{equation*}
		By the embedding of $ H $ into $ L^{2,\alpha} $, the same holds for $ \Phi_t : H \to L^{2,\alpha} $.
	\end{lemma}

	For $ t > 0 $, set
	\begin{equation} \label{def:Ut}
		U_t := \lbrace v \in H : \Phi_t(v,\cdot) \in \mathcal{V} \rbrace.
	\end{equation}
	Then, by Lemma~\ref{lemma:Phi_t_Lipschitz}, $ U_t $ is an open subset of $ H $ for any $ t > 0 $.
	We then state the following.
	
	\begin{proposition} \label{prop:Deta_t}
		For any $ t > 0 $, $ \eta_t : U_t \subset H \to \R $ is twice continuously Fréchet differentiable, and its derivatives have the representation
		\begin{equation} \label{Deta_t_scalar}
			D \eta_t(v) h = \langle D\eta_t(v), h \rangle_\alpha, \quad \forall h \in H,
		\end{equation}
		and
		\begin{equation} \label{D2eta_t_scalar}
			D^2 \eta_t(v) (h_1, h_2) = \int_{\R^2} D^2 \eta_t(v,x,y) h_1(x) h_2(y) e^{\alpha(x + y)} dx dy, \qquad \forall h_1, h_2 \in H,
		\end{equation}
		for some functions $ D\eta_t(v,\cdot) \in L^{2,\alpha} $ and $ D^2\eta_t(v,\cdot,\cdot) \in L^{2,\alpha}(\R^2) $.
		Moreover, for any $ K > 0 $, $ \eta_t $, $ D \eta_t $ and $ D^2 \eta_t $ are Lipschitz continuous on $ \Phi_t^{-1}(\mathcal{V}_K) $.
	\end{proposition}
	
	Let us now introduce the fractional Sobolev spaces $ H^{\gamma,\alpha} $, $ \gamma \in \R \setminus \N $, following \citep[Section~5.8.9]{evans_partial_2010}.
	Let $ \mathcal F : L^2 \to L^2 $ denote the Fourier transform, i.e.~for $\phi\in L^2$,
	\begin{align*}
		\mathcal{F} \phi(\xi) =  \int_\R e^{-i \xi x} \phi(x) dx,
	\end{align*}
	and let $ \mathcal{J} : L^{2,\alpha} \to L^2 $ denote the multiplication by $ x \mapsto e^{\frac{\alpha}{2}x} $, \textit{i.e.} $ \mathcal{J}\phi(x) = \phi(x) e^{\frac{\alpha}{2}x} $.
	For $ \phi \in L^{2,\alpha} $ and $ \gamma \in \R_+ $, let us introduce the notation
	\begin{equation} \label{def:Hnorm_gamma}
		\Hnorm{\phi}{\gamma} = \| (1 + |\cdot|^\gamma) \mathcal{F} \mathcal{J} \phi \|_{L^2},
	\end{equation}
	and let $ H^{\gamma,\alpha} $ be the subset of functions $ \phi \in L^{2,\alpha} $ such that $ \| \phi \|_{\tilde{H}^{\gamma,\alpha}} < \infty $.
	For $ \gamma \in \N $, this norm is equivalent to the one defined in \eqref{def:Hnorm_n}, since
	\begin{equation} \label{deriv_fourier}
		(-i\xi)^n \mathcal{F} \mathcal{J} \phi(\xi) = \mathcal{F} \left( \partial_{x}^n \mathcal{J} \phi \right)(\xi),
	\end{equation}
	and by Leibniz's rule,
	\begin{equation*}
		\partial_{x}^n (\mathcal{J} \phi)(x) = \sum_{i=0}^{n} \binom{n}{i} \left(\frac{\alpha}{2}\right)^{n-i} \partial_{x}^i \phi(x) e^{\frac{\alpha}{2} x}.
	\end{equation*}
	The following lemma provides bounds and continuity estimates on $ D \eta_t $ and $ D^2 \eta_t $, which will be used in the proof of Theorem~\ref{thm:main_result}.
	
	\begin{lemma} \label{lemma:bounds_Detat_fixed_time}
		For any fixed $ t \geq 1 $, $ K > 0 $ and $ \gamma \in [0,2) $, there exists $ C > 0 $ such that, for all $ v \in \Phi_t^{-1}(\mathcal{V}_K) $,
		\begin{equation} \label{Detat_Hgamma}
			\Hnorm{D \eta_t(v)}{\gamma} \leq C,
		\end{equation}
		for any $ v, v' \in \Phi_t^{-1}(\mathcal{V}_K) $,
		\begin{equation} \label{diff_Detat_Hgamma}
			\Hnorm{D\eta_t(v) - D\eta_t(v')}{\gamma} \leq C \| v- v' \|_{2,\alpha}.
		\end{equation}
		Moreover, for any $ t \geq 1 $, $ K > 0 $, $ p \geq 1 $ and $ q \in \R $ with $ |q| < p $, there exists $ C > 0 $ such that, for all  $ v \in \Phi_t^{-1}(\mathcal{V}_K) $,
		\begin{equation} \label{Deta_t_Lpq}
			\| D \eta_t(v) \|_{p,q} \leq C,
		\end{equation}
		and, if in addition $  \abs{\frac{q}{p} - \frac{\alpha}{2}} < 1 $, there exists $ C > 0 $ such that, for all $ v, v' \in \Phi_t^{-1}(\mathcal{V}_K) $, 
		\begin{equation} \label{diff_Detat_Lpq}
			\| D \eta_t(v) - D\eta_t(v') \|_{p,q} \leq C \| v - v' \|_{2,\alpha}.
		\end{equation}
		In addition, for any $ t \geq 1 $, $ K > 0 $ and $ q \in \R $ with $ |q| < 2 $ and $ |q-\alpha| < 1 $, there exists $ C > 0 $ such that, for all $ v \in \Phi_t^{-1}(\mathcal{V}_K) $,
		\begin{equation} \label{D2eta_t_L1q}
			\int_\R | D^2 \eta_t(v,x,x) | e^{qx} dx \leq C,
		\end{equation}
		and, for $ q \in \R $ with $ |q|<2 $, $ |q-\alpha| < 1 $ and $ \abs{q - \frac{3 \alpha}{2}} < 1 $, there exists $ C > 0 $ such that for any $ v, v' \in \Phi_t^{-1}(\mathcal{V}_K) $,
		\begin{equation} \label{diff_D2eta_t_L1q}
			\int_\R | D^2 \eta_t(v,x,x) - D^2 \eta_t(v',x,x) | e^{qx} dx \leq C \| v- v' \|_{2,\alpha}.
		\end{equation}
		For any $ t > 0 $, $ q \in \R $ with $ |q-\alpha| < 1 $ and $ K > 0 $, there exists $ C > 0 $ such that, for all $ v \in \Phi_t^{-1}(\mathcal{V}_K) $,
		\begin{multline} \label{D2eta_t_L1q_continuity}
			\int_{[-1,1]^2 \times \R} | D^2 \eta_t(v,y + \delta z_1, y + \delta z_2) - D^2 \eta_t(v,y,y) | \rho(z_1) \rho(z_2) e^{q y} dz_1 dz_2 dy \leq C \delta | \log(\delta) |.
		\end{multline}
	\end{lemma}

	The following proposition shows that, for $v \in \mathcal{V}_{\Beta{zeta},K,\Epsilon{zeta}}^{(\lambda)} $, $ D \eta_t $ and $ D^2 \eta_t $ converge to some limits $ D \zeta $ and $ D^2 \zeta $ as $ t \to \infty $, and that these limits satisfy suitable bounds and continuity estimates.
	Note that, by Theorem~\ref{thm:zeta}, for any $ v \in \mathcal{V}_{\Beta{zeta},K,\Epsilon{zeta}}^{(\lambda)} $, $ v \in U_t $ for all $ t \geq 0 $.
	
	\begin{proposition} \label{prop:Dzeta}
		For any $ K > 0 $, there exist $ D \zeta : \mathcal{V}_{\Beta{zeta},K,\Epsilon{zeta}}^{(\lambda)} \to L^{2,\alpha} $ and $ D^2 \zeta : \mathcal{V}_{\Beta{zeta},K,\Epsilon{zeta}}^{(\lambda)} \to L^{1,\alpha} $ such that the following holds.
		For any $ K > 0 $, $ \gamma \in [0,2) $, $ p \geq 1 $ and $ q \in \R $ with $ |q| < p $, there exist constants $ C > 0 $, $ \genericc > 0 $ and $ \delta > 0 $ such that, for all $ v \in \mathcal{V}_{\Beta{zeta},K,\Epsilon{zeta}}^{(\lambda)} $ and all $ t \geq 1 $
		\begin{align}
			\Hnorm{D \eta_t(v) - D\zeta(v)}{\gamma} &\leq C e^{-\genericc t}, \label{Detat-Dzeta_Hgamma} \\
			\| D \eta_t(v) - D\zeta(v) \|_{p,q} &\leq C\, e^{-\genericc t}.  \label{Detat-Dzeta_Lpq}
		\end{align}
		Furthermore,
		\begin{align} 
			\| D \zeta(v) \|_{p,q} &\leq C, \label{Dzeta_Lpq}\\
			\Hnorm{D \zeta(v)}{\gamma} &\leq C, \label{Dzeta_Hgamma}
		\end{align}
		and
		\begin{equation} \label{Dzeta_v-m_Lpq}
			\| D \zeta(v) - D \zeta(m_{\zeta(v)}) \|_{p,q} \leq C \, \dist(v,M)^\delta.
		\end{equation}
		In addition, for $ p \geq 1 $ and $ q \in \R $ satisfying
		\begin{equation} \label{condition_p_q_D2zeta}
			\abs{\frac{q}{p} - \alpha} < 1 \wedge (2-\alpha),
		\end{equation}
		there exist constants $ C > 0 $, $ \genericc > 0 $ and $ \delta > 0 $ such that, for all $ v \in \mathcal{V}_{\Beta{zeta},K,\Epsilon{zeta}}^{(\lambda)} $ and all $ t \geq 1 $,
		\begin{equation}
			\left( \int_\R | D^2 \eta_t(v,y,y) - D^2\zeta(v,y) |^p e^{q y} dy \right)^{1/2} \leq C \, e^{-\genericc t}, \label{D2etat-D2zeta}
		\end{equation}
		and
		\begin{align} 
			\| D^2 \zeta(v) \|_{p,q} &\leq C, \label{D2zeta_Lpq}\\
			\| D^2 \zeta(v) - D^2 \zeta(m_{\zeta(v)}) \|_{p,q} &\leq C \, \dist(v,M)^\delta.  \label{D2zeta_v-m}
		\end{align}
	\end{proposition}
	
	Note that we use the notation $ D\zeta $ and $ D^2\zeta $ even though $ \zeta $ is not a priori defined on an open subset of $ L^{2,\alpha} $ (or $ H $).
	We shall be careful in what follows to not treat these as proper Fréchet derivatives in the sense of Definition~\ref{def:Frechet_derivative}.
	We prove Proposition~\ref{prop:Dzeta} in Subsection~\ref{subsec:Dzeta}, along with the following, which gives the expressions for $ D\zeta $ and $ D^2\zeta $ when $ v \in M $.
	
	\begin{lemma} \label{lemma:Dzeta_m}
		For all $ \eta \in \R $,
		\begin{equation} \label{Dzeta_m}
		D\zeta(m_{\eta}, y) = - \frac{\partial_x m_\eta (y)}{\| \partial_x m_\eta \|_{2,\alpha}^2},
		\end{equation}
		and
		\begin{multline} \label{D2zeta_m}
		D^2\zeta (m_\eta,y) = - \frac{\alpha}{2} \frac{\partial_x m_\eta(y)^2}{\| \partial_x m_\eta \|_{2,\alpha}^4} \\ - \frac{1}{\| \partial_x m_\eta \|_{2,\alpha}^2} \int_{0}^{\infty} \int_\R p_{0,t}(m_\eta,x,y)^2 f''(m_\eta(x)) \partial_x m_\eta (x) e^{\alpha x} dx dt.
		\end{multline}
	\end{lemma}
%
%
	Recall the definition of $A_1$ and $A_2$ in~\eqref{eq:A1A2defn}.
	As mentioned after the statement of Theorem~\ref{thm:main_result}, we can now show that $|A_2|<\infty$ for $\alpha \in (0,3/2)$.
	
	\begin{cor} \label{cor:A2}
		For $\alpha \in (0,3/2)$, $A_1\in (0,\infty)$ and $|A_2|<\infty$.
	\end{cor} 
	
	\begin{proof}
		As noted after the statement of Theorem~\ref{thm:main_result}, we have $A_1\in (0,\infty)$ for $\alpha \in (0,3/2)$.
		By Lemma~\ref{lemma:Dzeta_m} we have
		\begin{equation} \label{A2_finite}
			A_2 = - \int_\R D^2 \zeta(m, y) m(y) (1-m(y)) e^{2\alpha y} dy - \frac{\alpha}{2} A_1.
		\end{equation}
		Then, by the Cauchy-Schwarz inequality,
		\begin{equation*}
			\abs{ \int_\R D^2 \zeta(m, y) m(y) (1-m(y)) e^{2\alpha y} dy } \leq \| D^2 \zeta(m) \|_{2,\frac{8\alpha}{3}} \| m(1-m) \|_{2,\frac{4\alpha}{3}}.
		\end{equation*}
		The second factor on the right-hand side is finite for $ \alpha \in (0,3/2) $ while the first factor is finite thanks to \eqref{D2zeta_Lpq} if
		\begin{equation*}
			\abs{ \frac{4\alpha}{3} - \alpha } < 1 \wedge (2-\alpha),
		\end{equation*}
		which is equivalent to $ \alpha \in (0,3/2) $.
	\end{proof}	
	
	Since $ \zeta(\Phi_t(v,\cdot)) = \zeta(v) $ for all $ t \geq 0 $, differentiating with respect to $ t $, we obtain that, for any $ t > 0 $,
	\begin{equation} \label{cancellation_t>0}
		\langle D \zeta(\Phi_t(v,\cdot), \cdot), \partial_{xx} \Phi_t(v,\cdot) + f(\Phi_t(v,\cdot)) + \alpha \partial_x \Phi_t(v,\cdot) \rangle_{\alpha} = 0.
	\end{equation}
	We would like to extend this equality to $ t = 0 $, but the scalar product on the left-hand side may not be well defined for $ t = 0 $ if $ v $ is not differentiable.
	The following lemma, which is proved in Subsection~\ref{subsec:Dzeta}, shows that one may let $ t $ tend to zero in \eqref{cancellation_t>0} provided $ v \in H^{\gamma,\alpha} $ for some $ \gamma > 0 $.
	
	\begin{lemma} \label{lemma:katzenberger}
		For any $K>0$ and $ \gamma > 0 $, for all $ v \in \mathcal{V}_{\Beta{zeta},K,\Epsilon{zeta}}^{(\lambda)} \cap H^{\gamma,\alpha} $,
		\begin{align} \label{katzenberger_cancellation}
			\langle D \zeta (v,\cdot) , \partial_{xx} v  + \alpha \partial_x v + f(v) \rangle_{{\alpha}} = 0.
		\end{align}
	\end{lemma}
	
	\section{Proof of the main result} \label{sec:proof}

	The main idea in the proof of Theorem~\ref{thm:main_result} is to show that $ (\zeta(v^N_t), t \geq 0) $ converges in distribution to $ (\xi_t, t \geq 0) $, satisfying \eqref{eds_xi}.
	We then use \eqref{diff_eta_zeta} to compare $ \eta(v^N_t) $ and $ \zeta(v^N_t) $.
	One issue is that $ \zeta(v^N_t) $ might not be defined if $ v^N_t $ strays too far from the stable manifold.
	We thus need to make sure that, with probability going to 1 as $ N \to \infty $, $ v^N $ does not leave a suitable neighbourhood of $ M $ before time $ T $.

	Let $ K > K_0 $ be fixed (it will be chosen at the beginning of the proof of Theorem~\ref{thm:main_result}), take $ \Beta{zeta} $ and $ \Epsilon{zeta} $ such that the statement of Theorem~\ref{thm:zeta} is satisfied, and let $\newBeta{proof} \in (0,\Beta{zeta})$ and $\newEpsilon{proof} \in (0, \Epsilon{zeta})$ be arbitrary (they will also be chosen at the beginning of the proof of Theorem~\ref{thm:main_result}). 
	We also fix $\cinit{big} \in (0,1) $. Eventually, we will fix constants $c_i >0$, $i=1,\dots,10$. 
	For the time being, we take constants $\cexp{KN}, \cexp{deltaN}, \cexp{rN}, \cexp{thetaN} \in (0,\cinit{big})$, whose value will be fixed below, and define the following sequences:
	\begin{align} \label{def:kappa_delta_r_theta}
		\kappa_N := N^{1-\cexp{KN}}, && \delta_N := N^{-\cexp{deltaN}}, && r_N := N^{-\cexp{rN}}, && \vartheta_N := N^{-\cexp{thetaN}} \wedge \Beta{proof}.
	\end{align}
	We also fix a constant $\newggamma{weight} \in [0,1/2)$ satisfying
	\begin{equation} \label{eq:weightdefn}
	\ggamma{weight} = 0 \text{ if } \alpha \in (0,1), \quad \ggamma{weight} \in (\alpha - 1, 1/2) \text{ if } \alpha \in [1,3/2).
	\end{equation}
	
	We then introduce a smooth approximation of $ v^N_t $ as follows.
	Let $ \rho : \R \to \R $ be a smooth, non-negative and symmetric function such that $ \rho(x) = 0 $ for all $ \abs{x} \geq 1 $ and $ \int_\R \rho(x) dx = 1 $.
	For $ \delta > 0 $, set
	\begin{align} \label{eq:rhodeltadefn}
	\rho^\delta(x) := \frac{1}{\delta} \rho \left(\frac{x}{\delta} \right)\quad \text{for }x\in \R
	\end{align}
	and
	\begin{align} \label{def_vdelta}
	v^{\delta,N}_t := \rho^\delta \ast v^N_t\quad \text{for }t\ge 0.
	\end{align}
	Also for $t\ge 0$ and $x\in \R$, define
	\begin{align} \label{eq:Ddefn}
		D^{\delta,N}(t,x) := \int_\R \rho^\delta(x-y) \abs{v^N_t(y) - v^N_t(x)} dy.
	\end{align}

	We then define the following stopping times:
	\begin{subequations}
		\begin{align}
		\newsigmaN{dist} &:= \inf \lbrace t \geq 0 : \| v^N_t - m_{\eta(v^N_t)} \|_{2,\alpha} > \vartheta_N \rbrace, \label{def_sigma_dist} \\
		\newsigmaN{eta} &:= \inf \lbrace t \geq 0 : |\eta(v^N_t)| > K \rbrace, \label{def_sigma_eta} \\
		\newsigmaN{supnorm} &:= \inf \lbrace t \geq 0 : \| s(v^{N}_t) \|_{\infty,\ggamma{weight}} > \Epsilon{proof} \rbrace, \label{def_sigma_sup}   \\
		\newsigmaN{vdelta} &:= \inf \left\lbrace t \geq 0 : \| 	D^{\delta_N,N}(t,\cdot) \|_{{2,\alpha}} > r_N \text{ or } \| v^{\delta_N,N}_t - v^N_t \|_\infty > r_N \right\rbrace, \label{def_sigma_vdelta} \\
		\newsigmaN{tail} &:= \inf \left\lbrace t \geq 0 : \int_\R v^N_t(x) e^{\alpha x} dx  > \kappa_N \right\rbrace \label{def_sigma_tail}, \\
		\newsigmaN{tail_right} &:= \inf\{t\geq 0: \exists y\in \R \text{ s.t. }v^N_t(y)\geq N^{\cexp{tail_right_N}}e^{-(1-\cexp{tail_right_exp}) y}\}, \label{eq:sigmatailrightdefn}
		\end{align}
	\end{subequations}
	where $ \cexp{tail_right_N} $ and $ \cexp{tail_right_exp} $ are constants in $ (0,1) $ whose values will be fixed below.
	We will need one further stopping time, denoted $\newsigmaN{interface}$, which controls the position of the right interface, i.e. $R_t$ from Lemma~\ref{lem:compactinterface}.
	The definition of $\sigmaN{interface}$ is somewhat more involved, and we postpone it to Section~\ref{sec:interface} (see \eqref{eq:sigmainterface}), where we analyze the right interface.
	The only reason we mention $\sigmaN{interface}$ here is because it is necessary for the statement of Proposition~\ref{prop:holder}.
	
	Note that $ \eta(v^N_t) $ might not be defined up to $ \sigmaN{dist} $ or $ \sigmaN{eta} $, but, by Lemma~\ref{lemma:eta}, it certainly is defined up to $ \sigmaN{dist} \wedge \sigmaN{eta} $, and this will suffice for what follows.
	Moreover, for $ t < \sigmaN{dist} \wedge \sigmaN{eta} $, $ \dist(v^N_t,M) = \| v - m_{\eta(v^N_t)} \|_{2,\alpha} \leq \vartheta_N $.
	We further set
	\begin{align} \label{def_sigma_N}
		\sigma_N := \inf \lbrace \sigma_{i,N}: 1 \leq i \leq \thesigmaN \rbrace,
	\end{align}
	and, for $ 1 \leq i \leq \thesigmaN $,
	\begin{align} \label{def_hat_sigma}
		\hat{\sigma}_{i, N} = \inf \lbrace \sigma_{j,N}: 1 \leq j \leq \thesigmaN, j \neq i \rbrace.
	\end{align}
	Our assumptions on $u_0^N$ give us the following result. Some of the proof is postponed to later sections.
	\begin{lemma} \label{lem:tauNpos}
	Suppose $K>K_0$ 
	and suppose that $u^N_0$ satisfies Assumption~\ref{assumpt:v0}. For $N$ sufficiently large, $\sigma_{i,N}>0$ almost surely for all $1 \leq i \leq \thesigmaN$. In particular, $\sigma_N >0$ almost surely.
	\end{lemma}
	\begin{proof}
		Since $\cexp{thetaN}<\cinit{big}$, for $N$ sufficiently large, by Assumption~\ref{assumpt:v0}.\ref{v0:dist} we have $\sigmaN{dist}>0$. 
		Since $K>K_0$, $\sigmaN{eta}>0$ for sufficiently large $N$ also by Assumption~\ref{assumpt:v0}.\ref{v0:dist}.
		For the proof that $\sigmaN{interface}>0$ for large enough $N$, see Lemma~\ref{lem:R0r0}.
		For the proof that $\sigmaN{tail_right} > 0$, see Lemma~\ref{lem:stoptimepos}. 
		This is followed immediately by the proof that $\sigmaN{supnorm} > 0$, $\sigmaN{vdelta} >0$ and $\sigmaN{tail} > 0$.
	\end{proof}
	We remark that the above implies that, for sufficiently large $N$,
	\begin{equation} \label{eq:u0V}
	u_0^N \in \mathcal{V}_{\Beta{zeta},K,\Epsilon{zeta}}^{(\ggamma{weight})},
	\end{equation}
	as can be verified from the definition of the latter set in \eqref{eq:VbetaKep}.	
	
	In order to prove Theorem~\ref{thm:main_result}, we will need two intermediate results.
	
		\begin{proposition} \label{prop:holder}
		There exists $\gamma>0$ such that under the assumptions of Theorem~\ref{thm:main_result}, for any $ T > 0 $ and $ N \geq 1 $, 
		there exists a random variable $U_{N,T} \in (0,\infty)$ such that 
		\begin{equation} \label{bound_Hgamma_vN}
			\sup_{t \in [0,\sigma_N \wedge T]} \Hnorm{v^N_t}{\gamma} \leq U_{N,T},
		\end{equation}
		and $\E{U_{N,T}} < \infty$.
		\end{proposition}
		
	
	
	\begin{proposition} \label{prop:control_sigma}
		For any fixed $ T > 0 $, $ K > K_0 $, $ \ggamma{weight} $ (satisfying \eqref{eq:weightdefn}) and $ \cinit{big} \in (0,1) $, there exist $ c_i \in (0,1) $ for $ i \in \{ \ref{cexp:KN}, \ref{cexp:deltaN},\ref{cexp:rN}, \ref{cexp:thetaN}, \ref{cexp:tail_right_exp}, \ref{cexp:tail_right_N} \}$, $\Beta{proof} \in (0,\Beta{zeta})$, $\Epsilon{proof}\in (0,\Epsilon{zeta})$ and $ \cinit{small} \in (0,1) $ 
		such that
		\begin{align}
			\cexp{tail_right_exp} < 1 - \frac{2\alpha}{3},  && \frac{\alpha}{6} \cexp{tail_right_N} < \cexp{thetaN} \left( 1 - \frac{2\alpha}{3} - \cexp{tail_right_exp} \right), \label{bound_cexp_tail_right}
		\end{align}
		and, under the assumptions of Theorem~\ref{thm:main_result}, for $ 1 \leq i \leq \thesigmaN $, $ i \neq \ref{sigmaN:eta} $,
		\begin{equation}  \label{control_sigmaN_i}
			\lim_{N \to \infty} \P{ \sigma_{i,N} \leq T \wedge \hat{\sigma}_{i,N} } = 0.
		\end{equation}
	\end{proposition}
	
	The statement of Proposition~\ref{prop:control_sigma} for $ i = \ref{sigmaN:dist} $ is proved in Section~\ref{sec:stability}, based on the techniques used in~\cite{funaki_scaling_1995}, using an energy functional to control the distance of $v^N_t$ to the stable manifold, with substantial adaptations to our setting.
	In Section~\ref{sec:interface}, we define $\sigmaN{interface}$ and prove Proposition~\ref{prop:control_sigma} for $i = \ref{sigmaN:interface}$.
	The statement for $ i \in \lbrace \ref{sigmaN:supnorm}, \ref{sigmaN:vdelta}, \ref{sigmaN:tail}, \ref{sigmaN:tail_right} \rbrace $ is proved in Section~\ref{sec:holder}, using estimates on the H\"older continuity and tail estimates on $ v^N_t $.

	Proposition~\ref{prop:control_sigma} implies that, for any $ T > 0 $,
	\begin{align*}
		\limsup_{N \to \infty} \P{\sigma_N \leq T} = \limsup_{N \to \infty} \P{\sigmaN{eta} \leq T \wedge \hat{\sigma}_{N,\ref{sigmaN:eta}}}.
	\end{align*}
	Proposition~\ref{prop:holder} is proved in Section~\ref{sec:H_gamma} while Proposition~\ref{prop:control_sigma} is proved in several steps in Sections~\ref{sec:stability}, \ref{sec:interface} and \ref{sec:holder}, using the choice of constants detailed at the end of the present section.
	
	We will also need the following simple lemma.	
	
	\begin{lemma} \label{lemma:diff_m_eta}
	Suppose $\alpha \in (0,2)$.
		There exists a constant $ C > 0 $ such that, for any $ \eta_1, \eta_2 \in \R $,
		\begin{align*}
			\| m_{\eta_1} - m_{\eta_2} \|_{2,\alpha} \leq C \, |\eta_1 - \eta_2|\, e^{\frac{\alpha}{2} (\eta_1 \vee \eta_2)}.
		\end{align*}
		Moreover, for any $ n \in \N $, there exists a constant $ C > 0 $ such that, for any $ \eta_1, \eta_2 \in \R $,
		\begin{align*}
		\| m_{\eta_1} - m_{\eta_2} \|_{H^{n,\alpha}} \leq C \, |\eta_1 - \eta_2|\, e^{\frac{\alpha}{2} (\eta_1 \vee \eta_2)}.
		\end{align*}
	\end{lemma}

	\begin{proof}
		We begin with the first statement.
		By Taylor's formula,
		\begin{align*}
			\| m_{\eta_1} - m_{\eta_2} \|_{2,\alpha}^2 &= \int_{\R} \left( \int_{0}^{1} \partial_{x} m(x - \eta_1 - t (\eta_2 - \eta_1)) dt \right)^2 (\eta_1 - \eta_2)^2 e^{\alpha x} dx.
		\end{align*}
		Using Jensen's inequality and then a change of variables, we obtain
		\begin{align*}
			\| m_{\eta_1} - m_{\eta_2} \|_{2,\alpha}^2 & \leq |\eta_1- \eta_2|^2 \int_{0}^{1} \int_\R |\partial_{x} m(x - \eta_1 - t (\eta_2-\eta_1))|^2 e^{\alpha x} dx dt \\
			&= |\eta_1-\eta_2|^2 \|\partial_{x} m\|_{2,\alpha}^2 \int_{0}^{1} e^{\alpha (\eta_1 + t (\eta_2-\eta_1))} dt.
		\end{align*}
		The first statement follows, and the second statement follows by the same argument for each term in the definition of $\|m_{\eta_1}-m_{\eta_2}\|_{H^{n,\alpha}}$ in~\eqref{def:Hnorm_n}, since $\partial^n_x m\in L^{2,\alpha}$ $\forall n\in \N$.
	\end{proof}	
	
	We shall also use the following auxiliary result, which follows from a straightforward application of H\"older's inequality.

	\begin{lemma} \label{lemma:Lpq_L2alpha_delta}
		For any $ p > 0 $, $ q \in \R $, $\delta \in (0,1 \wedge \frac{2}{p}) $, and any $ \phi : \R \to \R $ such that the right-hand side is finite,
		\begin{equation*}
			\| \phi \|_{p,q} \leq \| \phi \|_{2,\alpha}^{\delta} \| \phi \|_{p(\delta),q(\delta)}^{1-\delta},
		\end{equation*}
		where
		\begin{align*}
			p(\delta) := \frac{p(1-\delta)}{1-\frac{p\delta}{2}}, && q(\delta) := \frac{q - \frac{\alpha p \delta}{2}}{1-\frac{p\delta}{2}}.
		\end{align*}
	\end{lemma}

	Note that $ (p(\delta), q(\delta)) \to (p,q) $ as $ \delta \to 0 $ and that
	\begin{equation} \label{qp_delta}
		\frac{q(\delta)}{p(\delta)} - \frac{\alpha}{2} = \frac{1}{1-\delta} \left( \frac{q}{p} - \frac{\alpha}{2} \right),
	\end{equation}
	so that both sides of the above equation share the same sign.
	The following lemma will also be used in the proof of Theorem~\ref{thm:main_result}.
	It is proved in Appendix~\ref{sec:sobolev}.

	\begin{lemma} \label{lemma:sobolev_average}
		Let $ g : [-1,1] \to \R $ be a Lipschitz continuous function with $ g(0) = 1 $.
		For any $ \gamma \in \R $ and $ \gamma' \in (0,1] $, there exists $ C > 0 $ such that, for all $ v \in H^{\gamma + \gamma', \alpha} $ and $ \delta \in (0,1) $,
		\begin{equation*}
			\Hnorm{(\rho^\delta \cdot g) \ast v - v}{\gamma} \leq C \delta^{\gamma'} \Hnorm{v}{\gamma+\gamma'}.
		\end{equation*}
	\end{lemma}
	
	Let us now prove our main result.
	Note that, from the definition of $ (v^N_t, t \geq 0) $ in \eqref{def_vN} and from~\eqref{spde_u_N}, $ (v^N_t, t \geq 0) $ is a mild solution of
		\begin{align} \label{spde_vN}
			d v^N_t = N \left( \partial_{xx} v^N_t + \alpha \partial_x v^N_t + f(v^N_t) \right) dt + \sqrt{v^N_t(1-v^N_t)} d W(t),
		\end{align}
		where $ (W(t), t \geq 0) $ is again a cylindrical generalised Wiener process on $ L^2(\R) $.
	
	\begin{proof}[Proof of Theorem~\ref{thm:main_result}]
		Let
		\begin{equation*}
			\mathcal V_\zeta := \{v\in L^{2,\alpha}:\exists \, \zeta(v)\text{ s.t. }\|\Phi_t(v,\cdot)-m_{\zeta(v)}\|_{2,\alpha}\to 0\text{ as }t\to \infty\}.
		\end{equation*}
		Recall the definition of $ \drift $ and $ \variance $ in \eqref{parameters_xi} in the statement of Theorem~\ref{thm:main_result}, and recall from Corollary~\ref{cor:A2} that $|\drift|<\infty$ and that $\variance<\infty$.
	
		Fix arbitrary $ \epsilon > 0 $ and $ T > 0 $, and let $ K > 0 $ be such that
		\begin{align} \label{bound_K_proof}
			K > \variance \sqrt{\frac{T}{\epsilon}} + \drift T + \Cst{sup_zeta}(K_0) + 2\epsilon,
		\end{align}
		where $ K_0 $ is given in Assumption~\ref{assumpt:v0}.
		We choose constants $c_i \in (0,1)$, $\cinit{small} > 0$, $\Beta{proof} \in (0,\Beta{zeta})$, and $\Epsilon{proof}\in (0,\Epsilon{zeta})$ such that Proposition~\ref{prop:control_sigma} is satisfied.
		We then define the stopping times $ \sigma_{i,N} $ and $ \sigma_N $ accordingly (as in~\eqref{def_sigma_dist}-\eqref{eq:sigmatailrightdefn} and~\eqref{def_sigma_N}), noting that $\sigma_N$ depends on $\epsilon$ through our choice of $K$. For $ t \in [0,T] $, we set
		\begin{align*}
			\xi^{N,\epsilon}_t := \zeta(v^N_{t \wedge \sigma_N}).
		\end{align*}
		To see that this is well defined, first note that, by Theorem~\ref{thm:zeta}by the definition of $ \sigma_N $, we have $ v^N_{t} \in \mathcal{V}_{\vartheta_N,K,\Epsilon{proof}} $ for all $ t < \sigma_N $ and $ N \geq 1 $, and since $ \vartheta_N \leq \Beta{zeta} $ and $ \Epsilon{proof} \leq \Epsilon{zeta} $.
		Moreover, the stopping time $\sigma_N$ is such that $v^N_t$ remains in a closed subset of $\mathcal{V}_{\vartheta_N,K,\Epsilon{proof}}  \subset L^{2,\alpha}$ while $t < \sigma_N$. 
		Since $t\mapsto v^N_t$ is continuous in $L^{2,\alpha}$, this implies that $v^N_{\sigma_N} \in \mathcal{V}_{\vartheta_N,K,\Epsilon{proof}}$.
		
		We shall prove that, for such a choice of $ K $, $ \Beta{proof} $ and $ \Epsilon{proof} $,
		\begin{align} \label{bound_sigmaN}
			\limsup_{N \to \infty}\, \P{\sigma_N \leq T} \leq \epsilon,
		\end{align}
		and that, for all $ N \geq 1 $, there exists a standard Brownian motion $ (B^N_t, t \geq 0) $ such that
		\begin{align} \label{bound_diff_xi_tilde}
			\limsup_{N \to \infty}\, \P{\sup_{t \in [0,T]} \abs{\xi^{N,\epsilon}_t - (\zeta(v^N_0) - \drift t + \variance B^N_t)} > \epsilon } \leq \epsilon.
		\end{align}
		The statement of Theorem~\ref{thm:main_result} will now follow from \eqref{bound_sigmaN} and \eqref{bound_diff_xi_tilde}. Indeed, for $ t < \sigma_N $, using the triangle inequality, Lemma~\ref{lemma:diff_m_eta} and Theorem~\ref{thm:zeta},
		\begin{align*}
			\| v^N_t - m_{\xi^{N}_t} \|_{2,\alpha} \leq \dist(v^N_t, M) +  C e^{\frac{\alpha}{2} (K+\Cst{eta} \vartheta_N) \vee \Cst{sup_zeta}} | \zeta(v^N_t) - \eta(v^N_t) |.
		\end{align*}
		Then, by \eqref{def_sigma_dist} and \eqref{diff_eta_zeta}, we obtain that there exists $ C > 0 $ such that,
		\begin{equation} \label{vN-m_xi_2alpha}
			\| v^N_t - m_{\xi^{N}_t} \|_{2,\alpha} \leq C \vartheta_N \quad \text{for all }t < \sigma_N .
		\end{equation}
		Hence for $\epsilon '>0$, for $N$ sufficiently large that $C\vartheta_N < \epsilon'$, we have
		\[
		\P{\sup_{t\in [0,T]}\|v^N_t-m_{\xi^N_t}\|_{2,\alpha}>\epsilon'}\le \P{\sigma_N\le T}\le 2\epsilon
		\]
		where the last inequality holds for $N$ sufficiently large by~\eqref{bound_sigmaN}.
		Since $\epsilon>0$ was arbitrary, it follows that
		$\sup_{t \in [0,T]} \| v^N_t - m_{\xi^N_t} \|_{2,\alpha} \to 0$
		in probability as $N\to \infty$.
		
		We now prove \eqref{bound_sigmaN} and \eqref{bound_diff_xi_tilde}.
		As we observed above, $ v^N_t \in \mathcal{V}_{\vartheta_N,K,\Epsilon{proof}}^{(\ggamma{weight})} $ for all $ t \leq \sigma_N $. 
		Hence by Theorem~\ref{thm:zeta}, $v^N_t \in \Phi_r^{-1}(\mathcal V_{K+1}) \subset U_r$ for all $ r \geq 0 $ and $ t \leq \sigma_N $, where $U_r$ was defined in \eqref{def:Ut}.
		For $ r \geq 0 $, recalling the notation $ \eta_r $ introduced in \eqref{def:eta_t}, set
		\begin{equation*}
			\xi^{N,\epsilon}_{r,t} := \eta_r(v^{N}_{t \wedge \sigma_N}).
		\end{equation*}
		By \eqref{diff_eta_zeta},
		\begin{equation} \label{lim_xi_r}
			| \xi^{N,\epsilon}_{r,t} - \xi^{N,\epsilon}_t | \leq \Cst{zeta} \vartheta_N e^{-\lambda r}.
		\end{equation}
		We would like to apply It\^o's formula to write the semi-martingale decomposition of $ \xi^{N,\epsilon}_{r,t} $, but we note that the drift term in \eqref{spde_vN} is not in $ L^{2,\alpha} $, and so not in $ H $.
		To circumvent this, we use $ v^{\delta,N}_t $ as defined in \eqref{def_vdelta}.
		Then, from \eqref{spde_vN}, we write, for $ t \geq 0 $,
		\begin{equation*}
		v^{\delta,N}_t = \rho^\delta \ast v^{N}_0 + \int_{0}^{t} A_{N,\delta}(s) ds + \int_{0}^{t} B_{N,\delta}(s) dW(s),
		\end{equation*}
		where $ A_{N,\delta} $ is an $ H $-valued process and $ B_{N,\delta} $ takes values in the space $ L_{HS}(L^2(\R), H) $ of Hilbert-Schmidt operators from $ L^2(\R) $ to $ H $, respectively defined by
		\begin{align*}
			A_{N,\delta}(s) &:= N \left( \partial_{xx} v^{\delta,N}_s + \alpha \partial_x v^{\delta,N}_s + \rho^\delta \ast f(v^N_s) \right), \\
			B_{N,\delta}(s) &: \phi \mapsto \rho^\delta \ast \left( \sqrt{v^N_s(1-v^N_s)} \phi \right).
		\end{align*}
		Then, using the fact that $ \partial_x v^{N,\delta}_s = (\partial_x \rho^\delta) \ast v^N_s $, $ \partial_{xx} v^{N,\delta}_s = (\partial_{xx} \rho^\delta) \ast v^N_s $ and $ |f(u)| \leq C |u| $, we deduce the existence of a constant $ C > 0 $ such that, for any $ s \geq 0 $,
		\begin{equation*}
			\| A_{N,\delta}(s) \|_H \leq C \delta^{-2} N \| v^N_s \|_H.
		\end{equation*}
		Now using the fact that $ 0 \leq v^N_s \leq 1 $, we obtain that 
		\[ \| v^N_s \|_H^2 \leq \frac{3}{\alpha} + \| v^N_s \|_{2,\alpha}^2 \leq \frac{3}{\alpha} + \Hnorm{v^N_s}{\gamma}^2,\]
		and by Proposition~\ref{prop:holder}, we deduce that for any $ T > 0 $,
		\begin{equation*}
			\sup_{s \in [0,T]} \| A_{N,\delta}(s) \|_H  < \infty
		\end{equation*}
		almost surely.
		The space $ L_{HS}(L^2(\R), H) $ is equipped with the norm
		\begin{equation*}
			\| B \|_{HS} := \left( \Tr{ (B) (B)^*}  \right)^{1/2}, \quad \text{ for } B \in L_{HS}(L^2(\R), H),
		\end{equation*}
		where $ (B)^* : H \to L^2(\R) $ is the adjoint of $ B $, see \citep[Section~I.4.2]{prato_stochastic_2014}.
		Then, by the definition of $ B_{N,\delta} $, we have
		\begin{equation*}
			\| B_{N,\delta}(s) \|_{HS} = \left( \int_{\R^2} \rho^\delta(x-y)^2 v^N_s(y) (1-v^N_s(y)) w(x) dx dy \right)^{1/2},
		\end{equation*}
		where $ w $ was defined in \eqref{def:w}.
		Then, using \eqref{def_sigma_tail} and the fact that $ 0 \leq v^N_s \leq 1 $, we deduce that there exists $ C > 0 $ such that, almost surely
		\begin{equation*}
			\sup_{s \leq \sigma_N} \| B_{N,\delta}(s) \|_{HS} \leq C \sqrt{\frac{\kappa_N}{\delta}}.
		\end{equation*}
		Let us set
		\begin{equation*}
			\sigma_{N,\delta} := \inf\{ t \geq 0 : v^{\delta,N}_t \not\in \Phi_r^{-1}(\mathcal V_{K+1}) \} \wedge \sigma_N.
		\end{equation*}
		Then, using Proposition~\ref{prop:Deta_t} and the infinite-dimensional It\^o formula \citep[Theorem~4.32]{prato_stochastic_2014}, we write
		\begin{multline} \label{ito_vNdelta}
			\eta_r(v^{\delta,N}_{t \wedge \sigma_{N,\delta}}) = \eta_r(v^{\delta,N}_0) + \int_{0}^{t \wedge \sigma_{N,\delta}} D\eta_r(v^{\delta,N}_s) A_{N,\delta}(s)  ds \\ + \int_{0}^{t \wedge \sigma_{N,\delta}} D\eta_r(v^{\delta,N}_s) B_{N,\delta}(s) dW(s) + \frac{1}{2} \int_{0}^{t \wedge \sigma_{N,\delta}} \Tr{ D^2 \eta_r(v^{\delta,N}_s) B_{N,\delta}(s) B_{N,\delta}(s)^* } ds.
		\end{multline}
		In the last term, $ D^2 \eta_r(v^N_s) $ stands for the linear operator $ L : H \to H $ such that, for all $ h_1, h_2 \in H $,
		\begin{equation*}
			D^2 \eta_r(v^N_s) (h_1, h_2) = \langle L h_1, h_2 \rangle_H.
		\end{equation*}
		The remainder of the proof consists in taking successive limits in the above equation as $ \delta \to 0 $, $ r \to \infty $ and $ N \to \infty $.

		\paragraph*{Limit as $\delta \to 0$}

		By Lemma~\ref{lemma:sobolev_average} and Proposition~\ref{prop:holder}, for any $ T > 0 $,
		\begin{equation} \label{cvg_vdeltaN}
			\sup_{s \in [0,T]} \| v^{\delta,N}_s - v^N_s \|_{2,\alpha} \leq C \delta^{\gamma} U_{N,T} \to 0 \quad \text{as }\delta \to 0,
		\end{equation}
		almost surely.
		Hence, using the fact that $ \Phi_r^{-1}(\mathcal V_{K+1}) $ is an open subset of $ H \subset L^{2,\alpha} $ and that $ v^N_t \in \Phi_r^{-1}(\mathcal V_{K+1}) $ for all $ t \leq \sigma_N $, we deduce that,
		\begin{equation*}
			\lim_{\delta \to 0} \sigma_{N,\delta} = \sigma_N \quad \text{almost surely}.
		\end{equation*}
		Furthermore, by the continuity of $ \eta_r : H \to \R $,
		\begin{equation*}
			\lim_{\delta \to 0} \eta_r(v^{\delta,N}_{t \wedge \sigma_{N,\delta}}) = \eta_r(v^N_{t \wedge \sigma_N}) \quad \text{almost surely}.
		\end{equation*}
		In addition, setting
		\begin{equation*}
			A_N(s) := N \left( \partial_{xx} v^N_s + \alpha \partial_x v^N_s + f(v^N_s) \right),
		\end{equation*}
		we have, for any $ s \geq 0 $ and $ \gamma' \in (0,2] $, using the representation \eqref{Deta_t_scalar},
		\begin{multline*}
			| D\eta_r(v^{\delta,N}_s) A_{N,\delta}(s) - D\eta_r(v^{\delta,N}_s) A_N(s) | \leq N \Hnorm{ D\eta_r(v^{\delta,N}_s) }{2-\gamma'}  \Hnorm{ v^{N,\delta}_s - v^N_s }{\gamma'} \\+ N \| D\eta_r(v^{\delta,N}_s) \|_{2,\alpha} \| f(v^{N,\delta}_s) - f(v^N_s) \|_{2,\alpha}.
		\end{multline*}
		Choosing $ \gamma' < \gamma $ where $ \gamma $ is given in Proposition~\ref{prop:holder} and using \eqref{Detat_Hgamma}, Lemma~\ref{lemma:sobolev_average}, we deduce that, for all $ t \geq 0 $,
		\begin{equation*}
			\sup_{s \in [0, t \wedge \sigma_{N,\delta}]} | D\eta_r(v^{\delta,N}_s) A_{N,\delta}(s) - D\eta_r(v^{\delta,N}_s) A_N(s) | \to 0 \quad \text{as }\delta \to 0,
		\end{equation*}
		almost surely.
		In addition, using $ |f(u)| \leq C |u| $, there exists $ C > 0 $ such that
		\begin{equation*}
			| D\eta_r(v^{\delta,N}_s) A_N(s) - D\eta_r(v^N_s) A_N(s) | \leq C N \Hnorm{ D\eta_r(v^{\delta,N}_s) - D\eta_r(v^N_s) }{2-\gamma} \Hnorm{v^N_s}{\gamma}.
		\end{equation*}
		By \eqref{diff_Detat_Hgamma} and \eqref{cvg_vdeltaN}, we deduce that, for all $ t \geq 0 $,
		\begin{equation*}
			\sup_{s \in [0, t \wedge \sigma_{N,\delta}]} | D\eta_r(v^{\delta,N}_s) A_N(s) - D\eta_r(v^N_s) A_N(s) | \to 0 \quad \text{as }\delta \to 0,
		\end{equation*}
		almost surely.
		Turning to the third term on the right-hand side of \eqref{ito_vNdelta}, let us define a linear operator $ B_N(s) : L^2(\R) \to H $ by
		\begin{equation*}
			B_N(s) : \phi \mapsto \sqrt{v^N_s(1-v^N_s)} \phi.
		\end{equation*}
		Then, for any $ T \geq 0 $,
		\begin{multline*}
			\E{ \sup_{t \in [0,T]} \left| \int_{0}^{t \wedge \sigma_{N}} \left( \1{s \leq \sigma_{N,\delta}} D\eta_r(v^{\delta,N}_s) B_{N,\delta}(s) - D\eta_r(v^N_s) B_N(s) \right) dW(s) \right|^2 } \\ \leq 4 \E{ \int_{0}^{T \wedge \sigma_N} \| \1{s \leq \sigma_{N,\delta}} D\eta_r(v^{\delta,N}_s) B_{N,\delta}(s) - D\eta_r(v^N_s) B_N(s) \|_{2}^2 ds },
		\end{multline*}
		where both $ D\eta_r(v^{\delta,N}_s) B_{N,\delta}(s) $ and $ D\eta_r(v^N_s) B_N(s) $ are viewed as elements of $ L^2(\R) $ by identifying the Hilbert-Schmidt operator $ D\eta_r(v^{\delta,N}_s) B_{N,\delta}(s) : L^2(\R) \to \R $ with the element $ \psi $ of $ L^2(\R) $ such that
		\begin{equation*}
			D\eta_r(v^{\delta,N}_s) B_{N,\delta}(s) \phi = \langle \psi,  \phi \rangle_{L^2(\R)}, \quad \forall \phi \in L^2(\R),
		\end{equation*}
		and similarly for $ D\eta_r(v^N_s) B_N(s) $.
		We then write
		\begin{multline} \label{diff_martingale_term}
			\| \1{s \leq \sigma_{N,\delta}} D\eta_r(v^{\delta,N}_s) B_{N,\delta}(s) - D\eta_r(v^N_s) B_N(s) \|_{2} \\ \leq \1{s \leq \sigma_{N,\delta}} \| D\eta_r(v^{\delta,N}_s) (B_{N,\delta}(s) - B_N(s)) \|_{2} \\ + \1{s \leq \sigma_{N,\delta}} \| D\eta_r(v^{\delta,N}_s) B_N(s) - D\eta_r(v^N_s) B_N(s) \|_{2} \\ + \1{s > \sigma_{N,\delta}} \| D\eta_r(v^N_s) B_N(s) \|_{2},
		\end{multline}
		and we bound each term separately.
		Using \eqref{Deta_t_scalar} and the definition of $ B_{N,\delta}(s) $ and $ B_N(s) $, we have
		\begin{multline*}
			\| D\eta_r(v^{\delta,N}_s) (B_{N,\delta}(s) - B_N(s)) \|_{2}^2 \\ = \int_{\R} \left( \int_{\R} D\eta_r(v^{\delta,N}_s, y) \rho^\delta(x-y) e^{\alpha (y-x)} dy - D\eta_r(v^{\delta,N}_s, x) \right)^2 v^N_s(x) (1-v^N_s(x)) e^{2\alpha x} dx.
		\end{multline*}
		By the Cauchy-Schwarz inequality, the right-hand side above is bounded by
		\begin{equation*}
			\| (\rho^\delta \cdot e^{-\alpha \cdot}) \ast D\eta_r(v^{\delta,N}_s) - D\eta_r(v^{\delta,N}_s) \|_{4,\frac{8\alpha}{3}}^2 \| v^N_s (1-v^N_s) \|_{2,\frac{4\alpha}{3}}.
		\end{equation*}
		By the definition of $ \sigmaN{tail_right} $, for any $p > 1$ and $q \in (0,p)$ and $t \leq \sigmaN{tail_right}$,
		\begin{equation*}
			\| v^N_t \|_{p,q} \leq \left( \int_{-\infty}^{\frac{\cexp{tail_right_N}}{1-\cexp{tail_right_exp}} \log(N)} e^{qx} dx + N^{p\cexp{tail_right_N}}  \int_{\frac{\cexp{tail_right_N}}{1-\cexp{tail_right_exp}} \log(N)}^{+\infty} e^{-(p(1-\cexp{tail_right_exp})-q) x} dx \right)^{1/p}.
		\end{equation*}
		If in addition $\cexp{tail_right_exp} < 1 - \frac{q}{p} $, then there exists $ C > 0 $ such that, for $ t \leq \sigmaN{tail_right} $,
		\begin{equation}
			\| v^N_t \|_{p,q} \leq C N^{\frac{\cexp{tail_right_N}q}{p(1-\cexp{tail_right_exp})}}. \label{bound_vN_Lpq}
		\end{equation}
		Since $\alpha \in (0,3/2) $, $ \frac{4\alpha}{3} \in (0,2) $ and, by \eqref{bound_cexp_tail_right}, $ \cexp{tail_right_exp} < 1 - \frac{2\alpha}{3} $, so we can use this to obtain that there exist $ C > 0 $ such that, for all $ s \leq t \wedge \sigma_N $,
		\begin{equation} \label{bound_vN_2_4a/3}
			\| v^N_s (1-v^N_s) \|_{2,\frac{4\alpha}{3}} \leq C N^{\frac{\cexp{tail_right_N} 4\alpha}{3(1-\cexp{tail_right_exp})}}.
		\end{equation}
		In addition, by Lemma~\ref{lemma:Lpq_L2alpha_delta}, for any $ \delta' \in (0, 1/2) $, there exist $ p(\delta') $ and $ q(\delta') $ such that
		\begin{multline*}
			\| (\rho^\delta \cdot e^{-\alpha \cdot}) \ast D\eta_r(v^{\delta,N}_s) - D\eta_r(v^{\delta,N}_s) \|_{4,\frac{8\alpha}{3}} \\
			\begin{aligned}
				&\leq \| (\rho^\delta \cdot e^{-\alpha \cdot}) \ast D\eta_r(v^{\delta,N}_s) - D\eta_r(v^{\delta,N}_s) \|_{2,\alpha}^{\delta'} \\ &\qquad \times \| (\rho^\delta \cdot e^{-\alpha \cdot}) \ast D\eta_r(v^{\delta,N}_s) - D\eta_r(v^{\delta,N}_s) \|_{p(\delta'),q(\delta')}^{1-\delta'} \\
			 	&\leq C (\delta \Hnorm{D\eta_r(v^{\delta,N}_s)}{1})^{\delta'} \| D\eta_r(v^{\delta,N}_s) \|_{p(\delta'),q(\delta')}^{1-\delta'},
			\end{aligned}
		\end{multline*}
		using Lemma~\ref{lemma:sobolev_average} and the fact that there exists $ C > 0 $ such that, for all $ v \in L^{p,q} $,
		\begin{equation}
			\| (\rho^\delta \cdot e^{-\alpha \cdot}) \ast v \|_{p,q} \leq C \| v \|_{p,q}. \label{bound_convolution_Lpq}
		\end{equation}
		Choosing $ \delta' > 0 $ small enough that $ q(\delta') < p(\delta') $ and using \eqref{Detat_Hgamma} and \eqref{Deta_t_Lpq}, we deduce that, for all $ T \geq 0 $,
		\begin{equation} \label{cvg_Deta_Bdelta}
			\sup_{s \in [0, T \wedge \sigma_{N,\delta}]} \| D\eta_r(v^{\delta,N}_s) (B_{N,\delta}(s) - B_N(s)) \|_{2}^2 \to 0 \quad \text{as }\delta \to 0,
		\end{equation}
		almost surely and in $ L^1 $ for all $ N \geq 1$.
		In addition, for the second term on the right-hand side of \eqref{diff_martingale_term}, we have
		\begin{multline*}
			\| D\eta_r(v^{\delta,N}_s) B_N(s) - D\eta_r(v^N_s) B_N(s) \|_{2}^2 \\ 
			\begin{aligned}
				&= \int_{\R} \left( D\eta_r(v^{\delta,N}_s,x) - D\eta_r(v^N_s, x) \right)^2 v^N_s(x) (1-v^N_s(x)) e^{2\alpha x} dx \\
				&\leq \| D\eta_r(v^{\delta,N}_s) - D\eta_r(v^N_s) \|_{4,\frac{8\alpha}{3}}^2 \| v^N_s (1-v^N_s) \|_{2,\frac{4\alpha}{3}} \\
				&\leq C N^{\frac{\cexp{tail_right_N} 4\alpha}{3(1-\cexp{tail_right_exp})}} \| v^{\delta,N}_s - v^N_s \|_{2,\alpha}^2,
			\end{aligned}
		\end{multline*}
		using the Cauchy-Schwarz inequality in the second line, and \eqref{bound_vN_2_4a/3} and \eqref{diff_Detat_Hgamma} in the third.
		Note from Proposition~\ref{prop:holder} that $\E{\sup_{s \in [0,T \wedge \sigma_N]} \Hnorm{v^N_s}{\gamma}^2} < \infty$.
		Combined with Lemma~\ref{lemma:sobolev_average} and $\sigma_{N,\delta} \leq \sigma_N$, we deduce that, for any $ T > 0 $,
		\begin{equation} \label{diff_eta_r_vdelta_v}
			\E{  \sup_{s \in [0, T \wedge \sigma_{N,\delta}]} \| D\eta_r(v^{\delta,N}_s) B_N(s) - D\eta_r(v^N_s) B_N(s) \|_{2}^2 } \to 0 \quad \text{as }\delta \to 0,
		\end{equation}
		Finally, for the last term on the right of \eqref{diff_martingale_term}, for all $ s \leq \sigma_N $,
		\begin{align*}
			\| D\eta_r(v^N_s) B_N(s) \|_{2}^2 &= \int_{\R} D\eta_r(v^N_s, x)^2 v^N_s(x) (1-v^N_s(x)) e^{2\alpha x} dx \\
			&\leq \| D\eta_r(v^N_s) \|_{4,\frac{8\alpha}{3}}^2 \| v^N_s (1-v^N_s) \|_{2,\frac{4\alpha}{3}} \\
			&\leq C N^{\frac{\cexp{tail_right_N} 4\alpha}{3(1-\cexp{tail_right_exp})}}, \numberthis \label{bound_Deta_BN}
		\end{align*}
		using the Cauchy-Schwarz inequality in the second line, and \eqref{Deta_t_Lpq} and \eqref{bound_vN_2_4a/3} in the third.
		Combining \eqref{cvg_Deta_Bdelta}, \eqref{diff_eta_r_vdelta_v} and \eqref{bound_Deta_BN}, we deduce that, for any $ T \geq 0 $,
		\begin{equation*}
			\sup_{t \in [0,T]} \left| \int_{0}^{t \wedge \sigma_{N}} \left( \1{s \leq \sigma_{N,\delta}} D\eta_r(v^{\delta,N}_s) B_{N,\delta}(s) - D\eta_r(v^N_s) B_N(s) \right) dW(s) \right| \to 0 \text{ as }\delta \to 0,
		\end{equation*}
		in $ L^2 $ and hence we can extract a subsequence along which this convergence holds almost surely.
		We now let $ \delta \to 0 $ in the last term on the right of \eqref{ito_vNdelta}.
		By the definition of $ B_{N,\delta}(s) $ and by \eqref{D2eta_t_scalar}, the operator $ D^2 \eta_r(v^{\delta,N}_s) B_{N,\delta}(s) B_{N,\delta}(s)^* : H \to H $ is such that, for all $ h \in H $,
		\begin{multline*}
			D^2 \eta_r(v^{\delta,N}_s) B_{N,\delta}(s) B_{N,\delta}(s)^* h(x) \\ = \frac{e^{\alpha x}}{w(x)} \int_{\R^3} D^2 \eta_r(v^{\delta,N}_s, x, y) \rho^\delta(y-z) \rho^\delta(z-u) v^N_s(z) (1-v^N_s(z)) h(u) w(u) e^{\alpha y} dy dz du.
		\end{multline*}
		As a result,
		\begin{multline*}
			\Tr{ D^2 \eta_r(v^{\delta,N}_s) B_{N,\delta}(s) B_{N,\delta}(s)^* } \\ 
			\begin{aligned}
				&= \int_{\R^3} D^2 \eta_r(v^{\delta,N}_s, x, y) \rho^\delta(y-z) \rho^\delta(z-x) v^N_s(z) (1-v^N_s(z)) e^{\alpha (x + y)} dx dy dz \\
				&= \int_{[-1,1]^2 \times \R} D^2 \eta_r(v^{\delta,N}_s, z + \delta y_1, z + \delta y_2) \rho(y_1) \rho(y_2) v^N_s(z) (1-v^N_s(z)) e^{\alpha (2z + \delta (y_1 + y_2))} dy_1 dy_2 dz.
			\end{aligned}
		\end{multline*}
		Then, for $ s \leq \sigma_{N,\delta} $,
		\begin{multline*}
			\int_{[-1,1]^2 \times \R} \left| D^2 \eta_r(v^{\delta,N}_s, z + \delta y_1, z + \delta y_2) e^{\alpha \delta (y_1 + y_2)} - D^2\eta_r(v^{\delta,N}_s,z,z) \right| \\ \times \rho(y_1) \rho(y_2) v^N_s(z) (1-v^N_s(z)) e^{2\alpha z} dy_1 dy_2 dz \\
			\leq e^{2\alpha \delta} N^{\cexp{tail_right_N}} \int_{[-1,1]^2 \times \R} \left| D^2 \eta_r(v^{\delta,N}_s, z + \delta y_1, z + \delta y_2) - D^2\eta_r(v^{\delta,N}_s,z,z) \right| \rho(y_1) \rho(y_2) e^{(2\alpha-1+\cexp{tail_right_exp}) z} dy_1 dy_2 dz \\ + C \delta N^{\cexp{tail_right_N}} \int_{\R} D^2 \eta_r(v^{\delta,N}_s,z,z) e^{(2\alpha-1+\cexp{tail_right_exp}) z} dy_1 dy_2 dz.
		\end{multline*}
		Note that, by \eqref{bound_cexp_tail_right},
		\begin{equation} \label{bound_exponent}
			2\alpha - 1 + \cexp{tail_right_exp} < \frac{4\alpha}{3}.
		\end{equation}
		Then, using \eqref{D2eta_t_L1q} and \eqref{D2eta_t_L1q_continuity}, we obtain that, for any $ T > 0 $,
		\begin{multline*}
			\sup_{s \in [0, T \wedge \sigma_{N,\delta}]} \int_{[-1,1]^2 \times \R} \left| D^2 \eta_r(v^{\delta,N}_s, z + \delta y_1, z + \delta y_2) e^{\alpha \delta (y_1 + y_2)} - D^2\eta_r(v^{\delta,N}_s,z,z) \right| \\ \times \rho(y_1) \rho(y_2) v^N_s(z) (1-v^N_s(z)) e^{2\alpha z} dy_1 dy_2 dz \to 0 \quad \text{as }\delta \to 0,
		\end{multline*}
		almost surely and in $ L^1 $ for all $ N \geq 1 $.
		Furthermore, for all $ s \leq \sigma_{N,\delta} $,
		\begin{multline*}
			\int_\R | D^2 \eta_r(v^{\delta,N}_s, z, z) - D^2 \eta_r(v^N_s, z, z) | v^N_s(z) (1-v^N_s(z)) e^{2\alpha z} dz \\ \leq C N^{\cexp{tail_right_N}} \int_\R | D^2 \eta_r(v^{\delta,N}_s, z, z) - D^2 \eta_r(v^N_s, z, z) | e^{(2\alpha-1+\cexp{tail_right_exp}) z} dz.
		\end{multline*}
		Using \eqref{bound_exponent} again, along with \eqref{diff_D2eta_t_L1q}, we obtain that the above also vanishes locally uniformly as $ \delta \to 0 $.
		We can thus let $ \delta \to 0 $ in \eqref{ito_vNdelta} and obtain that, for all $ t \geq 0 $,
		\begin{multline} \label{ito_Deta_r}
			\xi^{N,\epsilon}_{r,t} = \eta_r(v^N_0) + N \int_{0}^{t \wedge \sigma_N}  D\eta_r(v^N_s) A_N(s) ds \\ + \int_{0}^{t \wedge \sigma_N} D\eta_r(v^N_s, \cdot) B_N(s) dW(s) \\ + \frac{1}{2} \int_{0}^{t \wedge \sigma_N} \int_{\R} D^2 \eta_r(v^N_s, y, y) v^N_s(y) (1-v^N_s(y)) e^{2\alpha y} dy ds.
		\end{multline}

		\paragraph*{Limit as $r \to \infty$}

		By Lemma~\ref{lemma:katzenberger} and Proposition~\ref{prop:holder}, almost surely,
		\begin{equation*}
			\langle D \zeta(v^N_s, \cdot), \partial_{xx} v^N_s + \alpha \partial_x v^N_s + f(v^N_s) \rangle_\alpha = 0 \quad \forall s\le \sigma_N.
		\end{equation*}
		As a result,
		\begin{multline*}
			| \langle D\eta_r(v^N_s, \cdot), \partial_{xx} v^N_s + \alpha \partial_{x} v^N_s + f(v^N_s) \rangle_\alpha | \leq \| D \eta_r(v^N_s,\cdot) - D\zeta(v^N_s, \cdot) \|_{H^{2-\gamma,\alpha}} \| v^N_s \|_{H^{\gamma, \alpha}} \\ + \| D\eta_r(v^N_s, \cdot) - D\zeta(v^N_s, \cdot) \|_{2,\alpha} \| f(v^N_s) \|_{2,\alpha}.
		\end{multline*}
		Then, combining \eqref{Detat-Dzeta_Hgamma}, and Proposition~\ref{prop:holder}, along with the fact that $ |f(u)| \leq C |u| $, we obtain that, for any $T > 0$,
		\begin{equation*}
			\lim_{r \to \infty} \sup_{s \in [0, T \wedge \sigma_N]} | \langle D\eta_r(v^N_s), \partial_{xx} v^N_s + \alpha \partial_{x} v^N_s + f(v^N_s) \rangle_\alpha | = 0,
		\end{equation*}
		almost surely.
		For the second term in \eqref{ito_Deta_r}, we note that
		\begin{equation} \label{def:diff_martingale}
			t \mapsto  \int_{0}^{t \wedge \sigma_N} \langle D\eta_r(v^N_s,\cdot) - D\zeta(v^N_s,\cdot), \sqrt{v^N_s(1-v^N_s)} dW(s) \rangle_\alpha
		\end{equation}
		is a continuous martingale with quadratic variation
		\begin{equation*}
			t \mapsto \int_{0}^{t \wedge \sigma_N} \int_{\R} | D\eta_r(v^N_s, y) - D\zeta(v^N_s, y) |^2 v^N_s(y) (1-v^N_s(y)) e^{2\alpha y} dy ds.
		\end{equation*}
		By the Cauchy-Schwarz inequality, the integrand of the integral over $s$ above is bounded by
		\begin{equation*}
			\| D\eta_r(v^N_s, \cdot) - D\zeta(v^N_s, \cdot) \|_{4, \frac{8\alpha}{3}}^2 \| v^N_s (1-v^N_s) \|_{2,\frac{4\alpha}{3}}.
		\end{equation*}
 		Combining \eqref{bound_vN_2_4a/3} and \eqref{Detat-Dzeta_Lpq}, we deduce that there exist $ C > 0 $ and $ \genericc > 0 $ such that, for all $ s \leq \sigma_N $,
		\begin{equation*}
			\| D\eta_r(v^N_s, \cdot) - D\zeta(v^N_s, \cdot) \|_{4, \frac{8\alpha}{3}} \| v^N_s (1-v^N_s) \|_{2,\frac{4\alpha}{3}} \leq C N^{\frac{\cexp{tail_right_N} 2\alpha}{3(1-\cexp{tail_right_exp})}} e^{-\genericc r}.
		\end{equation*}
		As a result, the martingale defined in \eqref{def:diff_martingale} vanishes almost surely as $ r \to \infty $, locally uniformly in $ t $.
		Finally, turning to the last term in \eqref{ito_Deta_r}, we note that
		\begin{multline*}
			\int_\R | D^2 \eta_r(v^N_s, y, y) - D^2\zeta(v^N_s, y) | v^N_s(y) (1-v^N_s(y)) e^{2\alpha y} dy \\ \leq \left( \int_\R | D^2 \eta_r(v^N_s, y, y) - D^2\zeta(v^N_s, y) |^2 e^{\frac{8\alpha}{3} y} dy \right)^{1/2} \| v^N_s (1-v^N_s) \|_{2,\frac{4\alpha}{3}}.
		\end{multline*}
		Hence, using \eqref{bound_vN_2_4a/3} again and \eqref{D2etat-D2zeta} (it is straightforward to check that \eqref{condition_p_q} is satisfied with $q = 8\alpha/3$ and $p=2$), we obtain that the above vanishes almost surely as $ r \to \infty $, locally uniformly in $ s $.
		We can thus let $ r \to \infty $ in \eqref{ito_Deta_r}, also using \eqref{lim_xi_r}, yielding
		\begin{multline} \label{ito}
		\xi^{N,\epsilon}_t = \zeta(v^N_0) + \int_{0}^{t \wedge \sigma_N} \Big\langle D \zeta (v^N_s,\cdot), \sqrt{v^N_s(1-v^N_s)} dW(s) \Big\rangle_\alpha \\+ \frac{1}{2} \int_{0}^{t \wedge \sigma_N} \int_\R D^2\zeta(v^N_s,y,y) v^N_s(y)(1-v^N_s(y)) e^{2\alpha y} dy\, ds.
		\end{multline}

		\paragraph*{Limit as $N \to \infty$}

		We will now use the fact that $ \dist(v^N_t, M)\le \vartheta_N$ for all $ t \leq \sigma_N $ and the estimates on the Fr\'echet derivatives of $ \zeta $ in Proposition~\ref{prop:Dzeta} to show that we can replace $ v^N_s $ by $ m_{\xi^{N,\epsilon}_s} $ to compute the limits of the above integrals. 
		More precisely, let $ (\tilde{\xi}^N_t, t \geq 0) $ be defined as
		\begin{multline} \label{def_xi_tilde}
			\tilde{\xi}^N_t := \zeta(v^N_0) + \int_{0}^{t} \Big\langle D\zeta(m_{\xi^{N,\epsilon}_{s}},\cdot), \sqrt{m_{\xi^{N,\epsilon}_{s}} (1-m_{\xi^{N,\epsilon}_{s}})} dW(s) \Big\rangle_{\alpha} \\ + \frac{1}{2} \int_{0}^{t} \int_\R D^2 \zeta(m_{\xi^{N,\epsilon}_{s}}, y, y) m_{\xi^{N,\epsilon}_{s}}(y) (1-m_{\xi^{N,\epsilon}_{s}}(y)) e^{2\alpha y} dy\, ds.
		\end{multline}
		(We omit the dependence of $\tilde{\xi}^N_t$ on $\epsilon$.)
		We then note that, from \eqref{Dzeta_m} and \eqref{parameters_xi}, for any $ \eta \in \R $,
		\begin{align*}
		\int_\R D\zeta(m_{\eta}, y)^2 m_{\eta}(y)(1-m_{\eta}(y)) e^{2\alpha y} dy = \variance^2,
		\end{align*}
		and, from \eqref{D2zeta_m} and~\eqref{parameters_xi}, as noted in~\eqref{A2_finite},
		\begin{align*}
		\frac{1}{2} \int_\R D^2\zeta(m_{\eta},y,y) m_{\eta}(y)(1-m_{\eta}(y)) e^{2\alpha y} dy = -\drift.
		\end{align*}
		It then follows that, for each $ N \geq 1 $, there exists a standard Brownian motion $ (B^N_t, t \geq 0) $ such that $v_0^N$ and $B^N$ are independent and
		\begin{align} \label{def_B_N}
			\tilde{\xi}^N_t = \zeta(v^N_0) - \drift t + \variance B^N_t.
		\end{align}
		
		We now use this to prove \eqref{bound_sigmaN} and \eqref{bound_diff_xi_tilde}.
		Subtracting \eqref{def_xi_tilde} from \eqref{ito}, we see that
		\begin{align} \label{diff_xi_xi_tilde}
			\xi^{N,\epsilon}_t - \tilde{\xi}^N_t = \frac{1}{2} \int_{0}^{t} V^N_s ds + M^N_t,
		\end{align}
		where
		\begin{multline*}
			V^N_s = \1{s \leq \sigma_N} \int_\R D^2\zeta(v^N_s,y,y) v^N_s(y)(1-v^N_s(y)) e^{2\alpha y} dy \\ - \int_\R D^2 \zeta(m_{\xi^N_{s}}, y, y) m_{\xi^N_{s}}(y) (1-m_{\xi^N_{s}}(y)) e^{2\alpha y} dy,
		\end{multline*}
		and $ (M^N_t, t \geq 0) $ is a continuous local martingale with quadratic variation
		\begin{multline} \label{q_var_MN}
			\left[ M^N \right]_t = \int_{0}^{t} \int_\R \left( \1{s \leq \sigma_N} D\zeta(v^N_s,y) \sqrt{v^N_s(y)(1-v^N_s(y))} \right.\\ \left.- D\zeta(m_{\xi^N_{s}},y) \sqrt{m_{\xi^N_{s}}(y) (1-m_{\xi^N_{s}}(y)) \vphantom{v^N_s} } \right)^2 e^{2\alpha y} dy\, ds.
		\end{multline}
		Let us first deal with $ V^N_s $.
		By the Cauchy-Schwarz inequality and, noting that, for $ v, m \in [0,1] $, $ | v(1-v) - m(1-m) | \leq |v - m | $, for any $ s \leq \sigma_N $ we have
		\begin{align} \label{eq:absV}
		\abs{V^N_s} &=\bigg|\int_{\R} D^2 \zeta (v^N_s,y) (v^N_s(y)(1-v^N_s(y))-m_{\zeta(v^N_s)}(y)(1-m_{\zeta(v^N_s)}(y)))e^{2\alpha y}dy \notag\\
		&\qquad +\int_{\R}(D^2\zeta (v^N_s,y)-D^2 \zeta (m_{\zeta(v^N_s)},y))m_{\zeta(v^N_s)}(y)(1-m_{\zeta(v^N_s)}(y))e^{2\alpha y}dy \bigg| \notag\\
		&\leq \| v^N_s - m_{\zeta(v^N_s)} \|_{2,\frac{4\alpha}{3}}  \| D^2\zeta(v^N_s) \|_{2,\frac{8\alpha}{3}} \notag \\
		&\quad + \| m_{\zeta(v^N_s)}(1-m_{\zeta(v^N_s)}) \|_{2,\frac{4\alpha}{3}}  \| D^2\zeta(v^N_s)-D^2\zeta(m_{\zeta(v^N_s)}) \|_{2,\frac{8\alpha}{3}}.
		\end{align}
		Note that, by Theorem~\ref{thm:zeta} and~\eqref{def_sigma_dist}-\eqref{def_sigma_sup}, for $ s \leq \sigma_N $ we have $ | \zeta(v^N_s) | \leq \Cst{sup_zeta} $, and hence
		\begin{align} \label{bound_m_1-m}
			\| m_{\zeta(v^N_s)}(1-m_{\zeta(v^N_s)}) \|_{2,\frac{4\alpha}{3}} \leq \| m(1-m) \|_{2,\frac{4\alpha}{3}} e^{\frac{2\alpha}{3}\Cst{sup_zeta}},
		\end{align}
		which is finite since $\alpha \in (0,3/2)$.
		In addition, by Lemma~\ref{lemma:Lpq_L2alpha_delta}, for any $ \delta \in (0,1) $,
		\begin{equation*}
			\| v^N_s - m_{\zeta(v^N_s)} \|_{2,\frac{4\alpha}{3}} \leq \| v^N_s - m_{\zeta(v^N_s)} \|_{2,\alpha}^\delta \| v^N_s - m_{\zeta(v^N_s)} \|_{2,q(\delta)}^{1-\delta},
		\end{equation*}
		where
		\begin{align} \label{def:q_delta}
			q(\delta) = \frac{(4 - 3 \delta) \alpha}{3(1-\delta)}.
		\end{align}
		Note that, by \eqref{bound_cexp_tail_right}, 
		\begin{equation*}
			\frac{\frac{2\alpha}{3} \cexp{tail_right_N}}{\cexp{thetaN}(1-\cexp{tail_right_exp}) + \frac{\alpha}{2} \cexp{tail_right_N}} < \frac{1 - \frac{2\alpha}{3} - \cexp{tail_right_exp}}{1 - \frac{\alpha}{2} - \cexp{tail_right_exp}}.
		\end{equation*}
		We then choose $ \delta \in (0,1) $ such that
		\begin{equation}
			\frac{\frac{2\alpha}{3} \cexp{tail_right_N}}{\cexp{thetaN}(1-\cexp{tail_right_exp}) + \frac{\alpha}{2} \cexp{tail_right_N}} < \delta < \frac{1 - \frac{2\alpha}{3} - \cexp{tail_right_exp}}{1 - \frac{\alpha}{2} - \cexp{tail_right_exp}}. \label{condition_delta}
		\end{equation}
		Using \eqref{qp_delta} and plugging the upper bound on $\delta$ on the right-hand side, we obtain that
		\begin{equation*}
			\frac{q(\delta)}{2} < 1 - \cexp{tail_right_exp}.
		\end{equation*}
		It follows that $ \| m_{\zeta(v^N_s)} \|_{2,q(\delta)} \leq C $ for some $ C > 0 $ and, by \eqref{bound_vN_Lpq}, we obtain that there exists $ C > 0 $ such that, for all $ s \leq \sigma_N $,
		\begin{equation*}
			\| v^N_s - m_{\zeta(v^N_s)} \|_{2,q(\delta)} \leq  C N^{\frac{\cexp{tail_right_N} q(\delta)}{2(1-\cexp{tail_right_exp})}}.
		\end{equation*}
		Combining this with \eqref{vN-m_xi_2alpha}, we obtain that there exists $ C > 0 $ such that, for all $ s \leq \sigma_N $,
		\begin{align*}
			\| v^N_s - m_{\zeta(v^N_s)} \|_{2,\frac{4\alpha}{3}} &\leq C \vartheta_N^\delta N^{\frac{\cexp{tail_right_N} q(\delta)}{2(1-\cexp{tail_right_exp})} (1-\delta)} \\
			&= C N^{- \cexp{thetaN} \delta + \frac{\cexp{tail_right_N} q(\delta)}{2(1-\cexp{tail_right_exp})} (1-\delta)}.
		\end{align*}
		By \eqref{def:q_delta} and the lower bound on $\delta$ in \eqref{condition_delta}, we have
		\begin{equation*}
			 - \cexp{thetaN} \delta + \frac{\cexp{tail_right_N} q(\delta)}{2(1-\cexp{tail_right_exp})} (1-\delta) = \frac{1}{1-\cexp{tail_right_exp}} \left[ \frac{2\alpha}{3} \cexp{tail_right_N} - \delta \left( \cexp{thetaN}(1-\cexp{tail_right_exp}) + \frac{\alpha}{2} \cexp{tail_right_N} \right) \right] < 0,
		\end{equation*}
		and hence there exist $ C > 0 $ and $ \genericc > 0 $ such that, for all $ s \leq \sigma_N $,
		\begin{equation} \label{vN-m_zeta_Lpq}
			\| v^N_s - m_{\zeta(v^N_s)} \|_{2,\frac{4\alpha}{3}} \leq C N^{-\genericc}.
		\end{equation}
		Combining this with \eqref{bound_m_1-m}, \eqref{D2zeta_Lpq} and \eqref{D2zeta_v-m}, we obtain that there exist $ C > 0 $ and $ \genericc > 0 $ such that, 
		\begin{equation} \label{bound_VN}
			\sup_{s \geq 0}\abs{V^N_{s \wedge \sigma_N}} \leq C N^{-\genericc}.
		\end{equation}
		Turning to the martingale term $ M^N_t $, by \eqref{q_var_MN},
		\begin{align}
			&\left[ M^N \right]_{T \wedge \sigma_N}\notag \\ &\leq 2\int_{0}^{T \wedge \sigma_N} \int_\R m_{\xi^N_s}(y) (1-m_{\xi^N_s}(y)) \left( D\zeta(v^N_s, y) - D\zeta(m_{\zeta(v^N_s)},y) \right)^2 e^{2\alpha y} dy ds \notag \\ 
			&\quad +  2 \int_{0}^{T \wedge \sigma_N} \int_\R D\zeta(v^N_s,y)^2 \notag \\
			&\hspace{3.5cm} \cdot \left( \sqrt{v^N_s(y)(1-v^N_s(y))} - \sqrt{m_{\zeta(v^N_s)}(y) (1-m_{\zeta(v^N_s)}(y))} \right)^2 e^{2\alpha y} dy ds. \label{eq:MNquadvar}
		\end{align}
		Note that, for any $ a, b> 0 $,
		\begin{align*}
			\abs{\sqrt{a}-\sqrt{b}} = \abs{\frac{a-b}{\sqrt{a}+\sqrt{b}}} &\leq \sqrt{\abs{a-b}} \frac{\sqrt{a+b}}{\sqrt{a}+\sqrt{b}} \\
			&\leq \sqrt{\abs{a-b}}.
		\end{align*}
		Then, using the Cauchy-Schwarz inequality for each term on the right-hand side of~\eqref{eq:MNquadvar},
		and recalling the definition of $\|\cdot \|_{p,q}$ from~\eqref{eq:normpq},
		\begin{multline*}
			\left[ M^N \right]_{T \wedge \sigma_N} \leq 2 \int_{0}^{T \wedge \sigma_N} \| m_{\zeta(v^N_s)}(1-m_{\zeta(v^N_s)}) \|_{2,\frac{4\alpha}{3}} \| D\zeta(v^N_s, \cdot) - D\zeta(m_{\zeta(v^N_s)},\cdot) \|_{4,\frac{8\alpha}{3}}^2 ds \\
		    + 2 \int_{0}^{T \wedge \sigma_N} \| D\zeta(v^N_s, \cdot) \|_{4,\frac{8\alpha}{3}}^2 \| v^N_s - m_{\zeta(v^N_s)} \|_{2,\frac{4\alpha}{3}} ds.
		\end{multline*}
		Applying \eqref{bound_m_1-m}, \eqref{vN-m_zeta_Lpq}, \eqref{Dzeta_Lpq} and \eqref{Dzeta_v-m_Lpq}, we obtain that, for any $T \geq 0 $, there exist constants $ C > 0 $ and $\genericc > 0$ such that
		\begin{align*}
		\left[ M^N_{\cdot \wedge \sigma_N} \right]_T \leq C N^{-\genericc},
		\end{align*}
		almost surely, for all $ N $ large enough.
		Hence $ (M^N_{t\wedge \sigma_N}, t \in [0,T]) $ is a square-integrable martingale, and, applying Doob's inequality (see e.g. Theorem~II.1.7 in \cite{revuz_continuous_2013}), we write
		\begin{align*}
		\P{\sup_{t \in [0,T]} | M^N_{t \wedge \sigma_N} | > \frac{\epsilon}{2}} &\leq \frac{4}{\epsilon^2} \E{\left(M^N_{T\wedge \sigma_N}\right)^2} \\
		&= \frac{4}{\epsilon^2} \E{\left[ M^N \right]_{T \wedge \sigma_N}} \\
		&\leq \frac{4 C}{\epsilon^2} N^{-\genericc}.
		\end{align*}
		Combining this with \eqref{diff_xi_xi_tilde} and \eqref{bound_VN}, we obtain that
		\begin{align} \label{bound_sup_xi_xi_tilde}
			\lim_{N \to \infty} \P{\sup_{t \in [0,T\wedge\sigma_N]} \abs{\xi^{N,\epsilon}_t-\tilde{\xi}^N_t} > \epsilon} = 0.
		\end{align}
		We then note that
		\begin{align} \label{bound_sup_eta_tilde_xi}
			\P{\sup_{t \in [0,T]} \abs{\xi^{N,\epsilon}_t - \tilde{\xi}^N_t} > \epsilon} &\leq \P{\sup_{t \in [0,T\wedge \sigma_N]} \abs{\xi^{N,\epsilon}_t - \tilde{\xi}^N_t} > \epsilon}\notag \\ &\quad + \P{\left\lbrace \sigma_N \leq T \right\rbrace \cap \left\lbrace \sup_{t \in [0,T\wedge \sigma_N]} \abs{\xi^{N,\epsilon}_t - \tilde{\xi}^N_t} \leq \epsilon \right\rbrace}.
		\end{align}
		The first term on the right-hand side of~\eqref{bound_sup_eta_tilde_xi} vanishes as $ N \to \infty $ by \eqref{bound_sup_xi_xi_tilde}.
		For the second term, we write
		\begin{align*}
			\left\lbrace \sigma_N \leq T \right\rbrace \subseteq  \left\lbrace \sigmaN{dist} \wedge \sigmaN{supnorm} \leq T \wedge \sigmaN{eta} \right\rbrace \cup \left\lbrace \sigmaN{eta} \leq T \wedge \sigmaN{dist} \wedge \sigmaN{supnorm} \right\rbrace.
		\end{align*}
		Hence,
		\begin{multline*}
			\P{\left\lbrace \sigma_N \leq T \right\rbrace \cap \left\lbrace \sup_{t \in [0,T\wedge \sigma_N]} \abs{\xi^{N,\epsilon}_t - \tilde{\xi}^N_t}\leq \epsilon \right\rbrace} \notag \\  \begin{aligned}
				\leq \quad& \P{\sigmaN{dist} \wedge \sigmaN{supnorm} \leq T \wedge \sigmaN{eta}} \\ &+ \P{\left\lbrace \sigmaN{eta} \leq T \wedge \sigmaN{dist} \wedge \sigmaN{supnorm} \right\rbrace \cap \left\lbrace \sup_{t \in [0,T\wedge \sigma_N]} \abs{\xi^{N,\epsilon}_t - \tilde{\xi}^N_t} \leq \epsilon \right\rbrace}.
			\end{aligned} 
		\end{multline*}
		The first sum vanishes as $ N \to \infty $ by Proposition~\ref{prop:control_sigma}.
		For the last term, note that by~\eqref{diff_eta_zeta} we have that for $t\le \sigma_N$,
		\begin{equation*}
			\abs{\eta(v^N_t)-\xi^{N,\epsilon}_t}=\abs{\eta(v^N_t)-\zeta(v^N_t)}\le \Cst{zeta} \vartheta_N.
		\end{equation*}	
		Hence for $N$ sufficiently large that $\Cst{zeta} \vartheta_N<\epsilon$, by the definition of $ \sigmaN{eta} $ in \eqref{def_sigma_eta} we have
		\begin{align*}
			\P{\left\lbrace \sigmaN{eta} \leq T \wedge \sigmaN{dist} \wedge \sigmaN{supnorm} \right\rbrace \cap \left\lbrace \sup_{t \in [0,T\wedge \sigma_N]} \abs{\xi^{N,\epsilon}_t - \tilde{\xi}^N_t} \leq \epsilon \right\rbrace} \leq \P{\sup_{t \in [0,T]} \abs{\tilde{\xi}^N_t} \ge K - 2\epsilon}.
		\end{align*}
		By~\eqref{eq:u0V} and~\eqref{diff_eta_zeta}, we have $\abs{\zeta(v_0)}\le \Cst{sup_zeta}(K_0) $.
		Therefore, using \eqref{def_B_N}, we obtain
		\begin{align*}
			\P{\sup_{t \in [0,T]} \abs{\tilde{\xi}^N_t} \ge K - 2\epsilon} &\leq \P{\sup_{t \in [0,T]} (\variance \abs{B^N_t}) \ge K - 2\epsilon - |\drift |T - \Cst{sup_zeta}(K_0)} \\
			&\leq \P{\sup_{t \in [0,T]} \abs{B^N_t} \ge \sqrt{\frac{T}{\epsilon}} },
		\end{align*}
		where the second inequality follows by \eqref{bound_K_proof}.
		Using Doob's inequality, we see that
		\begin{align*}
			\P{\sup_{t \in [0,T]} \abs{B^N_t} \ge \sqrt{\frac{T}{\epsilon}} } \leq \epsilon.
		\end{align*}
		Hence, returning to \eqref{bound_sup_eta_tilde_xi} and by~\eqref{def_B_N}, we have shown that
		\begin{align*}
			\limsup_{N \to \infty}\, \P{\sup_{t \in [0,T]} \abs{\xi^{N,\epsilon}_t - (\zeta(v_0^N)- \drift t+\variance B^N_t)} > \epsilon} \leq \epsilon,
		\end{align*}
		as claimed in \eqref{bound_diff_xi_tilde}.
		The bound \eqref{bound_sigmaN} follows from exactly the same arguments since
		\begin{multline*}
			\P{\sigma_N \leq T} \leq  \P{\sup_{t \in [0,T\wedge \sigma_N]} \abs{\xi^{N,\epsilon}_t - \tilde{\xi}^N_t} > \epsilon}\\ + \P{\left\lbrace \sigma_N \leq T \right\rbrace \cap \left\lbrace \sup_{t \in [0,T\wedge \sigma_N]} \abs{\xi^{N,\epsilon}_t - \tilde{\xi}^N_t} \leq \epsilon \right\rbrace}.
		\end{multline*}
		This concludes the proof of our main result.
	\end{proof}

	\subsection*{Outline of the remainder of the article}

	The rest of the paper is devoted to the proofs of the various results which were used in the proof of Theorem~\ref{thm:main_result}.
	It is organised as follows.
	Proposition~\ref{prop:control_sigma} is proved in Sections~\ref{sec:stability}, \ref{sec:interface} and~\ref{sec:holder}.
	Proposition~\ref{prop:holder} is proved in Section~\ref{sec:H_gamma}.
	We then prove Theorem~\ref{thm:zeta} in Section~\ref{sec:det_flow}.
	All the estimates on the Fréchet derivatives of $ \zeta $ are proved in Section~\ref{sec:frechet}, using some intermediate results on the linearised semigroup associated to the deterministic PDE, which are proved in Section~\ref{sec:semigroup}.
	In Appendix~\ref{sec:fermi}, we prove Lemma~\ref{lemma:eta}.
	Appendix~\ref{sec:sobolev} contains the proof of Lemma~\ref{lemma:sobolev_average}, and Appendix~\ref{sec:Phit} proves some properties of the flow $ \Phi^t $ associated to the deterministic PDE.
	Each section is essentially independent of the others, with the following caveats: Section~\ref{sec:holder} uses a result from Section~\ref{sec:interface}, Section~\ref{sec:H_gamma} uses a result from Section~\ref{sec:holder}, Section~\ref{sec:det_flow} uses a result from Section~\ref{sec:stability} (Proposition~\ref{prop:energy_estimates}), and Section~\ref{sec:frechet} uses some results from Section~\ref{sec:det_flow}, as well as the results from Section~\ref{sec:semigroup}.
	
	The proof of Proposition~\ref{prop:control_sigma} is split across three sections, and involves the stopping times $\sigmaN{dist}$ – $\sigmaN{interface}$.
	These stopping times depend on some small constants, some of which have already been introduced.
	Presently we will select the remaining constants and specify relations between them, as promised above. In particular, we will define
	\[\newcexp{holder},  \newcexp{tail_left}, \newcexp{KN}, \newcexp{tail_right_exp}, \newcexp{deltaN}, \newcexp{rN}, \newcexp{thetaN}, \newcexp{tail_right_N},  \newcexp{interface}, \newcexp{inter_integ},  
	 \in (0,1).\]
	For a given value of $\cinit{big}$, the conditions we require on these constants will ultimately impose a constraint on $\cinit{small}$.
	We now proceed systematically with the imposition of conditions on $c_i$, $i =1,\dots,10$. First, we any $\cexp{holder} \in (0,1)$ satisfying
	\begin{equation} \label{eq:holderconstdef}
	\cexp{holder} < \tfrac{\cinit{big}}{2} \wedge \tfrac 1 4.
	\end{equation}
	Noting that $f'(1) < 0$, we choose $\cexp{tail_left}$ to be sufficiently small so that
	\begin{equation} \label{eq:cleftfprime}
	\cexp{tail_left} < \cinit{big}, \quad  \alpha \cexp{tail_left} + \cexp{tail_left}^2 < - \frac 1 2 f'(1).
	\end{equation}
	We next remark that because $\alpha < 2$, there exists $c(\alpha) \in (0,1/2)$ and $b(\alpha) \in (0,1)$ such that
	\begin{equation}\label{eq:csumsmallxyz}
	\sup_{x,y,z \in (0,c(\alpha)]} x + (1+ y)((\alpha + z - 1) \vee 0) = b(\alpha) < 1.
	\end{equation}
	In the sequel we fill choose three different constants to be at most $c(\alpha)$, which will guarantee they satisfy the above. In the meantime, let us choose
	\begin{equation} \label{eq:cKNcsum0}
	\cexp{KN} < 1 - b(\alpha).
	\end{equation}	

	We recall from \eqref{eq:weightdefn} the constant $\ggamma{weight}$,which satisfies $\ggamma{weight} =0$ if $\alpha < 1$ and $\ggamma{weight} \in (\alpha - 1,1/2]$ if $\alpha \in [1,3/2)$. 
	Next, we fix $\cexp{tail_right_exp}$ satisfying
	\begin{equation}\label{bounds_cexptailrightexp}
	\cexp{tail_right_exp} <  \left( 1-\frac{2 \alpha}{3} \right) \wedge c(\alpha), \quad \text{ and } \,  \cexp{tail_right_exp} < 1 - 2\ggamma{weight} \, \text{ if } \alpha \in (1,3/2).
	\end{equation}
	If $\alpha \in (0,1]$, we simply ignore the second condition. Next, we choose $\cexp{deltaN}$ to satisfy
	 \begin{equation} \label{bounds_cexpdeltaN}
		\cexp{deltaN} <  \left( \cinit{big} \wedge \frac{\cexp{KN}}{3} \right). 
	\end{equation}
	Since $\cexp{tail_right_exp} < 1 - \alpha/2$, we may take $\cexp{rN}$ such that 
	\begin{equation} \label{eq:cexprNsmall}
	\cexp{rN}  <\left( \frac{\cexp{holder} \cexp{deltaN}}{4}\right) \wedge \cexp{tail_left} \wedge \frac{\alpha}{4} \wedge (1-\cexp{tail_right_exp}) 
	\wedge \left(\left(\frac 1 2 \wedge \frac{\cexp{holder}\cexp{deltaN}}{4\alpha}\right)(2(1-\cexp{tail_right_exp}) - \alpha) \right).
	\end{equation}
	Due to \eqref{bounds_cexpdeltaN}, we may next take $ \cexp{thetaN} \in (0,1) $ satisfying
	\begin{equation} \label{bounds_cexpthetaN}
		\cexp{thetaN} < \cinit{big} \wedge \cexp{deltaN} \wedge (\cinit{big} - \cexp{deltaN}) \wedge \cexp{rN} \wedge \left( \frac{\cexp{KN} - 3\cexp{deltaN}}{2} \right).
	\end{equation}
	Using the conditions from \eqref{bounds_cexptailrightexp}, \eqref{eq:cexprNsmall} and \eqref{bounds_cexpthetaN}, we may next choose $\cexp{tail_right_N} > 0$ satisfying
	\begin{align} \label{bounds_cexptailrightN}
	\cexp{tail_right_N} & < \,(\cexp{holder}\cexp{thetaN}) \wedge c(\alpha) \wedge (1 - \cexp{tail_right_exp} - \cexp{rN})  
	\\ &\quad \wedge  \left(\left(\frac 1 2 \wedge \frac{\cexp{holder}\cexp{deltaN}}{4\alpha}\right) (2(1-\cexp{tail_right_exp}) - \alpha) - \cexp{rN} \right) 
	\wedge \frac{6 \cexp{thetaN}}{\alpha}\left( 1 - \frac{2\alpha}{3} - \cexp{tail_right_exp} \right). \notag
	\end{align}
	Finally, we take $\cexp{interface}$ satisfying
	\begin{equation} \label{bounds_cexpinterface}
	\cexp{interface} < \cexp{tail_right_N}  \wedge \left((1-\cexp{tail_right_exp})^{-1} -1\right) \wedge c(\alpha)
	\end{equation}
	and $\cexp{inter_integ}$,
	\begin{equation} \label{bounds_cexpinterinteg}
	\cexp{inter_integ} < \left(\frac{\cexp{interface}}{1+\cexp{interface}} \right) \wedge \cexp{tail_right_exp} \wedge c(\alpha) \wedge \left( 1 - \frac{\alpha}{2}\right).
	\end{equation}

	We have now specified values for all of $\cexp{holder}$ – $\cexp{thetaN}$, which are defined only in terms of $\alpha$ and the constant $\cinit{big}$ appearing in Assumption~\ref{assumpt:v0}.
	Given the assignments above, and noting that \eqref{bounds_cexpinterinteg} implies that $(1+\cexp{interface})(1-\cexp{inter_integ}) > 1$, we now take $\cinit{small} >0$ to satisfy
	\begin{equation} \label{bounds_cinitsmall}
	\cinit{small} < \cexp{tail_right_exp} \wedge \cexp{interface} \wedge \cexp{tail_right_N} \wedge \cexp{inter_integ} \wedge \left((1+ \cexp{interface})(1-\cexp{inter_integ})- 1\right).
	\end{equation}
	This is the precise version of the condition, from Theorem~\ref{thm:main_result}, that $\cinit{small} > 0$ is sufficiently small.
	
	We conclude with little more bookkeeping. Since $\cexp{tail_right_exp}, \cexp{interface}, \cexp{tail_right_N} < c(\alpha)$, it follows from \eqref{eq:csumsmallxyz} and \eqref{eq:cKNcsum0} that
	\begin{equation} \label{eq:cKNcsum}
	\cexp{KN} < 1 - \cexp{tail_right_N} - (1+ \cexp{interface})((\alpha + \cexp{tail_right_exp} - 1) \vee 0).
	\end{equation}
	Also, since $\cexp{tail_left} \in (0,1)$ satisfies \eqref{eq:cleftfprime}, we may introduce the auxiliary constant
	\[ \rho_0 := \left(\|f'\|_\infty + \frac 1 2 |f'(1)| \right) \left(\left(\frac 1 2 |f'(1)| - \alpha \cexp{tail_left} - \cexp{tail_left}^2 \right)^{-1} \vee 1\right).\]
	Given this, and recalling the constant $K > 0$, we fix $\newCst{tail_left} = \Cst{tail_left}(K)$ to be sufficiently large so that 
	\begin{equation} \label{eq:lefttailinitconst}
	\Cst{tail_left} \geq 4 \vee (4 \Cst{init}) \vee (4 \rho_0)
	\end{equation}
	and
	\begin{align} \label{eq:CtailleftKexp}
	&\sup_{x \in \R} \sup_{y \in [-K,K]} (1 - m(x+y))e^{\cexp{tail_left}x} \leq \frac{\Cst{tail_left}}{8 \rho_0}.
	 \end{align}

	\section{Stability of the invariant manifold} \label{sec:stability}
	
	The aim of this section is to prove Proposition~\ref{prop:control_sigma} in the case $ i = \ref{sigmaN:dist} $.
	In order to do this, we introduce an energy functional that serves as a potential for the deterministic flow. 
	First, we state a few properties concerning the asymptotic behaviour of $m$. 
	We recall $\lambda_+$ from \eqref{eq:lambda+def}, that we have scaled the equation so that $\lambda_+ = 1$ and $f'(0) = \alpha - 1$,
	for $\alpha \in (0,3/2)$. Just as one computes $\lambda_+$, i.e. the rate of exponential decay of $m$ as $x \to \infty$, by linearizing near $0$,
	linearizing the equation near $1$ and recalling that $f'(1) < 0$ allows us to compute the decay rate as $x \to -\infty$ as
	\begin{equation} \label{eq:lambda-def}
	\lambda_- = -\frac \alpha 2 + \sqrt{\frac{\alpha^2}{4} - f'(1)} >0.
	\end{equation}
	It follows that 
	\begin{equation} \label{eq:masympneg}
	m(x) \sim k e^{-x} \, \text{ and } \, 1-m(-x) \sim k' e^{-\lambda' x} \,\, \text{ as $x \to \infty$.}
	\end{equation}
	for some $k' > 0$. It is also straightforward to show that $\partial_x m(x)$ has the same asymptotic behaviour as $m(x)$ at $\pm \infty$. 
	In particular, by the above and because $m(x) \sim k e^{-x}$ as $x \to + \infty$, we have, for some constants $b,b' > 0$,
	\begin{equation*}
	m(x) \sim b e^{-x} \,  \text{ and } \, 1-m(-x) \sim b' e^{-\lambda' x} \,\, \text{ as $x \to +\infty$}.
	\end{equation*}

	\subsection{The energy functional} \label{subsec:energy_functional}
	
	For $ v \in \R $, let
	\begin{equation} \label{def_F}
		F(v) := - \int_{0}^{v} f(u) du,
	\end{equation}
	where for convenience, we extend the definition of $f$ from~\eqref{f} in such a way that $f$ is smooth and $f(x)=0$ for $|x|\ge 2$.
	For $ v \in H^{1,\alpha} $, define the energy functional
	\begin{equation} \label{eq:Hdefn}
	\mathcal{H}(v) := \int_\R \left( \frac{1}{2} \abs{\partial_x v (x)}^2 + F(v(x)) \right) e^{\alpha x} dx.
	\end{equation}
	This $\mathcal H (v)$ is the equivalent in our $L^{2,\alpha}$ space of the energy functional defined in~\cite[(3.2)]{funaki_scaling_1995}.
	
	\begin{lemma} \label{lemma:H_is_zero_on_M}
	Suppose $\alpha \in (0,2)$.
		For any $ \eta \in \R $, $ \mathcal{H}(m_\eta) = 0 $.
	\end{lemma}

	\begin{proof}
		Recall from~\eqref{scalar_product} that we have set, for $ \phi, \psi \in L^{2,\alpha} $,
		\begin{align*}
			\langle \phi, \psi \rangle_\alpha = \int_\R \phi(x) \psi(x) e^{\alpha x} dx.
		\end{align*}
		Note further that, for $ \phi \in H^{2,\alpha} $ and $ \psi \in H^{1,\alpha} $, integrating by parts,
		\begin{align} \label{IPP_alpha}
			\langle \partial_{x} \phi, \partial_{x} \psi \rangle_\alpha = - \langle \partial_{xx} \phi + \alpha \partial_{x} \phi, \psi \rangle_{\alpha}.
		\end{align}
		Applying~\eqref{IPP_alpha} with $ \phi = v $ and $ \psi = \partial_{x} v $, we obtain, for $ v \in H^{2,\alpha} $,
		\begin{align} \label{IPP2}
			\langle \partial_{xx} v + \alpha \partial_{x} v, \partial_{x} v \rangle_\alpha = - \langle \partial_{x} v, \partial_{xx} v \rangle_\alpha.
		\end{align}
		Noting that $ \partial_{x} v \, \partial_{xx} v = \frac{1}{2} \partial_{x} (\partial_{x} v)^2 $ and integrating by parts again, we obtain
		\begin{align} \label{IPP3}
			\langle \partial_{xx} v + \alpha \partial_{x} v, \partial_{x} v \rangle_\alpha = \frac{\alpha}{2} \int_\R |\partial_{x} v(x)|^2 e^{\alpha x} dx.
		\end{align}
		Another integration by parts yields
		\begin{align*}
			\langle f(v), \partial_{x} v \rangle_\alpha = \alpha \int_\R F(v(x)) e^{\alpha x} dx.
		\end{align*}
		Summing the last two equalities, we obtain
		\begin{align*}
			\mathcal{H}(v) = \frac{1}{\alpha} \langle \partial_{xx} v + \alpha \partial_{x} v + f(v), \partial_{x} v \rangle_\alpha.
		\end{align*}
		The fact that $ \mathcal{H}(m_\eta) = 0 $ for all $\eta\in \R$ then follows from \eqref{def:m},
		 and since $m(x) \sim k e^{-x}$ as $x \to \infty$ (recall \eqref{eq:lambda+def}) and $\alpha \in (0,2)$, we have $ m \in H^{2,\alpha} $.
	\end{proof}
	
	Using \eqref{IPP_alpha}, we note that, for $ v \in H^{2,\alpha} $ and $ h \in H^{2,\alpha} $,
	\begin{align*}
		&\int_\R \left( | \partial_{x} v(x) + \partial_{x} h(x) |^2 - | \partial_{x} v(x) |^2 \right) e^{\alpha x} dx\\ &\quad = - \langle (\partial_{xx} + \alpha \partial_{x}) (v+h), v+h \rangle_{ \alpha} + \langle (\partial_{xx} + \alpha \partial_{x}) v, v \rangle_{ \alpha} \\
		&\quad = - 2 \langle (\partial_{xx} + \alpha \partial_{x}) v, h \rangle_{ \alpha} - \langle (\partial_{xx} + \alpha \partial_{x}) h, h \rangle_\alpha.
	\end{align*}
	As a result, setting for $v,h\in \R$,
	\begin{align} \label{def_Uvh}
	\mathcal{U}(v,h) &:= F(v+h)-F(v) - F'(v) h - \frac{1}{2} F''(v) h^2 \notag \\
	&= F(v+h) - F(v) + f(v) h + \frac{1}{2} f'(v) h^2,
	\end{align}
	we obtain that for $v\in H^{2,\alpha}$ and $h\in H^{2,\alpha}$,
	\begin{multline} \label{expansion_H}
	\mathcal{H}(v + h) - \mathcal{H}(v) = - \langle \partial_{xx} v + \alpha \partial_x v + f(v), h \rangle_{\alpha} - \frac{1}{2} \langle \partial_{xx}h + \alpha \partial_x h + f'(v)h , h \rangle_{\alpha} \\ + \int_\R \mathcal{U}(v(x), h(x)) e^{\alpha x} dx.
	\end{multline}
	Since $ | \mathcal{U}(v,h) | \leq C \min(|h|^3, |h|^2) $ for all $v,h\in \R$ for some constant $ C > 0 $, we see that $ \mathcal{H} $ is Fréchet differentiable on $ H^{2,\alpha} $ and that the operator $ D\mathcal{H}(v) $ can be identified with the function
	\begin{align} \label{def_DH}
		D\mathcal{H}(v) := - (\partial_{xx} v + \alpha \partial_{x} v + f(v)).
	\end{align}
	The aim of this subsection is to prove the following result, which will allow us to use $ \mathcal{H} $ to control the distance between $ v^N_t $ and $ M $.
	Recall from~\eqref{eq:Fermidefn} in Section~\ref{sec:statement} that for $K>0$, for $\beta_0$ as defined in Lemma~\ref{lemma:eta}, the Fermi coordinate $ s(v) \in L^{2,\alpha} $ is defined for $ v \in \mathcal{V}_{K, \Beta{eta}} $ by setting $ s(v) = v - m_{\eta(v)} $.
	
	\begin{proposition} \label{prop:energy_estimates}
		Suppose $\alpha \in (0,2)$.
		There exist $ \newCst{H_below}, \newCst{H_above} > 0 $ and $ \newEpsilon{H} > 0 $ such that, for all $ v \in H^{1,\alpha} $ and $ \eta \in \R $ such that
		\begin{equation} \label{v_eta_uniform}
			\langle v - m_\eta, \partial_{x} m_\eta \rangle_\alpha = 0, \quad \text{ and } \quad \| v - m_\eta \|_\infty \leq \Epsilon{H},
 		\end{equation}
 		then
		\begin{align} \label{energy_estimate_H_1}
			\Cst{H_below} \| v - m_\eta \|_{H^{1,\alpha}}^2 \leq \mathcal{H}&(v) \leq \Cst{H_above} \| v - m_\eta \|_{H^{1,\alpha}}^2.
		\end{align}
		Moreover, if in addition $ v \in H^{2,\alpha} $,
		\begin{align} \label{energy_estimate_H_2}
			\Cst{H_below} \| v - m_\eta \|_{H^{2,\alpha}}^2 \leq \| D\mathcal{H}&(v,\cdot) \|_{{2,\alpha}}^2 \leq \Cst{H_above} \| v-m_\eta \|_{H^{2,\alpha}}^2.
		\end{align}
	\end{proposition}
	
	To prove this, let us define, for $ \eta \in \R $, an operator $ \mathcal{A}_\eta : H^{2,\alpha} \to L^{2,\alpha} $ as
	\begin{align} \label{def_A}
	\mathcal{A}_\eta := - \partial_{xx} - \alpha \partial_x - f'(m_\eta).
	\end{align}
	Then, by \eqref{IPP_alpha}, $ \mathcal{A}_\eta $ is self-adjoint in $ L^{2,\alpha} $ and moreover, differentiating \eqref{def:m} with respect to $ x $,
	\begin{align} \label{eigenfunction_A}
	\mathcal{A}_\eta \partial_{x} m_\eta = 0.
	\end{align}
	Since $ \partial_x m_\eta < 0 $ by~\eqref{def:m}, this implies that the spectrum of $ \mathcal{A} $ lies in $ [0, \infty) $ and that its continuous spectrum lies in $ [\lambda, \infty) $ for some $ \lambda > 0 $ \citep[case 8 of Theorem~10.12.1]{zettl_sturm-liouville_2005}, see also \cite{fife_approach_1977}.
	It follows that there exists $ \specgap > 0 $ such that, if we set
	\begin{align} \label{eq:phietadefn}
	\varphi_\eta := - \frac{\partial_{x} m_\eta}{\| \partial_{x} m_\eta \|_{{2,\alpha}}},
	\end{align}
	then, for all $ h \in H^{1,\alpha} $ such that $ \langle h, \varphi_\eta \rangle_{\alpha} = 0 $,
	\begin{align} \label{spectral_gap}
	\langle \mathcal{A}_\eta h , h \rangle_{{\alpha}} \geq \specgap \left\| h \right\|_{{2,\alpha}}^2,
	\end{align}
	and for $ h \in H^{2,\alpha} $ such that $ \langle h, \varphi_\eta \rangle_{\alpha} = 0 $,
	\begin{align}
		\| \A_\eta h \|_{{2,\alpha}}^2 \geq \specgap^2 \left\| h \right\|_{{2,\alpha}}^2. \label{spectral_gap_L2}
	\end{align}
	Moreover, using that $\langle h_1(\cdot +\eta),h_2(\cdot +\eta)\rangle_{\alpha}=e^{-\alpha \eta}\langle h_1,h_2\rangle_{\alpha}$, we have that $\specgap$ is a constant that does not depend on $\eta$.
	The following lemma, combined with \eqref{expansion_H}, will allow us to prove Proposition~\ref{prop:energy_estimates}.
	This lemma is the equivalent of \cite[Lemma~3.2]{funaki_scaling_1995}.
	
	\begin{lemma} \label{lemma:estimates_A}
		Suppose $\alpha \in (0,2) $. There exist $ \newCst{A_below}, \newCst{A_above} > 0 $ such that, for all $ s \in H^{1,\alpha} $ with $ \langle s, \varphi_\eta \rangle_{{\alpha}} = 0 $,
		\begin{align} \label{energy_estimate_A_1}
		\Cst{A_below} \| s \|_{H^{1,\alpha}}^2 \leq \langle \mathcal{A}_\eta s, s \rangle_{{\alpha}} \leq \Cst{A_above} \| s \|_{H^{1,\alpha}}^2.
		\end{align}
		Moreover, if $s\in H^{2,\alpha}$ with $ \langle s, \varphi_\eta \rangle_{{\alpha}} = 0 $ then
		\begin{align} \label{energy_estimate_A_2}
			\Cst{A_below} \| s \|_{H^{2,\alpha}}^2 \leq \| \mathcal{A}_\eta  s \|_{{2,\alpha}}^2 \leq \Cst{A_above} \| s \|_{H^{2,\alpha}}^2.
		\end{align}
	\end{lemma}

	\begin{proof}[Proof of Lemma~\ref{lemma:estimates_A}]
		Let $ s \in H^{1,\alpha} $ be such that $ \langle s, \varphi_\eta \rangle_{{\alpha}} = 0 $.
		By \eqref{IPP_alpha} and the definition of $ \A_\eta $ in \eqref{def_A},
		\begin{align} \label{quadratic_form}
		\langle \A_\eta s, s \rangle_{{\alpha}} = \int_{\R} \left( \abs{\partial_{x} s (x)}^2 - f'(m_\eta(x)) \abs{s(x)}^2 \right) e^{\alpha x} dx.
		\end{align}
		Hence there exists a constant $ \Cst{A_above} > 0 $ such that $ \langle \A_\eta s, s \rangle_{{\alpha}} \leq \Cst{A_above} \| s \|_{H^{1,\alpha}}^2 $.
		From \eqref{spectral_gap}, since $ \langle s, \varphi_\eta \rangle_{{\alpha}} = 0 $, 
		\begin{align} \label{borne_inf_1}
		\langle \A_\eta s, s \rangle_{{\alpha}} \geq \specgap \| s \|_{{2,\alpha}}^2.
		\end{align}
		Furthermore, from \eqref{quadratic_form}, letting $ f'_{\mathrm{max}} := \sup_{x \in [0,1]} f'(x) >0$,
		\begin{align} \label{borne_inf_2}
		\langle \A_\eta s, s \rangle_{{\alpha}} \geq \| s \|_{H^{1,\alpha}}^2 - (f'_{\mathrm{max}} + 1) \| s \|_{{2,\alpha}}^2.
		\end{align}
		Multiplying \eqref{borne_inf_1} by $ f'_{\mathrm{max}} + 1 $ and multiplying~\eqref{borne_inf_2} by $ \specgap $, then summing the two inequalities, we obtain
		\begin{align*}
		(f'_{\mathrm{max}} + 1 + \specgap) \langle \A_\eta s, s \rangle_{{\alpha}} \geq \specgap \| s \|_{H^{1,\alpha}}^2,
		\end{align*}
		thus proving the lower bound in \eqref{energy_estimate_A_1}.
		
		Now suppose $s\in H^{2,\alpha}$.
		To prove \eqref{energy_estimate_A_2}, let us first note that
		\begin{align*}
			\left\| \mathcal{A}_\eta s \right\|_{2,\alpha}^2 = \| \partial_{xx} s + \alpha \partial_{x} s \|_{2,\alpha}^2 + 2 \langle \partial_{xx} s + \alpha \partial_{x} s, f'(m_\eta) s \rangle_{\alpha} + \| f'(m_\eta) s \|_{2,\alpha}^2.
		\end{align*}
		In addition,
		\begin{align} \label{dxxsidentity}
			\| \partial_{xx} s + \alpha \partial_{x} s \|_{2,\alpha}^2 & = \| \partial_{xx} s \|_{2,\alpha}^2 + 2 \alpha \langle \partial_{xx} s, \partial_{x} s \rangle_{ \alpha} + \alpha^2 \langle \partial_{x} s, \partial_{x} s \rangle_{ \alpha} \notag \\
			&= \| \partial_{xx} s \|_{2,\alpha}^2 + 2 \alpha \langle \partial_{xx} s+\alpha \partial_x s, \partial_{x} s \rangle_{ \alpha} - \alpha^2 \langle \partial_{x} s, \partial_{x} s \rangle_{ \alpha}.
		\end{align}
		By \eqref{IPP3}, we see that the last two terms on the right-hand side of~\eqref{dxxsidentity} cancel.
		As a result,
		\begin{align*}
		\| \A_\eta s \|_{{2,\alpha}}^2 = \| \partial_{xx} s \|_{{2,\alpha}}^2 + 2 \langle \partial_{xx}s + \alpha \partial_x s, f'(m_\eta) s \rangle_{{\alpha}} + \langle f'(m_\eta)^2 s, s \rangle_{{\alpha}}.
		\end{align*}
		Using the fact that $f'$ is bounded and the Cauchy-Schwarz inequality, we obtain that there exists a constant $ \Cst{A_above} > 0 $ such that $ \| \mathcal{A}_\eta s \|_{2,\alpha}^2 \leq \Cst{A_above} \| s \|_{H^{2,\alpha}}^2 $.
		For the lower bound, using~\eqref{IPP_alpha},
		\begin{align} \label{eq:Aetas2alpha}
		\| \A_\eta s \|_{{2,\alpha}}^2 = \| \partial_{xx} s \|_{{2,\alpha}}^2 + \tfrac{1}{2} \| \partial_x s \|_{{2,\alpha}}^2 + \left\langle \partial_{xx}s + \alpha \partial_x s, \left( 2 f'(m_\eta) + \tfrac{1}{2}\right) s \right\rangle_{{\alpha}} + \langle f'(m_\eta)^2 s, s \rangle_{{\alpha}}.
		\end{align}
		In addition, since $\frac 12 \partial_x(s^2)=s\partial_x s$,
		\begin{align*}
		\left\langle \partial_x s, \left(2 f'(m_\eta) + \tfrac{1}{2}\right) s \right\rangle_{\alpha} = \tfrac{1}{2} \left\langle \partial_x (s^2), \left(2 f'(m_\eta) + \tfrac{1}{2}\right) \right\rangle_{\alpha} .
		\end{align*}
		Integrating by parts as in \eqref{IPP_alpha}, we obtain
		\begin{align*}
			\left\langle \partial_x s, \left(2 f'(m_\eta) + \tfrac{1}{2}\right) s \right\rangle_{\alpha} = - \tfrac{1}{2} \left\langle s^2, (\partial_{x} + \alpha) \left(2 f'(m_\eta) + \tfrac{1}{2} \right) \right\rangle_{\alpha}.
		\end{align*}
		Moreover, setting $ C_{f,1} = \sup_{x \in [0,1]} \abs{2 f'(x) + \tfrac{1}{2}} $, by the Cauchy-Schwarz inequality,
		\begin{align*}
		\left\langle \partial_{xx} s, \left(2 f'(m_\eta) + \tfrac{1}{2}\right) s \right\rangle_{\alpha} \geq - C_{f,1} \| \partial_{xx} s \|_{2,\alpha} \| s \|_{2,\alpha}.
		\end{align*}
		Let $ C_{f,2} \geq 0 $ be such that $ \inf_{x\in \R} \{f'(m(x))^2 - \alpha(\partial_{x} + \alpha) f'(m)(x) - \frac{\alpha^2}{4} \} \geq - C_{f,2} $. We then have by~\eqref{eq:Aetas2alpha} that
		\begin{align*}
		\| \A_\eta s \|_{2,\alpha}^2 \geq \| \partial_{xx} s \|_{2,\alpha}^2 + \tfrac{1}{2} \| \partial_x s \|_{2,\alpha}^2 - C_{f,2} \| s \|_{2,\alpha}^2 - C_{f,1} \| \partial_{xx}s \|_{2,\alpha} \| s \|_{2,\alpha}.
		\end{align*}
		Note that
		\begin{align*}
		\tfrac{1}{2} \| \partial_{xx}s \|_{2,\alpha}^2 - C_{f,1} \| \partial_{xx}s \|_{2,\alpha} \| s \|_{2,\alpha} &= \tfrac{1}{2} \left( \| \partial_{xx}s \|_{2,\alpha} - C_{f,1} \| s \|_{2,\alpha} \right)^2 - \frac{C_{f,1}^2}{2} \| s \|_{2,\alpha}^2 \\
		&\geq - \frac{C_{f,1}^2}{2} \| s \|_{2,\alpha}^2.
		\end{align*}
		As a result, there exists a positive constant $ C_{f,3} > 0 $ such that
		\begin{align*}
		\| \A_\eta s \|_{2,\alpha}^2 \geq \tfrac{1}{2} \left( \| \partial_{xx}s \|_{2,\alpha}^2 + \| \partial_x s \|_{2,\alpha}^2 \right) - C_{f,3} \| s \|_{2,\alpha}^2.
		\end{align*}
		Together with \eqref{spectral_gap_L2}, this implies \eqref{energy_estimate_A_2} by the same argument as used to prove the lower bound in~\eqref{energy_estimate_A_1} from~\eqref{borne_inf_1} and~\eqref{borne_inf_2}.
	\end{proof}

	Let us now use Lemma~\ref{lemma:estimates_A} to prove Proposition~\ref{prop:energy_estimates}.

	\begin{proof}[Proof of Proposition~\ref{prop:energy_estimates}]
		Take $ v \in H^{1,\alpha} $ and $ \eta \in \R $ satisfying \eqref{v_eta_uniform} and set $ s := v-m_{\eta}$.
		Then, using \eqref{expansion_H} with $v=m_\eta$ and $h = s$, and using Lemma~\ref{lemma:H_is_zero_on_M} and~\eqref{def:m},
		\begin{align*}
		\mathcal{H}(v) = \frac{1}{2} \langle \A_{\eta}s, s \rangle_{\alpha} + \int_\R \mathcal{U}(m_\eta(x), s(x)) e^{\alpha x} dx.
		\end{align*}
		By the definition of $\mathcal U$ in~\eqref{def_Uvh} and by Taylor's theorem, there exists a constant $ C > 0 $ such that
		\begin{align*}
		\abs{\int_\R \mathcal{U}(m_\eta(x), s(x)) e^{\alpha x} dx} \leq C \| s \|_{{3,\alpha}}^3 \leq C \| s \|_{\infty} \| s \|_{2,\alpha}^2 \leq C \| s \|_\infty \| s \|_{H^{1,\alpha}}^2.
		\end{align*}
		As a result, from \eqref{energy_estimate_A_1} in Lemma~\ref{lemma:estimates_A}, and since by \eqref{v_eta_uniform} we have $ \langle s, \varphi_{\eta}  \rangle_{ \alpha}= 0 $,
		\begin{align*}
		\| s \|_{H^{1,\alpha}}^2 \left( \frac{\Cst{A_below}}{2} - C \| s \|_{\infty} \right) \leq \mathcal{H}(v) \leq \left( \frac{\Cst{A_above}}{2} + C \| s \|_{\infty} \right) \| s \|_{H^{1,\alpha}}^2.
		\end{align*}
		Hence there exist $ \Epsilon{H} > 0 $, $ \Cst{H_below} > 0 $ and $ \Cst{H_above} > 0 $ such that, if  $ \| s \|_{\infty} \leq \Epsilon{H} $,
		\begin{align*}
		\Cst{H_below} \| s \|_{H^{1,\alpha}}^2 \leq \mathcal{H}(v) \leq \Cst{H_above} \| s \|_{H^{1,\alpha}}^2.
		\end{align*}
		This proves \eqref{energy_estimate_H_1}.
		To prove \eqref{energy_estimate_H_2}, take $ v \in H^{2,\alpha} $ and use \eqref{def:m},~\eqref{def_DH} and~\eqref{def_A} to write
		\begin{equation} \label{DHexpression}
		D \mathcal{H}(v,\cdot) = \A_\eta s - \mathcal{V}(\cdot,s,\eta),
		\end{equation}
		where
		\begin{align*}
		\mathcal{V}(x,s,\eta) = f(v(x)) - f(m_\eta(x)) - f'(m_\eta(x)) s(x).
		\end{align*}
		Using Taylor's theorem, there exists a constant $ C > 0 $ such that
		\begin{align*}
		\| \mathcal{V}(\cdot,s,\eta) \|_{2,\alpha} \leq C \| s \|_{{4,\alpha}}^2 \leq C \| s \|_{\infty} \| s \|_{2,\alpha} \leq C \| s \|_\infty \| s \|_{H^{2,\alpha}}.
		\end{align*}
		Then, using~\eqref{DHexpression} and the Cauchy-Schwarz inequality and then using the right-hand side of~\eqref{energy_estimate_A_2} in the second line, we can write
		\begin{align*}
		\abs{ \| D \mathcal{H} (v,\cdot) \|_{2,\alpha}^2 - \| \A_\eta s \|_{2,\alpha}^2 } &\leq 2 \| \A_\eta s \|_{2,\alpha} \| \mathcal{V}(\cdot,s,\eta) \|_{2,\alpha} + \| \mathcal{V}(\cdot,s,\eta) \|_{2,\alpha}^2 \\
		&\leq C' \| s \|_{H^{2,\alpha}}^2 \left( \| s \|_\infty + \| s \|_\infty^2 \right),
		\end{align*}
		for some constant $ C' > 0 $.
		As a result, from \eqref{energy_estimate_A_2}, we can choose $ \Epsilon{H} > 0 $, $ \Cst{H_below} > 0 $ and $ \Cst{H_above} > 0 $ such that, if $ \| s \|_\infty \leq \Epsilon{H} $,
		\begin{align*}
		\Cst{H_below} \| s \|_{H^{2,\alpha}}^2 \leq \| D \mathcal{H}(v,\cdot) \|_{2,\alpha}^2 \leq \Cst{H_above} \| s \|_{H^{2,\alpha}}^2.
		\end{align*}
		This concludes the proof of Proposition~\ref{prop:energy_estimates}.
	\end{proof}

	\subsection{Stability of the manifold under the stochastic flow} \label{subsec:stability}
	
	We now turn to the proof of Proposition~\ref{prop:control_sigma} in the case $ i = \ref{sigmaN:dist} $. 
	In view of Proposition~\ref{prop:energy_estimates}, we would like to control $ \| s(v^N_t) \|_{2,\alpha} $ with the help of the energy functional $ \mathcal{H} $.
	The reason behind this is that by~\eqref{def_DH}, the SPDE \eqref{spde_vN} can be reformulated as
	\begin{align*}
		d v^N_t = - N D\mathcal{H}(v^N_t,\cdot) dt + \sqrt{v^N_t (1-v^N_t)} d W(t).
	\end{align*}
	Hence the drift term in this equation is pushing $ v^N_t $ towards the set of minimisers of $ \mathcal{H} $, i.e.~the manifold $ M $.
	
	However, one cannot hope to control $ \mathcal{H}(v^N_t) $ directly, because at positive times $t$, we have $ v^N_t \notin H^{1,\alpha} $.
	Instead, recalling \eqref{def_vdelta} and~\eqref{def:kappa_delta_r_theta}, we can use the smooth approximation $ v^{\delta_N, N}_t $ and control $ \mathcal{H}(v^{\delta_N,N}_t) $ to show that $ v^N_t $ stays in a suitable neighbourhood of $ M $.
	
	In order to apply Proposition~\ref{prop:energy_estimates} to $ v^{\delta_N,N}_{t \wedge \sigma_N} $, we need the following lemma.
	Recall the definitions of $\Epsilon{H}$ in Proposition~\ref{prop:energy_estimates} and $\kappa_N,\delta_N$ in~\eqref{def:kappa_delta_r_theta}.
	Recall that we could choose $\varepsilon_2$ in the statement of Proposition~\ref{prop:control_sigma}.
	From now on, we fix
	\begin{equation} \label{eq:eps2choice}
	\varepsilon_2\in (0,\varepsilon_1\wedge \varepsilon_3).
	\end{equation}
	
	\begin{lemma} \label{lemma:v_delta}
	Suppose $\alpha \in (0,2)$.
		For all $ N $ large enough, almost surely for all $t\ge 0$, $ v^{\delta_N,N}_{t \wedge \sigma_N} \in \mathcal{V}_{K,\Beta{eta}} $ and
		\begin{equation}\label{bound_s_v_delta}
			\dist(v^{\delta_N,N}_{t \wedge \sigma_N},M) \leq  \vartheta_N+r_N,
		\end{equation}
		and
		\begin{equation} \label{bound_eta_v_delta}
			\abs{ \eta(v^{\delta_N, N}_{t\wedge \sigma_N}) } < K + 1.
		\end{equation}
		In addition, for all $ N $ large enough, almost surely for all $ t \geq 0 $,
		\begin{align} \label{eq:s2infbounds}
			\| s(v^{\delta_N,N}_{t \wedge \sigma_N}) \|_\infty \leq \Epsilon{H}.
		\end{align}
		Moreover, there exists a constant $ \newCst{H2_v_delta} > 0 $ such that for all $ N $ large enough, almost surely for all $ t \geq 0 $, $ v^{\delta_N, N}_{t \wedge \sigma_N} \in H^{2,\alpha} $, 
		\begin{align} 
			\| s(v^{\delta_N,N}_{t \wedge \sigma_N}) \|_{H^{2,\alpha}} &\leq \Cst{H2_v_delta} (\delta_N+r_N+\delta_N^{-2}\vartheta_N),\label{bound_H2_v_delta}\\
			\text{and}\quad \| s(v^{\delta_N,N}_{t \wedge \sigma_N}) \|_{H^{1,\alpha}}
			&\le \Cst{H2_v_delta} (\delta_N+r_N+\delta_N^{-1}\dist(v^N_{t\wedge \sigma_N},M)). \label{bound_H1_v_delta}
		\end{align}
		
	\end{lemma}
	
	\begin{proof}[Proof of Lemma~\ref{lemma:v_delta}]
		We begin by proving the first part of the result.
		Note that for $v,v^\delta :\R\to \R$, if $ v \in \mathcal{V}_{K} $ then
		\begin{align}
		\| v^\delta - m_{\eta(v)} \|_{2,\alpha} \leq \dist(v,M) + \| v^\delta - v\|_{2,\alpha}. \label{dvdelta}
		\end{align}
		Recall from Lemma~\ref{lem:tauNpos} that $\sigma_N>0$ almost surely.
		Moreover, by~\eqref{def_vdelta} and \eqref{eq:Ddefn}, and then by~\eqref{def_sigma_vdelta}, for $t\ge 0$,
		\begin{equation} \label{vNtdeltadiff}
			\| v^N_{t \wedge \sigma_N} - v^{\delta_N,N}_{t \wedge \sigma_N} \|_{2,\alpha} \leq \| D^{\delta_N,N}(t \wedge \sigma_N,\cdot) \|_{2,\alpha} \leq r_N.
		\end{equation}
		Thus, by~\eqref{dvdelta} and \eqref{def_sigma_dist},
		\begin{align} \label{bound_L2_sv_delta}
			\| v^{\delta_N,N}_{t \wedge \sigma_N} - m_{\eta(v^N_{t \wedge \sigma_N})} \|_{2,\alpha} \leq \vartheta_N + r_N.
		\end{align}
		Hence if $ N $ is large enough that $ \vartheta_N + r_N < \Beta{eta} $, then $ v^{\delta_N,N}_{t \wedge \sigma_N} \in \mathcal{V}_{K} $ for any $t\ge 0$ (since $ | \eta(v^N_{t \wedge \sigma_N}) | \leq K $) and $ \dist(v^{\delta_N,N}_{t \wedge \sigma_N},M) \leq \vartheta_N + r_N $; this establishes~\eqref{bound_s_v_delta}.
		
		Then, by \eqref{bound_eta_eta_0} and \eqref{bound_L2_sv_delta}, we have
		\begin{equation*}
			| \eta(v^{\delta_N, N}_{t \wedge \sigma_N}) - \eta(v^N_{t\wedge \sigma_N}) | \leq \Cst{eta} (\vartheta_N + r_N).
		\end{equation*}
		This yields \eqref{bound_eta_v_delta} for $ N $ large enough.
		
		
		Next consider the following identity for $v,v^\delta \in \mathcal{V}_{K}$,
		\begin{align} \label{sv_svdelta}
		s(v^\delta) = s(v) + v^\delta - v - (m_{\eta(v^\delta)} - m_{\eta(v)}).
		\end{align}
		Using the fact that $ m $ is a Lipschitz function and using \eqref{eta_Lipschitz} (combined with~\eqref{def_sigma_dist},~\eqref{def_sigma_eta},~\eqref{bound_L2_sv_delta} and \eqref{bound_eta_v_delta}), for some constant $ C > 0 $, taking $ N $ large enough that $ \vartheta_N + r_N \leq \Beta{eta} $,
		\begin{align*}
			\| m_{\eta(v^{\delta_N,N}_{t \wedge \sigma_N})} - m_{\eta(v^N_{t\wedge \sigma_N})} \|_{\infty} \leq C \| v^{\delta_N,N}_{t\wedge \sigma_N} - v^N_{t\wedge\sigma_N} \|_{2,\alpha}.
		\end{align*}
		Hence, using~\eqref{sv_svdelta} and then \eqref{def_sigma_sup},~\eqref{vNtdeltadiff} and~\eqref{def_sigma_vdelta}, we obtain, for $ t \geq 0 $,
		\begin{align*}
		\| s(v^{\delta_N,N}_{t \wedge \sigma_N}) \|_{\infty} &\leq \| s(v^N_{t \wedge \sigma_N}) \|_\infty + \| v^{\delta_N,N}_{t \wedge \sigma_N} - v^N_{t \wedge \sigma_N} \|_\infty + C \| v^{\delta_N,N}_{t \wedge \sigma_N} - v^N_{t \wedge \sigma_N} \|_{2,\alpha} \\
		& \leq \Epsilon{proof} + (1+C)\, r_N.
		\end{align*}
		Thus since in~\eqref{eq:eps2choice} we chose $ \Epsilon{proof} < \Epsilon{H} $, we have $ \| s(v^{\delta_N,N}_{t \wedge \sigma_N}) \|_\infty \leq \Epsilon{H} $ for $ N $ large enough.
		
		Let us now turn to the bounds on $ \| s(v^{\delta_N,N}_{t \wedge \sigma_N}) \|_{H^{2,\alpha}} $ and $ \| s(v^{\delta_N,N}_{t \wedge \sigma_N}) \|_{H^{1,\alpha}} $.
		Note that since $m\in H^{2,\alpha}$, once we have established that $s(v^{\delta_N,N}_{t \wedge \sigma_N}) \in H^{2,\alpha} $, it will follow that $v^{\delta_N,N}_{t \wedge \sigma_N}\in H^{2,\alpha} $.
		Suppose $N$ is large enough that~\eqref{bound_s_v_delta} and~\eqref{bound_eta_v_delta} hold, and $\vartheta_N+r_N\le \beta_3$.
		Write $\delta:=\delta_N$.
		By~\eqref{def_vdelta}, for $t\le \sigma_N$ and $k\in \lbrace 0, 1, 2 \rbrace$,
		\begin{align} \label{deriv_v_delta}
		(\partial_{x})^k s(v^{\delta,N}_t) &= (\partial_{x})^k \rho^{\delta} \ast v^N_t - (\partial_{x})^k m_{\eta(v^{\delta,N}_t)} \notag \\
		&=(\partial_{x})^k \rho^{\delta} \ast (v^N_t -m_{\eta(v^{N}_t)}) +((\partial_{x})^k \rho^{\delta} \ast m_{\eta(v^{N}_t)} - (\partial_{x})^k m_{\eta(v^{\delta,N}_t)}).
		\end{align}
		For the second term on the right-hand side, for $x\in \R$, we can write
		\begin{align*}
		&\abs{(\partial_{x})^k \rho^{\delta} \ast m_{\eta(v^{N}_t)}(x) - (\partial_{x})^k m_{\eta(v^{\delta,N}_t)}(x)
		}\\
		&\quad = \abs{
		\int_{\R} \rho^{\delta}(y)((\partial_{x})^k m_{\eta(v^N_t)}(x-y)-(\partial_{x})^k m_{\eta(v^{\delta,N}_t)}(x))dy
		}\\
		&\quad \le 
		\int_{\R} \rho^{\delta}(y)\sup_{|z-x|\le K+1+\delta}|(\partial_{x})^{k+1}m(z)|(\delta+| \eta(v^N_t)-\eta(v^{\delta,N}_t)|)dy
		\\
		&\quad = 
		(\delta+| \eta(v^N_t)-\eta(v^{\delta,N}_t)|) \sup_{|z-x|\le K+1+\delta}|(\partial_{x})^{k+1}m(z)|,
		\end{align*}
		where in the inequality we used~\eqref{def_sigma_eta} and~\eqref{bound_eta_v_delta}, and the fact that $\rho^\delta$ is supported on $(-\rho,\rho)$ by~\eqref{eq:rhodeltadefn}.
		Then by~\eqref{def_sigma_dist},~\eqref{def_sigma_eta},~\eqref{bound_L2_sv_delta} and~\eqref{bound_eta_v_delta}, we can apply \eqref{eta_Lipschitz}, yielding
		\[
		| \eta(v^N_t)-\eta(v^{\delta,N}_t)|\le C \|v^N_t-v^{\delta,N}_t\|_{2,\alpha}\le Cr_N,
		\]
		where the second inequality follows by~\eqref{vNtdeltadiff}.
		Hence
		\begin{align} \label{eq:dkmbound}
		&\|(\partial_{x})^k \rho^{\delta} \ast m_{\eta(v^{N}_t)} - (\partial_{x})^k m_{\eta(v^{\delta,N}_t)}\|_{2,\alpha} \notag \\
		&\quad \le (\delta_N +Cr_N)\left( \int_{\R} \sup_{|z-x|\le K+1+\delta}|(\partial_{x})^{k+1}m(z)| ^2 e^{\alpha x}dx \right)^{1/2} \notag \\
		&\quad \le C' (\delta_N +Cr_N)
		\end{align}
		for some constant $C'>0$.
		For the first term on the right-hand side of~\eqref{deriv_v_delta}, by the Cauchy-Schwarz inequality,
		\begin{align*}
		&\| (\partial_{x})^k \rho^\delta \ast (v^N_t-m_{\eta(v^N_t)}) \|_{2,\alpha} \\
		&\leq \left( \int_\R \int_{\R} | (\partial_{x})^k \rho^\delta(x-y_1) | dy_1 \int_{\R} | (\partial_{x})^k \rho^\delta (x-y_2) | \, (v^N_t(y_2)-m_{\eta(v^N_t)}(y_2))^2 dy_2\, e^{\alpha x} dx \right)^{1/2}.
		\end{align*}
		Letting $z_1=x-y_1$ and $z_2=x-y_2$, the three integrals factorise, yielding
		\begin{multline*}
		\| (\partial_{x})^k \rho^\delta \ast (v^N_t-m_{\eta(v^N_t)}) \|_{2,\alpha} \\ \leq \left( \int_{\R} | (\partial_{x})^k \rho^\delta (z_1) | dz_1 \right)^{1/2}  \left( \int_{\R} | (\partial_{x})^k \rho^\delta (z_2) | e^{\alpha z_2} dz_2 \right)^{1/2} \\ \times \left( \int_{\R} (v^N_t(y_2)-m_{\eta(v^N_t)}(y_2))^2 e^{\alpha y_2} dy_2 \right)^{1/2}.
		\end{multline*}
		Since $ \rho^{\delta} $ is supported on $(-\rho,\rho)$, we can write $ e^{\alpha z_2} \leq e^{\alpha \delta} $ in the second integral.
		Moreover, by the definition of $\rho^\delta$ in~\eqref{eq:rhodeltadefn},
		\begin{equation*}
		\int_{\R} | (\partial_{x})^k \rho^\delta (x) | dx = \frac{1}{\delta^k} \| (\partial_{x})^k \rho \|_1,
		\end{equation*}
		and so
		\begin{align*}
		\| (\partial_{x})^k \rho^\delta \ast (v^N_t-m_{\eta(v^N_t)}) \|_{2,\alpha} \leq \frac{1}{\delta^k} e^{\frac{\alpha}{2} \delta} \| (\partial_{x})^k \rho \|_1 \|v^N_t-m_{\eta(v^N_t)}\|_{2,\alpha}.
		\end{align*}
By~\eqref{deriv_v_delta},~\eqref{eq:dkmbound} and~\eqref{def_sigma_dist}, this completes the proof.			
	\end{proof}
	
	Recall from \eqref{eq:eps2choice} that $ \Epsilon{proof} < \Epsilon{H} $.
	Lemma~\ref{lemma:v_delta} allows us to prove the following lemma.
	
	\begin{lemma} \label{lemma:bound_Hp}
	Suppose $\alpha \in (0,2)$.
		For $ \lambda > 0 $ and $t\ge 0$, set
		\begin{equation} \label{eq:FNlambdadefn}
			F^N_\lambda(t) := \exp\left( \lambda \mathcal{H}(v^{\delta_N,N}_t) \right)
		\end{equation}
		and let
		\begin{equation} \label{eq:dNdefn}
			d_N := r_N^2 + \frac{\kappa_N}{N \delta_N^3}.
		\end{equation}
		Then there exist constants $ \newCst{lambda} > 0 $ and $ \newCst{expH} >0 $ such that, for all $ \lambda > 0 $ such that
		\begin{equation} \label{upper_bound_lambda}
			\lambda \leq \Cst{lambda} \frac{N \delta_N}{\kappa_N},
		\end{equation}
		and for all $ t \geq 0 $,
		\begin{equation*}
			\E{ F^N_\lambda(t \wedge \sigma_N) } \leq \E{F^N_\lambda(0)} + \Cst{expH} \lambda \, d_N \left( N t + \E{F^N_\lambda(0)-1} \right) \exp\left( \Cst{expH} \, d_N \lambda \right).
		\end{equation*}
	\end{lemma}

	Before proving Lemma~\ref{lemma:bound_Hp}, let us see how it implies Proposition~\ref{prop:control_sigma} in the case $ i = 1 $.
	
	\begin{proof}[Proof of Proposition~\ref{prop:control_sigma} in the case $ i = \ref{sigmaN:dist} $]
	Suppose $\alpha \in (0,2)$.
		First note that, by~\eqref{def_sigma_dist}, and by~\eqref{bound_s_v_delta} in Lemma~\ref{lemma:v_delta} and the identity in~\eqref{sv_svdelta}, for $t\le \sigma_N$,
		\begin{align} \label{eq:s2alphalower}
			\| s(v^{\delta_N,N}_t) \|_{2,\alpha} \geq \| s(v^N_t) \|_{2,\alpha} - \| v^{\delta_N,N}_t - v^N_t \|_{2,\alpha} - \| m_{\eta(v^{\delta_N,N}_t)} - m_{\eta(v^N_t)} \|_{2,\alpha}.
		\end{align}
		Then taking $ N $ large enough that $ \vartheta_N + r_N \leq \Beta{eta} $ and applying Lemma~\ref{lemma:diff_m_eta}, Lemma~\ref{lemma:v_delta}, \eqref{def_sigma_dist}, \eqref{def_sigma_eta}, and \eqref{eta_Lipschitz}, we see that there exists a constant $ C > 0 $ such that for all $t\ge 0$,
		\begin{align} \label{bound_below_s_vdelta}
			\| s(v^{\delta_N,N}_{t \wedge \sigma_N}) \|_{2,\alpha} &\geq \| s(v^N_{t \wedge \sigma_N}) \|_{2,\alpha} - C \|v^{\delta_N,N}_{t\wedge \sigma_N}-v^{N}_{t\wedge \sigma_N}\|_{2,\alpha} \notag \\
			&\ge \| s(v^N_{t \wedge \sigma_N}) \|_{2,\alpha} - C r_N,
		\end{align}
		where the second inequality follows by~\eqref{vNtdeltadiff}.
		In addition, by the definitions of $ \sigmaN{dist} $ and $\sigmaN{eta}$ in~\eqref{def_sigma_dist} and~\eqref{def_sigma_eta}, and since $\sigma_N>0$ almost surely (by Lemma~\ref{lem:tauNpos}), for $t>0$,
		\begin{align*}
			\{\sigmaN{dist} \leq t \wedge \hat{\sigma}_{\ref{sigmaN:dist},N}\} &\subseteq \{ \| s(v^{N}_{t \wedge \sigma_N}) \|_{2,\alpha} \geq \vartheta_N\}.
		\end{align*}
		By \eqref{bound_below_s_vdelta}, this implies
		\begin{align*}
			\{\sigmaN{dist} \leq t \wedge \hat{\sigma}_{\ref{sigmaN:dist},N}\} &\subseteq \{ \| s(v^{\delta_N,N}_{t \wedge \sigma_N}) \|_{2,\alpha} \geq \vartheta_N - Cr_N\},
		\end{align*}
		and hence by Proposition~\ref{prop:energy_estimates} and Lemma~\ref{lemma:v_delta},
		\begin{align*}
			\{\sigmaN{dist} \leq t \wedge \hat{\sigma}_{\ref{sigmaN:dist},N}\} &\subseteq \{ \mathcal{H}(v^{\delta_N,N}_{t \wedge \sigma_N}) \geq \Cst{H_below} (\vartheta_N - C r_N)^2\}.
		\end{align*}
		As a result, for any $ \lambda > 0 $ and $ T > 0 $,
		by the definition of $F^N_\lambda$ in~\eqref{eq:FNlambdadefn} and by Markov's inequality,
		\begin{equation} \label{eq:sigma1Nmarkov}
			\P{ \sigmaN{dist} \leq T \wedge \hat{\sigma}_{\ref{sigmaN:dist},N} } \leq \frac{\E{F^N_\lambda(T \wedge \sigma_N)}}{\exp\left( \lambda \Cst{H_below} (\vartheta_N - C r_N)^2 \right)}.
		\end{equation}
		By \eqref{def:kappa_delta_r_theta}, for $ N $ large enough we have
		\begin{equation*}
			\frac{N \delta_N}{\kappa_N} \vartheta_N^2 = N^{\cexp{KN} - \cexp{deltaN} - 2 \cexp{thetaN}}.
		\end{equation*}
		Since, by~\eqref{bounds_cexpthetaN}, $ \cexp{thetaN} < \frac{\cexp{KN}- 3 \cexp{deltaN}}{2} $, we have $ \cexp{KN}-\cexp{deltaN}-2\cexp{thetaN} > 2 \cexp{deltaN} > 0 $.
		As a result, taking $ \lambda = \lambda_N = \Cst{lambda} \frac{N \delta_N}{\kappa_N} $, we have
		that there exists $c>0$ such that for $N$ sufficiently large,
		\begin{equation} \label{eq:lambdatheta2}
			\lambda_N \vartheta_N^2 =\Cst{lambda} N^c.
		\end{equation}
		By Lemma~\ref{lemma:bound_Hp}, it follows from~\eqref{eq:sigma1Nmarkov} with $\lambda=\lambda_N$ that
		\begin{align} \label{bound_proba_sigma1}
			&\P{ \sigmaN{dist} \leq T \wedge \hat{\sigma}_{\ref{sigmaN:dist},N} }\notag \\
			 &\quad \leq \mathbb E \Bigg[ \exp\left( \lambda_N \left( \mathcal{H}(v^{\delta_N,N}_0) - \Cst{H_below} (\vartheta_N - C r_N)^2 \right) \right)  \notag \\&\qquad   + \Cst{expH} \lambda_N \, d_N \left( N T + F^N_{\lambda_N}(0) - 1 \right) \exp\left( - \lambda_N \vartheta_N^2 \left( \Cst{H_below} \left(1 - C \frac{r_N}{\vartheta_N}\right)^2 - \Cst{expH} \frac{d_N}{\vartheta_N^2} \right) \right) \Bigg].
		\end{align}
		By~\eqref{eq:dNdefn} and then by \eqref{def:kappa_delta_r_theta}, for $ N $ large enough,
		\begin{equation} \label{eq:dnthetaN-2}
			\frac{d_N}{\vartheta_N^2} = \frac{1}{\vartheta_N^2} \left( r_N^2 + \frac{\kappa_N}{N \delta_N^3} \right) \leq 2 N^{2 \cexp{thetaN} -(2 \cexp{rN} \wedge (\cexp{KN}-3\cexp{deltaN}))}.
		\end{equation}
		Since, by~\eqref{bounds_cexpthetaN}, we have $ \cexp{thetaN} < \cexp{rN} \wedge \frac{\cexp{KN}-3\cexp{deltaN}}{2} $, the right-hand side of~\eqref{eq:dnthetaN-2} vanishes as $ N \to \infty $.
		In addition, by~\eqref{def:kappa_delta_r_theta} again, for $ N $ large enough,
		\begin{equation*}
			\frac{r_N}{\vartheta_N} = N^{\cexp{thetaN}-\cexp{rN}},
		\end{equation*}
		which also vanishes as $ N \to \infty $ since $ \cexp{thetaN} < \cexp{rN} $ by~\eqref{bounds_cexpthetaN}.
		It follows that, for all $ N $ large enough,
		\begin{equation*}
			\Cst{H_below} \left(1 - C \frac{r_N}{\vartheta_N}\right)^2 - \Cst{expH} \frac{d_N}{\vartheta_N^2} \geq \frac{\Cst{H_below}}{2}.
		\end{equation*}
		Coming back to \eqref{bound_proba_sigma1}, we have established that for all $ N $ large enough,
		\begin{multline} \label{eq:probsigma1N}
			\P{ \sigmaN{dist} \leq T \wedge \hat{\sigma}_{\ref{sigmaN:dist},N} } \leq \mathbb E \Bigg[ \exp\left( - \lambda_N \vartheta_N^2 \left( \frac{\Cst{H_below}}{2} - \frac{\mathcal{H}(v^{\delta_N,N}_0)}{\vartheta_N^2} \right) \right) \\ + \Cst{expH} \lambda_N d_N \left( NT + \exp\left( \lambda_N \mathcal{H}(v^{\delta_N,N}_0) \right) \right) \exp\left( - \frac{\Cst{H_below}}{2} \lambda_N \vartheta_N^2 \right) \Bigg].
		\end{multline}

		It remains to bound $\mathcal H(v_0^{\delta_N,N})$.
		By Proposition~\ref{prop:energy_estimates} and Lemma~\ref{lemma:v_delta}, and then by~\eqref{bound_H1_v_delta} in Lemma~\ref{lemma:v_delta},
		\begin{align*}
			\mathcal{H}(v^{\delta_N,N}_0) &\leq \Cst{H_above} \| s(v^{\delta_N,N}_0) \|_{H^{1,\alpha}}^2 
			\le \Cst{H_above}\Cst{H2_v_delta}^2 (\delta_N+r_N+\delta_N^{-1}\dist(v_0^N,M))^2.
		\end{align*}
		By Assumption~\ref{assumpt:v0}.\ref{v0:dist} and \eqref{def:kappa_delta_r_theta}, it follows that for $N$ sufficiently large,
		\begin{align*} 
			\vartheta_N^{-2}\mathcal{H}(v^{\delta_N,N}_0) \le  \Cst{H_above} \Cst{H2_v_delta}^{2}N^{2\cexp{thetaN}} (N^{-\cexp{deltaN}}+N^{-\cexp{rN}} + N^{\cexp{deltaN}-\cinit{big}})^{2}.
		\end{align*}
		By \eqref{bounds_cexpthetaN} we see that
		\begin{equation*}
		\cexp{thetaN}<\cexp{deltaN}\wedge \cexp{rN} \quad \text{and} \quad \cexp{deltaN}+\cexp{thetaN}<\cinit{big},
		\end{equation*}
		and hence there exists $c'>0$ such that for $N$ sufficiently large,
		\begin{align} \label{to0_2}
			\vartheta_N^{-2}\mathcal{H}(v^{\delta_N,N}_0) \le N^{-c'}.
		\end{align}
		Take $N$ sufficiently large that~\eqref{bound_proba_sigma1} and~\eqref{eq:lambdatheta2} hold, and that $N^{-c'}\le \frac 14 \Cst{H_below}$.
		Then by~\eqref{to0_2} and~\eqref{bound_proba_sigma1}, we have 
		\begin{align*}
			&\P{ \sigmaN{dist} \leq T \wedge \hat{\sigma}_{\ref{sigmaN:dist},N} }\notag \\
			 &\quad \leq \exp\left( -\lambda_N \vartheta_N^2 \cdot \tfrac 14 \Cst{H_below}  \right)  + \Cst{expH} \lambda_N \, d_N \left( N T + \exp\left( \lambda_N \vartheta_N^2 \cdot \tfrac 14 \Cst{H_below}  \right)\right) \exp\left( -\lambda_N \vartheta_N^2 \cdot \tfrac 12 \Cst{H_below}  \right)\\
			 &\quad = \exp\left( -\tfrac 14 \Cst{H_below} \Cst{lambda}N^c \right)  + \Cst{expH} \lambda_N \, d_N \left( N T\exp\left( -\tfrac 12 \Cst{H_below} \Cst{lambda}N^c \right)+\exp\left( -\tfrac 14 \Cst{H_below} \Cst{lambda}N^c \right) \right),
		\end{align*}
		where the last line follows by~\eqref{eq:lambdatheta2}.
		Since $c>0$, and by the definitions of $d_N$ in~\eqref{eq:dNdefn} and $\lambda_N$ before~\eqref{eq:lambdatheta2}, it follows that
		\begin{align*}
			\lim_{N \to \infty} \P{\sigmaN{dist} \leq T \wedge \hat{\sigma}_{\ref{sigmaN:dist},N}} = 0,
		\end{align*}
		which completes the proof of Proposition~\ref{prop:control_sigma} in the case $ i = \ref{sigmaN:dist} $.
	\end{proof}

	In the remainder of this subsection, we will prove Lemma~\ref{lemma:bound_Hp}; some of the ideas will be similar to those in Section~5 in \cite{funaki_scaling_1995}.
	Let us abuse notation by letting $ W(dx dt) $ denote the martingale measure defined by
	\begin{equation} \label{def_martingale_measure}
		\int_{[0,t] \times \R } \phi(x) W(dx ds) = \langle W(t), \phi \rangle_{L^2},
	\end{equation}
	where $(W(t),t\ge 0)$ is a cylindrical Wiener process on $L^2(\R)$ as in Definition~\ref{def:u}.
	Note that, by \eqref{spde_vN} and \eqref{def_vdelta}, for $\delta>0$ and $t>0$,
	\begin{align} \label{eq:dvdelta}
		d v^{\delta,N}_t = N \left( \partial_{xx} v^{\delta,N}_t + \alpha \partial_x v^{\delta,N}_t + \rho^\delta \ast f(v^N_t) \right) dt + \int_\R \rho^\delta(\cdot - y) \sqrt{v^N_t(y) (1-v^N_t(y)) } W(dy dt).
	\end{align}
	Moreover, by the integration by parts formula and \eqref{IPP_alpha},
	\begin{multline*} 
		d \langle \partial_{x} v^{\delta,N}_t, \partial_{x} v^{\delta,N}_t \rangle_{\alpha} = -2 \langle (\partial_{xx} + \alpha \partial_x) v^{\delta,N}_t, d v^{\delta,N}_t \rangle_{\alpha} \\ - \int_\R \langle (\partial_{xx} + \alpha \partial_x) \rho^\delta(\cdot-y), \rho^\delta(\cdot-y) \rangle_{\alpha} v^N_t(y) (1-v^N_t(y)) dy dt.
	\end{multline*}
	Also, by Itô's lemma, for $\delta>0$, $t>0$ and $x\in \R$,
	\begin{align*}
		d F(v^{\delta,N}_t(x)) = F'(v^{\delta,N}_t(x)) d v^{\delta,N}_t(x) + \frac{1}{2} F''(v^{\delta,N}_t(x)) d \left[ v^{\delta,N}(x) \right]_t.
	\end{align*}
	Hence by~\eqref{def_F} and~\eqref{eq:dvdelta},
	\begin{multline*}
		d \int_\R F(v^{\delta,N}_t(x)) e^{\alpha x} dx = - \langle f(v^{\delta,N}_t), d v^{\delta,N}_t \rangle_{\alpha} \\ - \frac{1}{2} \int_\R \left\langle f'(v^{\delta,N}_t) \rho^\delta(\cdot-y), \rho^\delta(\cdot-y) \right\rangle_{\alpha} v^N_t(y) (1-v^N_t(y)) dy dt.
	\end{multline*}
	As a result, in view of the definition of $\mathcal{H}$ in~\eqref{eq:Hdefn} and using \eqref{def_DH}, for $\delta>0$ and $t>0$,
	\begin{align*}
		d \mathcal{H}(v^{\delta,N}_t) &= \langle D \mathcal{H}(v^{\delta,N}_t,\cdot),dv_t^{\delta,N}  \rangle_{\alpha}  + A^{\delta,N}(t) dt\\
		&= N \langle D \mathcal{H}(v^{\delta,N}_t,\cdot), (\partial_{xx} + \alpha \partial_x) v^{\delta,N}_t + \rho^\delta \ast f(v^N_t) \rangle_{\alpha} dt + A^{\delta,N}(t) dt + d M^{\delta,N}(t),
	\end{align*}
	where
	\begin{align} \label{eq:Adeltadefn}
		A^{\delta,N}(t) = -\frac{1}{2} \int_\R \left\langle (\partial_{xx} + \alpha \partial_x + f'(v^{\delta,N}_t)) \rho^\delta(\cdot-y), \rho^\delta(\cdot-y) \right\rangle_{\alpha} v^N_t(y) (1-v^N_t(y)) dy
	\end{align}
	and
	\begin{align} \label{eq:Mdeltadefn}
		M^{\delta,N}(t) = \int_0^t \int_\R \langle D \mathcal{H}(v^{\delta,N}_s,\cdot), \rho^\delta(\cdot-y) \rangle_{\alpha} \sqrt{v^N_s(y) (1-v^N_s(y))} W(dy ds).
	\end{align}
	Setting
	\begin{align} \label{def_R}
		R^{\delta,N}(t) := \rho^\delta \ast f(v^N_t) - f(v^{\delta,N}_t),
	\end{align}
	we can write this as
	\begin{align} \label{dH(vdelta)}
		d \mathcal{H}(v^{\delta,N}_t) = - N \| D \mathcal{H}(v^{\delta,N}_t) \|_{2,\alpha}^2 dt + N \langle D \mathcal{H}(v^{\delta,N}_t), R^{\delta,N}(t) \rangle_{\alpha} dt + A^{\delta,N}(t) dt + d M^{\delta,N}(t).
	\end{align}
	From Proposition~\ref{prop:energy_estimates}, we see that taking $\delta=\delta_N$, if $\mathcal{H}(v^{\delta,N}_t)$ grows, then the first term on the right-hand side of \eqref{dH(vdelta)} grows in absolute value.
	We will thus be able to use~\eqref{dH(vdelta)} to show that $ \mathcal{H}(v^{\delta,N}_t) $ cannot grow too much, controlling the other terms on the right-hand side.
	
	Recall from~\eqref{eq:FNlambdadefn} that for $\lambda>0$ and $t\ge 0$, we set
	\begin{equation*}
		F^N_\lambda(t) = \exp\left( \lambda \mathcal{H}(v^{\delta_N,N}_t) \right).
	\end{equation*}
	By Itô's lemma,
	setting $\delta=\delta_N$, we have
	\begin{align*}
		d F^N_\lambda(t) = \lambda F^N_\lambda(t) d \mathcal{H}(v^{\delta, N}_t) + \tfrac{1}{2}\lambda^2 F^N_\lambda(t)  d \left[ \mathcal{H}(v^{\delta, N}_\cdot) \right]_t.
	\end{align*}
	If we set
	\begin{align} \label{def_BN}
		B^{\delta, N}(t) = \int_\R \langle D \mathcal{H}(v^{\delta, N}_t,\cdot), \rho^\delta(\cdot
	-y) \rangle_{ \alpha}^2 v^N_t(y) (1-v^N_t(y)) dy,
	\end{align}
	then by~\eqref{dH(vdelta)} we can write
	\begin{multline} \label{dFN_lambda}
		d F^N_\lambda(t) = - \lambda N  \| D\mathcal{H}(v^{\delta,N}_t) \|_{2,\alpha}^2 F^N_\lambda(t) dt + \lambda N \langle D\mathcal{H}(v^{\delta,N}_t), R^{\delta,N}(t) \rangle_\alpha F^N_\lambda(t) dt \\+ \lambda A^{\delta,N}(t) F^N_\lambda(t) dt + \lambda F^N_\lambda(t) d M^{\delta,N}(t) + \frac{\lambda^2}{2} B^{\delta,N}(t) F^N_\lambda(t) dt.
	\end{multline}
	The following lemma will allow us to prove Lemma~\ref{lemma:bound_Hp}.
	
	\begin{lemma} \label{lemma:true_martingale}
		Suppose $\alpha \in (0,2)$.
		For $p\ge 1$, the process
		\begin{align} \label{martingale}
			t \mapsto \int_{0}^{t \wedge \sigma_N} F^N_\lambda(s) d M^{\delta_N, N}(s)
		\end{align}
		is a square integrable martingale with respect to the natural filtration of $(v^N_t)_{t\geq 0}$.
		Moreover, there exists a constant $ \newCst{AB} > 0 $ such that, for all $ t \leq \sigma_N $,
		\begin{align} \label{bound_A}
			\abs{A^{\delta_N,N}(t)} \leq \Cst{AB} \frac{\kappa_N}{\delta_N^3} && \text{ and } && \abs{B^{\delta_N, N}(t)} \leq \Cst{AB} \frac{\kappa_N}{\delta_N} \| D \mathcal{H}(v^{\delta_N, N}_t) \|_{2,\alpha}^2.
		\end{align}
	\end{lemma}

	We will prove Lemma~\ref{lemma:true_martingale} at the end of this subsection.
	Let us now prove Lemma~\ref{lemma:bound_Hp}.
	
	\begin{proof}[Proof of Lemma~\ref{lemma:bound_Hp}]
		Fix $t\ge 0$. To simplify notation, let $\delta:=\delta_N$ and $\sigma:=\sigma_N$.
		For $\lambda>0$, integrating \eqref{dFN_lambda} up to time $ t \wedge \sigma_N $ and taking expectations on both sides, by~\eqref{martingale} in Lemma~\ref{lemma:true_martingale} we obtain
		\begin{multline*} 
			\E{ F^N_\lambda(t \wedge \sigma) } + \lambda N \E{ \int_{0}^{t \wedge \sigma} \| D\mathcal{H}(v^{\delta,N}_s) \|_{2,\alpha}^2 F^N_\lambda(s) ds } \\= \E{F^N_\lambda(0)} + \lambda N \E{ \int_{0}^{t \wedge \sigma} \langle D\mathcal{H}(v^{\delta,N}_s), R^{\delta,N}(s) \rangle_\alpha F^N_\lambda(s) ds } \\+ \lambda \E{ \int_{0}^{t \wedge \sigma} A^{\delta,N}(s) F^N_\lambda(s) ds } + \frac{\lambda^2}{2} \E{ \int_{0}^{t \wedge \sigma} B^{\delta,N}(s) F^N_\lambda(s) ds }.
		\end{multline*}
		Now we use that for $s\le \sigma_N$, since $ab\le \frac 12 (a^2+b^2)$ $\forall a,b\in \R$,
		\begin{align*}
			\langle D \mathcal{H}(v^{\delta, N}_s), R^{\delta, N}(s) \rangle_{\alpha} \leq \frac{1}{2} \left( \| D \mathcal{H}(v^{\delta,N}_s) \|_{2,\alpha}^2 + \| R^{\delta,N}(s) \|_{2,\alpha}^2 \right),
		\end{align*}
		and we note that by~\eqref{def_R} and then by~\eqref{eq:Ddefn}, for $ s \leq \sigma_N $,
		\begin{align*}
			\| R^{\delta,N}(s) \|_{2,\alpha} &\leq \| \rho^\delta \ast f(v^N_s) - f(v^N_s) \|_{2,\alpha} + \| f(v^N_s) - f(\rho^\delta \ast v^N_s) \|_{2,\alpha} \\
			&\leq 2\, \sup_{x\in [0,1]} |f'(x)|\, \| D^{\delta,N}(s,\cdot) \|_{2,\alpha} \\
			&\leq 2\, \sup_{x\in [0,1]} |f'(x)|\, r_N,
		\end{align*}
		where the last inequality follows by~\eqref{def_sigma_vdelta}.
		Together with~\eqref{bound_A} in Lemma~\ref{lemma:true_martingale}, and since $F^N_\lambda(s)\ge 0$, this yields
		\begin{multline} \label{bound_E_FNlambda}
			\E{ F^N_\lambda(t\wedge \sigma) } + \frac{\lambda N}{2} \left( 1 - \Cst{AB} \frac{\lambda \kappa_N}{N \delta_N} \right) \E{ \int_{0}^{t \wedge \sigma} \| D\mathcal{H}(v^{\delta,N}_s) \|_{2,\alpha}^2 F^N_\lambda(s) ds } \\ \leq \E{F^N_\lambda(0)} + C \lambda \left( N r_N^2 + \frac{\kappa_N}{\delta_N^3} \right) \E{ \int_{0}^{t\wedge \sigma} F^N_\lambda(s) ds },
		\end{multline}
		where $C:=\Cst{AB} \vee (2\sup_{x\in [0,1]}|f'(x)|^2)>0$ is a constant.
		For $\lambda\ge 0$, we set
		\begin{equation*}
			A_N(\lambda) := \E{ \int_{0}^{t \wedge \sigma} F^N_\lambda(s)ds }.
		\end{equation*}
		By Proposition~\ref{prop:energy_estimates} and Lemma~\ref{lemma:v_delta},
there exists deterministic $K_N>0$ such that
$\mathcal H(v^{\delta_N,N}_{s\wedge \sigma_N})\in [0,K_N]$ $\forall s\ge 0$, and so in particular,		
		 $ \lambda \mapsto A_N(\lambda) $ is locally bounded.
		Moreover, 
		for $h\in \R$ and $\lambda>0$ with $\lambda+h>0$, for $s\le \sigma_N$ we have
		\begin{align*}
		\left|\frac 1 h (F^N_{\lambda+h}(s)-F^N_\lambda(s))\right| 
		&= \exp(\lambda \mathcal H(v^{\delta_N,N}_s))\left|\frac 1h (\exp(h\mathcal H(v^{\delta_N,N}_s))-1)\right|\\
		&\le \exp(\lambda K_N)\cdot \exp (|h| K_N)K_N.
		\end{align*}
		Therefore, by dominated convergence, $A_N(\lambda)$ is differentiable, and
		\begin{align} \label{eq:AN'}
			A_N'(\lambda) &= \E{ \int_{0}^{t \wedge \sigma} \mathcal{H}(v^{\delta,N}_s) F^N_\lambda(s) ds }\notag  \\
			&\leq \Cst{H_above}\Cst{H_below}^{-1} \E{ \int_{0}^{t\wedge \sigma} \| D\mathcal{H}(v^{\delta,N}_s) \|_{2,\alpha}^2 F^N_\lambda(s) ds },
		\end{align}
		where the inequality follows by Proposition~\ref{prop:energy_estimates} and Lemma~\ref{lemma:v_delta}.
		Set $ \Cst{lambda} = (2\Cst{AB})^{-1}>0 $ and take $ \lambda \leq \Cst{lambda} \frac{N \delta_N}{\kappa_N} $.
		Then $1 - \Cst{AB} \frac{\lambda \kappa_N}{N \delta_N}\ge \frac 12$, and so by~\eqref{bound_E_FNlambda} and~\eqref{eq:AN'} we have
		\begin{align*}
			A_N'(\lambda) &\leq \Cst{H_above}\Cst{H_below}^{-1}\cdot \frac 4 {\lambda N} \left( \E{F^N_\lambda(0)- F^N_\lambda(t\wedge \sigma) } + C \lambda \left( N r_N^2 + \frac{\kappa_N}{\delta_N^3}\right) A_N(\lambda)\right)\\
			&\leq 4 \Cst{H_above}\Cst{H_below}^{-1} \left( \frac{\E{F^N_\lambda(0)-1}}{\lambda N} + C \, d_N A_N(\lambda) \right),
		\end{align*}
		where the second inequality follows since $F^N_\lambda(t\wedge \sigma)\ge 1$ and by the definition of $d_N$ in~\eqref{eq:dNdefn}.
		Integrating this inequality, and using the fact that $ A_N(0) \leq t $, we obtain
		that for $ \lambda \leq \Cst{lambda} \frac{N \delta_N}{\kappa_N} $,
		\begin{equation*}
			A_N(\lambda) \leq t + \frac{4 \Cst{H_above}\Cst{H_below}^{-1}}{N} \int_{0}^{\lambda} \frac{\E{F^N_u(0)-1}}{u} du + 4\Cst{H_above}\Cst{H_below}^{-1}C \, d_N \int_{0}^{\lambda} A_N(u) du.
		\end{equation*}
		Let $\tilde C:=4\Cst{H_above}\Cst{H_below}^{-1}C$; then
		by Gr\"onwall's inequality (using that $F^N_u(0)\ge 1$ $\forall u\ge 0$), this yields that for $ \lambda \leq \Cst{lambda} \frac{N \delta_N}{\kappa_N} $,
		\begin{equation*}
			A_N(\lambda) \leq \left( t + \frac{4 \Cst{H_above}\Cst{H_below}^{-1}}{N} \int_{0}^{\lambda} \frac{\E{F^N_u(0)-1}}{u} du \right) \exp\left( \tilde C \,	d_N \lambda \right).
		\end{equation*}
		Now note that for $a\ge 0$ fixed, the function $u\mapsto \frac 1 u (e^{au}-1)$ is non-decreasing on $(0,\infty)$, and so
		\[
		u\mapsto \frac{\E{F^N_u(0)-1}}{u}=\frac 1u \left(\exp(u\mathcal H(v_0^{\delta_N,N}))-1\right)
		\]
		is non-decreasing on $(0,\infty)$. Therefore we can also write
		\begin{equation*}
			A_N(\lambda) \leq \left( t + \frac{4 \Cst{H_above}\Cst{H_below}^{-1}}{N} \E{F^N_\lambda(0)-1} \right) \exp\left( \tilde C \, d_N \lambda \right).
		\end{equation*}
		Plugging this into the right-hand side of \eqref{bound_E_FNlambda}, and using the fact that the second term on the left-hand side of~\eqref{bound_E_FNlambda} is non-negative, along with the definition of $d_N$ in~\eqref{eq:dNdefn},  we obtain
		\begin{align*}
			\E{ F^N_\lambda(t \wedge \sigma) } &\leq \E{F^N_\lambda(0)} + C \lambda N \, d_N A_N(\lambda) \\
			&\leq \E{F^N_\lambda(0)} + C \lambda \, d_N \left( N t + 4 \Cst{H_above}\Cst{H_below}^{-1}\E{F^N_\lambda(0)-1} \right) \exp\left( \tilde C \, d_N \lambda \right),
		\end{align*}
		which completes the proof, since $F^N_\lambda(0)\ge 1$ and we can take $\Cst{expH}=\tilde C \vee C \vee (4 \Cst{H_above}\Cst{H_below}^{-1}C)$. 
	\end{proof}	
	We now finish this section by proving Lemma~\ref{lemma:true_martingale}.
	
	\begin{proof}[Proof of Lemma~\ref{lemma:true_martingale}]
	To simplify notation, write $\delta :=\delta_N$.
		We start with the second part of the statement.
		Recall the definition of $A^{\delta,N}(t)$ in~\eqref{eq:Adeltadefn}.
		Note that, for all $y \in \R$ and $t\ge 0$, using \eqref{IPP_alpha},
		\begin{align*}
		\abs{ \langle -(\partial_{xx} + \alpha \partial_x + f'(v^{\delta,N}_t)) \rho^\delta(\cdot - y), \rho^\delta(\cdot-y) \rangle_{\alpha} } &\leq \| \partial_{x} \rho^\delta(\cdot-y) \|_{2,\alpha}^2 + \sup_{x\in [0,1]}| f' (x)| \| \rho^\delta(\cdot-y) \|_{2,\alpha}^2 \\
		&\leq \frac{C}{\delta^3} e^{\alpha y},
		\end{align*}
		for some constant $ C > 0 $, where the second inequality follows by the definition of $\rho^\delta$ in~\eqref{eq:rhodeltadefn}.
		Hence, by \eqref{def_sigma_tail} and the fact that $ v^N_t \in [0,1] $, there exists a constant $ C > 0 $ such that, for all $ t \leq \sigma_N $,
		\begin{align*}
		|A^{\delta,N}(t)| \leq C \frac{\kappa_N}{\delta_N^3}.
		\end{align*}
		For the bound on $ B^{\delta, N}(t) $, note that for $s\le \sigma_N$ and $y\in \R$, by the Cauchy-Schwarz inequality,
		\begin{align*}
		\langle D \mathcal{H}(v^{\delta,N}_s,\cdot ), \rho^\delta(\cdot-y) \rangle_{\alpha}^2 &\leq \| D \mathcal{H}(v^{\delta,N}_s) \|_{2,\alpha}^2 \| \rho^\delta(\cdot -y) \|_{2,\alpha}^2 \\
		&\leq \frac{C}{\delta} e^{\alpha y} \| D \mathcal{H}(v^{\delta,N}_s) \|_{2,\alpha}^2
		\end{align*}
		for some constant $ C > 0 $.
		Plugging this into \eqref{def_BN} and using \eqref{def_sigma_tail}, we obtain, for $ t \leq \sigma_N $,
		\begin{align} \label{eq:Bdeltabound}
		|B^{\delta, N}(t)| \leq C \dfrac{\kappa_N}{\delta_N} \| D \mathcal{H}(v^{\delta,N}_t) \|_{2,\alpha}^2.
		\end{align}
		It remains to show that \eqref{martingale} is a square integrable martingale. 
		For this, we note that by the definition of $M^{\delta,N}(t)$ in~\eqref{eq:Mdeltadefn}, and then by~\eqref{eq:Bdeltabound} and since $F^N_\lambda(s)\ge 0$ for $s\le \sigma_N$ (by~\eqref{eq:FNlambdadefn}),
		\begin{align*}
		\left[ \int_{0}^{\cdot \wedge \sigma_N} F^N_\lambda(s) d M^{\delta,N}(s) \right]_t &= \int_0^{t\wedge \sigma_N} F^N_\lambda(s)^2 B^{\delta, N}(s) ds \\
		&\leq C\frac{\kappa_N}{\delta_N} \int_{0}^{t \wedge \sigma_N} F^N_\lambda(s)^2 \| D \mathcal{H}(v^{\delta,N}_s) \|_{2,\alpha}^2 ds. \numberthis \label{bound_qvar}
		\end{align*}
		By~\eqref{eq:FNlambdadefn} and then by Proposition~\ref{prop:energy_estimates} and Lemma~\ref{lemma:v_delta}, for $ s \leq \sigma_N $,
		\begin{align} \label{bound_H2}
			F^N_\lambda(s)^2 \| D \mathcal{H}(v^{\delta,N}_s) \|_{2,\alpha}^2 
= e^{2\lambda \mathcal{H}(v^{\delta, N}_s)}	\| D \mathcal{H}(v^{\delta,N}_s) \|_{2,\alpha}^2 		
			\leq e^{2\lambda \Cst{H_above} \| s(v^{\delta,N}_s) \|_{H^{2,\alpha}}^2}\Cst{H_above} \| s(v^{\delta,N}_s) \|_{H^{2,\alpha}}^2.
		\end{align}
		Therefore, using~\eqref{bound_H2_v_delta} in Lemma~\ref{lemma:v_delta}, there exists a deterministic $ K_N > 0 $ depending on $ N $ such that
		\begin{align*}
		\left[ \int_{0}^{\cdot \wedge \sigma_N} \mathcal{H}^{p-1}(v^{\delta, N}_s) d M^{\delta,N}(s) \right]_t \leq K_N (t \wedge \sigma_N)
		\end{align*}
		almost surely, and the proof is complete.
	\end{proof}

	\section{The right interface} \label{sec:interface}

	This section is devoted to the analysis of the position of the right interface of the solution, i.e. the right endpoint of its support, and a related process related to the integrability of the solution.
	The upper bounds obtained here for these quantities will be used in Section~\ref{sec:holder} to obtain further tail estimates for solutions, 
	which combined with Hölder estimates will allow us to complete the proof of Proposition~\ref{prop:control_sigma}.
	Given $v \in C(\R,[0,1])$, we define
\begin{equation} \label{def_Rv}
R(v) := \sup \{x \in \R : v(x) > 0 \}.
\end{equation}
We recall the constants $\cexp{interface}$ and $\cexp{inter_integ}$ defined in \eqref{bounds_cexpinterface} and \eqref{bounds_cexpinterinteg}.
We then define
\begin{equation} \label{def_rN}
r_{N}(v) := \inf \left\{ r \in \R : N \int_\R v(x) e^{(1-\cexp{inter_integ})(x-r)} dx  \leq C_X \right\}.
\end{equation}
In the above, 
\begin{equation} \label{eq:CX_prel}
C_X \geq 1 \vee ( (1-\cexp{inter_integ})^{-1} + (\cexp{inter_integ} - \cinit{small})^{-1})
\end{equation} is a constant whose value is independent of $N$. A precise value is specified later on, see \eqref{def_CX}.



In Section~\ref{sec:proof}, we alluded to a stopping time $\sigmaN{interface}$ related to the right interface of $u^N_t$. 
We now define $\sigmaN{interface}$ to be
\begin{equation}  \label{eq:sigmainterface}
	\sigmaN{interface} := \inf \left\{t \geq 0 : R(v^N_t) \wedge r_N(v^N_t) > (1+ \cexp{interface}) \log N \right\}.
\end{equation}
This section is devoted to the proof of the following, which proves Proposition~\ref{prop:control_sigma} for $i = \ref{sigmaN:interface}$.
\begin{proposition} \label{prop_R_main}
	There exists $a>0$ such that for any $T >0$, for sufficiently large $N$,
	\begin{equation*}
		\bP \left(\sigmaN{interface} \leq \sigma_{2,N} \wedge \sigma_{3,N}  \wedge T \right) < N^{-a}.
	\end{equation*}
\end{proposition}


The proof of Proposition~\ref{prop_R_main} uses a novel argument which compares the dynamics of $R(v^N_t) \wedge r_N(v^N_t)$ to a random walk with negative drift. 
To control $R(v^N_t)$, we use a variant of the method of Krylov \cite{krylov_result_1997} to study the compact support property for SPDEs. 
The process $t\mapsto r_N(v^N_t)$ is controlled via an auxiliary process $t \mapsto X^{N,S}_t$ which we introduce shortly.
The key idea is that, while our estimates are insufficient to obtain the desired control of $R(v^N_t)$ and $X^{N,S}_t$, they can be controlled jointly using the random walk formalism alluded to above.

We now collect some definitions and notation which will be used in the proof of Proposition~\ref{prop_R_main}. For the proof of this result, it will be convenient to work with the equation on its original time-scale. Within this section, we will primarily work with the equation 
\begin{equation}\label{eq_spde_movingframe} 
d w^N_t = (\partial_{xx}w^N_t + \alpha \partial_x w^N_t + f(w^N_t))\, dt  + \sqrt{\frac{w^N_t(1-w^N_t)}{N} } dW(t).
\end{equation}
In particular, on our original probability space, the relationship between $v^N_t$ and $w^N_t$ is summarized by $\bP(w^N_{Nt} = v^N_t \, \forall t \geq 0)  = 1$.

We will also work with solutions to this equation started from various initial conditions, which we also denote by $(w^N_t)_{t \geq 0}$. For $w\in C(\R,[0,1])$, we denote the law of the solution to \eqref{eq_spde_movingframe} started from $w_0^N = w$ by $\bP_w$, and write $\mathbb{E}_w$ for the associated expectation. Given a solution to \eqref{eq_spde_movingframe}, we introduce the stopping times
\begin{align*}
\hat{\sigma}_{2,N} &:= \inf \{t \geq 0 :  | \eta(w^N_t)| \geq K\}, \\
\hat{\sigma}_{3,N} &:= \inf \{t \geq 0 : \|s(w^N_t)\|_\infty > \varepsilon_2 \}.
\end{align*}
These are simply the images of $\sigma_{2,N}$ and $\sigma_{3,N}$ under the map that sends $v^N_t$ to $w^N_t$. In particular, $\hat{\sigma}_{i,N} = N \sigma_{i,N}$ for $i=2,3$.

Working under $\bP_w$ (or $\bP$), for $t \geq 0$ we define
\begin{equation} \label{def_RNrN}
R^N_t := R(w^N_t), \quad r^N_t := r_N(w^N_t).
\end{equation}
In order to facilitate the analysis of $r^N_t$, we introduce the following auxiliary process defined under $\bP$ and $\bP_w$. For $S \in \R$, for $t \geq 0$ we define
\begin{equation} \label{def_XNT}
X^{N,S}_t := N \int_{\R} e^{(1-\cexp{inter_integ})(x-S)} w^N_t(x) dx.
\end{equation}
It then follows from the definition of $r^N_t$ that 
\begin{equation} \label{eq_r_X_link}
r^N_t = \inf \big\{S \in \R : X^{N,S}_t \leq C_X \big\}.
\end{equation}
We introduce a set of initial conditions for the SPDE involving simultaneous control of $R$ and $r_N$. For $S \in \R$, we define
\begin{align} \label{def_WS}
\cW_S := \left\{ w \in C(\R,[0,1]) : R(w) \leq S \, \text{ and }\, r_N(w) \leq S \right\}.
\end{align}
The initial states of the processes $R^N$, $r^N$ and $X^{N,S}$ under $\bP_w$ for $w \in \cW_S$ are summarized below: 
\begin{equation} \label{eq_WS_timezero}
\text{For any $w \in \cW_S$, $R^N_0 \vee r^N_0 \leq S$ and $X^{N,S}_0 \leq C_X$, $\bP_w$-a.s.}
\end{equation}
Next, we introduce several constants related to the non-linearity $f$ which appear only within this section.
 Recalling that $\cexp{inter_integ} < 1-\alpha/2$, we fix parameters $q_1$ and $q_2$ satisfying
	\begin{equation} \label{def_a2a3}
		q_1 \in (0, \cexp{inter_integ}(2-\alpha - \cexp{inter_integ}) \wedge (\alpha^2/4+1-\alpha)), \quad q_2:=  \cexp{inter_integ}(2 - \alpha - \cexp{inter_integ}) - q_1  > 0.
	\end{equation}
	Since $f'(0) = \alpha - 1$ and $f$ is $C^2$, we have $f(u) = (\alpha - 1)u + O(u^2)$ as $u\to 0$. Hence there exists $u_* \in (0,1)$ such that
\begin{equation} \label{eq_f_small1}
\text{For all } u \in [0,u_*], \, f(u) \leq (\alpha - 1 + q_1)u.
\end{equation}
Remarking that
\[ (1-\cexp{inter_integ})^2 - \alpha(1-\cexp{inter_integ}) = (1-\alpha) - \cexp{inter_integ}(2 - \alpha - \cexp{inter_integ}), \]
we conclude that 
\begin{equation} \label{eq_f_small2}
\text{For all $u \in [0,u_*]$, } ((1-\cexp{inter_integ})^2 - \alpha(1-\cexp{inter_integ}))u + f(u) \leq - q_2 u.
\end{equation}

Next, we derive a related estimate for $w^N$ which will be used several times in the sequel. First, we recall the constant $\varepsilon_2>0$, which appeared not long ago in the definition of $\hat{\sigma}_{3,N}$. We will hereafter assume that $\varepsilon_2 \leq u_*/2$, which we can achieve by making $\varepsilon_2$ smaller if necessary. This does not interfere with its other dependencies. Next, since the stationary wave profile $m$ satisfies $\lim_{x \to \infty} m(x) = 0$, we may choose $x_0$ to be sufficiently large so that 
\begin{equation}
	m(x) \leq u_*/2 \, \text{for all $x \geq x_0$.}
\end{equation}
From the definitions of $\hat{\sigma}_{2,N}$ and $\hat{\sigma}_{3,N}$, it then follows that 
\begin{align} \label{eq_wNsmall_tail}
	\text{For all $t \leq \hat{\sigma}_{2,N} \wedge \hat{\sigma}_{3,N}$ and $x \geq K + x_0$,} \quad w^N_t(x) \leq \sup_{|y|\leq K}m(x - y) + \varepsilon_2 \leq u_*. 
\end{align}

	Before carrying on the the main argument, we discuss some elementary properties of solutions $w^N$ to \eqref{eq_spde_movingframe}. 
	For $t>0$ and $x\in \R$, let us set
	\begin{equation} \label{eq:Gdefn}
		G_t(x) = \frac{1}{\sqrt{4\pi t}} \exp \left( - \frac{x^2}{4 t} \right),
	\end{equation} \label{eq:Qdefn}
	and let $ (Q(t), t \geq 0) $ be the operator semigroup defined by
	\begin{equation}
		Q(t)w(x) = \int_\R G_t(x + \alpha t - y) w(y) dy,
	\end{equation}
	which sends $ L^{2,\alpha} $ into $ L^{2,\alpha} $, and admits $ \partial_{xx} + \alpha \partial_{x} $ as its infinitesimal generator.
	A solution $w^N$ to \eqref{eq_spde_movingframe} with initial condition $w \in C(\R,[0,1])$ can be written in mild form as 
	\begin{equation*}
		w^N_t = Q(t) w^N_0 + N \int_{0}^{t} Q(t-s) f(w^N_s) ds + \int_{0}^{t} Q(t-s) \left( \sqrt{\frac{w^N_s(1-w^N_s)}{N}} dW(s) \right).
	\end{equation*}
	Furthermore, for any $ a \in \R $, 
	\begin{multline} \label{mild_spde_wN}
		w^N_t = e^{a t} Q(t) w + N \int_{0}^{t} e^{a (t-s)} Q(N(t-s)) \left( f(w^N_s) - a w^N_s \right) ds \\ + \int_{0}^{t} e^{a (t-s)} Q(t-s) \left(\sqrt{\frac{w^N_s(1-w^N_s)}{N}}  dW(s) \right).
	\end{multline}
	We shall sometimes abuse notation by writing $W(dxdt)$ for the martingale measure as in~\eqref{def_martingale_measure}.
	The representation \eqref{mild_spde_wN} will be used variously, with different values of $a$, throughout Sections~\ref{sec:interface}-\ref{sec:holder}.
	Its first application is the following easy bound on the tail of $\E[w]{w^N_t(x)}$ as $x \to \infty$.
	
	\begin{lemma} \label{lem:wNtailbd}
	Let $w \in C(\R,[0,1])$ satisfy $R(w) \leq N$.
	Then for $x\geq 2N$ and $t> 0$,
	$$
	\E[w]{w^N_t(x)}
	\leq e^{\|f'\|_\infty t} e^{-x^2/(16 t)}.
	$$
	\end{lemma}
	\begin{proof}
	We observe that $f(w) - \|f'\|_\infty w \leq 0$ for all $ w\in[0,1]$. Hence, by \eqref{mild_spde_wN} with $ a = \|f'\|_\infty $, for $t> 0$ and $x\in \R$,
	\begin{align*}
	w^N_t(x)&=
	\int_{\R} G_{t}(x+\alpha  t-y)e^{\|f'\|_\infty t}w(y)dy\\
	&\qquad +\int_0^t \int_{\R} G_{t-s}(x+\alpha  (t-s)-y)e^{\|f'\|_\infty(t-s)}
	(f(w^N_s(y))-\|f'\|_\infty w^N_s(y))dyds\\
	&\qquad + N^{-1/2} \int_0^t \int_{\R} G_{t-s}(x+\alpha N(t-s)-y)e^{(1+\alpha)N(t-s)}
	\sqrt{w^N_s(y)(1-w^N_s(y))}W(dy ds)\\
	&\leq \int_{\R} G_{t}(x+\alpha  t-y)e^{\|f'\|_\infty t}w(y)dy\\
	&\qquad + \int_0^t \int_{\R} G_{t-s}(x+\alpha (t-s)-y)e^{(1+\alpha)(t-s)}
	\sqrt{w^N_s(y)(1-w^N_s(y))}W(dy ds).
	\end{align*}
	Recall that $R(w) \leq N$ by assumption. 
	Therefore for $x\geq 2N$ and $t>0$, letting $(B_s)_{s\ge 0}$ denote a Brownian motion in the second line, and using a Gaussian tail bound in the third line,
	\begin{align*} 
	\E{w^N_t(x)}
	&\leq e^{\|f'\|_\infty t}\int_{\R} G_{t}(x+\alpha t-y)w(y)dy \\
	&\leq e^{\|f'\|_\infty t} \P[x+\alpha  t]{B_{2 t}\leq N}  \\
	&\leq e^{\|f'\|_\infty t} e^{-(x+\alpha  t-N)^2/(4 t)}  \\
	&\leq e^{\|f'\|_\infty t} e^{-x^2/(16 t)},  
	\end{align*}
	where the last inequality follows
	since $(x-N)^2 \ge \frac 14 x^2$ by our choice of $x$.
	\end{proof}


We conclude our preliminary discussion by collecting some elementary properties of $t\mapsto R^N_t$. A stronger version of the compact interface property for \eqref{eq_spde_movingframe} is finite speed of propagation of the support, which entails, among other things, that $R^N_t$ cannot increase by jumps, and so can only increase continuously. While this is undoubtedly true, it would be cumbersome to give a full proof, so instead we will make do with two weaker properties which we now state. For $y > R^N_0$, let $\tau^R_y = \inf\{ t \geq 0 : R^N_t \geq y\}$. 
\begin{lemma} \label{lemma_Rsimple} Let $w \in C(\R,[0,1])$ satisfy $R(w) <\infty$. Then \\
 (a) For any $y > R(w)$,  $\bP_w(\tau^R_y > 0) = 1$. 
\\ (b) For any $y > R(w)$, $\bP_w(R^N_{\tau^R_y} = y) = 1$. 
\\ (c) For fixed $N>0$ and any $\Lambda >0$, $\lim_{x \to \infty} \bP_w(R^N_t > x \text{ for some } t \in [0,\Lambda]) = 0$.
\end{lemma}
Part (a) effectively states that $R^N_t$ propagates with finite speed over a fixed interval. Part (b) states that its passage over a fixed level is continuous. Both claims are standard for SPDEs like the one we consider, and can be shown for example using the method of Krylov \cite{krylov_result_1997}. In Section~\ref{s_krylov}, we adapt Krylov's general approach to obtain precise, tailor-made estimates concerning $R^N_t$ in the present setting. Since the claims above can be proved using a much simpler version of the same arguments along with more elementary moment estimates, we omit the proof. Part (c) can likewise be proved using a standard application of Krylov's method.

Finally, we state a lemma concerning the initial value of $R^N_0 \vee r^N_0$ under $\bP$, i.e. with initial conditions $u^N_0$.
\begin{lemma} \label{lem:R0r0}
Under Assumption~\ref{assumpt:v0}, under $\bP$ we have $R^N_0 \vee r^N_0 <  (1- \cexp{inter_integ})^{-1}(1+\cinit{small}) \log N$, and in particular $u^N_0 \in \cW_{ (1- \cexp{inter_integ})^{-1}(1+\cinit{small}) \log N}$.
Furthermore, $\sigmaN{interface} > 0$ almost surely.
\end{lemma}
\begin{proof}
From Assumption~\ref{assumpt:v0}.\ref{v0:compact_support} and \eqref{bounds_cexpinterface}, $R^N_0 < (1+\cinit{small}) \log N < (1+\cexp{interface})\log N$.  
Let $x_0 = (1- \cexp{inter_integ})^{-1}(1+\cinit{small}) \log N$. Then $e^{-(1-\cexp{inter_integ}) x_0} = N^{-1-\cinit{small}}$. Hence, by Assumption~\ref{assumpt:v0}.\ref{v0:tail}, and using $u^N_0(x) \leq 1$ for $x \leq 0$,
\begin{align*}
\int_{\R} u^N_0(x) e^{(1-\cexp{inter_integ})(x - x_0)}dx &\leq N^{-1-\cinit{small}} \left( (1-\cexp{inter_integ})^{-1} + N^{\cinit{small}} \int_0^{\infty} e^{-(\cexp{inter_integ} - \cinit{small})x} dx\right)
\\ &< N^{-1} \left( (1-\cexp{inter_integ})^{-1} + (\cexp{inter_integ} - \cinit{small})^{-1} \right).
\end{align*}
and hence, by \eqref{def_rN} and \eqref{eq:CX_prel}, we have $r^N_0 < x_0$. Hence, we have shown that $R^N_0 \vee r^N_0 <  (1- \cexp{inter_integ})^{-1}(1+\cinit{small}) \log N$ as desired.

To see that $\sigmaN{interface} > 0$, first note that $(1- \cexp{inter_integ})^{-1}(1+\cinit{small}) \log N < (1+\cexp{interface}) \log N$ by \eqref{bounds_cinitsmall}.
It follows from Lemma~\ref{lemma_Rsimple}(a) that $\inf\{t \geq 0 : R^N_0 > (1+\cexp{interface})\log N\} >0$. 
Similarly, the functional $t\mapsto \int_\R w^N_t(x) e^{(1-\cexp{inter_integ})x} dx$ is a.s. continuous. (We will prove this fact independently later on, see \eqref{eq_XN_weak0}.) 
This implies that $t \mapsto r^N_t$ is continuous, and hence $\inf\{t \geq 0 : r^N_0 > (1+\cexp{interface})\log N\} >0$ as well. The proof is complete.
\end{proof}

The rest of this section is organized as follows. In the next section, we prove Proposition~\ref{prop_R_main} assuming a key result, Proposition~\ref{prop_SN_inc}. This result is then proved in Section~\ref{s_markovpropproof} using intermediate results derived in Subsections~\ref{s_krylov}-\ref{s_Xanalysis}.

\subsection{The walk process and proof of Proposition~\ref{prop_R_main}} \label{s_walkproof}

%
%
%

The proof of Proposition~\ref{prop_R_main} is based on the analysis of a discrete-time process whose dynamics simultaneously incorporate the behaviour of $R^N_t$ and $r^N_t$. In this section, we introduce this walk process, which we denote $S^N = (S^N_n)_{n\in \N_0}$, and use it to prove Proposition~\ref{prop_R_main}. The proof uses Proposition~\ref{prop_SN_inc}, whose proof is the subject of Subsections~\ref{s_krylov}-\ref{s_markovpropproof}.

We define a walk process by observing our solution at discrete times spaced by time increments of size $\Lambda > 0$, where $\Lambda$ will be chosen to be large enough (independently of $N$) so that certain estimates hold. Heuristically, the increment of the walk process over an interval of length $\Lambda$ corresponds to the larger increment of the functionals $R^N$ and $r^N$ over the interval, and thus the process serves as an upper bound for the values of $R^N$ and $r^N$ of our solution. Thus, at an intuitive level, we are interested in the process
\begin{equation*}
S^{N,\Lambda,*}_n   = r^N_{n\Lambda} \vee R^N_{n \Lambda} =  r_{N}(w^N_{n\Lambda}) \vee R(w^N_{n\Lambda}),
\end{equation*}
where we recall that $w^N_t = v^N_{Nt}$ is a solution to \eqref{eq_spde_movingframe}.

The relatively simple definition above does not work for a few technical reasons. First, the estimates allowing us to control $R^N$ and $r^N$ require, roughly speaking, that the starting values of these quantities are not too small. To handle this, we incorporate a minimum with a deterministic lower bound into the definition of the process. Second, we require uniform control of $R^N_t$ and $r^N_t$, but $S^{N,*}_n$ as written above only describes their values along a discrete set of times. Therefore, keeping the above in mind to inform our intuition, we make our definition. First, for $n \in \N$ define the events
\begin{equation} \label{def_FNn}
	F^{N,\Lambda}_n = \{n\Lambda < \hat{\sigma}_{2,N}  \wedge \hat{\sigma}_{3,N} \}
\end{equation}
and
\begin{align*}
G^{N,\Lambda}_n = F^{\Lambda,N}_n \cap \bigg\{ \sup_{t \in [(n-1)\Lambda , n\Lambda]} R^N_t < S^N_{n-1} + y_0 \bigg\} \cap \bigg\{ \sup_{t \in [(n-1)\Lambda,n\Lambda]} X^{N,S^N_{n-1}}_t \leq 4C_X \bigg\}.
\end{align*}
In the above, $y_0$ is a positive constant whose value is independent of $N$, for which a value is specified later in \eqref{def_y0Lambda}. Then, we define the (random) interval $I^{N,\Lambda}_n$ by
\begin{equation} \label{def_INn}
	I^{N,\Lambda}_n := [((n-1)\Lambda) \wedge \hat{\sigma}_{2,N} \wedge \hat{\sigma}_{3,N}, (n\Lambda) \wedge \hat{\sigma}_{2,N} \wedge \hat{\sigma}_{3,N}],
\end{equation}
noting that $I^{N,\Lambda}_n = \{\hat{\sigma}_{2,N} \wedge \hat{\sigma}_{3,N}\}$ if $\hat{\sigma}_{2,N} \wedge \hat{\sigma}_{3,N} \leq (n-1)\Lambda$, and define
\begin{equation}  \label{def_calRn}
\mathcal{R}^{N,\Lambda}_n :=  \indc_{G_n}  R^N_{n\Lambda } + \indc_{G_n^c}  \sup_{t \in I^{N,\Lambda}_n} R^N_t. 
\end{equation}
$\mathcal{R}^N_n$ is a modification of the term $R^N_{n\Lambda}$ in the definition of $S^{N,\Lambda,*}_n$ which affords uniform control of the right end point while retaining the dynamics we need in order to prove our result. For a similarly uniformized version of $r^N_{n\Lambda}$, we define
\begin{equation}  \label{def_frakRn}
\mathfrak{R}^{N,\Lambda}_n := \indc_{G_n} r^N_{n \Lambda} + \indc_{G_n^c}  \sup_{t \in I^{N,\Lambda}_n} r^N_t.
\end{equation}
We then define the walk process $S^{N,\Lambda} = (S^{N,\Lambda}_n)_{n \in \N_0}$ by 
\begin{equation} \label{def_Sn}
S^{N,\Lambda}_n = \left((1-\cexp{inter_integ})^{-1} \log N\right) \vee \mathcal{R}^{N,\Lambda}_n \vee \mathfrak{R}^{N,\Lambda}_n.
\end{equation}
Next, for $n \in \N$ define
\begin{equation} \label{def_xin}
\xi^{N,\Lambda}_n = \left(\mathcal{R}^{N,\Lambda}_{n} \vee \mathfrak{R}^{N,\Lambda}_n \right) - S^{N,\Lambda}_{n-1}.
\end{equation}
Then it can be readily verified that
\begin{equation} \label{def_Sinc_xi}
S^{N,\Lambda}_n = \left(S^{N,\Lambda}_{n-1} + \xi^{N,\Lambda}_n \right) \vee \left((1-\cexp{inter_integ})^{-1} \log N \right)
\end{equation}
for all $n \in \N$. In particular, $S^{N,\Lambda}$ has increments given by $\{\xi^{N,\Lambda}_n : n \in \N\}$ but has ``reflection'' at $(1-\cexp{inter_integ})^{-1} \log N$.

We define a discrete filtration $(\mathcal{G}^{N,\Lambda})_{n \in \N_0}$ by $\mathcal{G}^{N,\Lambda}_n = \mathcal{F}_{n\Lambda }$, where $(\mathcal{F}_t)_{t \geq 0}$ is the filtration on the probability space on which our solution $w^N$ is constructed. 
The following is our key estimate concerning the increments of $S^{N,\Lambda}$.

\begin{proposition} \label{prop_SN_inc}
There exists a constant $\rho > 0$ such that the following holds: for any $\delta > 0$ there exists $\bar{\Lambda}(\delta) \geq 1$ such that if $\Lambda = \bar{\Lambda}(\delta)$, for sufficiently large $N$, for every $n \in \N$,
\begin{align*} 
\bP\left(\{\xi^{N,\Lambda}_n  > -\rho \Lambda\} \cap F^{N,\Lambda}_n \, | \, \mathcal{G}^{N,\Lambda}_{n-1} \right) \leq \frac 1 2 \delta \, \text{ a.s.}
\end{align*}
and
\begin{align*}
  \indc_{F^{N,\Lambda}_{n-1}} \,\,\bP\left(\xi^{N,\Lambda}_n  > k \rho \Lambda\, | \, \mathcal{G}^{N,\Lambda}_{n-1} \right)  \leq \frac 1 2 e^{-c \Lambda k^2}  \quad \text{ a.s.  for all } k \geq M,
 \end{align*}
 where $c>0$ and $M \geq 1$ are constants which do not depend on $\delta$.
\end{proposition}
This result is proved in Section~\ref{s_markovpropproof} using estimates derived in Sections~\ref{s_krylov}-\ref{s_Xanalysis}. For the rest of this subsection, we proceed using $\Lambda = \bar{\Lambda}(\delta)$ and omit dependence on $\Lambda$ in our notation from now on. For the time being we continue with arbitrary $\delta>0$, but it will be chosen to take a sufficiently small value later on.


%
%
Given $\delta >0$, we fix $k_\delta \geq M$ such that
\begin{equation} \label{eq_kdeltabd}
	e^{-c\Lambda k_\delta^2} \leq \frac 1 2 \delta.
\end{equation}
Next, we define
\begin{equation*}
\tilde{\xi}^N_n := \xi^N_n \indc_{F^N_{n-1}}  - \rho \Lambda \indc_{(F^N_n)^c} - k_\delta \rho \Lambda \indc_{ F^N_{n-1} \cap (F^N_{n})^c }, 
\end{equation*}
as well as the process $\tilde{S}^N = (\tilde{S}^N_n)_{n \in \N_0}$, with $\tilde{S}^N_0 = S^N_0$, defined recursively by
\begin{equation*}
\tilde{S}^N_{n} := (\tilde{S}^N_{n-1} + \tilde{\xi}^N_n) \vee \left((1-\cexp{inter_integ})^{-1} \log N\right).
\end{equation*}
It is then immediate from the definition of $\tilde{\xi}^N_n$ that for $n \in \N$,
\begin{equation} \label{eq_StildeS}
	S^N_k = \tilde{S}^N_k \text{ for all $k \leq n$ on $F^N_n$} \,\, \text{ and } \,\,  S^N_n = \tilde{S}^N_n + (1+k_\delta)\rho\Lambda  \text{  on $F^N_{n-1} \cap (F^N_n)^c$}.
\end{equation}
%
Using Proposition~\ref{prop_SN_inc}, we can easily obtain conditional estimates on the distribution of the increments $\{\tilde{\xi}^N_n, n \in \N\}$. Since $F^N_n \subseteq F^N_{n-1}$, on $(F^N_{n-1})^c$ we have $\tilde{\xi}^N_n = -\rho \Lambda$, and hence
\begin{equation*}
 \indc_{(F^N_{n-1})^c}  \bP(\tilde{\xi}^N_n > -\rho \Lambda \, | \, \mathcal{G}^N_{n-1}) = 0 \,\, \text{a.s.}
\end{equation*} 
Similarly, it a.s. holds that
\begin{align*} 
 &\indc_{F_{n-1}^N} \bP(\tilde{\xi}^N_n > -\rho \Lambda \, | \, \mathcal{G}^N_{n-1}) 
 \\ &\hspace{1 cm} =  \indc_{F^N_{n-1}} \bP( \{\tilde{\xi}^N_n > -\rho \Lambda \} \cap F^N_n  | \, \mathcal{G}^N_{n-1}) +   \indc_{F^N_{n-1}} \bP( \{\tilde{\xi}^N_n > -\rho \Lambda \} \cap (F_n^N)^c  | \, \mathcal{G}^N_{n-1})
\\ &\hspace{1 cm} =  \indc_{F^N_{n-1}} \bP( \{\xi^N_n > -\rho \Lambda\} \cap F^N_n  | \, \mathcal{G}^N_{n-1}) +   \indc_{F^N_{n-1}} \bP( \{\xi^N_n -(1+k_\delta)\rho \Lambda > -\rho \Lambda \} \cap (F_n^N)^c  | \, \mathcal{G}^N_{n-1})
\\ &\hspace{1 cm}\leq \indc_{F^N_{n-1}} \delta,
\end{align*}
where we apply Proposition~\ref{prop_SN_inc} and \eqref{eq_kdeltabd} in the last line. Combining the two previous statements yields 
\begin{align} \label{eq_tildeS_p1}
\bP(\tilde{\xi}^N_n > - \rho \Lambda \,  | \, \mathcal{G}^N_{n-1}) \leq  \delta \, \text{ a.s.},
\end{align}
and repetition of the same argument with $-\rho\Lambda$ replaced by $k \rho \Lambda$, and using the other bound from Proposition~\ref{prop_SN_inc}, leads to 
\begin{align} \label{eq_tildeS_pk}
\bP(\tilde{\xi}^N_n \geq  k \rho \Lambda \,| \, \mathcal{G}_{n-1}) \leq  \frac 1 2 e^{-c\Lambda k^2} \, \text{ a.s.} \, \text{ for all } k \geq M.
\end{align}

The next step is to couple $\tilde{S}^N$ with a random walk in which the steps are identically distributed. To do so, for $\delta \in(0,1)$ we define a probability measure $\mu_\delta$ on $\R$. First, let $\mu_\delta(\{-\rho \Lambda\}) = 1-\delta$, and then set $\lambda = \rho^2 \Lambda c^{-1}$ and define
\[ g(x) = \lambda^{-1/2} e^{-x^2/(2\lambda)} \,\text{ for $x \in \R$.}\]
We note that $\lambda$ and hence $g$ depend on $\delta$ through $\Lambda = \bar{\Lambda}(\delta)$, but we omit this dependence as it is not important. Next, if $b(\delta) = \inf \{ x \in \R: \int_x^\infty g(y) dy < \delta  \}$, we set $\mu_\delta(dx) =g(x)dx$ for $x \geq b(\delta)$.
Thus, $\mu_\delta$ has an atom of mass $1-\delta$ at $-\rho \Lambda$ and density $g(x)$ for $x \geq b(\delta)$. In other words, $\mu_\delta$ is the probability measure
\begin{equation} \label{def_mudelta} 
	\mu_\delta(dx) := (1-\delta)\cdot \delta_{-\rho \Lambda}(dx) + \indc_{x \geq b(\delta)}g(x) dx,\end{equation}
where $ \delta_{-\rho \Lambda} $ denotes a Dirac mass at $-\rho \Lambda$.

\begin{lemma} \label{lemma_stochdom}
For sufficiently small $\delta$, the following hold: if $\Lambda = \bar{\Lambda}(\delta)$ and $N$ is sufficiently large, for any $x \in \R$ and $n \in \N$, with probability one
\begin{equation}\label{eq_stochdomlemma} \bP(\tilde{\xi}^N_n >  x \,| \, \mathcal{G}_{n-1}) \leq \mu_\delta((x,\infty)).\end{equation}
In particular, the regular conditional distribution of $\tilde{\xi}^N_n$ given $\mathcal{G}_{n-1}$ a.s. satisfies
\[ \bP(\tilde{\xi}^N_n \in \cdot \,| \, \mathcal{G}_{n-1}) \preceq \mu_\delta,\]
where $\preceq$ denotes stochastic domination.
\end{lemma}


\begin{proof} First, we remark that for $x < -\rho \Lambda$, $\mu_\delta((x,\infty)) = 1$, so \eqref{eq_stochdomlemma} holds trivially for such $x$. For $x \in [-\rho \Lambda, b(\delta)]$, by \eqref{eq_tildeS_p1} we a.s. have 
\begin{align*}
\bP(\tilde{\xi}^N_n > x \,  | \, \mathcal{G}_{n-1}) \leq \bP(\tilde{\xi}^N_n > - \rho \Lambda \,  | \, \mathcal{G}_{n-1}) \leq \delta = \mu_\delta((x,\infty)).
\end{align*}
It remains to handle $x > b(\delta)$. We change variables to $k$ by setting $x = k\rho \Lambda$, so we need to prove the desired domination for $k \geq b(\delta)/(\rho \Lambda)$. 
We recall the Gaussian tail estimate $\int_x^\infty e^{-y^2/2} dy \geq (x^{-1} - x^{-3})e^{-x^2/2}$ for $x >0$, which implies that $\int_x^\infty e^{-y^2/2} dy \geq \tfrac{1}{2x}e^{-x^2/2}$ for $x \geq \sqrt 2$.
Hence, using scaling and the definition of $\lambda$, we have
\begin{align*}
\mu_\delta((k\rho \Lambda,\infty)) = \int_{\lambda^{-1/2} k \rho\Lambda}^\infty e^{-x^2/2}dx &\geq \frac{1}{2\lambda^{-1/2} k \rho \Lambda} e^{-\lambda^{-1} k^2 \rho^2 \Lambda^2/2}=  \frac{1}{2 k c^{1/2} \Lambda^{1/2} }e^{-c \Lambda k^2/2}
\end{align*}
provided $k \geq b(\delta)/(\rho \Lambda)$, which ensures that the first equality in the above is true, and $k \geq\lambda^{1/2} /( \rho\Lambda)$, which allows us to use the tail estimate. 
We remark that $\lim_{\delta \downarrow 0} b(\delta) = \infty$, which along with the definition of $\lambda$ implies that the first condition implies that second for sufficiently small $\delta>0$.
Finally, we note that the inequality above implies that for some $M'$, if $k\Lambda^{1/2} \geq M'$ and $k \geq b(\delta)/(\rho \Lambda)$,
\[ \mu_\delta((k\rho \Lambda,\infty))  \geq \frac 12  e^{-c \Lambda k^2}.\]
Again using $\lim_{\delta \downarrow 0} b(\delta) = \infty$ and $\Lambda \geq 1$, the above holds for $k \geq b(\delta)/(\rho \Lambda)$ for small enough $\delta$. 
Finally, we may choose $\delta$ to be small enough so that $M \leq b(\delta)/(\rho \Lambda)$. In this case, by \eqref{eq_tildeS_pk} we a.s. have
\begin{align*}
\bP(\tilde{\xi}^N_n \geq  k \rho \Lambda \,| \, \mathcal{G}_{n-1}) \leq  \frac 1 2 e^{-c\Lambda k^2} \leq \mu_\delta((k\rho \Lambda,\infty))
\end{align*} 
for all $k > b(\delta)/(\rho\Lambda)$. Thus the desired bound a.s. holds for $x > b(\delta)$, and hence for all $x \in \R$, which completes the proof.
\end{proof}

By the tower property, \eqref{eq_stochdomlemma} holds a.s. when $\mathcal{G}_{n-1}$ is replaced with $\sigma(\tilde{S}^N_0,\dots,\tilde{S}^N_{n-1})$, and the same is also true for the stochastic domination of the regular conditional distribution. One may then use standard coupling techniques (for example using \cite[Theorem IV.5.8]{lindvall_lectures_2002}) to construct a probability space on which there exists a monotone coupling between a copy of $\tilde{S}^N$ and a random walk whose increments are distributed according to $\mu_\delta$. We abuse notation and denote the law on this probability space by $\bP$, even though it does not necessarily still contain our original process $w^N$. On this probability space, $(\tilde{\xi}^N_n)_ {n \in \N}$ and $(\tilde{S}^N_n)_{n \geq 0}$ have the same distribution as on the original space and are adapted to a filtration $(\mathcal{H}_n)_{n \geq 0}$, and furthermore there exist a family of random variables $(J^\delta_n)_{n \in \N}$ such that the following hold for all $n \in \N$:
\begin{itemize}
	\item $J^\delta_n$ is $\mathcal{H}_n$-measurable.
	\item Conditional on $\mathcal{H}_{n-1}$, $J_n$ is distributed according to $\mu_\delta$. (More precisely, the regular conditional distribution of $J^\delta_n$ given $\mathcal{H}_{n-1}$ is a.s. equal to $\mu_\delta$.)
	\item With probability one, $\tilde{\xi}^N_n \leq J^\delta_n$.
\end{itemize}
Since the conditional distribution of $J^\delta_n$ given $\mathcal{H}_{n-1}$ is a.s. constant (i.e. $\mu_\delta$), we may without loss of generality assume that the $(J^\delta_n)_{n \in \N}$ are independent. Let $(\mathcal{S}^\delta_n)_{n \in \N_0}$ be the reflected random walk defined by $\mathcal{S}^\delta_0 = S^N_0 -  (1-\cexp{inter_integ})^{-1} \log N$ and
\[ \mathcal{S}^\delta_n = \left( \mathcal{S}^\delta_{n-1} + J^\delta_n \right) \vee 0, \quad n \in \N.\]
We have shifted the process by $(1-\cexp{inter_integ})^{-1} \log N$, so the reflection now occurs at $0$, which will simplify the notation of some upcoming arguments. Since $J^\delta_n \geq \tilde{\xi}^N_n$ a.s. for every $n$ and $\mathcal{S}^N_0 = (1-\cexp{inter_integ})^{-1} \log N + \tilde{S}^N_0$, it follows easily that with probability one,
\[ \tilde{S}^N_n \leq  (1-\cexp{inter_integ})^{-1} \log N + \mathcal{S}^\delta_n \quad \text{ for all } n \in \N_0.\]
By \eqref{eq_StildeS}, we have $S^N_n \leq \tilde{S}^N_n + (1+k_\delta)\rho\Lambda$ for all $n \leq \lceil \Lambda^{-1} ( \hat{\sigma}_{2,N} \wedge \hat{\sigma}_{3,N}) \rceil$, 
so the following result is now immediate.

\begin{lemma} \label{lemma_coupling} For sufficiently small $\delta > 0$ and $\Lambda = \bar{\Lambda}(\delta)$, for sufficiently large $N$,
\[\bP\bigg(\max_{n \leq M \wedge \lceil \Lambda^{-1}(\hat{\sigma}_{2,N} \wedge \hat{\sigma}_{3,N} )\rceil} \,S^{N,\Lambda}_n >  r + (1+k_\delta)\rho\Lambda\bigg)  \leq \bP\bigg(\max_{n \leq M } \, \mathcal{S}^\delta_n > r - (1-\cexp{inter_integ})^{-1} \log N \bigg)\]
for all $M \in \N$ and $r \in \R$. 
\end{lemma}

Thus, proving Proposition~\ref{prop_R_main} can be reduced to proving an appropriate statement about $\mathcal{S}^\delta$, which we address in the next lemma. 

\begin{lemma} \label{lemma_RW_bound} There exists $a>0$ such that the following holds: for sufficiently small $\delta$, $\Lambda = \bar{\Lambda}(\delta)$, and any $\ell >0$ and $r \in \R$,
for sufficiently large $N$,
\[ \bP \left(\max_{n \leq \lceil \ell N \rceil }\mathcal{S}^\delta_n > (1+ \cexp{interface}- (1-\cexp{inter_integ})^{-1} ) \log N +r \right) < N^{-a}. \]
\end{lemma}

In order to prove the above, we will use some relatively standard arguments on the random walk $\mathcal{S}^\delta$. Of course, under $\bP$, $\mathcal{S}^\delta$ is a random walk with reflection at $0$; we will make use of an excursion decomposition of $\mathcal{S}^\delta$, and for this it will be convenient to consider a version of the random walk without reflection started from any initial point. In the sequel, for $x \geq 0$ we write $\bP^{\mathcal{S}^\delta}_x$ to denote the law of a random walk $\mathcal{S}$ without reflection started from $x$. That is, under $\bP^{\mathcal{S}^\delta}_x$, $\mathcal{S}^\delta$ is a standard random walk with increments distributed according to $\mu_\delta$.

For a discrete path $s : \N_0 \to \R$ and $a \in \R$, define
\[\tau^-_a(s) := \inf \{n \in \N_0 : s(n) \leq a\}, \quad \tau^+_a(s) := \inf \{n \in \N_0 : s(n) \geq a\} .\]
By \eqref{bounds_cexpinterinteg}, we fix a parameter $q \in (0,1/2)$ defined implicitly by the equation
\begin{equation}\label{defeq_q} \left( \frac{1-q}{q} \right)^{(1 + \cexp{interface} - (1-\cexp{inter_integ})^{-1}) / (\rho\Lambda)}= e^3. \end{equation}

\begin{lemma} \label{lemma_gamblerbd} For sufficiently small $\delta>0$, 
\begin{equation*}
\bP^{\mathcal{S}^\delta}_x( \tau^+_{b}(\mathcal{S}^\delta) <  \tau^-_{a}(\mathcal{S}^\delta) ) \leq  \frac{\left( \frac{1-q}{q}\right)^{(x-a)/(\rho\Lambda)+1} - 1 }{\left(  \frac{1-q}{q}\right)^{(b - a)/(\rho\Lambda)} - 1}
\end{equation*}
for all real numbers $a \leq x \leq b$.
\end{lemma}

We postpone the proof of Lemma~\ref{lemma_gamblerbd} until the end of the section. First, we use it to prove Lemma~\ref{lemma_RW_bound}, then complete the proof of Proposition~\ref{prop_R_main}. 

\begin{proof}[Proof of Lemma~\ref{lemma_RW_bound}]
Within this proof, we suppress dependence on $\delta$ in our notation. We begin by fixing $\ell >0 $ and $r \in \R$ and defining some stopping times. We set $\gamma^+_0 = 0$, and for $k \in \N$ recursively define
\[\gamma^0_k := \inf \{ n > \gamma^+_{k-1} : \mathcal{S}_n = 0\}, \quad \gamma_k^+ = \inf \{ n > \gamma_k^0 : \mathcal{S}_n > 0\}.  \] 
That is, $\gamma^0_k$ and $\gamma^+_k$ denote the times of subsequent visits to, and departures from $0$, and (for sufficiently small $\delta>0$) they are clearly finite for all $k$ due to the negative drift of $\mathcal{S}$. For convenience, we will assume without loss of generality that $\mathcal{S}_0 >0$, or equivalently that $S^N_0 > (1-\cexp{inter_integ})^{-1} \log N$. The proof is identical in the complimentary case; one simply sets $\gamma^0_0 = 0$ instead of $\gamma^+_0 = 0$ and then defines the stopping times in the same way.

We remark that $\gamma^+_{\lceil \ell N \rceil} \geq \lceil \ell N \rceil$, which can be seen immediately from its definition. Thus, it suffices to consider the probability that $\mathcal{S}_n$ reaches the level $\epsilon \log N + r$ before time $\gamma^+_{ \lceil \ell N \rceil}$, where we define
\begin{equation} \label{eq_epsa0a1}
	\epsilon = 1+\cexp{interface} - (1-\cexp{inter_integ})^{-1} > 0.
\end{equation}
The probability that it does so is bounded above by the sum of the probability that $\mathcal{S}_n$ crosses $\epsilon \log N + r$ before its first visit to $0$ and the probability that one of the first  $\lceil  \ell N \rceil$ excursions above zero reaches $\epsilon \log N + r$.

Let us first consider the excursions started from $0$. By the strong Markov property, the excursions are independent. The first step of each excursion is a copy of $J$ (which has distribution $\mu_\delta$) conditioned to be positive. Let $\mathbf{J}$ denote a random variable with this distribution. Then the law of excursions of $\mathcal{S}$ away from zero is given by $\bP^{\mathcal{S}}_{\mathbf{J}}(\cdot)$. In particular, by the strong Markov property, the probability of the walk exceeding $\epsilon \log N + r$ on the $k$th excursion away from zero can be written
\begin{align} \label{eq_walk_excursion_markov}
\bP( \mathcal{S}_n >  \epsilon \log N \text{ for some } n \in \N \text{ satisfying } \gamma^+_k \leq n < \gamma_{k+1}^0 )  = \mathbb{E}[\bP^\mathcal{S}_{\mathbf{J}}(\tau^+_{\epsilon \log N + r}(\mathcal{S}) < \tau^-_0(\mathcal{S}))].
\end{align}
Since $J$ has distribution $\mu_\delta$, whose definition we recall from \eqref{def_mudelta}, we have
\begin{equation*}
\bP(\mathbf{J} > x ) = \bP(J > x \, | \, J >0) \leq \delta^{-1} \int_x^\infty \lambda^{-1/2} e^{-y^2/(2\lambda)}dy \leq \delta^{-1} \sqrt{2\pi} e^{-x^2/(2\lambda)}
\end{equation*}
for all $x \geq 0$, where $\lambda > 0$ is the constant introduced below \eqref{eq_tildeS_pk}, and in particular we have
\begin{equation} \label{eq_boldJ_bd}
\bP(\mathbf{J} > (\log N)^{2/3}) \leq \delta^{-1} \sqrt{2\pi} N^{-(\log N)^{1/3}/(2\lambda )} \leq N^{-2}, 
\end{equation}
where the second inequality holds for sufficiently large $N$. Next, we remark that by Lemma~\ref{lemma_gamblerbd}, for $x \leq \epsilon \log N + r$ we have
\begin{align*}
\bP^{\mathcal{S}}_x( \tau^+_{\epsilon \log N +r}(\mathcal{S}) < \tau^-_0(\mathcal{S})) & \leq \frac{\left( \frac{1-q}{q}\right)^{x/(\rho\Lambda)+1} - 1 }{\left( \frac{1-q}{q}\right)^{(\epsilon \log N + r)/(\rho\Lambda)} - 1}.
\end{align*}
It is easy to show that for sufficiently large $N$, the numerator is at most $N^{1/2}$ for all $x \leq (\log N)^{2/3}$. Simplifying the denominator using \eqref{eq_epsa0a1} and the definition of $q$ (see \eqref{defeq_q}), we obtain that for large $N$,
\begin{align} \label{eq_uniformgamblerbd}
&\text{For all $x \leq (\log N)^{2/3}$, } \quad \bP^{\mathcal{S}}_x( \tau^+_{\epsilon \log N +r}(\mathcal{S}) < \tau^-_0(\mathcal{S}))  \leq \frac{N^{1/2} }{(\frac{1-q}{q})^{r/(\rho \Lambda)} N^{3}  - 1} \leq N^{-2}.
\end{align}
We now return to \eqref{eq_walk_excursion_markov}. Conditioning on $\mathbf{J}$ and partitioning over $\{\mathbf{J} \leq (\log N)^{2/3}\}$, \eqref{eq_boldJ_bd} and \eqref{eq_uniformgamblerbd} now imply that for sufficiently large $N$, for every $k \in \N$,
\begin{align*}
\bP( \mathcal{S}_n >  \epsilon \log N + r \text{ for some } n \in \N \text{ satisfying } \gamma^+_k \leq n < \gamma_{k+1}^0 ) \leq 2N^{-2}.
\end{align*}

Recall that if $\gamma^0_{k} \leq n < \gamma^+_k$ for any $k$, then $\mathcal{S}_n = 0$, so that the probability that $\mathcal{S}$ attains a value above $\epsilon \log N + r$ for such an $n$ is zero. Hence, recalling our observation that $\gamma^+_{\lceil \ell N \rceil} \geq \lceil \ell N \rceil$, we obtain from the above that 
\begin{align*}
&\bP( \mathcal{S}_n >  \epsilon \log N + r \text{ for some } n \in \N \text{ satisfying } \gamma^0_1 \leq n \leq \lceil \ell N \rceil ) \notag
\\ &\hspace{1 cm}\leq  \bP( \mathcal{S}_n >  \epsilon \log N + r \text{ for some } n \in \N \text{ satisfying } \gamma^0_1 \leq n \leq \gamma_{ \lceil \ell N \rceil}^0 ) \notag
\\ &\hspace{1 cm}\leq  \sum_{k=1}^{ \lceil \ell N \rceil - 1} \bP( \mathcal{S}_n >  \epsilon \log N  +r \text{ for some } n \in \N \text{ satisfying } \gamma^+_k \leq n < \gamma^0_{k+1} ) \notag
\\ &\hspace{1 cm}\leq  2 \lceil  \ell N \rceil  N^{-2} ,
\end{align*}
%
and hence the probability vanishes as $N\to \infty$. It remains to bound the probability that $\mathcal{S}$ exceeds $\epsilon \log N + r$ before $\gamma^0_1$. We remark that 
\begin{align*}
\bP( \mathcal{S}_n >  \epsilon \log N + r \text{ for some } n\leq \gamma^0_1) &=  \bP^\mathcal{S}_{\mathcal{S}_0}(\tau^+_{\epsilon \log N + r}(\mathcal{S}) < \tau^-_0(\mathcal{S})).
\end{align*}
We recall that $\mathcal{S}_0 = S^N_0 - (1-\cexp{inter_integ})^{-1} \log N$. By Lemma~\ref{lem:R0r0}, we have $S^N_0 \leq (1-\cexp{inter_integ})^{-1}(1+\cinit{small}) \log N$. 
Noting that $ \cinit{small} (1-\cexp{inter_integ})^{-1} < \epsilon = 1+\cexp{interface} - (1-\cexp{inter_integ})^{-1}$ by \eqref{bounds_cinitsmall}, it now follows from Lemma~\ref{lemma_gamblerbd}
that there exists some $a>0$ such that for any $r \in \R$, for sufficiently large $N$,
\[ \bP^\mathcal{S}_{\mathcal{S}_0}(\tau^+_{\epsilon \log N + r}(\mathcal{S}) < \tau^-_0(\mathcal{S}))  < N^{-a}.\]
%
%
This completes the proof. \end{proof}

We can now give the proof of the main result of this section. 

\begin{proof}[Proof of Proposition~\ref{prop_R_main}]
The statement concerns the solution $v^N_t$. We consider instead the corresponding solution $w^N_t$ to \eqref{eq_spde_movingframe}, which satisfies $v^N_t = w^N_{Nt}$ for all $t \geq 0$. Since $R^N_t = R(w^N_t)$ and $r^N_t = r_N(w^N_t)$, we have
\begin{align*}
	\bP \left(\sigmaN{interface} \leq \sigma_{2,N} \wedge \sigma_{3,N} \wedge T \right)\leq \bP \bigg(\sup_{t \in [0,NT \wedge \hat{\sigma}_{2,N} \wedge \hat{\sigma}_{3,N}]} (R^N_t \vee r^N_t) \geq (1+ \cexp{interface}) \log N  \bigg).
\end{align*}
We will show that the probability above vanishes as $N \to \infty$. Let $\delta > 0$ be small enough so that Lemmas~\ref{lemma_coupling} and \ref{lemma_RW_bound} hold and fix $\Lambda = \bar{\Lambda}(\delta)$. For $n \in \N$, recall the interval 
\[I^N_t = [((n-1)\Lambda) \wedge \hat{\sigma}_{2,N} \wedge \hat{\sigma}_{3,N},(n\Lambda) \wedge \hat{\sigma}_{2,N} \wedge \hat{\sigma}_{3,N}]\] 
from \eqref{def_INn}, and the definitions of $\mathcal{R}^N_n$ and $\mathfrak{R}^N_n$ from \eqref{def_calRn} and \eqref{def_frakRn} as well as the event $G^N_n$ defined immediately preceding. It is immediate that for all $n \in \N$,
\begin{equation*}
	\sup_{t \in I^N_n} R^N_t  \leq \mathcal{R}^N_n \indc_{(G^N_n)^c} + (S^N_{n-1}+y_0) \indc_{G^N_n} \leq \left( S^N_{n-1} \vee S^N_n \right) + y_0,
\end{equation*}
where we use \eqref{def_Sn} in the last inequality. For $r^N_t$, we remark that on $G_n^c$ we have 
\begin{equation*}
	\sup_{t \in I^N_n}r^N_t =  \mathfrak{R}_n \leq S^N_n,
\end{equation*}
whereas on $G_n$, the relation \eqref{eq_r_X_link} implies that
\begin{equation*}
	\sup_{t \in [(n-1)\Lambda, n\Lambda]}r^N_t \leq S^N_{n-1} + (1-\cexp{inter_integ})^{-1} \log 4.
\end{equation*}
Combining the three previous inequalities, we obtain that for all $n \in \N$,
\begin{equation*}
	\sup_{t \in I^N_n} \left(R^N_t \vee r^N_t \right) \leq \left( S^N_{n-1} \vee S^N_n\right)  + L,
\end{equation*}
where $L = y_0 \vee ((1-\cexp{inter_integ})^{-1} \log 4)$. 
Hence, if $S^N_n < (1+\cexp{interface}) \log N - L$ for all $n \leq \lceil NT / \Lambda \rceil \wedge \lceil\Lambda^{-1} (\hat{\sigma}_{2,N} \wedge \hat{\sigma}_{3,N})\rceil$, the desired event occurs. This implies that
\begin{align*}
&\bP \bigg(\sup_{t \in [0,NT \wedge \hat{\sigma}_{2,N} \wedge \hat{\sigma}_{3,N} ]} \left(R^N_t \vee r^N_t\right)  \geq (1+ \cexp{interface}) \log N   \bigg)
\\ &\hspace{.8 cm} \leq \bP \bigg( \max_{n \leq \lceil NT / \Lambda \rceil \wedge \lceil\Lambda^{-1} (\hat{\sigma}_{2,N} \wedge \hat{\sigma}_{3,N})\rceil } S^N_n \geq (1+\cexp{interface}) \log N - L \bigg)
\\ &\hspace{.8 cm} \leq \bP \bigg( \max_{n \leq\lceil NT / \Lambda \rceil} \mathcal{S}^\delta_n > ((1+\cexp{interface}) - (1-\cexp{inter_integ})^{-1}) \log N -L - (1+k_\delta)\rho \Lambda \bigg).
\end{align*}
The last inequality uses Lemma~\ref{lemma_coupling}. The result now follows from Lemma~\ref{lemma_RW_bound}. 
\end{proof}

We now return to give the proof of Lemma~\ref{lemma_gamblerbd}, which requires the following elementary result. We recall the constant $q$ defined in \eqref{defeq_q}.

\begin{lemma} \label{lemma_geometricRV}
For sufficiently small $\delta > 0$ and $\Lambda = \bar{\Lambda}(\delta)$,
\[ \mu_\delta(((k-1)\rho\Lambda,\infty)) \leq q^{k+1} \]
for all $k \in \N_0$.
\end{lemma}
\begin{proof} We begin by recalling the definitions of $\mu_\delta$, as well as $\lambda$ and $b(\delta)$, from \eqref{def_mudelta} and the preceding discussion.  
Using Gaussian scaling, we obtain that for any $k \in \N$,
\begin{align*} \mu_\delta(((k-1)\rho \Lambda,\infty)) =  \int_{(\lambda^{-1/2} ((k-1)) \rho\Lambda \vee b(\delta))}^\infty  e^{-x^2/2}dx &\leq  \int_{\lambda^{-1/2} (k-1) \rho\Lambda}^\infty e^{-x^2/2}dx  \\&\leq \sqrt{2\pi}  e^{- c\Lambda (k-1)^2}.
\end{align*}
Since $\Lambda = \bar{\Lambda}(\delta) \geq 1$, clearly there exists $M \in \N$, which is independent of $\delta$, such that 
\begin{equation*}
\mu_\delta(((k-1)\rho \Lambda,\infty)) \leq q^{k+1}, \quad k \geq M.
\end{equation*}
We then choose $\delta >0$ to be small enough so that $\delta < q^{M+1}$, and recall that $\mu_\delta((-\rho\Lambda,\infty)) = \delta$. Then for $k \in \N_0$ satisfying $k \leq M$,
\begin{align*}
\mu_\delta(((k-1)\rho \Lambda,\infty)) \leq \delta < q^{M+1} \leq q^{k+1}.
\end{align*}
Combining the above with the previous condition for $k \geq M$, we obtain the desired inequality for all $k \in \N_0$, which completes the proof.\end{proof}

\begin{proof}[Proof of Lemma~\ref{lemma_gamblerbd}]
Let $\hat{J}$ be a geometric-like random variable with distribution given by
\[\bP(\hat{J} = -\rho \Lambda) = 1-q\]
and 
\[ \bP(\hat{J} = k \rho \Lambda) = (1-q)\cdot q^{k+1} \, \,\text{ for } k\in \N_0.\]
Then for all $k \in \N_0$, $\hat{J}$ satisfies
\[ \bP(\hat{J} \geq x) = q^{k+1} \, \text{ for all } \, x \in ((k-1)\rho\Lambda, k\rho\Lambda].\] 
We now prove that for small enough $\delta$, $\mu_\delta([x,\infty)) \leq \bP(\hat{J} \geq x)$ for all $x \in \R$. Suppose that $\delta \in (0, q]$ is sufficiently small so that Lemma~\ref{lemma_geometricRV} holds. Then for any $k \in \N_0$ and $x \in ((k-1)\rho \Lambda, k\rho \Lambda]$,
\[ \mu_\delta([x,\infty)) \leq \mu_\delta(((k-1)\rho\Lambda,\infty)) \leq q^{k+1} = \bP(\hat{J} \geq x).\]
On the other hand, for $x \leq -\rho \Lambda$, $\mu_\delta([x,\infty)) = \bP(\hat{J} \geq x) = 1$. Thus, we have shown that $\mu_\delta([x,\infty)) \leq \bP(\hat{J} \geq x)$ for all $x \in \R$. In particular, if $J$ is a random variable with distribution $\mu_\delta$, then $J \preceq \hat{J}$. Under $\bP^\mathcal{S}_x$, $\mathcal{S}$ is a random walk whose increments are independent copies of $J$. It is therefore possible to construct a monotone coupling of $\mathcal{S}$ with a random walk whose increments are independent copies of $\hat{J}$. Let $\{\hat{J}_k : k \in \N\}$ be an independent family of random variables with the same law as $\hat{J}$, and define $\mathcal{S}^{\text{Geom}} = (\mathcal{S}^{\text{Geom}}_n)_{n \in \N_0}$ by
\[\mathcal{S}^{\text{Geom}}_n = x + \sum_{k=1}^n \hat{J}_k.\]
If $\bP^{\text{Geom}}_x$ denotes the law of the above process, the existence of a monotone coupling of $\mathcal{S}$ and $\mathcal{S}^{\text{Geom}}$ implies that, for $a < x < b$, 
\begin{equation} \label{eq_couplingprob1}
\bP^\mathcal{S}_x( \tau^+_{b}(\mathcal{S}) <  \tau^-_{a}(\mathcal{S}) ) \leq \bP^{\hat{\mathcal{S}}}_x( \tau^+_{b}(\mathcal{S}^{\text{Geom}}) <  \tau^-_{a}(\mathcal{S}^{\text{Geom}})),
\end{equation}
and thus it suffices to obtain an upper bound for the right-hand side. To do so, we will make use of another coupling which represents $\mathcal{S}^{\text{Geom}}$ via a time-changed biased simple random walk. Let $\mathcal{S}^{\text{SRW}}= (\mathcal{S}^{\text{SRW}}_n)_{n \in \N_0}$ be a biased simple random walk on (a translate of) $\rho \Lambda \Z$ defined by
\[ \mathcal{S}^{\text{SRW}}_n = \mathcal{S}^{\text{SRW}}_0 + \sum_{k=1}^n I_k,\]
with independent and identically distributed increments with law
\[ \bP(I_k = \rho \Lambda) =q, \quad \bP(I_k = - \rho \Lambda) = 1-q\]
for all $k \in \N$. Next, define a non-decreasing functional $(\gamma_n)_{n \in \N_0}$ associated to $\mathcal{S}^{\text{SRW}}$. Set $\gamma_0 = 0$. Then, for each $n \in \N$, we define
\begin{align*}
\gamma_n = \begin{cases} \gamma_{n-1} + 1 & \text{ if } I_n = -\rho \Lambda,
\\ \gamma_{n-1} & \text{ if } I_n = \rho \Lambda. \end{cases}
\end{align*}
In particular, $\gamma_n$ increases by one at a given step if the step is negative, and remains constant if it is positive. We then define the inverse
\[ \Gamma_n := \inf \{k \in \N_0 : \gamma_k = n\}, \quad n \in \N_0.\]
We now claim that $\mathcal{S}^{\text{SRW}}$ time-changed by $\Gamma_n$ is a copy of $\mathcal{S}^{\text{Geom}}$. That is, the process $(\mathcal{S}^{\text{SRW}}_{\Gamma_n})_{n \in \N_0}$ has the dynamics of $(\mathcal{S}^{\text{Geom}}_n)_{n \in \N_0}$, started from $\mathcal{S}^{\text{SRW}}_0$. This can be seen by a straightforward analysis of the transition probabilities. Conditioning on $\Gamma_n$, we see that for any $n,m \in \N_0$,
\begin{align*}
\bP(\mathcal{S}^{\text{SRW}}_{\Gamma_{n+1}} - \mathcal{S}^{\text{SRW}}_{\Gamma_n}= -\rho \Lambda \, | \, \Gamma_n = m) &= \bP(I_{m+1} = -\rho \Lambda ) = 1- q.
\end{align*}
Similarly, for $k \in \N_0$,
\begin{align*}
\bP(\mathcal{S}^{\text{SRW}}_{\Gamma_{n+1}} - \mathcal{S}^{\text{SRW}}_{\Gamma_n}= k \rho \Lambda \, | \, \Gamma_n = m) &= \bP(I_{m+1} = \dots = I_{m+k+1} = \rho \Lambda, I_{k+m+2} = -\rho \Lambda)
\\ & =(1-q) \cdot q^{k+1},
\end{align*}
and hence $(\mathcal{S}^{\text{SRW}}_{\Gamma_n})_{n \in \N_0}$ has the same transition probabilities as $(\mathcal{S}^{\text{Geom}}_n)_{n \in \N_0}$. 


We define $\mathcal{S}^{\text{Geom}}$ via the coupling introduced above, that is $\mathcal{S}^{\text{Geom}}_n := \mathcal{S}^{\text{SRW}}_{\Gamma_n}$ for all $n \in \N_0$. We then claim that, if $\mathcal{S}^{\text{Geom}}_0= \mathcal{S}^{\text{SRW}}_0= x$ and $a < x < b$, then 
\begin{equation} \label{eq_walkcoupleevents}
\big\{ \tau^+_b(\mathcal{S}^{\text{Geom}}) <  \tau^-_a(\mathcal{S}^{\text{Geom}})  \big\} = \big\{ \tau^+_{b + \rho \Lambda}(\mathcal{S}^{\text{SRW}}) <  \tau^-_a(\mathcal{S}^{\text{SRW}} )\big\}.
\end{equation}
Indeed, since $\gamma_n$ increases with each down-step of $\mathcal{S}^{\text{SRW}}$, it is easy to argue that \[\tau^-_a(\mathcal{S}^{\text{Geom}}) = \gamma_{\tau^-_a(\mathcal{S}^{\text{SRW}})}.\]
On the other hand, if we define $\tau^*_b(\mathcal{S}^{\text{SRW}}) = \inf\{n \in \N_0 : b \leq \mathcal{S}^{\text{SRW}}_n < \mathcal{S}^{\text{SRW}}_{n-1} \}$,
we then claim that
\[ \tau^+_b(\mathcal{S}^{\text{Geom}}) = \gamma_{\tau^*_b(\mathcal{S}^{\text{SRW}})}.\]
Indeed, $\tau^*_b(\mathcal{S}^{\text{SRW}})$ is the first time that $\mathcal{S}^{\text{SRW}}$ is above $b$ after completing a down-step, and the above then follows from the definition of $\gamma_n$. Since $n \mapsto \gamma_n$ increases on down-steps of $\mathcal{S}^{\text{SRW}}$, and both $\tau^*_b(\mathcal{S}^{\text{SRW}})$ and $ \tau^-_a(\mathcal{S}^{\text{SRW}})$ occur on down-steps (and the two stopping times are not equal), the two previous equations imply that
\begin{equation*}
\big\{ \tau^+_b(\mathcal{S}^{\text{Geom}}) <  \tau^-_a(\mathcal{S}^{\text{Geom}})  \big\} = \big\{ \tau^*_{b}(\mathcal{S}^{\text{SRW}}) <  \tau^-_a(\mathcal{S}^{\text{SRW}})  \big\}.
\end{equation*}
Finally, we remark that, since the increments of $\mathcal{S}^{\text{SRW}}$ have magnitude identically $\rho \Lambda$, $\tau^*_{b}(\mathcal{S}^{\text{SRW}})$ is exactly the time of the first down-step made by $\mathcal{S}^{\text{SRW}}$ after $\tau^+_{b+\rho\Lambda}(\mathcal{S}^{\text{SRW}})$. Since $\mathcal{S}^{\text{SRW}}$ remains above $b$, and hence $a$, between $\tau^+_{b+\rho \Lambda}(\mathcal{S}^{\text{SRW}})$ and $\tau^*_{b}(\mathcal{S}^{\text{SRW}})$, it follows that
\begin{equation*}
 \big\{ \tau^+_{b+\rho \Lambda }(\mathcal{S}^{\text{SRW}}) <  \tau^-_a(\mathcal{S}^{\text{SRW}})  \big\} = \big\{ \tau^*_{b}(\mathcal{S}^{\text{SRW}}) <  \tau^-_a(\mathcal{S}^{\text{SRW}})  \big\}.
\end{equation*}
Combining this with the previous equality of events proves \eqref{eq_walkcoupleevents}, which along with \eqref{eq_couplingprob1} implies that
\begin{equation} \label{eq_couplingprob2}
\bP^\mathcal{S}_x( \tau^+_{b}(\mathcal{S}) <  \tau^-_{a}(\mathcal{S}) ) \leq \bP^{\text{SRW}}_x( \tau^+_{b+\rho\Lambda}(\mathcal{S}^{\text{SRW}}) <  \tau^-_{a}(\mathcal{S}^{\text{SRW}}) ),
\end{equation}
where $\bP^{SRW}_x$ is the law of $\mathcal{S}^{\text{SRW}}$ started from $x$. We now conclude the proof by proving the desired bound for $\tilde{\mathcal{S}}$, which is a biased simple random walk on $x + \rho \Lambda \Z$. First, assuming for simplicity that that $a, b \in x + \rho \Lambda \Z$, we can compute the probability above using the classical gambler's ruin formula for biased random walks. 
Adapting these to the setting of the lattice spaced by $\rho \Lambda$, we obtain
\begin{align} \label{eq_gamblerlattice}
\bP^{\mathcal{S}^{\text{SRW}}}_x( \tau^+_{b+\rho\Lambda}(\mathcal{S}^{\text{SRW}}) <  \tau^-_{a}(\tilde{\mathcal{S}}) ) = \frac{\left( \frac{1-q}{q}\right)^{(x-a)/(\rho\Lambda)} - 1 }{\left(  \frac{1-q}{q}\right)^{ (b - a)/(\rho\Lambda)+1} - 1 }.
\end{align}
For general $a, b \in \R$, define $a_-$ and $b_-$ to be the smallest points in $x + \rho\Lambda \Z$ such that $a_- \leq a$ and $b_- \leq b$. Then starting from $x \in (a,b)$, we have 
\[ \tau^+_{a }(\mathcal{S}^{\text{SRW}}) \leq \tau^+_{a_- }(\mathcal{S}^{\text{SRW}}) \,\, \text{ and } \,\, \tau^+_{b_-+\rho\Lambda}(\mathcal{S}^{\text{SRW}}) < \tau^+_{b +\rho\Lambda}(\mathcal{S}^{\text{SRW}}).\]
Hence, by \eqref{eq_gamblerlattice} and the fact that $|a_- - a| \leq \rho \Lambda$ and $|b_- - b| \leq \rho \Lambda$, we obtain 
\begin{align*}
\bP^{\mathcal{S}^{\text{SRW}}}_x( \tau^+_{b+\rho\Lambda}(\mathcal{S}^{\text{SRW}}) <  \tau^-_{a}(\mathcal{S}^{\text{SRW}}) )  
\leq \bP^{\mathcal{S}^{\text{SRW}}}_x( \tau^+_{b_- +\rho\Lambda}(\mathcal{S}^{\text{SRW}}) <  \tau^-_{a_-}(\mathcal{S}^{\text{SRW}}) ) \leq \frac{\left( \frac{1-q}{q}\right)^{(x-a)/(\rho\Lambda)+1} - 1 }{\left(  \frac{1-q}{q}\right)^{(b - a)/(\rho\Lambda)} - 1}.
\end{align*}
Substituting this inequality into \eqref{eq_couplingprob2} completes the proof. \end{proof}

\subsection{Interface estimates} \label{s_krylov}
In this subsection, we obtain probabilistic estimates on the position of the right endpoint of the support over a given time interval for solutions to \eqref{eq_spde_movingframe} with initial data compactly supported to the right. We adapt the method of Krylov \cite{krylov_result_1997} to include moving boundaries, which will be crucial for the proof of Proposition~\ref{prop_SN_inc}. The main result is Lemma~\ref{lemma_krylov2}, which along with Corollary~\ref{corollary_endpoint} will allow us to control the right endpoint. Throughout this section, we work under $\bP_w$, where $w \in C(\R,[0,1])$ satisfies $R(w) < \infty$.

The argument makes use of the weak form of the SPDE, which for \eqref{eq_spde_movingframe} is as follows: for any $\phi \in C^{1,2}_c([0,\infty) \times \R)$, with probability one, for all $t \geq 0$ we have
\begin{align} \label{eq_weakform_movingframe}
\langle w^N_t,\phi_t \rangle = \langle w, \phi_0 \rangle + \int_0^t \langle  w^N_s, \partial_{xx} \phi_s -\alpha \partial_x \phi_s + \dot{\phi}_s \rangle ds + \int_0^t \langle f(w^N_s), \phi_s \rangle ds +  M_t^\phi,
\end{align}
where we denote $\phi_s = \phi(s,\cdot)$ and $\dot{\phi}_s = \partial_s \phi_s$, and define
\begin{equation} \label{phi_mart_right}
M_t^\phi = \frac{1}{\sqrt N}\int_0^t \int_\R \phi(s,x) \sqrt{w^N_s(x)(1-w^N_s(x))} \dot{W}(dsdx).
\end{equation}
Let $z_0 \in \R$ and $V > 0$. We define a linear curve $z(t) = z_0 - V t$ for $t \geq 0$ and corresponding functional
\begin{equation} \label{eq_Adef}
A_t(z_0,V) := \int_0^t w^N_s(z(s))ds, \quad t \geq 0.
\end{equation}
We then define the test function
\begin{align} \label{def_phimoving}
     \phi_{z_0,V}(t,x) := e^{(\alpha/2 - V) (x - z(t))}(x-z(t))_+,\quad t\geq0, \,x\in \R.
\end{align}
To simplify notation in the ensuing calculations, we write $\phi(t,x) = \phi_{z_0,V}(t,x)$. Elementary calculations lead to the following expressions for the (distributional) derivatives of $\phi$: 
\begin{align*}
    \partial_x \phi(t,x) &= \left[ e^{ (\alpha/2 - V) (x-z(t))} + \left(\frac \alpha 2 -V\right)e^{ (\alpha/2 - V) (x-z(t))}(x-z(t))  \right] 1_{(z(t),\infty)}(x), \\
    \partial_{xx} \phi(t,x) &= \Big[( \alpha - 2V) e^{ (\alpha/2 - V) (x-z(t))} + \Big(\frac{\alpha^2}{4} + V^2 - \alpha V\Big)(x-z(t)) e^{ (\alpha/2 - V) (x-z(t))} \Big] 1_{(z(t),\infty)}(x) 
    \\ & \quad + \delta_0 (x-z(t)),
    \\ \dot{\phi}(t,x) &= \left[ V e^{(\alpha/2 - V) (x - z(t))} + V\left(\frac \alpha 2 - V\right)e^{ (\alpha/2 - V) (x-z(t))} (x-z(t)) \right]1_{(z(t),\infty)}(x) .
  \end{align*}
In the expression for $\partial_{xx} \phi(t,x)$, $\delta_0$ denotes the standard Dirac delta function. Using the above, for $x > z(t)$ we have
\begin{align} \label{eq:phiVderivatives}
    &\partial_{xx} \phi (t,x) - \alpha \partial_x  \phi(t,x) + \dot{\phi}(t,x) \notag 
    \\ &\hspace{1 cm}= e^{(\alpha/2 - V) (x-z(t))} \left[\alpha - 2V - \alpha + V\right] \notag
    \\ &\hspace{1.5 cm} + (x-z(t))e^{(\alpha/2 - V)(x-z(t))} \left[\frac{\alpha^2}{4} + V^2 - \alpha V  - \frac{\alpha^2}{2} + \alpha V + \frac{V \alpha}{2} - V^2\right] \notag
    \\ &\hspace{1 cm}= -V e^{(\alpha/2 - V) (x-z(t))} +  \left[- \frac{\alpha^2}{4} + \frac{V\alpha}{2}\right]\phi(t,x).
\end{align}
Using the weak form with $\phi = \phi_{z_0,V}$ therefore leads to the following.

\begin{lemma} \label{lemma_weak_phiV} Let $w \in C(\R,[0,1])$ and $z_0 \in \R$ satisfy $R(w) < z_0$. Then for any $V>0$, $\bP_w$-a.s.,
\begin{align} \label{e_ibp_A}
&\int_\R \phi_{z_0,V}(t,x) w^N_t(x) dx \leq A_t(z_0,V) + M_t^{\phi_{z_0,V}}   \notag
\\ &\qquad \qquad + \int_0^t \int_{z(s)}^\infty  \left[ f(w^N_s(x)) - \left( \frac{\alpha^2}{4} - \frac{V\alpha}{2} \right)w^N_s(x)\right] \phi_{z_0,V}(s,x) ds dx \quad  \forall \,t  \geq 0,
\end{align} 
where $M_t^{\phi_{z_0,V}}$ is defined as in \eqref{phi_mart_right}. \\
\end{lemma}
\begin{proof}
The weak form \eqref{eq_weakform_movingframe} cannot be applied directly with $\phi = \phi_{z_0,V}$, so a limiting argument is required. This is relatively routine, and we only sketch the details.
Let $(\phi^k)_{k \in \N}$ be a family of compactly supported test functions such that $\phi^k(t,x) = \phi(t,x)$ on $(-\infty, z(t) + k]$. (These can be chosen because $\phi(t,x) = 0$ for $x \leq z(t)$.)
Since the second derivative of $\phi^k$ is a measure, we still cannot apply the weak form directly. 
It is routine, however, to approximate $\phi^k$ over $[0,\Lambda]$ with a family of functions $(\phi^k_n)_{n \in \N} \subset C^{1,2}_c([0,\Lambda] \times \R)$, and pass to the limit for any $\Lambda>0$.
For the details, see \cite[Lemma~3.1]{krylov_result_1997}. This gives the weak form for $\phi^k$ for every $k$, including the distributional term corresponding to $A_{t}(z_0,V)$, i.e. 
\begin{align} \label{eq:weakapproxk}
\int_\R \phi^k(t,x) w^N_t(x) dx &=  \int_0^t \int_{(z(s),z(s) + k]} (\partial_{xx} \phi^k (s,x) - \alpha \partial_x  \phi^k(s,x) + \dot{\phi}^k(s,x)) w^N_s(x) dxds  \notag
\\ &\qquad + A_t(z_0,V) + M_t^{\phi^k}  + \int_0^t \int_{z(s)}^\infty f(w^N_s(x))  \phi^k(s,x) dxds. 
\end{align}
In order to let $k \to \infty$, we use a soft argument with the version of the compact interface property which we already have. 
By Lemma~\ref{lem:compactinterface}, $\sup_{t \in [0,\Lambda]} R(w^N_t) < \infty$ a.s. for any $\Lambda > 0$.
We recall that $z(t) = z_0 - Vt$. Then for any $\Lambda>0$, since $\phi^k(t,x) = \phi(t,x)$ for $x \in (-\infty, z(t) + k]$, on $\{\sup_{[0,\Lambda]} R(w^N_t) \leq z_0 - V\Lambda + k\}$, \eqref{eq:weakapproxk}
is equivalent to 
\begin{align*} 
\int_\R \phi(t,x) w^N_t(x) dx &=  \int_0^t \int_{(z(s)}^\infty (\partial_{xx} \phi (s,x) - \alpha \partial_x  \phi(s,x) + \dot{\phi}(s,x)) w^N_s(x) dxds  \notag
\\ &\qquad + A_t(z_0,V) + M_t^{\phi}  + \int_0^t \int_{z(s)}^\infty f(w^N_s(x))  \phi(s,x) dxds,
\end{align*}
for all $t \in [0,\Lambda]$, because the integrand in each term vanishes outside of the region where $\phi^k = \phi$. The fact that $\sup_{t \in[0,\Lambda]} R(w^N_t) < \infty$ a.s. therefore implies that
the above holds for all $t \in [0,\Lambda]$ a.s. The result now follows from \eqref{eq:phiVderivatives}.
\end{proof}


The integration by parts formula \eqref{e_ibp_A} will be particularly useful when we can show that the last integral term is non-positive. As we will see, for certain values of $V$ this can be established up to a natural stopping time. A short calculation linearizing $f$ at $0$ shows that this can be achieved for all $V < \frac \alpha 2 + \frac 2 \alpha -2$. Here, in order to more easily connect our results with arguments involving the stopping time $\hat{\sigma}_{2,N} \wedge \hat{\sigma}_{3,N}$, and to reduce the proliferation of definitions, we will use a slightly more restricted set of admissible velocities $V$. We define
\begin{equation} \label{def_Vstar}
	V_{*} = \frac \alpha 2 + \frac 2 \alpha (1-q_1) - 2>0,
\end{equation}
where the positivity can be seen as a consequence of \eqref{def_a2a3}. The above is equivalent to 
\[-\left(\frac{\alpha^2}{4} - \frac{ V_*\alpha}{2}\right) = 1- \alpha - q_1.\]
%
In particular, by \eqref{eq_f_small1} this implies that for any $V \in [0,V_*]$,
\begin{equation} \label{eq_Vcancel}
	\text{For all $u \in [0,u_*]$,} \quad f(u) - \left( \frac{\alpha^2}{4} - \frac{V\alpha}{2} \right)u \leq 0.
\end{equation}
For $x \in \R$, we then define the stopping time
\begin{equation} \label{def_GammaN}
	\Gamma_N(x) =  \inf \{t \geq 0 : \, \exists \,y \geq x \, \text{ such that } w^N_t(x) > u_* \}.
\end{equation}
We make an observation for later use, which is
\begin{align} \label{eq_gammaN_sigmaN}
	&\text{If $\lim_{N \to \infty} x(N) = \infty$, then for sufficiently large $N$, } \\
	&\qquad \hat{\sigma}_{2,N} \wedge \hat{\sigma}_{3,N} \leq  \Gamma_N(x(N)) \text{ $\bP_w$-a.s. for any $w \in C(\R,[0,1])$}. \notag
\end{align}
This follows from \eqref{eq_wNsmall_tail}.

The following is an immediate consequence of \eqref{eq_Vcancel}, \eqref{def_GammaN} and Lemma~\ref{lemma_weak_phiV}.

\begin{lemma} \label{lemma_weak_phiVbd} Let $w \in C(\R,[0,1])$, $z_0 > R(w)$, $V \in [0,V_*]$ and $\Lambda>0$. Then $\bP_w$-a.s.,
\begin{align} \label{e_ibp_A_eq_rightmove}
\int_\R \phi_{z_0,V}(t,x) w^N_t(x) dx  &\leq A_t(z_0,V)+  M_t^{\phi_{z_0,V}} \quad \text{for all }\, t  \in [0, \Lambda \wedge \Gamma_N(z_0 - V \Lambda) ],
\end{align} 
where $M_t^{\phi_{z_0,V}}$ is defined as in \eqref{phi_mart_right}. \\
\end{lemma}

The next result establishes that if $A_\Lambda(z_0,V) = 0$, then the right endpoint of the support $R^N_t$ (recall \eqref{def_RNrN}) remains below $z(t)$ for $t \in [0,\Lambda]$.

\begin{cor} \label{corollary_endpoint} Let $w \in C(\R,[0,1])$, $z_0 > R(w)$, $V \in [0,V_*]$ and $\Lambda>0$. Then
\begin{equation*}
\bP_w\left(R^N_t > z_0 - Vt \text{ for some } t \in [0,\Lambda \wedge \Gamma ], A_{\Lambda \wedge \Gamma}(z_0,V) = 0 \right) = 0,
\end{equation*}
where $\Gamma = \Gamma_N(z_0 - V \Lambda)$.
\end{cor}
\begin{proof}
Let $w, z_0, V$ and $\Lambda$ be as in the statement and write $\Gamma = \Gamma_N(z_0 - \Lambda V)$. Let $\gamma = \inf\{t \geq 0 :  A_t(z_0,V)  > 0\}$. Then by \eqref{e_ibp_A_eq_rightmove} and continuity of $t\mapsto A_t(z_0,V)$, we have
\[ M^{\phi_{z_0,V}}_{t \wedge \gamma \wedge \Gamma} \geq - A_{t\wedge \gamma \wedge \Gamma}(z_0,V) = 0\]
for all $t \in [0,\Lambda]$. Hence $t \mapsto M^{\phi_{z_0,V}}_{t \wedge \gamma \wedge \Gamma}$ is a non-negative local martingale, and therefore a supermartingale. A supermartingale started at $0$ with vanishing negative part must a.s.~be identically zero, whence it follows that with probability one,
\begin{equation*}
\int_\R \phi_{z_0,V}(t,x) w^N_t(x) dx  \leq A_t(z_0,V)+  M_t^{\phi_{z_0,V}} = 0  \quad \text{ for } t \in [0, \Lambda \wedge \gamma \wedge \Gamma].
\end{equation*}
Since the left-hand side is non-negative, this proves that it is identically zero for $t \leq \Lambda \wedge \gamma \wedge \Gamma$. Finally, we note that, because $w^N$ is continuous and $\phi_{z_0,V}(t,\cdot)$ is strictly positive on $(z(t),\infty)$,  the left-hand side equals zero if and only if $w^N_t(x) = 0$ for all $x > z(t)$, and hence it is identically zero until time $\Lambda \wedge \gamma \wedge \Gamma$ if and only if $R^N_t \leq z(t)$ for all $t \leq  \Lambda \wedge \gamma \wedge \Gamma$. This proves the result. \end{proof}

We now prove a key estimate.

\begin{lemma} \label{lemma_krylov1}
For every $\beta \in (0,1)$, there exists $C_1 = C_1(\beta)$, such that the following holds: suppose that $w \in C(\R,[0,1])$ and $z_0 \in \R$ satisfy $R(w) < z_0$, let $V \in [0, V_* \wedge (\alpha/2)]$, $\Lambda > 0$, and $a,b,r \in (0,1]$, and let $\sigma_0$ be a stopping time satisfying $\sigma_0 \leq \Gamma_N(z_0-V\Lambda)$ a.s.; then there exists $z_1 \in [z_0 + r, z_0 + 2r]$ such that
\begin{equation*}
\bP_w \left(A_{\Lambda \wedge \sigma_0}(z_1,V) \geq a, \, \Lambda < \sigma_0 \right) \leq \bP_w\left(A_{\Lambda \wedge \sigma_0}(z_0,V) \geq b, \, \Lambda < \sigma_0 \right) + C_1  N^{\beta/2}  r^{-3\beta/2}  \left(\frac{b}{a^{1/2}} \right)^\beta.
\end{equation*}
\end{lemma}

\begin{proof}
We fix $w, z_0, V, \Lambda$ and $\sigma_0$ as in the statement, and we write $M_t$ to denote $M^{\phi_{z_0,V}}_t$ throughout the proof. By Lemma~\ref{lemma_weak_phiVbd} and the assumption that $\sigma_0 \leq \Gamma_N(z_0-V\Lambda)$ a.s., with probability one we have
\begin{equation} \label{e_krylemma1_0}
  M_{t \wedge \sigma_0} \geq - A_{t \wedge \sigma_0}(z_0,V), \quad t \in [0,\Lambda].
  \end{equation}
Now fix $\beta \in (0,1)$, and let $a,b \in (0,1]$. Define $\tau = \inf \{t > 0 : A_t(z_0,V) \geq b\}$. Let $z > z_0$. Partitioning on $\{\tau \leq \Lambda  \wedge \sigma_0 \}$ and using the fact that $t \mapsto A_t(z,V)$ is non-decreasing, we have
\begin{align}  \label{eq_krylemma1_1}
\bP_w(A_{\Lambda \wedge \sigma_0}(z,V) \geq a)  &= \bP_w(A_{\Lambda \wedge \sigma_0}(z,V)\geq a, \tau \leq {\Lambda \wedge \sigma_0}) +  \bP_w(A_{\Lambda \wedge \sigma_0}(z,V) \geq a, {\Lambda \wedge \sigma_0} < \tau) \notag
\\ & \leq \bP_w(A_{\Lambda \wedge \sigma_0}(z_0,V)  \geq b ) +  \bP_w(A_{\Lambda \wedge \tau \wedge \sigma_0}(z,V) \geq a). 
\end{align}
By Markov's inequality,
\begin{align*} \bP_w(A_{\Lambda \wedge \tau \wedge \sigma_0}(z,V) \geq a) &\leq a^{-\beta /2} \mathbb{E}_w [A_{\Lambda \wedge \tau \wedge \sigma_0}(z,V)^{\beta /2}]
\\ & =a^{-\beta /2} \mathbb{E}_w \bigg[ \bigg(\int_0^{\Lambda \wedge \sigma_0 \wedge \tau} w^N_t(z - Vt) dt \bigg)^{\beta /2 }\bigg].
\end{align*}
Let $r \in (0,1]$. To handle the second term in \eqref{eq_krylemma1_1}, we integrate the above with respect to $z$ over $[z_0 + r, z_0+ 2r]$ to obtain
\begin{align} \label{e_kry_tech1}
\frac 1 r \int_{ z_0 + r}^{z_0 + 2r} &\bP_w(A_{\Lambda \wedge \tau \wedge \sigma_0}(z,V) \geq a) dz  \notag
\\ &\leq  \frac 1 r \int_{ z_0 + r}^{z_0 + 2r} a^{-\beta /2}  \mathbb{E}_w \bigg[ \bigg(\int_0^{\Lambda \wedge \sigma_0 \wedge \tau} w^N_t(z - Vt) dt \bigg)^{\beta /2} \bigg] dz \notag
\\ &= a^{-\beta /2 } \mathbb{E}_w \bigg[ \frac 1 r \int_{ z_0 + r}^{z_0 + 2r}  \bigg( \int_0^{\Lambda \wedge \sigma_0 \wedge \tau} w^N_t(z - Vt) dt \bigg)^{\beta /2} dz \bigg] \notag
\\ &\leq a^{-\beta /2 } \mathbb{E}_w \bigg[  \bigg(\frac 1 r \int_{ z_0 + r}^{z_0 + 2r}   \int_0^{\Lambda \wedge \sigma_0 \wedge \tau} w^N_t(z - Vt) dt dz \bigg)^{\beta  /2}  \bigg]\notag
\\ &= a^{-\beta /2 } \mathbb{E}_w \bigg[  \bigg(\frac 1 r  \int_0^{\Lambda \wedge \sigma_0 \wedge \tau}  \int_{z_0 - Vt + r}^{z_0 - Vt + 2r}w^N_t(x) dxdt \bigg)^{\beta  /2}  \bigg].
\end{align}
The second inequality follows from Jensen's inequality with concave function $u \mapsto u^{\beta / 2}$. 
We note that if $x \in [z_0 - Vt + r,z_0 - Vt + 2r]$, then by \eqref{def_phimoving},
\[ 1\leq \frac{(x - z_0 + Vt)^2}{r^2} \leq \frac{(x - z_0 + Vt)^2}{r^2} e^{(\alpha - 2V)(x-z_0 + Vt)} =  r^{-2} \phi_{z_0,V}(t,x)^2.\]
Here we use $V \leq \alpha/2$). Note that for all $t \leq \Lambda$ and $x \in [z_0 - Vt + r, z_0 - Vt + 2r]$, we have $x \geq z_0 - V\Lambda$.
Since $\sigma_0 \leq \Gamma_N(z_0-V\Lambda)$, by \eqref{def_GammaN} we have $w^N_t(x) \leq u_* \leq \frac 1 2$, and hence $w^N_t(x) \leq 2 w^N_t(x)(1-w^N_t(x))$.
Using the two bounds above, we conclude that
\begin{align*}
&\text{For all }\, (t,x) \in [0, \Lambda \wedge \sigma_0 \wedge \tau] \times [z_0 - Vt +r, z_0 - Vt + 2r].\\
& \qquad \qquad w^N_t(x) \leq 2  r^{-2} \phi_{z_0,V}(t,x)^2w^N_t(x)(1-w^N_t(x)).
\end{align*}
Substituting the above into \eqref{e_kry_tech1}, we obtain
\begin{align*}
&\frac 1 r \int_{ z_0 + r}^{z_0 + 2r} \bP_w(A_{\Lambda \wedge \tau \wedge \sigma_0}(z,V)\geq a)dz
\\ &\hspace{.5 cm} \leq a^{-\beta /2} \mathbb{E}_w \bigg[  \bigg( \frac{2}{r^3}  \int_0^{\Lambda \wedge \sigma_0 \wedge \tau}  \int_{z_0 - Vt + r}^{z_0 - Vt + 2r} \phi_{z_0,V}(t,x)^2 w^N_t(x)(1-w^N_t(x)) dx dt \bigg)^{\beta/2}  \bigg]
\\ &\hspace{.5 cm}=a^{-\beta /2 } 2^{\beta/2} r^{-3\beta/2} \mathbb{E}_w \bigg[  \bigg( \int_0^{\Lambda \wedge \sigma_0 \wedge \tau}  \int_{z_0 - Vt + r}^{z_0 - Vt + 2r}\phi_{z_0,V}(t,x)^2 w^N_t(x)(1-w^N_t(x)) dx dt\bigg)^{\beta/2}  \bigg].
\end{align*}
We now remark that, because $M_t = M^{\phi_{z_0,V}}_t$, as defined in \eqref{phi_mart_right}, the integrand in the above is the same integrand that appears in $\langle M\rangle_t$, the quadratic variation of $M_t$. In particular, we have
\begin{equation*}
\langle M \rangle_{\Lambda \wedge \sigma_0 \wedge \tau} = \frac 1 N \int_0^{\Lambda \wedge \sigma_0 \wedge \tau}  \int_\R \phi_{z_0,V}(t,x)^2 w^N_t(x)(1-w^N_t(x)) dxdt.
\end{equation*}
Hence, it follows that 
\begin{align} \label{eq_krylemma1_2}
&\frac 1 r \int_{z_0+ r}^{z_0+2r} \bP_w(A_{\Lambda \wedge \tau \wedge \sigma_0}(z,V) \geq a)dz \leq  a^{-\beta /2 }2^{\beta/2} r^{-3\beta/2} N^{\beta /2} \mathbb{E}_w [( \langle M \rangle_{\Lambda \wedge \sigma_0 \wedge \tau} )^{\beta/2} ].
\end{align}
To complete the proof, it suffices to bound the term $\mathbb{E}_w [( \langle M\rangle_{\Lambda \wedge \sigma_0 \wedge \tau} )^{\beta/2} ]$. We first remark that, since $M_t$ is a continuous local martingale, by the Burkholder-Davis-Gundy inequality there is a constant $C_\beta >0$ such that
\begin{equation*}
\mathbb{E}_w [( \langle M\rangle_{\Lambda \wedge \tau \wedge \sigma_0} )^{\beta/2} ]\leq C_\beta \mathbb{E}_w \bigg[ \sup_{t\in [0, \Lambda \wedge \tau \wedge \sigma_0]} |M_t|^\beta \bigg] = C_\beta \mathbb{E}_w \bigg[ \sup_{t \in [0,\Lambda]} |M_{t \wedge \sigma_0 \wedge \tau}|^\beta \bigg].
\end{equation*}
We will bound the right-hand side by applying \cite[Theorem III.6.8]{krylov_introduction_1994}, the statement of which we reproduce here for convenience. We then use its implication to complete the proof, before returning to verify that it can be applied. The result from \cite{krylov_introduction_1994} is as follows.

	\begin{theorem*}
 	Let $\eta_t$, $\zeta_t$ be continuous $\cF_t$-adapted processes. Suppose that $\eta_t \geq 0$, $\zeta_t \geq 0$, $\zeta_0 = 0$, and for all bounded stopping times $\sigma'$,
 	\[ \mathbb{E} [\eta_{\sigma'}] \leq \mathbb{E} [\zeta_{\sigma'}] .\]
 	Then for any $\beta \in (0,1)$ and any stopping time $\tau'$,
 	\[ \mathbb{E} \bigg[ \sup_{t} \eta^\beta_{t \wedge \tau'} \bigg] \leq \frac{2-\beta}{1-\beta} \mathbb{E} \bigg[ \sup_t \zeta^\beta_{t \wedge \tau'} \bigg].\]
 \end{theorem*}

We apply the above with $\eta_t = |M_{t \wedge \Lambda \wedge \sigma_0}|$ and $\zeta_t = 2 A_{t \wedge \Lambda \wedge \sigma_0}(z_0,V)$ at the stopping time $\tau' = \tau =\inf \{t \geq 0 : A_{t \wedge \sigma_0}(z_0,V) \geq b\}$. Since $\zeta_t$ is increasing, the supremum on the right hand side is $2A_{\Lambda \wedge \sigma_0 \wedge \tau}(z_0,V)$, which is at most $2b$ by definition of $\tau$. Thus, our application of the theorem yields
\begin{equation} \label{e_kry_uniformbd}
\mathbb{E}_w \bigg[ \sup_{t\in [0,\Lambda]} |M_{t \wedge \sigma_0 \wedge \tau}|^\beta \bigg] \leq \frac{2^\beta (2-\beta)}{1-\beta} b^\beta.
\end{equation}
Substituting this into the previous inequality and using this bound in \eqref{eq_krylemma1_2}, we obtain that for some $C_1 = C_1(\beta)$, 
\begin{align*} 
&\frac 1 r \int_{z_0+r}^{z_0+2r} \bP_w(A_{\Lambda \wedge \tau \wedge \sigma_0}(z,V) \geq a)dz  \leq C_1 r^{-3\beta/2}  N^{\beta /2} \left( \frac{b}{a^{1/2}}\right)^{\beta}.
\end{align*}
Returning to \eqref{eq_krylemma1_1} and integrating it with respect to $z$ over $[z_0 + r, z_0 + 2r]$, we may conclude that there exists $z_1 \in [z_0 +r, z_0 +2r]$ where the claimed bound is satisfied, which proves the result.

We now justify the application of \cite[Theorem III.6.8]{krylov_introduction_1994} used to obtain \eqref{e_kry_uniformbd}. It suffices to prove that for every bounded stopping time $\sigma'$,
\begin{equation} \label{e_kry_martbd0}
\mathbb{E}_w[|M_{\sigma' \wedge \sigma_0 \wedge \Lambda }|]\leq \mathbb{E}_w[2A_{ \sigma' \wedge \sigma_0\wedge \Lambda}(z_0,V)].
\end{equation}
Let $M^-_t := - (M_t \wedge 0) \geq 0$ and $M^+_t := M_t \vee 0$ denote the negative and positive parts of $M_t$. By \eqref{e_krylemma1_0}, we have $M^{-}_t \leq A_t(z_0,V)$ for all $t \in [0, \sigma_0 \wedge \Lambda]$ a.s., and hence 
\begin{equation} \label{e_kry_martbd1} 
\mathbb{E}_w[M^{-}_{\sigma' \wedge \sigma_0 \wedge \Lambda }] \leq \mathbb{E}_w[A_{ \sigma' \wedge \sigma_0 \wedge \Lambda}(z_0,V)]
\end{equation} 
for every bounded stopping time $\sigma'$. To see that we can obtain the same bound for $\mathbb{E}_w[M^+_{\sigma' \wedge \sigma_0 \wedge \Lambda }]$, we note that because $(M_t)_{t \geq 0}$ is a local martingale we can take an increasing sequence of stopping times $(\sigma_i)_{i \in \N}$ such that $\sigma_i \to \infty$ a.s. and $(M_{t \wedge \sigma_i})_{t \geq 0}$ is a martingale for each $i \in \N$. By the latter property, and because $\sigma' \wedge \sigma_0 \wedge \Lambda$ is bounded, we have
\[ \mathbb{E}_w[M_{\sigma' \wedge \sigma_0 \wedge \Lambda \wedge \sigma_i}^+] = \mathbb{E}_w[M_{\sigma' \wedge \sigma_0 \wedge \Lambda  \wedge \sigma_i}^-] \]
for each $i \in \N$. We then argue using Fatou's lemma, \eqref{e_kry_martbd1}, and the fact that $t \mapsto A_t(z_0,V)$ is non-decreasing to obtain that $\mathbb{E}_w[M^{+}_{\sigma' \wedge \sigma_0 \wedge \Lambda }] \leq \mathbb{E}_w[A_{\sigma' \wedge \sigma_0 \wedge \Lambda}(z_0,V)]$. Indeed,	
\begin{align*}
\mathbb{E}_w[M^+_{\sigma' \wedge \sigma_0 \wedge \Lambda }] = \mathbb{E}_w \left[\lim_{i \to \infty} M^+_{\sigma' \wedge \sigma_0 \wedge \Lambda  \wedge \sigma_i}\right] &\leq \liminf_{i \to \infty} \mathbb{E}_w[M^+_{\sigma' \wedge \sigma_0 \wedge \Lambda  \wedge \sigma_i}]
\\ &=  \liminf_{i \to \infty} \mathbb{E}_w[M^-_{\sigma' \wedge \sigma_0 \wedge \Lambda  \wedge \sigma_i}]
\\ &\leq  \liminf_{i \to \infty} \mathbb{E}_w[A_{\sigma' \wedge \sigma_0 \wedge \Lambda \wedge \sigma_i}(z_0,V)]
 \\ &= \mathbb{E}_w[A_{\sigma' \wedge \sigma_0 \wedge \Lambda}(z_0,V)], 
\end{align*}
where the final equality is by monotone convergence. The above and \eqref{e_kry_martbd1} imply that \eqref{e_kry_martbd0} holds for every bounded stopping time $\sigma'$. 
Thus the theorem is applicable as claimed. This completes the proof.
%
%
%
\end{proof}

We now state and prove the main result of the subsection.

\begin{lemma} \label{lemma_krylov2} For every $\beta \in (0,1)$, there exists $C_2 = C_2(\beta)$, such that the following holds: suppose that $w \in C(\R,[0,1])$ and $z_0 \in \R$ satisfy $R(w) < z_0$, let $V \in [0, V_* \wedge (\alpha/2)]$, $\Lambda > 0$, and let $\sigma_0$ be a stopping time satisfying $\sigma_0 \leq \Gamma_N(z_0-V\Lambda)$ a.s.; then for any $\zeta > 0$, there exists $y \in [z_0 + 1/2, z_0 + 1]$ such that
\begin{equation*}
\bP_w(A_{\Lambda \wedge \sigma_0}(y,V) > 0 ) \leq \bP_w(A_{\Lambda \wedge \sigma_0}(z_0,V) > \zeta) + C_2  ( N\zeta )^{\beta/2}.
\end{equation*}
\end{lemma} 
\begin{proof}
Let $w$, $z_0$ and $\Lambda$ be as above, and let $\beta \in (0,1)$ and $\zeta  >0$. For $n \in \N_0$, let $a_n = \zeta e^{-n}$, and for $n \in \N$, let $r_n = \rho(1+n)^{-2}$, where $\rho$ is chosen so that $\sum_{n=1}^\infty r_n = 1/2$. For $n \in \N$, we iteratively apply Lemma~\ref{lemma_krylov1} with $a = a_n$, $b = a_{n-1}$ and $r = r_n$ to construct a sequence of points $(z_n)_{n \in \N_0}$ for which $A_{\Lambda \wedge \sigma_0}(z_n, V)$ satisfies an estimate in terms of $A_{\Lambda \wedge \sigma_0}(z_{n-1},V)$. Starting from $z_0$, the first application gives a point $z_1 \in [z_0 + r_1,z_0 +  2 r_1]$ such that the estimate from Lemma~\ref{lemma_krylov1} holds with $A_{\Lambda \wedge \sigma_0}(z_{0},V)$ and $A_{\Lambda \wedge \sigma_0}(z_{1},V)$. For the second iteration, starting from $z_1$ we obtain a point $z_2 \in [z_1 + r_2, z_1+ 2r_2]$ such that the estimate holds for $A_{\Lambda \wedge \sigma_0}(z_{1},V)$ and $A_{\Lambda \wedge \sigma_0}(z_{2},V)$. Proceeding in this fashion, the $n$th iteration produces the estimate
\begin{align*}
\bP_w( A_{\Lambda \wedge \sigma_0}(z_{n},V)\geq a_n) &\leq \bP_w(A_{\Lambda \wedge \sigma_0}(z_{n-1},V) \geq a_{n-1})  + C_1N^{\beta/2} \rho^{-3\beta/2}(1+n)^{3\beta}\left(\frac{a_{n-1}}{a_n^{1/2}} \right)^\beta ,
\end{align*}
where $C_1 = C_1(\beta)$ as in the statement of Lemma~\ref{lemma_krylov1}, and $z_n \in [z_{n-1} + r_n, z_{n-1} + 2r_n]$. We remark that $a_{n-1} / a_n^{1/2} = \zeta^{1/2} e^{1 - n/2}$. Substituting the estimate for $A_{\Lambda \wedge \sigma_0}(z_{n-1},V) $ into the above and iterating, we obtain that 
\begin{align*}
&\bP_w(A_{\Lambda \wedge \sigma_0}(z_{n},V)  \geq a_n)  \leq \bP_w(A_{\Lambda \wedge \sigma_0}(z_{0},V) \geq \zeta) + C_1    N^{\beta/2} \zeta^{\beta/2}  \rho^{-3\beta/2}\sum_{k=1}^n  (1+k)^{3\beta}  e^{(1 - k/2)\beta}.
\end{align*}
The sequence $(z_n)_{n \in \N_0}$ is increasing and bounded above by $z_0 + \sum_{n=1}^\infty 2r_n = z_0 + 1$. In particular, there exists $y \in [z_0 + 1/2, z_0 + 1]$ such that $\lim_{n \to \infty}z_n = y$. Because $w^N$ is continuous, $z \mapsto A_{\Lambda \wedge \sigma_0}(z,V)$ is also continuous. In particular, for any $\epsilon >0$, by Fatou's lemma and the estimate above we have
\begin{align*}
\bP_w(A_{\Lambda \wedge \sigma_0}(y,V)> \epsilon) &= \bP_w \left( \liminf_{n \to \infty} A_{\Lambda \wedge \sigma_0}(z_{n},V)> \epsilon \right) 
\\ &\leq \mathbb{E}_w \left( \liminf_{n \to \infty} \indc_{\{A_{\Lambda \wedge \sigma_0}(z_{n},V)> \epsilon/2 \}} \right) 
\\ &\leq \liminf_{n \to \infty}  \bP_w(A_{\Lambda \wedge \sigma_0}(z_{n},V) > \epsilon/2)
\\ & \leq\liminf_{n \to \infty}  \bP_w(A_{\Lambda \wedge \sigma_0}(z_{n},V)) > a_n)
\\ &\leq \bP_w(A_{\Lambda \wedge \sigma_0}(z_{0},V) \geq \zeta) 
\\ &\qquad + C_1   N^{\beta/2} \zeta^{\beta/2} \rho^{-3\beta/2} \sum_{k=1}^\infty  (1+k)^{3\beta}  e^{(1 - k/2)\beta}.
\end{align*}
Since the above holds for all $\epsilon >0$, we may conclude that the same inequality holds when $\epsilon$ is taken to be zero. Since the sum appearing in the right-hand side is convergent, this completes the proof. \end{proof}

\subsection{Some moment bounds}
In this subsection we obtain upper bounds on the expectation of $w^N_t(x)$ under certain assumptions on the initial condition. These bounds will later be used alongside Lemma~\ref{lemma_krylov2} to control $R^N_t$. 
We recall $\cexp{inter_integ}$ from \eqref{bounds_cexpinterinteg} and $q_1$ from \eqref{def_a2a3}. We also note that, under $\bP_w$, the initial state $X^{N,S}_0$ of the process $X^{N,S}$ (from \eqref{def_XNT}) is determined by $w$. The main result of this subsection is the following.

\begin{proposition} \label{prop_rightmoment} There exists $\epsilon'>0$ such that the following holds: for any $\Lambda > 0$, for sufficiently large $N$, for all $S \in \R$ and $w \in C(\R,[0,1])$ such that $R(w) \leq S$ and $X := X^{N,S}_0 \leq N^{1/3}$ under $\bP_w$,
\begin{align*}
&\mathbb{E}_w[w^N_t(x)] 
\\ &\hspace{.5 cm} \leq \min \bigg( N^{-1} X e^{(1-\cexp{inter_integ})(S+2-x)}e^{-(\cexp{inter_integ}(2-\alpha-\cexp{inter_integ}) - q_1)t}  \bP^B_{x+\alpha t - 1 -2(1-\cexp{inter_integ})t}(B_{2t} \leq S) + N^{-1-\epsilon'}, 
\\ &\hspace{1.5 cm} N^{-1} X e^{(1-\cexp{inter_integ})(S+2 - x)}e^{(\|f'\|_\infty +1)t}  \bP^B_{x+\alpha t - 1 -2(1-\cexp{inter_integ})t}(B_{2t} \leq S)  \bigg)
\end{align*}
for all $t \in [0,\Lambda]$ and $x \geq S+1-\alpha t$.
\end{proposition}

The rest of this subsection proves this proposition. The proofs in this section will make use of several different integral representations for $w^N_t(x)$. To this end, we introduce 
\begin{equation*}
I_t w (x) := \int_{\R} G_t(x + \alpha t - y) w(y) dy, \,\, \text{ for } t \geq 0, \, x\in \R.
\end{equation*}
That is, $I_t$ is a shifted heat semigroup. Within this section, we will write 
\[ p= \alpha -1 +q_1\]
and remark that $|p| \leq 2$. Hence by \eqref{eq_f_small1}, we have 
\begin{equation} \label{eq_f_small3}
\text{For all } u \in [0,u_*], f(u) \leq p u.
\end{equation}

We consider the solution with law $\bP_w$ for $w \in C(\R,[0,1])$. By \eqref{mild_spde_wN} with $a = p$, we note that for every $(t,x) \in (0,\infty) \times \R$, it a.s. holds that
\begin{align} \label{eq_wN_pmild}
w^N_t(x) &= e^{pt}I_t(x) + \int_0^t\int_{\R} G_{t-s}(x + \alpha (t-s) - y) e^{p(t-s)}\left( f(w^N_s(y)) - p w^N_s(y)\right) dy ds\notag 
\\ &\quad +N^{-1/2} \int_0^t\int_{\R} G_{t-s}(x + \alpha (t-s) - y) e^{p(t-s)} \sqrt{w^N_s(y)(1-w^N_s(y))} \dot{W}(dsdy).
\end{align}
Using instead $a = \|f'\|_\infty$,
\begin{align*} 
w^N_t(x) &=e^{\|f'\|_\infty t} I_t(x) + \int_0^t\int_{\R} G_{t-s}(x + \alpha (t-s) - y) e^{\|f'\|_\infty(t-s)}\left( f(w^N_s(y)) - \|f'\|_\infty w^N_s(y)\right) dy ds 
\\ &\quad + N^{-1/2}\int_0^t\int_{\R} G_{t-s}(x + \alpha (t-s) - y) e^{\|f'\|_\infty(t-s)} \sqrt{w^N_s(y)(1-w^N_s(y))} \dot{W}(dsdy) .
\end{align*}
Noting that $ f(u) - \|f'\|_\infty u \leq 0$ for all $u \in [0,1]$, we obtain from the above that
\begin{align} \label{eq_wN_fnormmild}
w^N_t(x)  \leq  e^{\|f'\|_\infty t} I_t(x)  + N^{-1/2} \int_0^t\int_{\R} G_{t-s}(x + \alpha (t-s) - y) e^{\|f'\|_\infty(t-s)} \sqrt{w^N_s(y)(1-w^N_s(y))} \dot{W}(dsdy).
\end{align}
As in the statement of Proposition~\ref{prop_rightmoment}, throughout this section we write $X = X^{N,S}_0$. Because $S$ will always be fixed for the duration of a given proof, this will not cause any confusion. We will frequently use the following estimate: for $w \in \cW_S$ and $z \in \R$, under $\bP_w$,
\begin{equation}\label{eq_intervalX_bd}
\int_z^{z+1} w(y)dy \leq e^{-(1-\cexp{inter_integ})z} \int_{z}^{z+1} w(y) e^{(1-\cexp{inter_integ})y}dy  \leq  N^{-1} X e^{(1-\cexp{inter_integ})(S-z)}.
\end{equation}
In order to prove Proposition~\ref{prop_rightmoment}, we need the two following lemmas. 

\begin{lemma} \label{lemma_rightmoment1} Suppose $w \in C(\R,[0,1])$, $S \in \R$, $x \in \R$ and $t \geq 0$ satisfy $R(w) \leq S$ and $x \geq S+1 - \alpha t$. Then under $\bP_w$,
\begin{equation*} I_tw(x) \leq N^{-1} X_0^{N,S} e^{(1-\cexp{inter_integ})(S+2)}e^{(1-\cexp{inter_integ})^2 t - (1-\cexp{inter_integ})(x + \alpha t)}  \bP^B_{x+\alpha t - 1 -2(1-\cexp{inter_integ})t}(B_{2t} \leq S).
\end{equation*}
\end{lemma}

\begin{lemma} \label{lemma_rightmoment2} There exists $\epsilon>0$ such that the following holds: for any $\Lambda >0$, for sufficiently large $N$, for all $S \in \R$ and $w \in C(\R,[0,1])$ satisfying $R(w) \leq S$ and, under $\bP_w$, $X_0^{N,S} \leq N^{1/3}$, 
\begin{equation*}
\bP_w(w^N_t(x) > u^*) \leq N^{-2} 
\end{equation*}
for all $(t,x) \in [N^{-1-\epsilon},\Lambda] \times [S - \epsilon \log N,\infty)$. 
\end{lemma}

Before proving these lemmas, we use them to prove Proposition~\ref{prop_rightmoment}.

\begin{proof}[Proof of Proposition~\ref{prop_rightmoment}]
Suppose that $S$, $w$ and $X = X_0^{N,S}$ are as in the statement and let $\Lambda > 0$. We begin by taking the expectation of \eqref{eq_wN_pmild}. 
For any $x$, a short argument using the Burkholder-Davis-Gundy inequality and Lemma~\ref{lem:wNtailbd} implies that
the martingale term is square integrable over $t \in [0,\Lambda]$, and hence for any $t \in [0,\Lambda]$,
\begin{align*}
\mathbb{E}_w[w^N_t(x)] &= e^{pt} I_t(x) +  \int_0^t\int_{\R} G_{t-s}(x + \alpha (t-s) - y) e^{p(t-s)} \mathbb{E}_w[ f(w^N_s(y)) - p w^N_s(y)] dy ds.
\end{align*}
If $x \geq S+1 - \alpha t$, we may bound the first term in our expression for $\mathbb{E}_w[w^N_t(x)]$ using Lemma~\ref{lemma_rightmoment1}. We also note that if $w^N_s(y) \in [0,u^*]$, then $f(w^N_s(y)) - p w^N_s(y) \leq 0$ by \eqref{eq_f_small3}. It follows (recall that $|p| \leq 2$) that for any $s \geq 0$ and $y \in \R$,
\begin{equation*} \mathbb{E}_w [ f(w^N_s(y)) - p w^N_s(y)] \leq (2+\|f\|_\infty) \bP_w(w^N_s(y) > u^*).\end{equation*}
We may therefore bound this expectation using Lemma~\ref{lemma_rightmoment2} for values of $(s,y)$ where it is applicable, and by $2+\|f\|_\infty$ where it is not. Combining these estimates, we obtain that for sufficiently large $N$, for $t \in [0,\Lambda]$ and $x \geq S+1-\alpha t$,
\begin{align} \label{eq_wNt_momentbdprop1}
\mathbb{E}_w[w^N_t(x)] &\leq e^{pt} N^{-1} X e^{(1-\cexp{inter_integ})(S+2)}e^{(1-\cexp{inter_integ})^2 t - (1-\cexp{inter_integ})(x + \alpha t)}  \bP^B_{x+\alpha t - 1 -2(1-\cexp{inter_integ})t}(B_{2t} \leq S)\notag 
\\ &\quad + (2+\|f \|_\infty)  N^{-2} \int_{N^{-1-\epsilon} \wedge t}^t \int_{S-\epsilon \log N}^\infty G_{t-s}(x + \alpha (t-s) - y) e^{p(t-s)}  dy ds \notag 
\\ &\quad + (2+\|f \|_\infty) \int_{N^{-1-\epsilon} \wedge t}^t \int_{-\infty}^{S - \epsilon \log N} G_{t-s}(x + \alpha (t-s) - y) e^{p(t-s)} dyds \notag
\\ &\quad + (2+\|f \|_\infty) \int_0^{N^{-1-\epsilon} \wedge t} \int_{-\infty}^\infty G_{t-s}(x + \alpha (t-s) - y) e^{p(t-s)} dyds.
\end{align}
We remark that since $t \leq \Lambda$,
\begin{align*}
\int_0^{N^{-1-\epsilon} \wedge t} \int_{-\infty}^\infty G_{t-s}(x + \alpha (t-s) - y) e^{p(t-s)} dyds \leq N^{-1-\epsilon} e^{(p \vee 0)\Lambda}.
\end{align*}
To bound the penultimate term in \eqref{eq_wNt_momentbdprop1}, we remark that we may make $N$ large enough so that $x \geq S+1 - \alpha t$ implies that $x \geq S - (\epsilon/2)\log N$ for all $t \in [0,\Lambda]$. Hence, for large $N$, for $t \in [0,\Lambda]$ and $x \geq S+1 - \alpha t$ we have
\begin{align*}
&\int_{N^{-1-\epsilon} \wedge t}^t \int_{-\infty}^{S - \epsilon \log N} G_{t-s}(x + \alpha (t-s) - y) e^{p(t-s)} dyds 
\\ &\hspace{1 cm}\leq e^{(p \vee 0)\Lambda}  \int_0^t \bP^B_{x+\alpha(t-s)} (B_{2(t-s)} \leq S - \epsilon \log N)ds
\\ &\hspace{1 cm}\leq  e^{(p \vee 0)\Lambda} \Lambda \bP_{0}^B(B_{2\Lambda} \geq (\epsilon /2) \log N)
\\ &\hspace{1 cm}\leq e^{(p \vee 0)\Lambda} \Lambda N^{-\epsilon^2 \log N/(16\Lambda)},
\end{align*}
which is bounded above by $e^{(p \vee 0)\Lambda} \Lambda N^{-2}$ for large enough $N$. We conclude from these estimates and a trivial bound on the first integral in \eqref{eq_wNt_momentbdprop1} that for that for sufficiently large $N$, for $t \in [0,\Lambda]$ and $x \geq S+1 -\alpha t$, 
\begin{align} \label{eq_wNt_momentbdprop2}
\mathbb{E}_w[w^N_t(x)] &\leq  N^{-1} X e^{(1-\cexp{inter_integ})(S+2)}e^{pt+(1-\cexp{inter_integ})^2 t - (1-\cexp{inter_integ})(x + \alpha t)}  \bP^B_{x+\alpha t - 1 -2(1-\cexp{inter_integ})t}(B_{2t} \leq S) \notag
\\ &\quad + e^{p\Lambda}(1+\Lambda)(2+\|f \|_\infty)  (2N^{-2} + N^{-1-\epsilon}).
\end{align}
Recalling from before \eqref{eq_f_small3} that $p = \alpha - 1 + q_1$, we obtain the simplification
\[e^{pt+(1-\cexp{inter_integ})^2 t - (1-\cexp{inter_integ})( \alpha t)}  = e^{-(\cexp{inter_integ}(2-\alpha-\cexp{inter_integ}) - q_1)t}.\]
Choosing $\epsilon' \in (0,\epsilon)$, the second term on the right-hand side of \eqref{eq_wNt_momentbdprop2} is at most $N^{-1-\epsilon'}$ for sufficiently large $N$. This completes the proof that $\mathbb{E}_w[w^N_t(x)]$ is bounded above by the first term in the minimum from the statement. We now show that it is also bounded above by the second term. Taking the expectation of \eqref{eq_wN_fnormmild} and applying Lemma~\ref{lemma_rightmoment1} and the integrability of the martingale term yields that for $t \in [0,\Lambda]$ and $x \geq S+1 - \alpha t$,
\begin{align*}
\mathbb{E}_w[w^N_t(x)] &\leq e^{\|f'\|_\infty t} I_t(x)
\\ &\leq N^{-1} X e^{(1-\cexp{inter_integ})(S+2)}e^{\|f'\|_\infty t + (1-\cexp{inter_integ})^2 t - (1-\cexp{inter_integ})(x + \alpha t)}  \bP^B_{x+\alpha t - 1 -2(1-\cexp{inter_integ})t}(B_{2t} \leq S)
\\ &\leq N^{-1} X e^{(1-\cexp{inter_integ})(S+2 - x)}e^{(\|f'\|_\infty +1)t}  \bP^B_{x+\alpha t - 1 -2(1-\cexp{inter_integ})t}(B_{2t} \leq S),
\end{align*}
which is the desired inequality. The proof is complete.
 \end{proof}

We now give the proofs of the lemmas.

\begin{proof}[Proof of Lemma~\ref{lemma_rightmoment1}]
Suppose that $R(w) \leq S$. Let $z \in \R$. By \eqref{eq_intervalX_bd} and monotonicity of the heat kernel, it follows that for $t \geq 0$ and $x\in \R$ satisfying $x + \alpha t \geq z + 2$, 
\begin{align*}
\int_z^{z+1} G_t(x+\alpha t - y) w(y) dy &\leq G_t(x+\alpha t - (z+1))  N^{-1} X e^{(1-\cexp{inter_integ})(S-z)}
\\ &\leq N^{-1} X e^{(1-\cexp{inter_integ})(S+1)}  \int_{z}^{z+1} G_t(x+ \alpha t - (y+1)) e^{-(1-\cexp{inter_integ})y}dy.
\end{align*}
Since $R(w) \leq S$, we may conclude that if $x +\alpha t \geq S+1$, under $\bP_w$ we have
\begin{align} \label{eq_Itx_bd1}
I_t w (x)& \leq   \int_{-\infty}^{S} G_t(x + \alpha t - y)w(y)dy  \notag
\\ &\leq N^{-1} X e^{(1-\cexp{inter_integ})(S+1)} \int_{-\infty}^{S} G_t(x+ \alpha t - 1-y) e^{-(1-\cexp{inter_integ})y}dy. 
\end{align}
We remark that
\begin{align*}
G_t(x+ \alpha t - 1-y) e^{-(1-\cexp{inter_integ})y} &= \frac{1}{\sqrt{4\pi t}} e^{-\frac{1}{4t}(y - (x + \alpha t -1))^2 - (1-\cexp{inter_integ})y}
\\ &=  \frac{1}{\sqrt{4\pi t}} e^{-\frac{1}{4t}(y - (x + \alpha t -1) +2 (1-\cexp{inter_integ}) t )^2 + (1-\cexp{inter_integ})^2 t- (1-\cexp{inter_integ})(x + \alpha t - 1)}
\\ &= G_t((x+\alpha t -1 -2(1-\cexp{inter_integ})t) - y) e^{(1-\cexp{inter_integ})^2 t - (1-\cexp{inter_integ})(x + \alpha t - 1)}.
\end{align*}
This implies that 
\begin{align*}
\int_{-\infty}^{S} G_t(x+ \alpha t - 1-y) e^{-(1-\cexp{inter_integ})y}dy &= e^{(1-\cexp{inter_integ})^2 t - (1-\cexp{inter_integ})(x + \alpha t - 1)} \bP^B_{x+\alpha t - 1 -2(1-\cexp{inter_integ})t}(B_{2t} \leq S).
\end{align*}
The above combined with \eqref{eq_Itx_bd1} completes the proof.
\end{proof}

\begin{proof}[Proof of Lemma~\ref{lemma_rightmoment2}]
Suppose $R(w) \leq S$ and $X = X^{N,S}_0 \leq N^{1/3}$ under $\bP_w$, and let $\Lambda>0$. Next, suppose that $t \in [0,\Lambda]$ and $(t,x)$ is such that
\begin{equation} \label{eq_Itx_bd}
I_t w(x) \leq u^* e^{-\|f'\|_\infty \Lambda}/2.
\end{equation}
If the above holds, then it follows from \eqref{eq_wN_fnormmild} that
\begin{align} \label{eq_vmin_prob}
&\bP_w(w^N_t(x) > u^*) 
\\ &\hspace{.5 cm} \leq \bP_w \left( \int_0^t\int_{\R} G_{t-s}(x + \alpha (t-s) - y) e^{\|f'\|_\infty(t-s)} \sqrt{w^N_s(y)(1-w^N_s(y))} \dot{W}(ds,dy) \geq  u^*  N^{1/2}/2 \right).\notag
\end{align}
Using $w^N_s(y) \in [0,1]$, we obtain from the Burkholder-Davis-Gundy inequality, for a universal constant $C>0$, 
\begin{align*}
&\mathbb{E}_w \left[ \left(\int_0^t\int_{\R} G_{t-s}(x + \alpha (t-s) - y) e^{\|f'\|_\infty(t-s)} w^N_s(y)(1-w^N_s(y))) \dot{W}(dsdy)\right)^6 \right]
\\ &\hspace{1 cm}\leq C\left(\int_0^t\int_{\R} G_{t-s}(x + \alpha (t-s) - y)^2 e^{2 \|f'\|_\infty(t-s)} dyds\right)^3
\\ &\hspace{1 cm}\leq  C \left( \int_0^t  \,e^{2 \|f'\|_\infty(t-s)} \int_\R G_{t-s}(y)^2 dy ds \right)^3
\\ &\hspace{1 cm} = C (8\pi)^{-3/2} \left(\int_0^t  e^{2\|f'\|_\infty(t-s)} (t-s)^{-1/2} ds\right)^3
\\&\hspace{1 cm} \leq 8C e^{6 \|f'\|_\infty \Lambda} \Lambda^{3/2},
\end{align*}
where the final constant $C$ depends only on $\sigma$, $f$ and universal quantities, and we have used $t \leq \Lambda$. In particular, by Markov's inequality the probability on the right-hand side of \eqref{eq_vmin_prob} is bounded above by
\begin{equation*}
\frac{2^9 C e^{6 \|f'\|_\infty \Lambda} \Lambda^{3/2} }{u_*^6 N^3} \leq N^{-2},
\end{equation*}
with the inequality holding for sufficiently large $N$ depending on $\Lambda$. To complete the proof, it therefore suffices to establish that, for some $\epsilon >0$, for $\Lambda>0$, \eqref{eq_Itx_bd} holds for $(t,x) \in [N^{-1-\epsilon},\Lambda] \times [S - \epsilon \log N,\infty)$ for sufficiently large $N$.

Let $\epsilon > 0$, with the value to be specified later on. Then for $t \geq 0$ and $x \in \R$,
\begin{align*}
I_t w(x) &=  \int_{-\infty}^{x + \alpha t - \epsilon \log N} G_t(x+\alpha t-y) w(y) dy +  \int_{x + \alpha t - \epsilon \log N}^{\infty} G_t(x+\alpha t-y) w(y) dy
\\ &\leq \bP_0^B(B_{2t} \geq \epsilon \log N) + \frac{1}{\sqrt{4\pi t}} \int_{x + \alpha t - \epsilon \log N}^\infty w(y) dy.
\end{align*}
Using \eqref{eq_intervalX_bd}, we compute 
\begin{align*}
 \int_{x + \alpha t - \epsilon \log N}^\infty w(y) dy = \sum_{z=0}^\infty \int_{x+z + \alpha t - \epsilon \log N}^{x+z+1+\alpha t - \epsilon \log N} w(y)dy &\leq N^{-1} X \sum_{z=0}^\infty e^{(1-\cexp{inter_integ})(S + \epsilon \log N - x - \alpha t - z)}
 \\ & \leq C N^{-1 + \epsilon(1-\cexp{inter_integ})} X e^{(1-\cexp{inter_integ})(S- x - \alpha t)},
\end{align*}
where $C\geq 1$ depends only on $\cexp{inter_integ}$. Returning to the previous inequality, and assuming further that $x \geq S - \epsilon \log N$, we obtain
\begin{align*}
I_t w(x) & \leq \bP_0^B(B_{2t} \geq \epsilon \log N) + C t^{-1/2} N^{-1 + 2\epsilon(1-\cexp{inter_integ})} X.
\end{align*}
Note that $\bP_0^B(B_{2t} \geq \epsilon \log N) \leq e^{-(\epsilon \log N)^2/(4t)} \leq N^{-\epsilon^2 \log N/(4\Lambda)}$ for all $t \in [0,\Lambda]$. Using this, and the assumption that $X \leq N^{1/3}$, we conclude that for $x \geq S - \epsilon \log N$ and $t \in [N^{-1-\epsilon},\Lambda]$ that
\begin{align*}
I_t w(x) \leq N^{-\epsilon^2 \log N/(4\Lambda)} +  CN^{-1 + \frac 1 2 + \frac 1 3 + \frac \epsilon 2 + 2\epsilon(1-\cexp{inter_integ})} =   N^{-\epsilon^2 \log N/(4\Lambda)} +  CN^{-\frac 1 6 + \epsilon (2(1-\cexp{inter_integ}) + \frac 1 2)}.
\end{align*}
If $0<\epsilon < (2(1-\cexp{inter_integ}) + \frac 1 2)^{-1}/6$, the power on the second term is negative. Choose any such $\epsilon$. Then for any $\Lambda>0$, we may take $N$ large enough so that the above is arbitrarily small. Since this bound is uniform over $(t,x) \in [N^{-1-\epsilon},\Lambda] \times [S - \epsilon \log N,\infty)$, this implies \eqref{eq_Itx_bd} for all such $(t,x)$ for sufficiently large $N$, which completes the proof.
\end{proof}

\subsection{Analysis of $X^{N,S}_t$} \label{s_Xanalysis}
This subsection derives various estimates for the process $X^{N,S} = (X^{N,S}_t)_{t \geq 0}$, which is defined in \eqref{def_XNT}. These will be used to control $r^N_t$ in our proof of Proposition~\ref{prop_SN_inc}.
Throughout the subsection we work under $\bP_w$ for $w \in C(\R,[0,1])$ satisfying $R(w) < \infty$ (and hence $X^{N,S}_0 < \infty$) with later assumptions imposed as necessary. The parameter $S$ is in general fixed throughout our analysis, and to the extent possible we will omit dependence on $S$ in notations which are internal to this subsection.

We may express $X^{N,S}_t$ using the weak form \eqref{eq_weakform_movingframe} by taking $\phi(x) = N e^{(1-\cexp{inter_integ})(x-S)}$. 
We can justify this choice of $\phi$ in \eqref{eq_weakform_movingframe} with the same argument we used in Lemma~\ref{lemma_weak_phiV}.
Thus, with probability one, for any $S \in \R$ and all $t \geq 0$, 
\begin{align} \label{eq_XN_weak0}
X^{N,S}_t - X^{N,S}_0 &= \int_0^t \int_\R N e^{(1-\cexp{inter_integ})(x - S)}  \left(((1-\cexp{inter_integ})^2 -\alpha(1-\cexp{inter_integ})) w^N_s(x) + f(w^N_s(x)) \right) dx ds \notag
\\ & \hspace{1 cm} + \sqrt{N} e^{-(1-\cexp{inter_integ})S} \int_0^t \int_\R e^{(1-\cexp{inter_integ})x} w^N_s(x)(1-w^N_s(x)) \dot{W}(ds dx). \notag
\\ &=: \int_0^t \mu^N_s ds + M^N_t,  
\end{align}
where we define 
\begin{equation} \label{def_muN}
\mu_t^N = N \int_{\R} e^{(1-\cexp{inter_integ})(x-S)} \left(((1-\cexp{inter_integ})^2 -\alpha(1-\cexp{inter_integ})) w^N_t(x) + f(w^N_t(x)) \right) dx.
\end{equation}
We also introduce 
\begin{equation*}
(\sigma^N_t)^2 = N \int_\R w^N_t(x)(1-w^N_t(x)) e^{2(1-\cexp{inter_integ})x} e^{-2(1-\cexp{inter_integ})S}  dx,
\end{equation*}
so that 
\[\langle M^N \rangle_t = \int_0^t (\sigma^N_u)^2 du.\]
In particular, there is a Brownian motion $B^N$ such that
\begin{align} \label{eq_XN_weak}
X^{N,S}_t =  X^{N,S}_0 + \int_0^t \mu^N_s ds + \int_0^t \sigma^N_s dB^N_s.
\end{align}
We begin by obtaining a bound on $\mu^N_t$. 

We recall the constant $q_2$, see \eqref{eq_f_small2} and the preceding discussion, as well as the constants $K$.
\begin{lemma} \label{lemma_mu_bd} There exists a constant $C_{K} \geq 1$ such that for any $S \in \R$ and $w \in C(\R,[0,1])$ with $R(w) < \infty$, $\bP_w$-a.s.,
\begin{align*}
\mu^N_t \leq -q_2 X^{N,S}_t + N e^{-(1-\cexp{inter_integ})S} C_K \quad \text{for all $t \leq \hat{\sigma}_{2,N} \wedge \hat{\sigma}_{3,N}$}.
\end{align*}
\end{lemma}
\begin{proof}

We recall the constant $x_0$, which depends only on $m$, appearing in \eqref{eq_wNsmall_tail}. From \eqref{eq_f_small2} and \eqref{eq_wNsmall_tail} we obtain that if $t \leq \hat{\sigma}_{2,N} \wedge \hat{\sigma}_{3,N}$, then
\begin{align*}
&  \int_{\R} e^{(1-\cexp{inter_integ})x}  \left(((1-\cexp{inter_integ})^2 - \alpha(1-\cexp{inter_integ}))w^N_t(x) + f(w^N_t(x)) \right)dx
\\ &\hspace{1 cm} \leq -q_2 \int_{K + x_0}^\infty e^{(1-\cexp{inter_integ})x} w^N_t(x) dx + (1 + \|f\|_\infty) \int_{-\infty}^{K + x_0} e^{(1-\cexp{inter_integ})x}dx
\\ &\hspace{1 cm} \leq -q_2 \int_{\R}  e^{(1-\cexp{inter_integ})x} w^N_t(x) dx + (1 + \|f\|_\infty + q_2) \int_{-\infty}^{K + x_0} e^{(1-\cexp{inter_integ})x} dx
\\ &\hspace{1 cm} \leq  -q_2 \int_{\R}  e^{(1-\cexp{inter_integ})x} w^N_t(x) dx + C_K,
\end{align*}
where we let
\[ C_K := \frac{(1 + \|f\|_\infty + q_2) e^{(1-\cexp{inter_integ})(K + x_0)} }{1-\cexp{inter_integ}}.\]
From \eqref{def_XNT} and \eqref{def_muN}, this completes the proof.
%
\end{proof} 

To simplify notation, for the remainder of this subsection we will write $R_t := R^N_t$. We next introduce a variant of $X^{N,S}_t$ as follows. We write $R^*_t := \sup_{s \in [0,t]} R_s$ for $t \geq 0$ and will use the same notation for the maximal process associated to other time-indexed processes. If $R(w) < \infty$, then $R^*_t < \infty$ $\bP_w$-a.s. for any $t>0$ by Lemma~\ref{lemma_Rsimple}(c). We fix a constant $y_0 > 0$ whose value will be specified later on (see \eqref{def_y0Lambda}) and for $t\geq 0$ define
\begin{equation} \label{def_XNTtilde}
\tilde{X}^{N,S}_t := \int_{\R} e^{(1-\cexp{inter_integ})(x-S)} w^N_t(x) N e^{-(1-\cexp{inter_integ})\left\lfloor \left(\frac{R_t^* - S}{y_0} \right) \vee \,0 \right\rfloor y_0} dx.
\end{equation}
In particular, if $R_t^* < S + y_0$, then $\tilde{X}^{N,S}_s = X^{N,S}_s$ for all $s \in [0,t]$. 

Let $\tau^0_{y_0} = 0$. For $k \in \N$, we define $\tau^k_{y_0}$ by
\[\tau^k_{y_0} = \inf \left\{ t \geq 0 : R_t^* - S \geq k y_0 \right\}.\]
Note that by Lemma~\ref{lemma_Rsimple}(c), for any $\Lambda>0$,
\begin{equation} \label{eq_tauky0_fin}
\lim_{k \to \infty} \bP_w(\tau^k_{y_0} < \Lambda) = 0.
\end{equation}
Next, suppose that $t \leq \tau^k_{y_0}$. Then $w^N_s(x) = 0$ for $s \leq t$ and $x \geq S + ky_0$, and hence
\begin{align*}
(\sigma^N_t)^2 &\leq  N  \int_{-\infty}^{S + ky_0} w^N_t(x) e^{2(1-\cexp{inter_integ})(x-S)} dx
\\ &\leq  e^{(1-\cexp{inter_integ})ky_0}  N \int_{-\infty}^{S + ky_0} w^N_t(x) e^{(1-\cexp{inter_integ})(x-S)} dx
\\ &\leq e^{(1-\cexp{inter_integ})k y_0} X^{N,S}_t.
\end{align*}
Thus, we have
\begin{equation} \label{eq_sigmaN_kbd}
\text{For all $t \leq \tau^{k}_{y_0}$,} \quad (\sigma^N_t)^2 \leq  e^{(1-\cexp{inter_integ})k y_0} X^{N,S}_t.
\end{equation}

The next lemma gives an upper bound for $\tilde{X}^{N,S}_t$.
\begin{lemma} \label{lemma_tildeXNbd} Suppose that $w \in C(\R,[0,1])$ and $S \in \R$ satisfy $R(w) < \infty$ and $S \geq (1-\cexp{inter_integ})^{-1} \log N$. Then $\bP_w$-a.s., \begin{align}\label{eq_lemma_tildeXNbd}
\tilde{X}^{N,S}_{t} &\leq e^{-q_2 t} X^{N,S}_0 +  q_2^{-1} C_K + e^{-q_2 t} \tilde{M}_t \quad \text{for all $t \leq \hat{\sigma}_{2,N} \wedge \hat{\sigma}_{3,N}$,} 
\end{align}
where $\tilde{M}$ is a continuous local martingale satisfying
\[ \langle \tilde{M}\rangle_t \leq C_0 \int_0^t e^{(1-\cexp{inter_integ})y_0}  e^{2 q_2 s} \tilde{X}^{N,S}_s ds.\]
\end{lemma}
\begin{proof}
Let $w$ and $S$ satisfy the stated assumptions. Within this proof, we omit dependence on $N$ and $S$ of $X^{N,S}_t$ and $\tilde{X}^{N,S}_t$, simply writing $X_t$ and $\tilde{X}_t$. To begin, we note that for every $k \in \N_0$, on $\{0<\tau^k_{y_0} < \infty\}$ we have $\tau^k_{y_0} < \tau^{k+1}_{y_0}$ a.s. This can be seen by stopping at $\tau^k_{y_0}$, then applying the strong Markov property and Lemma~\ref{lemma_Rsimple}(a). We also have 
\begin{equation} \label{eq_tildeXN_condrep}
\tilde{X}_t = e^{-k(1-\cexp{inter_integ})y_0} X_t, \quad t \in [\tau^k_{y_0}, \tau^{k+1}_{y_0}).
\end{equation}
Since $X_t$ is continuous, e.g. by \eqref{eq_XN_weak}, this implies that the only jumps of $\tilde{X}^N_t$ occur at the times $\{\tau^k_{y_0}, k \in \N\}$. Furthermore, we remark that for each $k \in \N$,
\begin{equation}
\tilde{X}_{\tau^k_{y_0}} = e^{-(1-\cexp{inter_integ})y_0} \tilde{X}_{\tau^k_{y_0}-}.
\end{equation}
In particular, 
\begin{align*}
\Delta \tilde{X}_{\tau^k_{y_0}} :=  \tilde{X}_{\tau^k_{y_0}} -  \tilde{X}_{\tau^k_{y_0}-} = -(1-e^{-(1-\cexp{inter_integ})y_0}) \tilde{X}_{\tau^k_{y_0}-}.
\end{align*}
Thus, if for $t \geq 0$ we define the non-decreasing jump process
\begin{equation*}
J_t := \sum_{k = 1}^\infty 1(\tau^{k}_{y_0} \leq t) (1-e^{-(1-\cexp{inter_integ})y_0}) \tilde{X}_{\tau^k_{y_0}-}, 
\end{equation*}
then for a continuous process $(H_t)_{t \geq 0}$ with $H_0 = 0$ we have
\begin{equation*}
\tilde{X}_t = \tilde{X}_0 + H_t - J_t.
\end{equation*}
In order to give a representation for $H_t$, we define, for $t \geq 0$,
\[ \xi_t := \sum_{k=0}^\infty 1_{(\tau^k_{y_0}, \tau^{k+1}_{y_0}]}(t) e^{-k(1-\cexp{inter_integ})y_0}.\]
Then $\xi_t$ is bounded and predictable and it is then immediate from \eqref{eq_tildeXN_condrep} that, with probability one,
\[ H_t = \int_0^t \xi_s dX_s, \quad t \geq 0.\]
We note that, in view of \eqref{eq_tildeXN_condrep}, $\xi_t$ has the ``wrong'' value on $\{\tau^k_{y_0}, k \in \N\}$, but changing the diffusion coefficient on this countable set of times a.s. has no effect on $H_t$.

Recall from \eqref{eq_XN_weak} that $dX_s = \mu_s ds + \sigma_s dB_s$ for some Brownian motion $B$, with simplified notation $\mu_s = \mu^N_s$ and $\sigma_s = \sigma^N_s$. Hence, $dH_s = \xi_s \mu_s ds + \xi_s \sigma_s dB_s$. Let $\bar{\sigma}_s =  \xi_s \sigma_s$. To control this diffusion coefficient, we remark from \eqref{eq_sigmaN_kbd}, \eqref{eq_tildeXN_condrep} and the definition of $\xi_s$ that if $s \in (\tau^{k}_{y_0}, \tau^{k+1}_{y_0})$ for some $k \in \N_0$, then
 \begin{align*}
\bar{\sigma}_s^2  = (\xi_s \sigma_s)^2 \leq  e^{-2k(1-\cexp{inter_integ})y_0}  e^{(k+1)(1-\cexp{inter_integ}) y_0} X_s =  e^{(1-\cexp{inter_integ})y_0} \tilde{X}_s.
\end{align*}
This bound holds if $s \in (\tau^{k}_{y_0}, \tau^{k+1}_{y_0})$ for any $k \in \N_0$. As a consequence of \eqref{eq_tauky0_fin} we must have, for any $\Lambda>0$, $\tau^k_{y_0} > \Lambda$ for some (random) $k \in \N$, and thus it follows that with probability one, 
\begin{equation} \label{eq_barsigma_bd}
\bar{\sigma}_s^2 \leq e^{(1-\cexp{inter_integ}) y_0} \tilde{X}_s \text{ for a.e. $s \geq 0$} .
\end{equation}
To conclude, we apply Itô's lemma to $e^{q_2 t} \tilde{X}_t$. We remark that
\[d\tilde{X}_t =  dH_t - dJ_t = \xi_t \mu_t dt + \bar{\sigma}_t dB_t - dJ_t.\]
In particular, our application of Itô's lemma gives 
\begin{align*}
e^{q_2 t} \tilde{X}_t &= \tilde{X}_0 + \int_0^{t} e^{q_2 s} \left( \xi_s \mu_s + q_2 \tilde{X}_s \right) ds + \int_0^{t} e^{q_2 s}  \bar{\sigma}_s dB_s - \int_0^t e^{q_2 s} dJ_s 
\\ &\leq \tilde{X}_0 + \int_0^{t} e^{q_2 s} \left( \xi_s \mu_s + q_2 \tilde{X}_s \right) ds + \int_0^{t} e^{q_2 s}  \bar{\sigma}_s dB_s,
\end{align*}
where the inequality uses the fact that $J$ has positive jumps only. By Lemma~\ref{lemma_mu_bd} and the assumption that $S \geq (1-\cexp{inter_integ})^{-1} \log N$, we have $\mu_s \leq -q_2 X_s + C_K$ for all $s \leq \hat{\sigma}_{2,N} \wedge \hat{\sigma}_{3,N}$. Hence, assuming $s \leq \hat{\sigma}_{2,N} \wedge \hat{\sigma}_{3,N}$, for $k \in \N_0$, if $s \in (\tau^k_{y_0}, \tau^{k+1}_{y_0})$, then by the definition of $\xi_s$ and \eqref{def_XNTtilde} we have
\begin{align*}
\xi_s \mu_s + q_2 \tilde{X}_s&= e^{-k(1-\cexp{inter_integ})y_0 } \mu_s + q_2 e^{-k(1-\cexp{inter_integ})y_0 }  X_s 
\\&\leq e^{-k(1-\cexp{inter_integ})y_0  } (-q_2 X_s + C_K) + q_2 e^{-k(1-\cexp{inter_integ})y_0 }  X_s 
\\ &=  e^{-k(1-\cexp{inter_integ})y_0} C_K
\\ &\leq C_K.
\end{align*}
This bound is uniform in $k \in \N_0$ and in particular it a.s.~holds for a.e.~$s \in [0, \hat{\sigma}_{2,N} \wedge \hat{\sigma}_{3,N}]$. Hence, if $t \leq  \hat{\sigma}_{2,N} \wedge \hat{\sigma}_{3,N}$, then
\begin{align*}
\tilde{X}_{t} &\leq e^{-q_2 t} \tilde{X}_0 +  C_K \int_0^{t } e^{-q_2 (t - s)} ds + e^{-q_2 t} \int_0^{t } e^{q_2 s} \bar{ \sigma}_s dB_s.
\end{align*}
The expression above, \eqref{eq_barsigma_bd}, and the fact that $\tilde{X}_0 \leq X_0$ now imply the result.
\end{proof}

We now prove hitting time bounds for $\tilde{X}^{N,S}_t$. For $\lambda > 0$, define
\begin{equation} \label{def_tautildeX}
	\tau^{\tilde{X}}_\lambda = \inf \{t \geq 0 : \tilde{X}^{N,S}_t > \lambda \}.
\end{equation}
We also recall the constant $C_X$ introduced in \eqref{def_CX} (but whose value is assigned in the next subsection). Here we impose the condition
\begin{equation} \label{eq_CX_lwrbd}
C_X \geq C_K q_2^{-1}.
\end{equation}


\begin{lemma} \label{lemma_Xhit2}
Suppose that $w \in C(\R,[0,1])$ and $S \in \R$ with $R(w) <\infty$, $S \geq (1-\cexp{inter_integ})^{-1} \log N$, and $X_0^{N,S} \leq C_X$ under $\bP_w$. Then for any $\ell \geq 4$ and $\Lambda >0$,
\begin{equation}
\bP_w\left( \tau^{\tilde{X}}_{\ell C_X} \leq \Lambda \wedge \hat{\sigma}_{2,N} \wedge \hat{\sigma}_{3,N} \right) \leq 4(\Lambda+1)  \exp\left(-\frac{\ell C_X }{4q_2^{-1} e^{y_0(1-\cexp{inter_integ}) + 2q_2}}\right).
\end{equation}
\end{lemma}

\begin{proof}
We again suppress dependence of all processes on $N$ and $S$. Furthermore, within this proof we write $\tau_\lambda := \tau^{\tilde{X}}_\lambda$. We first prove a result for general $\lambda$, then specialize to $\lambda = \ell C_X$. Let $\lambda \geq 2(C_Kq_2^{-1} + X_0)$. Our starting point is the upper bound for $\tilde{X}_t$ obtained in Lemma~\ref{lemma_tildeXNbd}. Since $\tilde{X}_0 \leq X_0$ and $X_0 < \lambda$, and recalling that $\tilde{X}_t$ only has negative jumps, as established in the proof of Lemma~\ref{lemma_tildeXNbd}, it follows that $\tilde{X}_{\tau_\lambda} = \lambda$. Hence, if $\tau_\lambda \leq \hat{\sigma}_{2,N} \wedge \hat{\sigma}_{3,N}$, by \eqref{eq_lemma_tildeXNbd} (from Lemma~\ref{lemma_tildeXNbd}) at $t = \tau_{\lambda}$ we have
\begin{align*}
\lambda \leq e^{-q_2 \tau_\lambda} X_0 +  q_2^{-1} C_K + e^{-q_2 \tau_\lambda} \tilde{M}_{\tau_\lambda}.
\end{align*}
By Dubins-Schwarz, we can represent $\tilde{M}$ as a time-changed Brownian motion $\tilde{B}$, i.e. $\tilde{M}_t = \tilde{B}_{\langle \tilde{M} \rangle_t}$ for all $t\geq 0$. We may then rewrite the above as
\begin{align*}
\tilde{B}_{\langle \tilde{M} \rangle_{\tau_\lambda}} \geq e^{q_2 \tau_\lambda}(\lambda - q_2^{-1} C_K) - X_0.
\end{align*}
Furthermore, Lemma~\ref{lemma_tildeXNbd} implies that, if $\tau_\lambda \leq  \hat{\sigma}_{2,N} \wedge \hat{\sigma}_{3,N}$,
\[ \langle \tilde{M} \rangle_{\tau_\lambda} \leq \lambda e^{(1-\cexp{inter_integ})y_0} e^{2q_2 \tau_\lambda }/(2 q_2).\]
Finally, we remark that for $0 \leq t_1 \leq t_2$ and $c>0$, $\tilde{B}_{t_1} \geq c$ implies that $ |\tilde{B}_{t_2}|^* \geq c$, where $|\tilde{B}_{t}|^*= \sup_{s \in [0,t]} |\tilde{B}_{s}|$ for $t \geq 0$. The two inequalities above then imply that
\begin{align*}
&\{\tau_\lambda = t, \hat{\sigma}_{2,N} \wedge \hat{\sigma}_{3,N}\}
\\ &\hspace{1 cm} \subseteq \left\{\tilde{B}_{\langle \tilde{M} \rangle_{t}} \geq e^{q_2 t}(\lambda - q_2^{-1} C_K) - X_0, \tau_\lambda =  t \leq  \hat{\sigma}_{2,N} \wedge \hat{\sigma}_{3,N}\right\}
\\ &\hspace{1 cm}\subseteq \left\{ |\tilde{B}_{\lambda C e^{2 q_2 t }}|^* \geq e^{q_2 t}(\lambda - q_2^{-1} C_K) - X_0\right\},
\end{align*}
where we set $C =  e^{(1-\cexp{inter_integ})y_0}/(2q_2)$. In particular, for $\Lambda>0$,
\begin{align} \label{eq_taulambdapf1}
&\bP_w(\tau_\lambda \leq \Lambda \wedge  \hat{\sigma}_{2,N} \wedge \hat{\sigma}_{3,N} ) \notag
\\ &\hspace{1 cm} \leq \bP_w\left( \exists t \in [0,\Lambda] : |\tilde{B}_{\lambda C e^{2q_2 t }}|^* \geq e^{q_2 t}(\lambda - q_2^{-1} C_K) - X_0\right) \notag
\\ &\hspace{1 cm} \leq \sum_{k=0}^{\lfloor \Lambda \rfloor} \bP_w\left( \exists t \in [k,k+1] : |\tilde{B}_{\lambda C e^{2q_2 t }}|^* \geq e^{q_2 t}(\lambda - q_2^{-1} C_K) - X_0\right).
\end{align}
For a given $k$ in the sum above, we remark that 
\begin{align*}
 &\bP_w\left( \exists t \in [k,k+1] : |\tilde{B}_{\lambda C e^{2q_2 t }}|^* \geq e^{q_2 t}(\lambda - q_2^{-1} C_K) - X_0\right)
 \\ & \hspace{1 cm} \leq  \bP_w\left( |\tilde{B}_{ \lambda C e^{2 q_2 (k+1)}}|^* \geq e^{q_2 k}(\lambda - q_2^{-1} C_K) - X_0\right)
 \\ &\hspace{1 cm} \leq 4 \cdot \bP_w \left(\tilde{B}_{\lambda C e^{2q_2 (k+1)}} \geq  e^{q_2 k}(\lambda - q_2^{-1} C_K) - X_0 \right),
\end{align*}
where the last line uses the reflection principle and our assumption $\lambda \geq 2 (C_K q_2^{-1} + X_0)$. Now, using this assumption again, it follows that for each $k\geq 0$, 
\begin{align*}
\bP_w \left(\tilde{B}_{\lambda C e^{2q_2 (k+1)}} \geq  e^{q_2 k }(\lambda - q_2^{-1} C_K) - X_0 \right) & \leq \bP\left(\tilde{B}_{\lambda C e^{2q_2 (k+1)}} \geq  e^{q_2 k}\lambda/2 \right)
\\ &= \bP_w\left( \tilde{B}_1 \geq \lambda^{1/2} C^{-1/2} e^{-q_2} / 2\right) 
\\ &\leq \exp(-\lambda/(8C e^{2q_2})),
\end{align*}
where we used Brownian scaling and then a standard Gaussian tail estimate. Using the above in the previous inequality, then substituting it into \eqref{eq_taulambdapf1} yields
\begin{align*}
\bP_w(\tau_\lambda \leq \Lambda \wedge \hat{\sigma}_{2,N} \wedge \hat{\sigma}_{3,N} ) \leq 4(\Lambda+1) \exp(-\lambda/(8C e^{2q_2})).
\end{align*}
Substituting the value of $C$ into the inequality above, we obtain that
\begin{align} \label{eq_Xsmall1lambda}
\bP_w \left( \tau_\lambda \leq \Lambda \wedge \hat{\sigma}_{2,N} \wedge \hat{\sigma}_{3,N} \right) \leq 4 (\Lambda +1)   \exp\left(-\frac{\lambda}{4q_2^{-1} e^{y_0(1-\cexp{inter_integ}) + 2q_2}}\right).
\end{align}
This holds under our initial assumption that $\lambda \geq 2(C_{K}q_2^{-1} + X_0)$. Our assumption $X_0 \leq C_X$ and \eqref{eq_CX_lwrbd} implies that this condition is satisfied by $\lambda = \ell C_X$ for $\ell \geq 4$, and substituting this value into \eqref{eq_Xsmall1lambda} completes the proof. \end{proof}

We end this subsection by proving that if $\Lambda \leq \tau^{\tilde{X}}_{4 C_X} \wedge \hat{\sigma}_{2,N} \wedge \hat{\sigma}_{3,N} $, then $\tilde{X}^{N,S}_\Lambda$ is small with high probability. 

\begin{lemma} \label{lemma_Xsmall}
Suppose that $w \in C(\R,[0,1])$ and $S \in \R$ with $R(w) <\infty$, $S \geq (1-\cexp{inter_integ})^{-1} \log N$, and $X_0^{N,S} \leq C_X$ under $\bP_w$. Then if $\epsilon>0$ and $\Lambda >0$ are such that  $C_X (\epsilon - e^{-q_2 \Lambda}) \geq 2q_2^{-1}C_K$,
\begin{equation}
\bP_w \left( \tilde{X}^{N,S}_\Lambda > \epsilon C_X , \Lambda \leq \tau^{\tilde{X}}_{4C_X } \wedge \hat{\sigma}_{2,N} \wedge \hat{\sigma}_{3,N} \right) \leq  2 \exp \left( -\frac{C_X(\epsilon - e^{-q_2 \Lambda})^2}{16e^{(1-\cexp{inter_integ})y_0}q_2^{-1}} \right).
\end{equation}
\end{lemma}
\begin{proof}
We continue to suppress dependence of processes on $N$ and $S$, and write $\tau_{4C_X} = \tau^{\tilde{X}}_{4C_X }$. Let $\epsilon, \Lambda >0$ and suppose that $\epsilon$, $\Lambda$ and $X_0$ satisfy the stated conditions. From \eqref{eq_lemma_tildeXNbd} in Lemma~\ref{lemma_tildeXNbd}, $\bP_w$-a.s. on $\{\Lambda \leq \tau_{4X^N_0} \wedge \hat{\sigma}_{2,N} \wedge \hat{\sigma}_{3,N}\}$ we have
\[ \tilde{X}_\Lambda > \epsilon C_X \Rightarrow e^{-q_2 \Lambda} X_0 +  q_2^{-1} C_K + e^{-q_2 \Lambda} \tilde{M}_\Lambda > \epsilon C_X.\]
Using this and rearranging terms gives
\begin{align*}
&\bP_w(\tilde{X}_\Lambda > \epsilon C_X, \Lambda \leq \tau_{4C_X} \wedge \hat{\sigma}_{2,N} \wedge \hat{\sigma}_{3,N})
\\ &\hspace{1 cm}\leq \bP_w( e^{-q_2 \Lambda} \tilde{M}_\Lambda >  \epsilon C_X - e^{-q_2 \Lambda } X_0  - q_2^{-1} C_K, \Lambda \leq \tau_{4C_X} \wedge \hat{\sigma}_{2,N} \wedge \hat{\sigma}_{3,N})
\\ &\hspace{1 cm}\leq \bP_w\left( e^{-q_2 \Lambda} \tilde{M}_\Lambda^* >  C_X (\epsilon - e^{-q_2 \Lambda}) - q_2^{-1} C_K , \Lambda \leq \tau_{4C_X} \wedge \hat{\sigma}_{2,N} \wedge \hat{\sigma}_{3,N} \right),
\end{align*}
where we set $\tilde{M}^*_\Lambda = \sup_{s\in[0,\Lambda]} \tilde{M}_s$ and use $X_0 \leq C_X$ in the second inequality. Then if $\Lambda \leq \tau_{4C_X} \wedge \hat{\sigma}_{2,N} \wedge \hat{\sigma}_{3,N} $, by Lemma~\ref{lemma_tildeXNbd} we have $ \langle \tilde{M} \rangle_\Lambda \leq 2 C_X e^{(1-\cexp{inter_integ})y_0}e^{2q_2 \Lambda}/q_2$. Again using Dubins-Schwarz to write $\tilde{M}_t = \tilde{B}_{\langle \tilde{M} \rangle_t}$ for $t \geq 0$, and letting $\tilde{B}^*_t = \sup_{s \in [0,t]} \tilde{B}_s$, it follows that 
\[\tilde{M}_\Lambda^* \leq \tilde{B}_{2C_X e^{(1-\cexp{inter_integ})y_0}e^{2q_2 \Lambda}/q_2}^*\]
on $\{\Lambda \leq \tau_{4C_X} \wedge \hat{\sigma}_{2,N} \wedge \hat{\sigma}_{3,N}\}$. Applying this in our previous inequality yields
\begin{align*}
&\bP_w(\tilde{X}_\Lambda > \epsilon C_X,  \Lambda \leq \tau_{4C_X} \wedge \hat{\sigma}_{2,N} \wedge \hat{\sigma}_{3,N})
\\ &\hspace{1 cm} \leq \bP_w\left(  e^{-q_2 \Lambda} \tilde{B}_{2C_X e^{(1-\cexp{inter_integ})y_0}e^{2q_2 \Lambda}/q_2}^* >  C_X (\epsilon - e^{-q_2 \Lambda}) -  q_2^{-1} C_K\right)
\\ &\hspace{1 cm} \leq \bP_w \left( \tilde{B}_{2C_X e^{(1-\cexp{inter_integ})y_0}/q_2}^* >  C_X (\epsilon - e^{-q_2 \Lambda}) - q_2^{-1} C_K\right),
\end{align*}
where the final inequality uses Brownian scaling. Applying scaling again followed by the reflection principle yields the bound
\begin{align*}
&\bP_w(\tilde{X}_\Lambda > \epsilon C_X, \Lambda \leq\tau_{4C_X} \wedge \hat{\sigma}_{2,N} \wedge \hat{\sigma}_{3,N})
\\ &\hspace{1 cm} \leq 2\cdot \bP_w \left((2e^{(1-\cexp{inter_integ})y_0}/q_2)^{1/2} C_X^{1/2} \tilde{B}_{1} > C_X (\epsilon - e^{-q_2 \Lambda}) -  q_2^{-1}C_K  \right).
\end{align*}
Then, grouping terms and using our assumption that $C_X (\epsilon - e^{-q_2 \Lambda}) \geq 2q_2^{-1}C_K$, we obtain from the above that
\begin{align*}
\bP_w (\tilde{X}_\Lambda > \epsilon C_X, \Lambda \leq\tau_{4C_X} \wedge \hat{\sigma}_{2,N} \wedge \hat{\sigma}_{3,N}) &\leq 2 \bP_w\left(\tilde{B}_{1} > \frac 1 2C_X^{1/2} (\epsilon - e^{-q_2 \Lambda})(2e^{(1-\cexp{inter_integ})y_0}/q_2)^{-1/2} \right)
\\ &\leq 2 \exp \left( -\frac{C_X(\epsilon - e^{-q_2\Lambda})^2}{16e^{(1-\cexp{inter_integ})y_0}q_2^{-1}} \right),
\end{align*}
which completes the proof.
\end{proof}

\subsection{Proof of Proposition~\ref{prop_SN_inc}} \label{s_markovpropproof}
We are now ready to combine the results of the three previous subsections in order to prove Proposition~\ref{prop_SN_inc}. Essentially, this means obtaining simultaneous control, over some time interval, of $R^N_t$ and $X^{N,S}_t$ under $\bP_w$ for $w \in \cW_S$. 

We begin by fixing several parameters. Recall  \eqref{eq_Adef} and \eqref{def_Vstar}. We now fix a value of $V$ satisfying
\begin{equation} \label{def_Vchoice}
	V \in (0, V_* \wedge (\alpha/2) \wedge (\cexp{inter_integ}(2-\alpha - \cexp{inter_integ}) - q_1)).
	\end{equation}
We will obtain estimates which hold over time intervals of length $\Lambda$; we now fix values of the constants $C_X$ and $y_0$ which depend on $\Lambda$. In particular, we set
\begin{equation} \label{def_y0Lambda}
	y_0 = y_0(\Lambda) := \frac 1 2 V \Lambda+ 1.\end{equation}
Next, we fix a constant $\epsilon_1>0$. For the time being, we simply specify that $\epsilon_1 < q_2 \wedge 1/4$. We will require further constraints on the value of $\epsilon_1$ that are imposed implicitly in the proof of Lemma~\ref{lemma_RNT_step}, but these are all universal in terms of $V$, $\alpha$ and $\cexp{inter_integ}$, and instead of stating the cumbersome upper bound explicitly we refer to the proof of that result.

For the constant $C_X$ introduced earlier, we now specify the value
\begin{equation}\label{def_CX} C_X = C_X(\Lambda) :=  e^{(1-\cexp{inter_integ})y_0(\Lambda)} e^{4\epsilon_1 \Lambda},\end{equation}
where $y_0 = y_0(\Lambda)$ as above. For sufficiently large $\Lambda$, this is consistent with the assumed lower bound \eqref{eq_CX_lwrbd}. We assume henceforth that $\Lambda \geq \Lambda_0 \geq (4V^{-1}) \vee 4$, where $\Lambda_0$ is large enough so that \eqref{eq_CX_lwrbd} holds.

\begin{lemma} \label{lemma_RNT_step} (a) Let $\delta > 0$. There exists $\Lambda_1(\delta) \geq \Lambda_0$ such that if $\Lambda \geq \Lambda_1(\delta)$, then for sufficiently large $N$, for any $S \geq (1-\cexp{inter_integ})^{-1} \log N$ and $w \in \cW_S$,
\begin{equation*}
\bP_w \left(R^N_t > S + y_0 - V t \text{ for some } t \in[0,\Lambda \wedge \hat{\sigma}_{2,N} \wedge \hat{\sigma}_{3,N}] \right) \leq \delta.
\end{equation*}
(b) There exists $\Lambda_2 \geq \Lambda_0$ such that if $\Lambda \geq \Lambda_2$, for sufficiently large $N$, for any $S \geq (1-\cexp{inter_integ})^{-1} \log N$ and $w \in \cW_S$, 
\begin{equation*}
	\bP_w \left( \sup_{t \in [0,\Lambda \wedge \hat{\sigma}_{2,N} \wedge \hat{\sigma}_{3,N}]} R^N_t >  S + (5+k)\Lambda + 1\right) \leq e^{-\Lambda k^2/32}
\end{equation*}
for all $k \geq 1$.
\end{lemma}
\begin{proof}
The proofs of both parts rely on applications of Proposition~\ref{prop_rightmoment}, which holds for sufficiently large $N$ depending on $\Lambda$. Within this proof, we variously impose size constraints on $\Lambda$; it is always understood that we are choosing $N$ to be sufficiently large so that Proposition~\ref{prop_rightmoment} holds for whichever value of $\Lambda$ we are using. This causes no problem, because at no point do we assume that $\Lambda$ is large in any way which depends upon $N$.

We begin with the proof of part (a). Let $z_0 = S + \frac 1 2 V\Lambda$. Since $V < \alpha$, for sufficiently large $\Lambda$ we have $z_0 - Vt > S + 1 - \alpha t$ for all $t \geq 0$, and hence we can apply Proposition~\ref{prop_rightmoment} at $x = z_0 - Vt$ for all $t \in [0,\Lambda]$. Applying this result and using the first upper bound (i.e. the first term in the minimum), and using the fact that $X^{N,S}_0 \leq C_X$ from \eqref{eq_WS_timezero}, we obtain
\begin{align} \label{eq_ATmomentbd1}
&\mathbb{E}_w [ A_\Lambda(z_0,V) ] \notag
\\ &\hspace{.5 cm}= \int_0^\Lambda \mathbb{E}_w [w^N_t(z_0 - Vt)]dt \notag
\\ &\hspace{.5 cm}\leq  \int_0^\Lambda N^{-1} C_X e^{(1-\cexp{inter_integ})(2-\frac 1 2 V\Lambda)}e^{-(\cexp{inter_integ}(2-\alpha-\cexp{inter_integ}) - q_1 - V)t} \bP^B_{z_0-Vt+\alpha t - 1 -2(1-\cexp{inter_integ})t}(B_{2t} \leq S) dt \notag
\\ &\hspace{1 cm}+ \Lambda N^{-1-\epsilon'},
\end{align}
We will break the integral over $[0,\Lambda]$ into two parts, $[0,b\Lambda]$ and $[b\Lambda,\Lambda]$, for a constant $b \in (0,1)$ which we specify shortly.

By definition of $z_0$, we have
\begin{align*}
z_0-Vt+\alpha t - 1 -2(1-\cexp{inter_integ})t = S + \frac{1}{2} V\Lambda - 1 - (2(1-\cexp{inter_integ}) +V - \alpha)t.
\end{align*}
Note that $(2(1-\cexp{inter_integ}) +V - \alpha) > 0$ due to \eqref{bounds_cexpinterinteg} and $V>0$. Hence, for $t \leq \frac 1 6 V\Lambda (2(1-\cexp{inter_integ})+V - \alpha)^{-1}$, it follows that
\begin{align*}
\bP^B_{z_0-Vt+\alpha t - 1 -2(1-\cexp{inter_integ})t}(B_{2t} \leq S) \leq \bP^B_{S+ \frac 1 3 V\Lambda -1}(B_{2t} \leq S)&\leq \bP^B_{\frac 1 4 V\Lambda}(B_{2t} \leq 0)
\\ &= \bP^B_0(B_{2t} \geq V\Lambda/4)
\\ &\leq e^{-(V\Lambda)^2/(64t)},
\end{align*}
where the second inequality assumes that $\Lambda$ is large enough so that $V\Lambda \geq 12$. Hence, if $b = (\frac V 6 (2(1-\cexp{inter_integ})+V - \alpha)^{-1}) \vee 1$, then 
\begin{align} \label{eq_bTbd1}
& \int_0^{b\Lambda} N^{-1} C_X e^{(1-\cexp{inter_integ})(2-\frac 1 2 V\Lambda)}e^{-(\cexp{inter_integ}(2-\alpha-\cexp{inter_integ}) - q_1 - V)t} \bP^B_{z_0-Vt+\alpha t - 1 -2(1-\cexp{inter_integ})t}(B_{2t} \leq S) dt\notag
 \\ &\hspace{1cm}\leq N^{-1} C_X e^{(1-\cexp{inter_integ})(2-\frac 1 2 V\Lambda )} b\Lambda e^{-(V^2/(64b))\Lambda}, 
 \end{align}
 where we have used, from \eqref{def_Vchoice}, that $\cexp{inter_integ}(2-\alpha-\cexp{inter_integ})-q_1 -V>0$. To handle the rest of the integral, we bound the probability above by $1$ and compute
\begin{align*}
 &\int_{b\Lambda}^\Lambda N^{-1} C_X e^{(1-\cexp{inter_integ})(2-\frac 1 2 V\Lambda )}e^{-(\cexp{inter_integ}(2-\alpha-\cexp{inter_integ}) - q_1 - V)t} \bP^B_{z_0-Vt+\alpha t - 1 -2(1-\cexp{inter_integ})t}(B_{2t} \leq S) dt
\\ &\hspace{1cm}\leq N^{-1} C_X e^{(1-\cexp{inter_integ})(2-\frac 1 2 V\Lambda)} \int_{b\Lambda}^\Lambda e^{-(\cexp{inter_integ}(2-\alpha-\cexp{inter_integ}) - q_1 - V)t} dt
\\ &\hspace{1cm}\leq N^{-1} C_X e^{(1-\cexp{inter_integ})(2-\frac 1 2 V\Lambda)} e^{-(\cexp{inter_integ}(2-\alpha-\cexp{inter_integ}) - q_1 - V)b\Lambda} (\cexp{inter_integ}(2-\alpha-\cexp{inter_integ}) - q_1 - V)^{-1}.
\end{align*}
Recall from \eqref{def_y0Lambda} and \eqref{def_CX} that $C_X = e^{(1-\cexp{inter_integ})y_0 + 4\epsilon_1 \Lambda} = e^{(1-\cexp{inter_integ}) + (1-\cexp{inter_integ}) \frac 1 2 V \Lambda + 4\epsilon_1 \Lambda}$. Hence, by the above and \eqref{eq_bTbd1}, we obtain from \eqref{eq_ATmomentbd1} that
\begin{align} \label{eq:momentbdALambda}
&\mathbb{E}_w [ A_\Lambda(z_0,V) ]  \notag
\\ &\quad \leq N^{-1} e^{3(1-\cexp{inter_integ}) + 4\epsilon_1 \Lambda} \left[  b\Lambda e^{-(V^2/(64b))\Lambda} +e^{-(\cexp{inter_integ}(2-\alpha-\cexp{inter_integ}) - q_1 - V)b\Lambda} (\cexp{inter_integ}(2-\alpha-\cexp{inter_integ}) - q_1 - V)^{-1} \right]  \notag
\\ &\qquad+ \Lambda N^{-1-\epsilon'}  \notag
\\ &\quad \leq N^{-1} e^{-b' \Lambda} + \Lambda N^{-1-\epsilon'},
\end{align}
where $ b' \in  (0, ((\cexp{inter_integ}(2-\alpha-\cexp{inter_integ}) - q_1 - V)b)  \wedge (V^2/(64b)) - 4\epsilon_1)$, and the final inequality holds for sufficiently large $\Lambda$. (We note that this implicitly imposes an upper bound on the value of $\epsilon_1$, as we have previously remarked.) 

	We now conclude using Lemma~\ref{lemma_krylov2}. First, we recall that $z_0 = S + \frac 1 2 V\Lambda$. The application of the lemma therefore requires that 
\begin{equation} \label{eq_gammaN_sigmaN2}
	\hat{\sigma}_{2,N} \wedge \hat{\sigma}_{3,N} \leq \Gamma_N(S - V\Lambda / 2)\, \text{a.s.},
\end{equation}
where $\Gamma_N(S-V\Lambda / 2)$ is defined in \eqref{def_GammaN}.  Since $S \geq (1-\cexp{inter_integ})^{-1} \log N$, this is true for sufficiently large $N$ by \eqref{eq_gammaN_sigmaN}. Applying Lemma~\ref{lemma_krylov2} at $z_0$, with $\beta = 1/2$ and $\zeta \in (0,1]$, followed by the fact that $t \mapsto A_t(z_0,V)$ is non-decreasing, and then Markov's inequality using \eqref{eq:momentbdALambda}, we obtain that for some $x \in [1/2,1]$,
\begin{align}\label{eq_krylovmarkov1}
\bP_w(A_{\Lambda \wedge \hat{\sigma}_{2,N} \wedge \hat{\sigma}_{3,N}}(z_0+x, V) > 0) &\leq \bP_w(A_{\Lambda \wedge \hat{\sigma}_{2,N} \wedge \hat{\sigma}_{3,N}}(z_0,V) > \zeta) + C_2 \big( N \zeta\big)^{1/4} \notag
\\ &\leq \bP_w(A_{\Lambda 	}(z_0,V) > \zeta) + C_2 \big( N \zeta\big)^{1/4} \notag
\\ &\leq \zeta^{-1} \left(N^{-1} e^{-b' \Lambda} + \Lambda N^{-1-\epsilon'}\right) + C_2 \big( N \zeta\big)^{1/4}.  
\end{align}
Finally, we choose $\zeta = N^{-1} \delta^4 / (3 C_2)^4$ and let $\Lambda$ be sufficiently large so that $e^{-b' \Lambda} \leq \delta^5 (3C_2)^{-4}/3$. We thus obtain
\begin{align*}
\bP_w(A_{\Lambda \wedge \hat{\sigma}_{2,N} \wedge \hat{\sigma}_{3,N}}(z_0+x, V) > 0) \leq \delta/3 + (3C_2)^4 \delta^{-4} \Lambda N^{-\epsilon'} + \delta /3 \leq \delta,
\end{align*}
where the last inequality holds for large enough $N$. Hence, by Corollary~\ref{corollary_endpoint} and \eqref{eq_gammaN_sigmaN2}, we have
\[\bP_w \left(R^N_t > z_0 +x- V t \text{ for some } t \in[0,\Lambda \wedge \hat{\sigma}_{2,N} \wedge \hat{\sigma}_{3,N}] \right) \leq \delta.\]
Finally, recall that $z_0 = S + \frac 1 2 V \Lambda$. Since $x \leq 1$, by \eqref{def_y0Lambda} we have $z_0 + x \leq S + y_0$, and hence the above implies the claim from part (a).

The proof of part (b) is similar but with $A_\Lambda(S + (5+k)\Lambda ,0)$ replacing $A_\Lambda(z_0,V)$. We begin by obtaining moment bounds for this quantity. Note that for any $k \geq 1$, Proposition~\ref{prop_rightmoment} can be applied at $(t, S+ (5+k)\Lambda)$ for all $t \in [0,\Lambda]$. Arguing as in part (a), but this time using the second upper bound from Proposition~\ref{prop_rightmoment} instead, for $k \geq 1$ we obtain
\begin{align*}
&\mathbb{E}_w [ A_\Lambda(S	 + (5+k)\Lambda,0) ]
\\ &\hspace{1 cm} \leq \int_0^\Lambda N^{-1} C_X e^{(1-\cexp{inter_integ})(2 - (5+k)\Lambda)}e^{(\|f'\|_\infty +1)t}  \bP^B_{S+ (5+k)\Lambda +\alpha t - 1 -2(1-\cexp{inter_integ})t}(B_{2t} \leq S) dt
\\ &\hspace{1 cm} \leq N^{-1} \int_0^\Lambda e^{(\|f'\|_\infty + 1) t} \bP^B_{k\Lambda}(B_{2t} \leq 0) dt
\\ & \hspace{1 cm} \leq N^{-1} e^{(\|f'\|_\infty+1) \Lambda} \Lambda e^{-\Lambda k^2/4}.
\end{align*}
In the second inequality, we have used the fact that $4\epsilon_1 + (1-\cexp{inter_integ})\frac 1 2 V  \leq 2$ and assumption $\Lambda \geq 4$, which implies (using \eqref{def_y0Lambda} and \eqref{def_CX}) that
\[ C_X e^{(1-\cexp{inter_integ})(2 - (5+k)\Lambda)} \leq e^{2\Lambda + 1 + (1-\cexp{inter_integ})(2 - (5+k)\Lambda)} \leq 1.  \] 
Similarly, $\Lambda \geq 4$ and $t \in [0,\Lambda]$ imply that  
\[ (5+k)\Lambda +\alpha t - 1 -2(1-\cexp{inter_integ})t \geq k \Lambda,\]
which implies the bound used on the Brownian probability. From our previous inequality, it follows that for sufficiently large $\Lambda$, for $k \geq 1$,
\[\mathbb{E}_w [ A_\Lambda(S + (5+k)\Lambda,0)]  \leq N^{-1}e^{-\Lambda k^2/5}. \]
We now argue as in \eqref{eq_krylovmarkov1}, using Lemma~\ref{lemma_krylov2} with $\zeta = N^{-1} e^{-\Lambda k^2/6} $ and Markov's inequality, to obtain that there exists $x \in [1/2,1]$ such that 
\begin{align*}
\bP(A_{\Lambda \wedge \hat{\sigma}_{2,N} \wedge \hat{\sigma}_{3,N}} (S+ (5+k)\Lambda + x,0) > 0) &\leq  e^{-\Lambda k^2(\frac 1 5 - \frac{1}{6})} + C_2 e^{-\Lambda k^2/24}
\\ &\leq e^{-\Lambda k^2/32},
\end{align*}
where the second inequality holds for all $k \geq 1$ for large enough $\Lambda$. The desired result again follows from Corollary~\ref{corollary_endpoint} and \eqref{eq_gammaN_sigmaN2}.
 \end{proof}

\begin{lemma} \label{lemma_XN_step}
(a) Let $\delta >0$. There exists $\Lambda_3(\delta) \geq \Lambda_0$ such that if $\Lambda \geq \Lambda_3(\delta)$, for any $S \geq (1-\cexp{inter_integ})^{-1} \log N$ and $w \in \cW_S$,
\[  \bP_w\left( \left( \bigg\{ \tilde{X}^{N,S}_\Lambda > e^{-\epsilon_1\Lambda} C_X \bigg\} \cup \bigg\{ \sup_{t \in [0,\Lambda]}  \tilde{X}^{N,S}_t > 4C_X \bigg\} \right) \cap  \{\Lambda \leq \hat{\sigma}_{2,N} \wedge \hat{\sigma}_{3,N}\} \right) \leq \delta. \]
(b) There exists $\Lambda_4 \geq \Lambda_0$ such that if $\Lambda \geq \Lambda_4$, for any $S \geq (1-\cexp{inter_integ})^{-1} \log N$, $w \in \cW_S$ and $\ell \geq 4$, 
\[  \bP_w\left( \sup_{t \in [0,\Lambda \wedge \hat{\sigma}_{2,N} \wedge \hat{\sigma}_{3,N}]} \tilde{X}^{N,S}_t > \ell C_X \right) \leq \exp(- \ell \Lambda). \]
\end{lemma}
\begin{proof}
We assume throughout that $S \geq (1-\cexp{inter_integ})^{-1} \log N$ and $w \in \cW_S$. First we will prove part (a). Recalling the definition of $\tau_{4C_X}^{\tilde{X}}$ from \eqref{def_tautildeX}, by a union bound we obtain
\begin{align} \label{eq_Xcontrol1}
& \bP_w \left( \left( \bigg\{ \tilde{X}^{N,S}_\Lambda > e^{-\epsilon_1\Lambda} C_X \bigg\} \cup \bigg\{ \sup_{t \in [0,\Lambda]}  \tilde{X}^{N,S}_t > 4C_X \bigg\} \right) \cap  \{\Lambda \leq \hat{\sigma}_{2,N} \wedge \hat{\sigma}_{3,N}\} \right) 
\\& \hspace{0.5 cm}\leq \bP_w \left(\tau^{\tilde{X}}_{4C_X} < \Lambda, \, \Lambda \leq  \hat{\sigma}_{2,N} \wedge \hat{\sigma}_{3,N}\right) + \bP_w \left(\tilde{X}^{N,S}_\Lambda > e^{-\epsilon_1 \Lambda} C_X, \,\Lambda \leq \tau^{\tilde{X}}_{4C_X} \wedge \hat{\sigma}_{2,N} \wedge \hat{\sigma}_{3,N}\right). \notag
\end{align}
By Lemma~\ref{lemma_Xhit2} (which applies since $X^{N,S}_0 \leq C_X$ by \eqref{eq_WS_timezero}) and the definition of $C_X$ in \eqref{def_CX}, 
\begin{align*}
 \bP_w \left(\tau^{\tilde{X}}_{4C_X} < \Lambda, \, \Lambda \leq \hat{\sigma}_{2,N} \wedge \hat{\sigma}_{3,N}\right) \leq 4(\Lambda+1) \exp \left( -\frac{e^{4 \epsilon_1 \Lambda} }{e^{2q_2} q_2^{-1}} \right).
\end{align*}
We will bound the other probability in \eqref{eq_Xcontrol1} using Lemma~\ref{lemma_Xsmall} with $\epsilon = e^{-\epsilon_1 \Lambda}$. In order for the lemma to apply, we require $C_X(e^{-\epsilon_1 \Lambda} - e^{-q_2 \Lambda}) \geq 2 q_2^{-1} C_{K}$. Since $\epsilon_1 < q_2$ and $C_X = e^{4\epsilon_1 \Lambda +(1-\cexp{inter_integ})y_0}$, this clearly holds for sufficiently large $\Lambda$, and hence 
\begin{align}
\bP_w \left(\tilde{X}^{N,S}_\Lambda > e^{-\epsilon_1 \Lambda} C_X, \Lambda  \leq \tau^{\tilde{X}}_{4C_X} \wedge \hat{\sigma}_{2,N} \wedge \hat{\sigma}_{3,N}\right) \leq  2 \exp \left(- \frac{e^{4\epsilon_1 \Lambda}(e^{-\epsilon_1 \Lambda} - e^{-q_2\Lambda})^2}{16 q_2^{-1}} \right),
\end{align}
where we have used the form of $C_X$ to simplify the expression. Again using $\epsilon_1 < q_2$, for large enough $\Lambda$ the numerator in the exponential above is at least $e^{2\epsilon_1 \Lambda}/2$. Using the above, combined with the previous estimate, we can conclude from \eqref{eq_Xcontrol1} that 
\begin{align*}
 \bP_w \left(\tilde{X}^N_\Lambda > e^{-\epsilon_1\Lambda} C_X, \Lambda < \hat{\sigma}_{2,N} \wedge \hat{\sigma}_{3,N}\right) & \leq 4(\Lambda+1) \exp \left( -\frac{e^{4 \epsilon_1 \Lambda} }{e^{2q_2} q_2^{-1}} \right) + 2 \exp \left(- \frac{e^{2\epsilon_1 \Lambda}}{32 q_2^{-1}} \right). \notag
\end{align*}
Given $\delta >0$, the above is bounded above by $\delta$ for sufficiently large $\Lambda$, proving part (a).

To prove part (b), for $\ell \geq 4$ we apply Lemma~\ref{lemma_Xhit2} to obtain
\begin{align}
\bP_w\left( \sup_{t \in [0,\Lambda \wedge \hat{\sigma}_{2,N} \wedge \hat{\sigma}_{3,N}]} \tilde{X}^{N,S}_t > \ell C_X \right) &\leq \bP_w \left(\tau^{\tilde{X}}_{\ell C_X} \leq \Lambda \wedge \hat{\sigma}_{2,N} \wedge \hat{\sigma}_{3,N}\right)  \notag
\\ &\leq  4(\Lambda +1) \exp \left( - \frac{\ell e^{4 \epsilon_1 \Lambda} }{4 e^{2q_2} q_2^{-1}} \right) \notag
\\ &\leq   \exp \left( - \Lambda \ell \right),
\end{align}
where the last inequality holds for all $\ell \geq 4$ for sufficiently large $\Lambda$. \end{proof}

Let $w \in C(\R,[0,1])$ such that $R(w) < \infty$, and note that, by \eqref{def_rN}, this also implies that $r_N(w) < \infty$. In particular, we may define
\begin{equation} \label{def_Sw}
S(w) := R(w) \vee r_	N(w) \vee \left((1-\cexp{inter_integ})^{-1} \log N\right),
\end{equation}
and we remark that, by \eqref{def_RNrN}, under $\bP_w$ this is equal to the random variable $R^N_0 \vee r^N_0 \vee ((1-\cexp{inter_integ})^{-1} \log N)$. By the definition of $\cW_S$ in \eqref{def_WS}, we have $w \in \cW_{S(w)}$. Take $S \geq S(w)$. We define an event $G(S)$ by
\begin{align*}
G(S) = \left\{\Lambda < \hat{\sigma}_{2,N} \wedge \hat{\sigma}_{3,N} \right\} \cap \bigg\{ \sup_{t \in [0,\Lambda]} R^N_t < S + y_0 \bigg\} \cap \bigg\{ \sup_{t \in [0,\Lambda]} \tilde{X}^{N,S}_t \leq 4C_X \bigg\}
\end{align*}
and then let
\begin{equation} \label{def_calR} \mathcal{R}^N(S) :=  \indc_{G(S)}  R^N_{\Lambda } + \indc_{G(S)^c}  \sup_{t \in [0,  \Lambda \wedge \hat{\sigma}_{2,N} \wedge\hat{\sigma}_{3,N}]} R^N_t
\end{equation}
and
\begin{equation} \label{def_frakR} \mathfrak{R}^N(S) := \indc_{G(S)} r^N_\Lambda + \indc_{G(S)^c}  \sup_{t \in [0,  \Lambda \wedge \hat{\sigma}_{2,N} \wedge\hat{\sigma}_{3,N}]} r^N_t. 
\end{equation}
These of course depend implicitly on a time parameter $\Lambda$, but the omission of this dependence will not cause any confusion. Finally, we define
\begin{equation} \label{def_xi}
\xi^N(S) := \left( \mathcal{R}^N(S) \vee \mathfrak{R}^N(S)\right) - S.
\end{equation}
We now define the constant $\rho>0$ appearing in Proposition~\ref{prop_SN_inc} by
\[ \rho :=  \left(\epsilon_1 (1-\cexp{inter_integ})^{-1}\right) \wedge (V/4).\]

Proposition~\ref{prop_SN_inc} will follow from the next lemma and the Markov property.
\begin{lemma} \label{lemma_incbd_final}
There exist positive constants $M$ and $c$ such that for $\delta >0$ there exists $\bar{\Lambda}(\delta) \geq 1$, such that for sufficiently large $N >0$, the following holds: for any $w \in C(\R,[0,1])$ with $R(w) < \infty$, for any $S \geq S(w)$, and for $\Lambda = \bar{\Lambda}(\delta)$, 
\begin{align*} 
\bP_w \left(\xi^N(S) > -\rho \Lambda , \Lambda \leq \hat{\sigma}_{2,N} \wedge \hat{\sigma}_{3,N} \right) \leq \delta,
\end{align*}
and for all $ k \geq M$,
\begin{align*}
\bP_w (\xi^N (S)\geq k \rho \Lambda) \leq \frac 1 2 e^{-c \Lambda k^2}.
 \end{align*}
\end{lemma}
\begin{proof}
Let $\delta>0$, and let $\Lambda =  \Lambda_1(\delta/2) \vee \Lambda_2 \vee \Lambda_3(\delta/2) \vee \Lambda_4$, as introduced in Lemmas~\ref{lemma_RNT_step} and \ref{lemma_XN_step}. Then for sufficiently large $N$, the bounds from Lemmas~\ref{lemma_RNT_step} and \ref{lemma_XN_step} hold (with parameter $\delta/2$, where applicable) for our choice of $\Lambda$ for any $(w,S) \in (C(\R,[0,1]),\R)$ with $S \geq S(w)$, where we assume $R(w) < \infty$.

We begin with the proof of the first inequality. We define events
\begin{align*}
E_1 &= \bigg\{R_{\Lambda}^N \leq S - V\Lambda /2 +1, \sup_{s \in [0,\Lambda]} R^N_s < S + y_0  \bigg\},
\\ E_2 &= \bigg\{\tilde{X}^{N,S}_{\Lambda} \leq e^{-\epsilon_1 \Lambda} C_X, \sup_{s \in [0,\Lambda]} \tilde{X}^{N,S}_s \leq 4 C_X \bigg\},
\end{align*}
and set $E = \{\Lambda \leq \hat{\sigma}_{2,N} \wedge \hat{\sigma}_{3,N} \} \cap E_1 \cap E_2$. We observe that, by the remark after \eqref{def_XNTtilde}, on $E_1$ we have $X^{N,S}_{\Lambda} = \tilde{X}^{N,S}_{\Lambda}$, and hence the first condition in $E_2$ implies that 
\[ E \subseteq \big\{X^{N,S}_{\Lambda } \leq e^{-\epsilon_1 \Lambda} C_X \big\} .\]
Consequently, on $E$ we have
\begin{align*}
X^{N,S}_\Lambda \leq e^{-\epsilon_1 \Lambda} C_X.
\end{align*}
It follows from \eqref{def_XNT} and \eqref{eq_r_X_link} that on $E$,
\[ r^N_\Lambda \leq S - \epsilon_1 (1-\cexp{inter_integ})^{-1} \Lambda. \]
Recall the definition of $\Lambda_0$ from before the statement of Lemma~\ref{lemma_RNT_step}. Since $\Lambda \geq \Lambda_0$, we also have $- V\Lambda /2 + 1 \leq -V\Lambda/4$, and so the definition of $E_1$ implies that on $E$ we have
\[ R^N_{\Lambda} \leq S - V\Lambda/4. \]
Since $E \subseteq G(S)$, it follows that $\mathcal{R}^N(S) = R^N_\Lambda$ and $\mathfrak{R}^N(S) = r^N_\Lambda$ on $E$. This observation combined with the last two inequalities which hold on $E$ imply that, on $E$,
\begin{align*}
\xi^N(S) = \left(R^N_\Lambda \vee r^N_\Lambda\right) - S  \leq - \left((\epsilon_1(1-\cexp{inter_integ})^{-1} \right) \wedge (V/4)) \Lambda = - \rho \Lambda.
\end{align*}
To prove part (a), it therefore suffices to prove that 
\begin{align} \label{eq_Ec_probdelta}
\bP_w(E^c \cap \{ \Lambda \leq \hat{\sigma}_{2,N} \wedge \hat{\sigma}_{3,N}\}) \leq \delta.
\end{align}
To prove this, we observe that
\begin{align*}
\bP_w(E^c \cap \{\Lambda \leq \hat{\sigma}_{2,N} \wedge \hat{\sigma}_{3,N}\}) &\leq \bP_w(E_1^c \cap \{\Lambda \leq \hat{\sigma}_{2,N} \wedge \hat{\sigma}_{3,N}\}) +  \bP_w(E_2^c \cap \{\Lambda \leq \hat{\sigma}_{2,N} \wedge \hat{\sigma}_{3,N}\}).
\end{align*}
We remark that $E_2^c \cap \{\Lambda \leq \hat{\sigma}_{2,N} \wedge \hat{\sigma}_{3,N}\}$ is precisely the event from Lemma~\ref{lemma_XN_step}(a). Since $\Lambda \geq \Lambda_3(\delta/2)$, it follows that the second probability in the above is at most $\delta/2$ for sufficiently large $N$. To handle the first probability, we remark that since $y_0 = 1 + \frac 1 2 V \Lambda$ by \eqref{def_y0Lambda}, 
\begin{align*}
E_1^c \subseteq \left\{R^N_t >  S + y_0 - Vt \text{ for some } t \in [0,\Lambda]  \right\} \cup \left\{ \tau^R_{S + y_0/2} = 0 \right\}.
\end{align*}
Since $R(w) \leq S < S + y_0/2$, the second event has probability zero by Lemma~\ref{lemma_Rsimple}(a). Since $\Lambda \geq \Lambda_1(\delta/2)$, the intersection of the first event with $\{\Lambda \leq \hat{\sigma}_{2,N} \wedge \hat{\sigma}_{3,N}\}$ has probability at most $\delta/2$ for sufficiently large $N$ by Lemma~\ref{lemma_RNT_step}(a). Hence \eqref{eq_Ec_probdelta} holds for sufficiently large $N$, and the first bound is proved.

To prove the second inequality, for $k\geq 1$ we define 
\[ F_k = \{\xi^N(S) \geq k\rho \Lambda\} = \{ \mathcal{R}^N(S) \vee \mathfrak{R}^N(S) > S + k \rho \Lambda \}.\]
Note that the definitions of $\mathcal{R}^N(S)$ and $\mathfrak{R}^N(S)$ imply that 
\[ \mathcal{R}^N(S) \vee \mathfrak{R}^N(S) \leq (R^N_\Lambda)^* \vee (r^N_{\Lambda})^*,\]
where $(R^N_\Lambda)^* = \sup_{t\in [0,\Lambda \wedge \hat{\sigma}_{2,N} \wedge \hat{\sigma}_{3,N}]} R^N_t$ and $(r^N_{\Lambda})^*$ is defined similarly, and hence $F_k \subseteq F_k'$, where
\[F_k' = \{ (R^N_\Lambda)^* \vee (r^N_{\Lambda})^* \geq S + k \rho \Lambda \}.\]
Hereafter we can and do consider $F_k'$ instead of $F_k$. Next, we define the events
\begin{align*}
B_{1,k} = \{(R^N_\Lambda)^* \geq S + k \rho \Lambda /2 \}, \quad B_{2,k} = \{(\tilde{X}^{N,S}_{\Lambda})^* > k^3 C_X \},
\end{align*}
where $(\tilde{X}^{N,S}_{\Lambda})^* = \sup_{t\in [0,\Lambda \wedge \hat{\sigma}_{2,N} \wedge \hat{\sigma}_{3,N}]} \tilde{X}^{N,S}_t$. We now show that for sufficiently large $k$, $F_k' \subseteq B_{1,k} \cup B_{2,k}$, and so it suffices to bound the probability of the latter event. Clearly on $B_{1,k}^c$ we have
\[ (R^N_\Lambda)^* < S + k\rho \Lambda.\]
On the other hand, on $B_{1,k}^c \cap B_{2,k}^c$,  
\begin{align*}
\sup_{t \in [0,\Lambda \wedge \hat{\sigma}_{2,N} \wedge \hat{\sigma}_{3,N}]} X^{N,S}_t \leq  e^{(1-\cexp{inter_integ})\lfloor \frac{k \rho \Lambda / 2}{y_0}\rfloor  y_0}  (\tilde{X}^{N,S}_\Lambda)^* &\leq e^{(1-\cexp{inter_integ}) k \rho \Lambda / 2} k^3 C_X
\\& =e^{(1-\cexp{inter_integ}) k \rho \Lambda / 2 + 3 \log k} C_X.
\end{align*}
By \eqref{eq_r_X_link}, this implies that
\begin{align*}
(r^N_\Lambda)^* \leq S + k \rho \Lambda / 2 + 3(1-\cexp{inter_integ})^{-1} \log k < S + k \rho \Lambda,
\end{align*}
where the latter inequality holds for $k \geq M_0 \in \N$ for some $M_0$ which is independent of $N$ and $\Lambda$ (since $\Lambda \geq 1$). Combining the above bounds on $(R^N_\Lambda)^*$ and $(r^N_\Lambda)^*$, we have shown that for $k \geq M_0$, on $ B_{1,k}^c \cap B_{2,k}^c$, $(R^N_\Lambda)^* \vee (r^N_\Lambda)^* < S + k\rho \Lambda$. In particular, $  B_{1,k}^c \cap B_{2,k}^c \subseteq (F_k')^c$, and hence, for $k \geq M_0$,
\begin{equation}\label{eq_Fk_probbd}
\bP_w(F_k') \leq  \bP_w(B_{1,k}) + \bP_w(B_{2,k} ).
\end{equation} 
To complete the proof, we need only bound the two probabilities on the right-hand side. 

Let $M_1 = 24/\rho$ and suppose that $k \geq M_1$. If $m_k =  k \rho / 2 - 6 $, then
\[B_{1,k} = \{ (R^N_\Lambda)^*  \geq S + (6+m_k)\Lambda \} \subseteq \{ (R^N_\Lambda)^*  \geq S + (5+m_k)\Lambda +1 \}\]
Moreover, the condition $k \geq M_1$ implies that $m_k \geq (k \rho)/4$. In particular, using Lemma~\ref{lemma_RNT_step}(b), we obtain
\begin{align}\label{eq_B1k_probbd}
\bP_w( B_{1,k}) \leq  \bP_w((R^N_\Lambda)^* \geq S + (5+m_k)\Lambda + 1)  \leq e^{-\Lambda m_k^2/32} \leq e^{-c \Lambda k^2},  
\end{align}
where $c = \rho^2 /(16\cdot 32)$. On the other hand, by Lemma~\ref{lemma_XN_step}(b), for all $k \geq 4^{1/3}$ we have
\begin{align*}
\bP_w(B_{2,k} ) \leq e^{-\Lambda k^3}.
\end{align*}
Combining \eqref{eq_Fk_probbd}, \eqref{eq_B1k_probbd} and the above, and using the fact that $F_k \subset F_k'$, we may conclude that for $k \geq M := M_0 \vee M_1 \vee 4^{1/3}$, 
\[ \bP_w(F_k )  \leq  e^{-c \Lambda k^2} +e^{-\Lambda k^3}. \] 
Since $\Lambda \geq 1$, we may modify the value of the constants $c$ and $M$, so that the above is bounded above by $\frac 1 2 e^{-c\Lambda k^2}$ for all $k \geq M$. This completes the proof.
\end{proof}

We are nearly ready to give the proof of Proposition~\ref{prop_SN_inc}. First, we briefly discuss the Markov property. As a consequence of uniqueness-in-law for solutions to \eqref{eq_spde_movingframe} (which we have assumed), solutions are Markov. 
In order to facilitate arguments using the Markov property, for the time being we will assume that the solution is defined on a canonical space equipped with the canonical shift operators $(\theta_t : t \geq 0)$, i.e. for $s, t \geq 0$, $(\theta_t(\omega))(s) = \omega(t+s)$. We keep technical details to a minimum here, but in view of the compact interface property, but we observe that a suitable canonical state space is given by $\Omega = C(\R_+,\mathcal{C})$,
where
\[ \mathcal{C} := \left\{ u \in C(\R,[0,1]) : \lim_{ x \to - \infty} u(x) = 1 \text{ and } \lim_{x \to \infty} u(x) = 0\right\},\]
which is separable with respect to the uniform topology. Then the Markov property for $(w^N_t)_{t \geq 0}$ can be expressed as
	\begin{equation*}
		\mathbb{E}_w[Z \circ \theta_t \, | \, \cF_t](\omega) = \mathbb{E}_{w^N_t(\omega)}[Z] \, \text{ a.s.}
	\end{equation*}
	for any bounded, measurable random variable $Z : C(\R_+,\mathcal{C}) \to \R$, where $w^N_t(\omega)= \omega(t)$ is the coordinate map,
	and the $N$ is included simply to remind us which solution we are working with.



\begin{proof}[Proof of Proposition~\ref{prop_SN_inc}] 
Let $\delta >0$ and $\Lambda = \bar{\Lambda}(\delta)$ as in Lemma~\ref{lemma_incbd_final}, and let $N$ be large enough so that the conclusions of that result hold. We recall the event $F^N_n$ from \eqref{def_FNn}, the definitions \eqref{def_calRn}-\eqref{def_xin} of the quantities $\mathcal{R}^N_n$, $\mathfrak{R}^N_n$, $S^N_n$ and $\xi^N_n$ (we omit dependence on $\Lambda$) defined under the  probability measure $\bP$, as well as the definitions \eqref{def_Sw}-\eqref{def_xi} of $S(w)$, $\mathcal{R}^N$, $\mathfrak{R}^N$, and $\xi^N$, defined under $\bP_w$ for any $w \in C(\R,[0,1])$ with $R(w) < \infty$. Next, we claim that for $x \in \R$,
\begin{equation} \label{eq_xi_markov}
	\bP(\xi^N_n > x \, | \, \mathcal{G}^N_{n-1}) = \bP_{w^N_t}(\xi^N(S^N_{n-1}) > x) \quad \text{$\bP$-a.s. on $F^N_{n-1}$} 
\end{equation}
and 
\begin{equation} \label{eq_xi_markov2}
	\bP(\{ \xi^N_n > x\} \cap F^N_n   \, | \, \mathcal{G}^N_{n-1}) = \bP_{w^N_t}(\xi^N(S^N_{n-1}) > x, \Lambda < \hat{\sigma}_{2,N}\wedge \hat{\sigma}_{3,N}) \quad \text{$\bP$-a.s. on $F^N_{n-1}$}.
\end{equation}
Both claims essentially follow from the Markov property and the fact that 
\begin{equation} \label{eq_RN_shift}
	\left(\mathcal{R}^N_n  \vee \mathfrak{R}^N_n\right) = \left(\mathcal{R}^N \vee \mathfrak{R}^N \right) \circ \theta_{\Lambda(n-1)} \,\, \text{ on $F^N_{n-1}$.}
\end{equation}
and hence $\xi^N_n \circ \circ \theta_{\Lambda(n-1)} = \xi^N(S^N_{n-1})$ on $F^N_{n-1}$. 
The application of the Markov property is complicated by the fact that $\xi^N(S^N_{n-1}$) depends on $\mathcal{G}_{n-1}$ through $S^N_{n-1}$; 
however, it does so only parametrically, in that $S^N_{n-1}$ ``selects'' which random variable we compute, but given this selection the random variable depends only on the future. 
That the Markov property is applicable in this situation can be justified via a monotone class argument which we omit. Thus, \eqref{eq_xi_markov} and \eqref{eq_xi_markov2} are proved.

We remark from the definition of $S^N_{n-1}$ in \eqref{def_Sn} and from \eqref{def_Sw} that 
\begin{equation*}
	S^N_{n-1} \geq S(w^N_{\Lambda(n-1)}) \,\, \text{ a.s. on $F^N_{n-1}$.}
\end{equation*}
In particular, Lemma~\ref{lemma_incbd_final} can a.s.~be applied to the right-hand sides of \eqref{eq_xi_markov} and \eqref{eq_xi_markov2}, respectively with $x=-\rho \Lambda$ and $x = k\rho\Lambda$ for $k \geq M$. This completes the proof. \end{proof}

	\section{Hölder and tail estimates: proof of Proposition~\ref{prop:control_sigma}}  \label{sec:holder}
		
	In this section, we complete the proof of Proposition~\ref{prop:control_sigma}.
	We begin by stating four intermediate results that we use to prove Proposition~\ref{prop:control_sigma}; then in Sections~\ref{subsec:holder}-\ref{subsec:techG-G}, we prove these intermediate results.
	We work as usual with $u^N_0$ satsifying Assumption~\ref{assumpt:v0}, with $K>0$, $\beta_2 \in (0,\beta_1)$ and $\varepsilon_2\in (0,\varepsilon_1)$ fixed to their previous values.
	
	We begin with an estimate on the H\"older continuity of $v^N_t$, which will be proved in Section~\ref{subsec:holder}.
	For $x,y\in \R$ and $s,t\ge 0$, we let
	\begin{equation} \label{eq:2normdefn}
	\| (x,t) - (y,s) \|_2 := \left( |x-y|^2 + |t-s| \right)^{1/2}.
	\end{equation}
	We further recall from \eqref{eq:holderconstdef} that $\cexp{holder}$ satisfies
	\begin{equation*} 
		\cexp{holder} \in (0, \tfrac 14 \wedge \tfrac{\cinit{big}}{2}).
	\end{equation*}
	\begin{proposition} \label{prop:tildeudiffunif}
	For each $k \in \N$ there exists $A_k<\infty$ such that the following holds. 
	For any $x\in \R$ and $t\geq 0$, there exists a random variable $U^N_{t,x,k}$ such that
	\begin{equation*} 
	|v^N_{t_1}(x_1)-v^N_{t_2}(x_2)|
	\leq U^N_{t,x,k} \|(x_1,Nt_1)-(x_2,N t_2)\|_2^{\cexp{holder}}
	\quad \forall x_1,x_2 \in [x,x+2], t_1,t_2\in [t,t+2N^{-1}]
	\end{equation*}
	and
	 $\E{(U^N_{t,x,k})^{2k}}\le A_k $.
	\end{proposition}

	Recall the definition of $\sigmaN{tail_right}$ from \eqref{eq:sigmatailrightdefn}.
	We also define a stopping time $\newsigmaN{tail_left}$ that controls the tail behaviour of $v^N_t(x)$ as $x \to -\infty$ as
	\begin{equation} \label{eq:sigmalefttaildefn}
		\sigmaN{tail_left}:=\inf\{t\geq 0: \exists y\in \R \text{ s.t. }1-v^N_t(y)\geq \Cst{tail_left}(e^{\cexp{tail_left} y}+N^{-\cexp{tail_left}})\}.
	\end{equation}		

	The following result, which we remark implies Proposition~\ref{prop:control_sigma} for $i = \ref{sigmaN:tail_right}$, is proved in Section~\ref{subsec:utail}.
	\begin{proposition} \label{prop:sigmaNrighttail}
		For $T>0$, for $N$ sufficiently large,
		\begin{align*}
			\P{\sigmaN{tail_right}\leq  \sigmaN{interface} \wedge T} < N^{-2}.
		\end{align*}
	\end{proposition}
%
	We shall prove the following result in Section~\ref{subsec:1-utail}.
	
	\begin{proposition} \label{prop:sigmatailleft}
		For $T>0$, as $N\to \infty$,
		\begin{align*}
			\P{\sigmaN{tail_left}\leq \sigmaN{dist}\wedge \sigmaN{eta}  \wedge T} \leq 3 N^{-2}.
		\end{align*}
	\end{proposition}
	
	We will also require a variant of our Hölder continuity estimate with spatial decay in the constants. The following result is proved in Section~\ref{subsec:expholder}.
	\begin{proposition} \label{prop:holderexp}
	There exists a constant $A_{(0)}$ such that 
	\begin{align*}
	&\mathbb{P} \bigg( \exists t_1,t_2 \in [0, \sigmaN{interface} \wedge \sigmaN{tail_right} \wedge N], \, x_1, x_2 \geq 0 \text{ with } |t_1 - t_2| \leq N^{-1} \text{ and }|x_1-x_2| \leq 1 :
	\\ &\hspace{1 cm} |v^N_{t_1}(x_1) - v^N_{t_2}(x_2)| > A_{(0)} N^{\cexp{tail_right_N}} e^{-((1-\cexp{tail_right_exp})/2)(x_1 \wedge x_2)} (|x_1 - x_2|^{\cexp{holder}} +|t_1-t_2|^{\cexp{holder}/2}) \bigg) < N^{-2}.
	\end{align*}
	\end{proposition}
	
%
	
	Note that our assumptions on $u^N_0$ give us the following result.
	\begin{lemma} \label{lem:stoptimepos}
	For $N$ sufficiently large, almost surely, 
	\[
	\sigmaN{tail_right}>0 \quad \text{ and } \quad \sigmaN{tail_left}>0  
	\]
	\end{lemma}
	\begin{proof}
	From Assumption~\ref{assumpt:v0}.\ref{v0:tail}, \eqref{eq:cleftfprime}, \eqref{bounds_cinitsmall} and \eqref{eq:lefttailinitconst}, we have
	\[u_0^N(x) \leq N^{\cexp{tail_right_N}} e^{-(1-\cexp{tail_right_exp})x}, \quad 1-u^N_0(x) \leq \Cst{init}(e^{\cexp{tail_left}x} + N^{-\cexp{tail_left}})  \]
	for all $x \geq 0$, and so the desired conditions are satisfied at time $t=0$. 
	By Lemma~\ref{lem:compactinterface}, solutions are continuous and a.s. constant outside of a (random) bounded interval for $t \in [0,1]$. 
	This implies that the inequalities above hold a.s. with $u^N_t(x)$ replacing $u^N_0(x)$ for sufficiently small $t$, which completes the proof.
	\end{proof}

Note that the above implies Lemma~\ref{lem:tauNpos} for $i = \ref{sigmaN:tail_right}$. We now complete the proof of Lemma~\ref{lem:tauNpos}.

\begin{proof}[Proof of Lemma~\ref{lem:tauNpos} for $i = \ref{sigmaN:supnorm}, \ref{sigmaN:vdelta}, \ref{sigmaN:tail}$]
	First we consider $\sigmaN{vdelta}$ and $\sigmaN{tail}$. We begin by showing proving that for sufficiently large $N$,
	\begin{equation} \label{eq:initbdsigma56}
	\int_\R u_0^N(x) e^{\alpha x} dx < \kappa_N, \quad \|D^{\delta_N,N}(t,\cdot) \|_{2,\alpha} < r_N, \,\,\, \text{ and } \,\,\, \|v^{\delta_N,N}_0 - v^N_0\|_\infty < r_N.
	\end{equation}
	We begin with the first condition. By Assumption~\ref{assumpt:v0}.\ref{v0:compact_support} and~\ref{assumpt:v0}.\ref{v0:tail},
	assuming first that $\alpha - 1 + \cinit{small} > 0$, we obtain
	\begin{align} \label{eq:sigma5Npos}
		\int_{\R}u^N_0(x)e^{\alpha x}dx
		&\le \int_{-\infty}^0 e^{\alpha x}dx + N^{\cinit{small}} \int_{0}^{(1+\cinit{small})\log N} e^{\alpha - (1-\cinit{small})x}dx \notag\\
		&\le \alpha^{-1} +  (\alpha - 1 + \cinit{small}))^{-1} N^{\cinit{small} + (1+\cinit{small}) (\alpha - 1+ \cinit{small})) }\notag \\
		&< N^{1-\cexp{KN}}.
	\end{align}
	The last line holds 
	where the last line holds for sufficiently large $N$ due to \eqref{bounds_cinitsmall} and \eqref{eq:cKNcsum}. Since $\kappa_N = N^{1-\cexp{KN}}$ by \eqref{def:kappa_delta_r_theta}, 
	it follows that $\sigmaN{tail}>0$. 
	If $\alpha - 1 + \cinit{small} \leq 0$, the second integral in the first line is at most $(1+\cinit{small}) \log N$, and again the result follows for sufficiently large $N$ by \eqref{bounds_cinitsmall}.

	To avoid repeating arguments, we will omit the proof of the second and third inequalities in \eqref{eq:initbdsigma56}. 
	Indeed, the proof is essentially contained in the proof of Proposition~\ref{prop:control_sigma} in the case $i=\ref{sigmaN:vdelta}$, which we give shortly.
	That proof uses Proposition~\ref{prop:tildeudiffunif} to obtain adequate Hölder continuity and uses the control provided by the stopping times $\sigmaN{tail_right}$ and $\sigmaN{tail_left}$ for integrability.
	Since $\cinit{big} > \cexp{holder} \vee \cexp{tail_left}$ and $\cinit{small} < \cexp{tail_right_exp} \wedge \cexp{tail_right_N}$ from \eqref{eq:holderconstdef}, \eqref{eq:cleftfprime} and \eqref{bounds_cexptailrightexp},
	it follows that Assumption~\ref{assumpt:v0}.\ref{v0:holder} and~\ref{assumpt:v0}.\ref{v0:tail} offer stronger Hölder continuity and tail decay than
	those provided by Proposition~\ref{prop:tildeudiffunif}, $\sigmaN{tail_right}$ and $\sigmaN{tail_left}$, and hence result can be proved in the same way with essentially no modifications. We omit the details. 
	
	Given \eqref{eq:initbdsigma56}, it will follow that $\sigmaN{vdelta} >0$ and $\sigmaN{tail}>0$ if the maps
	\[ t \mapsto \int_\R u_0^N(x)e^{\alpha x} dx, \quad t\mapsto \|D^{\delta_N,N}(t,\cdot)\|_{2,\alpha}, \,\,\, \text{ and } \,\,\, t \mapsto \|v^{\delta_N,N}_t - v^N_t \|_\infty\]
	are continuous. For the first two processes, this can be seen directly from the compact interface property (c.f. Lemma~\ref{lem:compactinterface}) and dominated convergence.
	For the third process, continuity follows as a result of (space-time) continuity of $v^{\delta_N,N}_t(x) - v^N_t(x)$ and the compact interface property. This completes the proof for $i = \ref{sigmaN:vdelta}, \ref{sigmaN:tail}$.
	
	Now consider $i = \ref{sigmaN:supnorm}$. We will again defer to a coming argument to avoid repeating proofs here. 
	In the proof of Proposition~\ref{prop:control_sigma} for $i = \ref{sigmaN:supnorm}$, we obtain control of $\sigmaN{supnorm}$ via the conditions defining $\sigmaN{dist}$ and $\sigmaN{eta}$,
	and Hölder regularity from Propositions~\ref{prop:tildeudiffunif} and~\ref{prop:holderexp}. To prove that 
	\[ \| u^N_0 \|_{\infty, \ggamma{weight}} < \Epsilon{proof},\]
	these ingredients above may be substituted with Assumptions~\ref{assumpt:v0}.\ref{v0:dist} and~\ref{assumpt:v0}.\ref{v0:holder}, 
	using the facts that $\cexp{thetaN}, \cexp{holder} < \cinit{big}$ and $\cexp{tail_right_exp}, \cexp{tail_right_N} > \cinit{small}$ from \eqref{eq:holderconstdef}, \eqref{bounds_cexptailrightexp}, \eqref{bounds_cexpthetaN}, \eqref{bounds_cexptailrightN}, and the proof is then identical. Continuity of $t \mapsto \| v^N_t \|_{\infty, \ggamma{weight}}$ then follows from the compact interface property and continuity of solutions,
	and it follows that $\sigmaN{supnorm} > 0$.
	\end{proof}

	We now use Propositions~\ref{prop:tildeudiffunif}-\ref{prop:holderexp} to complete the proof of Proposition~\ref{prop:control_sigma}. 
	Recall that in Section~\ref{subsec:stability} we proved Proposition~\ref{prop:control_sigma} in the case $i=\ref{sigmaN:dist}$; 
	it remains to consider the cases $i=\ref{sigmaN:supnorm}$, \ref{sigmaN:vdelta} and \ref{sigmaN:tail}.
		

	We begin with the case $i= \ref{sigmaN:supnorm}$.
	First, we state and prove a lemma which we will re-use during the proof of Proposition~\ref{prop:sigmatailleft}.
	The lemma's proof only uses Proposition~\ref{prop:tildeudiffunif}, the proof of which is independent of Proposition~\ref{prop:sigmatailleft}, and hence there is no issue using the lemma later on in that setting.
	\begin{lemma} \label{lem:scontrol}
	For any $\epsilon >0$ and $C>0$, for $N$ sufficiently large,
	\begin{equation*}
	\P{\exists x \in [-C,C], t \in [0, \sigmaN{dist} \wedge \sigmaN{eta} \wedge N] : |s(v^N_t)(x)| > \epsilon} < N^{-2}.
	\end{equation*}
	\end{lemma}
	\begin{proof}
		Recall from~\eqref{def_sigma_dist}-\eqref{def_sigma_sup} and Lemma~\ref{lemma:eta} that
		$\eta(v^N_t)$ and $s(v^N_t)$ are well defined for $t\le \sigmaN{dist}\wedge \sigmaN{eta}$.		
		We recall the definition of $\sigmaN{dist}$ from~\eqref{def_sigma_dist}, and recall from~\eqref{def:kappa_delta_r_theta} that $\vartheta_N=N^{-\cexp{thetaN}}\wedge \Beta{proof}$, 
		Let $\epsilon > 0$ and $C >0$ and take $a\in (0,\cexp{thetaN} \cexp{holder})$. We fix $k\in \N$ with $ak> 2$; recall the definition of the random variables $U^N_{t,x,k}$ in Proposition~\ref{prop:tildeudiffunif}.
		By a union bound and Markov's inequality, and then by Proposition~\ref{prop:tildeudiffunif},
		we have
		\begin{multline} \label{eq:UNK'bd}
		\P{\exists x\in \Z\cap [-C-1,+C],\, t\in N^{-1} \Z\cap [0,N]:U^N_{t,x,k}\geq N^a} \\
		\begin{aligned}
		&\quad \leq \sum_{x\in \Z \cap [-C-1,C]} \: \sum_{t \in N^{-1} \Z \cap [0,N]} N^{-2ak}\E{(U^N_{t,x,k})^{2k}} \\
		&\quad \leq (2C+2) (N^2+1)N^{-2ak}A_k \\
		&\quad < N^{-2}
		\end{aligned}
		\end{multline}
		for $N$ sufficiently large, where the last inequality follows since $ak> 2$.
		Now suppose $N$ is sufficiently large that 
		\begin{align} \label{choice_N}
			N^{-a}+N^{-2a/\cexp{holder}}\|\partial_{x} m\|_\infty<\tfrac \epsilon 2 && \text{ and } && \tfrac 1 {4} e^{-\alpha C} \epsilon^2 N^{-2a/\cexp{holder}}>\vartheta_N^2,
		\end{align}
		which is possible since $ a $ has been chosen such that $a\cexp{holder}^{-1}<\cexp{thetaN}$.
		Suppose, aiming for a contradiction, that there exist $y\in [-C,C]$ and $t\in [0,N\wedge \sigmaN{dist}\wedge \sigmaN{eta}]$ with $|s(v^N_t)(y)|> \epsilon$, and suppose
		$U^N_{N^{-1}\lfloor tN \rfloor, \lfloor y \rfloor,k}\le N^a$.
		Then by Proposition~\ref{prop:tildeudiffunif} and the triangle inequality, 
		for $y'\in [\lfloor y \rfloor ,\lfloor y \rfloor +2]$,
		$$
		|s(v^N_t)(y')-s(v^N_t)(y)|\le N^a |y'-y|^{\cexp{holder}} +|y'-y|\|\partial_{x} m\|_\infty.
		$$
		By our choice of $N$ in~\eqref{choice_N}, for $y'\in [y,y+N^{-2a/\cexp{holder}}]$, the right-hand side is bounded by $ \tfrac{1}{2} \epsilon$, and we have
		$|s(v^N_t)(y')|>\tfrac 12 \epsilon$.
		Hence
		\begin{align*}
			\|s(v^N_t)\|^2_{2,\alpha} &\ge \int_{y}^{y + N^{-2a/\cexp{holder}}} | s(v^N_t)(y') |^2 e^{\alpha y'} dy' \\
			&\ge e^{-\alpha C}(\tfrac 12 \epsilon)^2 N^{-2a/\cexp{holder}} \\
			&> \vartheta_N^2,
		\end{align*}
		where the last inequality follows by \eqref{choice_N}.
		By~\eqref{def_sigma_dist}, this contradicts our assumption that $t\le \sigmaN{dist}\wedge \sigmaN{eta}$.
		We have now established that on the event 
		\[\{U^N_{s,x,k}\le N^a \, \forall x\in \Z\cap [-C-1,C],s\in N^{-1}\Z\cap [0,N]\},\]
		 for $t\le \sigmaN{dist}\wedge \sigmaN{eta} \wedge N$ we have
		$\sup_{y\in [-C,C]}|s(v^N_t)(y)|\le \epsilon$.
		Thus, by \eqref{eq:UNK'bd}, the proof is complete.
	\end{proof}

	\begin{proof}[Proof of Proposition~\ref{prop:control_sigma} in the case $i= \ref{sigmaN:supnorm} $]
	 	As in the proof of Lemma~\ref{lem:scontrol}, recall from~\eqref{def_sigma_dist}-\eqref{def_sigma_sup} and Lemma~\ref{lemma:eta} that
		$\eta(v^N_t)$ and $s(v^N_t)$ are well defined for $t\le \sigmaN{dist}\wedge \sigmaN{eta}$. 
		We recall the constant $\ggamma{weight} < 1$ from \eqref{eq:weightdefn}. By \eqref{eq:masympneg}, we may take $K'>0$ sufficiently large so that
		\begin{align} \label{eq:K'choice1}
			&\Cst{tail_left}e^{-\cexp{tail_left} K'}<\tfrac 16 \Epsilon{proof}
		\end{align}
		and
		\begin{align} \label{eq:K'choice2}
			\sup_{x \leq -K', |y| \leq K} (1-m(x+y))<\tfrac 16 \Epsilon{proof} \quad \text{ and } \quad \sup_{x \geq K', |y| \leq K} e^{-\ggamma{weight} x} m(x+y) <\tfrac 16 \Epsilon{proof}.
		\end{align}
		Let $N$ be sufficiently large so that $(1+\cexp{interface}) \log N \geq K'$ and $N^{-\cexp{tail_left}} < \tfrac 16 \Epsilon{proof}$. Then if $t \leq \sigmaN{dist} \wedge \sigmaN{eta} \wedge \sigmaN{interface}$, 
		$v^N_t$ vanishes above $(1+\cexp{interface})\log N$, and hence for all $x \geq (1+\cexp{interface})\log N$,
		\begin{align}\label{eq:svtNepsbd}
		|s(v^N_t)(x)| \leq m_{\eta(v^N_t)}(x) \leq \sup_{y \in [-K,K]} m(x+y) \leq \frac 12 e^{-\ggamma{weight} x} \Epsilon{proof},
		\end{align}
		where the last inequality uses \eqref{eq:K'choice2}. 
		Next, we note that for $t \leq \sigmaN{dist} \wedge \sigmaN{eta} \wedge \sigmaN{tail_left}$ and $x \leq -K'$, by \eqref{eq:K'choice1}, \eqref{eq:K'choice2} and our choice of $N$,
		\begin{align} \label{eq:svtNepsbd2}
		|s(v^N_t)(x)|&=|(1-v^N_t(x))-(1-m_{\eta(v^N_t)}(x))| \notag
		\\& \le \Cst{tail_left}(e^{\cexp{tail_left} x}+N^{-\cexp{tail_left}})+\sup_{y \in [-K,K]} (1-m(x+y)) < \tfrac 12\Epsilon{proof},
		\end{align}
		Finally, by Lemma~\ref{lem:scontrol} applied with $\epsilon = e^{-\ggamma{weight}} \tfrac 12 \Epsilon{proof}$ and $C = K'$, we have that
		\begin{equation} \label{eq:svtNepsbd3}
		\P{\exists x \in [-K',1], t \in [0, \sigmaN{dist} \wedge \sigmaN{eta} \wedge N] : \frac{|s(v^N_t)(x))|}{e^{-\ggamma{weight} x} \wedge 1} > \frac{\Epsilon{proof}}{2} }  < N^{-2}
		\end{equation}
		for sufficiently large $N$. It remains to handle the range $x \in [1,(1+\cexp{interface})\log N]$. The argument is similar to that used to prove Lemma~\ref{lem:scontrol},
		but uses Proposition~\ref{prop:holderexp} instead of Proposition~\ref{prop:tildeudiffunif}.
		
		We again will use $\ggamma{weight}$ from \eqref{eq:weightdefn}, and recall $A_{(0)}$ from Proposition~\ref{prop:holderexp}.
		Let $t\geq0$ and suppose that for all $x_1, x_2 \in [0, (1+\cexp{interface})\log N]$ with $|x_1 -x_2| \leq 1$, 
		\begin{equation} \label{eq:scontrolproofexpholder}
		|v^N_{t}(x_1) - v^N_{t}(x_2)| \leq A_{(0)} N^{\cexp{tail_right_N}} e^{-\ggamma{weight}x_1 } |x_1 - x_2|^{\cexp{holder}}.
		\end{equation}
		Now suppose that there exists $y \in  [0, (1+\cexp{interface})\log N]$ such that 
		\begin{equation*}
		s(v^N_t)(y) > \tfrac{1}{2}\Epsilon{proof} e^{-\ggamma{weight} y}.
		\end{equation*}
		By \eqref{eq:scontrolproofexpholder}, for $y \in [y',y+1]$
		\begin{align*}
		&|s(v^N_t)(y) - s(v^N_t)(y')|  \notag
		\\& \hspace{ 1cm }\leq  e^{-\ggamma{weight} y} \left( A_{(0)}N^{\cexp{tail_right_N}}  |y' -y|^{\cexp{holder}} + |y' - y| \|\partial_x m\|_{\infty,\ggamma{weight}}\right)\notag
		\\ &\hspace{ 1cm } \leq 2e^{-\ggamma{weight} y} A_{(0)}N^{\cexp{tail_right_N}} |y' -y|^{\cexp{holder}},
		\end{align*}
		where the second inequality holds for sufficiently large $N$. We remark that because $\ggamma{weight} < 1$, $\|\partial_x m\|_{\infty,\ggamma{weight}} < \infty$ by\eqref{eq:masympneg}.
		Now let 
		\[ k_N(y) := \left(\frac{\Epsilon{proof}}{4A_{(0)}}\right)^{1/\cexp{holder}} N^{-\cexp{tail_right_N}/\cexp{holder}}.\]
		The previous inequality then implies that $|s(v^N_t)(y) - s(v^N_t)(y')| \geq \tfrac{1}{4} e^{-\ggamma{weight} y} \Epsilon{proof}$ for all $y' \in [y,y + k_N(y)]$. In particular,
		\begin{align*}
			\|s(v^N_t)\|^2_{2,\alpha} &\ge \int_{y}^{y + k_N(y)} | s(v^N_t)(y') |^2 e^{\alpha y'} dy' \\
			&\ge k_N(y) \left(\tfrac 14 \Epsilon{proof}\right)^2 e^{(\alpha - 2\lambda)y}\\
			& = \frac{N^{-\cexp{tail_right_N}/\cexp{holder}}}{16(4A_{(0)})^{1/\cexp{holder}}} \Epsilon{proof}^{2 + 1/\cexp{holder}} e^{(\alpha - 2\ggamma{weight})y}\\
			& \geq \frac{N^{-\cexp{tail_right_N}/\cexp{holder}}}{16(4A_{(0)})^{1/\cexp{holder}}} \Epsilon{proof}^{2 + 1/\cexp{holder}}\\
			& > \vartheta_N^{2}.
		\end{align*}
		In the second last line, we have used the fact that $2\alpha - 2\ggamma{weight} > 0$ from \eqref{eq:weightdefn} and $y \geq 0$.
		The last line, which holds for sufficiently large $N$, follows because $\cexp{tail_right_N} < \cexp{holder} \cexp{thetaN}$, from \eqref{bounds_cexptailrightN},
		and $\vartheta_N=N^{-\cexp{thetaN}}\wedge \Beta{proof}$.
		Thus, we have shown that if \eqref{eq:scontrolproofexpholder} holds at some $t \geq 0$, 
		then either $t > \sigmaN{dist}$ or $s(v^N_t)(y) \leq \tfrac{1}{2}\Epsilon{proof} e^{-\ggamma{weight} y}$ for all $y \in  [1, (1+\cexp{interface})\log N]$.
		By \eqref{bounds_cexptailrightexp}, we have $\ggamma{weight} < (1-\cexp{tail_right_exp})/2$; it therefore follows from the above and Proposition~\ref{prop:holderexp} that
		\begin{align*} 
		&\P{\exists x \in [0,(1+\cexp{interface})\log N], t \in [0, \sigmaN{dist} \wedge \sigmaN{eta} \wedge \sigmaN{interface} \wedge \sigmaN{tail_right} \wedge N] : 
		\frac{|s(v^N_t)(x))|}{e^{-\ggamma{weight} x} \wedge 1} > \frac{\Epsilon{proof}}{2} }  
		\\ &\hspace{1 cm}< N^{-1}.
		\end{align*}
		Combining this with~\eqref{eq:svtNepsbd},~\eqref{eq:svtNepsbd2} and~\eqref{eq:svtNepsbd3}, it follows that for $N$ sufficiently large,
		\begin{equation*}
			\P{\sigmaN{supnorm}\le \sigmaN{dist}\wedge \sigmaN{eta}\wedge  \sigmaN{tail_right}\wedge \sigmaN{interface} \wedge \sigmaN{tail_left}\wedge N}\le 2N^{-1}.
		\end{equation*}
		By a union bound, and then by Propositions~\ref{prop_R_main}, \ref{prop:sigmaNrighttail} and~\ref{prop:sigmatailleft}, it follows that for $T>0$, for $N$ sufficiently large,
		\begin{align*}
			&\P{\sigmaN{supnorm}\le \sigmaN{dist}\wedge \sigmaN{eta}\wedge T}\\
				&\quad \le 2N^{-1}+ \P{\sigmaN{tail_right}\le  \sigmaN{interface}\wedge T}
			\\	&\qquad + \P{\sigmaN{tail_left}\le  \sigmaN{dist} \wedge \sigmaN{eta} \wedge \wedge T} + \P{\sigmaN{interface} \leq \sigmaN{eta} \wedge \sigmaN{supnorm} \wedge T} 
			\\	&\quad = o(1),
		\end{align*}
		which completes the proof.
	\end{proof}
	
	\begin{proof}[Proof of Proposition~\ref{prop:control_sigma} in the case $i=\ref{sigmaN:vdelta}$]
		Recall from \eqref{def_sigma_vdelta} that
		\begin{equation} \label{eq:sigma4Nremind}
			\sigmaN{vdelta}= \inf \left\lbrace t \geq 0 : \| D^{\delta_N,N}(t,\cdot) \|_{{2,\alpha}} > r_N \text{ or } \| v^{\delta_N,N}_t - v^N_t \|_\infty > r_N \right\rbrace.
		\end{equation}
		By the definition of $D^{\delta,N}(t,x)$ in~\eqref{eq:Ddefn}, and then by Jensen's inequality and since $\int_{\R} \rho^{\delta_N}(z)dz=1$, for $t\ge 0$ we have
		\begin{align} \label{eq:Dnormbound}
			\|D^{\delta_N,N}(t,\cdot)\|^2_{2,\alpha}
			&= \int_{\R} e^{\alpha x} \left(\int_{\R} \rho^{\delta_N}(x-y)|v^N_t(y)-v^N_t(x)| dy
			\right)^2 dx \notag \\
			&\leq \int_{\R} e^{\alpha x}\int_{\R} \rho^{\delta_N}(x-y)(v^N_t(y)-v^N_t(x))^2 dy
			\, dx \notag \\
			&\le \frac{\|\rho\|_\infty}{\delta_N}\int_{\R} e^{\alpha x}\int_{-\delta_N}^{\delta_N} (v^N_t(x+z)-v^N_t(x))^2 dz\,
			dx,
		\end{align}
		where the last line follows by the definition of $\rho^{\delta_N}$ in~\eqref{eq:rhodeltadefn}.
		Recall from~\eqref{bounds_cexpdeltaN} that $\cexp{deltaN}>0$, and fix $k\in \N$ sufficiently large that $\frac 12 \cexp{deltaN}\cexp{holder}k>2$.
		By Proposition~\ref{prop:tildeudiffunif}, for $t\in [0,N]$ and $x,y \in [-\log N-1,\log N+1]$ with $|x-y|\le \delta_N$ we have
		\begin{equation} \label{eq:vxvyclose}
		|v^N_t(x)-v^N_t(y)|\le (\delta_N)^{\cexp{holder}} \max_{x'\in \Z \cap [-\log N-2,\log N], \, t'\in N^{-1} \Z\cap [0,N]}U^N_{t',x',k}.
		\end{equation}
		Now recall from~\eqref{def:kappa_delta_r_theta} that $\delta_N=N^{-\cexp{deltaN}}$,
		and suppose the event 
		\begin{equation} \label{eq:eventENdefn}
			E_N:=\left\{\max_{x'\in \Z \cap [-\log N-2,\log N],\, t'\in N^{-1} \Z\cap [0,N]}U^N_{t',x',k} \le N^{\frac 12 \cexp{deltaN} \cexp{holder}}=(\delta_N)^{-\cexp{holder}/2}\right\}
		\end{equation}
		occurs.
		Take 
		$t\in [0,\sigmaN{tail_right}\wedge N]$; by~\eqref{eq:sigmatailrightdefn} we have $v^N_t(x)\le N^{\cexp{tail_right_N}} e^{-(1-\cexp{tail_right_exp}) x}$ for $x\ge a \log N -\delta_N$ for any $a >0$.
		We fix $a \in (0,1)$, and will choose a certain value shortly. By \eqref{eq:vxvyclose} and since $v_t^N\in [0,1]$, we have
		\begin{align} \label{eq:v-v2bd}
			&\int_{\R} \int_{-\delta_N}^{\delta_N}
			(v_t^N(x+y)-v^N_t(x))^2 dy \, e^{\alpha x}dx  \notag
			\\ &\hspace{1.5 cm} \leq \int_{-\infty}^{-a \log N} \int_{-\delta_N}^{\delta_N} e^{\alpha x} dy \, dx + \int_{-a \log N}^{a \log N} \int_{-\delta_N}^{\delta_N} (\delta_N)^{\cexp{holder}}dy \, e^{\alpha x}dx \notag
			\\ &\hspace{1.5 cm} \quad  + \int_{a \log N}^\infty \int_{-\delta_N}^{\delta_N} e^{2\delta_N} N^{2\cexp{tail_right_N}} e^{-2(1-\cexp{tail_right_exp})x}dy \, e^{\alpha x} dx \notag
			\\ &\hspace{1.5 cm}\leq 2\delta_N \alpha^{-1} N^{-\alpha a} + (\delta_N)^{\cexp{holder}} 2\delta_N \alpha^{-1} N^{\alpha a} \notag
			\\ &\hspace{1.5 cm} \quad + 2\delta_Ne^{2\delta_N} N^{2\cexp{tail_right_N}}\cdot (2(1-\cexp{tail_right_exp}) - \alpha)^{-1} N^{-a(2(1-\cexp{tail_right_exp}) - \alpha)}.
		\end{align}

		Now choose $a =  1 \wedge ( \cexp{deltaN} \cexp{holder}/ (2\alpha)) > 0$. Then
		\[ N^{2\cexp{tail_right_N} -a(2(1-\cexp{tail_right_exp}) - \alpha)} = N^{- 2 (( \frac 1 2 \wedge (\frac{\cexp{deltaN}\cexp{holder}}{4\alpha} ))(2(1-\cexp{tail_right_exp}) -\alpha) - \cexp{tail_right_N})}.\]
		Furthermore, by \eqref{def:kappa_delta_r_theta} we have
		\[N^{-\alpha a}  \leq N^{-( (\frac{\cexp{deltaN}\cexp{holder}}{2} )\wedge  \frac \alpha 2 )} \quad \text{ and } \quad  \delta_N^{ \cexp{holder}} N^{\alpha a}  \leq N^{-\frac{\cexp{deltaN}\cexp{holder}}{2}}.\]
		It follows from \eqref{eq:cexprNsmall} and \eqref{bounds_cexptailrightN} that the exponent in each term above is less than $-2\cexp{rN}$. 
		In particular, returning to \eqref{eq:v-v2bd}, we obtain that there exists $c > 2\cexp{rN}$ such that for $N$ sufficiently large,	
		$$
		\int_{\R} \int_{-\delta_N}^{\delta_N}
		(v_t^N(x+y)-v^N_t(x))^2 dy \, e^{\alpha x}dx
		\le \delta_N N^{-c}.
		$$
		Therefore by~\eqref{eq:Dnormbound}, for $N$ sufficiently large, on the event $E_N$, for $t\in [0, \sigmaN{tail_right}\wedge N]$ we have
		\begin{equation} \label{eq:Ddeltabound}
		\|D^{\delta_N,N}(t,\cdot)\|^2_{2,\alpha}\le \|\rho\|_\infty N^{-c}.
		\end{equation}
		By~\eqref{eq:eventENdefn}, a union bound, Markov's inequality and Proposition~\ref{prop:tildeudiffunif}, and then since we chose $k\in \N$ with $\frac 12 \cexp{deltaN} \cexp{holder}k>2$ we have that
		\begin{align} \label{eq:UNnotbig}
		\P{(E_N)^c}
			&\leq (N^2+1)(2\log N +3)N^{-2k\cdot \frac 12 \cexp{deltaN} \cexp{holder}} A_k \notag \\
			&\leq N^{-1}
		\end{align}
		for $N$ sufficiently large.
		Therefore, since $r_N=N^{-\cexp{rN}}$ by~\eqref{def:kappa_delta_r_theta} and since we chose $c>2\cexp{rN}$, by~\eqref{eq:Ddeltabound}, for $N$ sufficiently large,
		\begin{equation} \label{eq:DNbound}
		\P{\exists\, t\leq \sigmaN{tail_right}\wedge N : \|D^{\delta_N,N}(t,\cdot)\|_{2,\alpha} \ge r_N } \leq N^{-1}.
		\end{equation}
		Suppose
		again that the event $E_N$ occurs, as defined in~\eqref{eq:eventENdefn}.
		Then by the definition of $v^{\delta_N,N}_t$ in~\eqref{def_vdelta} and then by~\eqref{eq:vxvyclose} and since $\|\rho^{\delta_N}\|_1=1$,
		for $t\in [0,N]$ and $x\in [-\log N,\log N]$, we have that
		\begin{align} \label{eq:vdelta-v}
		|v^{\delta_N,N}_t(x)-v^N_t(x)|&\le 
		\int_{-\delta_N}^{\delta_N} \rho^{\delta_N}(y) |v^N_t(x-y)-v^N_t(x)| dy
		\le (\delta_N)^{\cexp{holder}/2}.
		\end{align}
		For $t\le \sigmaN{tail_right}\wedge \sigmaN{tail_left}$,
		for $x\le -\log N$, by the definition of $\sigmaN{tail_left}$ in~\eqref{eq:sigmalefttaildefn} we have
		\begin{equation} \label{eq:vdelta-v2}
		|v^{\delta_N,N}_t(x)-v^N_t(x)|\leq (1-v^{\delta_N,N}_t(x))+(1-v^N_t(x))
		\le   \Cst{tail_left} (e^{\cexp{tail_left}\delta_N}+3)N^{-\cexp{tail_left}},
		\end{equation}
		and for $x\ge \log N$, by the definition of $\sigmaN{tail_right}$ in~\eqref{eq:sigmatailrightdefn} we have
		\begin{equation} \label{eq:vdelta-v3}
		|v^{\delta_N,N}_t(x)-v^N_t(x)|\leq v^{\delta_N,N}_t(x)+v^N_t(x)
		\le N^{\cexp{tail_right_N}} (e^{(1-\cexp{tail_right_exp}) \delta_N}+1) N^{-(1-\cexp{tail_right_exp})}.
		\end{equation}
		Recall from \eqref{eq:cexprNsmall} and \eqref{bounds_cexptailrightN} that $\cexp{rN} < (\frac 12 \cexp{holder}\cexp{deltaN})\wedge (1- \cexp{tail_right_exp} - \cexp{tail_right_N}) \wedge \cexp{tail_left}$, 
		and that $\delta_N=N^{-\cexp{deltaN}}$, $r_N=N^{-\cexp{rN}}$ by~\eqref{def:kappa_delta_r_theta}.
		Hence combining~\eqref{eq:vdelta-v},~\eqref{eq:vdelta-v2} and~\eqref{eq:vdelta-v3}, and then using~\eqref{eq:UNnotbig}, for $N$ sufficiently large
		\begin{equation} \label{eq:vdiffrN}
		\P{\exists\, t\le \sigmaN{tail_right}\wedge \sigmaN{tail_left}\wedge N : \|v^{\delta_N,N}_t-v^N_t\|_{\infty} \ge r_N}\le \P{(E_N)^c}\le N^{-1}.
		\end{equation}
		By a union bound, and then, recalling~\eqref{eq:sigma4Nremind},
		by~\eqref{eq:DNbound} and~\eqref{eq:vdiffrN}, and by Propositions~\ref{prop_R_main}, \ref{prop:sigmaNrighttail} and \ref{prop:sigmatailleft},
		we now have that for $T>0$, for $N$ sufficiently large,
		\begin{align*}
		&\P{\sigmaN{vdelta}\le \sigmaN{dist}\wedge \sigmaN{eta}\wedge \sigmaN{supnorm}\wedge T}\\
		&\le \P{\sigmaN{vdelta}\le  \sigmaN{tail_right}\wedge \sigmaN{tail_left}\wedge T}
		+\P{\sigmaN{tail_right}\le \sigmaN{dist}\wedge \sigmaN{eta}  \wedge \sigmaN{supnorm} \wedge \sigmaN{interface} \wedge T}\\
		&\qquad+\P{\sigmaN{tail_left}\le \sigmaN{dist}\wedge \sigmaN{eta}  \wedge T} 
		+  \P{\sigmaN{interface}\le  \sigmaN{eta} \wedge \sigmaN{supnorm} \wedge T }\\
		&\le 2N^{-1}+o(1),
		\end{align*}
		which completes the proof.
	\end{proof}
	
	\begin{proof}[Proof of Proposition~\ref{prop:control_sigma} with $i=\ref{sigmaN:tail}$]
	Recall from \eqref{def_sigma_tail} that
	\begin{equation*}
		\sigmaN{tail}=  \inf \left\lbrace t \geq 0 : \int_\R v^N_t(x) e^{\alpha x} dx  > \kappa_N \right\rbrace.
	\end{equation*}
	Take $t\in [0,\sigmaN{interface} \wedge \sigmaN{tail_right} ]$;
	by \eqref{eq:sigmatailrightdefn} and \eqref{eq:sigmainterface} we have that
	\begin{align} \label{eq:intvNsig68}
	\int_{\R} v^N_t(x)e^{\alpha x}dx &\leq \int_{-\infty}^0 e^{\alpha x}dx +\int_{0}^{(1+\cexp{interface})\log N} N^{\cexp{tail_right_N}} e^{-(1-\cexp{tail_right_exp})x} e^{\alpha x}dx 
	 \\ &\leq \alpha^{-1} + \frac{1}{\alpha + \cexp{tail_right_exp} - 1} N^{\cexp{tail_right_N}} N^{(1+ \cexp{interface})((\alpha + \cexp{tail_right_exp} - 1) \vee 0)},	\notag \end{align}
	provided $\alpha + \cexp{tail_right_exp} - 1 > 0$. If $\alpha  + \cexp{tail_right_exp} -1 \leq 0$, $N^{(1+ \cexp{interface})((\alpha + \cexp{tail_right_exp} - 1) \vee 0)}$ is replaced
	with either $\log N$ in the case of equality, or otherwise a negative power of $N$. In particular, the case shown gives the worst case (largest) behaviour in $N$. 
	Hence, by \eqref{eq:cKNcsum} and the bound above, for $N$ sufficiently large, for $t\in [0, \sigmaN{interface} \wedge \sigmaN{tail_right} ]$,
	\begin{align*}
	\int_{\R} v^N_t(x)e^{\alpha x}dx  <\kappa_N,
	\end{align*}
	since $\kappa_N =N^{1-\cexp{KN}}$ by \eqref{def:kappa_delta_r_theta}.
	It follows that for $N$ sufficiently large we have $\sigmaN{tail}>  \sigmaN{interface}\wedge \sigmaN{tail_right}$. Therefore
	for $T>0$, for $N$ sufficiently large, by a union bound,
	\begin{align*}
	&\P{\sigmaN{tail}\le \sigmaN{dist}\wedge \sigmaN{eta}\wedge \sigmaN{supnorm} \wedge T}\\
	&\hspace{1 cm}\le \P{\sigmaN{interface}\le \sigmaN{eta} \wedge \sigmaN{supnorm} \wedge T} +\P{\sigmaN{tail_right}\le \ \sigmaN{interface} \wedge T} = o(1),
	\end{align*}
	where to conclude we have used Propositions~\ref{prop_R_main} and~\ref{prop:sigmaNrighttail}. \end{proof}

	It remains to prove Propositions~\ref{prop:tildeudiffunif}-\ref{prop:holderexp}; we will prove these results in Sections~\ref{subsec:holder}-\ref{subsec:expholder} below.
	The proof of a technical lemma in Section~\ref{subsec:holder} is postponed to Section~\ref{subsec:techG-G}.

	\subsection{Proof of Proposition~\ref{prop:tildeudiffunif}} \label{subsec:holder}
	
	It will be convenient in the next few subsections to work with the following rescaled process:
	recalling~\eqref{def_vN},
	for $x\in \R$ and $t\ge 0$, let
	\begin{equation} \label{eq:tildeudefn}
	 w^N_t(x):=u^N_t(x+\alpha t)=v^N_{t/N}(x).
	\end{equation}
	(We have already used this notation for a generic solution to \eqref{eq_spde_movingframe}; hereafter, $w^N_t$ refers to the solution specified by \eqref{eq:tildeudefn}.)
	By \eqref{mild_spde_wN}, we obtain that for $a\in \R$, for $t> 0$ and $x\in \R$,
	\begin{multline} \label{eq:utildeformula}
	w^N_t = e^{a t} Q(t) u^N_0 + \int_{0}^{t} e^{a(t-s)} Q(t-s) \left( f(w^N_s) - a w^N_s \right) ds \\ + N^{-1/2} \int_{0}^{t} e^{a(t-s)} Q(t-s) \left( \sqrt{w^N_s(1-w^N_s)} dW(s) \right),
	\end{multline}
	where $(W(s),s\ge 0)$ is a cylindrical Wiener process.
	The following technical lemma will be used in the proofs of Propositions~\ref{prop:tildeudiffunif},~\ref{prop:sigmaNrighttail} and~\ref{prop:holderexp}; its proof is postponed to Section~\ref{subsec:techG-G}.
	Recall the definition of $G_t(x)$ in~\eqref{eq:Gdefn}.
	\begin{lemma} \label{lem:Gdiffintbound}
	For $0\le s< t_1\le t_2$ and $x\in \R$, let
	\begin{equation} \label{eq:Idefn}
	 I_{t_1,t_2}(s,x):=\int_{\R} (G_{t_1-s}(y-x)-G_{t_2-s}(y))^2 dy.
	\end{equation}
	Then for $a>0$, there exists a constant $\newCst{I}=\Cst{I}(a)<\infty$ such that for $x\in \R$ and $0\le t_1 \le t_2$,
	\begin{align} \label{eq:Gdiff1}
	\int_0^{t_1} e^{-a(t_1-s)} I_{t_1,t_2}(s,x) \,ds
	&\leq \Cst{I}((t_2-t_1)^{1/4}+|x|^{1/2}),
	\end{align}
	and there exists a constant $\newCst{I2}<\infty$ such that for $t_0\in (0,\infty)$, $0\le t_1 \le t_2$ and $x\in \R$, if $0\le s\le t_1-t_0$ then
	\begin{equation} \label{eq:Gdiff2}
	I_{t_1,t_2}(s,x)
	\leq \Cst{I2} \, t_0^{-3/2} (x^2+(t_2-t_1)).
	\end{equation}
	\end{lemma}
	We can use Lemma~\ref{lem:Gdiffintbound} to prove the following moment bound. We recall the constant $\cinit{big}>0$ as appearing in Assumption~\ref{assumpt:v0}.\ref{v0:holder}.
	\begin{lemma} \label{lem:kolmest}
	For $k\in \N$ there exists a constant $A_k^{(1)}<\infty$ such that for $0\le t_1\le t_2$ and $x_1,x_2\in \R$ with $t_2-t_1\leq 1$ and $|x_1-x_2|\leq 1$,
	\begin{align*}
	\E{(w^N_{t_1}(x_1)-w^N_{t_2}(x_2))^{2k}}
	&\leq A^{(1)}_k (|x_1-x_2|^{(\frac 1 4 \wedge \cinit{big})2k}+(t_2-t_1)^{(\frac 1 4 \wedge \cinit{big}) k}).
	\end{align*}
	\end{lemma}
	Before proving Lemma~\ref{lem:kolmest}, we show that Proposition~\ref{prop:tildeudiffunif} follows easily from Lemma~\ref{lem:kolmest}.
	\begin{proof}[Proof of Proposition~\ref{prop:tildeudiffunif}]
	Recall the definition of $\|\cdot \|_2$ in~\eqref{eq:2normdefn}.
	By Lemma~\ref{lem:kolmest} and Kolmogorov's continuity criterion (see Corollary~1.2 in~\cite{walsh_introduction_1986}), for $k\in \N$ with $k>8\vee (2\cinit{big}^{-1})$ there exist constants $A^{(2)}_k<\infty$ and $A^{(3)}_k<\infty$ such that for $x\in \R$ and $t\ge 0$, there exists a random variable $\tilde U^N_{t,x,k}$ such that for any $x_1$, $x_2\in [x,x+2]$ and $t_1,t_2 \in [t,t+2N^{-1}]$,
	\begin{align*}
	&|w^N_{Nt_1}(x_1)-w^N_{Nt_2}(x_2)|\\
	&\quad \leq \tilde U^N_{t,x,k} \|(x_1,Nt_1)-(x_2,Nt_2)\|_2^{(1/4\wedge \cinit{big})-1/k}(\log (A_k^{(2)}/\|(x_1,Nt_1)-(x_2,Nt_2)\|_2))^{1/k},
	\end{align*}
	and moreover $\E{(\tilde U^N_{t,x,k})^{2k}}\leq A^{(3)}_k$.
	Since $v^N_t(x)= w^N_{Nt}(x)$ by~\eqref{eq:tildeudefn}, and $\cexp{holder} \in (0, \tfrac 1 4 \wedge \cinit{big})$, the result follows.
	\end{proof}
	\begin{proof}[Proof of Lemma~\ref{lem:kolmest}]
	Take $0\leq t_1\leq t_2$ with $t_2-t_1\le 1$ and $x_1,x_2\in \R$ with $|x_1-x_2|\le 1$ such that $(x_1,t_1)\neq (x_2,t_2)$.
	By~\eqref{eq:Qdefn} and~\eqref{eq:utildeformula} with $a=-1$,
	\begin{align} \label{eq:(G)}
	&w^N_{t_1}(x_1)-w^N_{t_2}(x_2) \notag\\
	&=\int_{\R} \left(G_{t_1}(x_1+\alpha t_1-y)e^{-t_1}
	-G_{t_2}(x_2+\alpha t_2-y)e^{-t_2}\right)
	u_0^N(y)dy \notag\\
	&\qquad +\int_0^{t_1} \int_{\R} \bigg(G_{t_1-s}(x_1+\alpha (t_1-s)-y)e^{-(t_1-s)}-G_{t_2-s}(x_2+\alpha (t_2-s)-y)e^{-(t_2-s)}\bigg) \notag\\
	&\hspace{11cm}(f(w^N_s(y))+w^N_s(y))dyds \notag\\
	&\qquad -\int_{t_1}^{t_2} \int_{\R} G_{t_2-s}(x_2+\alpha (t_2-s)-y)e^{-(t_2-s)}
	(f(w^N_s(y))+ w^N_s(y))dyds \notag\\
	&\qquad + N^{-1/2}\int_0^{t_1} \int_{\R} \bigg(G_{t_1-s}(x_1+\alpha (t_1-s)-y)e^{-(t_1-s)} \notag \\
	&\hspace{4cm}-G_{t_2-s}(x_2+\alpha (t_2-s)-y)e^{-(t_2-s)}\bigg)
	\sqrt{w^N_s(y)(1- w^N_s(y))}W(dy ds)\notag \\
	&\qquad - N^{-1/2}\int_{t_1}^{t_2} \int_{\R} G_{t_2-s}(x_2+\alpha (t_2-s)-y)e^{-(t_2-s)}
	\sqrt{w^N_s(y)(1- w^N_s(y))}W(dy ds).
	\end{align}
	For the first term on the right-hand side of~\eqref{eq:(G)}, by the triangle inequality,
	\begin{align} \label{eq:1sttermholder0}
	&\left| \int_{\R} \left(G_{t_1}(x_1+\alpha t_1-y)e^{-t_1}
	-G_{t_2}(x_2+\alpha t_2-y)e^{-t_2}\right)
	u_0^N(y)dy\right| \notag \\
	&\leq \left| \int_{\R} (G_{t_1}(x_1+\alpha t_1-y)
	-G_{t_2}(x_2+\alpha t_2-y))e^{-t_1}
	u_0^N(y)dy\right| \notag \\
	&\quad +\left| \int_{\R} (e^{-t_1}
	-e^{-t_2})G_{t_2}(x_2+\alpha t_2-y)
	u_0^N(y)dy\right| \notag \\
	&\leq e^{-t_1}\left| \int_{\R} G_{t_1}(x_1+\alpha t_1-y)
	\int_{\R} G_{t_2-t_1}(y-z)(u_0^N(y)-u_0^N(z+x_2-x_1+\alpha(t_2-t_1)))dz dy\right| \notag \\
	&\quad + (e^{-t_1}
	-e^{-t_2}),
	\end{align}
	where the second inequality follows
	since $|u_0^N(y)|\le 1$ $\forall y\in \R$ and since $t_1\le t_2$.
	Recall from Assumption~\ref{assumpt:v0}.\ref{v0:holder} that 
	for all $ x, y \in \R $ we have
		$| u^N_0(x)-u^N_0(y) | \leq \Cst{init} |x-y|^{\cinit{big}}$.
	We can use this to bound the first term on the right-hand side above and write
	\begin{align} \label{eq:1sttermholder}
	&\left| \int_{\R} \left(G_{t_1}(x_1+\alpha t_1-y)e^{-t_1}
	-G_{t_2}(x_2+\alpha t_2-y)e^{-t_2}\right)
	u_0^N(y)dy\right| \notag \\
	&\leq \int_{\R} G_{t_1}(x_1+\alpha t_1-y)
	\int_{\R}  G_{t_2-t_1}(y-z)\Cst{init}|z+x_2-x_1+\alpha (t_2-t_1)-y|^{\cinit{big}} dz
	dy +(t_2-t_1) \notag \\
	&\leq 
	\Cst{init}|x_2-x_1+\alpha (t_2-t_1)|^{\cinit{big}}+\Cst{init} \mathbb E_0 \left[|B_{2(t_2-t_1)}|^{\cinit{big}} \right]
	+ (t_2-t_1)\notag \\
	&=
	\Cst{init}|x_2-x_1+\alpha (t_2-t_1)|^{\cinit{big}}+\Cst{init} (t_2-t_1)^{\cinit{big}/2}\mathbb E_0 \left[|B_{2}|^{\cinit{big}} \right]
	+ (t_2-t_1),
	\end{align}
	where $(B_t)_{t\ge 0}$ is a Brownian motion, the second line follows since $\cinit{big}<1$ and so $|a+b|^{\cinit{big}}\le |a|^{\cinit{big}}+|b|^{\cinit{big}}$ $\forall a,b\in \R$, and the third line follows by Brownian scaling.
	For the second term on the right-hand side of~\eqref{eq:(G)},
	since $|f(u)+u|\le 2+\alpha$ for $u\in [0,1]$ (recall the definition of $f$ in~\eqref{f}), and then by the triangle inequality,
	\begin{align} \label{eq:(H)}
	&\bigg|\int_0^{t_1} \int_{\R} \bigg(G_{t_1-s}(x_1+\alpha (t_1-s)-y)e^{-(t_1-s)}-G_{t_2-s}(x_2+\alpha (t_2-s)-y)e^{-(t_2-s)}\bigg) \notag\\
	&\hspace{11cm}(f(w^N_s(y))+ w^N_s(y))dyds\bigg| \notag\\
	&\leq (2+\alpha) \int_0^{t_1} e^{-(t_1-s)}
	\int_{\R} \left|
	G_{t_1-s}(x_1+\alpha (t_1-s)-y)-G_{t_2-s}(x_2+\alpha (t_2-s)-y)e^{-(t_2-t_1)}
	\right|dy ds \notag \\
	&\leq (2+\alpha) \int_0^{t_1} e^{-(t_1-s)}
	\int_{\R} \left|
	G_{t_1-s}(x_1+\alpha (t_1-s)-y)-G_{t_2-s}(x_2+\alpha (t_2-s)-y)
	\right|dy ds \notag \\
	&\quad + (2+\alpha) \int_0^{t_1} e^{-(t_1-s)}(1-e^{-(t_2-t_1)})
	\int_{\R} 
	G_{t_2-s}(x_2+\alpha (t_2-s)-y)
	dy ds \notag \\
	&\leq (2+\alpha) \int_0^{t_1} e^{-(t_1-s)}
	\int_{\R} \left|
	G_{t_1-s}(x_1+\alpha (t_1-s)-y)-G_{t_2-s}(x_2+\alpha (t_2-s)-y)
	\right|dy ds \notag \\
	&\quad + (2+\alpha) (t_2-t_1).
	\end{align}
	We now bound the first term on the right-hand side of~\eqref{eq:(H)}.
	By H\"older's inequality, for $s\in [0,t_1)$,
	\begin{align} \label{eq:GintIbd}
	&\int_{\R} \left|
	G_{t_1-s}(x_1+\alpha (t_1-s)-y)-G_{t_2-s}(x_2+\alpha (t_2-s)-y)
	\right|dy \notag \\
	&\leq \bigg(\int_{\R} \left|
	G_{t_1-s}(x_1+\alpha (t_1-s)-y)-G_{t_2-s}(x_2+\alpha (t_2-s)-y)
	\right|^{1/2} dy \bigg)^{2/3} \notag \\
	& \qquad \cdot \bigg(\int_{\R} \left(
	G_{t_1-s}(x_1+\alpha (t_1-s)-y)-G_{t_2-s}(x_2+\alpha (t_2-s)-y)
	\right)^2 dy\bigg)^{1/3} \notag \\
	&\leq \bigg(\int_{\R}
	G_{t_1-s}(x_1+\alpha (t_1-s)-y)^{1/2}dy +\int_{\R} G_{t_2-s}(x_2+\alpha (t_2-s)-y)
	^{1/2} dy \bigg)^{2/3} \notag \\
	& \qquad \cdot \big( I_{t_1,t_2}(s,x_1-x_2+\alpha(t_1-t_2))\big)^{1/3},
	\end{align}
	where the second inequality follows by the definition of $I_{t_1,t_2}(s,x)$ in~\eqref{eq:Idefn} and since $|a+b|^{1/2}\le |a|^{1/2}+|b|^{1/2}$ $\forall a,b\in \R$.
	For $t>0$, we can calculate
	\begin{align*}
	\int_{\R}
	G_{t}(y)^{1/2}dy
	&= \int_{\R} (4\pi t)^{-1/4} e^{-y^2/(8t)}dy
	=2\pi^{1/4}t^{1/4}.
	\end{align*}
	Therefore, by~\eqref{eq:GintIbd} and using~\eqref{eq:Gdiff2} in Lemma~\ref{lem:Gdiffintbound}, for $t_0\in (0,\infty)$, if $0\le s\le t_1-t_0$ then
	\begin{align*}
	&\int_{\R} \left|
	G_{t_1-s}(x_1+\alpha (t_1-s)-y)-G_{t_2-s}(x_2+\alpha (t_2-s)-y)
	\right|dy \\
	&\leq 
	(2\pi^{1/4}((t_1-s)^{1/4}+(t_2-s)^{1/4}))^{2/3}(\Cst{I2} t_0^{-3/2}((x_1-x_2+\alpha(t_1-t_2))^2 +(t_2-t_1)))^{1/3}.
	\end{align*}
	Also, for any $s\in [0,t_1)$,
	$$
	\int_{\R} \left|
	G_{t_1-s}(x_1+\alpha (t_1-s)-y)-G_{t_2-s}(x_2+\alpha (t_2-s)-y)
	\right|dy \le 2.
	$$
	Hence by~\eqref{eq:(H)}, for any choice of $t_0\in (0,t_1]$, we can write
	\begin{align}
	&\bigg|\int_0^{t_1} \int_{\R} \bigg(G_{t_1-s}(x_1+\alpha (t_1-s)-y)e^{-(t_1-s)}-G_{t_2-s}(x_2+\alpha (t_2-s)-y)e^{-(t_2-s)}\bigg) \notag\\
	&\hspace{11cm}(f(w^N_s(y))+ w^N_s(y))dyds\bigg|\notag \\
	&\leq  (2+\alpha) \int_0^{t_1-t_0}e^{-(t_1-s)}4^{2/3}\pi^{1/6}(t_2-s)^{1/6}
	\Cst{I2}^{1/3}t_0^{-1/2}((x_1-x_2+\alpha(t_1-t_2))^2+(t_2-t_1))^{1/3}
	ds\notag \\
	&\quad +(2+\alpha)\int_{t_1-t_0}^{t_1}2 ds +(2+\alpha)(t_2-t_1) \label{eq:GGint1}\\
	&\leq (2+\alpha)4^{2/3}\pi^{1/6}
	\Cst{I2}^{1/3}t_0^{-1/2}((x_1-x_2+\alpha(t_1-t_2))^2+(t_2-t_1))^{1/3}\int_0^\infty (s'+t_2-t_1)^{1/6}e^{-s'}ds' \notag \\
	&\quad +2(2+\alpha)t_0+(2+\alpha)(t_2-t_1). \label{eq:GGint2}
	\end{align}
	Therefore, setting $t_0:=\min(t_1,((x_1-x_2+\alpha(t_1-t_2))^2+(t_2-t_1))^{2/9})$, by~\eqref{eq:GGint2} (and using~\eqref{eq:GGint1} for the case $t_0=t_1$) we have that there exists a constant $C<\infty$ such that
	\begin{align} \label{eq:2ndtermholder}
	&\bigg|\int_0^{t_1} \int_{\R} \bigg(G_{t_1-s}(x_1+\alpha (t_1-s)-y)e^{-(t_1-s)}-G_{t_2-s}(x_2+\alpha (t_2-s)-y)e^{-(t_2-s)}\bigg) \notag \\
	&\hspace{11cm}(f(w^N_s(y))+ w^N_s(y))dyds\bigg| \notag \\
	 &\leq C (|x_1-x_2|^{4/9}+(t_2-t_1)^{2/9}).
	\end{align}
	For the third term on the right-hand side of~\eqref{eq:(G)}, since $|f(u)+u|\le 2+\alpha$ for $u\in [0,1]$,
	\begin{align} \label{eq:3rdtermholder}
	\left|\int_{t_1}^{t_2} \int_{\R} G_{t_2-s}(x_2+\alpha (t_2-s)-y)e^{-(t_2-s)}
	(f(w^N_s(y))+ w^N_s(y))dyds \right|
	&\leq (2+\alpha)(t_2-t_1).
	\end{align}
	For the fourth term on the right-hand side of~\eqref{eq:(G)}, recall that $w^N_s(y)\in [0,1]$ $\forall s\ge 0$ and $y\in \R$. Hence, by the Burkholder-Davis-Gundy inequality, for $k\in \N$, there exists a constant $C(k)<\infty$ such that
	\begin{align} \label{eq:GGWBDGineq}
	&\mathbb E\bigg[\bigg(\int_0^{t_1} \int_{\R} \bigg(G_{t_1-s}(x_1+\alpha (t_1-s)-y)e^{-(t_1-s)}-G_{t_2-s}(x_2+\alpha (t_2-s)-y)e^{-(t_2-s)}\bigg) \notag \\
	&\hspace{9cm}\sqrt{w^N_s(y)(1-w^N_s(y))}W(dy ds)\bigg)^{2k}\bigg]\notag \\
	&\leq C(k)
	\bigg(\int_0^{t_1} \int_{\R} \bigg(G_{t_1-s}(x_1+\alpha (t_1-s)-y)e^{-(t_1-s)} \notag \\
	&\hspace{5cm}-G_{t_2-s}(x_2+\alpha (t_2-s)-y)e^{-(t_2-s)}\bigg)^2 dy ds\bigg)^{k}.
	\end{align}
	Note that for $t>0$,
	\begin{equation} \label{eq:G2int}
	\int_{\R} G_t(y)^2 dy =(4\pi t)^{-1}(2\pi t)^{1/2}=2^{-3/2}\pi^{-1/2}t^{-1/2}.
	\end{equation}
	Since $(a+b)^2 \le 2a^2 +2b^2$ $\forall a,b\in \R$, and then by recalling the definition of $I_{t_1,t_2}(s,x)$ in~\eqref{eq:Idefn} and by~\eqref{eq:G2int}, we have
	\begin{align*}
	&\int_0^{t_1} \int_{\R} \bigg(G_{t_1-s}(x_1+\alpha (t_1-s)-y)e^{-(t_1-s)}
	-G_{t_2-s}(x_2+\alpha (t_2-s)-y)e^{-(t_2-s)}\bigg)^2 dy ds\\
	&\leq 2 \int_0^{t_1} e^{-2(t_1-s)}\int_{\R} (G_{t_1-s}(x_1+\alpha (t_1-s)-y)
	-G_{t_2-s}(x_2+\alpha (t_2-s)-y))^2 dy ds\\
	&\qquad + 2 \int_0^{t_1} e^{-2(t_1-s)}\int_{\R} G_{t_2-s}(x_2+\alpha (t_2-s)-y)^2(1-e^{-(t_2-t_1)})^2 dy ds\\
	&\le 2 \int_0^{t_1} e^{-2(t_1-s)}I_{t_1,t_2}(s,x_1-x_2+\alpha(t_1-t_2))ds \\
	&\qquad +2(t_2-t_1)^2 \int_0^\infty e^{-2s'}2^{-3/2}\pi^{-1/2}(s'+t_2-t_1)^{-1/2}ds'.
	\end{align*}
	By~\eqref{eq:Gdiff1} in Lemma~\ref{lem:Gdiffintbound} and~\eqref{eq:GGWBDGineq},
	it follows that there exists a constant $C'(k)$ such that
	\begin{align} \label{eq:4thtermholder}
	&\mathbb E\bigg[\bigg(\int_0^{t_1} \int_{\R} \bigg(G_{t_1-s}(x_1+\alpha (t_1-s)-y)e^{-(t_1-s)}-G_{t_2-s}(x_2+\alpha (t_2-s)-y)e^{-(t_2-s)}\bigg) \notag \\
	&\hspace{9cm}\sqrt{w^N_s(y)(1- w^N_s(y))}W(dy ds)\bigg)^{2k}\bigg] \notag \\
	&\leq C'(k)
	(|x_1-x_2|^{k/2}+(t_2-t_1)^{k/4}).
	\end{align}
	Finally, for the fifth term on the right-hand side of~\eqref{eq:(G)}, by the Burkholder-Davis-Gundy inequality, for $k\in \N$, there exists a constant $C(k)<\infty$ such that
	\begin{align} \label{eq:5thtermholder}
	&\mathbb E\bigg[\bigg(\int_{t_1}^{t_2} \int_{\R} G_{t_2-s}(x_2+\alpha (t_2-s)-y)e^{-(t_2-s)}
	\sqrt{w^N_s(y)(1- w^N_s(y))}W(dy ds)\bigg)^{2k}\bigg] \notag \\
	&\leq C(k)\bigg(\int_{t_1}^{t_2} \int_{\R} G_{t_2-s}(x_2+\alpha (t_2-s)-y)^2 e^{-2(t_2-s)}dy ds\bigg)^{k} \notag \\
	&\leq  C(k) \bigg(\int_{t_1}^{t_2} 2^{-3/2}\pi^{-1/2}(t_2-s)^{-1/2} ds\bigg)^{k} \notag \\
	&\leq (2\pi)^{-k/2}C(k) (t_2-t_1)^{k/2},
	\end{align}
	where the second inequality follows by~\eqref{eq:G2int}.
	Hence by combining~\eqref{eq:1sttermholder},~\eqref{eq:2ndtermholder},~\eqref{eq:3rdtermholder},~\eqref{eq:4thtermholder}, and~\eqref{eq:5thtermholder} with~\eqref{eq:(G)}, the result follows.
	\end{proof}
	
	\subsection{Proof of Proposition~\ref{prop:sigmaNrighttail}} \label{subsec:utail}
	
	Recall the definition of $w^N$ in~\eqref{eq:tildeudefn}.
	We start by obtaining a representation for $w^N_t(x)$ which is used in this section as well as in Section~\ref{subsec:expholder}.
	By~\eqref{eq:utildeformula} with $a=\|f'\|_\infty$, for $0 \leq t_0 \leq t$ and $x\in \R$,
	\begin{align} \label{eq:wNt0t}
	w^N_t(x)&=
	\int_{\R} G_{t-t_0}(x+\alpha (t-t_0) -y)e^{\|f'\|_\infty(t-t_0)} w^N_{t_0}(y)dy\\
	&\quad +\int_{t_0}^t \int_{\R} G_{t - s}(x+\alpha (t-s)-y)e^{\|f'\|_\infty(t-s)}
	(f(w^N_s(y)) - \|f'\|_\infty  w^N_s(y))dyds \notag\\
	&\quad + N^{-1/2}\int_{t_0}^t \int_{\R} G_{t-s}(x+\alpha (t-s)-y)e^{\|f'\|_\infty(t-s)}
	\sqrt{w^N_s(y)(1-w^N_s(y))}W(dy ds). \notag
	\end{align}
	For $0 \leq t_0 \leq t$ and $x\in \R$, let 
	\begin{align} \label{eq:Xtxdefn}
		&X_{t_0,t}(x)
		\\ &=N^{-1/2}\int_{t_0}^t \int_{\R} G_{t-s}(x+\alpha (t-s)-y)e^{\|f'\|_\infty(t-s)}
		\1{s\leq N\sigmaN{tail_right}}\sqrt{w^N_s(y)(1-w^N_s(y))}W(dy ds). \notag
	\end{align}
	We then remark that, if $0 \leq t_0 \leq t \leq N \sigmaN{tail_right}$, then by \eqref{eq:wNt0t} we have
	\begin{align} \label{eq:wNtXt0t}
	w^N_t(x)&=
	\int_{\R} G_{t-t_0}(x+\alpha (t-t_0) -y)e^{\|f'\|_\infty(t-t_0)} w^N_{t_0}(y)dy +  X_{t_0,t}(x) \\
	&\quad +\int_{t_0}^t \int_{\R} G_{t - s}(x+\alpha (t-s)-y)e^{\|f'\|_\infty(t-s)} \notag
	(f(w^N_s(y)) - \|f'\|_\infty  w^N_s(y))dyds.
	\end{align}
	We now continue with the proof of Proposition~\ref{prop:sigmaNrighttail}, which will require several lemmas; first, we prove an upper bound on $w^N_t$ for $t\le N \sigmaN{tail_right}$

	
	
	\begin{lemma} \label{lem:tildeuboundX}
		For $N$ sufficiently large, for $t\le N ( \sigmaN{interface}\wedge \sigmaN{tail_right} )$ and $x\in \R$, letting $t_0 = (\lfloor t \rfloor -1) \vee 0$,
		$$
		w^N_t(x)\le 
		\frac 1 2 N^{\cexp{tail_right_N}} e^{-(1-\cexp{tail_right_exp})x}
		+X_{t_0,t}(x).
		$$
	\end{lemma}
	\begin{proof}
	For $x \in \R$ and $0 \leq t_0 \leq t \leq N\sigmaN{tail_right}$, by \eqref{eq:wNtXt0t} and the fact that $f(w) \leq \|f'\|_\infty w$ for all $w \in [0,1]$ implies that 
	\begin{align} \label{eq:(M)}
	w^N_t(x)&\leq 
	\int_{\R} G_{t-t_0}(x+\alpha (t-t_0)-y)e^{\|f'\|_\infty(t-t_0)}w^N_{t_0} (y)dy +  X_{t_0,t}(x).
	\end{align}
	Let $x \geq 0$ and suppose that $t \geq 1$, so that $t_0 = \lfloor t \rfloor - 1 \geq 0$. 
	If $t_0 \leq N \sigmaN{interface}$, we have $r^N(w^N_{t_0}) \leq (1+ \cexp{interface})$, where we recall that $r^N(\cdot)$ is defined in \eqref{def_rN}.
	Hence, for $t_0 \leq N \sigmaN{interface}$, 
	\begin{align*}
	\int_\R w^N_{t_0}(y) e^{(1-\cexp{inter_integ})y} dy \leq N^{-1} C_X e^{(1-\cexp{inter_integ}) (1 + \cexp{interface} \log N)} = C_X N^{(1-\cexp{inter_integ})(1+\cexp{interface}) - 1}.
	\end{align*}
	This implies that for every $k \in \N_0$,
	\begin{align*}
	\int_k^{k+1} w^N_{t_0}(y) dy \leq C_X N^{(1-\cexp{inter_integ})(1+\cexp{interface}) - 1} e^{-(1-\cexp{inter_integ})k}.
	\end{align*}
	Let $(B_t)_{t\ge 0}$ denote a Brownian motion. 
	Then, using the above and the fact that $t - t_0 \in [1,2]$, if $t \leq N \sigmaN{interface}$ we have
	\begin{align} \label{eq:wNexpbdprel1}
	&\int_{\R} G_{t-t_0}(x+\alpha (t-t_0)-y) w^N_{t_0} (y) dy  
	\\& \hspace{.5 cm} \leq \P[x]{B_4 \leq 0} +\sum_{k \in \N_0} \int_k^{k+1} G_{t-t_0}(x+\alpha (t-t_0)-y) w^N_{t_0} (y)dy \notag
	\\ &\hspace{.5 cm} \leq \P[x]{B_4 \leq 0} +\sum_{k \in \N_0} \bigg( \sup_{y \in [k, k+1]} G_{t-t_0}(x+\alpha(t-t_0)-y)\bigg) \int_k^{k+1} w^N_{t_0} (y)dy \notag
	\\ &\hspace{.5 cm} \leq \P[x]{B_4 \leq 0} + C_X N^{(1-\cexp{inter_integ})(1+\cexp{interface}) - 1} \sum_{k \in \N_0} \bigg( \sup_{y \in [k, k+1]} G_{t-t_0}(x+\alpha(t-t_0)-y)\bigg)e^{-(1-\cexp{inter_integ})k}. \notag
	\end{align}
	Let $z = x + \alpha(t-t_0)$. Since $G_{t-t_0}$ is radially decreasing, 
	\begin{align*}
	\sum_{k \in \N_0} &\bigg( \sup_{y \in [k, k+1]} G_{t-t_0}(z-y)\bigg) e^{-(1-\cexp{inter_integ})k} 
	\\ & \leq \sum_{k \in \N_0 : k \leq z - 2} \int_{k}^{k+1} G_{t-t_0}(z-1-y)e^{-(1-\cexp{inter_integ})(y-1)}dy
	\\ &\quad +  \sum_{k \in \N_0 : k \geq z+2} \int_{k}^{k+1}G_{t-t_0}(z+1-y)  e^{-(1-\cexp{inter_integ})(y-1)}dy  
	\\ &\hspace{1 cm} + 4\left(\frac{1}{\sqrt{4\pi (t - t_0)}} e^{-(1-\cexp{inter_integ})(z - 2)} \right)
	\\ & \leq e^{2(1-\cexp{inter_integ})} \left( \E[z-1]{ e^{-(1-\cexp{inter_integ})B_{2(t-t_0)}}} + \E[z+1]{e^{-(1-\cexp{inter_integ})B_{2(t-t_0)}}} + 2 e^{-(1-\cexp{inter_integ})z} \right)
	\\ &\leq K_1 e^{-(1-\cexp{inter_integ})z}.
	\end{align*} 
	for some constant $K_1\geq 1$. Combined with \eqref{eq:wNexpbdprel1} and recalling that $z = x + \alpha(t - t_0)$ and $t - t_0 \in [1,2]$, we obtain that
	\begin{align*} 
	\int_{\R} G_{t-t_0}(x+\alpha (t-t_0)-y) w^N_{t_0} (y) dy  &\leq e^{-x^2/8}+ K_1 C_X N^{(1-\cexp{inter_integ})(1+\cexp{interface}) - 1} e^{-(1-\cexp{inter_integ})(x + \alpha(t-t_0))}
	\\ &\leq \frac 1 2 e^{-2\|f'\|_\infty} N^{\cexp{tail_right_N}} e^{-(1-\cexp{tail_right_exp})x},
	\end{align*}
	The second inequality holds for sufficiently $N$, because, by \eqref{bounds_cexpinterface} and \eqref{bounds_cinitsmall}, $(1-\cexp{inter_integ})(1+\cexp{interface}) - 1 < \cexp{interface} < \cexp{tail_right_N}$
	and $\cexp{inter_integ} < \cexp{tail_right_exp}$. 
	The desired inequality for $w^N_t(x)$ for $t \leq N(\sigmaN{interface} \wedge N\sigmaN{tail_right})$ and $x \geq 0$ now follows from the above and \eqref{eq:(M)}.
	If $x < 0$, the inequality is trivial.
		
	It remains to consider the case when $t \in [0,1)$, in which case $t_0 = 0$. Using Assumption~\ref{assumpt:v0}.\ref{v0:tail}, we have
	\begin{align*}
	\int_\R G_t(x+\alpha t-y) u_0^N(y) dy \leq N^{\cinit{small}} \E[x+\alpha t]{e^{-(1-\cinit{small})B_{2t}}} \leq N^{\cinit{small}} e^{-(1-\cinit{small})(x+\alpha t) + (1-\cinit{small})^2 t},
	\end{align*}
	and the result follows from \eqref{eq:(M)} and \eqref{bounds_cinitsmall}.
	\end{proof}

	In order to apply Kolmogorov's continuity criterion to the process $(t,x) \mapsto X_{t_0,t}(x)$, we prove the following moment bound on its increments.
	\begin{lemma} \label{lem:XtxdiffE}
	For $k\in \N$, there exists a constant $A^{(4)}_k<\infty$ such that for any $t_0 \geq 0$, $t_1,t_2 \in [t_0, t_0 + 2]$ and $x_1,x_2\in \R$ with $|x_1-x_2|\leq 1$,
	\begin{align*}
	\E{(X_{t_1}(x_1)-X_{t_2}(x_2))^{2k}}
	&\leq A^{(4)}_k N^{-k} (N^{-\cexp{tail_right_N}} e^{-(1-\cexp{tail_right_exp}) x_1})^k (|x_1-x_2|^{1/2}+|t_1-t_2|^{1/4})^k .
	\end{align*}
	\end{lemma}
	\begin{proof}
	Recall the definition of $X_{t_0,t}(x)$ in~\eqref{eq:Xtxdefn}.
	Let $t_0 \geq 0$, $t_1,t_2 \in [t_0, t_0 +2]$ and $x_1,x_2\in \R$ and assume without loss of generality that $t_1 \leq t_2$. Then
	\begin{align*}
	&X_{t_0,t_1}(x_1)-X_{t_0,t_2}(x_2)\\
	&=N^{-1/2}
	\bigg( \int_{t_0}^{t_1} \int_{\R} \bigg(G_{t_1-s}(x_1+\alpha (t_1-s)-y)e^{\|f'\|_{\infty})(t_1-s)}\\
	&\hspace{3.6cm} -G_{t_2-s}(x_2+\alpha (t_2-s)-y)e^{\|f'\|_\infty(t_2-s)}\bigg)\\
	&\hspace{6cm}
	\1{s\leq N\sigmaN{tail_right}}\sqrt{\sigma^2(w^N_s(y))}W(dy ds)\\
	&\quad -\int_{t_1}^{t_2} \int_{\R} G_{t_2-s}(x_2+\alpha (t_2-s)-y)e^{-\|f'\|_\infty(t_2-s)}
	\1{s\leq N\sigmaN{tail_right}}\sqrt{w^N_s(y)(1-w^N_s(y))}W(dy ds)\bigg).
	\end{align*}
	For $k\in \N$,
	by the Burkholder-Davis-Gundy inequality it follows that there exists a constant $C(k)<\infty$ such that
	\begin{align} \label{eq:dagger}
	&\E{(X_{t_0,t_1}(x_1)-X_{t_0,t_2}(x_2))^{2k}} \\
	&\leq N^{-k}2^{2k-1}  C(k) \notag\\
	&\quad  \bigg( \mathbb E \bigg[
	\bigg( \int_{t_0}^{t_1} \int_{\R} \bigg(G_{t_1-s}(x_1+\alpha (t_1-s)-y)e^{\|f'\|_\infty(t_1-s)} \notag\\
	&\hspace{2cm} -G_{t_2-s}(x_2+\alpha (t_2-s)-y)e^{\|f'\|_\infty(t_2-s)}\bigg)^2
	\1{s\leq N\sigmaN{tail_right}}w^N_s(y)(1-w^N_s(y)) dy ds \bigg)^{k}\bigg] \notag\\
	&\quad +\mathbb E\bigg[ \bigg(\int_{t_1}^{t_2} \int_{\R} G_{t_2-s}(x_2+\alpha (t_2-s)-y)^2 e^{2\|f'\|_\infty(t_2-s)}
	\1{s\leq N\sigmaN{tail_right}} w^N_s(y)(1-w^N_s(y))dy ds\bigg)^{k}
	\bigg]\bigg). \notag
	\end{align}
	Note that for $0<s<t$ and $x,y\in \R$, 
	\begin{align} \label{eq:Gshift}
	&G_{t-s}(x+\alpha (t-s)-y)
	\\&\quad=G_{t-s}(x-y+(\alpha-(1-\cexp{tail_right_exp}))(t-s))e^{-\frac 12 (1-\cexp{tail_right_exp})( x-y+(\alpha -(1-\cexp{tail_right_exp}))(t-s))}e^{-\frac 14 (1-\cexp{tail_right_exp})^2(t-s)} \notag
	\\&\quad=G_{t-s}(x-y+(\alpha-(1-\cexp{tail_right_exp}))(t-s))e^{-\frac 12 (1-\cexp{tail_right_exp}) x} e^{\frac 12 (1-\cexp{tail_right_exp}) y} \notag
	e^{(\frac 14 (1-\cexp{tail_right_exp})^2 -\frac 12 \alpha (1-\cexp{tail_right_exp}))(t-s)}.
	\end{align}
	For $s\in [0,t_1]$, let
	\begin{equation} \label{eq:z1sz2sdefn}
	 z_1(s):=x_1+(\alpha -(1-\cexp{tail_right_exp}))(t_1-s) \quad \text{ and }\quad z_2(s):=x_2+(\alpha -(1-\cexp{tail_right_exp}))(t_2-s).
	\end{equation}
	Then
	for the first term on the right-hand side of~\eqref{eq:dagger},
	by the definition of $\sigmaN{tail_right}$ in~\eqref{eq:sigmatailrightdefn} and then by~\eqref{eq:Gshift}, and $w^N_s(y)(1-w^N_s(y)) \leq w^N_s(y)$,
	 we have that
	\begin{align} \label{eq:star}
	&\int_{t_0}^{t_1} \int_{\R} \bigg(G_{t_1-s}(x_1+\alpha (t_1-s)-y)e^{\|f'\|_\infty)(t_1-s)} \notag\\
	&\hspace{2cm} -G_{t_2-s}(x_2+\alpha (t_2-s)-y)e^{\|f'\|_\infty)(t_2-s)}\bigg)^2
	\1{s\leq N\sigmaN{tail_right}} w^N_s(y)(1-w^N_s(y)) dy ds \notag \\
	&\leq \int_{t_0}^{t_1} \int_{\R} \bigg(G_{t_1-s}(x_1+\alpha (t_1-s)-y)e^{\|f'\|_\infty(t_1-s)} \notag\\
	&\hspace{3cm} -G_{t_2-s}(x_2+\alpha (t_2-s)-y)e^{\|f'\|_\infty(t_2-s)}\bigg)^2
	 N^{\cexp{tail_right_N}}  e^{-(1-\cexp{tail_right_exp})y} dy ds \notag \\
	&\leq \int_{t_0}^{t_1} \int_{\R} \bigg(G_{t_1-s}(z_1(s)-y)e^{-\frac 12 (1-\cexp{tail_right_exp}) x_1}e^{(\frac 14 (1-\cexp{tail_right_exp})^2 -\frac 12 (1-\cexp{tail_right_exp}) \alpha +\|f'\|_\infty)(t_1-s)} \notag \\
	&\hspace{3cm} -G_{t_2-s}(z_2(s)-y)e^{-\frac 12 (1-\cexp{tail_right_exp}) x_2}e^{(\frac 14 (1-\cexp{tail_right_exp})^2 -\frac 12 (1-\cexp{tail_right_exp}) \alpha + \|f'\|_\infty)(t_2-s)}\bigg)^2  \notag
	\\ &\hspace{3.5 cm} \times e^{(1-\cexp{tail_right_exp}) y} \cdot  N^{\cexp{tail_right_N}} e^{-(1-\cexp{tail_right_exp})y}
	dy ds  .
	\end{align}
	Let 
	\begin{equation*} \label{eq:rhodefn}
	\rho :=\tfrac 14 (1-\cexp{tail_right_exp})^2 +\tfrac 12 (1-\cexp{tail_right_exp}) \alpha + \|f'\|_\infty.
	\end{equation*}
	In order to bound the right-hand side of~\eqref{eq:star}, we write
	\begin{align*} 
	&\int_{t_0}^{t_1} \int_{\R} \bigg(G_{t_1-s}(z_1(s)-y)e^{-\frac 12 (1-\cexp{tail_right_exp}) x_1}e^{\rho(t_1-s)}
	 -G_{t_2-s}(z_2(s)-y)e^{-\frac 12(1-\cexp{tail_right_exp}) x_2}e^{\rho(t_2-s)}\bigg)^2 dy ds \notag \\
	&=\int_{t_0}^{t_1} e^{2\rho (t_1-s)}e^{-(1-\cexp{tail_right_exp}) x_1} \notag \\
	&\hspace{1cm}\int_{\R} \left(G_{t_1-s}(z_1(s)-y) -G_{t_2-s}(z_2(s)-y)e^{-\frac 12 (1-\cexp{tail_right_exp})(x_2-x_1)}e^{\rho (t_2-t_1)}\right)^2 dy ds \notag \\
	&\leq e^{-(1-\cexp{tail_right_exp}) x_1} \int_{t_0}^{t_1} e^{2\rho (t_1-s)} 
	\int_{\R} \bigg(2(G_{t_1-s}(z_1(s)-y) -G_{t_2-s}(z_2(s)-y))^2 \notag \\
	&\hspace{3cm}+2G_{t_2-s}(z_2(s)-y)^2 \left(1-e^{-\frac 12 (1-\cexp{tail_right_exp}) (x_2-x_1)}e^{\rho(t_2-t_1)}\right)^2\bigg) dy ds,
	\end{align*}
	where the last line follows since $(a+b)^2\le 2a^2 +2b^2$ $\forall a,b\in \R$.
	Hence using~\eqref{eq:Idefn} and~\eqref{eq:G2int}, if $t_1,t_2 \in [t_0,t_0+2]$ and $|x_1-x_2|\leq 1$,
	\begin{align} \label{eq:Gdiffz1z2}
	& \int_{t_0}^{t_1}\int_{\R} \left(G_{t_1-s}(z_1(s)-y)e^{-\frac 12 (1-\cexp{tail_right_exp}) x_1} e^{\rho(t_1-s)}-
	G_{t_2-s}(z_2(s)-y)e^{-\frac 12 (1-\cexp{tail_right_exp}) x_2} e^{\rho(t_2-s)}\right)^2 dy ds \notag \\
	&\leq 2 e^{-(1-\cexp{tail_right_exp})x_1}\int_{t_0}^{t_1}e^{2\rho(t_1-s)}I_{t_1,t_2}(s,z_1(s)-z_2(s))ds \notag \\
	&\quad +2e^{-(1-\cexp{tail_right_exp})x_1}\int_{t_0}^{t_1} e^{2\rho(t_1-s)} 2^{-3/2}\pi^{-1/2}
	(t_2-s)^{-1/2} \left(1-e^{-\frac 12 (1-\cexp{tail_right_exp}) (x_2-x_1)}e^{\rho(t_2-t_1)}\right)^2\, ds \notag \\
	&\leq K_1 e^{-(1-\cexp{tail_right_exp}) x_1}((t_2-t_1)^{1/4}+|x_1-x_2|^{1/2})
	\end{align}
	for some constant $K_1<\infty$. To obtain the second inequality, we bound the term on the second line using~\eqref{eq:Gdiff1} in Lemma~\ref{lem:Gdiffintbound} 
	combined with the definition of $z_1(s)$ and $z_2(s)$ in~\eqref{eq:z1sz2sdefn}, and use the fact that
	$\int_0^2 s^{-1/2}<\infty$. Since $t_1 - t_0 \leq 2$, the contribution of the term $e^{2\rho(t_1-s)}$ can be absorbed into $K_1$.
	For the term in the third line, the $s$-dependent coefficient in the integral is handled similarly, and an elementary calculation shows that 
	\[ \left(1-e^{-\frac 12 (1-\cexp{tail_right_exp}) (x_2-x_1)}e^{\rho(t_2-t_1)}\right) 
	\leq \tfrac 12 (1-\cexp{tail_right_exp})e^{\frac 12 (1-\cexp{tail_right_exp})} |x_2-x_1| + e^{\rho(t_2-t_1) + \frac 12(1-\cexp{tail_right_exp})} |t_2 - t_1|,\]
	and $t_2 - t_1 \in [0,2]$ again implies the desired bound.
	Hence, for $t_1,t_2 \in [t_0, t_0 + 2]$ and $|x_1-x_2|\leq 1$, by~\eqref{eq:star} and \eqref{eq:Gdiffz1z2},
	\begin{align} \label{eq:Gdiffbeta}
	&\int_0^{t_1}\int_{\R} \bigg(G_{t_1-s}(x_1+\alpha(t_1-s)-y)e^{\|f'\|_\infty(t_1-s)} \notag \\
	&\hspace{1cm}-
	G_{t_2-s}(x_2+\alpha(t_2-s)-y) e^{\|f'\|_\infty(t_2-s)}\bigg)^2 \1{s\leq N\sigmaN{tail_right}} \sigma^2(w^N_s(y)) dy ds \notag \\
	&\leq  K_1 e^{-(1-\cexp{tail_right_exp}) x_1}((t_2-t_1)^{1/4}+|x_1-x_2|^{1/2}).	\end{align}
	For the second term on the right-hand side of~\eqref{eq:dagger}, by the definition of $\sigmaN{tail_right}$ in~\eqref{eq:sigmatailrightdefn}, 
	and then by~\eqref{eq:Gshift} and the definition of $\rho$ in~\eqref{eq:rhodefn}, we have that
	\begin{align} \label{eq:G2square}
	& \int_{t_1}^{t_2} \int_{\R} G_{t_2-s}(x_2+\alpha (t_2-s)-y)^2 e^{2\|f'\|_\infty(t_2-s)} \notag
	\1{s\leq N\sigmaN{tail_right}} w^N_s(y)(1-w^N_s(y)) dy ds \\ \notag
	&\leq N^{\cexp{tail_right_N}} \int_{t_1}^{t_2}e^{2\|f'\|_\infty(t_2-s)} \int_{\R} 
	G_{t_2-s}(x_2+\alpha(t_2-s)-y)^2 e^{-(1-\cexp{tail_right_exp}) y} dy ds\\ \notag
	&\leq N^{\cexp{tail_right_N}} \int_{t_1}^{t_2} \int_{\R} G_{t_2-s}(x_2-y+(\alpha-(1-\cexp{tail_right_exp})(t_2-s))^2 \\ \notag
	&\hspace{3 cm} \times e^{-(1-\cexp{tail_right_exp}) x_2}e^{(1-\cexp{tail_right_exp}) y}e^{2\rho (t_2-s)} e^{-(1-\cexp{tail_right_exp}) y}dy ds\\ 
	&\leq e^{(2\rho \vee 0)(t_2-t_1)} e^{-(1-\cexp{tail_right_exp}) x_2} N^{\cexp{tail_right_N}}  \int_{t_1}^{t_2}  \int_{\R} G_{t_2-s}(y)^2 dy ds.
	\end{align}
	Hence by~\eqref{eq:G2int},
	\begin{align*}
	&\int_{t_1}^{t_2} \int_{\R} G_{t_2-s}(x_2+\alpha (t_2-s)-y)^2 e^{2\|f'\|_\infty(t_2-s)}
	\1{s\leq N\sigmaN{tail_right}} \sigma^2(w^N_s(y)) dy ds \\
	&\leq e^{(2\rho \vee 0)(t_2-t_1)} e^{-(1-\cexp{tail_right_exp}) x_2} N^{\cexp{tail_right_N}}  \int_{t_1}^{t_2} 2^{-3/2} \pi^{-1/2}(t_2-s)^{-1/2}ds\\
	&\leq e^{(2\rho \vee 0)(t_2-t_1)} e^{-(1-\cexp{tail_right_exp}) x_2} N^{\cexp{tail_right_N}} (t_2-t_1)^{1/2}.
	\end{align*}
	Therefore, by~\eqref{eq:dagger} and~\eqref{eq:Gdiffbeta}, for $t_0 \geq 0$ with $t_1,t_2 \in [t_0,t_0+2]$ and $x_1,x_2\in \R$ with $|x_1-x_2|\leq 1$, for $k\in \N$, 
	\begin{align*}
	&\E{(X_{t_0,t_1}(x_1)-X_{t_0,t_2}(x_2))^{2k}}\\
	&\leq N^{-k}2^{2k-1}C(k)
	\left( \left(K_1 e^{-(1-\cexp{tail_right_exp}) x_1}\right)^k(|t_2-t_1|^{1/4}+|x_1-x_2|^{1/2})^k \right.\\
	&\hspace{7cm}\left.+ \left( N^{\cexp{tail_right_N}} e^{4\rho \vee 0} e^{-(1-\cexp{tail_right_exp}) x_2}\right)^k  |t_1- t_2|^{k/2}\right),
	\end{align*}
	which completes the proof.
	\end{proof}
	
	We may now use Lemma~\ref{lem:XtxdiffE}  to obtain a modulus of continuity for $X_{t_0,t}(x)$.
	\begin{lemma} \label{lem:Xtholder}
	For each $k \in \N$ there exists $A_k' < \infty$ such that the following holds: for any $x \in \R$ and $t_0 \geq 0$, there exists a random variable $Y^N_{t_0,x,k}$ such that for all $x_1,x_2\in [x,x+2]$ and $t_1,t_2\in [t_0,t_0+2]$,
	\begin{align} \label{eq:(F1)}
	|X_{t_0,t_1}(x_1)-X_{t_0, t_2}(x_2)| \leq Y^N_{t_0,x,k}(|x_1 -x_2|^{\frac 1 4} + |t_1 - t_2|^{\frac 1 8})
	\end{align}
	and
	\begin{equation} \label{eq:(E1)}
	\E{(Y^N_{t_0,x,k})^{2k}}\leq A_k' N^{-k}( N^{\cexp{tail_right_N}} e^{-(1-\cexp{tail_right_exp}) x})^k.
	\end{equation}
	\end{lemma}
	\begin{proof}
	Similarly to the proof of Proposition~\ref{prop:tildeudiffunif}, this result follows immediately from Lemma~\ref{lem:XtxdiffE} 
	and an application of Kolmogorov's continuity criterion (as stated in Corollary 1.2 in~\cite{walsh_introduction_1986}).
	\end{proof}

	We can now use Lemma~\ref{lem:Xtholder} to prove a uniform upper bound on $X_{t_0,t}(x)$ that holds with high probability.
	\begin{lemma} \label{lem:Xtxunifbound}
	For $a>0$, for $N$ sufficiently large, for any $t_0 \geq 0$,
	\begin{align*}
	\P{\exists t\in [t_0,t_0+2], \,x\in [0,N^3] :X_{t_0,t}(x)\geq N^{-1/2+a}(N^{\cexp{tail_right_N}}e^{-(1-\cexp{tail_right_exp})x})^{1/2}}
	&\leq N^{-5}.
	\end{align*}
	\end{lemma}
	\begin{proof}
	Fix $a>0$ and let $\rho = 2^{\frac 1 4} + 2^{\frac 1 8}$.
	Recalling the random variable $Y^N_{t_0,x,k}$ from Lemma~\ref{lem:Xtholder}, 
	and because $X_{t_0,t_0}(x) = 0$ for any $x$, it follows from \eqref{eq:(F1)} that
	\begin{equation} \label{eq:Xt0tunifholder}
	\sup_{y \in [x,x+2], t \in [t_0,t_0+2]} X_{t_0,t}(y) \leq \rho \cdot Y^N_{t_0,x,k}.
	\end{equation}
	By Lemma~\ref{lem:Xtholder} and applying Markov's inequality to \eqref{eq:(E1)}, we obtain
	\begin{align*} 
	&\P{Y_{t_0,x,k}\geq \tfrac{1}{\rho} N^{-1/2+a/2}(N^{\cexp{tail_right_N}}e^{-(1-\cexp{tail_right_exp})(x+2)})^{1/2}} \notag
	\\ &\hspace{1.5 cm} \leq \rho^{-2k} N^{k-ka}(N^{\cexp{tail_right_N}}e^{-(1-\cexp{tail_right_exp})(x+2)})^{-k}\E{(Y_{t_0,x,k})^{2k}} \notag
	\\ & \hspace{1.5 cm}\leq \rho^{-2k} A_k' e^{(1-\cexp{tail_right_exp})2k} N^{-ak}.
	\end{align*}
	Hence, choosing $k \in \N$ sufficiently large so that that $ka > 8$, by \eqref{eq:Xt0tunifholder} and a union bound,
	\begin{align*}
	&\P{\exists \,x\in \Z\cap [0,N^3] :\sup_{y\in [x,x+1],t\in [t_0,t_0+2]}X_{t_0,t}(y)\geq N^{-1/2+a}(N^{\cexp{tail_right_N}}e^{-(1-\cexp{tail_right_exp})x})^{1/2}}\\
	&\leq (N^3+1) \rho^{-2k} A_k' e^{(1-\cexp{tail_right_exp})2k}  N^{-ka}\\
	&\leq N^{-5}
	\end{align*}
	for $N$ sufficiently large. The result follows.
	\end{proof}
	
	We can now combine Lemma~\ref{lem:tildeuboundX} and Lemma~\ref{lem:Xtxunifbound} to prove the following upper bound on $w^N_t$.
	\begin{cor} \label{cor:tildeuunifbound}
	For $N$ sufficiently large,
	\begin{align*}
	&\P{\exists t\in [0,(N ( \sigmaN{interface} \wedge \sigmaN{tail_right}))\wedge N^2], \,x\in [0, (1-\cexp{tail_right_exp})^{-1}\log N ] : w^N_t(x)\geq N^{\cexp{tail_right_N}}e^{-(1-\cexp{tail_right_exp})x}}
	\\ &\hspace{1 cm}< N^{-2}.
	\end{align*}
	\end{cor}
	\begin{proof}
	Take $a\in (0,\cexp{tail_right_N} /2)$, and suppose $N$ is sufficiently large that Lemma~\ref{lem:tildeuboundX} and Lemma~\ref{lem:Xtxunifbound} hold, and that $N^{a-\cexp{tail_right_N}/2}<\frac 12$.
	Suppose that for all $t_0 \in \Z \cap [0,N^2]$,
	\[X_{t_0,t}(x) <  N^{- 1/2 + a}(N^{\cexp{tail_right_N}}(e^{-(1-\cexp{tail_right_exp})x})^{1/2}\]
	for all $t\in [t_0,t_0+2]$ and $x\in [0,N^3]$.
	Then by Lemma~\ref{lem:tildeuboundX}, for $t\in [0,(N ( \sigmaN{interface} \wedge \sigmaN{tail_right}))\wedge N^2]$ and $x\in [0,(1-\cexp{tail_right_exp})^{-1} \log N]$,
	\begin{align*}
	w^N_t(x)
	&\leq < 12 N^{\cexp{tail_right_N}}e^{-(1-\cexp{tail_right_exp})x}  + N^{-1/2+a + \cexp{tail_right_N}/2}e^{-(1-\cexp{tail_right_exp})x/2} <  N^{\cexp{tail_right_N}}e^{-(1-\cexp{tail_right_exp})x} 
	\end{align*}
	for sufficiently large $N$. The second inequality holds because we have assumed that $N^{a-\cexp{tail_right_N}/2}<\frac 12$
	which implies $ N^{-1/2+a + \cexp{tail_right_N}/2} < \frac 12 N^{-1/2 + \cexp{tail_right_N}}$. 
	Then, we use the fact that $N^{-1/2} \leq e^{-(1-\cexp{tail_right_exp})x/2}$, which follows because $x \leq (1-\cexp{tail_right_exp})^{-1} \log N$.
	The result follows by Lemma~\ref{lem:Xtxunifbound} and a union bound over $t_0 \in \Z \cap [0,N^2]$.
	\end{proof}

	\begin{proof}[Proof of Proposition~\ref{prop:sigmaNrighttail}]
	By the definition of $\sigmaN{interface}$ in~\eqref{eq:sigmainterface}, and since $1+ \cexp{interface} < (1- \cexp{tail_right_exp})^{-1}$ from \eqref{bounds_cexpinterface}, for $t \leq N \sigmaN{interface}$ we have
	\[ w^N_t(x) = 0 \,\, \forall \, x \geq (1-\cexp{tail_right_exp})^{-1}.\]
	By the definition of $\sigmaN{tail_right}$ in~\eqref{eq:sigmatailrightdefn}, and taking $N>T$,
	the result follows directly from Corollary~\ref{cor:tildeuunifbound}. 
	\end{proof}

	\subsection{Proof of Proposition~\ref{prop:sigmatailleft}} \label{subsec:1-utail}
	
	For $t\ge 0$ and $x\in \R$,
	let $w^N_t(x)=1-u^N_t(x)$.
	For $w \in [0,1]$, let $g(w) = -f(1-w)$.  
	Then $\tilde{u}^N_t$ is a mild solution of the SPDE
	$$
	d \tilde{u}^N_t = (\partial_{xx} \tilde{u}^N_t+g(\tilde{u}^N_t))dt-\sqrt{\frac 1 N \tilde{u}^N_t(1- \tilde{u}^N_t)} dW(t),
	$$
	where $(W(t),t\ge 0)$ is a cylindrical Wiener process.
	By~\eqref{eq:utildeformula}, for $a\in \R$, for $t\ge 0$ and $x\in \R$,
	\begin{align*} 
	&\tilde{u}_t^N(-x+\alpha t)\\
	&=1- w^N_t(-x)\\
	&=
	\int_{\R} G_t(-x+\alpha t-y)e^{at}(1-u^N_0(y))dy \notag \\
	&\qquad +\int_0^t \int_{\R} G_{t-s}(-x+\alpha (t-s)-y)e^{a(t-s)}
	(g(\tilde{u}^N_s(y+\alpha s))-a \tilde{u}^N_s(y+\alpha s))dyds \notag \\
	&\qquad - N^{-1/2}\int_0^t \int_{\R} G_{t-s}(-x+\alpha (t-s)-y)e^{a(t-s)}
	\sqrt{\tilde{u}_s(y+\alpha s)(1-\tilde{u}_s(y+\alpha s))}W(dy ds).
	\end{align*}
	For $t\ge 0$ and $x\in \R$, let
	$\tilde w^N_t(x)=\tilde{u}^N_t(-x+\alpha t)$. Then letting $y'=-y$ inside the integrals, we can write
	\begin{align} \label{eq:wNformula}
	\tilde w^N_t(x)&= 
	\int_{\R} G_t(x-\alpha t-y')e^{a t}(1-u^N_0(-y'))dy'  \\
	&\qquad +\int_0^t \int_{\R} G_{t-s}(x-\alpha (t-s)-y')e^{a(t-s)}
	(g(\tilde w^N_s(y'))-a\tilde w^N_s(y'))dy' ds \notag \\
	&\qquad - N^{-1/2}\int_0^t \int_{\R} G_{t-s}(x-\alpha (t-s)-y')e^{a(t-s)}
	\sqrt{\tilde{w}^N_s(y')(1-\tilde{w}^N_s(y'))}W(dy' ds).\notag 
	\end{align} 
	Note that 
	\begin{equation} \label{eq:tildewdefn}
	\tilde w^N_t(x)=1-v^N_{t/N}(-x)=1-w^N_t(-x)
	\end{equation}
	by~\eqref{def_vN} and~\eqref{eq:tildeudefn}.
	For $t\geq 0$ and $x\in \R$, let 
	$$
	X'_t(x)=N^{-1/2}\int_0^t \int_{\R} G_{t-s}(x-\alpha (t-s)-y)e^{\frac{1}{2}g'(0)(t-s)}
	 \sqrt{\tilde{w}^N_s(y)(1- \tilde{w}^N_s(y))}W(dy ds).
	$$
	
	\begin{lemma} \label{lem:tildewboundY}
	For $N$ sufficiently large,
	\begin{align*}
	&\P{ \exists t \in [0, N( \sigmaN{dist} \wedge \sigmaN{eta} \wedge \sigmaN{tail_left}) \wedge N^2], x \in \R : \tilde{w}^N_t(x) > \frac 1 2 \Cst{tail_left}(e^{-\cexp{tail_left}x} + N^{-\cexp{tail_left}}) - X'_t(x)}
	\\ &\hspace{1 cm} < N^{-2}.
	\end{align*}
	\end{lemma}

	\begin{proof}
	We observe that $g'(0) = f'(1) < 0$. We introduce constants
	\begin{equation} \label{eq:ell1ell2defn}
	\ell_1 := \|g'\|_\infty + \frac 1 2 |g'(0)| \quad \text{ and } \quad \ell_2 := \left( \frac 1 2 |g'(0)| - \alpha \cexp{tail_left} - \cexp{tail_left}^2 \right)^{-1},
	\end{equation}	
	and remark that $\ell_2 >0$ by \eqref{eq:cleftfprime}. 
	Furthermore, from these definitions, \eqref{eq:CtailleftKexp} and the fact that $g'(0) = f'(1)$, we obtain that
	\begin{equation} \label{eq:CtailleftKexpapp}
	\forall \, x \in \R, \,\, \sup_{y \in [-K,K]} 1 - m(x+y) \leq \frac{1}{8\ell_1 (\ell_2 \vee 1)} \Cst{tail_left} e^{\cexp{tail_left}x}.
	\end{equation}
	Next, recall from \eqref{eq:sigmalefttaildefn} that
	\begin{equation*}
	\sigmaN{tail_left} = \inf \left\{ t \geq 0 : \exists x \in \R \text{ such that } 1 -v^N_t(x) \geq \Cst{tail_left}(e^{\cexp{tail_left}x} + N^{-\cexp{tail_left}}) \right\}.
	\end{equation*}
	Let $u_- \in (0,1)$ be sufficiently small so that $g(u) - \tfrac 12 g'(0) u \leq 0$ for all $u \in [u,u_-]$. Then choose $R_->0$ and let $N$ be sufficiently large so that
	\[\Cst{tail_left} e^{-\cexp{tail_left}y} < \frac 1 2 u_- \,\, \forall y \geq R_- \quad \text{ and } \quad \Cst{tail_left} N^{-\cexp{tail_left}} < \frac 1 2 u_-.\]
	It is then immediate that if $s \leq N \sigmaN{tail_left}$ and $x \geq R_-$, then $\tilde{w}^N_s(y) < u_-$, and hence
	\begin{equation} \label{eq:gtildewbd1}
	g(\tilde{w}^N_s(y)) - \frac 1 2 g'(0)\tilde{w}^N_s(y) \leq 0 \,\,\quad \forall \, s  \leq \sigmaN{tail_left}, \, x \leq R_-.
	\end{equation}
	by our choice of $u_-$. Now let 
	\[ \epsilon = \frac{\Cst{tail_left}}{8\ell_1(\ell_2 \vee 1)} e^{-\cexp{tail_left}R_-}\]
	and define
	\begin{equation*}
	E_N := \left\{ \exists x \in [-R_-,0], t \in [0,N(\sigmaN{dist} \wedge \sigmaN{eta}) \wedge N^2] : |s(w^N_t)(x)| > \epsilon \right\}.
	\end{equation*}
	Then by Lemma~\ref{lem:scontrol}, for sufficiently large $N$,
	\begin{equation}
	\P{E_N} < N^{-2}.
	\end{equation}
	By our choice of $\epsilon$ and \eqref{eq:CtailleftKexpapp}, on $E_N$, for $t \in [0,N(\sigmaN{dist} \wedge \sigmaN{eta}) \wedge N^2]$, we have
	\begin{equation*}
	\tilde{w}^N_t(x) \leq \frac{\Cst{tail_left}}{4 \ell_1 (\ell_2 \vee 1)} e^{-\cexp{tail_left}x} \,\, \forall x \leq R_-.
	\end{equation*}
	(The case $x<0$ in the above holds trivially because $\Cst{tail_left} \geq 4\ell_1(\ell_2 \vee 1)$, see the discussion above \eqref{eq:CtailleftKexp}.)
	Hence, by our choice of $\ell_1$ in \eqref{eq:ell1ell2defn}, we conclude from the above that on $E_N$,
	\begin{equation} \label{eq:gtildewbd2}
	g(\tilde{w}^N_t(x)) - \frac 1 2 g'(0) \tilde{w}^N_t(x) \leq \frac{\Cst{tail_left}}{ 4(\ell_2 \vee 1)} e^{-\cexp{tail_left}x} \quad \,\, \forall t \in [0,N(\sigmaN{dist} \wedge \sigmaN{eta}) \wedge N^2],  x \leq R_-.
	\end{equation}
	We recall Assumption~\ref{assumpt:v0}.\ref{v0:tail}, which along with the relations \eqref{eq:cleftfprime} and \eqref{eq:lefttailinitconst} implies that 
	\[ 1-u^N_0(y) \leq \frac 1 4 \Cst{tail_left}(e^{-\cexp{tail_left}y} + N^{-\cexp{tail_left}}) \quad \forall \, y \in \R.\]
	(For $y \leq 0$, we simply use the bound $1-u^N_0(y) \leq 1$ and $\Cst{tail_left} \geq 4$ from \eqref{eq:lefttailinitconst}.)
	Hence, by~\eqref{eq:wNformula} with $a=\tfrac 1 2 g'(0) < 0$, using \eqref{eq:gtildewbd1} and \eqref{eq:gtildewbd2},
	for $t\leq N(\sigmaN{dist} \wedge \sigmaN{eta} \wedge \sigmaN{tail_left})$ and $x \in \R$,
	\begin{align} \label{eq:(M)2}
	\tilde w^N_t(x)&\leq 
	\int_{\R} G_t(x-\alpha t-y)e^{\frac 1 2 g'(0)t}\frac 1 4 \Cst{tail_left}(e^{-\cexp{tail_left}y}+N^{-\cexp{tail_left}}) dy   \\
	&\qquad +\int_0^t \int_{\R} G_{t-s}(x-\alpha (t-s)-y)e^{\frac 1 2 g'(0)(t-s)}
	\frac{\Cst{tail_left}}{4(\ell_2 \vee 1)} e^{-\cexp{tail_left}y} dyds  - X'_t(x).\notag
	\end{align}
	If $(B_t)_{t \geq 0}$ denotes a Brownian motion, for any $s \geq 0$,
	\begin{equation*}
	\int_{\R} G_s(x-\alpha s-y) e^{-\cexp{tail_left}y} dy = \E[x- \alpha s]{e^{-\cexp{tail_left} B_{2s}}} = e^{-\cexp{tail_left}(x - \alpha s)} e^{\cexp{tail_left}^2 s}.
	\end{equation*}
	Hence, using \eqref{eq:cleftfprime} and $g'(0) = f'(1)$ and the definition of $\ell_2>0 $ in \eqref{eq:ell1ell2defn}, we obtain from \eqref{eq:(M)2} that 
	\begin{align*} 
	&\tilde w^N_t(x) \leq 
	\Cst{tail_left}( e^{-\cexp{tail_left}x} + N^{-\cexp{tail_left}}) \bigg( \frac 1 4 e^{-\ell_2^{-1} t} +  \frac{1}{4(\ell_2 \vee 1)} \int_0^t   e^{-\ell_2^{-1} s} ds \bigg)  - X'_t(x).\notag
	\end{align*}
	for all $t\leq N(\sigmaN{dist} \wedge \sigmaN{eta} \wedge \sigmaN{tail_left})$ and $x \in \R$. 
	The result follows.
	\end{proof}
	We now prove moment bounds on $X'_t(x)$ and its increments. 
	\begin{lemma} \label{lem:X'txmoments}
	For $k\in \N$, there exist constants $A^{(6)}_k<\infty$ and $A^{(7)}_k<\infty$ such that for $0\leq t_1\leq t_2$ and $x_1,x_2\in \R$ with $t_2-t_1\leq 1$ and $|x_1-x_2|\leq 1$,
	\begin{align*}
	\E{(X'_{t_1}(x_1)-X'_{t_2}(x_2))^{2k}}
	&\leq A^{(6)}_k N^{-k} (|x_1-x_2|^{1/2}+(t_2-t_1)^{1/4})^k,
	\end{align*}
	and for $t\ge 0$ and $x\in \R$,
	\begin{align*}
	\E{(X'_t(x))^{2k}}
	&\leq A^{(7)}_k N^{-k}.
	\end{align*}
	\end{lemma}
	 \begin{proof}
	Since $\tfrac 12 g'(0) < 0$, the first inequality (i.e. the increment bound) follows from the exact same argument used to prove 
	\eqref{eq:4thtermholder} and \eqref{eq:5thtermholder} in the proof of Lemma~\ref{lem:kolmest}. 
	To prove the second inequality, from the Burkholder-Davis-Gundy inequality and the fact that $\tilde{w}^N_s$ takes values in $[0,1$], for some $C(k) \geq 1$, 
	\begin{align*}
	\E{(X_t'(x))^{2k}} \leq C(k) N^{-k} \left (\int_0^t \int_\R e^{g'(0)(t-s)} G_{t-s}(x-\alpha(t-s) - y)^2 dy ds \right)^k.
	\end{align*}
	The desired inequality now follows using \eqref{eq:G2int} and the fact that $\int_0^\infty s^{-1/2} e^{\frac 1 2 g'(0)s} ds < \infty$ as in the proof of Lemma~\ref{lem:kolmest}.
%
%
	 \end{proof}
	 
	 Next we can state a result that corresponds to Corollary~\ref{cor:tildeuunifbound}.
	 \begin{cor} \label{cor:wtildeunifmid}
	For $N$ sufficiently large,
	\begin{align*}
	\P{\exists t\in [0,N(\sigmaN{dist} \wedge \sigmaN{eta} \wedge \sigmaN{tail_left})\wedge N^2], \,x\in [0,N^3] :\tilde w^N_t(x) \geq \Cst{tail_left}(e^{-\cexp{tail_left} x}+N^{-\cexp{tail_left}})}
	&<  2 N^{-2}.
	\end{align*}
	\end{cor}
	\begin{proof}
	Let $k \in \N$ be sufficiently large so that $k(1-2\cexp{tail_left}) > 8$. Then for $x\in \R$ and $t \geq 0$,
	by Markov's inequality and then by Lemma~\ref{lem:X'txmoments},
	\begin{align} \label{eq:X'txMarkov}
	\P{|X_{t}'(x)|\geq \tfrac 12 \Cst{tail_left} N^{-\cexp{tail_left}}} \leq  2^{k} N^{2\cexp{tail_left}} \E{(X_{t_0,t}(x))^{2k}}  &\leq 2^{k}  A^{(7)}_k N^{-k(1 - 2\cexp{tail_left})}.
	\end{align}
	Similarly to the proofs of Proposition~\ref{prop:tildeudiffunif} and Lemma~\ref{lem:Xtholder}, we may use the increment moment bounds from Lemma~\ref{lem:X'txmoments}
	to conclude that for any $k \in \N$, $t \geq 0$ and $x \in \R$, there exist random variables $\tilde{Y}^N_{t,x,k}$ such that for all $t_1,t_2 \in [t,t+2]$ and $x_1, x_2 \in [x,x+2]$,

	\begin{equation} \label{eq:Xt'inc}
	|X'_{t_1}(x_1)-X'_{t_2}(x_2)| \leq  \tilde{Y}^N_{t,x,k} \|(x_1,t_1) - (x_2,t_2)\|_2^{\cexp{holder}},
	\end{equation}
	and $\E{(\tilde{Y}^N_{t,x,k})^{2k}} \leq A_k'' N^{-k}$ for some constant $A_k'' \geq 1$.
	In particular, 
	\begin{align*}
	\P{\tilde{Y}^N_{t,x,k}\geq \tfrac 14 N^{-\cexp{tail_left}}} \leq 4^{k} A_k'' N^{2k\cexp{tail_left}} \E{(Y_{t_0,x,k})^{2k}}  
	&\leq 4^{k} A_k'' N^{-k(1 - 2\cexp{tail_left})}.
	\end{align*}
	By~\eqref{eq:X'txMarkov},~\eqref{eq:Xt'inc}, the above, and a union bound, 
	\begin{align*}
	\P{\sup_{y\in [x,x+2],s\in [t,t+2]}X_{s}'(y)\geq \tfrac 12   \Cst{tail_left} N^{-\cexp{tail_left}}}
	&\leq 4^{k}(A^{(7)} +  A_k'') N^{-k(1 - 2\cexp{tail_left})}.
	\end{align*}
	Hence, by another union bound,
	\begin{align*}
	&\P{\exists \, t\in [0,N^2-1], x\in \Z\cap [0,N^3] :\sup_{y\in [x,x+1],s\in [t,t+2]}X_{s}'(y)\geq  \tfrac 12 \Cst{tail_left}N^{-\cexp{tail_left}} } \\
	&\leq N^2 (N^3+1) 4^{k}(A^{(7)} +  A_k'') N^{-k(1 - 2\cexp{tail_left})} \\
	&\leq N^{-2}
	\end{align*}
	for $N$ sufficiently large, where the last line follows because $k(1-2\cexp{tail_left}) > 8$. The now follows from the above and Lemma~\ref{lem:tildewboundY}.
	\end{proof}
	
	The following result handles the tail of $\tilde{w}_t^N(x)$ for $x \geq N^3$.
	\begin{lemma} \label{lem:wtildeuniftail}
	For $N$ sufficiently large,
	\begin{align*}
	\P{\exists x\geq N^3,\, t\in [0,N^2] : \tilde w^N_t(x)\geq N^{-\cexp{tail_left}}}
	&\leq  e^{-N}.
	\end{align*}
	\end{lemma}
	\begin{proof}
	Note that $g(w)\le \|g'\|_\infty w$ for $w\in [0,1]$ and recall that  $u^N_0(y)=1$ $\forall y\le -N$ by Assumption~\ref{assumpt:v0}.\ref{v0:compact_support}.
	Therefore, by~\eqref{eq:wNformula} with $a=\|f'\|_\infty$ and the same argument as for the proof of Lemma~\ref{lem:wNtailbd},
	for $x\in \R$ and $t\geq 0$ with $x-\alpha t\ge N$,
	$$
	\E{\tilde w^N_t(x)}
	\le e^{\|g'\|_\infty t} \int_{-\infty}^N G_t(x-\alpha t-y)dy
	\le e^{\|g'\|_\infty t} e^{-(x-\alpha t-N)^2/(4t)}.
	$$
	Hence for $N$ sufficiently large, for $t\in [0,N^2]$ and $x\ge N^3$,
	we have $x-\alpha t-N \ge x/2$ and so
	\begin{equation} \label{eq:wmomentlargex}
	\E{\tilde w^N_t(x)}
	\leq e^{\|g'\|_\infty t} e^{-x^2/(16t)}.
	\end{equation}
	Now take $\cexp{holder}>0$ as in Proposition~\ref{prop:tildeudiffunif}
	and let $a\in (0,\cexp{holder})$ and $k\in \N$ with $k>1/a$.
	Take $N$ sufficiently large that $N^3-1\ge 2N$ and
	\begin{equation} \label{eq:NlargeN3lem}
	 2^{(a-\cexp{holder})(N^3-1)}2^{\cexp{holder}/2}\le \tfrac 12 N^{-\cexp{tail_left}},
	\end{equation}
	and take $(U^N_{t,x,k})_{t\ge 0,x\in \Z}$ as in Proposition~\ref{prop:tildeudiffunif}.
	By a union bound,
	\begin{align} \label{eq:utailunionbd}
	&\P{\exists x\in \Z \cap[ N^3-1,\infty),\,t\in \Z\cap [0,N^2] \text{ s.t. }
	\tilde w^N_s(y)\geq N^{-\cexp{tail_left}}\text{ for some }s\in [t,t+1), y\in [x,x+1)} \notag \\
	&\leq \sum_{x\in \Z \cap[ N^3-1,\infty),t\in \Z\cap [0,N^2]}
	\bigg(
	\P{\{U^N_{N^{-1}t,-x,k}\leq 2^{ax}\}\cap \{\exists s\in [t,t+1),y\in [x,x+1):\tilde w^N_s(y)\geq N^{-\cexp{tail_left}}\}} \notag \\
	&\qquad \qquad \qquad \qquad \qquad \qquad  +\P{U^N_{N^{-1} t,-x,k}\geq 2^{a x}}
	\bigg).
	\end{align}
	Suppose for some $x\in \Z\cap [N^3-1,\infty)$ and $t\in \Z\cap [0,N^2]$, we have that $U^N_{N^{-1}t,-x,k}\leq 2^{ax}$ and
	$$
	\tilde w^N_{s'}(y')\leq \tfrac 12 N^{-\cexp{tail_left}}\,\,\, \forall s'\in 2^{-x}\Z\cap [t,t+1), \, y'\in 2^{-x}\Z\cap [x,x+1).
	$$
	Then
	for $s\in [t,t+1)$ and $y\in [x,x+1)$, 
	letting $s'=2^{-x}\lfloor 2^x s\rfloor$ and $y'=2^{-x}\lfloor 2^x y\rfloor$ we have that
	 $\tilde w^N_{s'}(y')\leq \frac 12 N^{-\cexp{tail_left}}$ and $|s-s'|\vee |y-y'|\leq 2^{-x}$.
	Hence by Proposition~\ref{prop:tildeudiffunif},
	and recalling how $\tilde{w}^N$ relates to $v^N$ by~\eqref{eq:tildewdefn},
	\begin{align*}
	\tilde w^N_s(y)&\leq \tfrac 12 N^{-\cexp{tail_left}}+2^{ax}\|(y,s)-(y',s')\|_2^{\cexp{holder}}
	\leq N^{-\cexp{tail_left}},
	\end{align*}  
	where the second inequality follows by our choice of $N$ in~\eqref{eq:NlargeN3lem} and since
	$x\ge N^3-1$.
	Therefore, for  $x\in \Z\cap [N^3-1,\infty)$ and $t\in \Z\cap [0,N^2]$,
	by a union bound and then by Markov's inequality,
	\begin{align*}
	&\P{\{U^N_{N^{-1} t,-x,k}\leq 2^{ax}\}\cap \{\exists s\in [t,t+1),y\in [x,x+1):\tilde w^N_s(y)\geq N^{-\cexp{tail_left}}\}}\\
	&\leq \sum_{y'\in 2^{-x}\Z\cap [x,x+1), \, s'\in 2^{-x}\Z \cap [t,t+1)}
	\P{\tilde w^N_{s'}(y')\geq \tfrac 12 N^{-\cexp{tail_left}}}\\
	&\leq 
	\sum_{y'\in 2^{-x}\Z\cap [x,x+1), \, s'\in 2^{-x}\Z \cap [t,t+1)}
	2N^{\cexp{tail_left}}\E{\tilde w^N_{s'}(y')}\\
	&\leq 
	\sum_{y'\in 2^{-x}\Z\cap [x,x+1), \, s'\in 2^{-x}\Z \cap [t,t+1)}
	2N^{\cexp{tail_left}}e^{(1+\alpha)(t+1)}e^{-x^2/(16(t+1))},
	\end{align*}
	where the last inequality follows
	by~\eqref{eq:wmomentlargex}.
	Hence by~\eqref{eq:utailunionbd} and Markov's inequality, and then using Proposition~\ref{prop:tildeudiffunif},
	\begin{align*}
	&\P{\exists x\geq N^3,\, t\in [0,N^2] : \tilde w^N_t(x)\geq N^{-\cexp{tail_left}}}\\
	&\leq \sum_{x\in \Z \cap[ N^3-1,\infty),t\in \Z\cap [0,N^2]}
	\bigg(
	2^{2x}\cdot 2 N^{\cexp{tail_left}}e^{\|g'\|_\infty(N^2+1)}e^{-x^2/(16(N^2+1))}+2^{-2a kx}\E{(U^N_{N^{-1}t,-x,k})^{2k}}
	\bigg)\\
	&\leq \sum_{x\in \Z \cap[ N^3-1,\infty),t\in \Z\cap [0,N^2]}
	\bigg(
	2 N^{\cexp{tail_left}}e^{\|g'\|_\infty(N^2+1)}e^{-(N^3-1)(16(N^2+1))^{-1}x+2x\log 2}+2^{-2akx}A_k
	\bigg)\\
	&\leq (N^2+1)
	\bigg(
	2 N^{\cexp{tail_left}}e^{\|g'\|_\infty(N^2+1)}\sum_{x=\lceil N^3-1\rceil}^\infty e^{(2\log 2-(N^3-1)(16(N^2+1))^{-1})x}+A_k \sum_{x=\lceil N^3-1 \rceil}^\infty 2^{-2x}
	\bigg)\\
	&\leq e^{-N}
	\end{align*}
	for $N$ sufficiently large, where the third inequality follows since we chose $k>1/a$.
	\end{proof}
	\begin{proof}[Proof of Proposition~\ref{prop:sigmatailleft}]
	Since $\Cst{tail_left} \geq 1$, the result follows directly from
	Corollary~\ref{cor:wtildeunifmid} and Lemma~\ref{lem:wtildeuniftail}.	\end{proof}

	\subsection{Proof of Proposition~\ref{prop:holderexp}} \label{subsec:expholder}
	Our starting point is the representation for $w^N_t(x)$ for $t \leq N \sigmaN{tail_right}$ given in \eqref{eq:wNtXt0t}, where we recall that for $0 \leq t_0 \leq t$ and $x \in \R$, $X_{t_0,t}(x)$ is defined in \eqref{eq:Xtxdefn}. 
	There are several constants appearing only in this subsection which we denote by $K_A$, $K_B$, etc. 
	\begin{lemma} \label{lem:holderexpbd}
	There exists a constant $K_A \geq 1$ such that the following holds: for any $t_0 \geq 0$ and $x \in \R$, a.s. on $\{t_0 \leq N \sigmaN{tail_right}\}$,
	\begin{align*}
	&|w^N_{t_1}(x_1) - w^N_{t_2}(x_2)| \leq  | X_{t_0,t_1}(x_1) - X_{t_0,t_2}(x_2)| + K_A  N^{\cexp{tail_right_N}} e^{-(1-\cexp{tail_right_exp})x_1 }(|x_1 - x_2|^{\cexp{holder}} + |t_1 - t_2|^{\cexp{holder}/2}).
	\\ &\hspace{1 cm}+  e^{\|f'\|_\infty(t_1-t_0)} \left| \int_\R(G_{t_1 - t_0} (x_1 + \alpha(t_1 - t_0) - y)  - G_{t_2 - t_0} (x_2 + \alpha(t_2 - t_0) - y))w^N_{t_0}(y) dy \right|.
	\end{align*}
	for all $t_1,t_2 \in [t_0,(t_0+2)\wedge N \sigmaN{tail_right}]$ and $x_1, x_2 \in \R$ with $|x_1 - x_2 | \leq 1$.
	\end{lemma}
	\begin{proof}
	Our starting point is the representation for $w^N_t(x)$ for $t \leq N \sigmaN{tail_right}$ given in \eqref{eq:wNtXt0t}. Let $t_0 \geq 0$ and suppose that $t_0 \leq N \sigmaN{tail_right}$. 
	Let $t_1, t_2 \in [t_0, (t_0 +2)\wedge N \sigmaN{tail_right}]$ and $x_1, x_2 \in [x,x+1]$. 
	We also assume without loss of generality that $t_1 \leq t_2$. We also define
	\begin{equation}\label{eq:z1z2expholder}
	z_1(s) := x_1 + \alpha(t_1 - s), \quad z_2(s) := x_2 + \alpha(t_2-s).
	\end{equation}
	Then from \eqref{eq:wNtXt0t},
	\begin{align} \label{eq:wNinc1}
	&w^N_{t_1}(x_1) - w^N_{t_2}(x_2) 
	\\ &\hspace{.6 cm} = \int_\R(G_{t_1 - t_0} (z_1(t_0) - y) e^{\|f'\|_\infty(t_1-t_0)} - G_{t_2 - t_0} (z_2(t_0) - y) e^{\|f'\|_\infty(t_2-t_0)}) w^N_{t_0}(y) dy \notag
	\\ &\hspace{.6 cm} \quad+ \int_{t_0}^{t_1}  \int_\R (G_{t_1 - s} (z_1(s) - y) e^{\|f'\|_\infty(t_1-s)} - G_{t_2 - s} (z_2(s) - y) e^{\|f'\|_\infty(t_2-s)}) \notag
	\\ &\hspace{2 cm} \quad \times (f(w^N_{s}(y)) - \|f'\|_\infty w^N_{s}(y)) dy ds \notag
	\\ &\hspace{.6 cm} \quad - \int_{t_1}^{t_2}  \int_\R G_{t_2 - s} (z_2(s) - y) e^{\|f'\|_\infty(t_2-s)} (f(w^N_{s}(y)) - \|f'\|_\infty w^N_{s}(y)) dy ds \notag
	\\&\hspace{.6 cm} \quad +X_{t_0,t_1}(x_1) - X_{t_0,t_2}(x_2). \notag
	\\ &\hspace{0.6 cm} =: D_1(t_1,t_2,x_1,x_2) + D_2(t_1,t_2,x_1,x_2) + D_3(t_1,t_2,x_2)  + X_{t_0,t_1}(x_1) - X_{t_0,t_2}(x_2). \notag
	\end{align}
	First, we consider $D_2(t_1,t_2,x_1,x_2)$. Let
	\[J(s,y) =f(w^N_{s}(y)) - \|f'\|_\infty w^N_{s}(y) .\]
	By definition we have $J(s,y) \leq 0$. We also remark that 
	\begin{equation} \label{eq:Jlipbd}
	|J(s,y)|  \leq 2 \|f'\|_\infty N^{\cexp{tail_right_N}}  e^{-(1-\cexp{tail_right_exp}) y} \quad \forall \, s \in [0, N \sigmaN{tail_right}] \text{ and } y \in \R.
	\end{equation}	
	Next, we observe that
	\begin{align} \label{eq:D2split}
	&D_2(t_1,t_2,x_1,x_2)  
	\\ & \hspace{ 1cm}=  \int_{t_0}^{t_1}  \int_\R \left(e^{\|f'\|_\infty (t_1 - s)} - e^{\|f'\|_\infty(t_2-s)} \right)(G_{t_2 - s} (z_2(s) - y) )  J(s,y)  dy ds \notag
	\\ &\hspace{ 1cm}  \quad +  \int_{t_0}^{t_1} e^{\|f'\|_\infty(t_1-s)} \int_\R (G_{t_1 - s} (z_1(s) - y)  - G_{t_2 - s} (z_2(s) - y)) J(s,y) dy ds.  \notag
	\end{align}
	Using the fact that $t_1, t_2 \in [t_0,t_0+2]$, and then \eqref{eq:Jlipbd}, we have
	\begin{align} \label{eq:D2firstbd}
	&\left| \int_{t_0}^{t_1}  \int_\R \left(e^{\|f'\|_\infty (t_1 - s)} - e^{\|f'\|_\infty(t_2-s)} \right)(G_{t_2 - s} (z_2(s) - y) )  J(s,y)  dy ds \right| \notag
	\\ &\hspace{1 cm} \leq 2\|f'\|_\infty e^{2 \|f'\|_\infty} (t_2 -t_1)  N^{\cexp{tail_right_N}}  \int_{t_0}^{t_1} \int_\R (G_{t_2 - s} (z_2(s) - y) ) e^{-(1-\cexp{tail_right_exp}) y}  dy ds \notag 
	\\ &\hspace{1 cm} = 2\|f'\|_\infty e^{2 \|f'\|_\infty} (t_2 -t_1)  N^{\cexp{tail_right_N}}  \int_{t_0}^{t_1} \mathbb{E}_{z_2(s)}\left[e^{-(1-\cexp{tail_right_exp})B_{2(t_1 - s)}} \right] dy ds \notag
	\\ &\hspace{1 cm} = 2\|f'\|_\infty e^{2 \|f'\|_\infty} (t_2 -t_1)  N^{\cexp{tail_right_N}}  \int_{t_0}^{t_1} e^{-(1-\cexp{tail_right_exp})z_2(s) + (1-\cexp{tail_right_exp})^2(t_1 - s) }  dy ds \notag
	\\ &\hspace{1 cm} \leq 4\|f'\|_\infty e^{2 \|f'\|_\infty + 2} (t_2 -t_1)  N^{\cexp{tail_right_N}}  e^{-(1-\cexp{tail_right_exp})x_2 }.
	\end{align}
	In the third line, $(B_t)_{t \geq 0}$ is a Brownian motion, and the fourth line is a standard identity. 
	The last line uses $t_1 - t_0 \leq 2$ twice, as well as $z_2(s) \geq x_2$ for all $s \in [t_0,t_1]$, which follows from \eqref{eq:z1z2expholder}. 
	In the sequel, simplifications from $t_1,t_2 \in [t_0,t_0 + 2]$ and $z_1(s), z_2(s) \geq x$ will be made without explicit reference.
	To handle the other term in \eqref{eq:D2split}, we first consider the integrand for a fixed $s \in [t_0,t_1]$. We use \eqref{eq:Jlipbd}, then argue using Hölder's inequality as in \eqref{eq:GintIbd} to obtain
	\begin{align} \label{eq:D2holder}
	&\left|  e^{\|f'\|_\infty(t_1-s)} \int_\R (G_{t_1 - s} (z_1(s) - y)  - G_{t_2 - s} (z_2(s) - y)) J(s,y) dy  \right|  \notag 
	\\ &\leq e^{2\|f'\|_{\infty}} N^{\cexp{tail_right_N}} \int_{t_0}^{t_1} \int_\R |G_{t_1 - s} (z_1(s) - y)  - G_{t_2 - s} (z_2(s) - y)|^{1/2} e^{-(1-\cexp{tail_right_exp})y } dy   \notag
	\\ &\leq e^{2\|f'\|_{\infty}} N^{\cexp{tail_right_N}} \left( \int_\R |G_{t_1 - s} (z_1(s) - y)  - G_{t_2 - s} (z_2(s) - y)|^{1/2} e^{- \frac 3 2 (1-\cexp{tail_right_exp})y } dy \right)^{2/3} \notag
	\\ &\hspace{2 cm} \times \left(\int_\R (G_{t_1 - s} (z_1(s) - y)  - G_{t_2 - s} (z_2(s) - y))^2 dy ds \right)^{1/3}.
	\end{align}
	Observing that for any $t>0$ and $y \in R$, $G_t(y)^{1/2} = 2 (\pi t)^{1/4} G_{2t}(y)$, we have
	\begin{align}
	& \int_\R |G_{t_1 - s} (z_1(s) - y)  - G_{t_2 - s} (z_2(s) - y)|^{1/2} e^{-3 (1-\cexp{tail_right_exp})y } dy\notag
	\\ &\leq \int_\R (G_{t_1 - s} (z_1(s) - y)^{1/2}  + G_{t_2 - s} (z_2(s) - y)^{1/2}) e^{-\frac 3 2 (1-\cexp{tail_right_exp})y } dy \notag
	\\ &\leq 2 (\pi (t_1-s))^{1/4} \mathbb{E}_{z_1(s)} \left[e^{-\frac 32(1-\cexp{tail_right_exp})B_{4(t_1-s)}} \right] +  2 ( \pi (t_2-s))^{1/4} \mathbb{E}_{z_2(s)} \left[e^{-\frac 32(1-\cexp{tail_right_exp})B_{4(t_2-s)}} \right] \notag
	\\ &\leq 4 (t_1-s)^{1/4} e^{-\frac 32(1-\cexp{tail_right_exp})z_1(s) + \frac 92(1-\cexp{tail_right_exp})^2(t_1 - s) } + 4 (t_2-s)^{1/4} e^{-\frac 32(1-\cexp{tail_right_exp})z_2(s) + \frac 92(1-\cexp{tail_right_exp})^2(t_2 - s) } \notag
	\\ &\leq 4 e^{9} e^{-\frac 32 (1-\cexp{tail_right_exp})(x_1 \wedge x_2)} \left(  (t_1-s)^{1/4}  +  (t_2-s)^{1/4}  \right). \notag
	\end{align} 
	We recall the definition of $I_{t_1,t_2}(s,x)$ from \eqref{eq:Idefn} and remark that the last line of \eqref{eq:D2holder} is equal to $I_{t_1,t_2}(z_1(s) - z_2(s),s)^{1/3}$. 
	Hence, by the above and \eqref{eq:D2holder}, 
	\begin{align*} 
	&\left| \int_{t_0}^{t_1} e^{\|f'\|_\infty(t_1-s)} \int_\R (G_{t_1 - s} (z_1(s) - y)  - G_{t_2 - s} (z_2(s) - y)) J(s,y) dy ds \right|  \notag 
	\\ &\leq (4e^{9})^{2/3}e^{ 2\|f'\|_{\infty}} N^{\cexp{tail_right_N}} e^{-(1-\cexp{tail_right_exp})x} \int_{t_0}^{t_1} \left(  (t_1-s)^{1/4}  +  (t_2-s)^{1/4}  \right) I_{t_1,t_2}(z_1(s) - z_2(s),s)^{1/3} ds \notag
	\\ & \leq 8(4e^{9})^{2/3}e^{ 2\|f'\|_{\infty}} N^{\cexp{tail_right_N}} e^{-(1-\cexp{tail_right_exp})(x_1 \wedge x_2)} \int_{t_0}^{t_1} I_{t_1,t_2}(z_1(s) - z_2(s),s)^{1/3} ds. \notag
	\end{align*}
	To control the remaining integral term, we argue exactly as in \eqref{eq:GGint1} and the preceding calculation, from which we conclude that for some constant $K_1  \geq 1$,
	\begin{align*} 
	&\left| \int_{t_0}^{t_1} e^{\|f'\|_\infty(t_1-s)} \int_\R (G_{t_1 - s} (z_1(s) - y)  - G_{t_2 - s} (z_2(s) - y)) J(s,y) dy ds \right|  
	\\ &\hspace{1 cm}\leq K_1 N^{\cexp{tail_right_N}} e^{-(1-\cexp{tail_right_exp})(x_1 \wedge x_2)} (|x_1-x_2|^{4/9} + |t_1 - t_2|^{2/9}). 
	\end{align*}
	Combined with \eqref{eq:D2firstbd}, and modifying the value of $K_1$ as necessary, we obtain that
	\begin{align} \label{eq:D2control}
	&\text{On $\{ t_0 + 2 \leq 6\sigmaN{tail_right}\}$, for all $t_1,t_2 \in [t_0,t_0+2]$ and $x_1, x_2 \in [x,x+1]$,} \notag
	\\ &\hspace{2 cm} |D_2(t_1,t_2,x_1,x_2)| \leq K_1 N^{\cexp{tail_right_N}} e^{-(1-\cexp{tail_right_exp})x} (|x_1-x_2|^{4/9} + |t_1 - t_2|^{2/9}).
	\end{align}
	We can handle $D_3(t_1,t_2,x_2)$ using essentially the same argument. We have
	\begin{align*}
	|D_3(t_1,t_2,x_2)| &\leq e^{2\|f'\|_\infty}  \int_{t_1}^{t_2}  \int_\R G_{t_2 - s} (z_2(s) - y) |J(s,y)| dy ds  
	\\ &\leq 2\|f'\|_\infty e^{2\|f'\|_\infty} N^{\cexp{tail_right_N}}  \int_{t_1}^{t_2}  \int_\R G_{t_2 - s} (z_2(s) - y) e^{-(1-\cexp{tail_right_exp})y} dy ds.
	\end{align*}
	To obtain an upper bound for the above, we argue as in \eqref{eq:D2firstbd}, which implies that for some $K_2 \geq 1$,
	\begin{align} \label{eq:D3control}
	&\text{On $\{ t_0 + 2 \leq 6\sigmaN{tail_right}\}$, for all $t_1,t_2 \in [t_0,t_0+2]$ and $x_1, x_2 \in [x,x+1]$,} \notag
	\\ &\hspace{2 cm} |D_3(t_1,t_2,x_2)| \leq K_2 N^{\cexp{tail_right_N}} e^{-(1-\cexp{tail_right_exp})x_2} |t_1 - t_2|.
	\end{align}
	Finally, we consider the term $D_1(t_1,t_2,x_1,x_2)$ from \eqref{eq:wNinc1}. We observe that
	\begin{align} \label{eq:D1split}
	&D_1(t_1,t_2,x_1,x_2) \notag 
	\\ &\leq  e^{\|f'\|_\infty(t_1-t_0)} \int_\R(G_{t_1 - t_0} (z_1(t_0) - y)  - G_{t_2 - t_0} (z_2(t_0) - y) e^{\|f'\|_\infty(t_2-t_0)}) w^N_{t_0}(y) dy \notag
	\\ &\quad + \left(  e^{\|f'\|_\infty(t_1-t_0)}  -  e^{\|f'\|_\infty(t_2-t_0)} \right) \int_\R G_{t_2 - t_0} (z_2(t_0) - y) w^N_{t_0}(y) dy.
	\end{align}
	The second term above can be handled by a slight modification of the calculation in \eqref{eq:D2firstbd}. In particular, there exists a constant $K_3 \geq 1$ such that
	\begin{align*}
	&\text{On $\{ t_0 + 2 \leq 6\sigmaN{tail_right}\}$, for all $t_1,t_2 \in [t_0,t_0+2]$ and $x_1, x_2 \in [x,x+1]$,} 
	\\ &\hspace{1 cm} |D_1(t_1,t_2,x_1,x_2)| \leq K_3 N^{\cexp{tail_right_N}} e^{-(1-\cexp{tail_right_exp})x_2} |t_1 - t_2| \notag
	\\ &\hspace{1 cm} \quad + e^{\|f'\|_\infty(t_1-t_0)} \int_\R(G_{t_1 - t_0} (z_1(t_0) - y)  - G_{t_2 - t_0} (z_2(t_0) - y) e^{\|f'\|_\infty(t_2-t_0)}) w^N_{t_0}(y) dy. \notag
	\end{align*}
	We recall that $\cexp{holder} \leq 1/4 < 4/9$.
	Using the fact that $|x_1 - x_2| \leq 1$, the result now follows from \eqref{eq:wNinc1}, \eqref{eq:D2control}, \eqref{eq:D3control} and the above. 
	\end{proof}

	It will be useful to codify the integrability offered by $\sigmaN{tail_right}$; the next two lemmas concern functions $w \in C(\R,[0,1])$ satisfying
	\begin{equation} \label{eq:wintegright}
	w(y) \leq N^{\cexp{tail_right_N}} e^{-(1-\cexp{tail_right_exp}) y} \,\, \forall \, y \in \R.
	\end{equation}

	\begin{lemma} \label{lem:Gexpsmooth}
	There is a constant $K_B \geq 1$, such that if $w : \R \to [0,1]$ satisfies \eqref{eq:wintegright}, then
	\begin{equation*} 
	|G_{t_1} * w(x_1) - G_{t_2} * w(x_2)| \leq  K_B N^{\cexp{tail_right_N}} e^{-((1-\cexp{tail_right_exp})/2)(x_1 \wedge x_2)} (|x_1 - x_2| + |t_1 - t_2|)
	\end{equation*}
	for all $x_1,x_2 \in \R$ and $t_1, t_2 \in [1,2]$
	\end{lemma}
	\begin{proof}
	Write $W(t,x) = G_t * w(x)$. For $t \in [1,2]$ and $x \in \R$, taking the derivative under the integral and using \eqref{eq:wintegright},
	\begin{align*}
	|\partial_x W(t,x)| &= \left| \int_\R \frac{-(x-y)}{2\sqrt{4\pi t}} e^{-(x-y)^2/(4t)} w(y) dy \right|
	\\ &\leq  N^{\cexp{tail_right_N}} e^{-(1-\cexp{tail_right_exp})x} \int_\R \frac{|x-y|}{2\sqrt{4\pi t}} e^{-(x-y)^2/(4t)} e^{-(1-\cexp{tail_right_exp})(y-x)} dy. 
	\end{align*}
	By a variation of the argument in \eqref{eq:Gincexpholder2} using Cauchy-Schwarz, the integral in the above is bounded above uniformly for $t \in [1,2]$ by some finite constant. In particular,
	Thus, for some $K_1 \geq 1$, for all $x \in \R$ and $t \in [1,2]$,
	\begin{equation} \label{eq:dxWtx}
	|\partial_x W(t,x)| \leq K_1 N^{\cexp{tail_right_N}} e^{-(1-\cexp{tail_right_exp})x}.
	\end{equation}
	Computing $\partial_t W(t,x)$ using the same argument, we obtain
	\begin{align*}
	|\partial_t W(t,x)| &\leq N^{\cexp{tail_right_N}} e^{-(1-\cexp{tail_right_exp})x} \int_\R  \left(\frac{|x-y|^2}{4t^2} - \frac{1}{2t} \right) G_t(x-y) e^{-(1-\cexp{tail_right_exp})(y-x)} dy,
	\end{align*}
	and we conclude in the same way that for some $K_2 \geq 1$, for all $x \in \R$ and $t \in [1,2]$,
	\begin{equation} \label{eq:dtWtx}
	|\partial_t W(t,x)| \leq K_2 N^{\cexp{tail_right_N}} e^{-(1-\cexp{tail_right_exp})x}.
	\end{equation}
	Let $t_1 ,t_2 \in [1,2]$ and $x_1, x_2 \in \R$. Suppose first that $x_1 ,x_2 \geq 0$. Then the desired result follows from \eqref{eq:dxWtx} and \eqref{eq:dtWtx} and the triangle inequality.
	If $x_1, x_2 \leq 0$, a simpler version of the arguments above show that, since $w$ takes values in $[0,1]$, 
	$|\partial_x W(t,x)|$ and $\partial_t W(t,x)|$ are bounded above uniformly on $(t,x) \in [1,2] \times (-\infty,0]$. The result then follows for $x_1, x_2 \leq 0$. 
	Finally, if $x_1 > 0$ and $x_2 < 0$, we write $|x_1 - x_2| = |x_1| + |x_2|$, and the inequality
	follows from the first two cases. This completes the proof.
	\end{proof}
	
	We will need the following lemma to handle the contribution of the initial conditions $u_0^N$ for small values of $t$.
	\begin{lemma} \label{lem:expholderinitial} 
	There exists a constant $K_C \geq 1$ such that the following holds: if $w \in C(\R,[0,1])$ satisfies \eqref{eq:wintegright} and
	\begin{equation} \label{eq:wexpholder}
	|w(x_1) - w(x_2)|  \leq   N^{\cexp{tail_right_N}} e^{-((1-\cexp{tail_right_exp})/2)x_1} |x_1 - x_2|^{\cexp{holder}} \,\, \forall \,x_1,x_2 \in \R \text{ with } |x_1 - x_2| \leq 1,\end{equation}
	then
	\begin{align*}
	&|G_{t_1} * w (x_1 + \alpha t_1) - G_{t_2} * w(x_2 + \alpha t_2)| 
	\\ &\hspace{ 1 cm} \leq  K_C \cdot   N^{\cexp{tail_right_N}} e^{-((1-\cexp{tail_right_exp})/2)x_1} (|t_1 - t_2|^{ {\cexp{holder}}/2} + |x_1 - x_2 |^{{\cexp{holder}}})
	\end{align*}
	for all $t_1,t_2 \in [0,2]$ and $x_1,x_2 \in \R$ satisfying $|x_1 - x_2| \leq 1$.
	\end{lemma}
	\begin{proof}
	Let $w \in C(\R,[0,1])$ satisfy \eqref{eq:wintegright} and \eqref{eq:wexpholder}. We first remark that 
	\begin{equation*}
	w(y) \leq N^{\cexp{tail_right_N}} e^{-((1-\cexp{tail_right_exp})/2)y} \,\, \forall \, y \in \R,
	\end{equation*}
	where we use \eqref{eq:wintegright} for $y >0$ and $w(y) \leq 1 \leq N$ for $y \leq 0$. Hence, for $x_1, x_2 \in \R$ with $|x_1 - x_2| > 1$,
	\begin{align*}
	|w(x_1) - w(x_2)| \leq w(x_1) + w(x_2) &\leq 2 N^{\cexp{tail_right_N}} e^{-((1-\cexp{tail_right_exp})/2)(x_1 \wedge x_2)} 
	\\ &\leq 2 N^{\cexp{tail_right_N}} e^{-((1-\cexp{tail_right_exp})/2)(x_1 \wedge x_2)} |x_1-x_2|^{\cexp{holder}}.
	\end{align*}
	Combined with \eqref{eq:wexpholder}, this implies that
	\begin{align} \label{eq:wincrement1}
	|w(x_1) - w(x_2)| &\leq 2e^{(1-\cexp{tail_right_exp})} \cdot N^{\cexp{tail_right_N}} e^{-((1-\cexp{tail_right_exp})/2)(x_1 \wedge x_2)} |x_1-x_2|^{\cexp{holder}} \,\, \forall \, x_1, x_2 \in \R.
	\end{align}
	Let $t_1 ,t_2 \in [0,2]$ and $x_1,x_2 \in \R$ satisfy $|x_1 - x_2 | \leq 1$. For $i=1$ and $2$, let $z_i = x_i + \alpha t_i$. Arguing as in \eqref{eq:1sttermholder0}, we remark that
	\begin{align} 
	&|G_{t_1} * w (z_1) - G_{t_2} * w(z_2)| \leq \int_\R G_{t_1}(z_1 - y) \int_\R G_{t_2 - t_1} (y-z) |w(y) - w(z +z_2 - z_1)| dz dy. \notag
	\end{align}
	By \eqref{eq:wincrement1}, 
	\begin{align*}
	 &|w(y) - w(z +z_2 - z_1)|  
	 \\ &\hspace{ 1 cm} \leq  2e^{(1-\cexp{tail_right_exp})} \cdot N^{\cexp{tail_right_N}}e^{-((1-\cexp{tail_right_exp})/2)(y \wedge (z - z_2 + z_1))} |(y-z) + (z_2-z_1)|^{{\cexp{holder}}}
	 \\ &\hspace{1 cm}\leq   2e^{((1-\cexp{tail_right_exp})/2)(3+2\alpha)}  \cdot  N^{\cexp{tail_right_N}} e^{-((1-\cexp{tail_right_exp})/2)(y \wedge z)} (|y-z|^{{\cexp{holder}}} + |z_2-z_1|^{{\cexp{holder}}}),
	\end{align*}
	where we have used $z_2 - z_1 \geq -1-2\alpha$. Substituting the above into the previous equation, we obtain that for some constant $K_1$,
	\begin{align} \label{eq:Gincexpholder1} 
	&|G_{t_1} * w (z_1) - G_{t_2} * w(z_2)| 
	\\ &\hspace{1 cm} \leq K_1  \cdot N^{\cexp{tail_right_N}} \int_\R G_{t_1}(z_1 - y) \int_\R G_{t_2 - t_1} (y-z) \notag
	\\ &\hspace{ 3cm} \times \left( e^{-((1-\cexp{tail_right_exp})/2)y} + e^{-((1-\cexp{tail_right_exp})/2)z}\right) (|y-z|^{{\cexp{holder}}} + |z_2-z_1|^{{\cexp{holder}}}) dz dy. \notag
	\end{align}
	We handle the above term by term. First observe that
	\begin{align} \label{eq:Gincexpholder2} 
	e^{-((1-\cexp{tail_right_exp})/2)y} \int_\R G_{t_2 - t_1} (y-z) |y-z|^{{\cexp{holder}}} dz = e^{-((1-\cexp{tail_right_exp})/2)y} |t_2 -t_1|^{{\cexp{holder}}/2} \mathbb{E}_0[ |B_2|^{{\cexp{holder}}}].
	\end{align} 
	Using Cauchy-Schwarz, we also compute
	\begin{align} \label{eq:Gincexpholder3}
	&\int_\R G_{t_2 - t_1} (y-z) e^{-((1-\cexp{tail_right_exp})/2)z} |y-z|^{{\cexp{holder}}} dz  \notag
	\\ &\hspace{1 cm}= e^{-(1-\cexp{tail_right_exp})y}\int_\R G_{t_2 - t_1} (y-z) e^{-((1-\cexp{tail_right_exp})/2)(z-y)} |y-z|^{{\cexp{holder}}} dz  \notag
	\\ &\hspace{1 cm}= e^{-((1-\cexp{tail_right_exp})/2)y} \mathbb{E}_0[ |B_{2(t_2-t_1)}|^{{\cexp{holder}}} e^{-((1-\cexp{tail_right_exp})/2)B_{2(t_2-t_1)}}] \notag
	\\ &\hspace{1 cm}\leq e^{-((1-\cexp{tail_right_exp})/2)y} \mathbb{E}_0[ |B_{2(t_2-t_1)}|^{2 {\cexp{holder}}}]^{1/2} \mathbb{E}_0[ e^{-(1-\cexp{tail_right_exp})B_{2(t_2-t_1)}}]^{1/2} \notag 
	\\ &\hspace{1 cm}= e^{-((1-\cexp{tail_right_exp})/2)y} (t_2-t_1)^{{\cexp{holder}}/2} \mathbb{E}_0[|B_2|^{2 {\cexp{holder}}}]^{1/2} e^{2(1-\cexp{tail_right_exp})^2},
	\end{align}
	where in the last line we have used $t_2 -t_1 \leq 2$.
	Similarly, one can argue that
	\begin{align*} 
	& \int_\R G_{t_2 - t_1} (y-z) \left( e^{-((1-\cexp{tail_right_exp})/2)y} + e^{-((1-\cexp{tail_right_exp})/2)z}\right)|z_2-z_1|^{{\cexp{holder}}} dz
	\\&\hspace{1 cm} \leq \left( 1+e^{2(1-\cexp{tail_right_exp})^2}\right) e^{-((1-\cexp{tail_right_exp})/2)y}  |z_1 - z_2|^{{\cexp{holder}}}.
	\end{align*}
	By the above, \eqref{eq:Gincexpholder1}, \eqref{eq:Gincexpholder2} and \eqref{eq:Gincexpholder3}, it follows that, for an increased value of the constant $K_1$, 
	\begin{align*}
	&|G_{t_1} * w (z_1) - G_{t_2} * w(z_2)| 
	\\ &\leq  K_1 \cdot N^{\cexp{tail_right_N}} \int_\R G_{t_1}(z_1 - y) e^{-((1-\cexp{tail_right_exp})/2)y} \left( (t_2-t_1)^{{\cexp{holder}}/2}  + |z_1 - z_2|^{{\cexp{holder}}} \right)dy
	\\ &\leq  K_1 \cdot \left( (t_2-t_1)^{{\cexp{holder}}/2}  + |z_1 - z_2|^{{\cexp{holder}}} \right)N^{\cexp{tail_right_N}} e^{-((1-\cexp{tail_right_exp})/2)z_1} e^{(1-\cexp{tail_right_exp})^2}.
	 \end{align*}
	We recall that $z_i = x_i + \alpha t_i$. Using this and its implication $z_1 \geq x_1$ completes the proof.
	\end{proof}

		
	\begin{proof}[Proof of Proposition~\ref{prop:holderexp}]
	Let $N \geq 1$. We begin by recalling the random variables $Y^N_{t_0,x,k}$ from Lemma~\ref{lem:Xtholder}. 
	Let $k \in \N$ be large enough so that $k \cexp{tail_right_N} >4$.
	For $n \in \N_0$ define
	\begin{align*}
	E^N_n := &\left\{ Y_{n,x,k} \leq N^{ \cexp{tail_right_N}/2} e^{-((1-\cexp{tail_right_exp})/2)(x+2)} \,\,\, \forall \, x \in \Z \cap [1,1+(1+\cexp{interface})\log N] \right\}.
	\end{align*}
	Next, we remark that 
	\begin{align} \label{eq:holderXnbd}
	&\text{On $E^N_n$, for all $t_1,t_2 \in [n,n+2]$ and $x_1, x_2 \in [0, (1+\cexp{interface})\log N]$ with $|x_1 - x_2| \leq 1$,  } \notag
	\\ &\qquad  |X_{n,t_1}(x_1) - X_{n,t_2 }(x_2)| \leq N^{\cexp{tail_right_N}/2} e^{-((1-\cexp{tail_right_exp})/2)x_1} (|x_1 - x_2|^{\cexp{holder}} + |t_1 - t_2|^{\cexp{holder}/2} )
	\end{align}
	On $\{n \leq N\sigmaN{tail_right}\}$, $w^N_n$ satisfies \eqref{eq:wintegright}.
	Thus, by Lemma~\ref{lem:Gexpsmooth}, on $ \{n \leq N\sigmaN{tail_right}\}$ we have
	\begin{align*} 
	&|G_{t_1-n}*w^N_n (x_1) - G_{t_2-n}* w^N_n(x_2) | \leq K_B e^1 \cdot N^{\cexp{tail_right_N}}e^{-((1-\cexp{tail_right_exp})/2)x_1} (|x_1 - x_2|+ |t_1 -t_2|) \notag
	\end{align*}
	for all $t_1,t_2 \in [n+1,n+2]$ and $x_1,x_2 \in \R$ satisfying $|x_1 - x_2| \leq 1$. 
	By Lemma~\ref{lem:holderexpbd} combined with \eqref{eq:holderXnbd} and the above, we conclude that, on $E^N_n \cap \{n + 1 \leq N \sigmaN{tail_right}\}$, 
	for all $t_1,t_2 \in [(n+1) \wedge N\sigmaN{tail_right},(n+2) \wedge N\sigmaN{tail_right}]$ and $ x_1, x_2 \in [0,1+ (1+\cexp{interface})\log N]$ such that $|x_1 - x_2 | \leq 1$,
	\begin{align} \label{eq:expholderfinal}
	&|w^N_{t_1}(x_1) - w^N_{t_2}(x_2)|  
	\\ &\hspace{ 1cm}\leq (1 + K_A + e^{1+2\|f'\|_\infty}K_B) N^{\cexp{tail_right_N}} e^{-((1-\cexp{tail_right_exp})/2)(x_1 \wedge x_2)}  (|x_1 - x_2|^{\cexp{holder}} + |t_1 -t_2|^{\cexp{holder}/2}). \notag
	\end{align}
	where we assume $N\geq 1$ to guarantee $N^{\cexp{tail_right_N}/2}  \leq N^{\cexp{tail_right_N}}$.
	Recalling the definition of $\sigmaN{interface}$ from \eqref{eq:sigmainterface}, we remark that that if $t_1, t_2 \leq N \sigmaN{interface}$ and $x_1, x_2 \geq (1+\cexp{interface})\log N$,
	then $w^N_{t_1}(x_1) - w^N_{t_2}(x_2) = 0$, and hence \eqref{eq:expholderfinal} holds trivially for such $t_1$, $t_2$, $x_1$ and $x_2$. 
	(If $t_i = N \sigmaN{interface}$ for $i=1$ or $2$, we additionally appeal to Lemma~\ref{lemma_Rsimple}(b) to see that $w^N_{t_i}(\cdot)$ vanishes above $ (1+\cexp{interface})\log N$.)
	This implies that on $ E^N_n \cap \{n +1 \leq N(\sigmaN{interface} \wedge \sigmaN{tail_right})\}$, 
	\eqref{eq:expholderfinal} holds for all $t_1 , t_2 \in [(n+1) \wedge N(\sigmaN{interface} \wedge \sigmaN{tail_right})], (n+2) \wedge N(\sigmaN{interface} \wedge \sigmaN{tail_right})]$ and $x_1, x_2 \geq 0$
	satisfying $|t_1 - t_2| \leq 1$ and $|x_1 - x_2| \leq 1$. 
	In particular, we have shown that on the event $E = \cap_{n=0}^{N^2-2} E^N_n$, 
	\eqref{eq:expholderfinal} holds for all $t_1,t_2 \in [1, N(\sigmaN{interface} \wedge \sigmaN{tail_right}) \wedge N^2]$ and $x_1, x_2 \geq 0$ (which satisfy $|t_1 - t_2| \leq 1$ and $|x_1 - x_2| \leq 1$).
	
	It remains to handle the case that either $t_1$ or $t_2 \in [0,1]$. This will be satisfied if we can prove that \eqref{eq:expholderfinal} holds for all $t_1,t_2 \in [0,2]$. 
	We will apply Lemma~\ref{lem:holderexpbd} with $t_0 = 0$, this time using Lemma~\ref{lem:expholderinitial} instead of Lemma~\ref{lem:Gexpsmooth} to handle the term from the initial data. 
	We remark that by Assumption2~\ref{assumpt:v0}.\ref{v0:holder} and ~\ref{assumpt:v0}.\ref{v0:tail}, as well as \eqref{eq:holderconstdef} and \eqref{bounds_cinitsmall}, 
	that $u_0^N$ satisfies \eqref{eq:wintegright} and \eqref{eq:wexpholder}.
	Hence, we can apply Lemma~\ref{lem:expholderinitial} with $u^N_0$, and we obtain that 
	\begin{equation*}
	|G_{t_1} * u^N_0(x_1 + \alpha t_1) - G_{t_2} * u^N_0(x_2+\alpha t_2)| \leq K_C \cdot N^{\cexp{tail_right_N}} e^{-((1-\cexp{tail_right_exp})/2)x_1}(|x_1 - x_2|^{\cexp{holder}} + |t_1 - t_2|^{\cexp{holder}/2})
	\end{equation*}
	for all $x_1,x_2 \in \R$ with $|x_1 - x_2| \leq 1$. We then apply Lemma~\ref{lem:holderexpbd} with $t_0 = 0$ and use \eqref{eq:holderXnbd}
	to conclude that on $E^N_0$, \eqref{eq:expholderfinal} holds, with $K_C$ replacing $K_B$, 
	for all $t_1, t_2 \in [0, \sigmaN{tail_right} \wedge 2]$ and $x_1,x_2 \in \R$ with $|x_1 - x_2| \leq 1$.

	To conclude, we obtain an upper bound on the probability of $E^c$, where $E = \cap_{n=0}^{N^2-1} E^N_n$. 
	By a union bound, Markov's inequality, and Lemma~\ref{lem:Xtholder},
	\begin{align*}
	\bP(E^c) & \leq \sum_{n = 0}^{N^2 - 1} \sum_{x  \in \Z \cap [0, (1+\cexp{interface})\log N]} \P{ Y_{n,x,k} \geq N^{a' + \cexp{tail_right_N}/2} e^{-((1-\cexp{tail_right_exp})/2)(x+1)}}
	\\ &\leq  A'_k e^{k(1-\cexp{tail_right_exp})}N^2 ((1+\cexp{interface})\log N + 1)  N^{-k\cexp{tail_right_N}} 
	\end{align*}
	Since $k\cexp{tail_right_N} > 4$, the above is at most $N^{-2}$ for sufficiently large $N$. This completes the proof.	
	\end{proof}

	\subsection{Proof of Lemma~\ref{lem:Gdiffintbound}} \label{subsec:techG-G}
	We end Section~\ref{sec:holder} by proving the technical result stated in Section~\ref{subsec:holder}.
	\begin{proof}[Proof of Lemma~\ref{lem:Gdiffintbound}]
	For $0\le s < t_1 \le t_2$ and $x\in \R$, by~\eqref{eq:Idefn},
	\begin{align} \label{eq:Ibound}
	&I_{t_1,t_2}(s,x) \notag \\
	&=\int_{\R} (G_{t_1-s}(y-x) -G_{t_2-s}(y))^2 dy \notag \\
	&=\int_{\R} \left((4\pi (t_1-s))^{-1/2}e^{-(4(t_1-s))^{-1}(y-x)^2} -(4\pi (t_2-s))^{-1/2}e^{-(4(t_2-s))^{-1}y^2}\right)^2 dy \notag \\
	&=(4\pi (t_1-s))^{-1}\int_{\R} e^{-(2(t_1-s))^{-1}(y-x)^2}dy \notag \\
	&\quad -2(4\pi)^{-1}((t_1-s)(t_2-s))^{-1/2}
	\int_{\R} e^{-(4(t_1-s))^{-1}(y-x)^2-(4(t_2-s))^{-1}y^2}dy \notag \\
	&\quad +(4\pi (t_2-s))^{-1}\int_{\R} e^{-(2(t_2-s))^{-1}y^2}dy.
	\end{align}
	By completing the square, we have that
	\begin{align*}
	&e^{-(4(t_1-s))^{-1}(y-x)^2-(4(t_2-s))^{-1}y^2}\\
	&=
	e^{-\frac 14 ((t_1-s)^{-1}+(t_2-s)^{-1})
	\left(y-\frac{(t_1-s)^{-1}}{(t_1-s)^{-1}+(t_2-s)^{-1}}x
	\right)^2}
	e^{\frac 14 \frac{(t_1-s)^{-2}}{(t_1-s)^{-1}+(t_2-s)^{-1}}x^2}
	e^{-\frac 14 (t_1-s)^{-1}x^2}.
	\end{align*}
	Note that
	$$
	\frac{(t_1-s)^{-2}}{(t_1-s)^{-1}+(t_2-s)^{-1}}-(t_1-s)^{-1}=-(t_1+t_2-2s)^{-1},
	$$
	and so
	$$
	e^{-(4(t_1-s))^{-1}(y-x)^2-(4(t_2-s))^{-1}y^2}
	=e^{-\frac 14 ((t_1-s)^{-1}+(t_2-s)^{-1})
	\left(y-\frac{(t_1-s)^{-1}}{(t_1-s)^{-1}+(t_2-s)^{-1}}x
	\right)^2}
	e^{-\frac 14 x^2 (t_1+t_2-2s)^{-1}}.
	$$
	Hence by~\eqref{eq:Ibound},
	\begin{align} \label{eq:(A)}
	&I_{t_1,t_2}(s,x) \notag \\
	&=(4\pi (t_1-s))^{-1}(2\pi (t_1-s))^{1/2} +(4\pi (t_2-s))^{-1}(2\pi (t_2-s))^{1/2}  \notag \\
	&\qquad -(2\pi)^{-1}((t_1-s)(t_2-s))^{-1/2}(4 \pi ((t_1-s)^{-1}+(t_2-s)^{-1})^{-1})^{1/2}
	 e^{-\frac 14 x^2 (t_1+t_2-2s)^{-1}} \notag \\
	&=(8\pi (t_1-s))^{-1/2}+(8\pi (t_2-s))^{-1/2}
	-(\pi (t_1+t_2-2s))^{-1/2}e^{-\frac 14 x^2 (t_1+t_2-2s)^{-1}} .
	\end{align}
	It follows that for $t_0\in (0,\infty)$, if $0 \le s\le t_1-t_0$ then since $ t_2 \geq t_1 $ we can write
	\begin{align} \label{eq:IboundB}
	I_{t_1,t_2}(s,x)
	&=(4\pi)^{-1/2} ((2t_1-2s)^{-1/2}-(t_1+t_2-2s)^{-1/2}) \notag \\
	&\qquad +(4\pi)^{-1/2} ((2t_2-2s)^{-1/2}-(t_1+t_2-2s)^{-1/2}) \notag \\
	&\qquad +(\pi(t_1+t_2-2s))^{-1/2}(1- e^{-\frac 14 x^2 (t_1+t_2-2s)^{-1}}) \notag \\
	&\le (4\pi)^{-1/2} (t_2-t_1)\cdot \tfrac 12 (2t_1-2(t_1-t_0))^{-3/2} \notag \\
	&\qquad +(\pi (2t_1-2(t_1-t_0)))^{-1/2} \cdot \tfrac 14 x^2 (2t_1-2(t_1-t_0))^{-1} \notag \\
	&=2^{-7/2}\pi^{-1/2}t_0^{-3/2} ((t_2-t_1)+x^2),
	\end{align}
	which completes the proof of~\eqref{eq:Gdiff2}.
	
	It remains to prove~\eqref{eq:Gdiff1}.
	We have by~\eqref{eq:(A)} that in particular,
	$I_{t_1,t_2}(s,x)\le (2\pi (t_1-s))^{-1/2}$.
	Hence by~\eqref{eq:IboundB}, for $a>0$ and $t_0\in (0,t_1]$,
	\begin{align*}
	&\int_0^{t_1} e^{-a(t_1-s)}I_{t_1,t_2}(s,x)ds\\
	&\le \int_0^{t_1-t_0}e^{-a(t_1-s)} 2^{-7/2}\pi^{-1/2}t_0^{-3/2}((t_2-t_1)+x^2)ds
	+\int_{t_1-t_0}^{t_1}(2\pi (t_1-s))^{-1/2}ds.
	\end{align*}
	By taking $t_0=\min(t_1,((t_2-t_1)+x^2)^{1/2})$, it follows that
	\begin{align*}
	\int_0^{t_1} e^{-a(t_1-s)}I_{t_1,t_2}(s,x)ds 
	&\le 2^{-7/2}\pi^{-1/2}((t_2-t_1)+x^2)^{1/4}a^{-1}
	+2^{1/2}\pi^{-1/2}((t_2-t_1)+x^2)^{1/4},
	\end{align*}
	which completes the proof.
	\end{proof}

	\section{Additional H\"older estimates} \label{sec:H_gamma}
	
	In this section, we prove Proposition~\ref{prop:holder}, which states that there exists $\gamma>0$ such that for $T>0$ we have $ v^N_t \in H^{\gamma,\alpha} $ for all $t\in [0,T]$ almost surely and that $ t \mapsto \Hnorm{v^N_t}{\gamma} $ is almost surely locally bounded.
	
	We start by defining, for $ \gamma \in (0,1) $,
	\begin{equation*}
		\| v \|_{\hat{H}^{\gamma,\alpha}} := \left( \| v \|_{2,\alpha}^2 + \int_{\R} \int_{[-1,1]} \frac{|v(x + z) - v(x)|^2}{|z|^{1+2\gamma}} e^{\alpha x} dz dx \right)^{1/2},
	\end{equation*}
	for all $ v \in L^{2,\alpha} $ such that the above integral is finite.
	
	\begin{lemma} \label{lemma:norm_equivalence_Hgamma}
		For any $ \gamma \in (0,1) $, the set $ \lbrace v \in L^{2,\alpha} : \| v \|_{\hat{H}^{\gamma,\alpha}} < + \infty \rbrace $ coincides with $ H^{\gamma,\alpha} $ and the two norms $ \| \cdot \|_{\hat{H}^{\gamma, \alpha}} $ and $ \Hnorm{\cdot}{\gamma} $ are equivalent.
	\end{lemma}
	
	\begin{proof}
		By \citep[Proposition~3.4]{di_nezza_hitchhikers_2012}, for $ \gamma \in (0,1) $, the norm $ \Hnorm{\cdot}{\gamma} $ is equivalent to
		\begin{equation*}
			\left( \| v \|_{2,\alpha}^2 + \int_{\R^2} \frac{| v(y) e^{\frac{\alpha}{2} y} - v(x) e^{\frac{\alpha}{2} x} |^2}{|x-y|^{1 + 2 \gamma}} dx dy \right)^{1/2}.
		\end{equation*}
		Hence it suffices to show that the above is equivalent to $ \| \cdot \|_{\hat{H}^{\gamma,\alpha}} $.
		We start by writing
		\begin{equation*}
			\int_{\R^2} \frac{| v(y) e^{\frac{\alpha}{2} y} - v(x) e^{\frac{\alpha}{2} x} |^2}{|x-y|^{1 + 2 \gamma}} dx dy = \int_{\R^2} \frac{| v(x+z) e^{\frac{\alpha}{2} z} - v(x)|^2}{| z |^{1 + 2\gamma}} e^{\alpha x} dz dx.
		\end{equation*}
		We then split the integral over $ z $ on the right-hand side to write
		\begin{multline} \label{split_int_z}
			\int_{\R^2} \frac{| v(x+z) e^{\frac{\alpha}{2} z} - v(x)|^2}{| z |^{1 + 2\gamma}} e^{\alpha x} dz dx \leq \int_{\R} \int_{[-1,1]} \frac{| v(x+z) e^{\frac{\alpha}{2} z} - v(x)|^2}{| z |^{1 + 2\gamma}} e^{\alpha x} dz dx\\  + \int_{\R \times \lbrace | z | > 1 \rbrace} \frac{| v(x+z) e^{\frac{\alpha}{2} z} - v(x)|^2}{| z |^{1 + 2\gamma}} e^{\alpha x} dz dx.
		\end{multline}
		In the first integral on the right, we use that, for $ z \in [-1,1] $,
		\begin{equation*}
			| v(x+z) e^{\frac{\alpha}{2} z} - v(x)| \leq | v(x+z) - v(x) | e^{\frac{\alpha}{2}} + | v(x) | \, | e^{\frac{\alpha}{2} z} - 1 |,
		\end{equation*}
		and $ | e^{\frac{\alpha}{2} z} - 1 | \leq \frac{\alpha}{2} e^{\frac{\alpha}{2}} |z| $, to obtain
		\begin{multline*}
			\left( \int_{\R} \int_{[-1,1]} \frac{| v(x+z) e^{\frac{\alpha}{2} z} - v(x)|^2}{| z |^{1 + 2\gamma}} e^{\alpha x} dz dx \right)^{1/2} \\ \leq e^{\frac{\alpha}{2}} \left( \int_{\R} \int_{[-1,1]} \frac{| v(x+z) - v(x)|^2}{| z |^{1 + 2\gamma}} e^{\alpha x} dz dx \right)^{1/2} + \frac{\alpha}{2} e^{\frac{\alpha}{2}} \| v \|_{2,\alpha} \left( \int_{[-1,1]} |z|^{1-2\gamma} dz \right)^{1/2},
		\end{multline*}
		and we note that the last integral over $ z $ on the right is finite for $ \gamma \in (0,1) $.
		For the second term on the right of \eqref{split_int_z}, we write
		\begin{equation*}
			| v(x+z) e^{\frac{\alpha}{2} z} - v(x) | e^{\frac{\alpha}{2} x} \leq | v(x+z) | e^{\frac{\alpha}{2} (x+z)} + |v(x)| e^{\frac{\alpha}{2} x}.
		\end{equation*}
		As a result,
		\begin{multline*}
			\left( \int_{\R \times \lbrace | z | > 1 \rbrace} \frac{| v(x+z) e^{\frac{\alpha}{2} z} - v(x)|^2}{| z |^{1 + 2\gamma}} e^{\alpha x} dz dx \right)^{1/2} \leq \left( \int_{\R \times \lbrace | z | > 1 \rbrace} \frac{| v(x+z) |^2}{| z |^{1 + 2\gamma}} e^{\alpha (x+z)} dz dx \right)^{1/2} \\ + \left( \int_{\R \times \lbrace | z | > 1 \rbrace} \frac{| v(x)|^2}{| z |^{1 + 2\gamma}} e^{\alpha x} dz dx \right)^{1/2}.
		\end{multline*}
		By a change of variables in the first integral over $ x $, we see that the two integrals on the right coincide, and factorise to yield
		\begin{equation*}
			\left( \int_{\R \times \lbrace | z | > 1 \rbrace} \frac{| v(x+z) e^{\frac{\alpha}{2} z} - v(x)|^2}{| z |^{1 + 2\gamma}} e^{\alpha x} dz dx \right)^{1/2} \leq \| v \|_{2,\alpha} \left( \int_{\lbrace | z | > 1 \rbrace} \frac{1}{| z |^{1+2\gamma}} dz \right)^{1/2},
		\end{equation*}
		and we note that the integral over $ z $ on the right converges for $ \gamma > 0 $.
		This proves that there exists $ C > 0 $ such that, for all $ v \in L^{2,\alpha} $ such that $ \| v \|_{\hat{H}^{\gamma,\alpha}} < + \infty $,
		\begin{equation*}
			\left( \| v \|_{2,\alpha}^2 + \int_{\R^2} \frac{| v(y) e^{\frac{\alpha}{2} y} - v(x) e^{\frac{\alpha}{2} x} |^2}{|x-y|^{1 + 2 \gamma}} dx dy \right)^{1/2} \leq C \| v \|_{\hat{H}^{\gamma,\alpha}}.
		\end{equation*}
		Now take $ v \in H^{\gamma,\alpha} $ and, writing $ v(x+z) = v(x+z) e^{\frac{\alpha}{2} z} e^{-\frac{\alpha}{2} z} $, we obtain
		\begin{multline*}
			\left( \int_{\R} \int_{[-1,1]} \frac{| v(x+z) - v(x)|^2}{| z |^{1 + 2\gamma}} e^{\alpha x} dz dx \right)^{1/2} \\ \leq \left( \int_{\R} \int_{[-1,1]} \frac{| v(x+z) e^{\frac{\alpha}{2}(x+z)} - v(x) e^{\frac{\alpha}{2} x}|^2}{| z |^{1 + 2\gamma}} dz dx \right)^{1/2} \\ + \left( \int_{\R} \int_{[-1,1]} \frac{| v(x+z) |^2 e^{\alpha (x+z)}}{| z |^{1 + 2\gamma}} | e^{-\frac{\alpha}{2} z} - 1 |^2 dz dx \right)^{1/2}.
		\end{multline*}
		By a change of variables, the first term on the right-hand side is bounded from above by
		\begin{equation*}
			\left( \int_{\R^2} \frac{| v(y)e^{\frac{\alpha}{2} y} - v(x) e^{\frac{\alpha}{2} x} |^2}{|z|^{1+2\gamma}} dx dy \right)^{1/2},
		\end{equation*}
		and, by a change of varibles in the integral over $ x $, the second term on the right is equal to
		\begin{equation*}
			\| v \|_{2,\alpha} \left( \int_{[-1,1]} \frac{| e^{-\frac{\alpha}{2} z} - 1|^2}{|z|^{1+2\gamma}} dz \right)^{1/2},
		\end{equation*}
		and we note that the integral over $ z $ above is finite for $ \gamma \in (0,1) $ since $ | e^{-\frac{\alpha}{2} z} - 1 | \leq \frac{\alpha}{2} e^{\frac{\alpha}{2}} |z| $.
		This shows that there exists another constant $ C > 0 $ such that, for all $ v \in H^{\gamma,\alpha} $,
		\begin{equation*}
			\| v \|_{\hat{H}^{\gamma,\alpha}} \leq C \left( \| v \|_{2,\alpha}^2 + \int_{\R^2} \frac{| v(y) e^{\frac{\alpha}{2} y} - v(x) e^{\frac{\alpha}{2} x} |^2}{|x-y|^{1 + 2 \gamma}} dx dy \right)^{1/2},
 		\end{equation*}
 		and thus concludes the proof of Lemma~\ref{lemma:norm_equivalence_Hgamma}.
 	\end{proof}
 	
 	We now turn to the proof of Proposition~\ref{prop:holder}.

	\begin{proof}[Proof of Proposition~\ref{prop:holder}]
		Recall the definition of $\cexp{holder}$ in Proposition~\ref{prop:tildeudiffunif}, and take $\gamma\in (0,\cexp{holder}\wedge 1)$.
		Take $n\in \Z$ and $t\ge 0$.
		For $x\in [n,n+1]$ and $z\in [-1,1]$, we have either $\{x,x+z\}\subset [n,n+2]$ or $\{x,x+z\}\subset [n-1,n+1]$.
		Hence by Proposition~\ref{prop:tildeudiffunif}, for any $x\in [n,n+1]$ and $z\in [-1,1]$, we have
		\begin{equation} \label{eq:vxzvx_holder}
			|v^N_t(x+z)-v^N_t(x)|\le \left(U^N_{N^{-1}\lfloor Nt\rfloor ,n,1}\vee U^N_{N^{-1}\lfloor Nt\rfloor ,n-1,1}\right)|z|^{\cexp{holder}}.
		\end{equation}
		Moreover, by Proposition~\ref{prop:tildeudiffunif}, almost surely $U^N_{N^{-1}k_1,k_2,1}<\infty$ for all $k_1\in \N_0$, $k_2\in \Z$.
		
		By the compact interface result in Lemma~\ref{lem:compactinterface}, for $T>0$ there exist random variables $R_T$ and $L_T$ with $R_T<\infty$ and $L_T>-\infty$ almost surely, such that for every $t\in [0,T]$,
		\begin{equation} \label{eq:compact_interface}
			v^N_t(x)=\begin{cases}
				0 \quad &\text{ for }x\ge R_T,\\
				1 \quad &\text{ for }x\le L_T.
			\end{cases}
		\end{equation}
		
		Define the event
		\begin{equation*}
			E=\left\{
				U^N_{N^{-1}k_1,k_2,1}<\infty \; \forall k_1\in \N_0,\, k_2\in \Z
				\right\}
			\cap \left\{ R_T<\infty, L_T>-\infty\right\}.
		\end{equation*}
		By~\eqref{eq:vxzvx_holder} and~\eqref{eq:compact_interface}, it follows that for $t\in [0,T]$ we have 
		\begin{align*}
			&\| v^N_t \|_{2,\alpha}^2 + \int_\R \int_{[-1,1]} \frac{|v^N_t(x + z) - v^N_t(x)|^2}{|z|^{1+2\gamma}} e^{\alpha x} dz \, dx\\
			&\le \int_{-\infty}^{R_T}e^{\alpha x} dx\\
			&+\sum_{n\in \Z\cap [L_T-2,R_T+1]}\int_{[n,n+1]}\int_{[-1,1]}
			\left(U^N_{N^{-1}\lfloor Nt\rfloor ,n,1}\vee U^N_{N^{-1}\lfloor Nt\rfloor ,n-1,1}\right)^2 |z|^{2\cexp{holder} - 2\gamma - 1 } e^{\alpha x} dz \, dx,
		\end{align*}
		which is finite on the event $E$, since we chose $\gamma<\cexp{holder}$.
		Hence there exists a random variable $ U_{N,T} $ taking values in $ (0,\infty) $ such that, on $ E $,
		\begin{equation*}
			\sup_{t \in [0,T]} \| v^N_t \|_{\hat{H}^{\gamma,\alpha}} \leq U_{N,T}.
		\end{equation*}
		Since $E$ occurs almost surely, the result follows.
	\end{proof}

	\begin{proof}[Proof of Proposition~\ref{prop:holder}]
	Recall the definition of $\cexp{holder}$ in Proposition~\ref{prop:tildeudiffunif}, and take $\gamma\in (0,\cexp{holder}\wedge 1)$.
		Suppose $n\in \Z$ and $t\ge 0$.
		For $x\in [n,n+1]$ and $z\in [-1,1]$, we have either $\{x,x+z\}\subset [n,n+2]$ or $\{x,x+z\}\subset [n-1,n+1]$.
		Hence by Proposition~\ref{prop:tildeudiffunif}, for any $x\in [n,n+1]$, $z\in [-1,1]$, we have
		\begin{equation} \label{eq:vxzvx_holder}
			|v^N_t(x+z)-v^N_t(x)|\le \left(U^N_{N^{-1}\lfloor Nt\rfloor ,n,2}\vee U^N_{N^{-1}\lfloor Nt\rfloor ,n-1,2}\right)|z|^{\cexp{holder}}.
		\end{equation}
		Recalling the definition of $\sigmaN{interface}$ from \eqref{eq:sigmainterface}, and using the bound $\cexp{interface} < 1$, it follows from the above that for $t \in [0, \sigmaN{interface}]$,
		\begin{align} \label{eq:vsobobd1}
			&\| v^N_t \|_{2,\alpha}^2 + \int_\R \int_{[-1,1]} \frac{|v^N_t(x + z) - v^N_t(x)|^2}{|z|^{1+2\gamma}} e^{\alpha x} dz \, dx \\
			&\le \int_{-\infty}^{2\log N }e^{\alpha x} dx \notag \\
			&+\sum_{n\in \Z\cap (-\infty, \lfloor 2 \log N \rfloor]}\int_n^{n+1}\int_{-1}^1
			\left(U^N_{N^{-1}\lfloor Nt\rfloor ,n,2}\vee U^N_{N^{-1}\lfloor Nt\rfloor ,n-1,2}\right)^2 |z|^{2\cexp{holder} - 2\gamma - 1 } e^{\alpha x} dz \, dx. \notag
		\end{align}
		Since $\gamma < \cexp{holder}$, we can integrate $z$ out of the second term, which contributes a multiplicative constant $C>0$. 
		Next, let $s \in \N_0 \cap [0,NT]$ and $\ell \in \N$, 
		\begin{equation} \label{eq:holdereventdef}
		E^\ell_s := \cap_{n\in \Z\cap (-\infty, \lfloor 2 \log N \rfloor]} \left\{ \left(U^N_{s ,n,2}\vee U^N_{s ,n-1,2}\right) \leq \ell \cdot e^{-\frac{\alpha n}{4} - \frac 1 2}  \right\}.
		\end{equation} 
		Then from \eqref{eq:vsobobd1} and the subsequent remark, we obtain that if $t \leq \sigmaN{interface}$, then on $E^{\ell}_{N^{-1} \lfloor  Nt \rfloor}$,
		\begin{align*} 
			&\sum_{n\in \Z\cap (-\infty, \lfloor 2 \log N \rfloor]}\int_n^{n+1}\int_{-1}^1
			\left(U^N_{N^{-1}\lfloor Nt\rfloor ,n,2}\vee U^N_{N^{-1}\lfloor Nt\rfloor ,n-1,2}\right)^2 |z|^{2\cexp{holder} - 2\gamma - 1 } e^{\alpha x} dz \, dx \notag
			\\ &\leq C \ell^2 \cdot  \sum_{n\in \Z\cap (-\infty, \lfloor 2 \log N \rfloor]} e^{\alpha n / 2} 
			\\ &\leq N^{2\alpha} \ell^2
		\end{align*}
		where the final inequality holds for sufficiently large $N$. In particular, from \eqref{eq:vsobobd1}, we obtain
		\begin{equation}
		\text{On $\{t \leq \sigmaN{interface} \} \cap E^{\ell}_{N^{-1} \lfloor  Nt \rfloor}$, } \|v^N_t\|_{\hat{H}^{\gamma,\alpha}} \leq N^{2\alpha} +  N^{2\alpha} \ell^2.
		\end{equation}
		We may therefore conclude that on $\cap_{s = 0}^{\lfloor NT \rfloor} E^\ell_s$, we have
		\[ \sup_{ t \in [0,\sigmaN{interface} \wedge T]}  \|v^N_t\|_{\hat{H}^{\gamma,\alpha}} \leq N^{2\alpha} (1+ \ell^2).\]
		In particular, for $N \geq T$, 
		\begin{align} \label{eq:soboholderbadevent}
		&\E{ \sup_{ t \in [0,\sigmaN{interface} \wedge T]}  \|v^N_t\|_{\hat{H}^{\gamma,\alpha}}}  \notag
		\\ &\qquad = \int_0^\infty \P{\sup_{ t \in [0,\sigmaN{interface} \wedge T]}  \|v^N_t\|_{\hat{H}^{\gamma,\alpha}} > \lambda }d\lambda \notag
		\\ &\qquad= N^{2\alpha} + 2N^{2\alpha}  \int_0^\infty \ell \cdot \P{\sup_{ t \in [0,\sigmaN{interface} \wedge T]}  \|v^N_t\|_{\hat{H}^{\gamma,\alpha}} > N^{2\alpha} (1+ \ell^2)  }  d\ell \notag
		\\ &\qquad\leq N^{2\alpha} + 2N^{2\alpha} \left(2 + \sum_{\ell=1}^\infty  \ell \cdot \P{ \sup_{ t \in [0,\sigmaN{interface} \wedge T]}  \|v^N_t\|_{\hat{H}^{\gamma,\alpha}} \leq N^{2\alpha}  (1+\ell^2)} \right)\notag 
		\\ &\qquad\leq N^{2\alpha} \left(3 + \sum_{\ell=1}^\infty  \ell \cdot \P{ \cup_{s = 0}^{\lfloor NT \rfloor} (E^\ell_s)^c } \right).
		\end{align}
		Recall the definition of $E^\ell_s$ from \eqref{eq:holdereventdef}. 
		To bound the probability appearing above, we use a union bound, then Proposition~\ref{prop:tildeudiffunif} with Markov's inequality, 
		\begin{align*}
		\P{ \cup_{s = 0}^{\lfloor NT \rfloor} (E^\ell_s)^c} &\leq \sum_{s=0}^{N^2 - 1} \bP{ (E^\ell_s)^c}
		\\ &\leq \sum_{s=0}^{\lfloor NT \rfloor} \sum_{n\in \Z\cap (-\infty, \lfloor 2 \log N \rfloor]} \P{ \left(U^N_{s ,n,2}\vee U^N_{s ,n-1,2}\right)  > e^{-n\alpha/2} \ell}
		\\ &\leq  2 (NT+1) \sum_{n\in \Z\cap (-\infty, \lfloor 2 \log N \rfloor]} A_2 e^{2\alpha n } \ell^{-4}
		\\ &\leq C (NT+1)  \ell^{-4}
		\end{align*}
		for some constant $C>0$. Substituting the above into \eqref{eq:soboholderbadevent}, we obtain that for another constant $C>0$,
		\[\E{ \sup_{ t \in [0,\sigmaN{interface} \wedge T]}  \|v^N_t\|_{\hat{H}^{\gamma,\alpha}}} \leq C (NT+1)N^{2\alpha}, \]
		which proves result.
%
%
%
%
%
%
	\end{proof}
	
	
	\section{Analysis of the deterministic flow} \label{sec:det_flow}
	
	The aim of this section is to prove Theorem~\ref{thm:zeta}.
	We first use estimates on the energy functional along the deterministic flow to show that $\mathrm{ dist}(\Phi_t(v,\cdot),M) $ tends to zero as $ t \to \infty $ and we conclude by showing that $ \eta(\Phi_t(v,\cdot)) $ also converges as $ t \to \infty $.
	Most of this section is loosely adapted from Section~9 of \cite{funaki_scaling_1995}.
	
	We start by stating some basic facts about the deterministic flow $ \Phi_t $. The proof of the following can be found in Appendix~\ref{sec:Phit}.
	
	\begin{proposition} \label{prop:Phit}
		The function $ (t,v) \mapsto \Phi_t(v) $ is continuous from $ \R_+ \times L^{2,\alpha} $ to $ L^{2,\alpha} $.
		In addition, for all $ T > 0 $ and $ \gamma \in [0,2) $, there exists a constant $ C_{T,\gamma} > 0 $ such that, for all $ v, v' \in L^{2,\alpha} $ and $ t \in (0,T] $,
		\begin{align} \label{bounds_Phit_Hgamma}
			\Hnorm{\Phi_t(v,\cdot)}{\gamma} \leq \frac{C_{T,\gamma}}{t^{\gamma/2}} \| v \|_{2,\alpha}, && \Hnorm{\Phi_t(v,\cdot) - \Phi_t(v',\cdot)}{\gamma} \leq \frac{C_{T,\gamma}}{t^{\gamma/2}} \| v - v' \|_{2,\alpha}.
		\end{align}
		Moreover, for any $ 0 \leq \gamma' < \gamma < 2 $ and $ T > 0 $, there exists $ C > 0 $ such that, for all $ v \in H^{\gamma,\alpha} $ and $ t \in [0,T] $,
		\begin{equation} \label{continuity_Phit_Hgamma}
			\Hnorm{\Phi_t(v,\cdot) - v}{\gamma'} \leq C t^{\frac{\gamma-\gamma'}{2}} \Hnorm{v}{\gamma}.
		\end{equation}
	\end{proposition}
	
	The proof of the following is straightforward, and we omit it.
	\begin{lemma} \label{lemma:Phit_Linfty}
		For any $\lambda \geq 0$ satisfying $\lambda \in (\alpha - 1, \alpha/2)$, for any $ v \in L^{2,\alpha}$ with $\|v\|_{\infty,\lambda} < \infty$, $ t \mapsto \Phi_t(v,\cdot) $ is continuous from with respect to $\|\cdot\|_{\infty,\lambda}$.
	\end{lemma}
	
	\subsection{The energy functional along the deterministic flow} \label{subsec:H_Phi}
	
	Let us start by recalling the following result from \cite{fife_approach_1977} and \cite{rothe_convergence_1981}.
	Recall the definition of the norm $ \| \cdot \|_{\infty,\lambda} $ in \eqref{eq:weightdefn}.
	
	\begin{lemma}[\cite{fife_approach_1977}, Lemma~4.1, \cite{rothe_convergence_1981}, Lemma~2] \label{lemma:sub_sup_solutions}
		For any $ K > 0 $ and $ \lambda \geq 0$ such that $\alpha-1 < \lambda < 1$, there exist $ \newEpsilon{rothe} > 0 $, $ \genericc > 0 $ and $\newCst{rothe} > 0$ such that, if $ \| v - m_\eta \|_{\infty,\lambda} \leq \Epsilon{rothe} $ for some $ \eta \in [-K,K] $, then there exist $ \eta_-, \eta_+ \in \R $ such that, for all $ t \geq 0 $ and $ x \in \R $,
		\begin{align*}
			m_{\eta_-(v)}(x) - \| v-m_\eta \|_{\infty,\lambda} e^{-\genericc t} (1 \vee e^{-\lambda x}) \leq \Phi_t(v,x) \leq m_{\eta_+(v)}(x) + \| v-m_\eta \|_{\infty,\lambda} e^{-\genericc t} (1 \vee e^{-\lambda x}),
		\end{align*}
		and $ | \eta_\pm(v) - \eta | \leq \Cst{rothe} \| v - m_{\eta} \|_{\infty,\lambda} $.
	\end{lemma}

	We now fix $ K > 0 $ and $\lambda \geq 0$ such that $\alpha-1 < \lambda < 1$, and take $ v \in \mathcal{V}_{\Beta{zeta}, K, \Epsilon{zeta}}^{(\lambda)} $ for $ \Beta{zeta} \in (0,\Beta{eta}(K)) $ and $ \Epsilon{zeta} \in (0,\Epsilon{rothe}) $ to be chosen later.
	Let us then define
	\begin{equation*}
		\newstime{dist} := \inf \lbrace t \geq 0 : \Phi_t(v,\cdot) \notin \mathcal{V}_{K+1} \rbrace.
	\end{equation*}
	Then, by Lemma~\ref{lemma:eta}, for all $ t \in [0, \stime{dist}) $, $ \eta(\Phi_t(v, \cdot)) $ is well defined.
	We then set, for $ t \in [0, \stime{dist}) $,
	\begin{align*}
		\eta_t(v) := \eta(\Phi_t(v,\cdot)), && s_t(v) := \Phi_t(v,\cdot) - m_{\eta_t(v)}.
	\end{align*}
	Recall the definition of $ \Epsilon{H} $ in Proposition~\ref{prop:energy_estimates} and define
	\begin{equation*}
		\newstime{supnorm} := \inf \lbrace t \in [0,\stime{dist}) : \| s_t(v) \|_{\infty,\lambda} > \Epsilon{H} \rbrace,
	\end{equation*}
	using the convention $ \inf \emptyset = + \infty $.
	For $ \newEpsilon{eta_t} > 0 $, also set
	\begin{align} \label{def:t_eta}
		\newstime{eta} := \inf \lbrace t \in [0,\stime{dist}) : | \eta_t(v) - \eta(v) | > \Epsilon{eta_t} \rbrace.
	\end{align}
	Finally, let us write
	\begin{equation*}
		t_{min} := \stime{dist} \wedge \stime{supnorm} \wedge \stime{eta}.
	\end{equation*}
	
	The following result then follows directly from Lemma~\ref{lemma:sub_sup_solutions}.
	
	\begin{cor} \label{cor:sup_norm}
		There exist $ \Epsilon{zeta} > 0 $ and $ \Epsilon{eta_t} > 0 $ such that, for all $ v \in \mathcal{V}_{\Beta{zeta}, K, \Epsilon{zeta}}^{(\lambda)} $ and $ t \geq 0 $,
		\begin{equation*}
			t_{min} \wedge t < \stime{supnorm}.
		\end{equation*}
	\end{cor}

	\begin{proof}
		From Lemma~\ref{lemma:sub_sup_solutions}, we deduce that there exists $ \newCst{rothe2} > 0 $ such that, if $ \| v- m_\eta \|_{\infty,\lambda} \leq \Epsilon{rothe} $ for some $ \eta \in [-K,K] $, then, for all $ t \geq 0 $,
		\begin{equation*}
			\| \Phi_t(v,\cdot) - m_{\eta(v)} \|_{\infty,\lambda} \leq \Cst{rothe2} \| v - m_\eta \|_{\infty,\lambda}.
		\end{equation*}
		(See for example Lemma~2 in \cite{rothe_convergence_1981}.)
		We thus choose $ \Epsilon{zeta} $ such that $ \Epsilon{zeta} < \Epsilon{rothe} $, so that we can apply Lemma~\ref{lemma:sub_sup_solutions} to $ v $.
		We obtain that, for all $ t \geq 0 $,
		\begin{equation*}
			\| \Phi_t(v,\cdot) - m_{\eta(v)} \|_{\infty,\lambda} \leq \Cst{rothe2}\, \Epsilon{zeta}.
		\end{equation*}
		Then, for $ t \in [0, \stime{dist} \wedge \stime{eta}) $,
		\begin{align*}
			\| s_t(v) \|_{\infty,\lambda} &\leq \| \Phi_t(v,\cdot) - m_{\eta(v)} \|_{\infty,\lambda} + \| m_{\eta(v)} - m_{\eta_t(v)} \|_{\infty,\lambda} \\
			&\leq \Cst{rothe2}\, \Epsilon{zeta} + \| \partial_x m \|_{\infty,\lambda} \, e^{\lambda | \eta(v) - \eta_t(v) |} | \eta_t(v) - \eta(v) | \\
			&\leq  \Cst{rothe2}\,\Epsilon{zeta} + \| \partial_x m \|_{\infty,\lambda}\, e^{\lambda \Epsilon{eta_t}} \Epsilon{eta_t}. \numberthis \label{bound_st_infty}
		\end{align*}
		We then choose $ \Epsilon{zeta} $ and $ \Epsilon{eta_t} $ small enough that
		\begin{equation} \label{choice_epsilons}
			 \Cst{rothe2}\, \Epsilon{zeta} + \| \partial_x m \|_{\infty,\lambda}\, e^{\lambda \Epsilon{eta_t}} \Epsilon{eta_t} < \Epsilon{H},
		\end{equation}
		By Proposition~\ref{prop:Phit}, Lemma~\ref{lemma:Phit_Linfty} and \eqref{eta_Lipschitz}, $ t \mapsto \| s_t(v) \|_{\infty,\lambda} $ is continuous on $ [0,\stime{dist}) $.
		Assume, aiming for a contradiction, that $ \stime{supnorm} \leq t_{min} \wedge t $.
		This is only possible if $ \stime{supnorm} < \stime{dist} $, so, in particular, $ t \mapsto \| s_t(v) \|_{\infty,\lambda} $ is continuous at $ t = \stime{supnorm} $.
		By \eqref{bound_st_infty} and \eqref{choice_epsilons} and the fact that $ \stime{supnorm} \leq \stime{dist} \wedge \stime{eta} $, this implies
		\begin{equation*}
			\| s_{\stime{supnorm}}(v) \|_{\infty,\lambda} < \Epsilon{H}.
		\end{equation*}
		On the other hand, since $ \stime{supnorm} \leq t < +\infty $, there exists a sequence $ (t_n, n \geq 1) $ approaching $ \stime{supnorm} $ from above such that $ \| s_{t_n}(v) \|_{\infty,\lambda} > \Epsilon{H} $, yielding a contradiction.
	\end{proof}
	
	The following result then provides a bound on $ \mathcal{H}(\Phi_t(v,\cdot)) $ for $ t > 0 $, even when $ v \notin H^{1,\alpha} $.
	
	\begin{lemma} \label{lemma:bounds_H}
		For all $ T > 0 $, there exists $ C_\mathcal{H}(T) > 0 $ such that for all $ v \in \mathcal{V} $,
		\begin{equation} \label{bound_H_small_times}
			\mathcal{H}(\Phi_t(v,\cdot)) \leq \frac{C_\mathcal{H}(T)}{t} \dist(v,M)^2, \quad \forall\, t \in (0,T].
		\end{equation}
		Moreover, there exists $ \genericc > 0 $ such that, for any $ v \in \mathcal{V}_{\Beta{zeta},K,\Epsilon{zeta}} $ and any $ t, t' \in (0, \stime{dist} \wedge \stime{supnorm}] $ with $ t \leq t' $,
		\begin{equation} \label{bound_H_large_times}
			\mathcal{H}(\Phi_{t'}(v,\cdot)) \leq e^{-\genericc(t'-t)} \mathcal{H}(\Phi_t(v,\cdot)).
		\end{equation}
	\end{lemma}
	
	\begin{proof}
		For the first bound, apply \eqref{expansion_H} with $ v = m_\eta $ and $ h = \Phi_t(v,\cdot) - m_\eta $.
		Then, by \eqref{lemma:H_is_zero_on_M} and \eqref{def:m},
		\begin{equation*}
			\mathcal{H}(\Phi_t(v,\cdot)) = -\frac{1}{2} \langle \mathcal{A}_\eta (\Phi_t(v,\cdot) - m_\eta), \Phi_t(v,\cdot) - m_\eta \rangle_\alpha + \int_\R \mathcal{U}(m_\eta(x), \Phi_t(v,x) - m_\eta(x)) e^{\alpha x} dx.
		\end{equation*}
		By the definition of $ \mathcal{A}_\eta $ and \eqref{IPP_alpha}, we obtain
		\begin{multline*}
			\mathcal{H}(\Phi_t(v,\cdot)) = \frac{1}{2} \| \partial_{x} (\Phi_t(v,\cdot)-m_\eta) \|_{2,\alpha}^2 - \frac{1}{2} \int_\R f'(m_\eta(x)) (\Phi_t(v,x) - m_\eta(x))^2 e^{\alpha x} dx \\ + \int_\R \mathcal{U}(m_\eta(x), \Phi_t(v,x) - m_\eta(x)) e^{\alpha x} dx.
		\end{multline*}
		Then, using the fact that $ f' $ is bounded and that $ |U(v,h)| \leq C |h|^2 $, we see that there exists a constant $ C > 0 $ such that, for any $ \eta \in \R $,
		\begin{equation*}
			\mathcal{H}(\Phi_t(v,\cdot)) \leq C \| \Phi_t(v,\cdot) - m_\eta \|_{H^{1,\alpha}}^2.
		\end{equation*}
		The bound \eqref{bound_H_small_times} then follows from \eqref{bounds_Phit_Hgamma} with $ \gamma = 1 $, $ v = \Phi_t(v,\cdot) $ and $ v' = m_{\eta(v)} $, using the fact that $ m_\eta $ is a fixed point of $ \Phi_t $.
		
		Let us now prove \eqref{bound_H_large_times}.
		Note that by~\eqref{def_DH}, for $t > 0$,
		\begin{align} \label{eq:dtH2alpha}
			\partial_t \mathcal{H}(\Phi_t(v,\cdot)) &= \langle D \mathcal{H}(\Phi_t(v,\cdot),\cdot), \partial_t \Phi_t(v,\cdot) \rangle_{\alpha} 
				= - \| \partial_t \Phi_t(v,\cdot) \|_{2,\alpha}^2.
		\end{align}
		Now take $v\in \mathcal{V}_{\Beta{eta},K,\Epsilon{zeta}} $, and $ t, t' \in (0, \stime{dist} \wedge \stime{supnorm}] $.
		Then, by \eqref{eq:dtH2alpha}, for all $ s \in (t, t') $,
		\begin{align*}
			\partial_s \mathcal{H}(\Phi_s(v,\cdot)) = - \| D \mathcal{H}(\Phi_s(v,\cdot),\cdot) \|_{2,\alpha}^2.
		\end{align*}
		Since $ s < \stime{dist} \wedge \stime{supnorm} $ and $\| \cdot \|_{\infty} \leq \| \cdot \|_{\infty,\lambda} $, $ \Phi_s(v, \cdot) $ satisfies the assumptions of Proposition~\ref{prop:energy_estimates} with $ \eta = \eta_s(v) $.
		As a result, for all $ s \in (t,t') $,
		\begin{align*}
			\| D\mathcal{H}(\Phi_{s}(v,\cdot),\cdot) \|_{2,\alpha}^2 \geq \frac{\Cst{H_below}}{\Cst{H_above}} \mathcal{H}(\Phi_{s}(v,\cdot)).
		\end{align*}
		Integrating between $ t $ and $ t' $ yields
		\begin{align} \label{eq:HPhi}
			\mathcal{H}(\Phi_{t'}(v,\cdot)) \leq e^{-\frac{\Cst{H_below}}{\Cst{H_above}}(t'-t)} \mathcal{H}(\Phi_t(v,\cdot)),
		\end{align}
		which proves \eqref{bound_H_large_times}.
	\end{proof}

	\subsection{Convergence to the stable manifold, proof of Theorem~\ref{thm:zeta}} \label{subsec:proof_zeta}

	The aim of this section is to prove Theorem~\ref{thm:zeta}, namely the exponential convergence  of $ \Phi_t(v,\cdot) $ to some $ m_{\zeta} $ as $ t \to \infty $ in $ L^{2,\alpha} $.
	Following Section~9 of \cite{funaki_scaling_1995}, we first prove that $ \| s_t(v) \|_{2,\alpha} \to 0 $ as $ t \to \infty $ (Lemma~\ref{lemma:bounds_st} below).
	Since
	\begin{align*}
		\Phi_t(v,\cdot) - m_{\zeta} = s_t(v) + m_{\eta_t(v)} - m_{\zeta},
	\end{align*}
	it then remains to show that $ \eta_t(v) \to \zeta(v) $ for some $ \zeta(v) \in \R $, which is done with the help of Proposition~\ref{prop:D_eta} below.

	\begin{lemma} \label{lemma:bounds_st}
		For any $ K > 0 $, there exist  $ \newCst{cvg_st} > 0 $, $ \genericc > 0 $, $ \Epsilon{zeta} > 0 $, $ \Epsilon{eta_t} > 0 $ and $ \Beta{zeta} \in (0,\Beta{eta}) $ such that, for any $ v \in \mathcal{V}_{\Beta{zeta}, K, \Epsilon{zeta}} $ and $ t \geq 1 $,
		\begin{equation*}
			1 \leq \stime{eta} \wedge t < \stime{dist} \wedge \stime{supnorm},
		\end{equation*}
		and
		\begin{align} \label{bound_st_large_times}
			\| s_t(v) \|_{H^{1,\alpha}} \leq \Cst{cvg_st} \, e^{-\genericc t} \dist(v,M), \quad \forall\, t \in [1,t_{min}).
		\end{align}
		and
		\begin{align} \label{bound_st_large_times_L2}
			\| s_t(v) \|_{2,\alpha} \leq \Cst{cvg_st}\, e^{-\genericc t} \dist(v,M), \quad \forall t \in [0,t_{min}).
		\end{align}
	\end{lemma}

	\begin{proof}
		We start by showing that $ t_{min} \geq 1 $.
		By \eqref{bounds_Phit_Hgamma} in Proposition~\ref{prop:Phit} with $ \gamma = 0 $, for $ t \in [0,1] $,
		\begin{equation} \label{bound_Phi_t_m_eta_v}
			\| \Phi_t(v,\cdot) - m_{\eta(v)} \|_{2,\alpha} \leq C_{1,0} \| v - m_{\eta(v)} \|_{2,\alpha} \leq C_{1,0} \Beta{zeta}.
		\end{equation}
		In addition, by \eqref{bound_eta_eta_0},
		\begin{equation*}
			|\eta(v)| \leq K + \Cst{eta}(K) \Beta{zeta}.
		\end{equation*}
		Choosing $ \Beta{zeta} $ such that $ C_{1,0} \Beta{zeta} < \Beta{eta}(K+1) $ and $ \Cst{eta}(K) \Beta{zeta} < 1 $ ensures that, for all $ t \in [0,1] $, $ \Phi_t(v,\cdot) \in \mathcal{V}_{K+1} $, so $ \stime{dist} \geq 1 $.
		Moreover, by \eqref{bound_eta_eta_0}, for $ t \in [0,1] $,
		\begin{equation*}
			| \eta_t(v) - \eta(v) | \leq \Cst{eta}(K+1) \| \Phi_t(v,\cdot) - m_{\eta(v)} \|_{2,\alpha} \leq \Cst{eta}(K+1) C_{1,0} \Beta{zeta}.
		\end{equation*}
		Thus, choosing $ \Beta{zeta} $ such that $ \Cst{eta}(K+1) C_{1,0} \Beta{zeta} < \Epsilon{eta_t} $, we ensure that $ \stime{eta} \geq 1 $. By Corollary~\ref{cor:sup_norm}, we conclude that $ t_{min} \geq 1 $.
		
		Then, by \eqref{bound_H_small_times} in Lemma~\ref{lemma:bounds_H} with $ T = 1 $, we have
		\begin{align} \label{bound_H_Phi_1}
			\mathcal{H}(\Phi_1(v,\cdot)) \leq C_\mathcal{H}(1)\, \dist(v,M)^2 \leq C_\mathcal{H}(1) \Beta{zeta}^2.
		\end{align}
		In addition, by \eqref{bound_H_large_times}, for all $ t \in [1, \stime{dist} \wedge \stime{supnorm}] $,
		\begin{equation} \label{bound_H_Phit}
			\mathcal{H}(\Phi_t(v,\cdot)) \leq e^{-\genericc (t-1)} \mathcal{H}(\Phi_1(v,\cdot)) \leq C_\mathcal{H}(1) \Beta{zeta}^2.
		\end{equation}
		Assume now, aiming for a contradiction, that $ \stime{dist} \leq t \wedge \stime{supnorm} \wedge \stime{eta} $.
		Then, for all $ t \in [1, \stime{dist}) $, we can apply Proposition~\ref{prop:energy_estimates} to $ \Phi_t(v,\cdot) $ with $ \eta = \eta_t(v) $, yielding
		\begin{equation} \label{bound_st_H}
			\| s_t(v) \|_{H^{1,\alpha}}^2 \leq \frac{1}{\Cst{H_below}} \mathcal{H}(\Phi_t(v,\cdot)) \leq \frac{C_K(1) \Beta{zeta}^2}{\Cst{H_below}}.
		\end{equation}
		In addition, for all $ t \in [1, \stime{dist}) $, by \eqref{bound_eta_eta_0},
		\begin{equation*}
			| \eta_t(v) | \leq K + \Cst{eta}(K) \Beta{zeta} + \Epsilon{eta_t}.
		\end{equation*}
		Now if $ |t - \stime{dist}| \leq 1 $, by the flow property and by \eqref{bounds_Phit_Hgamma} with $ \gamma = 0 $,
		\begin{equation*}
			\| \Phi_{\stime{dist}}(v,\cdot) - m_{\eta_t(v)} \|_{2,\alpha} \leq C_{1,0} \| \Phi_{t}(v,\cdot) - m_{\eta_t(v)} \|_{2,\alpha}.
		\end{equation*}
		By \eqref{bound_st_H}, the right-hand side is bounded from above by $ C_{1,0} \sqrt{C_\mathcal{H}(1) / \Cst{H_below}} \Beta{zeta} $.
		Hence, if we choose $ \Beta{zeta} $ and $ \Epsilon{eta_t} $ such that $ \Cst{eta}(K) \Beta{zeta} + \Epsilon{eta_t} \leq 1 $ and $ C_{1,0} \sqrt{C_\mathcal{H}(1) / \Cst{H_below}} \Beta{zeta} < \Beta{eta}(K+1) $, we ensure that $ \Phi_{\stime{dist}}(v,\cdot) \in \mathcal{V}_{K+1} $. Since $ \mathcal{V}_{K+1} $ is an open subset of $ L^{2,\alpha} $ and $ t \mapsto \Phi_t(v,\cdot) $ is continuous in $ L^{2,\alpha} $, this is a contradiction.
		We have thus obtained that, for any $ t \geq 1 $, $ t \wedge \stime{supnorm} \wedge \stime{eta} < \stime{dist} $.
		By Corollary~\ref{cor:sup_norm}, this implies $ t \wedge \stime{eta} < \stime{dist} \wedge \stime{supnorm} $.
		
		To conclude, combining \eqref{bound_H_Phi_1}, \eqref{bound_H_Phit} and \eqref{bound_st_H}, we obtain that there exist $ C > 0 $ and $ \genericc > 0 $, inedependent of $ v $ such that, for all $ t \in [1, t_{min}) $,
		\begin{equation*}
			\| s_t(v) \|_{H^{1,\alpha}} \leq C e^{-\genericc t} \dist(v,M),
		\end{equation*}
		yielding \eqref{bound_st_large_times}.
		Combined with \eqref{bound_Phi_t_m_eta_v}, this yields \eqref{bound_st_large_times_L2}.
	\end{proof}
	
	To prove Theorem~\ref{thm:zeta}, it remains to show that $ \eta_t(v) $ converges as $ t \to \infty $ and that $ \stime{eta} = + \infty $.
	This is done with the help of the following lemma.
	We define, for $ v \in \mathcal{V} $,
	\begin{align}
		\Vphi{1}(v) &:= - \| \partial_{x} m_{\eta(v)} \|_{2,\alpha}^2 + \langle v - m_{\eta(v)}, \partial_{xx} m_{\eta(v)} \rangle_\alpha \label{def_phi_tilde} \\
		\Vphi{2}(v) &:= \frac{3\alpha}{2} \| \partial_{x} m_{\eta(v)} \|_{2,\alpha}^2 + \langle v - m_{\eta(v)}, \partial_{x}^3 m_{\eta(v)} \rangle_\alpha.
	\end{align}

	\begin{proposition} \label{prop:D_eta}
		The map $ \eta : \mathcal{V} \subset L^{2,\alpha} \to \R $ is twice Fréchet differentiable and its derivatives are given by
		\begin{align*}
			D \eta(v) h = \langle D \eta(v), h \rangle_\alpha, && D^2 \eta(v) (h,h) = \int_{\R} D^2\eta(v,y_1,y_2) h(y_1) h(y_2) e^{\alpha(y_1 + y_2)} dy_1 dy_2,
		\end{align*}
		with
		\begin{align} \label{Deta}
			D \eta (v,y) = \frac{\partial_{x} m_{\eta(v)}(y)}{\Vphi{1}(v)},
		\end{align}
		and
		\begin{multline*}
			D^2 \eta(v,y_1,y_2) = \frac{\Vphi{2}(v)}{\Vphi{1}(v)^3} \partial_{x} m_{\eta(v)}(y_1) \partial_{x} m_{\eta(v)}(y_2) \\ - \frac{1}{\Vphi{1}(v)^2} \left( \partial_{x}^2 m_{\eta(v)}(y_1) \partial_{x} m_{\eta(v)}(y_2) + \partial_{x} m_{\eta(v)}(y_1) \partial_{x}^2 m_{\eta(v)}(y_2) \right).
		\end{multline*}
		For any $ K > 0 $, there exist $ \newCst{varphi_sup}(K) > 0 $ and $ \newCst{varphi_inf}(K) > 0 $ such that, for all $ v \in \mathcal{V}_{K} $,
		\begin{equation} \label{bound_varphi_tilde} 
			|\Vphi{k}(v) | \leq \Cst{varphi_sup}(K), \quad k \in \lbrace 1, 2 \rbrace,
		\end{equation}
		and 
		\begin{equation}
			|\Vphi{1}(v) | \geq \Cst{varphi_inf}(K). \label{bound_varphi_tilde_below}
		\end{equation}
		Moreover, for all $ v, v' \in \mathcal{V}_{K} $ and for $ k \in \lbrace 1, 2 \rbrace $,
		\begin{equation} \label{varphi_Lipschitz}
			\abs{\Vphi{k}(v) - \Vphi{k}(v')} \leq \Cst{varphi_sup} \| v - v' \|_{2,\alpha}.
		\end{equation}
		In addition, for any $ K > 0 $, $ \eta : \mathcal{V}_K \to \R $, $ D\eta : \mathcal{V}_K \to L^{2,\alpha} $ and $ D^2 \eta : \mathcal{V}_K \to L(L^{2,\alpha}, L^{2,\alpha}) $ are Lipschitz continuous.
	\end{proposition}												

	We prove Proposition~\ref{prop:D_eta} in Appendix~\ref{sec:fermi}.
	For now, let us conclude the proof of Theorem~\ref{thm:zeta}.
	
	\begin{proof}[Proof of Theorem~\ref{thm:zeta}]
		By Proposition~\ref{prop:D_eta}, for all $ t \in [0,t_{min}) $,
		\begin{align*}
			\deriv{\eta_t(v)}{t} &= \langle D\eta(\Phi_t(v,\cdot),\cdot), \partial_t \Phi_t(v,\cdot) \rangle_{\alpha} \\
			&= \frac{\langle \partial_{xx} \Phi_t(v,\cdot) + \alpha \partial_{x} \Phi_t(v,\cdot) + f(\Phi_t(v,\cdot)), \partial_{x} m_{\eta_t(v)} \rangle_\alpha}{\Vphi{1}(\Phi_t(v,\cdot))}. \numberthis \label{deta_t/dt}
		\end{align*}
%
		Since $ \partial_{xx} m_\eta + \alpha \partial_x m_\eta + f(m_\eta) = 0 $ for any $\eta\in \R$ and by \eqref{IPP_alpha}, we can write, for $ t < \stime{dist} $
		\begin{multline}\label{eq:phikbound}
			| \langle \partial_{xx} \Phi_t(v,\cdot) + \alpha \partial_{x} \Phi_t(v,\cdot) + f(\Phi_t(v,\cdot)), \partial_{x} m_{\eta_t(v)} \rangle_\alpha | \\
			\begin{aligned}
				&\leq \abs{ \langle \Phi_t(v,\cdot) - m_{\eta_t(v)}, (\partial_{xx} + \alpha \partial_x) \partial_{x} m_{\eta_t(v)} \rangle_{\alpha} } + \abs{ \langle f(\Phi_t(v,\cdot)) - f(m_{\eta_t(v)}), \partial_{x} m_{\eta_t(v)} \rangle_{\alpha} } \\
				&\leq C e^{\frac{\alpha}{2} \eta_t(v)} \| s_t(v) \|_{2,\alpha},
			\end{aligned}
		\end{multline}
		for some constant $ C > 0 $, where the second line follows by the Cauchy-Schwarz inequality, noting that, for $k \in \N$,
		\begin{align} \label{nabla_k}
			\| \partial_{x}^{k} m_\eta \|_{2,\alpha} = e^{\frac{\alpha}{2} \eta} \| \partial_{x}^{k} m \|_{2,\alpha}.
		\end{align}
		In addition, for $ t \in [0,t_{min}) $, 
		\begin{align} \label{bound_eta}
			\eta_t(v) \leq \eta(v) + \Epsilon{eta_t} \leq K + \Cst{eta}(K) \Beta{zeta} + \Epsilon{eta_t}.
		\end{align}
		Hence there exists $ C > 0 $ and $ \genericc > 0 $ (depending on $ K $) such that, for all $ v \in \mathcal{V}_{\Beta{zeta},K,\Epsilon{zeta}}^{(\lambda)} $, by \eqref{bound_st_large_times_L2} in Lemma~\ref{lemma:bounds_st},  for $ t \in [0,t_{min}) $,
		\begin{align*}
			\abs{\langle \partial_{xx} \Phi_t(v,\cdot) + \alpha \partial_{x} \Phi_t(v,\cdot) + f(\Phi_t(v,\cdot)), \partial_{x} m_{\eta_t(v)} \rangle_\alpha} \leq C \, \dist(v,M) \, e^{-\genericc t}.
		\end{align*}
		Plugging this in \eqref{deta_t/dt} and using \eqref{bound_varphi_tilde_below}, we obtain, for all $ t \in [0,t_{min}) $,
		\begin{align} \label{eq:detabound}
			\abs{ \deriv{\eta_t(v)}{t} } \leq \frac{C}{\Cst{varphi_inf}(K+1)} \dist(v,M) \, e^{-\genericc t}.
		\end{align}
		Assume now, aiming for a contradiciton, that $ \stime{eta} < +\infty $.
		By Lemma~\ref{lemma:bounds_st}, $ \stime{eta} < \stime{dist} \wedge \stime{supnorm} $.
		The above yields
		\begin{equation*}
			|\eta_{\stime{eta}}(v) - \eta(v)| \leq \frac{C}{\Cst{varphi_inf}(K+1) \genericc} \dist(v,M).
		\end{equation*}
		Choosing $ \Beta{zeta} $ such that
		\begin{equation*}
			\frac{C}{\Cst{varphi_inf}(K+1) \genericc} \Beta{zeta} < \Epsilon{eta_t},
		\end{equation*}
		we obtain a contradiciton (using the fact that $ t \mapsto \eta(\Phi_t(v, \cdot)) $ is continuous).
		Hence $ \stime{eta} = +\infty $ and, by Lemma~\ref{lemma:bounds_st}, $ \stime{dist} = \stime{supnorm} = +\infty $.
		As a result,~\eqref{eq:detabound} holds for all $ t \geq 0 $.
		We can thus set
		\begin{align*}
			\zeta(v): = \eta(v) + \int_{0}^{\infty} \deriv{\eta_t(v)}{t} dt.
		\end{align*}
		Then, for all $ t \geq 0 $,
		\begin{align} \label{cvg_eta}
			\abs{\eta_t(v) - \zeta(v)}=\abs{\int_t^{\infty} \frac{d\eta_s(v)}{ds}ds} \leq \frac{C}{\Cst{varphi_inf}(K+1) \genericc} \dist(v,M) e^{-\genericc t}.
		\end{align}
		Hence we obtain \eqref{diff_eta_zeta} and $ \eta_t(v) \to \zeta(v) $ as $ t \to \infty $.
		
		Let us now conclude by proving that $ \Phi_t(v,\cdot) \to m_{\zeta(v)} $.
		By the definition of $ s_t(v) $ and Lemma~\ref{lemma:diff_m_eta},
		\begin{align*}
		\| \Phi_t(v,\cdot) - m_{\zeta(v)} \|_{2,\alpha} &\leq \| s_t(v) \|_{2,\alpha} + \| m_{\eta_t(v)} - m_{\zeta(v)} \|_{2,\alpha} \\
		&\leq \| s_t(v) \|_{2,\alpha} + C \abs{\eta_t(v) - \zeta(v)} e^{\frac{\alpha}{2} (\eta_t(v) \vee \zeta(v))}. \numberthis \label{bound_Phi_t}
		\end{align*}
		Note that $\eta_t(v) \vee \zeta(v) \leq K + \Cst{eta}(K) \Beta{zeta} + \Epsilon{eta_t}$.
		Hence, plugging~\eqref{cvg_eta} into \eqref{bound_Phi_t} and using~\eqref{bound_st_large_times_L2} for the first term, we obtain that there exists $C>0$ such that, for all $ v \in \mathcal{V}_{\Beta{zeta},K,\Epsilon{zeta}} $ and $t\ge 0$,
		\begin{align*}
			\| \Phi_t(v,\cdot) - m_{\zeta(v)} \|_{2,\alpha} \leq C \dist(v,M) e^{-\genericc t}.
		\end{align*}
		Using \eqref{bound_st_large_times} and the same reasoning, we obtain the bound in the $ H^{1,\alpha} $ norm for $ t \geq 1 $ (using the second part of Lemma~\ref{lemma:diff_m_eta} for the second term in the first line of \eqref{bound_Phi_t}).
	\end{proof}

	Combining this result with Lemma~\ref{lemma:sub_sup_solutions}, we obtain the convergence in the uniform topology, which will be used in Section~\ref{sec:frechet}.

	\begin{cor} \label{cor:uniform_cvg}
		For any $\lambda \geq 0$ such that $\lambda \in (\alpha - 1,\alpha/2)$, there exist constants $ C > 0 $ and $ \genericc > 0 $ such that, for all $ v \in \mathcal{V}_{\Beta{zeta},K,\Epsilon{zeta}} $ and $ t \geq 0 $,
		\begin{align*}
			\| \Phi_t(v,\cdot) - m_{\zeta(v)} \|_{\infty,\lambda} \leq C \left(\| s(v) \|_{\infty,\lambda} + \dist(v,M) \right) e^{-\genericc t}.
		\end{align*}
	\end{cor}
	\begin{proof}
	We start by proving the result for $ t \in [0,1] $.
	By Lemma~\ref{lemma:sub_sup_solutions} (with $ \eta = \eta(v) $), for $ t \geq 0 $,
	\begin{multline} \label{bound_Phi_m_zeta}
		| \Phi_t(v,x) - m_{\zeta(v)}(x) | \leq | m_{\eta(v)}(x) - m_{\zeta(v)}(x) | +(1 \vee e^{-\lambda x}) \| v - m_{\eta(v)} \|_{\infty,\lambda} e^{-\genericc t} 
		\\ + \left( \| m_{\eta(v)}(x) - m_{\eta_-(v)}(x) | \vee | m_{\eta(v)}(x) - m_{\eta_+(v)}(x) | \right).
	\end{multline}
	From Lemma~\ref{lemma:sub_sup_solutions} we have $|\eta(v) - \eta_{\pm}(v)| \leq C \|v - m_{\eta(v)}\|_{\infty,\lambda}$, and we also recall that
	\begin{equation} \label{inequalities_zetaeta}
	|\zeta(v)  - \eta(v)| \leq C \dist(v,M), \qquad |m_\eta - m_{\eta'}| \leq C |\eta - \eta'| . 
	\end{equation}
	Hence, from \eqref{bound_Phi_m_zeta} we obtain
	\begin{align*}
			| \Phi_t(v,x) - m_{\zeta(v)}(x) | \leq C( \dist(v,M) + \| s(v)\|_{\infty,\lambda}).
	\end{align*}
	This proves the result for $t \in [0,1]$.
	
	Next, by the Sobolev inequality (see 	\citep[Theorem~8.8]{brezis_functional_2011}), there exists a constant $ C > 0 $ such that, for all $ s \in H^{1,\alpha} $,
		\begin{align*}
			\sup_{x \in \R} \abs{ s(x) e^{\frac{\alpha}{2} x} } \leq C \| s \|_{H^{1,\alpha}}.
		\end{align*}
		As a result, by Theorem~\ref{thm:zeta}, for $ t \geq 1 $ and any $x \in \R$,
		\begin{align*}
			| \Phi_t(v,x) - m_{\zeta(v)}(x) | \leq C \dist(v,M) e^{-\genericc t} e^{-\frac{\alpha}{2} x}.
		\end{align*}
	For $t\geq 1$ we fix $A_t \geq 0$, which will be assigned a value shortly. Then by the above, for $x \geq -A_t$, 
		\begin{align} \label{eq:belowAt}
			| \Phi_t(v,x) - m_{\zeta(v)}(x) | (1 \vee e^{\lambda x})  \leq C(1 \vee e^{\lambda x}) \dist(v,M) e^{-\genericc t} e^{-\frac{\alpha}{2} x} \notag
			\\  \leq C \dist(v,M) e^{-\genericc t} e^{-\frac{\alpha}{2} A_t}.
		\end{align}
	From \eqref{eq:masympneg}, we have $\sup_{x \in \R} |\partial_x m(x)| e^{-\lambda_- x} < \infty$, where $\lambda_- > 0$ is as in \eqref{eq:lambda-def}.
	We may then argue using \eqref{bound_Phi_m_zeta} (with the fact that $\zeta(v) = \zeta(\Phi_t(v,\cdot))$), \eqref{inequalities_zetaeta} and \eqref{diff_eta_zeta} that for all $ x \in \R $,
		\begin{align*}
			| \Phi_t(v,x) - m_{\zeta(v)}(x) | \leq Ce^{\lambda_- x}\dist(v,M) + C(e^{\lambda_- x} + e^{-\genericc t}) \|s(v)\|_{\infty,\lambda}.
		\end{align*}
	Hence, for $x \leq -A_t$, 
	\begin{align} \label{eq:aboveAt}
			| \Phi_t(v,x) - m_{\zeta(v)}(x) | \leq Ce^{- \lambda_- A_t}\dist(v,M) + C(e^{- \lambda_- A_t} + e^{-\genericc t}) \|s(v)\|_{\infty,\lambda}.
	\end{align}
	Now let $A_t:= \tfrac{c}{\lambda_- + \alpha/2} t$. It then follows from \eqref{eq:belowAt} and \eqref{eq:aboveAt} that the desired claim holds (with some smaller constant $c' < c$).
	\end{proof}

	\section{Fr\'echet derivatives of the Katzenberger coordinate} \label{sec:frechet}
	
	The aim of this section is to prove Proposition~\ref{prop:Deta_t}, Lemma~\ref{lemma:bounds_Detat_fixed_time}, Proposition~\ref{prop:Dzeta}, Lemma~\ref{lemma:Dzeta_m}, and Lemma~\ref{lemma:katzenberger}, which are all concerned with the properties of the Fréchet derivatives of $ v \mapsto \eta_t(v) $ and their limit as $ t \to \infty $.
	These estimates are all based on properties of the semigroup generated by the linearised equation \eqref{def_pst} and its fundamental solution $ p_{s,t}(v,\cdot,\cdot) $.
	Subsection~\ref{subsec:semigroup} introduce this (two parameter) semigroup and states various bounds on this semigroup wihch are valid over a fixed time horizon.
	These estimates are used in Subsection~\ref{subsec:Deta_t} to prove Proposition~\ref{prop:Deta_t} and Lemma~\ref{lemma:bounds_Detat_fixed_time}.
	Subsection~\ref{subsec:Dzeta} starts by stating several technical lemmas which are shown to imply Proposition~\ref{prop:Dzeta}, Lemma~\ref{lemma:Dzeta_m} and Lemma~\ref{lemma:katzenberger}.
	The proof of these technical lemmas requires several results on the long time behaviour of the linearised semigroup, which are stated in Subsection~\ref{subsec:semigroup_long_time}.
	Subsections~\ref{subsec:psi_k} and \ref{subsec:Ut_Upsilon} then show how the latter imply the technical lemmas stated at the beginning of Subsection~\ref{subsec:Dzeta}.
	The proofs of the all the results concerning the linearised semigroup are gathered in Section~\ref{sec:semigroup}.
	
	Although the present section is loosely based on what is done in Section~9 in \cite{funaki_scaling_1995}, our proofs depart from Funaki's in several places and the exponential weights that appear in the norm $ \| \cdot \|_{2,\alpha} $ also make them somewhat more involved.

	In this whole section, we consider that $ K > 0 $ is fixed but arbitrary, and that $ \Beta{zeta} > 0 $ and $ \Epsilon{zeta} > 0 $ are chosen so that the statement of Theorem~\ref{thm:zeta} is satisfied.
	
	\subsection{The linearised semigroup} \label{subsec:semigroup}
	
	Recall from~\eqref{def_pst} that $ p_{s,t}(v,x,\cdot) $ denotes the fundamental solution of $ \partial_t - \partial_{xx} - \alpha \partial_x - f'(\Phi_t(v,\cdot)) $, \textit{i.e.} $p_{s,t}$ is such that for any $s\ge 0$ and $h\in L^{2,\alpha}$,
	\begin{equation} \label{eq:utdefn}
	u_t(x) = \langle p_{s,t}(v,x,\cdot), h \rangle_{\alpha}, \quad t > s,\, x\in \R ,
	\end{equation}
	solves the equation
	\begin{equation} \label{linear_equation}
	\left\lbrace
	\begin{aligned}
	& \partial_t u_t = \partial_{xx} u_t + \alpha \partial_x u_t + f'(\Phi_t(v,\cdot)) u_t, \quad t > s, \\
	& \lim_{t \downarrow s} u_t = h.
	\end{aligned}
	\right.
	\end{equation}
	For $v\in L^{2,\alpha}$, define a two-parameter semigroup $ T_{s,t,v} $ on $ L^{2,\alpha} $ by setting
	\begin{equation} \label{eq:Tsemigpdefn}
	T_{s,t,v} \phi (y) := \langle p_{s,t}(v,\cdot,y), \phi \rangle_\alpha, \quad 0\le s \le t, \, y\in \R.
	\end{equation}
	This operator satisfies the semigroup property by \eqref{linear_equation}, i.e., noting that $u_t = T_{s,t,v}^* h$, for any $ 0 \leq s \leq r \leq t $,
	\begin{equation*}
		T_{r,t,v}^*\, T_{s,r,v}^* = T_{s,t,v}^*.
	\end{equation*}
	Taking the adjoint of the above, it follows that
	\begin{equation} \label{semigroup_ppty}
	T_{s,r,v}\, T_{r,t,v}=T_{s,t,v}
	\end{equation}
	Recall the definition of $\mathcal A_{\eta}$ from~\eqref{def_A}.
	When $ v = m_{\eta} $ for some $ \eta \in \R $, we also write $ T_{s,t,m_{\eta}} = e^{-(t-s) \mathcal{A}_{\eta}} $.
	Also recall from~\eqref{eq:phietadefn} that we set
	\begin{align*}
	\varphi_\eta = - \frac{\partial_{x} m_\eta}{\| \partial_{x} m_\eta \|_{2,\alpha}},
	\end{align*}
	and let $ P_{\lbrace \varphi_\eta \rbrace} = \langle \cdot, \varphi_\eta \rangle_\alpha \varphi_\eta $ denote the projection in $ L^{2,\alpha} $ onto the linear subspace spanned by $ \varphi_\eta $ and $ P_{\lbrace \varphi_\eta \rbrace^{\perp}} = I - P_{\lbrace \varphi_\eta \rbrace} $ the projection onto the orthogonal subspace.
	We denote the operator norm of a linear operator on $ X $ by $\|\cdot \|_{X \to X}$.
	Let us introduce the following norm
	\begin{equation} \label{def:norm_exp}
		\enorm{\phi} := \sup_{x \in \R} |\phi(x)| e^{|x|},
	\end{equation}
	and let $ L^{\infty,e} $ denote the subset of $ L^\infty $ on which this norm is finite.
	
	\begin{lemma} \label{lemma:projection_varphi}
		For all $ \zeta \in \R $,
		\begin{align*}
			\left\| P_{\lbrace \varphi_{\zeta} \rbrace^\perp} \right\|_{L^{2,\alpha} \to L^{2,\alpha}} \leq 2,
		\end{align*}
		and, for any $ K > 0 $, there exists $ C > 0 $ such that, for all $ |\zeta| \leq \Cst{sup_zeta}(K) $,
		\begin{equation*}
			\left\| P_{\lbrace \varphi_{\zeta} \rbrace^\perp} \right\|_{L^{\infty,e} \to L^{\infty,e}} \leq C.
		\end{equation*}
	\end{lemma}
	
	\begin{proof}
		By definition,
		\begin{equation*}
			P_{\lbrace \varphi_{\zeta} \rbrace^\perp} \phi = \phi - \langle \phi, \varphi_{\zeta} \rangle_\alpha \varphi_{\zeta}.
		\end{equation*}
		By the Cauchy-Schwarz inequality and the fact that $ \| \varphi_{\zeta} \|_{2,\alpha} = 1 $, we easily obtain that
		\begin{equation*}
			\left\| P_{\lbrace \varphi_{\zeta} \rbrace^\perp} \phi \right\|_{2,\alpha} \leq 2 \| \phi \|_{2,\alpha}.
		\end{equation*}
		In addition, by the definition of $ \enorm{\cdot} $,
		\begin{equation*}
			\enorm{P_{\lbrace \varphi_{\zeta} \rbrace^\perp} \phi} \leq \enorm{\phi} \left( 1 + \enorm{\varphi_{\zeta}} \int_\R | \varphi_{\zeta}(x) | e^{-|x|} dx \right).
		\end{equation*}
		The result then follows from the fact that $ | \varphi_{\zeta}(x) | \leq C e^{-|x-\zeta|} $ for some $ C > 0 $.
	\end{proof}
	
	We now state some results on $ T_{s,t,v} $ which hold for $ 0 \leq s \leq t \leq T $ for some arbitrary $ T $.
	They are all proved in Subsection~\ref{subsec:Tst_short_times}.
		
	\begin{lemma} \label{lemma:Tstv_Hgamma}
		For any $ T > 0 $ and $ \gamma \in [0,2) $, there exists $ C > 0 $ such that, for all $ \phi \in L^{2,\alpha} $, $ v \in L^{2,\alpha} $ and all $ 0 \leq s < t \leq T $,
		\begin{equation} \label{Tstv_Hgamma}
			\Hnorm{T_{s,t,v} \phi}{\gamma} \leq \frac{C}{((t-s)\wedge 1)^{\gamma/2}} \| \phi \|_{2,\alpha}.
		\end{equation}
		Moreover, for any $ T > 0 $ and $ \gamma \in [0,2) $, there exists $ C > 0 $ such that, for any $ \phi \in L^{\infty} $, $ v, v' \in L^{2,\alpha} $ and for all $ 0 \leq s \leq t \leq T $,
		\begin{equation} \label{Tstv-v'_Hgamma}
			\Hnorm{(T_{s,t,v} - T_{s,t,v'}) \phi}{\gamma} \leq C \| \phi \|_\infty \| v - v' \|_{2,\alpha}.
		\end{equation}
	\end{lemma}
	
	\begin{lemma} \label{lemma:Tstv_Lpq_short_times}
		For any $ T > 0 $, $ p \geq 1 $ and $ q \in \R $, there exists $ C > 0 $ such that, for all $ \phi \in L^{p,q} $, $ v \in L^{2,\alpha} $ and $ 0 \leq s \leq t \leq T $,
		\begin{equation} \label{Tstv_Lpq_short_time}
			\| T_{s,t,v} \phi \|_{p,q} \leq C \| \phi \|_{p,q},
		\end{equation}
		and for any $ \phi \in L^{p,q} \cap L^{2p,2q-p\alpha} $, $ v, v' \in L^{2,\alpha} $ and $ 0 \leq s \leq t \leq T $,
		\begin{equation} \label{diff_Tstv_Lpq}
			\| (T_{s,t,v} - T_{s,t,v'}) \phi \|_{p,q} \leq C \| v - v' \|_{2,\alpha} \| \phi \|_{2p,2q-p\alpha}.
		\end{equation}
	\end{lemma}
	
	\begin{lemma} \label{lemma:diff_pst2_short_times}
		For any $ T > 0 $ and $ q \in \R $, there exists $ C > 0 $ such that, for all $ \phi \in L^{1,q-\alpha} \cap L^{2,2q-3\alpha} $, $ v, v' \in L^{2,\alpha} $ and $ t \in [0,T] $,
		\begin{equation} \label{diff_p0t_v_v'}
			\int_\R \langle |\phi|, |p_{0,t}(v,\cdot,y) - p_{0,t}(v',\cdot,y) | p_{0,t}(v,\cdot,y) \rangle_\alpha e^{qy} dy \leq C t^{1/4} \| v - v' \|_{2,\alpha} \| \phi \|_{2,2q-3\alpha}.
		\end{equation}
		Moreover, for any $ T > 0 $ and $ q \in \R $, there exists $ C > 0 $ such that, for all $ \phi \in L^{1,q-\alpha} $, $ v \in L^{2,\alpha} $, $ \delta \in (0,1) $ and $ t \in [0,T] $,
		\begin{multline} \label{phi_p0t_delta}
			\int_{0}^{t} \int_{[-1,1]^2 \times \R} | \langle \phi, p_{0,s}(v,\cdot,y + \delta z_1) p_{0,s}(v,\cdot,y + \delta z_2) \rangle_\alpha \\ - \langle \phi, p_{0,s}(v,\cdot,y)^2 \rangle_\alpha | \rho(z_1) \rho(z_2) e^{qy} dz_1 dz_2 dy ds \\ \leq C \delta | \log(\delta) | \| \phi \|_{1,q-\alpha},
		\end{multline}
		where $ \rho : [-1,1] \to \R_+ $ is smooth and bounded.
	\end{lemma}

	\subsection{Fr\'echet derivatives of the finite time coordinate on the stable manifold} \label{subsec:Deta_t}
	
	The following result gives the Fr\'echet derivatives of $ v \mapsto \Phi_t(v,\cdot) $, which allows us to obtain those of $ \eta_t $ and prove Proposition~\ref{prop:Deta_t}.
	We then also prove Lemma~\ref{lemma:bounds_Detat_fixed_time}.
	Note that since $ T_{s,t,v}^* $ acts on $ L^{2,\alpha} $ and $ H $ is embedded in $ L^{2,\alpha} $, $ T_{s,t,v}^* $ also acts on $ H $.

	\begin{proposition} \label{prop:deriv_Phi}
		For any $t>0$, the map $ \Phi_t : H \to L^{2,\alpha} $ is twice Fréchet differentiable and, for $ h, h' \in H $,
		\begin{align}
		& D\Phi_t(v) h = T_{0,t,v}^* h, \label{DPhi} \\
		& D^2 \Phi_t(v)(h,h') = \int_{0}^{t} T_{s,t,v}^* \lbrace f''(\Phi_s(v,\cdot)) (T_{0,s,v}^* h) (T_{0,s,v}^* h') \rbrace ds. \label{D2Phi}
		\end{align}
		Moreover, $ \Phi_t : H \to L^{2,\alpha} $, $ D \Phi_t : H \to L(H,L^{2,\alpha}) $ and $ D^2 \Phi_t : H \to L^2(H, L^{2,\alpha}) $ are Lipschitz continuous.
	\end{proposition}

	We prove Proposition~\ref{prop:deriv_Phi} in Appendix~\ref{sec:Phit}.
	A similar result is proved in \cite{funaki_scaling_1995} (see in particular Lemma~9.7).
	Combining Proposition~\ref{prop:deriv_Phi} and Proposition~\ref{prop:D_eta}, we obtain Proposition~\ref{prop:Deta_t}, namely that the set $ U_t = \Phi_t^{-1}(\mathcal{V}) $ is an open subset of $ H $ and that $ v \mapsto \eta_t(v) $ is twice Fréchet differentiable on $ U_t $ with Lipschitz continuous derivatives on each $ \Phi_t^{-1}(\mathcal{V}_K) $.
	We can also compute the derivatives of $ \eta_t $ on $ U_t $ as follows.
	For $ k \geq 1 $, set
	\begin{equation} \label{def:psi_k}
		\Vpsi{k}(s,t,v,y) := T_{s,t,v} \partial_{x}^k m_{\eta_t(v)}(y).
	\end{equation}
	Also, set, for $ k \in \lbrace 1, 2 \rbrace $,
	\begin{equation*}
		\Vphi{k}(t,v) := \Vphi{k}(\Phi_t(v,\cdot)),
	\end{equation*}
	where $ \Vphi{k}(\cdot) $ was defined in \eqref{def_phi_tilde}.
	By the chain rule, for $ v \in U_t $ and $ h \in H $,
	\begin{equation*}
		D\eta_t(v)(h) = D\eta(\Phi_t(v))(D\Phi_t(v)(h)).
	\end{equation*}
	Substituting $ D\Phi_t(v)(h) = T_{0,t,v}^* h $ from Proposition~\ref{prop:deriv_Phi} and \eqref{Deta} from Proposition~\ref{prop:D_eta}, we obtain the representation
	\begin{equation*}
		D\eta_t(v) h = \langle D\eta_t(v), h \rangle_\alpha,
	\end{equation*}
	where
	\begin{equation} \label{Deta_t}
		D\eta_t(v,y) = \frac{\Vpsi{1}(0,t,v,y)}{\Vphi{1}(t,v)}.
	\end{equation}
	In the same way, for $ v \in U_t $ and $ h_1, h_2 \in H $,
	\begin{equation*}
		D^2 \eta_t(v)(h_1, h_2) = D^2 \eta(\Phi_t(v))(D\Phi_t(v,\cdot)(h_1), D\Phi_t(v,\cdot)(h_2)) + D\eta(\Phi_t(v))(D^2\Phi_t(v)(h_1, h_2)).
	\end{equation*}
	Using Proposition~\ref{prop:D_eta} and Proposition~\ref{prop:deriv_Phi}, this yields the representation
	\begin{equation*}
		D^2 \eta_t(v)(h_1, h_2) = \int_{\R^2} D^2 \eta_t(v,x,y) h_1(x) h_2(y) e^{\alpha(x + y)} dx dy,
	\end{equation*}
	with
	\begin{multline} \label{D2eta_t}
		D^2\eta_t(v, y_1, y_2) = \frac{\Vphi{2}(t,v)}{\Vphi{1}(t,v)^3} \Vpsi{1}(0,t,v,y_1) \Vpsi{1}(0,t,v,y_2) \\ - \frac{1}{\Vphi{1}(t,v)^2} \left( \Vpsi{1}(0,t,v,y_1) \Vpsi{2}(0,t,v,y_3) + \Vpsi{2}(0,t,v,y_1) \Vpsi{1}(0,t,v,y_2) \right) \\ + \frac{1}{\Vphi{1}(t,v)} \int_{0}^{t} \langle \Vpsi{1}(s,t,v,\cdot), f''(\Phi_s(v,\cdot)) p_{0,s}(v,\cdot,y_1) p_{0,s}(v,\cdot,y_2) \rangle_{\alpha} ds.
	\end{multline}
	We can now proceed with the proof of Lemma~\ref{lemma:bounds_Detat_fixed_time}.
	
	\begin{proof}[Proof of Lemma~\ref{lemma:bounds_Detat_fixed_time}]
		By Lemma~\ref{lemma:Tstv_Hgamma}, for any fixed $ t \geq 1 $ and $ \gamma \in [0,2) $, there exists $ C > 0 $ such that, for all $ v \in \Phi_t^{-1}(\mathcal{V}_K) $,
		\begin{equation} \label{Hnorm_psi1}
			\Hnorm{\psi_1(0,t,v,\cdot)}{\gamma} \leq C \| \partial_x m_{\eta_t(v)} \|_{2,\alpha} \leq C',
		\end{equation}
		for some $ C' > 0 $ (depending on $ t $) since $ \eta_t $ is bounded on $ \Phi_t^{-1}(\mathcal{V}_K) $.
		The bound \eqref{Detat_Hgamma} then follows from \eqref{bound_varphi_tilde_below}.
		In addition, using \eqref{bound_varphi_tilde_below} and \eqref{varphi_Lipschitz}, for any $ v, v' \in \Phi_t^{-1}(\mathcal{V}_K) $,
		\begin{equation} \label{diff_Detat}
			\Hnorm{D\eta_t(v) - D\eta_t(v')}{\gamma} \leq C \Hnorm{\psi_1(0,t,v,\cdot) - \psi_1(0,t,v',\cdot)}{\gamma} + C \| v - v' \|_{2,\alpha} \Hnorm{\psi_1(0,t,v',\cdot)}{\gamma}.
		\end{equation}
		We then write
		\begin{multline*}
			\Hnorm{\psi_1(0,t,v,\cdot) - \psi_1(0,t,v',\cdot)}{\gamma} \leq \Hnorm{(T_{0,t,v} - T_{0,t,v'}) \partial_{x} m_{\eta_t(v)}}{\gamma} \\ + \Hnorm{T_{0,t,v'} (\partial_{x} m_{\eta_t(v)} - \partial_{x} m_{\eta_t(v')})}{\gamma}.
		\end{multline*}
		Using Lemma~\ref{lemma:Tstv_Hgamma}, we obtain,
		\begin{equation*}
			\Hnorm{\psi_1(0,t,v,\cdot) - \psi_1(0,t,v',\cdot)}{\gamma} \leq C \| \partial_{x} m \|_\infty \| v-v' \|_{2,\alpha} + C \| \partial_{x} m_{\eta_t(v)} - \partial_{x} m_{\eta_t(v')} \|_{2,\alpha}.
		\end{equation*}
		Since $ \eta_t $ is Lipschitz continuous on $ \Phi_t^{-1}(\mathcal{V}_K) $, we conclude that there exists $ C >0 $ such that
		\begin{equation*}
			\Hnorm{\psi_1(0,t,v,\cdot) - \psi_1(0,t,v',\cdot)}{\gamma} \leq C \| v - v' \|_{2,\alpha}.
		\end{equation*}
		Plugging this and \eqref{Hnorm_psi1} in \eqref{diff_Detat} yields \eqref{diff_Detat_Hgamma}.
		In addition, by \eqref{Tstv_Lpq_short_time}, for any $ p \geq 1 $ and $ q \in \R $, there exist $ C > 0 $ such that, for all $ v \in U_t $ and $ 0 \leq s \leq t \leq T $,
		\begin{equation*}
			\| \psi_k(s,t,v,\cdot) \|_{p,q} \leq C \| \partial_{x}^k m_{\eta_t(v)} \|_{p,q}.
		\end{equation*}
		Thus, if $ |q| < p $, for any $ K > 0 $ there exists $ C > 0 $ such that, for all $ v \in \Phi_t^{-1}(\mathcal{V}_K) $,
		\begin{equation} \label{psi_k_Lpq}
			\| \psi_k(s,t,v,\cdot) \|_{p,q} \leq C,
		\end{equation}
		using the fact that $ \eta_t $ is bounded on $ \Phi_t^{-1}(\mathcal{V}_K) $.
		We thus obtain \eqref{Deta_t_Lpq} from \eqref{bound_varphi_tilde_below}.
		In addition, 
		\begin{multline*}
			\| \psi_k(s,t,v,\cdot) - \psi_k(s,t,v',\cdot) \|_{p,q} \leq \| (T_{s,t,v} - T_{s,t,v'}) \partial_{x} m_{\eta_t(v)} \|_{p,q} \\ + \| T_{s,t,v'} (\partial_{x} m_{\eta_t(v)} - \partial_{x} m_{\eta_t(v')}) \|_{p,q}.
		\end{multline*}
		Using the fact that $ \abs{\frac{q}{p} - \frac{\alpha}{2}} < 1 $ and $ |q| < p $, we see that $ \partial_{x} m \in L^{2p,2q-p\alpha} \cap L^{p,q} $.
		By Lemma~\ref{lemma:Tstv_Lpq_short_times}, we then obtain
		\begin{multline*}
			\| \psi_k(s,t,v,\cdot) - \psi_k(s,t,v',\cdot) \|_{p,q} \leq C \| v - v' \|_{2,\alpha} \| \partial_{x} m_{\eta_t(v)} \|_{2p,2q-p\alpha} \\ + C \| \partial_{x} m_{\eta_t(v)} - \partial_{x} m_{\eta_t(v')} \|_{p,q}.
		\end{multline*}
		Since $ \eta_t $ is bounded and Lipschitz continuous on $ \Phi_t^{-1}(\mathcal{V}_K) $ and $ \partial_{x x} m \in L^{p,q} $, we obtain that there exists $ C > 0 $ such that
		\begin{equation} \label{continuity_psi_k_Lpq}
			\| \psi_k(s,t,v,\cdot) - \psi_k(s,t,v',\cdot) \|_{p,q} \leq C \| v - v' \|_{2,\alpha}
		\end{equation}
		The bound \eqref{diff_Detat_Lpq} then follows from \eqref{continuity_psi_k_Lpq}, \eqref{psi_k_Lpq}, \eqref{bound_varphi_tilde}, \eqref{bound_varphi_tilde_below} and \eqref{varphi_Lipschitz}, using \eqref{Deta_t}.
		We now turn to the bounds on $ D^2\eta_t(v,\cdot,\cdot) $.
		The terms in the first two lines on the right-hand side of \eqref{D2eta_t} can be dealt with using \eqref{psi_k_Lpq} and \eqref{continuity_psi_k_Lpq}, \eqref{bound_varphi_tilde}, \eqref{bound_varphi_tilde_below} and \eqref{varphi_Lipschitz}, using the Cauchy-Schwarz inequality.
		For $ v \in \Phi_t^{-1}(\mathcal{V}) $, let us introduce the notation
		\begin{equation} \label{eq:Udefn}
			U_t(v,y_1,y_2) := \int_{0}^{t} \langle \Vpsi{1}(s,t,v,\cdot), f''(\Phi_s(v,\cdot)) p_{0,s}(v,\cdot,y_1) p_{0,s}(v,\cdot,y_2) \rangle_{\alpha} ds.
		\end{equation}
		We then aim to find $ C > 0 $ such that
		\begin{equation} \label{U_t_L1q}
			\int_\R | U_t(v,y,y) | e^{qy} dy \leq C,
		\end{equation}
		and
		\begin{equation} \label{continuity_Ut}
			\int_\R | U_t(v,y,y) - U_y(v',y,y) | e^{q y} dy \leq C \| v - v' \|_{2,\alpha}.
		\end{equation}
		The bounds \eqref{D2eta_t_L1q} and \eqref{diff_D2eta_t_L1q} then follow using \eqref{bound_varphi_tilde_below} and \eqref{varphi_Lipschitz}.
		From Lemma~\ref{lemma:bound_pst}, we obtain that, for $ s \in [0,T] $,
		\begin{align*}
			\int_\R \langle |\phi|, p_{0,s}(v,\cdot,y)^2 \rangle_\alpha e^{qy} dy &\leq \int_{\R^2} |\phi(x)| G_{t}(x-y)^2 e^{(q-\alpha)y} dy \\
			&\leq \frac{C}{\sqrt{s}} \| \phi \|_{1,q-\alpha}. \numberthis \label{pst2_short_time}
		\end{align*}
 		Hence, for $ |q-\alpha|<1 $, we can apply \eqref{psi_k_Lpq} and we obtain
 		\begin{align*}
 			\int_\R |U_t(v,y,y)| e^{qy} dy &\leq \int_{0}^{t} \frac{C}{\sqrt{s}} \| \psi_1(s,t,v,\cdot) \|_{1,q-\alpha} ds \\
 			&\leq C',
 		\end{align*}
 		for some $ C' > 0 $ depending only on $ T $ and $ q $.
 		Let us now write
 		\begin{multline} \label{diff_Ut}
 			| U_t(v,y,y) - U_t(v',y,y) | \leq \int_{0}^{t} \| f'' \|_\infty \langle |\psi_1(s,t,v,\cdot) - \psi_1(s,t,v',\cdot)|, p_{0,s}(v,\cdot,y)^2 \rangle_\alpha ds \\ + \int_{0}^{t} \| f^{(3)} \|_\infty \langle |\psi_1(s,t,v',\cdot)|, |\Phi_s(v,\cdot) - \Phi_s(v',\cdot)| \, p_{0,s}(v,\cdot,y)^2 \rangle_\alpha ds \\ + \int_{0}^{t} \| f'' \|_\infty \langle |\psi_1(s,t,v',\cdot)|, |p_{0,s}(v,\cdot,y)^2 - p_{0,s}(v',\cdot,y)^2| \rangle_\alpha ds.
 		\end{multline}
 		By \eqref{pst2_short_time}, for $ | q - \frac{3\alpha}{2} | < 1 $ and $ | q - \alpha | < 1 $, using \eqref{continuity_psi_k_Lpq} in the second line,
 		\begin{align*}
 			\int_\R \langle |\psi_1(s,t,v,\cdot) - \psi_1(s,t,v',\cdot)|, p_{0,s}(v,\cdot,y)^2 \rangle_\alpha e^{qy} dy &\leq \frac{C}{\sqrt{s}} \| \psi_1(s,t,v,\cdot) - \psi_1(s,t,v',\cdot) \|_{1,q-\alpha} \\
 			&\leq \frac{C'}{\sqrt{s}} \| v - v' \|_{2,\alpha},
 		\end{align*}
 		for some $ C' > 0 $.
 		Likewise, by \eqref{pst2_short_time} and the Cauchy-Schwarz inequality,
 		\begin{multline*}
 			\int_\R \langle |\psi_1(s,t,v',\cdot)|, |\Phi_s(v,\cdot) - \Phi_s(v',\cdot)| \, p_{0,s}(v,\cdot,y)^2 \rangle_\alpha e^{qy} dy \\ \begin{aligned}
 			&\leq \frac{C}{\sqrt{s}} \| \Phi_s(v,\cdot) - \Phi_s(v',\cdot) \|_{2,\alpha} \| \psi_1(s,t,v',\cdot) \|_{2,2q - 2\alpha} \\
 			&\leq \frac{C'}{\sqrt{s}} \| v - v' \|_{2,\alpha},
 			\end{aligned}
 		\end{multline*}
 		for some $ C' > 0 $, using \eqref{bounds_Phit_Hgamma} (with $ \gamma = 0 $) and \eqref{psi_k_Lpq} (for $ |q-\alpha|<1 $) in the second line.
 		Finally, we write
 		\begin{multline*}
 			\langle |\psi_1(s,t,v',\cdot)|, |p_{0,s}(v,\cdot,y)^2 - p_{0,s}(v',\cdot,y)^2| \rangle_\alpha \\ \leq \langle |\psi_1(s,t,v',\cdot)|, |p_{0,s}(v,\cdot,y) - p_{0,s}(v',\cdot,y)| p_{0,s}(v,\cdot,y) \rangle_\alpha \\ + \langle |\psi_1(s,t,v',\cdot)|, |p_{0,s}(v,\cdot,y) - p_{0,s}(v',\cdot,y)| p_{0,s}(v',\cdot,y) \rangle_\alpha.
 		\end{multline*}
 		Hence, by Lemma~\ref{lemma:diff_pst2_short_times}, for $ |q-\alpha|<1 $ and $ |2q-3\alpha| < 2 $, using \eqref{psi_k_Lpq} in the second line,
 		\begin{multline*}
 			\int_\R \langle |\psi_1(s,t,v',\cdot)|, |p_{0,s}(v,\cdot,y)^2 - p_{0,s}(v',\cdot,y)^2| \rangle_\alpha e^{qy} dy \\ \begin{aligned}
 			&\leq C s^{1/4} \| v - v' \|_{2,\alpha} \| \psi_1(s,t,v',\cdot) \|_{2,2q-3\alpha} \\
 			&\leq C' s^{1/4} \| v - v' \|_{2,\alpha}.
 			\end{aligned}
 		\end{multline*}
 		Coming back to \eqref{diff_Ut}, we have obtained \eqref{continuity_Ut}, which concludes the proof of Lemma~\ref{lemma:bounds_Detat_fixed_time}.
		The bound \eqref{D2eta_t_L1q_continuity} is obtained in a similar way, using \eqref{phi_p0t_delta}.
 	\end{proof}
	
	We now wish to let $ t \to \infty $ in the expressions for $ D \eta_t $ and $ D^2 \eta_t $.
	
	\subsection{Fréchet derivatives of the Katzenberger coordinate} \label{subsec:Dzeta}
	
	The aim of this subsection is to prove Proposition~\ref{prop:Dzeta}, Lemma~\ref{lemma:Dzeta_m} and Lemma~\ref{lemma:katzenberger}.
	They will follow from the following results, whose proofs are deffered to subsections~\ref{subsec:psi_k} and \ref{subsec:Ut_Upsilon}.
	
	\begin{lemma} \label{lemma:Psi}
		\begin{subequations}
			For any $ s \geq 0 $ and $ v \in \mathcal{V}_{\Beta{zeta},K,\Epsilon{zeta}} $, there exists $ \Psi(s,v,\cdot) \in L^{2,\alpha} $ such that the following holds.
			For any $ p \geq 1 $, $ q \in \R $ with $ |q| < p $, there exist $ C > 0 $, $\genericc > 0 $ and $ \delta > 0 $ such that, for all $ t > s $ and all $ v \in \mathcal{V}_{\Beta{zeta},K,\Epsilon{zeta}} $,
			\begin{align}
				\| \psi_1(s,t,v,\cdot) - \Psi(s,v,\cdot) \|_{p,q} &\leq C \dist(v,M)^\delta e^{-\genericc t}, \label{psi_1-Psi} \\
				\| \psi_2(0,t,v,\cdot) + \frac{\alpha}{2} \Psi(0,v,\cdot) \|_{p,q} &\leq C e^{-\genericc t}, \label{psi_2-Psi}
			\end{align}
			and
			\begin{align} \label{Psi-gradm}
				\| \Psi(s,v,\cdot) - \partial_{x} m_{\zeta(v)} \|_{p,q} &\leq C \, \dist(v,M)^\delta e^{-\genericc s}, \\
				\| \Psi(s,v,\cdot) \|_{p,q} &\leq C. \label{Lpq_Psi}
			\end{align}
			Moreover, there exists $ C > 0 $ such that, for all $ v \in \mathcal{V}_{\Beta{zeta},K,\Epsilon{zeta}} $ and all $ 0 \leq s \leq t $,
			\begin{equation} \label{bound_psi1_m}
				\enorm{\psi_1(s,t,v,\cdot)} \leq C.
			\end{equation}
			and
			\begin{equation} \label{bound_Psi_m}
				\enorm{\Psi(s,v,\cdot)}  \leq C.
			\end{equation}
		\end{subequations}
	\end{lemma}

	\begin{lemma} \label{lemma:Psi_Hgamma}
		\begin{subequations}
			For any $ \gamma \in [0,2) $, there exist $ C > 0 $ and $ \genericc > 0 $ such that, for all $ v \in \mathcal{V}_{\Beta{zeta},K,\Epsilon{zeta}} $ and $ t \geq 1 $,
			\begin{equation} \label{psi_1_Psi_H_gamma}
				\Hnorm{\psi_1(0,t,v,\cdot) - \Psi(0,v,\cdot)}{\gamma} \leq C \dist(v,M) e^{-\genericc t},
			\end{equation}
			and, for all $ t \geq 0 $,
			\begin{equation} \label{Psi_H_gamma}
				\Hnorm{\Psi(t,v,\cdot)}{\gamma} \leq C.
			\end{equation}
			In addition, for any $ \gamma' \in (0,2-\gamma) $, there exists $ C > 0 $ such that, for all $ v \in \mathcal{V}_{\Beta{zeta},K,\Epsilon{zeta}} $ and $ s \leq 1 $,
			\begin{equation} \label{lim_Psi_H_gamma}
				\Hnorm{\Psi(s,v,\cdot) - \Psi(0,v,\cdot)}{\gamma} \leq C s^{\gamma'/2}.
			\end{equation}
		\end{subequations}
	\end{lemma}

	Let us also define, for $ t \geq 0 $, $ v \in \mathcal{V}_{\Beta{zeta}, K, \Epsilon{zeta}} $ and $ y \in \R $,
	\begin{equation} \label{def:Upsilon}
		\Upsilon(v,y) := \int_{0}^{\infty} \langle \Psi(s,v,\cdot), f''(\Phi_s(v,\cdot)) p_{0,s}(v,\cdot,y)^2 \rangle_{\alpha} ds.
	\end{equation}
	This allows us to state the following.
	
	\begin{lemma} \label{lemma:U_Upsilon}
		\begin{subequations}
			For any $ K > 0 $ and any $ p \geq 1 $ and $ q \in \R $ satisfying \eqref{condition_p_q}, there exist $ C > 0 $, $ \delta > 0 $ and $ \genericc > 0 $ such that, for all $ v \in \mathcal{V}_{\Beta{zeta},K,\Epsilon{zeta}} $ and $ t \geq 0 $,
			\begin{align}
				\| U_t(v,\cdot) - \Upsilon(v,\cdot) \|_{p,q} &\leq C e^{-\genericc t}, \label{Ut-Upsilon} \\
				\| \Upsilon(v,\cdot) - \Upsilon(m_{\zeta(v)},\cdot) \|_{p,q} &\leq C \, \dist(v,M)^\delta, \label{Upsilon_v-m}
			\end{align}
			and
			\begin{equation} \label{bound_Upsilon}
				\| \Upsilon(v,\cdot) \|_{p,q} \leq C.
			\end{equation}
		\end{subequations}
	\end{lemma}
	
	Let us now use the above lemmas to prove Proposition~\ref{prop:Dzeta}, Lemma~\ref{lemma:Dzeta_m} and Lemma~\ref{lemma:katzenberger}.
	Taking $ t \to \infty $ in \eqref{Deta_t} and \eqref{D2eta_t}, we can define, for $ v \in \mathcal{V}_{\Beta{zeta},K,\Epsilon{zeta}} $,
	\begin{equation} \label{def:Dzeta}
		D\zeta(v,y) := - \frac{\Psi(0,v,y)}{\| \partial_{x} m_{\zeta(v)} \|_{2,\alpha}^2},
	\end{equation}
	and
	\begin{equation} \label{def:D2zeta}
		D^2\zeta(v,y) := - \frac{\alpha}{2 \| \partial_x m_{\zeta(v)} \|_{2,\alpha}^4} \Psi(0,v,y)^2 - \frac{1}{\| \partial_{x} m_{\zeta(v)} \|_{2,\alpha}^2} \Upsilon(v,y).
	\end{equation}
	By \eqref{Psi-gradm} and \eqref{def:Upsilon}, we already obtain Lemma~\ref{lemma:Dzeta_m}.
	
	\begin{proof}[Proof of Proposition~\ref{prop:Dzeta}]
		Take $ \gamma \in [0,2) $, $ p \geq 1 $ and $ q \in \R $ with $ |q| < p $.
		Then
		\begin{multline} \label{Detat_Dzeta}
			D \eta_t(v,y) - D\zeta(v,y) = \frac{\psi_1(0,t,v,y) - \Psi(0,v,y)}{\Vphi{1}(t,v)} \\ - \frac{\Psi(0,v,y)}{ \Vphi{1}(t,v) \| \partial_{x} m_{\zeta(v)} \|_{2,\alpha}^2} \left( \Vphi{1}(t,v) + \| \partial_{x} m_{\zeta(v)} \|_{2,\alpha}^2 \right).
		\end{multline}
		Combining \eqref{bound_varphi_tilde_below}, \eqref{varphi_Lipschitz}, Theorem~\ref{thm:zeta}, \eqref{psi_1_Psi_H_gamma} and \eqref{Psi_H_gamma}, we obtain \eqref{Detat-Dzeta_Hgamma}.
		We also obtain \eqref{Dzeta_Hgamma} from \eqref{Psi_H_gamma} and the fact that $ |\zeta(v)| \leq \Cst{sup_zeta} $.
		From \eqref{psi_1-Psi} and \eqref{Lpq_Psi}, we also deduce \eqref{Detat-Dzeta_Lpq}.
		The bounds \eqref{Dzeta_Lpq} and \eqref{Dzeta_v-m_Lpq} follow from the fact that $ |\zeta(v)| \leq \Cst{sup_zeta} $ and from \eqref{Lpq_Psi} and \eqref{Psi-gradm}.
		We now turn to the estimates on $ D^2 \eta_t $ and $ D^2 \zeta $.
		We start by noting that, using the Cauchy-Schwarz inequality,
		\begin{multline*}
			\left( \int_\R \left| \psi_1(0,t,v,y) \psi_2(0,t,v,y) + \frac{\alpha}{2} \Psi(0,v,y)^2 \right|^p e^{qy} dy \right)^{1/p} \\ \leq \| \psi_1(0,t,v,\cdot) - \Psi(0,v,\cdot) \|_{2p,q} \| \psi_2(0,t,v,\cdot) \|_{2p,q} \\ + \| \Psi(0,v,\cdot) \|_{2p,q} \left\| \psi_2(0,t,v,\cdot) + \frac{\alpha}{2} \Psi(0,v,\cdot) \right\|_{2p,q}.
		\end{multline*}	
		We can then bound the right-hand side using Lemma~\ref{lemma:Psi}, since \eqref{condition_p_q_D2zeta} implies $ |q| < 2p $.
		Proceeding in this way and using \eqref{varphi_Lipschitz}, we can deal with the first terms appearing in \eqref{D2eta_t}.
		The last term is dealt with using \eqref{Ut-Upsilon} and we obtain \eqref{D2etat-D2zeta}.
		The last bounds \eqref{D2zeta_Lpq} and \eqref{D2zeta_v-m} are obtained in a similar way, using \eqref{varphi_Lipschitz}, Lemma~\ref{lemma:Psi} and Lemma~\ref{lemma:U_Upsilon}.	
	\end{proof}
	
	Before proving Lemma~\ref{lemma:katzenberger}, we need the following, which is proved in Subsection~\ref{subsec:psi_k}
	
	\begin{lemma} \label{lemma:Psi_time_shift}
		For any $ v \in \mathcal{V}_{\Beta{zeta},K,\Epsilon{zeta}} $ and $ t \geq 0 $,
		\begin{equation*}
			\Psi(0,\Phi_t(v,\cdot),\cdot) = \Psi(t,v,\cdot).
		\end{equation*}
	\end{lemma}
	
	Let us now prove Lemma~\ref{lemma:katzenberger}.
	
	\begin{proof}[Proof of Lemma~\ref{lemma:katzenberger}]
		We then note that, since $ \zeta(\Phi_t(v,\cdot)) = \zeta(v) $ for all $ t > 0 $,
		\begin{align*}
			\frac{d}{dt} \zeta(\Phi_t(v,\cdot)) = \langle D\zeta(\Phi_t(v,\cdot),\cdot), \partial_{xx} \Phi_t(v,\cdot) + \alpha \partial_{x} \Phi_t(v,\cdot) + f(\Phi_t(v,\cdot)) \rangle_{\alpha} = 0.
		\end{align*}
		The scalar product above is well defined for any $ t > 0 $ and $ v \in L^{2,\alpha} $ thanks to \eqref{bounds_Phit_Hgamma} and \eqref{Dzeta_Hgamma}.
		We now wish to take the limit as $ t \to 0 $, using the fact that $ v \in H^{\gamma,\alpha} $ for some $ \gamma > 0 $.
		To do this, note that, using the fact that
		\begin{align*}
			\langle u, (\partial_{x x} + \alpha \partial_{x}) v \rangle_\alpha &= \langle \mathcal{F} \mathcal{J}u, \mathcal{F} \mathcal{J} (\partial_{x x} + \alpha \partial_{x}) v \rangle_{L^2} \\
			&= \left\langle \mathcal{F}\mathcal{J} u, \mathcal{F} \left(\partial_{x x} - \frac{\alpha^2}{4}\right) \mathcal{J} v \right\rangle_{L^2},
		\end{align*}
		and \eqref{deriv_fourier}, we can write
		\begin{multline} \label{katz_fourier}
			\langle D\zeta(\Phi_t(v,\cdot),\cdot), \partial_{xx} \Phi_t(v,\cdot) + \alpha \partial_{x} \Phi_t(v,\cdot) + f(\Phi_t(v,\cdot)) \rangle_{\alpha} \\= \langle \mathcal{F J} D\zeta(\Phi_t(v,\cdot),\cdot), (-\xi^2) \mathcal{F J} \Phi_t(v,\cdot) \rangle_{L^2} + \langle \mathcal{F J} D\zeta(\Phi_t(v,\cdot),\cdot), \mathcal{F J} \tilde{f}(\Phi_t(v,\cdot)) \rangle_{L^{2}},
		\end{multline}
		with $ \tilde{f}(u) = f(u) - \frac{\alpha^2}{4} u $.
		Note that, by \eqref{Dzeta_Hgamma} and the Cauchy-Schwarz inequality, this can be extended to $ t = 0 $ to give a definition of the left-hand side for $ t = 0 $ when $ v \in \mathcal{V}_{\Beta{zeta},K,\Epsilon{zeta}} \cap H^{\gamma,\alpha} $ for some $ \gamma > 0 $ since (assuming without loss of generality that $ \gamma < 2 $)
		\begin{align*}
			| \langle \mathcal{F J} D\zeta(v,\cdot), (-\xi^2) \mathcal{F J} v \rangle_{L^2} | \leq \Hnorm{D\zeta(v,\cdot)}{2-\gamma} \Hnorm{v}{\gamma}.
		\end{align*}
		Using the Cauchy-Schwarz inequality again, we write, for some $ \gamma' \in (0,2) $,
		\begin{multline} \label{diff_katzenberger}
			\left| \langle \mathcal{F J} D\zeta(v,\cdot), (-\xi^2) \mathcal{F J} v \rangle_{L^2} - \langle \mathcal{F J} D\zeta(\Phi_t(v,\cdot),\cdot), (-\xi^2) \mathcal{F J} \Phi_t(v,\cdot) \rangle_{L^2} \right| \\ \leq \Hnorm{D\zeta(\Phi_t(v,\cdot),\cdot)}{2-\gamma'} \Hnorm{\Phi_t(v,\cdot) - v}{\gamma'} \\ + \Hnorm{D\zeta(\Phi_t(v,\cdot),\cdot) - D\zeta(v,\cdot)}{2-\gamma'} \Hnorm{v}{\gamma'}.
		\end{multline}
		Then, by Lemma~\ref{lemma:Psi_time_shift} and the fact that $ \zeta(\Phi_t(v,\cdot)) = \zeta(v) $,
		\begin{align*}
			D\zeta(\Phi_t(v,\cdot),y) = - \frac{\Psi(t,v,y)}{\| \partial_{x} m_{\zeta(v)} \|_{2,\alpha}^2}.
		\end{align*}
		Thus, choosing $ \gamma' \in (0,\gamma) $ and $ \gamma'' \in (0,\gamma') $ and applying \eqref{Psi_H_gamma}, \eqref{lim_Psi_H_gamma} and \eqref{continuity_Phit_Hgamma}, we obtain
		\begin{multline*}
			\left| \langle \mathcal{F J} D\zeta(v,\cdot), (-\xi^2) \mathcal{F J} v \rangle_{L^2} - \langle \mathcal{F J} D\zeta(\Phi_t(v,\cdot),\cdot), (-\xi^2) \mathcal{F J} \Phi_t(v,\cdot) \rangle_{L^2} \right| \\ \leq C t^{\frac{\gamma-\gamma'}{2}} \Hnorm{v}{\gamma} + C t^{\gamma''/2} \Hnorm{v}{\gamma}.
		\end{multline*}
		The second term on the right of \eqref{katz_fourier} is treated in a similar way (it is in fact easier to deal with since there is no factor $ \xi^2 $).
		We thus obtain that
		\begin{multline*}
			\langle D\zeta(v,\cdot), \partial_{x x} v + \alpha \partial_{x} v + f(v) \rangle_\alpha \\ = \lim_{t \downarrow 0} \langle D\zeta(\Phi_t(v,\cdot),\cdot), \partial_{x x} \Phi_t(v,\cdot) + \alpha \partial_{x} \Phi_t(v,\cdot) + f(\Phi_t(v,\cdot)) \rangle_\alpha  = 0,
		\end{multline*}
		as claimed.
	\end{proof}
	
	It remains to prove Lemma~\ref{lemma:Psi}, Lemma~\ref{lemma:Psi_Hgamma},  Lemma~\ref{lemma:Psi_time_shift} and Lemma~\ref{lemma:U_Upsilon}.
	To do so, we first state in the next subsection several results on the long time behaviour of $ T_{s,t,v} $.
	
	\subsection{Long-time behaviour of the linearised semigroup} \label{subsec:semigroup_long_time}
	
	The following results are proved in Section~\ref{sec:semigroup}.
	Recall the norm $ \enorm{\cdot} $ introduced in \eqref{def:norm_exp}.
	
	\begin{lemma} \label{lemma:T-m}
		For any $ K > 0 $, there exists $ C > 0 $ such that, for all $ v \in \mathcal{V}_{\Beta{zeta},K,\Epsilon{zeta}} $, $ \phi \in L^{\infty,e} $ and all $ 0 \leq s \leq t $
		\begin{equation*}
			\enorm{T_{s,t,v} \phi} \leq C \enorm{\phi}.
		\end{equation*}
	\end{lemma}

	\begin{lemma} \label{lemma:Tst_L2alpha}
		\begin{subequations}
			For all $ t \geq 0 $ and $ \eta \in \R $,
		\begin{equation}
			\| e^{-t \mathcal{A}_\eta} \|_{L^{2,\alpha} \to L^{2,\alpha}} \leq 1\label{L2_A}
		\end{equation}
		In addition, $ e^{t \mathcal{A}_\eta} P_{\lbrace \varphi_\eta \rbrace} = P_{\lbrace \varphi_\eta \rbrace} $ and there exists $ \genericc > 0 $ such that, for all $ t \geq 0 $ and $ \eta \in \R $,
		\begin{equation}
			\| e^{-t \mathcal{A}_\eta} P_{\lbrace \varphi_\eta \rbrace^{\perp}} \|_{L^{2,\alpha} \to L^{2,\alpha}} \leq e^{-\genericc t}. \label{L2_A-P}
		\end{equation}
		Moreover, for any $ K > 0 $, there exists $ C > 0 $ such that, for all $ v \in \mathcal{V}_{\Beta{zeta},K,\Epsilon{zeta}} $ and all $ 0 \leq s \leq t $,
		\begin{align}
			& \| T_{s,t,v} \|_{L^{2,\alpha} \to L^{2,\alpha}} \leq C, \label{L2_T} \\
			& \| (T_{s,t,v} - e^{-(t-s)\mathcal{A}_{\zeta(v)}}) P_{\lbrace \varphi_{\zeta(v)} \rbrace^\perp} \|_{L^{2,\alpha} \to L^{2,\alpha}} \leq C (\| s(v) \|_\infty + \sqrt{\dist(v,M)}) e^{-\genericc t}. \label{L2_T-A}
		\end{align}
		\end{subequations}
	\end{lemma}

	Combining the two previous results, we obtain the following.
	
	\begin{cor} \label{cor:Tst_Lpq}
		\begin{subequations}
			For all $ K > 0 $, $ p \geq 1 $ and $ q \in \R $ such that $ | q | < p $, there exist $ C > 0 $, $ \genericc > 0 $ and $ \delta \in (0,1) $ such that, for all $ v \in \mathcal{V}_{\Beta{zeta},K,\Epsilon{zeta}} $, all $ \phi \in L^{2,\alpha} \cap L^{\infty,e} $ and all $ 0 \leq s < t $,
			\begin{equation} \label{Lpq_T}
				\| T_{s,t,v} \phi \|_{p,q} \leq C \| \phi \|_{2,\alpha}^{\delta} \enorm{\phi}^{1-\delta}.
			\end{equation}
			Moreover, for all $ | \zeta | \leq \Cst{sup_zeta}(K) $ and all $ t > 0 $,
			\begin{equation} \label{Lpq_A}
				\| e^{-t \mathcal{A}_{\zeta}} P_{\lbrace \varphi_{\zeta} \rbrace^\perp} \phi \|_{p,q} \leq C e^{-\genericc t} \| \phi \|_{2,\alpha}^{\delta} \enorm{\phi}^{1-\delta},
			\end{equation}
			and, for all $ v \in \mathcal{V}_{\Beta{zeta},K,\Epsilon{zeta}} $ and $ 0 \leq s \leq t $,
			\begin{multline} \label{Lpq_T-A}
				\| (T_{s,t,v} - e^{-(t-s) \mathcal{A}_{\zeta(v)}}) P_{\lbrace \varphi_{\zeta(v)} \rbrace^\perp} \phi \|_{p,q} \\ \leq C  e^{-\genericc t} (\| s(v) \|_\infty + \sqrt{\dist(v,M)})^\delta \| \phi \|_{2,\alpha}^{\delta} \enorm{\phi}^{1-\delta}.
			\end{multline}
		\end{subequations}
	\end{cor}
	
	\begin{proof}
		First note that 
		By H\"older's inequality, for any $ \delta \in (0,2/p) $,
		\begin{equation*}
			\| \phi \|_{p,q} \leq \| \phi \|_{2,\alpha}^{\delta} \| \phi \|_{p(\delta), q(\delta)}^{1-\delta},
		\end{equation*}
		with
		\begin{align*}
			p(\delta) := \frac{p(1-\delta)}{1 - \frac{\delta p}{2}}, && q(\delta) := \frac{q - \frac{\alpha}{2} \delta p}{1 - \frac{\delta p}{2}}.
		\end{align*}
		We then choose $ \delta $ small enough that $ |q(\delta)| < p(\delta) $ (using the fact that $ |q| < p $).
		For such a $ \delta $, there exists $ C > 0 $ such that, for all $ \phi \in L^{\infty,e} $,
		\begin{equation*}
			\| \phi \|_{p(\delta), q(\delta)} \leq C \enorm{\phi}.
		\end{equation*}
		We then obtain \eqref{Lpq_T}, \eqref{Lpq_A} and \eqref{Lpq_T-A} from \eqref{L2_T}, \eqref{L2_A-P}, \eqref{L2_T-A}, Lemma~\ref{lemma:projection_varphi} and Lemma~\ref{lemma:T-m}.
	\end{proof}
	
	Furthermore, we can also extend the estimates from Lemma~\ref{lemma:Tstv_Hgamma} as follows.
	
	\begin{lemma} \label{lemma:Tstv_Hgamma_large_times}
		\begin{subequations}
			For any $ K > 0 $ and $ \gamma \in [0,2) $, there exists a constant $ C >0 $ such that, for all $ v \in \mathcal{V}_{\Beta{zeta}, K, \Epsilon{zeta}} $, $ \phi \in L^{2,\alpha} $ and all $ 0 \leq s \leq t $
			\begin{equation} \label{T_H_gamma}
				\Hnorm{T_{s,t,v} \phi }{\gamma} \leq \frac{C(1 + |t-s|)}{(|t-s| \wedge 1)^{\gamma / 2}} \| \phi \|_{2,\alpha},
			\end{equation}
			and, for all $ \gamma' \in (0,2-\gamma) $, there exists $ C > 0 $ such that, for all $ 0 \leq s \leq t $,
			\begin{equation} \label{diff_Tst-T0t_Hgamma}
				\Hnorm{(T_{s,t,v}-T_{0,t,v}) \phi}{\gamma} \leq \frac{C s^{\gamma'/2} (1+t-s)}{((t-s)\wedge 1)^{\frac{\gamma+\gamma'}{2}}} \| \phi \|_{2,\alpha}.
			\end{equation}
		\end{subequations}
	\end{lemma}
	
	Lemma~\ref{lemma:Tstv_Hgamma_large_times} is proved in Subsection~\ref{subsec:Tst_large_times}, along with the two following results.
	
	\begin{lemma} \label{lemma:pst_phi_delta}
		For any $ \delta_0 > 0 $, any $ K > 0 $ and any $ \delta \in (0, \delta_0 \wedge (1-\frac{\alpha}{2})) $, there exist $ C > 0 $ and $ \delta' \in (0,1) $ such that, for any $ \phi \in L^{2,\alpha} \cap L^\infty $, any $ v \in \mathcal{V}_{\Beta{zeta}, K, \Epsilon{zeta}} $ and all $ t > 0 $,
		\begin{equation*}
			\abs{ \int_\R \phi(x) p_{0,t}(v,x,y) e^{-\delta_0 |x| + \frac{\alpha}{2} x} dx } \leq \frac{C}{(t \wedge 1)^{\delta'/4}} \| \phi \|_{2,\alpha}^{\delta'} \| \phi \|_\infty^{1-\delta'} e^{-\delta |y| - \frac{\alpha}{2} y}.
		\end{equation*}
	\end{lemma}

	For $ s, t \geq 0 $, $ \zeta \in \R $, $ y, z \in \R $, let us set
	\begin{equation} \label{def:Xi0}
		\Xi_0(s,t,\zeta,y,z) := \langle \partial_{x} m_{\zeta} f''(m_{\zeta}), p_{0,s}(m_\zeta, \cdot, y) p_{0,t}(m_\zeta, \cdot, z) \rangle_\alpha.
	\end{equation}
	
	\begin{lemma} \label{lemma:bounds_pst2}
		For any $ K > 0 $ and $ \delta \in \left( 0, \frac{1}{2} \wedge \left( 1 - \frac{\alpha}{2} \right) \right) $, there exist $ C > 0 $ and $ \genericc > 0 $ such that, for $ | \zeta | \leq \Cst{sup_zeta}(K) $ and all $ s > 0 $, $ t > 0 $,
		\begin{equation} \label{bound_Xi0}
			| \Xi_0(s,t,\zeta,y,z) | \leq \frac{C}{((s\wedge1)(t\wedge1))^{1/4}} \left( e^{-\genericc t} + e^{-\genericc s} \right) e^{-\delta(|z| + |y|)-\frac{\alpha}{2} (y + z)}.
		\end{equation}
		In addition, for all $ p > 0 $ and $ q \in \R $ such that
		\begin{equation} \label{condition_p_q}
			\abs{\frac{q}{p} - \alpha} < 1 \wedge (2-\alpha),
		\end{equation}
		there exists $ C > 0 $, $ \delta \in (0,1) $ and $ \gamma \in (0,3/4) $ such that, for all $ v \in \mathcal{V}_{\Beta{zeta},K,\Epsilon{zeta}} $, all $ \phi \in L^{2,\alpha} \cap L^{\infty,e} $ and all $ t > 0 $,
		\begin{equation} \label{Lpq_p2}
			\left( \int_\R | \langle \phi, p_{0,t}(v,\cdot,y)^2 \rangle_\alpha |^p e^{qy} dy \right)^{1/p} \leq \frac{C}{(t\wedge 1)^\gamma} \| \phi \|_{2,\alpha}^{\delta} \enorm{\phi}^{1-\delta}.
		\end{equation}
	\end{lemma}

	In addition, we also have the following.
	
	\begin{lemma} \label{lemma:T-A}
		For any $ v \in \mathcal{V}_{\Beta{zeta},K,\Epsilon{zeta}} $,
		\begin{equation*}
			p_{s,t}(v,x,y) = p_{s,t}(m_{\zeta(v)},x,y) + \int_{s}^{t} \int_\R p_{s,r}(v,z,y) c_r(v,z) p_{r,t}(m_{\zeta(v)},x,z) e^{\alpha z} dz dr,
		\end{equation*}
		where
		\begin{equation} \label{eq:crdefn}
			c_r(v,z) = f'(\Phi_r(v,z))-f'(m_{\zeta(v)}(z)).
		\end{equation}
	\end{lemma}
	
	\begin{proof}
		Take $v\in \mathcal{V}_{\Beta{zeta},K,\Epsilon{zeta}}$ and $ \phi \in L^{2,\alpha} $.
		for $0\le s \le t$, let
		\begin{equation*}
			u_{s,t} := T_{s,t,v}^* \phi - e^{-(t-s) \mathcal{A}_{\zeta(v)}} \phi.
		\end{equation*}
		Then $ u_{s,s} = 0 $ and, for $ t \geq s $,
		\begin{equation*}
			\partial_t u_{s,t} = - \mathcal{A}_{\zeta(v)} u_{s,t} + (f'(\Phi_t(v,\cdot)) - f'(m_{\zeta(v)})) T_{s,t,v}^* \phi.
		\end{equation*}
		As a consequence,
		\begin{equation*}
			u_{s,t} = \int_{s}^{t} e^{-(t-r) \mathcal{A}_{\zeta(v)}} \left\lbrace c_r(v,\cdot) T_{s,r,v}^* \phi \right\rbrace dr.
		\end{equation*}
		Taking the adjoint on both sides, we see that
		\begin{equation*}
			T_{s,t,v} \phi = e^{-(t-s) \mathcal{A}_{\zeta(v)}} \phi + \int_{s}^{t} T_{s,r,v} \left\lbrace c_r(v,\cdot) e^{-(t-r) \mathcal{A}_{\zeta(v)}}\phi \right\rbrace dr,
		\end{equation*}
		from which the claim follows.
	\end{proof}

	\subsection{Estimates on $ \psi_k(s,t,v,\cdot) $} \label{subsec:psi_k}

	Here, we prove Lemma~\ref{lemma:Psi}, Lemma~\ref{lemma:Psi_Hgamma} and Lemma~\ref{lemma:Psi_time_shift}, which are all concerned with $ \psi_k(s,t,v,\cdot) $ and $ \Psi(s,v,\cdot) $.
	We first prove Lemma~\ref{lemma:Psi}.
	
	\begin{proof}[Proof of Lemma~\ref{lemma:Psi}]
		Recall from \eqref{def:psi_k} that
		\begin{equation*}
			\psi_k(s,t,v,\cdot) = T_{s,t,v} \partial_{x}^k m_{\eta_t(v)}.
		\end{equation*}
		By \eqref{Lpq_T}, for $ v \in \mathcal{V}_{\Beta{zeta},K,\Epsilon{zeta}} $ and $ |q| \leq p $,
		\begin{equation*}
			\| T_{s,t,v} \lbrace \partial_{x}^k m_{\eta_t(v)} - \partial_{x}^k m_{\zeta(v)} \rbrace \|_{p,q} \leq C \enorm{\partial_{x}^k m_{\eta_t(v)} - \partial_{x}^k m_{\zeta(v)}}^{1-\delta} \| \partial_{x}^k m_{\eta_t(v)} - \partial_{x}^k m_{\zeta(v)} \|_{2,\alpha}^{\delta}.
		\end{equation*}
		Using \eqref{diff_eta_zeta} and the fact that there exists $ C > 0 $ such that $ | \zeta(v) | \leq C $ for all $ v \in  \mathcal{V}_{\Beta{zeta},K,\Epsilon{zeta}} $, we obtain that there exists a constant $ C > 0 $ such that, for all  $ v \in  \mathcal{V}_{\Beta{zeta},K,\Epsilon{zeta}} $,
		\begin{equation*}
			\enorm{\partial_{x}^k m_{\eta_t(v)} - \partial_{x}^k m_{\zeta(v)}} \leq C.
		\end{equation*}
		In addition, proceeding as in the proof of Lemma~\ref{lemma:diff_m_eta}, we obtain that there exists a constant $ C > 0 $ such that
		\begin{equation*}
			\| \partial_{x}^k m_{\eta_t(v)} - \partial_{x}^k m_{\zeta(v)} \|_{2,\alpha} \leq C | \eta_t(v) - \zeta(v) |.
		\end{equation*}
		As a result, by Theorem~\ref{thm:zeta}, there exist $ C > 0 $, $ \delta \in [0,1) $ and $ \genericc > 0 $ such that
		\begin{equation} \label{bound_T_m_eta_zeta}
			\| T_{s,t,v} \lbrace \partial_{x}^k m_{\eta_t(v)} - \partial_{x}^k m_{\zeta(v)} \rbrace \|_{p,q} \leq C e^{-\genericc t} \dist(v,M)^{\delta}.
		\end{equation}
		We then note that, in view of the definition of $ T_{s,t,v} $,
		\begin{equation*}
			T_{s,t,v} \partial_{x} m_{\zeta(v)} = \partial_{x} m_{\zeta(v)} + \int_{s}^{t} T_{s,r,v} \left( (\partial_{xx} + \alpha \partial_{x} + f'(\Phi_r(v,\cdot))) \partial_{x} m_{\zeta(v)} \right) dr.
		\end{equation*}
		Since $ (\partial_{xx} + \alpha \partial_{x} + f'(m_\zeta)) \partial_{x} m_{\zeta} = 0 $, this is also
		\begin{equation} \label{T_grad_m}
			T_{s,t,v} \partial_{x} m_{\zeta(v)} = \partial_{x} m_{\zeta(v)} + \int_{s}^{t} T_{s,r,v} \left( (f'(\Phi_r(v,\cdot)) - f'(m_{\zeta(v)})) \partial_{x} m_{\zeta(v)} \right) dr.
		\end{equation}
		By \eqref{Lpq_T}, there exist $ C > 0 $ and $ \delta \in (0,1) $ such that
		\begin{multline*}
			\| T_{s,r,v} \left( (f'(\Phi_r(v,\cdot)) - f'(m_{\zeta(v)})) \partial_{x} m_{\zeta(v)} \right) \|_{p,q} \\ \leq C \enorm{\partial_x m_{\zeta(v)}}^{1-\delta} \| \Phi_r(v,\cdot) - m_{\zeta(v)} \|_{\infty}^{1-\delta} \| \Phi_r(v,\cdot) - m_{\zeta(v)} \|_{2,\alpha}^{\delta}.
		\end{multline*}
		By Corollary~\ref{cor:uniform_cvg} and Theorem~\ref{thm:zeta}, and because $\| \cdot \|_\infty \leq \|\cdot \|_{\infty,\lambda}$, this yields
		\begin{equation*}
			\| T_{s,r,v} \left( (f'(\Phi_r(v,\cdot)) - f'(m_{\zeta(v)})) \partial_{x} m_{\zeta(v)} \right) \|_{p,q} \leq C \dist(v,M)^{\delta} e^{-\genericc r}.
		\end{equation*}
		We can thus set
		\begin{equation} \label{def_Psi}
			\Psi(s,v,y) := \partial_{x} m_{\zeta(v)}(y) + \int_{s}^{\infty} T_{s,r,v} \left( (f'(\Phi_r(v,\cdot)) - f'(m_{\zeta(v)})) \partial_{x} m_{\zeta(v)} \right)(y) dr.
		\end{equation}
		Then
		\begin{align*}
			\| T_{s,t,v} \partial_{x} m_{\zeta(v)} - \Psi(s,v,\cdot) \|_{p,q} &\leq \int_{t}^{\infty} \left\| T_{s,r,v} \left( (f'(\Phi_r(v,\cdot)) - f'(m_{\zeta(v)})) \partial_{x} m_{\zeta(v)} \right) \right\|_{p,q} dr \\
			&\leq C \int_{t}^{\infty} e^{-\genericc r} dr \, \dist(v,M)^{\delta} \\
			&= \frac{C}{\genericc} e^{-\genericc t} \dist(v,M)^{\delta}. \numberthis \label{Tmzeta_Psi}
		\end{align*}
		Combined with \eqref{bound_T_m_eta_zeta} with $ k = 1 $, this proves \eqref{psi_1-Psi}.
		We also note that
		\begin{align*}
			\| \Psi(s,v,\cdot) - \partial_{x} m_{\zeta(v)} \|_{p,q} &\leq \int_{s}^{\infty} \left\| T_{s,r,v} \left( (f'(\Phi_r(v,\cdot)) - f'(m_{\zeta(v)})) \partial_{x} m_{\zeta(v)} \right) \right\|_{p,q} dr \\
			&\leq C \int_{s}^{\infty} e^{-\genericc r } dr \, \dist(v,M)^{\delta}, \\
			&= \frac{C}{\genericc} e^{-\genericc s} \dist(v,M)^{\delta},
		\end{align*}
		yielding \eqref{Psi-gradm}.
		We also obtain \eqref{Lpq_Psi} from the above and the fact that $ \| \partial_{x} m_{\zeta(v)} \|_{p,q} \leq C $ for $ |q| < p $ and $ v \in \mathcal{V}_{\Beta{zeta},K,\Epsilon{zeta}} $.
		We start by noting that, by Lemma~\ref{lemma:T-m},
		\begin{equation*}
			\enorm{\psi_1(s,t,v,\cdot)} \leq C \enorm{\partial_{x} m_{\eta_t(v)}}.
		\end{equation*}
		By Theorem~\ref{thm:zeta}, there exists a constant $ C > 0 $ such that, for all $ v \in \mathcal{V}_{\Beta{zeta},K,\Epsilon{zeta}} $, $ \enorm{\partial_{x} m_{\eta_t(v)}} \leq C $, yielding \eqref{bound_psi1_m}.
		In addition, using Lemma~\ref{lemma:T-m} in \eqref{def_Psi}, we have
		\begin{equation*}
			\enorm{\Psi(s,v,\cdot)} \leq \enorm{\partial_x m_{\zeta(v)}} \left( 1 + \int_{s}^{\infty} \| c_r(v,\cdot) \|_\infty dr \right).
		\end{equation*}
		Then, using Corollary~\ref{cor:sup_norm}, we obtain \eqref{bound_Psi_m}.
		We now turn to $ \psi_2(0,t,v,\cdot) $, and we write
		\begin{multline*}
			T_{0,t,v} \partial_{xx} m_{\zeta(v)} = \langle \varphi_{\zeta(v)}, \partial_{xx} m_{\zeta(v)} \rangle_\alpha T_{0,t,v} \varphi_{\zeta(v)} \\ + \left( T_{0,t,v} - e^{-t \mathcal{A}_{\zeta(v)}} \right) P_{\lbrace \varphi_{\zeta(v)} \rbrace^\perp} \partial_{xx} m_{\zeta(v)} + e^{-t \mathcal{A}_{\zeta(v)}} P_{\lbrace \varphi_{\zeta(v)} \rbrace^\perp} \partial_{xx} m_{\zeta(v)}.
		\end{multline*}
		We can then use Corollary~\ref{cor:Tst_Lpq} to bound the last two terms.
		By \eqref{Lpq_T-A}, there exist $ C > 0 $ and $ \genericc > 0 $ such that
		\begin{equation*}
			\| \left( T_{0,t,v} - e^{-t \mathcal{A}_{\zeta(v)}} \right) P_{\lbrace \varphi_{\zeta(v)} \rbrace^\perp} \partial_{xx} m_{\zeta(v)} \|_{p,q} \leq C e^{-\genericc t}.
		\end{equation*}
		In addition, by \eqref{Lpq_A}, there exist $ C > 0 $ and $ \genericc > 0 $ such that
		\begin{equation*}
			\|  e^{-t \mathcal{A}_{\zeta(v)}} P_{\lbrace \varphi_{\zeta(v)} \rbrace^\perp} \partial_{xx} m_{\zeta(v)} \|_{p,q} \leq C e^{-\genericc t}.
		\end{equation*}
		In addition, by \eqref{IPP2} and \eqref{IPP3},
		\begin{equation} \label{dxxm_dxm}
			\langle \partial_{x} m_\eta, \partial_{xx} m_\eta \rangle_\alpha = -\frac{\alpha}{2} \| \partial_{x} m_\eta \|_{2,\alpha}^2.
		\end{equation}
		Combining this with the definition of $ \varphi_{\zeta} $ in \eqref{eq:phietadefn}, we obtain
		\begin{align*}
			\langle \varphi_{\zeta(v)}, \partial_{xx} m_{\zeta(v)} \rangle_\alpha T_{0,t,v} \varphi_{\zeta(v)} &= \frac{\langle \partial_{xx} m_{\zeta(v)}, \partial_{x} m_{\zeta(v)} \rangle_\alpha}{\| \partial_{x} m_{\zeta(v)} \|_{2,\alpha}^2} T_{0,t,v} \partial_{x} m_{\zeta(v)} \\
			&= -\frac{\alpha}{2} T_{0,t,v} \partial_{x} m_{\zeta(v)}.
		\end{align*}
		As a result,
		\begin{equation*}
			\| T_{0,t,v} \partial_{xx} m_{\zeta(v)} + \frac{\alpha}{2} T_{0,t,v} \partial_{x} m_{\zeta(v)} \|_{p,q} \leq C e^{-\genericc t}.
		\end{equation*}
		The bound \eqref{psi_2-Psi} then follows from \eqref{bound_T_m_eta_zeta} and \eqref{Tmzeta_Psi}, which concludes the proof of the lemma.
	\end{proof}

	\begin{proof}[Proof of Lemma~\ref{lemma:Psi_Hgamma}]
		We start by noting that, using \eqref{T_H_gamma} and Theorem~\ref{thm:zeta}, for all $ t \geq 1 $,
		\begin{equation} \label{T_m_eta_zeta_H_gamma}
			\Hnorm{T_{0,t,v} \lbrace \partial_{x} m_{\eta_t(v)} - \partial_{x} m_{\zeta(v)} \rbrace}{\gamma} \leq C \dist(v,M) (1 + t) e^{-\genericc t}.
		\end{equation}
		Then, using \eqref{T_grad_m} and \eqref{def_Psi}, we obtain
		\begin{equation*}
			\Hnorm{T_{0,t,v} \partial_{x} m_{\zeta(v)} - \Psi(0,v,\cdot)}{\gamma} \leq \int_{t}^{\infty} \Hnorm{T_{0,r,v} \lbrace c_r(v,\cdot) \partial_{x} m_{\zeta(v)} \rbrace}{\gamma} dr.
		\end{equation*}
		Combining \eqref{T_H_gamma} and Theorem~\ref{thm:zeta} yields, for $ r \geq s \geq 0 $
		\begin{equation} \label{T_c_gradm_H_gamma}
			\Hnorm{T_{s,r,v} \lbrace c_r(v,\cdot) \partial_{x} m_{\zeta(v)} \rbrace}{\gamma} \leq \frac{C(1+|r-s|)}{(|r-s| \wedge 1)^{\gamma/2}} \dist(v,M) e^{-\genericc r},
		\end{equation}
		from which we deduce that there exists $ C > 0 $ and $ \genericc > 0 $ such that, for $ t \geq 1 $,
		\begin{equation} \label{T_m_Psi_H_gamma}
			\Hnorm{T_{0,t,v} \partial_{x} m_{\zeta(v)} - \Psi(0,v,\cdot)}{\gamma} \leq C \dist(v,M) e^{-\genericc t}.
		\end{equation}
		Combining \eqref{T_m_eta_zeta_H_gamma} and \eqref{T_m_Psi_H_gamma}, we obtain \eqref{psi_1_Psi_H_gamma}.
		Using \eqref{T_c_gradm_H_gamma} in \eqref{def_Psi}, we obtain
		\begin{equation*}
			\Hnorm{\Psi(s,v,\cdot)}{\gamma} \leq \Hnorm{\partial_{x} m_{\zeta(v)}}{\gamma} +  \int_{s}^{\infty} \frac{C(1+|r-s|)}{(|r-s|\wedge 1)^{\gamma/2}} \Beta{zeta} e^{-\genericc r} dr.
		\end{equation*}
		Since $ \gamma \in [0,2) $ and $ \partial_{x} m \in H^{2,\alpha} $, we obtain \eqref{Psi_H_gamma}.
		Finally, from \eqref{def_Psi}, we have
		\begin{multline} \label{Psi_s_Psi_0}
			\Psi(s,v,y) - \Psi(0,v,y) = \int_{s}^{\infty} (T_{s,r,v} - T_{0,r,v}) \lbrace c_r(v,\cdot) \partial_{x} m_{\zeta(v)} \rbrace(y) dr \\ - \int_{0}^{s} T_{0,r,v} \lbrace c_r(v,\cdot) \partial_{x} m_{\zeta(v)} \rbrace(y) dr.
		\end{multline}
		The second term on the right-hand side tends to zero as $ s \to 0 $ by \eqref{T_c_gradm_H_gamma} and the fact that $ \gamma \in [0,2) $.
		For the first term, by \eqref{diff_Tst-T0t_Hgamma} and Theorem~\ref{thm:zeta}, for $ \gamma' \in (0,2-\gamma) $,
		\begin{equation*}
			\Hnorm{(T_{s,r,v} - T_{0,r,v}) \lbrace c_r(v,\cdot) \partial_{x} m_{\zeta(v)} \rbrace}{\gamma} \leq \frac{C s^{\gamma'/2}(1+r-s)}{(|r-s|\wedge 1)^{\frac{\gamma+\gamma'}{2}}} \Beta{zeta} e^{-\genericc r}.
		\end{equation*}
		Plugging this in the first term in \eqref{Psi_s_Psi_0}, we obtain \eqref{lim_Psi_H_gamma}.
	\end{proof}
	
	To prove Lemma~\ref{lemma:Psi_time_shift}, we need the follwing, which is proved in Subsection~\ref{subsec:Tst_short_times}.
	
	\begin{lemma} \label{lemma:Tstv_time_shift}
		For any $ v \in L^{2,\alpha} $, $ 0 \leq s \leq t $ and $ \tau \geq 0 $,
		\begin{equation*}
			T_{s + \tau, t + \tau, v} = T_{s,t,\Phi_t(v,\cdot)}.
		\end{equation*}
	\end{lemma}
	
	This allows us to prove Lemma~\ref{lemma:Psi_time_shift}.
	
	\begin{proof}[Proof of Lemma~\ref{lemma:Psi_time_shift}]
		By the definition of $ \Psi(s,v,\cdot) $ in \eqref{def_Psi} and by the fact that $ \zeta(\Phi_t(v,\cdot)) = \zeta(v) $,
		\begin{equation*}
			\Psi(0,\Phi_t(v,\cdot),\cdot) = \partial_{x} m_{\zeta(v)} + \int_{0}^{t} T_{0,r,\Phi_t(v,\cdot)} \left( c_r(\Phi_t(v,\cdot),\cdot) \partial_{x} m_{\zeta(v)} \right) dr.
		\end{equation*}
		Combining Lemma~\ref{lemma:Tstv_time_shift} and the fact that $ \Phi_r(\Phi_r(v,\cdot),\cdot) = \Phi_{r+t}(v,\cdot) $, we obtain
		\begin{equation*}
			\Psi(0,\Phi_t(v,\cdot),\cdot) = \partial_{x} m_{\zeta(v)} + \int_{0}^{\infty} T_{t,r+t,v} \left( c_{r+t}(v,\cdot) \partial_{x} m_{\zeta(v)} \right) dr.
		\end{equation*}
		By a change of variables and using \eqref{def_Psi} again, we obtain the result.
	\end{proof}

	\subsection{Convergence of $ U_t(v,\cdot) $} \label{subsec:Ut_Upsilon}
	
	We now prove Lemma~\ref{lemma:U_Upsilon}.
	To do this, we define the following quantities, for $0\le s \le t$, $v \in \mathcal{V}_{\Beta{zeta},K,\Epsilon{zeta}}$, $ \zeta \in \R $ and $ y \in \R $:
	\begin{equation} \label{eq:Xidefs}
		\begin{aligned}
			&\Xi_1(s,t,v,y) = \langle \psi_1(s,t,v,\cdot) - \Psi(s, v,\cdot), f''(\Phi_s(v,\cdot)) p_{0,s}(v,\cdot,y)^2 \rangle_{{\alpha}}, \\
			&\Xi_2(s,v,y) = \langle \Psi(s,v,\cdot) f''(\Phi_s(v,\cdot)) - \partial_{x} m_{\zeta(v)} f''(m_{\zeta(v)}), p_{0,s}(v,\cdot, y)^2 \rangle_{{\alpha}}, \\
			&\Xi_3(s,v,y) = \langle \partial_{x} m_{\zeta(v)} f''(m_{\zeta(v)}), \left( p_{0,s}(v,\cdot,y)^2 - p_{0,s}(m_{\zeta(v)},\cdot,y)^2 \right) \rangle_{{\alpha}}, \\
			&\Xi_4(s,\zeta,y) = \langle \partial_{x} m_{\zeta} f''(m_{\zeta}), p_{0,s}(m_{\zeta},\cdot,y)^2 \rangle_{{\alpha}}.
		\end{aligned}
	\end{equation}
	With these definitions, and recalling the definitions of $U_t$ and $\Upsilon$ in~\eqref{eq:Udefn} and~\eqref{def:Upsilon}, we can write
	\begin{align} \label{decom_U-Upsilon}
		U_t(v,y) - \Upsilon(v,y) &= \int_{0}^{t} \Xi_1(s,t,v,y) ds - \int_{t}^{\infty} \Xi_2(s,v,y) ds \notag \\
		&\qquad\quad - \int_{t}^{\infty} \Xi_3(s,v,y) ds - \int_{t}^{\infty} \Xi_4(s,\zeta(v),y) ds, \\
		\Upsilon(m_{\zeta},y) &= \int_{0}^{\infty} \Xi_4(s,\zeta,y) ds, \label{eq:Upsilon_m} \\
		\intertext{and}
		\Upsilon(v,y) - \Upsilon(m_{\zeta(v)},y) &= \int_{0}^{\infty} \Xi_2(s,v,y) ds + \int_{0}^{\infty} \Xi_3(s,v,y) ds. \label{decomp_Upsilon_v-m}
	\end{align}
	
	The following lemma then allows us to prove Lemma~\ref{lemma:U_Upsilon}.
	
	\begin{lemma} \label{lemma:Xi}
		For any $ K > 0 $ and any $ p \geq 1 $, $ q \in \R $ satisfying \eqref{condition_p_q}, there exist $ C > 0 $, $ \delta \in (0,1) $, $ \genericc > 0 $ and $ \gamma \in (0,3/4) $ such that, for all $ v \in \mathcal{V}_{\Beta{zeta},K,\Epsilon{zeta}} $ and any $ 0 < s < t $,
		\begin{subequations}
			\begin{align}
				\| \Xi_1(s,t,v,\cdot) \|_{p,q} &\leq \frac{C}{(s \wedge 1)^\gamma} e^{-\genericc t} \dist(v,M)^\delta, \label{bound_Xi1} \\
				\| \Xi_2(s,v,\cdot) \|_{p,q} &\leq \frac{C}{(s \wedge 1)^\gamma} e^{-\genericc s} \dist(v,M)^\delta, \label{bound_Xi2} \\
				\| \Xi_3(s,v,\cdot) \|_{p,q} &\leq \frac{C}{(s \wedge 1)^{1/4}} e^{-\genericc s} \dist(v,M)^\delta, \label{bound_Xi3} \\
				\| \Xi_4(s,\zeta(v),\cdot) \|_{p,q} &\leq \frac{C}{\sqrt{s \wedge 1}} e^{-\genericc s}. \label{bound_Xi4}
			\end{align}
		\end{subequations}
	\end{lemma}
	
	Lemma~\ref{lemma:U_Upsilon} then follows directly from the above, using \eqref{decom_U-Upsilon}, \eqref{eq:Upsilon_m} and \eqref{decomp_Upsilon_v-m}.
	It thus remains to prove Lemma~\ref{lemma:Xi}.
	
	\begin{proof}[Proof of Lemma~\ref{lemma:Xi}]
		We bound each $ \Xi_i $ separately.
		
		\paragraph*{Bound on $ \Xi_1 $.}
		
		Using \eqref{Lpq_p2}, \eqref{bound_psi1_m} and \eqref{bound_Psi_m}, we obtain that there exists $ C > 0 $, $ \gamma \in (0,3/4) $ and $ \delta \in (0,1) $ such that
		\begin{equation*}
			\| \Xi_1(s,t,v,\cdot) \|_{p,q} \leq \frac{C}{(s \wedge 1)^\gamma} \| \psi_1(s,t,v,\cdot) - \Psi(s,v,\cdot) \|_{2,\alpha}^{\delta}.
		\end{equation*}
		The bound \eqref{bound_Xi1} then follows from \eqref{psi_1-Psi}.
		
		\paragraph*{Bound on $ \Xi_2 $.}
		
		Using \eqref{Lpq_p2} and \eqref{bound_Psi_m}, we obtain that there exists $ C > 0 $, $ \gamma \in (0,3/4) $ and $ \delta \in (0,1) $ such that
		\begin{equation*}
			\| \Xi_2(s,v,\cdot) \|_{p,q} \leq \frac{C}{(s \wedge 1)^{\gamma}} \| \Psi(s,v,\cdot) f''(\Phi_s(v,\cdot)) - \partial_{x} m_{\zeta(v)} f''(m_{\zeta(v)}) \|_{2,\alpha}^{\delta}.
		\end{equation*}
		We then write
		\begin{multline*}
			\| \Psi(s,v,\cdot) f''(\Phi_s(v,\cdot)) - \partial_{x} m_{\zeta(v)} f''(m_{\zeta(v)}) \|_{2,\alpha} \\ \leq C \left( \| \Psi(s,v,\cdot) - \partial_{x} m_{\zeta(v)} \|_{2,\alpha} + \| \Phi_s(v,\cdot) - m_{\zeta(v)} \|_{2,\alpha} \right),
		\end{multline*}
		and use Theorem~\ref{thm:zeta} and \eqref{Psi-gradm} to obtain \eqref{bound_Xi2}.
		
		\paragraph*{Bound on $ \Xi_4 $.}
		
		For the bound on $ \Xi_4 $, recalling the definition of $ \Xi_0 $ in \eqref{def:Xi0}, we note that
		\begin{equation*}
			\Xi_4(s,\zeta,y) = \Xi_0(s,s,\zeta,y,y).
		\end{equation*}
		Then, by \eqref{bound_Xi0} in Lemma~\ref{lemma:bounds_pst2}, for any $ \delta \in (0,\frac{1}{2} \wedge (1-\frac{\alpha}{2})) $, there exist $ C > 0 $ and $ \genericc > 0 $ such that, for $ | \zeta | \leq \Cst{sup_zeta}(K) $ and $ s > 0 $,
		\begin{equation*}
			\| \Xi_4(s,\zeta,y) \|_{p,q} \leq \frac{C}{(s\wedge 1)^{1/2}} e^{-\genericc s} \left( \int_\R e^{-2p\delta |y| + (q-p\alpha) y} dy \right)^{1/p}.
		\end{equation*}
		Since $ p $ and $ q $ satisfy \eqref{condition_p_q}, there exists $ \delta \in (0,\frac{1}{2} \wedge (1-\frac{\alpha}{2})) $ such that
		\begin{equation*}
			2 \delta p > \abs{q - p\alpha},
		\end{equation*}
		yielding \eqref{bound_Xi4}.
		
		\paragraph*{Bound on $ \Xi_3 $.}
		
		We start by noting that, if we set
    	\begin{equation*}
    		\Xi_{3,1}(s,v_1,v_2,y) = \langle \partial_{x} m_{\zeta(v_1)} f''(m_{\zeta(v_1)}), (p_{0,s}(v_1,\cdot,y)-p_{0,s}(m_{\zeta(v_1)},\cdot,y)) p_{0,s}(v_2,\cdot,y) \rangle_{{\alpha}},
    	\end{equation*}
    	we can write
    	\begin{equation}\label{decomp_Xi_3}
    		\Xi_3(s,v,y) = \Xi_{3,1}(s,v,v,y) + \Xi_{3,1}(s,v,m_{\zeta(v)},y).
    	\end{equation}
    	Then, by Lemma~\ref{lemma:T-A},
    	\begin{equation} \label{Xi31-Xi32}
    		\Xi_{3,1}(s,v_1,v_2,y) = \int_{0}^{s} \Xi_{3,2}(s,r,v_1,v_2,y) dr,
    	\end{equation}
    	with
    	\begin{multline} \label{Xi32}
    		\Xi_{3,2}(s,r,v_1,v_2,y) = \int_{\R^2} \partial_{x} m_{\zeta(v_1)}(x) f''(m_{\zeta(v_1)}(x)) p_{0,s}(v_2,x,y)\\ 
    		\times p_{0,r}(v_1,z,y) c_r(v_1,z) p_{r,s}(m_{\zeta(v_1)},x,z) e^{\alpha(x+z)} dx dz,
    	\end{multline}
    	Applying Lemma~\ref{lemma:T-A} again to $ p_{0,s}(v_2,x,y) $ in \eqref{Xi32}, we write
    	\begin{equation} \label{decomp_Xi_31}
    		\Xi_{3,2}(s,r,v_1, v_2, y) = \Xi_{3,3}(s,r, v_1, v_2, y) + \int_{0}^{s} \Xi_{3,4}(s,r,u,v_1, v_2, y) du,
    	\end{equation}
    	where
    	\begin{multline*}
    		\Xi_{3,3}(s,r,v_1, v_2, y) := \int_{\R^2} \partial_{x} m_{\zeta(v_1)}(x) f''(m_{\zeta(v_1)}(x)) p_{0,s}(m_{\zeta(v_2)},x,y) \\ \times p_{0,r}(v_1, z, y) c_r(v_1,z) p_{r,s}(m_{\zeta(v_1)},x,z) e^{\alpha (z+x)} dx dz,
    	\end{multline*}
    	and
    	\begin{multline*}
    		\Xi_{3,4}(s,r,u,v_1,v_2,y) \\ := \int_{\R^3} \partial_{x} m_{\zeta(v_1)}(x) f''(m_{\zeta(v_1)}(x)) p_{0,u}(v_2, z_2, y) c_u(v_2, z_2) p_{u,s}(m_{\zeta(v_2)}, x, z_2) \\ \times  p_{0,r}(v_1, z_1, y) c_r(v_1,z_1) p_{r,s}(m_{\zeta(v_1)},x,z_1) e^{\alpha (z_1 + z_2 + x)} dx dz_1 dz_2.
    	\end{multline*}
    	Recalling that $ \zeta(v_1) = \zeta(v_2) = \zeta(v) $ and that $ p_{r,s}(m_{\zeta(v_1)}, x, z) = p_{0,s-r}(m_{\zeta(v_1)}, x, z) $, by the definition of $ \Xi_0 $ in \eqref{def:Xi0}, we can write
    	\begin{equation} \label{Xi33-Xi0}
    		\Xi_{3,3}(s,r,v_1, v_2, y) = \int_\R \Xi_0(s, s-r, \zeta(v), y, z) p_{0,r}(v_1, z, y) c_r(v_1,z) e^{\alpha z } dz,
    	\end{equation}
    	and
    	\begin{multline*}
	    	\Xi_{3,4}(s,r,u,v_1,v_2,y) = \int_{\R^2} \Xi_0(s-u, s-r, \zeta(v), z_2, z_1) p_{0,u}(v_2, z_2, y) c_u(v_2, z_2) \\ \times p_{0,r}(v_1, z_1, y) c_r(v_1,z_1) e^{\alpha(z_1 + z_2)} dz_1 dz_2.
    	\end{multline*}
    	We can then write, using \eqref{Xi31-Xi32} and \eqref{decomp_Xi_31},
    	\begin{equation} \label{decomp_Xi_31_small_times}
    		\Xi_{3,1}(s,v_1, v_2, y) = \int_{0}^{s} \Xi_{3,3}(s,r,v_1,v_2,y) dr + \int_{0}^{s} \int_{0}^{s} \Xi_{3,4}(s,r,u,v_1,v_2,y) du dr.
    	\end{equation}
    	We then bound each term separately.
    	
    	By \eqref{bound_Xi0} in Lemma~\ref{lemma:bounds_pst2}, for any $ \delta \in (0,\frac{1}{2} \wedge (1-\frac{\alpha}{2})) $, there exist $ C > 0 $ and $ \genericc > 0 $ such that
    	\begin{multline} \label{1st_bound_Xi33}
    		| \Xi_{3,3}(s,r,v_1,v_2,y) | \leq \frac{C \left( e^{-\genericc s} + e^{-\genericc (s-r)} \right)}{((s \wedge 1) ((s-r) \wedge 1))^{1/4}}  e^{-\delta |y| - \frac{\alpha}{2} y} \\ \times  \int_\R |c_r(v_1,z)| p_{0,r}(v_1,z,y) e^{-\delta |z| + \frac{\alpha}{2} z} dz.
    	\end{multline}
    	Then, by Lemma~\ref{lemma:pst_phi_delta}, Theorem~\ref{thm:zeta}, Corollary~\ref{cor:uniform_cvg}, and $\|\cdot\|_\infty \leq \|\cdot\|_{\infty,\lambda}$, for any $ \delta_1 \in (0, \delta) $, there exists $ C > 0 $, $ \genericc > 0 $ and $ \delta' \in (0,1) $ such that
    	\begin{equation} \label{bound_cr_p0r}
    		\int_\R |c_r(v_1,z)| p_{0,r}(v_1,z,y) e^{-\delta |z| + \frac{\alpha}{2} z} dz \leq \frac{C e^{-\genericc r}}{(r \wedge 1)^{\delta'/4}} \dist(v_1,M)^{\delta'} e^{-\delta_1 |y| - \frac{\alpha}{2} y}.
    	\end{equation}
    	Plugging this in \eqref{1st_bound_Xi33}, we obtain
    	\begin{equation*}
    		| \Xi_{3,3}(s,r,v_1,v_2,y) | \leq \frac{C \left( e^{-\genericc s} + e^{-\genericc (s-r)} \right) e^{-\genericc r}}{((s \wedge 1) ((s-r) \wedge 1))^{1/4} (r \wedge 1)^{\delta'/4}}  \, \dist(v_1,M)^{\delta'} e^{-(\delta+\delta_1)|y| - \alpha y}.
    	\end{equation*}
    	We then note that
    	\begin{equation*}
    		\int_\R e^{-p(\delta + \delta_1)|y| + (q-\alpha) y} dy < \infty
    	\end{equation*}
    	if and only if
    	\begin{equation*}
    		\delta + \delta_1 > \abs{\frac{q}{p} - \alpha}.
    	\end{equation*}
    	Thanks to \eqref{condition_p_q}, it is possible to choose $ \delta $ and $ \delta_1 $ in such a way, yielding
    	\begin{equation} \label{bound_Xi33}
       		\| \Xi_{3,3}(s,r,v_1,v_2,\cdot) \|_{p,q} \leq  \frac{C \left( e^{-\genericc s} + e^{-\genericc (s-r)} \right) e^{-\genericc r}}{((s \wedge 1) ((s-r) \wedge 1))^{1/4} (r \wedge 1)^{\delta'/4}}  \, \dist(v_1,M)^{\delta'}.
    	\end{equation}
    	
    	We now turn to $ \Xi_{3,4} $.
    	By \eqref{bound_Xi0} in Lemma~\ref{lemma:bounds_pst2}, for any $ \delta \in (0,\frac{1}{2} \wedge (1-\frac{\alpha}{2})) $, there exist $ C > 0 $ and $ \genericc > 0 $ such that
    	\begin{multline*}
    		| \Xi_{3,4}(s,r,u,v_1,v_2,y) | \leq \frac{C \left( e^{-\genericc (s-u)} + e^{-\genericc (s-r)} \right)}{(((s-u)\wedge 1)((s-r)\wedge 1))^{1/4}} \int_\R |c_u(v_2,z)| p_{0,u}(v_2,z,y) e^{-\delta |z| + \frac{\alpha}{2} z} dz \\ \times \int_\R |c_r(v_1,z)| p_{0,r}(v_1,z,y) e^{-\delta |z| + \frac{\alpha}{2} z} dz.
    	\end{multline*}
    	Using \eqref{bound_cr_p0r}, for any $ \delta_1 \in (0,\delta) $, there exists $ C > 0 $, $ \genericc > 0 $ and $ \delta' \in (0,1) $ such that
    	\begin{multline*}
    	 	| \Xi_{3,4}(s,r,u,v_1,v_2,y) | \leq \frac{C \left( e^{-\genericc (s-u)} + e^{-\genericc (s-r)} \right) e^{-\genericc (u + r)}}{(((s-u)\wedge 1)((s-r)\wedge 1))^{1/4} ((u\wedge 1)(r\wedge 1))^{\delta'/4}} \\ \times \dist(v_1,M)^{\delta'} \dist(v_2,M)^{\delta'} e^{-2\delta_1 |y| - \alpha y}.
    	\end{multline*}
    	Thanks to \eqref{condition_p_q}, it is possible to choose $ \delta_1 $ such that
    	\begin{equation*}
    		2 \delta_1 > \abs{\frac{q}{p} - \alpha},
    	\end{equation*}
    	which ensures that
    	\begin{multline} \label{bound_Xi34}
    		\| \Xi_{3,4}(s,r,u,v_1,v_2,\cdot) \|_{p,q} \leq  \frac{C \left( e^{-\genericc (s-u)} + e^{-\genericc (s-r)} \right) e^{-\genericc (u + r)}}{(((s-u)\wedge 1)((s-r)\wedge 1))^{1/4} ((u\wedge 1)(r\wedge 1))^{\delta'/4}} \\ \times \dist(v_1,M)^{\delta'} \dist(v_2,M)^{\delta'}.
    	\end{multline}
    	
    	Combining \eqref{bound_Xi33} and \eqref{bound_Xi34} and returning to \eqref{decomp_Xi_31_small_times}, we obtain that there exist $ C > 0 $ and $ \genericc > 0 $ such that
    	\begin{equation*}
    		\| \Xi_{3,1}(s,v_1,v_2,\cdot) \|_{p,q} \leq \frac{C e^{-\genericc s}}{(s \wedge 1)^{1/4}} \dist(v_1,M)^{\delta'} + C e^{-\genericc s} \dist(v_1,M)^{\delta'} \dist(v_2,M)^{\delta'}.
    	\end{equation*}
    	The bound \eqref{bound_Xi3} then follows from \eqref{decomp_Xi_3}.
    \end{proof}

	\section{The linearised semigroup} \label{sec:semigroup}
	
	In this Section, we prove all the lemmas stated in Section~\ref{sec:frechet} concerning the linearised semigroup $ T_{s,t,v} $.
	We begin by proving some basic estimates on $ p_{s,t}(v,\cdot,\cdot) $, and stating some more delicate estimates which will be proved later.
	
	\subsection[First estimates]{First estimates on $ p_{s,t}(v,\cdot,\cdot) $} \label{subsec:basic_estimates}
	
	For $ t > 0 $ and $ x \in \R $, let $ G_t(x) := \frac{1}{\sqrt{4\pi t}} \exp\left( - \frac{x^2}{4t} \right) $.
	
	\begin{lemma} \label{lemma:bound_pst}
		For all $ 0 \leq s \leq t $ and $ v \in L^{2,\alpha} $,
		\begin{equation} \label{positivity_pst}
			p_{s,t}(v,z,y) \geq 0.
		\end{equation}
		Moreover, for any $ T > 0 $, there exists a constant $ C_T > 0 $ such that, for all $ 0 \leq s \leq t $ with $ t-s \leq T $ and for all $ v \in L^{2,\alpha} $, $ z, y \in \R $,
		\begin{align*}
			p_{s,t}(v,z,y) \leq C_T\, G_{t-s}(y-z) e^{- \frac{\alpha}{2} (y+z)}.
		\end{align*}
	\end{lemma}
	
	\begin{proof}
		By the Feynman-Kac formula, $ u_t $ in \eqref{linear_equation} can also be written
		\begin{align*}
			u_t(z) = \E{ \exp\left( \int_{s}^{t} f'(\Phi_{t-u}(v,X_u)) du \right) h(X_t) }{X_s = z},
		\end{align*}
		where $ X_u = X_s + \alpha (u-s) + \sqrt{2} (B_u - B_s) $ for $u\ge s$ and $ (B_t, t \geq 0) $ is standard Brownian motion.
		As a result, we can write
		\begin{align*}
			p_{s,t}(v,z,y) = G_{t-s}(y-z-\alpha (t-s)) e^{-\alpha y} \E{ \exp \left( \int_{s}^{t} f'\Big(\Phi_{t-u}\Big(v,z + \frac{u-s}{t-s} (y-z) + \sqrt{2} \bm{B}_u\Big)\Big) du \right) },
		\end{align*}
		where $ (\bm{B}_u, u \in [s, t]) $ is a standard Brownian bridge.
		This already yields the positivity of $ p_{s,t}(v,z,y) $.
		We can then bound $ f'(\Phi_{t-u}(v,\cdot)) $ by $ \sup_{u \in [0,1]} f'(u) $ and use the fact that
		\begin{equation} \label{identity_G_alpha}
			G_{t}(y-z-\alpha t) = e^{-\frac{\alpha^2}{4} t} G_t(y-z) e^{\frac{\alpha}{2}(y-z)}
		\end{equation}
		to obtain the result.
	\end{proof}

	\begin{lemma} \label{lemma:grad_pst}
		For any $ T > 0 $, there exists a constant $ C_T > 0 $ such that, for all $ 0 \leq s \leq t $ with $ t-s \leq T $ and all $ v \in L^{2,\alpha} $, $ x, y \in \R $,
		\begin{equation*}
			| \partial_y p_{s,t}(v,x,y) | \leq C_T \left( 1 + \frac{|y-x|}{t-s} \right) G_{t-s}(y-x) e^{-\frac{\alpha}{2}(x + y)}.
		\end{equation*}
	\end{lemma}
	
	\begin{proof}
		First note that, by the semigroup property,
		\begin{equation*}
			\partial_s T_{s,t,v} \phi = - (\partial_{xx} + \alpha \partial_{x} + f'(\Phi_s(v,\cdot))) T_{s,t,v} \phi.
		\end{equation*}
		As a result,
		\begin{equation*}
			T_{s,t,v} \phi = Q(t-s) \phi + \int_{s}^{t} Q(r-s) \lbrace f'(\Phi_r(v,\cdot)) T_{r,t,v} \phi \rbrace dr.
		\end{equation*}
		This translates to
		\begin{equation*}
			p_{s,t}(v,x,y) = G_{t-s}(y + \alpha (t-s) - x) e^{-\alpha x} + \int_{s}^{t} \int_\R G_{r-s}(y + \alpha (r-s) - z) f'(\Phi_r(v,z)) p_{r,t}(v,x,z) dz dr.
		\end{equation*}
		Using the fact that
		\begin{equation*}
			\partial_x G_t(x) = - \frac{x}{2t} G_t(x),
		\end{equation*}
		we obtain
		\begin{multline*}
			\partial_y p_{s,t}(v,x,y) = - \left( \frac{y-x}{2(t-s)} + \frac{\alpha}{2} \right) G_{t-s}(y + \alpha(t-s) - x) e^{-\alpha x} \\ - \int_{s}^{t} \int_\R \left( \frac{y-z}{2(r-s)} + \frac{\alpha}{2} \right) G_{r-s}(y+ \alpha(r-s) - z) f'(\Phi_r(v,z)) p_{r,t}(v,x,z) dz dr.
		\end{multline*}
		Then, using \eqref{identity_G_alpha} and Lemma~\ref{lemma:bound_pst}, we obtain that there exists $ C > 0 $ such that, for $ t-s \leq T $,
		\begin{multline} \label{bound_grad_pst}
			| \partial_y p_{s,t}(v,x,y) | \leq C \left( 1 + \frac{|y-x|}{t-s} \right) G_{t-s}(y-x) e^{-\frac{\alpha}{2}(y + x)} \\ + C \int_{s}^{t} \int_\R \left( 1 + \frac{|y-z|}{r-s} \right) G_{r-s}(y-z) G_{t-r}(z-x) e^{-\frac{\alpha}{2}(x + y)} dz dr.
		\end{multline}
		We then note that
		\begin{equation*}
			G_{r-s}(y-z) G_{t-r}(z-x) = G_{t-s}(y-x) G_{\frac{(r-s)(t-r)}{t-s}}\left( z - \frac{(t-r)y + (r-s) x}{t-s} \right).
		\end{equation*}
		By a change of variables,
		\begin{multline*}
			\int_\R \left( 1 + \frac{|y-z|}{r-s} \right) G_{\frac{(r-s)(t-r)}{t-s}}\left( z - \frac{(t-r)y + (r-s) x}{t-s} \right) dz \\
			\begin{aligned}
				&= \int_\R \left( 1 + \frac{1}{r-s} \left| \frac{r-s}{t-s} (y - x) - z \right| \right) G_{\frac{(r-s)(t-r)}{t-s}}(z) dz. \\
				&\leq 1 + \frac{|y-x|}{t-s} + C \sqrt{\frac{t-r}{(r-s)(t-s)}}.
			\end{aligned}
		\end{multline*}
		Plugging this in \eqref{bound_grad_pst} and integrating over $ r $ yields the existence of a constant $ C_T > 0 $ such that, for $ t-s \leq T $,
		\begin{equation*}
			| \partial_y p_{s,t}(v,x,y) | \leq C_T \left( 1 + \frac{|y-x|}{t-s} \right) G_{t-s}(y-x) e^{-\frac{\alpha}{2}(y + x)},
		\end{equation*}
		as stated.
	\end{proof}
	
	The following lemma is crucial in the proof of Lemma~\ref{lemma:bounds_pst2}.
	Its proof is fairly technical and relies on a connection between $ p_{s,t}(m,\cdot,\cdot) $ and the transition density of a one-dimensional diffusion.
	The proof is detailed in Subsection~\ref{subsec:diffusion}.
	
	\begin{lemma} \label{lemma:long_time_pstm}
		There exists $ C > 0 $ such that, for all $ v \in \mathcal{V}_{\Beta{zeta},K,\Epsilon{zeta}} $, for all $ 0 \leq s \leq t $ and $ x, y \in \R $,
		\begin{equation} \label{bound_pstv_pstm}
			p_{s,t}(v,x,y) \leq C p_{0,t-s}(m_{\zeta(v)},x,y).
		\end{equation}
		Furthermore, for any $ \delta_1 $ and $ \delta_2 $ such that
		\begin{equation} \label{condition_delta_1_2}
			0 < \delta_1 < \delta_2 < 1-\frac{\alpha}{2},
		\end{equation}
		and for any $ K > 0 $, there exist constants $ C > 0 $ and $ \genericc > 0 $ such that, for any $ \zeta \in \R $ with $ | \zeta | \leq \Cst{sup_zeta}(K) $, any $ t \geq 1 $ and any $ x, y \in \R $,
		\begin{equation} \label{bound_pst_varphi}
			| p_{0,t}(m_\zeta,x,y) - \varphi_{\zeta}(x) \varphi_{\zeta}(y) | \leq C e^{-\genericc t} e^{-\delta_1 |y| + \delta_2 |x| - \frac{\alpha}{2}(x + y)}.
		\end{equation}
	\end{lemma}

	\subsection{Short time bounds} \label{subsec:Tst_short_times}
	
	Let us start by proving Lemma~\ref{lemma:Tstv_Hgamma}.
	
	\begin{proof}[Proof of Lemma~\ref{lemma:Tstv_Hgamma}]
		First note that, by the definition of $ T_{s,t,v} $ in Subsection~\ref{subsec:semigroup}, using the semigroup property \eqref{semigroup_ppty},
		\begin{equation*}
			\partial_s T_{s,t,v} \phi = - (\partial_{xx} + \alpha \partial_{x} + f'(\Phi_s(v,\cdot))) T_{s,t,v} \phi.
		\end{equation*}
		Hence, by the variation of constants formula,
		\begin{equation} \label{Tstv_Q}
			T_{s,t,v} \phi = Q(t-s) \phi + \int_{s}^{t} Q(r-s) \left( f'(\Phi_r(v,\cdot)) T_{r,t,v} \phi \right) dr,
		\end{equation}
		where $ (Q(t), t \geq 0) $ denotes the semigroup generated by $ (\partial_{x x} + \alpha \partial_{x}) $, see \eqref{eq:Qdefn}.
		From \eqref{Tstv_Q} and Lemma~\ref{lemma:Qt}, we obtain
		\begin{equation} \label{Hnorm_Tstv}
			\Hnorm{T_{s,t,v} \phi}{\gamma} \leq \frac{C_\gamma}{((t-s)\wedge1)^{\gamma/2}} \| \phi \|_{2,\alpha} + \int_{s}^{t} \frac{C_\gamma}{((r-s)\wedge 1)^{\gamma/2}} \| f' \|_\infty \| T_{r,t,v} \phi \|_{2,\alpha} dr.
		\end{equation}
		In addition, by Lemma~\ref{lemma:bound_pst},
		\begin{equation} \label{bound_T_Q}
			| T_{r,t,v} \phi(y) | \leq C_T | Q(t-r) \phi(y) |,
		\end{equation}
		Combined with Lemma~\ref{lemma:Qt}, this yields
		\begin{equation*}
			\| T_{r,t,v} \phi \|_{2,\alpha} \leq C_T \| \phi \|_{2,\alpha}.
		\end{equation*}
		Plugging this in \eqref{Hnorm_Tstv} and integrating over $ r $ yields \eqref{Tstv_Hgamma}.
		Turning to \eqref{Tstv-v'_Hgamma}, we write, using \eqref{Tstv_Q},
		\begin{multline} \label{diff_Tstv}
			(T_{s,t,v} - T_{s,t,v'}) \phi = \int_{s}^{t} Q(r-s) \left( \left( f'(\Phi_r(v,\cdot)) - f'(\Phi_r(v',\cdot)) \right) T_{r,t,v} \phi \right) dr \\ + \int_{s}^{t} Q(r-s) \left( f'(\Phi_r(v,\cdot)) \left( T_{r,t,v}\phi - T_{r,t,v'} \phi \right) \right) dr.
		\end{multline}
		Hence, using Lemma~\ref{lemma:Qt},
		\begin{multline*}
			\Hnorm{(T_{s,t,v} - T_{s,t,v'}) \phi}{\gamma} \leq \int_{s}^{t} \frac{C_\gamma}{((r-s)\wedge 1)^{\gamma/2}} \| \left( f'(\Phi_r(v,\cdot)) - f'(\Phi_r(v',\cdot)) \right) T_{r,t,v} \phi \|_{2,\alpha} dr \\ + \int_{s}^{t} \frac{C_\gamma}{((r-s)\wedge 1)^{\gamma/2}} \| f' \|_\infty \| (T_{r,t,v} - T_{r,t,v'}) \phi \|_{2,\alpha} dr.
		\end{multline*}
		In addition, by Lemma~\ref{lemma:bound_pst},
		\begin{equation*}
			\| \left( f'(\Phi_r(v,\cdot)) - f'(\Phi_r(v',\cdot)) \right) T_{r,t,v} \phi \|_{2,\alpha} \leq C_T \| f'' \|_\infty \| \Phi_r(v,\cdot) - \Phi_r(v',\cdot) \|_{2,\alpha} \| Q(t-r) |\phi| \|_\infty.
		\end{equation*}
		Combined with \eqref{bounds_Phit_Hgamma} (with $ \gamma = 0 $) and the fact that $ \| Q(t)\phi \|_\infty \leq \| \phi \|_\infty $, we obtain that, for some constant $ C > 0 $ (depending only on $ T $ and $ \gamma $),
		\begin{equation} \label{Tstv-v'_pre_Gronwall}
			\Hnorm{(T_{s,t,v} - T_{s,t,v'}) \phi}{\gamma} \leq C \| v - v' \|_{2,\alpha} \| \phi \|_\infty + \int_{s}^{t} \frac{C_\gamma}{((r-s)\wedge 1)^{\gamma/2}} \| f' \|_\infty \| (T_{r,t,v} - T_{r,t,v'}) \phi \|_{2,\alpha} dr.
		\end{equation}
		Now take $ \gamma = 0 $ in the above equation to obtain, by Gronwall's inequality,
		\begin{equation*}
			\| (T_{s,t,v} - T_{s,t,v'}) \phi \|_{2,\alpha} \leq C \| v - v' \|_{2,\alpha} \| \phi \|_\infty,
		\end{equation*}
		for all $ 0 \leq s \leq t \leq T $ for some $ C > 0 $ depending only on $ T $ and $ \gamma $.
		Plugging this back into the right-hand side of \eqref{Tstv-v'_pre_Gronwall} yields \eqref{Tstv-v'_Hgamma}.
	\end{proof}
	
	We now prove Lemma~\ref{lemma:Tstv_Lpq_short_times}.

	\begin{proof}[Proof of Lemma~\ref{lemma:Tstv_Lpq_short_times}]
		The bound \eqref{Tstv_Lpq_short_time} follows from \eqref{Tstv_Q}, \eqref{bound_T_Q} and Lemma~\ref{lemma:Q_Lpalpha} (with $ r = 1 $).
		From \eqref{diff_Tstv} and Lemma~\ref{lemma:Q_Lpalpha}, we obtain
		\begin{multline*}
			\| (T_{s,t,v}-T_{s,t,v'}) \phi \|_{p,q} \leq C \int_{s}^{t} \| (f'(\Phi_r(v,\cdot))-f'(\Phi_r(v',\cdot))) T_{r,t,v} \phi \|_{p,q} dr \\ + C \| f' \|_\infty \int_{s}^{t} \| (T_{r,t,v}-T_{r,t,v}) \phi \|_{p,q} dr.
		\end{multline*}
		By the Cauchy-Schwarz inequality,
		\begin{align*}
			\| (f'(\Phi_r(v,\cdot))-f'(\Phi_r(v',\cdot))) T_{r,t,v} \phi \|_{p,q} &\leq \| f'' \|_\infty \| \Phi_r(v,\cdot) - \Phi_r(v',\cdot) \|_{2p,p\alpha} \| T_{r,t,v} \phi \|_{2p,2q-p\alpha} \\
			&\leq \frac{C}{r^{\frac{1}{4}\left( 1-\frac{1}{p} \right)}} \| v - v' \|_{2,\alpha} \| \phi \|_{2p,2q-p\alpha},
		\end{align*}
		using Corollary~\ref{cor:Phit_Lpalpha}, \eqref{bound_T_Q} and Lemma~\ref{lemma:Q_Lpalpha} in the second line.
		Note also that, since $ \phi \in L^{p,q} $, $ r \mapsto \| (T_{r,t,v}-T_{r,t,v}) \phi \|_{p,q} $ is locally bounded by Lemma~\ref{lemma:Q_Lpalpha}.
		Hence, by Gronwall's inequality, we obtain that there exists $ C > 0 $ such that, for all $ 0 \leq s \leq t \leq T $,
		\begin{equation*}
			\| (T_{s,t,v}-T_{s,t,v'}) \phi \|_{p,q} \leq C \| v-v' \|_{2,\alpha} \| \phi \|_{2p,2q-p\alpha},
		\end{equation*}
		yielding \eqref{diff_Tstv_Lpq}.
	\end{proof}
	
	Let us now prove Lemma~\ref{lemma:diff_pst2_short_times}.
	
	\begin{proof}[Proof of Lemma~\ref{lemma:diff_pst2_short_times}]
		We start by proving \eqref{diff_p0t_v_v'}.
		By Lemma~\ref{lemma:bound_pst}, there exists $ C $ such that, for $ t \in [0,T] $,
		\begin{equation} \label{bound_f0t}
			\int_\R \langle |\phi|, |p_{0,t}(v,\cdot,y) - p_{0,t}(v',\cdot,y)| p_{0,t}(v,\cdot,y) \rangle_\alpha e^{qy} dy \leq \frac{C}{\sqrt{t}} f(0,t),
		\end{equation}
		where
		\begin{equation*}
			f(s,t) := \int_{\R^2} |\phi(x)| |p_{s,t}(v,x,y) - p_{s,t}(v',x,y)| e^{(q-\frac{\alpha}{2}) y + \frac{\alpha}{2} x} dx dy.
		\end{equation*}
		Using Lemma~\ref{lemma:bound_pst} again, we note that
		\begin{align*}
			f(s,t) &\leq C \int_{\R^2} |\phi(x)| G_{t-s}(x-y) e^{(q-\alpha) y} dx dy \\
			&\leq C' \| \phi \|_{1,q-\alpha}.
		\end{align*}
		Hence $ s \mapsto f(s,t) $ is locally bounded on $ [0,t] $ for any $ t \geq 0 $ since $ \phi \in L^{1,q-\alpha} $.
		Then, from \eqref{diff_Tstv} and the definition of $ Q(t) $ in \eqref{eq:Qdefn} we deduce that
		\begin{multline*}
			p_{s,t}(v,x,y) - p_{s,t}(v',x,y) \\ = \int_{s}^{t} \int_\R G_{r-s}(y+\alpha (r-s) - z) \left( f'(\Phi_r(v,z)) - f'(\Phi_r(v',z)) \right) p_{r,t}(v,x,z) dz dr \\ + \int_{s}^{t} \int_\R G_{r-s}(y + \alpha (r-s) - z) f'(\Phi_r(v',z)) \left( p_{r,t}(v,x,z) - p_{r,t}(v',x,z) \right) dz dr.
		\end{multline*}
		Using Lemma~\ref{lemma:bound_pst} once again and \eqref{identity_G_alpha}, we obtain that there exists $ C > 0 $ such that, for $ 0 \leq s \leq t \leq T $,
		\begin{multline*}
			| p_{s,t}(v,x,y) - p_{s,t}(v',x,y) | \leq C \int_{s}^{t} \int_\R |\Phi_r(v,z) - \Phi_r(v',z) | G_{r-s}(y-z) G_{t-r}(z-x) e^{-\frac{\alpha}{2}(y + x)} dz dr \\ + C \int_{s}^{t} \int_\R | p_{r,t}(v,x,z) - p_{r,t}(v',x,z) | G_{r-s}(y-z) e^{\frac{\alpha}{2}(z-y)} dz dr.
		\end{multline*}
		Plugging this in the definition of $ f(s,t) $, we obtain
		\begin{multline} \label{bound_fst_pre_Gronwall}
			f(s,t) \leq C \int_{s}^{t} \int_{\R^3} |\phi(x)| \, |\Phi_r(v,z) - \Phi_r(v',z)| G_{r-s}(y-z) G_{t-r}(z-x) e^{(q-\alpha) y} dx dy dz dr \\ + C \int_{s}^{t} \int_{\R^3} |\phi(x)| |p_{r,t}(v,x,z) - p_{r,t}(v',x,z)| G_{r-s}(y-z) e^{(q-\alpha) y + \frac{\alpha}{2} (x-z)} dx dy dz dr.
		\end{multline}
		By the Cauchy-Schwarz inequality, we can bound the integrand in the first term as follows
		\begin{multline*}
			\int_{\R^3} |\phi(x)| \, |\Phi_r(v,z) - \Phi_r(v',z)| G_{r-s}(y-z) G_{t-r}(z-x) e^{(q-\alpha) y} dx dy dz \\ \leq \left( \int_{\R^2} | \phi(x) |^2 G_{t-s}(x-y) e^{(2q-3\alpha) y} dx dy \right)^{1/2} \\ \times \frac{1}{(4\pi (t-s))^{1/4}} \left( \int_{\R^2} | \Phi_r(v,z) - \Phi_r(v',z) |^2 G_{\frac{(t-r)(r-s)}{t-s}}(y-z) e^{\alpha y} dy dz \right)^{1/2},
		\end{multline*}
		where we have used the identity
		\begin{equation*}
			G_{t}(x-y) G_{s}(y-z) = G_{t+s}(x-z) G_{\frac{st}{s + t}}\left( y - \frac{sx + tz}{s + t} \right)
		\end{equation*}
		in each term.
		Integrating over $ y $ in each integral, we obtain that there exists $ C > 0 $ such that, for $ 0 \leq s \leq r \leq t \leq T $,
		\begin{multline*}
			\int_{\R^3} |\phi(x)| \, |\Phi_r(v,z) - \Phi_r(v',z)| G_{r-s}(y-z) G_{t-r}(z-x) e^{(q-\alpha) y} dx dy dz \\ \leq \frac{C}{(t-s)^{1/4}} \| \phi \|_{2,2q-3\alpha} \| \Phi_r(v,\cdot) - \Phi_r(v',\cdot) \|_{2,\alpha}.
		\end{multline*}
		Plugging this in \eqref{bound_fst_pre_Gronwall} and using \eqref{bounds_Phit_Hgamma} (with $ \gamma = 0 $), we obtain, also integrating over $ y $ in the second term,
		\begin{equation*}
			f(s,t) \leq C (t-s)^{3/4} \| \phi \|_{2,2q-3\alpha} \| v - v' \|_{2,\alpha} + C \int_{s}^{t} f(r,t) dr.
		\end{equation*}
		By Gronwall's lemma, we obtain that there exists $ C > 0 $ such that, for all $ 0 \leq s \leq t \leq T $,
		\begin{equation*}
			f(s,t) \leq C (t-s)^{3/4} \| \phi \|_{2,2q-3\alpha} \| v - v' \|_{2,\alpha}.
		\end{equation*}
		Plugging this in \eqref{bound_f0t} yields \eqref{diff_p0t_v_v'}.
		We now turn to the proof of \eqref{phi_p0t_delta} and we write
		\begin{multline*}
			p_{0,t}(v,x,y + \delta z_1) p_{0,t}(v,x,y + \delta z_2) - p_{0,t}(v,x,y)^2 \\ = \delta \int_{0}^{1} \Big\lbrace z_1 \partial_y p_{0,t}(v,x,y + u \delta z_1) p_{0,t}(v,x,y + u \delta z_2) \\ + z_2 p_{0,t}(v,x,y + u \delta z_1) \partial_y p_{0,t}(v,x,y + u \delta z_2) \Big\rbrace du.
		\end{multline*}
		Hence, by Lemma~\ref{lemma:grad_pst} and Lemma~\ref{lemma:bound_pst},
		\begin{multline*}
			| p_{0,t}(v,x,y + \delta z_1) p_{0,t}(v,x,y + \delta z_2) - p_{0,t}(v,x,y)^2 | \\ \leq C \delta \int_{0}^{1} \left\lbrace |z_1| \left( 1 + \frac{| y + u \delta z_1 - x |}{t} \right) + |z_2| \left( 1 + \frac{| y + u \delta z_2 - x |}{t} \right) \right\rbrace \\ \times  G_t(y + u \delta z_1 - x) G_t(y + u \delta z_2 - x) e^{-\alpha(y + x + \frac{u \delta}{2} (z_1 + z_2))} du.
		\end{multline*}
		Then, using the identity
		\begin{equation*}
			G_t(y + u \delta z_1 - x) G_t(y + u \delta z_2 - x) = G_{t/2}\left(y-x + \frac{u\delta}{2}(z_1+z_2)\right) G_{2t}\left( u \delta (z_1 - z_2) \right),
		\end{equation*}
		and a change of variables, we obtain
		\begin{multline*}
			\int_{[-1,1]^2 \times \R} | p_{0,t}(v,x,y + \delta z_1) p_{0,t}(v,x,y + \delta z_2) - p_{0,t}(v,x,y)^2 | \rho(z_1) \rho(z_2) e^{qy} dz_1 dz_2 dy \\ \leq C \delta e^{(q-2\alpha)x} \int_{0}^{1} \int_{[-1,1]^2 \times \R} \left\lbrace | z_1 | \left( 1 + \frac{|y + \frac{u \delta}{2} (z_1-z_2)|}{t} \right) + | z_2 | \left( 1 + \frac{| y + \frac{u \delta}{2} (z_2-z_1) |}{t} \right) \right\rbrace \\ \times G_{t/2}(y) G_{2t}(u \delta (z_1-z_2)) \rho( z_1) \rho(z_2) e^{(q-\alpha)y} e^{-\frac{q u \delta}{2} (z_1 + z_2)} dz_1 dz_2 dy du.
		\end{multline*}
		We then note that, for any $ T \geq 0 $, there exists $ C > 0 $ such that, for all $ t \in (0,T] $,
		\begin{align*}
			\int_\R G_{t/2}(y) e^{(q-\alpha)y} dy \leq C, && \int_\R \frac{|y|}{t} G_{t/2}(y) e^{(q-\alpha)y} dy \leq \frac{C}{\sqrt{t}}.
		\end{align*}
		In addition, there exists $ C > 0$ such that, for all $ t \in (0,T] $,
		\begin{equation*}
			\int_{[-1,1]^2} G_{2t}(u \delta (z_1-z_2)) \rho(z_1) \rho(z_2) dz_1 dz_2 \leq \frac{C}{\sqrt{t}}.
		\end{equation*}
		Moreover, by a change of variables, and using the fact that $ z \mapsto \rho(z) $ is bounded,
		\begin{equation*}
			\int_{[-1,1]^2} \frac{u\delta}{t} |z_1-z_2| G_{2t}(u \delta (z_1-z_2)) \rho(z_1) \rho(z_2) dz_1 dz_2 \leq \frac{C}{u \delta},
		\end{equation*}
		and
		\begin{equation*}
			\int_{[-1,1]^2} \frac{u\delta}{t} |z_1-z_2| G_{2t}(u \delta (z_1-z_2)) \rho(z_1) \rho(z_2) dz_1 dz_2 \leq \frac{C}{\sqrt{t}}.
		\end{equation*}
		As a result, there exists $ C > 0 $ such that, for all $ t \in (0,T] $,
		\begin{multline*}
			\int_{[-1,1]^2 \times \R} | p_{0,t}(v,x,y + \delta z_1) p_{0,t}(v,x,y + \delta z_2) - p_{0,t}(v,x,y)^2 | \rho(z_1) \rho(z_2) e^{qy} dz_1 dz_2 dy \\ 
			\begin{aligned}
				&\leq C \delta e^{(q-2\alpha)x} \int_{0}^{1} \frac{1}{\sqrt{t}} \left\lbrace  \left( \frac{1}{\sqrt{t}} \wedge \frac{1}{u \delta} \right) + 1 \right\rbrace du \\
				&\leq C \delta e^{(q-2\alpha)x} \frac{1}{\sqrt{t}} \left\lbrace \left(\frac{1}{\delta} \wedge \frac{1}{\sqrt{t}}\right) - \frac{1}{\delta} \log\left(\frac{\sqrt{t}}{\delta} \wedge 1\right) + 1 \right\rbrace,
			\end{aligned}
		\end{multline*}
		It follows that
		\begin{multline*}
			\int_{[-1,1]^2 \times \R} | \langle \phi, p_{0,t}(v,\cdot,y + \delta z_1) p_{0,t}(v,\cdot,y + \delta z_2) \rangle_\alpha - \langle \phi, p_{0,t}(v,\cdot,y)^2 \rangle_\alpha | \rho(z_1) \rho(z_2) e^{qy} dz_1 dz_2 dy \\ \leq C \delta \| \phi \|_{1,q-\alpha} \frac{1}{\sqrt{t}} \left\lbrace \left(\frac{1}{\delta} \wedge \frac{1}{\sqrt{t}}\right) - \frac{1}{\delta} \log\left(\frac{\sqrt{t}}{\delta} \wedge 1\right) + 1 \right\rbrace.
		\end{multline*}
		Integrating over $t$ yields the result.
	\end{proof}
	
	In passing, we can prove Lemma~\ref{lemma:Tstv_time_shift}.
	
	\begin{proof}[Proof of Lemma~\ref{lemma:Tstv_time_shift}]
		First note that, by \eqref{Tstv_Q}, for any $ 0 \leq s \leq t $ and $ \tau \geq 0 $, by a change of variables and the fact that $ \Phi_r(\Phi_\tau(v,\cdot),\cdot) = \Phi_{r + \tau}(v,\cdot) $,
		\begin{equation*}
			T_{s+\tau,t + \tau, v} \phi = Q(t-s) \phi + \int_{s}^{t} Q(r-s) \left( f'(\Phi_{r +\tau}(v,\cdot)) T_{r + \tau, t + \tau, v} \phi \right) dr.
		\end{equation*}
		Using \eqref{Tstv_Q} again with $ \Phi_\tau(v,\cdot) $ instead of $ v $, we obtain that
		\begin{equation*}
			T_{s+\tau,t+\tau,v} = T_{s,t,\Phi_t(v,\cdot)},
		\end{equation*}
		yielding the result.
	\end{proof}
	
	\subsection{Large time bounds} \label{subsec:Tst_large_times}
	
	The aim of this subsection is to prove Lemma~\ref{lemma:T-m}, Lemma~\ref{lemma:Tst_L2alpha}, Lemma~\ref{lemma:Tstv_Hgamma_large_times}, Lemma~\ref{lemma:pst_phi_delta} and Lemma~\ref{lemma:bounds_pst2}, which all deal with the large time behaviour of $ T_{s,t,v} $.
	Lemma~\ref{lemma:bound_pst}, combined with Lemma~\ref{lemma:T-A}, allows us to prove Lemma~\ref{lemma:T-m}.
	
	\begin{proof}[Proof of Lemma~\ref{lemma:T-m}]
		First note that, since
		\begin{equation*}
			\inf_{x \in \R} |\partial_{x} m(x)| e^{|x|} =: c > 0,
		\end{equation*} 
		for any $ \eta \in \R $ and $ \phi \in L^{\infty,e} $,
		\begin{equation} \label{upper_bound_grad_m}
			|\phi(x)| \leq c \, e^{|\eta|}  \enorm{\phi} |\partial_x m_{\eta} (x) |, \quad \forall x \in \R.
		\end{equation}
		By \eqref{upper_bound_grad_m} and the positivity of $ T_{s,t,v} $ and $ \partial_{x} m < 0 $, we can write
		\begin{equation*}
			|T_{s,t,v} \phi(y)| \leq c e^{| \zeta(v)|} \enorm{\phi} |T_{s,t,v} \partial_{x} m_{\zeta(v)}(y) |, \quad \forall y \in \R.
		\end{equation*}
		Let us then check that $ t \mapsto \enorm{T_{s,t,v} \partial_{x} m_{\zeta(v)}} $ is locally bounded for any $ y \in \R $.
		Using Lemma~\ref{lemma:bound_pst}, for $ t \in [s, s + T] $,
		\begin{equation*}
			| T_{s,t,v} \partial_{x} m_{\zeta(v)} (y) | \leq C_T e^{-\frac{\alpha}{2} y} \int_{\R} G_{t-s}(y-z) | \partial_{x} m_{\zeta(v)}(z) | e^{\frac{\alpha}{2} z} dz.
		\end{equation*}
		Since $ \enorm{\partial_{x} m_{\zeta(v)}} < \infty $, we have
		\begin{equation*}
			| T_{s,t,v} \partial_{x} m_{\zeta(v)} (y) | \leq C_T \enorm{\partial_{x} m_{\zeta(v)}} e^{-\frac{\alpha}{2}y} \int_{\R} G_{t-s}(y-z) e^{-|z| + \frac{\alpha}{2} z} dz.
		\end{equation*}
		By a change of variables and the fact that $ |y| \leq |y + z| + |z| $, we obtain that $ t \mapsto \enorm{T_{s,t,v} \partial_{x} m_{\zeta(v)}} $ is locally bounded.
		Finally, using Lemma~\ref{lemma:T-A} and noting that $ e^{-t\mathcal{A}_{\zeta(v)}} \partial_{x} m_{\zeta(v)} = \partial_{x} m_{\zeta(v)} $ for $t\ge 0$, for $y\in \R$,
		\begin{align*}
			T_{s,t,v} \partial_{x} m_{\zeta(v)} (y) = \partial_{x} m_{\zeta(v)}(y) + \int_{s}^{t} T_{s,r,v} \left\lbrace c_r(v,\cdot) \partial_{x} m_{\zeta(v)} \right\rbrace (y) dr.
		\end{align*}
		Taking absolute values on both sides and using the positivity of $ T_{s,t,v} $ (by Lemma~\ref{lemma:bound_pst}) and Corollary~\ref{cor:uniform_cvg}, we obtain
		\begin{align*}
			| T_{s,t,v} \partial_{x} m_{\zeta(v)} (y) | \leq | \partial_{x} m_{\zeta(v)} (y) | + C (\| s(v) \|_{\infty,\lambda} + \sqrt{\dist(v,M)}) \int_{s}^{t} e^{-\genericc r} | T_{s,r,v} \partial_{x} m_{\zeta(v)} | (y) dr.
		\end{align*}
		As a result,
		\begin{equation*}
			\enorm{ T_{s,t,v} \partial_{x} m_{\zeta(v)}} \leq \enorm{\partial_{x} m_{\zeta(v)}} + C (\| s(v) \|_{\infty,\lambda} + \sqrt{\dist(v,M)}) \int_{s}^{t} e^{-\genericc r} \enorm{T_{s,r,v} \partial_{x} m_{\zeta(v)} } dr
		\end{equation*}
		The result then follows by Gr\"onwall's inequality (using the fact that $ t \mapsto \enorm{T_{s,t,v} \partial_{x} m_{\zeta(v)}} $ is locally bounded).
	\end{proof}

	Next, we prove Lemma~\ref{lemma:Tst_L2alpha}.
	
	\begin{proof}[Proof of Lemma~\ref{lemma:Tst_L2alpha}]
		The first bound \eqref{L2_A} is a consequence of the spectral theorem already used in \eqref{spectral_gap}. Indeed, for $\phi \in L^{2,\alpha}$ and $t>0$,
		\begin{equation*}
			\frac{d}{dt} \| e^{-t \mathcal{A}_\eta} \phi \|_{2,\alpha}^2 = -2 \langle \mathcal{A}_\eta e^{-t \mathcal{A}_\eta} \phi, e^{-t \mathcal{A}_\eta} \phi \rangle_\alpha \leq 0.
		\end{equation*}
		We now prove~\eqref{L2_A-P}.
		Suppose $\psi \in L^{2,\alpha}$ with  $ \langle \psi, \varphi_\eta \rangle_\alpha = 0 $. Then using the fact that $ \mathcal{A}_\eta $ is self-adjoint and that (by~\eqref{eigenfunction_A}) $ \mathcal{A}_\eta \varphi_\eta = 0 $, we have that for $t>0$,
		\begin{equation*}
			\frac{d}{dt} \langle e^{-t \mathcal{A}_\eta} \psi, \varphi_\eta \rangle_\alpha = -\langle e^{-t\mathcal{A}_\eta} \psi, \mathcal{A}_\eta \varphi_\eta \rangle_\alpha 
			= 0.
		\end{equation*}
		As a result, $ \langle e^{-t \mathcal{A}_\eta} \psi, \varphi_\eta \rangle_\alpha = 0 $ for all $ t \geq 0 $.
		Hence, using \eqref{spectral_gap}, for $t>0$,
		\begin{equation*}
			\frac{d}{dt} \| e^{-t \mathcal{A}_\eta} \psi \|_{2,\alpha}^2 =-2 \langle \mathcal A_{\eta}e^{-t\mathcal A_{\eta}}\psi, e^{-t\mathcal A_{\eta}}\psi\rangle_{\alpha} \leq - 2\specgap \| e^{-t \mathcal{A}_\eta} \psi \|^2_{2,\alpha},
		\end{equation*}
		from which we deduce that
		\begin{equation} \label{eq:psibd}
			\| e^{-t \mathcal{A}_\eta} \psi \|_{2,\alpha} \leq e^{-\specgap t} \| \psi \|_{2,\alpha}.
		\end{equation}
		Then let $ \psi :=\phi - P_{\lbrace \varphi_\eta \rbrace}\phi $ for some $\phi \in L^{2,\alpha}$ and note that
		\begin{equation} \label{eq:psibd2}
			\frac{d}{dt} e^{-t \mathcal{A}_\eta} P_{\lbrace \varphi_\eta \rbrace} \phi = - \langle \phi, \varphi_\eta \rangle_\alpha e^{-t \mathcal{A}_\eta} \mathcal{A}_\eta \varphi_\eta 
			= 0
		\end{equation}
		since $\mathcal A_{\eta}\varphi_{\eta}=0$. Therefore
		$e^{-\mathcal A_{\eta}}\psi = (e^{-t\mathcal A_{\eta}}-P_{\{\varphi_{\eta}\}})\phi$, and since $\|\psi\|_{2,\alpha}\le \|\phi\|_{2,\alpha}$,
		we obtain \eqref{L2_A-P} from~\eqref{eq:psibd}.
		Let us now prove the second part of the lemma.
		Take $v\in \mathcal V_{\beta_1,K,\varepsilon_1}$ and $ \phi \in L^{2,\alpha} $, and let us first check that $ t \mapsto \| T_{s,t,v} \phi \|_{2,\alpha} $ is locally bounded.
		Using Lemma~\ref{lemma:bound_pst} and Lemma~\ref{lemma:Qt} (with $ \gamma = \gamma' = 0 $), for $ t \in [s, s + T] $,
		\begin{equation*}
			\| T_{s,t,v} \phi \|_{2,\alpha} \leq C_T \| Q(t) |\phi| \|_{2,\alpha} \leq C \| \phi \|_{2,\alpha},
		\end{equation*}
		for some $ C > 0 $ independent of $ \phi \in L^{2,\alpha} $.
		Hence $ t \mapsto \| T_{s,t,v}  \|_{L^{2,\alpha} \to L^{2,\alpha}} $ is locally bounded.
		From Lemma~\ref{lemma:T-A}, we then deduce that
		\begin{multline*}
			\| T_{s,t,v} \|_{L^{2,\alpha} \to L^{2,\alpha}} \leq \| e^{-(t-s) \mathcal{A}_{\zeta(v)}} \|_{L^{2,\alpha} \to L^{2,\alpha}}\\ + \int_{s}^{t} \| T_{s,r,v} \|_{L^{2,\alpha} \to L^{2,\alpha}} \| c_r(v,\cdot) \|_\infty \| e^{-(t-r) \mathcal{A}_{\zeta(v)}} \|_{L^{2,\alpha} \to L^{2, \alpha}} dr.
		\end{multline*}
		Using \eqref{L2_A} and Corollary~\ref{cor:uniform_cvg}, there exists a constant $C>0$ such that
		\begin{align*}
			\| T_{s,t,v} \|_{L^{2,\alpha} \to L^{2,\alpha}} \leq 1 + C \left(\| s(v) \|_{\infty,\lambda} + \sqrt{\dist(v,M)}\right) \int_{s}^{t}  e^{-\genericc r}  \| T_{s,r,v} \|_{L^{2,\alpha} \to L^{2,\alpha}} dr,
		\end{align*}
		from which \eqref{L2_T} follows by Grownall's inequality (using the fact that $ t \mapsto \| T_{s,t,v} \|_{L^{2,\alpha} \to L^{2,\alpha}} $ is locally bounded).
		Using Lemma~\ref{lemma:T-A}, we write
		\begin{multline*}
			\| (T_{s,t,v} - e^{-(t-s) \mathcal{A}_{\zeta(v)}}) P_{\lbrace \varphi_{\zeta(v)} \rbrace^\perp} \|_{L^{2,\alpha} \to L^{2,\alpha}} \\ \leq \int_{s}^{t} \| T_{s,r,v} \|_{L^{2,\alpha} \to L^{2,\alpha}} \| c_r(v,\cdot) \|_\infty \| e^{-(t-r) \mathcal{A}_{\zeta(v)}} P_{\lbrace \varphi_{\zeta(v)} \rbrace^\perp} \|_{L^{2,\alpha} \to L^{2,\alpha}} dr.
		\end{multline*}
		Inside the integral, we use \eqref{L2_T}, Corollary~\ref{cor:uniform_cvg} and \eqref{L2_A-P} (combined with~\eqref{eq:psibd2}) to obtain
		\begin{multline*}
			\| (T_{s,t,v} - e^{-(t-s) \mathcal{A}_{\zeta(v)}}) P_{\lbrace \varphi_{\zeta(v)} \rbrace^\perp} \|_{L^{2,\alpha} \to L^{2,\alpha}} \leq C \left(\| s(v) \|_{\infty,\lambda} + \sqrt{\dist(v,M)}\right) \int_{s}^{t} e^{-\genericc r} e^{-\genericc(t-r)} dr.
		\end{multline*}
		(We note that the constants denoted by $ \genericc $ in Corollary~\ref{cor:uniform_cvg} and in \eqref{L2_A-P} need not be the same, but we can always choose the smaller of the two.)
		Thus there exist constants $ \genericc' > 0 $ and $ C' > 0 $ such that
		\begin{align*}
			\| (T_{s,t,v} - e^{-(t-s) \mathcal{A}_{\zeta(v)}}) P_{\lbrace \varphi_{\zeta(v)} \rbrace^\perp} \|_{L^{2,\alpha} \to L^{2,\alpha}} \leq C' \left(\| s(v) \|_{\infty,\lambda} + \sqrt{\dist(v,M)}\right) e^{-\genericc' t},
		\end{align*}
		proving \eqref{L2_T-A}.
	\end{proof}
	
	Let us now prove Lemma~\ref{lemma:Tstv_Hgamma_large_times}.
	
	\begin{proof}[Proof of Lemma~\ref{lemma:Tstv_Hgamma_large_times}]
		Combining \eqref{Hnorm_Tstv} and \eqref{L2_T}, we obtain that there exists $ C > 0 $ such that, for all $ v \in \mathcal{V}_{\Beta{zeta},K,\Epsilon{zeta}} $,
		\begin{equation*}
			\Hnorm{T_{s,t,v} \phi}{\gamma} \leq \frac{C}{(|t-s| \wedge 1)^{\gamma/2}} \| \phi \|_{2,\alpha} + \int_{s}^{t} \frac{C}{(|r-s|\wedge 1)^{\gamma/2}} \| \phi \|_{2,\alpha} dr,
		\end{equation*}
		from which \eqref{T_H_gamma} follows
		By \eqref{Tstv_Q}, we write
		\begin{multline} \label{diff_Tstv_T0tv}
			T_{s,t,v}\phi - T_{0,t,v} \phi = (I-Q(s)) Q(t-s) \phi + \int_{s}^{t} (I - Q(s)) Q(r-s) \left( f'(\Phi_r(v,\cdot)) T_{r,t,v} \phi \right) dr \\ - \int_{0}^{s} Q(r) \left( f'(\Phi_r(v,\cdot)) T_{r,t,v} \phi  \right) dr.
		\end{multline}
		Then, by \eqref{continuity_Qt}, for any $ \gamma' > 0 $,
		\begin{align*}
			\Hnorm{(I-Q(s))Q(t-s) \phi}{\gamma} & \leq C s^{\gamma'/2} \Hnorm{Q(t-s)\phi}{\gamma+\gamma'} \\
			&\leq \frac{C s^{\gamma'/2}}{((t-s)\wedge1)^{\frac{\gamma+\gamma'}{2}}} \| \phi \|_{2,\alpha},
		\end{align*}
		where we have used \eqref{Qt_Hgamma} in the second line.
		Similarly, by \eqref{L2_T},
		\begin{equation*}
			\Hnorm{(I - Q(s)) Q(r-s) \left( f'(\Phi_r(v,\cdot)) T_{r,t,v} \phi \right)}{\gamma} \leq \frac{C s^{\gamma'/2}}{((r-s)\wedge 1)^{\frac{\gamma+\gamma'}{2}}} \| \phi \|_{2,\alpha}.
		\end{equation*}
		Finally, by \eqref{Qt_Hgamma},
		\begin{equation*}
			\Hnorm{Q(r) \left( f'(\Phi_r(v,\cdot)) T_{r,t,v} \phi  \right)}{\gamma} \leq \frac{C}{(r\wedge 1)^{\gamma/2}} \| \phi \|_{2,\alpha}.
		\end{equation*}
		Plugging the above estimates in \eqref{diff_Tstv_T0tv} and choosing $ \gamma' $ such that $ \gamma + \gamma' < 2 $, we obtain \eqref{diff_Tst-T0t_Hgamma} (using the fact that $ \frac{\gamma'}{2} < 1-\frac{\gamma}{2} $).
	\end{proof}
	
	We then prove the following result.
	\begin{lemma} \label{lemma:bounds_integral_phi}
		For any $ \delta \in [0,1] $, there exists $ C > 0 $ such that, for all $ \phi \in L^{2,\alpha} \cap L^{\infty,e} $ and all $ y \in \R $, $ t \in (0,1] $,
		\begin{equation} \label{bound_Gt_phi}
			\int_\R G_t(y-x) |\phi(x)| dx \leq \frac{C}{t^{\delta/4}} \| \phi \|_{2,\alpha}^{\delta} \enorm{\phi}^{1-\delta} e^{-(1-\delta) |y| + - \delta \frac{\alpha}{2} y}.
		\end{equation}
		Moreover, for any $ \delta_0 \in (-\infty,1) $ and any $ \delta \in [0,1] $ with
		\begin{equation} \label{condition_delta0}
			\delta < \frac{1 - \delta_0}{1 + \frac{\alpha}{2}},
		\end{equation}
		there exists $ C > 0 $ such that, for all $ \phi \in L^{2,\alpha} \cap L^{\infty,e} $,
		\begin{equation} \label{bound_phi_delta0}
			\int_\R | \phi(x) | e^{\delta_0 |x|} dx \leq C \| \phi \|_{2,\alpha}^{\delta} \enorm{\phi}^{1-\delta}.
		\end{equation}
	\end{lemma}

	\begin{proof}
		We write, for $ \delta \in [0,1] $,
		\begin{equation} \label{bound_phi_x}
			|\phi(x)| \leq \enorm{\phi}^{1-\delta} |\phi(x)|^\delta e^{-(1-\delta) |x|}.
		\end{equation}
		Then, by H\"older's inequality, setting $ r = (1-\frac{\delta}{2})^{-1} $,
		\begin{equation*}
			\int_\R G_{t}(x-y) |\phi(x)| dx \leq \| \phi \|_{2,\alpha}^{\delta} \enorm{\phi}^{1-\delta} \left( \int_\R G_{t}(x-y)^r e^{-r(1-\delta)|x| + \alpha(1-r) x} dx \right)^{1/r}.
		\end{equation*}
		We then use the fact that $ G_t(x) \leq \frac{1}{\sqrt{4\pi t}} $ to write that there exists $ C > 0 $ such that
		\begin{equation*}
			\int_\R G_{t}(x-y) |\phi(x)| dx \leq \frac{C}{t^{\delta/4}} \| \phi \|_{2,\alpha}^{\delta} \enorm{\phi}^{1-\delta} \left( \int_\R G_{t}(x-y) e^{-r(1-\delta) |x| +  \alpha(1-r) x} dx \right)^{1/r}.
		\end{equation*}
		By a change of variables and by $ |x+y| \geq |x| - |y| $, we obtain that there exists $ C > 0 $ such that, for $ t \leq 1 $,
		\begin{equation} \label{gaussian_convolution_bound}
			\int_\R G_{t}(x-y) e^{-r(1-\delta) |x| +  \alpha(1-r) x} dx \leq C e^{- r(1-\delta) |y| +  \alpha(1-r) y}.
		\end{equation}
		As a result, using the fact that $ r = (1-\frac{\delta}{2})^{-1} $, there exists $ C > 0 $ such that, for $ t \leq 1 $,
		\begin{equation*}
			\int_\R G_{t}(x-y) |\phi(x)| dx \leq \frac{C}{t^{\delta/4}} \| \phi \|_{2,\alpha}^{\delta} \enorm{\phi}^{1-\delta} e^{-(1-\delta) |y| - \delta \frac{\alpha}{2} y},
		\end{equation*}
		yielding \eqref{bound_Gt_phi}.
		We then turn to the proof of \eqref{bound_phi_delta0}.
		Using \eqref{bound_phi_x} and H\"older's inequality again, we obtain
		\begin{equation*}
			\int_\R | \phi(x) | e^{\delta_0 |x|} dx \leq \| \phi \|_{2,\alpha}^{\delta} \enorm{\phi}^{1-\delta} \left( \int_\R e^{-r(1-\delta) |x| + r \delta_0 |x | +  \alpha(1-r) x} dx \right)^{1/r},
		\end{equation*}
		with $ r = (1-\frac{\delta}{2})^{-1} $.
		Then the integral on the above right-hand side is finite if and only if
		\begin{equation*}
			1 - \delta - \delta_0 > \alpha \abs{ \frac{1}{r} - 1 } = \delta \frac{\alpha}{2}.
		\end{equation*}
		We then conclude by noting that the above is equivalent to \eqref{condition_delta0}.
	\end{proof}

	Let us continue by proving Lemma~\ref{lemma:pst_phi_delta}.
	
	\begin{proof}[Proof of Lemma~\ref{lemma:pst_phi_delta}]
		Note that, by \eqref{bound_pstv_pstm}, it is enough to treat the case $ v = m_\zeta $ with $ | \zeta | \leq \Cst{sup_zeta} $.
		Writing
		\begin{equation*}
			| \phi(x) | \leq |\phi(x)|^{\delta'} \| \phi \|_\infty^{1-\delta'}
		\end{equation*}
		and applying H\"older's inequality, we obtain that
		\begin{equation*}
			\abs{ \int_\R \phi(x) p_{0,t}(m_\zeta,x,y) e^{-\delta_0 |x| + \frac{\alpha}{2} x} dx } \leq \| \phi \|_{2,\alpha}^{\delta'} \| \phi \|_\infty^{1-\delta'} \left( \int_\R e^{-r \delta_0 |x| + \alpha(1-\frac{r}{2}) x} p_{0,t}(m_\zeta,x,y)^r dx \right)^{1/r},
		\end{equation*}
		with $ r = (1-\frac{\delta'}{2})^{-1} $.
		For $ t \in (0,1] $, we use Lemma~\ref{lemma:bound_pst} to write
		\begin{equation*}
			\left( \int_\R e^{-r \delta_0 |x| + \alpha(1-\frac{r}{2}) x} p_{0,t}(m_\zeta,x,y)^r dx \right)^{1/r} \leq C_1 e^{-\frac{\alpha}{2} y} \left( \int_\R G_t(x-y)^r e^{-r \delta_0 |x| + \alpha(1-r) x} dx \right)^{1/r}.
		\end{equation*}
		By the argument used in \eqref{gaussian_convolution_bound}, there exists $ C > 0 $ such that, for $ t \in (0,1] $,
		\begin{equation*}
			\left( \int_\R e^{-r \delta_0 |x| + \alpha(1-\frac{r}{2}) x} p_{0,t}(m_\zeta,x,y)^r dx \right)^{1/r} \leq \frac{C}{t^{\delta'/4}} e^{-\frac{\alpha}{2} y} e^{-\delta_0 |y| - \frac{\alpha \delta'}{2} y},
		\end{equation*}
		where we have used the fact that $ \frac{1}{r} = 1 - \frac{\delta'}{2} $.
		Choosing $ \delta' $ small enough that $ \delta_0 - \frac{\alpha \delta'}{2} > \delta $, we obtain the bound for $ t \in (0,1] $.

		For $ t > 1 $, using Minkowski's inequality and \eqref{bound_pst_varphi}, we obtain that, for any $ \delta_1 $, $ \delta_2 $ satisfying \eqref{condition_delta_1_2}, there exist $ C > 0 $ and $ \genericc > 0 $ such that, for all $ t \geq 1 $,
		\begin{multline*}
			\left( \int_\R e^{-r \delta_0 |x| + \alpha(1-\frac{r}{2}) x} p_{0,t}(m_\zeta,x,y)^r dx \right)^{1/r} \leq |\varphi_{\zeta}(y)| \left( \int_\R e^{-r \delta_0 |x| + \alpha(1-\frac{r}{2}) x} |\varphi_{\zeta}(x)|^r dx \right)^{1/r} \\ + C e^{-\genericc t} e^{-\delta_1 |y| - \frac{\alpha}{2} y} \left( \int_\R e^{-r(\delta_0 - \delta_2) |x| + \alpha(1-r) x} dx \right)^{1/r}.
		\end{multline*}
		Since $ |\varphi_{\zeta}(x)| \leq C e^{-|x|} $, the first integral on the right is bounded by a constant since
		\begin{equation*}
			1 + \delta_0 > \alpha \abs{\frac{1}{r} - \frac{1}{2}} = \alpha \frac{1-\delta'}{2}.
		\end{equation*}
		To ensure that the second integral on the right-hand side is finite, we need to choose $ \delta_2 < \delta_0 $ and $ \delta' > 0 $ small enough that
		\begin{equation*}
			\delta_0 - \delta_2 > \alpha \abs{\frac{1}{r} - 1} = \frac{\alpha \delta'}{2}.
		\end{equation*}
		For such a choice, using the fact that $ |\varphi_{\zeta}(y)| \leq C e^{-\delta_1 |y| - \frac{\alpha}{2} y} $, we obtain that there exists $ C > 0 $ such that, for $ t \geq 1 $,
		\begin{equation*}
			\left( \int_\R e^{-r \delta_0 |x| + \alpha(1-\frac{r}{2}) x} p_{0,t}(m_\zeta,x,y)^r dx \right)^{1/r} \leq C e^{-\delta_1 |y| - \frac{\alpha}{2} y}.
		\end{equation*}
		The result then follows by taking $ \delta_1 = \delta $.
	\end{proof}
	
	Before proving Lemma~\ref{lemma:bounds_pst2}, we establish the following.
	
	\begin{lemma} \label{lemma:phi_pst2}
		For any $ K > 0 $ and $ \delta \in [0,1] $, there exists $ C > 0 $ such that, for all $ v \in \mathcal{V}_{\Beta{zeta},K,\Epsilon{zeta}} $, $ t \in (0,1] $ and $ \phi \in L^{2,\alpha} \cap L^{\infty,e} $,
		\begin{equation} \label{bound_phi_pst2}
			\langle |\phi|, p_{0,t}(v,\cdot,y)^2 \rangle_\alpha \leq \frac{C}{t^{\frac{2+\delta}{4}}} \| \phi \|_{2,\alpha}^{\delta} \enorm{\phi}^{1-\delta} e^{-(1-\delta)|y| - \alpha(1+\frac{\delta}{2}) y}.
		\end{equation}
		In addition, for any $ K > 0 $ and $ \delta_1 \in (0, \frac{1}{2} \wedge (1-\frac{\alpha}{2})) $, there exist $ C > 0 $, $ \genericc > 0 $ and $ \delta \in (0,1) $ such that, for all $ \zeta \in \R $ with $ | \zeta | \leq \Cst{sup_zeta}(K) $, all $ y \in \R $ and $ t \geq 1 $,
		\begin{equation} \label{bound_phi_p-vphi}
			\langle |\phi|, | p_{0,t}(m_\zeta,\cdot,y) - \varphi_{\zeta}(y) \varphi_{\zeta} |^2 \rangle_\alpha \leq C e^{-\genericc t} \| \phi \|_{2,\alpha}^{\delta} \enorm{\phi}^{1-\delta} e^{-2\delta_1 |y| - \alpha y}.
		\end{equation}
	\end{lemma}

	\begin{proof}
		For $ t \in (0,1] $, we use Lemma~\ref{lemma:bound_pst} to write
		\begin{align*}
			\langle |\phi|, p_{0,t}(v,\cdot,y)^2 \rangle_\alpha &\leq C_1^2 e^{-\alpha y} \int_\R | \phi(x) | G_t(x-y)^2 dx \\
			&\leq \frac{C_1^2}{\sqrt{4\pi t}} e^{-\alpha y} \int_\R |\phi(x)| G_t(x-y) dx.
		\end{align*}
		We then use \eqref{bound_Gt_phi} from Lemma~\ref{lemma:bounds_integral_phi} to obtain \eqref{bound_phi_pst2}.
		To prove \eqref{bound_phi_p-vphi}, we apply \eqref{bound_pst_varphi} to write, for $ \delta_1 $, $ \delta_2 $ satisfying \eqref{condition_delta_1_2},
		\begin{equation*}
			\langle |\phi|, | p_{0,t}(m_\zeta,\cdot,y) - \varphi_{\zeta}(y) \varphi_{\zeta} |^2 \rangle_\alpha \leq C e^{-\genericc t} e^{-2 \delta_1 |y| - \alpha y} \int_\R |\phi(x)| e^{2\delta_2 |x|} dx.
		\end{equation*}
		The bound \eqref{bound_phi_p-vphi} then follows by \eqref{bound_phi_delta0} with $ q = 0 $, choosing $ \delta_2 < 1/2 $ and $ \delta $ such that
		\begin{equation*}
			\delta < \frac{1 - 2 \delta_2}{1 + \frac{\alpha}{2}}.
		\end{equation*}
	\end{proof}
	
	Let us now prove Lemma~\ref{lemma:bounds_pst2}.
	
	\begin{proof}[Proof of Lemma~\ref{lemma:bounds_pst2}]
		We start by proving \eqref{bound_Xi0} in the case $ s \leq 1 $ and $ t \leq 1 $.
		By the Cauchy-Schwarz inequality,
		\begin{equation*}
			| \Xi_0(s,t,\zeta,y,z) | \leq \langle | \partial_{x} m_\zeta f''(m_\zeta)|, p_{0,s}(m_\zeta,\cdot,y)^2 \rangle_\alpha^{1/2} \langle | \partial_{x} m_\zeta f''(m_\zeta) |, p_{0,t}(m_\zeta,\cdot,z)^2 \rangle_\alpha^{1/2}.
		\end{equation*}
		Then, applying \eqref{bound_phi_pst2} with $ \delta = 0 $ to each factor, we obtain that there exists $ C > 0 $ such that, for $ s \leq 1 $ and $ t \leq 1 $,
		\begin{equation*}
			| \Xi_0(s,t,\zeta,y,z) | \leq \frac{C}{(st)^{1/4}} e^{-\frac{1}{2}(|y| + |z|) - \frac{\alpha}{2}(y + z)}.
		\end{equation*}
		If now $ s \leq 1 $ and $ t \geq 1 $,
		\begin{multline} \label{bound_Xi0_small_s}
			| \Xi_0(s,t,\zeta,y,z) | \leq | \langle \partial_{x} m_\zeta f''(m_\zeta), p_{0,s}(m_\zeta,\cdot,y) (p_{0,t}(m_\zeta,\cdot,z) - \varphi_{\zeta}(z) \varphi_{\zeta}) \rangle_\alpha | \\ + | \varphi_{\zeta}(z) \langle \partial_{x} m_\zeta f''(m_\zeta), p_{0,s}(m_\zeta, \cdot,y) \varphi_{\zeta} \rangle_\alpha |.
		\end{multline}
		By the Cauchy-Schwarz inequality, the first term on the right-hand side is bounded from above by
		\begin{equation*}
			\langle | \partial_{x} m_\zeta f''(m_\zeta)|, p_{0,s}(m_\zeta,\cdot,y)^2 \rangle_\alpha^{1/2} \langle | \partial_{x} m_\zeta f''(m_\zeta) |, | p_{0,t}(m_\zeta,\cdot,z) - \varphi_{\zeta}(z) \varphi_{\zeta}|^2 \rangle_\alpha^{1/2}.
		\end{equation*}
		Applying \eqref{bound_phi_pst2} with $ \delta = 0 $ to the first factor and \eqref{bound_phi_p-vphi} to the second, we obtain that, for any $ \delta \in (0,\frac{1}{2} \wedge (1-\frac{\alpha}{2})) $, there exist $ C > 0 $ and $ \genericc > 0 $ such that this is bounded from above by
		\begin{equation*}
			\frac{C}{s^{1/4}} e^{-\genericc t} e^{-\frac{1}{2} |y| - \delta |z| - \frac{\alpha}{2}(y+z)}.
		\end{equation*}
		For the second term on the right-hand side of \eqref{bound_Xi0_small_s}, we apply the Cauchy-Schwarz inequality again to write
		\begin{multline*}
			| \varphi_{\zeta}(z) \langle \partial_{x} m_\zeta f''(m_\zeta), p_{0,s}(m_\zeta, \cdot,y) \varphi_{\zeta} \rangle_\alpha | \\ \leq | \varphi_{\zeta}(z) | \langle | \partial_{x} m_\zeta f''(m_\zeta) |, p_{0,s}(m_\zeta,\cdot,y)^2 \rangle_\alpha^{1/2} \langle | \partial_{x} m_\zeta f''(m_\zeta) |, \varphi_{\zeta}^2 \rangle_\alpha^{1/2}.
		\end{multline*}
		We then apply \eqref{bound_phi_pst2} with $ \delta = 0 $ along with the fact that $ | \varphi_{\zeta}(z) | \leq C e^{-|z|} $ for $ | \zeta | \leq \Cst{sup_zeta}(K) $ to obtain that the right-hand side is bounded from above by
		\begin{equation*}
			\frac{C}{s^{1/4}} e^{-|z| - \frac{1}{2} |y| + \frac{\alpha}{2}y}.
		\end{equation*}
		We then note that, for any $ \delta \leq 1-\frac{\alpha}{2} $, $ |z| \geq \delta |z| + \frac{\alpha}{2} z $ for all $ z \in \R $. This proves \eqref{bound_Xi0} for $ s \leq 1 $ and $ t \geq 1 $.
		The case $ s \geq 1 $ and $ t \leq 1 $ is similar.
		In the case $ s \geq 1 $ and $ t \geq 1 $, we start by noting that
		\begin{equation*}
			\langle \partial_{x} m_\zeta f''(m_\zeta), \varphi_{\zeta}^2 \rangle_\alpha = \frac{\langle (\partial_{x} m)^2 f''(m), \partial_{x} m \rangle_\alpha}{\| \partial_{x} m \|_{2,\alpha}^2}.
		\end{equation*}
		Recall that $ \mathcal{A} \partial_{x} m = 0 $.
		Differentiating this expression yields
		\begin{equation*}
			\mathcal{A} \partial_{xx} m = (\partial_{x} m)^2 f''(m),
		\end{equation*}
		from which it follows that
		\begin{equation*}
			\langle (\partial_x m)^2 f''(m), \partial_{x} m \rangle_\alpha = \langle \partial_{xx}m, \mathcal{A} \partial_{x} m \rangle_\alpha = 0.
		\end{equation*}
		As a result, we can write
		\begin{multline} \label{decomp_Xi0}
			| \Xi_0(s,t,\zeta,y,z) | \leq |\varphi_{\zeta}(z)| | \langle \partial_{x} m_\zeta f''(m_\zeta), \varphi_{\zeta} (p_{0,s}(m_\zeta,\cdot,y) - \varphi_{\zeta}(y) \varphi_{\zeta}) \rangle_\alpha | \\ + |\varphi_{\zeta}(y)| | \langle \partial_{x} m_\zeta f''(m_\zeta), \varphi_{\zeta} (p_{0,t}(m_\zeta,\cdot,z) - \varphi_{\zeta}(z) \varphi_{\zeta}) \rangle_\alpha \\ + | \langle \partial_{x} m_\zeta f''(m_\zeta), (p_{0,s}(m_\zeta,\cdot,y) - \varphi_{\zeta}(y) \varphi_{\zeta})(p_{0,t}(v,\cdot,z) - \varphi_{\zeta}(z) \varphi_{\zeta}) \rangle_\alpha |.
		\end{multline}
		We then bound each term on the right separately.
		By the Cauchy-Schwarz inequality, the first term is bounded from above by
		\begin{equation*}
			| \varphi_{\zeta}(z) | \langle |\partial_{x}m_\zeta f''(m_\zeta) |, \varphi_{\zeta}^2 \rangle_\alpha^{1/2} \langle | \partial_{x}m_\zeta f''(m_\zeta)| , | p_{0,s}(m_\zeta,\cdot,y) - \varphi_{\zeta}(y) \varphi_{\zeta} |^2 \rangle_\alpha^{1/2}.
		\end{equation*}
		By \eqref{bound_phi_p-vphi}, for any $ \delta \in (0, \frac{1}{2} \wedge (1-\frac{\alpha}{2})) $ this is bounded from above by
		\begin{equation*}
			C e^{-\genericc s} e^{-|z| - \delta |y| - \frac{\alpha}{2} y},
		\end{equation*}
		for some $ C > 0 $ and $ \genericc > 0$, for all $ s \geq 1 $.
		The second term on the right-hand side of \eqref{decomp_Xi0} is bounded in the same way.
		By the Cauchy-Schwarz inequality, the third term is bounded from above by
		\begin{equation*}
			\langle | \partial_{x}m_\zeta f''(m_\zeta)| , | p_{0,s}(m_\zeta,\cdot,y) - \varphi_{\zeta}(y) \varphi_{\zeta} |^2 \rangle_\alpha^{1/2} \langle | \partial_{x}m_\zeta f''(m_\zeta)| , | p_{0,t}(m_\zeta,\cdot,z) - \varphi_{\zeta}(z) \varphi_{\zeta} |^2 \rangle_\alpha^{1/2}.
		\end{equation*}
		By \eqref{bound_phi_p-vphi}, for any $ \delta \in (0,\frac{1}{2} \wedge (1-\frac{\alpha}{2})) $, this is bounded from above by
		\begin{equation*}
			C e^{-\genericc (t+s)} e^{- \delta (|y| + |z|) - \frac{\alpha}{2} (y + z)},
		\end{equation*}
		for some $ C > 0 $ and $ \genericc > 0 $, for all $ t \geq 1 $ and $ s \geq 1 $.
		As a result, using the fact that $ |z| \geq \delta |z| + \frac{\alpha}{2} z $ for any $ \delta \leq 1 - \frac{\alpha}{2} $, we obtain that, for any $ \delta \in (0,\frac{1}{2} \wedge (1-\frac{\alpha}{2})) $ there exist $ C > 0 $ and $ \genericc > 0 $ such that, for $ t \geq 1 $ and $ s \geq 1 $,
		\begin{equation*}
			| \Xi_0(s,t,\zeta,y,z) | \leq C (e^{-\genericc s} + e^{-\genericc t}) e^{-\delta (|y| + |z|) + \frac{\alpha}{2}(y + z)}.
		\end{equation*}
		This concludes the proof of \eqref{bound_Xi0}.
		
		We now turn to the proof of \eqref{Lpq_p2}.
		By \eqref{bound_phi_pst2}, for any $ \delta \in [0,1] $, there exists $ C > 0 $ such that, for $ t \in (0,1] $,
		\begin{equation*}
			\left( \int_\R | \langle \phi, p_{0,t}(v,\cdot,y)^2 \rangle_\alpha |^p e^{qy} dy \right)^{1/p} \leq \frac{C}{t^{\frac{2+\delta}{4}}} \| \phi \|_{2,\alpha}^{\delta} \enorm{\phi}^{1-\delta} \left( \int_\R e^{-p(1-\delta)|y| + (q-\alpha p(1+\frac{\delta}{2})) y} dy \right)^{1/p}.
		\end{equation*}
		We then choose $ \delta \in (0,1) $ such that
		\begin{equation*}
			1-\delta > \abs{\frac{q}{p} - \alpha\left(1+\frac{\delta}{2}\right)}.
		\end{equation*}
		This is equivalent to
		\begin{equation*}
			\delta < \left( \frac{1-\alpha + \frac{q}{p}}{1+\frac{\alpha}{2}} \right) \wedge \left( \frac{1 + \alpha - \frac{q}{p}}{1 - \frac{\alpha}{2}} \right),
		\end{equation*}
		which we can always satisfy thanks to \eqref{condition_p_q}.
		When $ t \geq 1 $, by \eqref{bound_pstv_pstm}, we have
		\begin{equation*}
			| \langle \phi, p_{0,t}(v,\cdot,y)^2 \rangle_\alpha | \leq C \langle |\phi|, p_{0,t}(m_{\zeta(v)}, \cdot,y)^2 \rangle_\alpha.
		\end{equation*}
		In addition, using the inequality $ a^2 \leq 2b^2 + 2(a-b)^2 $,
		\begin{equation*}
			\langle |\phi|, p_{0,t}(m_{\zeta}, \cdot,y)^2 \rangle_\alpha \leq 2|\varphi_{\zeta}(y)|^2 \langle |\phi|, \varphi_{\zeta}^2 \rangle_\alpha + 2 \langle |\phi|, |p_{0,t}(m_\zeta,\cdot,y) - \varphi_{\zeta}(y) \varphi_{\zeta}|^2 \rangle_\alpha.
		\end{equation*}
		We then note that, for $ |\zeta| \leq \Cst{sup_zeta}(K) $,
		\begin{equation*}
			\int_\R |\varphi_{\zeta}(y)|^{2p} e^{qy} dy \leq C,
		\end{equation*}
		for some $ C > 0 $, since
		\begin{equation*}
			-2 < \alpha-1 < \frac{q}{p} < 2.
		\end{equation*}
		On the other hand, by \eqref{bound_pst_varphi}, for any $ \delta_1 $ and $ \delta_2 $ satisfying \eqref{condition_delta_1_2}, there exist $ C > 0 $ and $ \genericc > 0 $ such that
		\begin{equation*}
			\langle |\phi|, |p_{0,t}(m_\zeta,\cdot,y) - \varphi_{\zeta}(y) \varphi_{\zeta}|^2 \rangle_\alpha \leq C e^{-\genericc t} e^{-2 \delta_1 |y| - \alpha y} \int_\R | \phi(x) | e^{2\delta_2 |x|} dx.
		\end{equation*}
		Then, by \eqref{bound_phi_delta0}, if $ \delta_2 < 1/2 $, there exists $ \delta \in (0,1) $ and $ C > 0 $ such that
		\begin{equation*}
			\int_\R | \phi(x) | e^{2\delta_2 |x|} dx \leq C \| \phi \|_{2,\alpha}^{\delta} \enorm{\phi}^{1-\delta}.
		\end{equation*}
		As a result, for any $ 0 < \delta_1 < \delta_2 < \frac{1}{2} \wedge (1-\frac{\alpha}{2}) $, there exist $ C > 0 $, $ \genericc > 0 $ and $ \delta \in (0,1) $ such that
		\begin{multline*}
			\left( \int_\R \langle |\phi|, |p_{0,t}(m_\zeta,\cdot,y) - \varphi_{\zeta}(y) \varphi_{\zeta}|^2 \rangle_\alpha^p e^{qy} dy \right)^{1/p} \\  \leq C e^{-\genericc t} \| \phi \|_{2,\alpha}^{\delta} \enorm{\phi}^{1-\delta} \left( \int_\R e^{-2p \delta_1 |y| + (q-\alpha p) y} dy \right)^{1/p}.
		\end{multline*}
		It remains to note that, thanks to \eqref{condition_p_q}, we can choose $ \delta_1 $ and $ \delta_2 $ such that
		\begin{equation*}
			2 \delta_1 > \abs{\frac{q}{p} - \alpha}.
		\end{equation*}
		This concludes the proof of \eqref{Lpq_p2}.
	\end{proof}

	\subsection{The linearised semigroup and a diffusion process} \label{subsec:diffusion}
	
	We now set out to prove Lemma~\ref{lemma:long_time_pstm}.
	The proof will exploit a connection between the semigroup $ e^{-t \mathcal{A}} $ and the Markov semigroup of a real-valued diffusion process.
	Let $ (X_t, t \geq 0) $ be a diffusion process with generator
	\begin{equation} \label{def_generator_diffusion}
		\mathcal{L} g(x) = \partial_{xx}g(x) + \left( \alpha + 2 \frac{\partial_{x}\varphi(x)}{\varphi(x)} \right) \partial_{x}g(x),
	\end{equation}
	and let $ \E[x]{\cdot} $ denote expectation with respect to the law of $ (X_t, t \geq 0) $ with initial condition $ X_0 = x $.
	Recall that
	\begin{equation*}
		c_t(v,x) = f'(\Phi_t(v,x)) - f'(m_{\zeta(v)}(x)).
	\end{equation*}
	
	\begin{proposition} \label{prop:duality_lineage}
		For any $ \phi \in L^{2,\alpha} $, $ x \in \R $ and $ t \geq 0 $, 
		\begin{equation*}
			e^{-t\mathcal{A}}\phi(x) = \varphi(x)\, \E[x]{\frac{\phi(X_t)}{\varphi(X_t)}}.
		\end{equation*}
		Moreover, for any $ v \in \mathcal{V}_{\Beta{zeta},K,\Epsilon{zeta}} $ with $ \zeta(v) = 0 $ and for any $ 0 \leq s \leq t $,
		\begin{equation*}
			T_{s,t,v}^*\phi(x) = \varphi(x)\, \E[x]{\frac{\phi(X_{t-s})}{ \varphi(X_{t-s})} \exp \left(\int_{0}^{t-s} c_{t-r}(v,X_r) dr \right) }.
		\end{equation*}
	\end{proposition}
	
	\begin{proof}
		Let $ (P_t, t \geq 0) $ denote the transition semigroup of $ (X_t, t \geq 0) $.
		The first claim can then be written
		\begin{equation} \label{etA0} 
			e^{-t \mathcal{A}} \phi = \varphi P_t \left( \frac{\phi}{\varphi} \right).
		\end{equation}
		This identity clearly holds for $ t = 0 $; 
		for general $t \geq 0$, as a consequence of \cite[Proposition~I.2.9]{ethier_markov_1986}  it suffices to prove that 
		\begin{equation} \label{link_generator}
			\varphi \mathcal{L} \left( \frac{\phi}{\varphi} \right) = - \mathcal{A} \phi.
		\end{equation}
		For the second identity, note that \eqref{def_pst} is equivalent to
		\begin{equation*}
			\frac{d}{dt}T_{s,t,v}^*\phi = (\partial_{xx} + \alpha \partial_{x} + f'(\Phi_t(v,\cdot)))T_{s,t,v}^*\phi.
		\end{equation*}
		As a result, setting $ g_t := \frac{T_{s,t,v}^*\phi}{\varphi} $, we have, for $ t \geq s $,
		\begin{align*}
			\partial_t g_t &= - \frac{\mathcal{A} T_{s,t,v}^* \phi}{\varphi} + (f'(\Phi_{t}(v,\cdot))-f'(m)) g_t \\
			&= \mathcal{L} g_t + c_t(v,\cdot) g_t.
		\end{align*}
		By Ito's lemma, it follows that
		\begin{align*}
			Y_u = \exp \left(\int_{0}^{u} c_{t-r}(v,X_r) dr \right) g_{t-u}(X_u), \quad 0 \leq u \leq t-s,
		\end{align*}
		is a martingale with respect to the natural filtration of the process $ (X_u, u \geq 0) $.
		As a result,
		\begin{equation*}
			g_t(x) = \E[x]{ \exp\left( \int_{0}^{t-s} c_{t-r}(v,X_r) dr \right) \frac{\phi(X_{t-s})}{\varphi(X_{t-s})} },
		\end{equation*}
		yielding the result.
	\end{proof}
	
	\begin{lemma} \label{lemma:translation_pst}
		For any $ v \in L^{2,\alpha} $ and any $ \eta \in \R $,
		\begin{equation*}
			p_{s,t}(v_{\eta}, x, y) = p_{s,t}(v, x-\eta, y-\eta) e^{-\alpha \eta}.
		\end{equation*}
	\end{lemma}
	
	\begin{proof}
		First note that, by definition,
		\begin{equation*}
			\Phi_t(v_\eta,x) = \Phi_t(v,x-\eta).
		\end{equation*}
		As a result,
		\begin{equation*}
			T_{s,t,v_\eta}^* (\phi_\eta) = \left( T_{s,t,v}^* \phi \right)_\eta.
		\end{equation*}
		The result then follows by evaluating this expression on $ \phi_{-\eta} $ and a change of variables.
	\end{proof}
	
	\begin{cor} \label{cor:bound_pv_pm}
		For any $ K > 0 $, there exists a constant $ C > 0 $ such that, for all for all $ v \in \mathcal{V}_{\Beta{zeta}, K, \Epsilon{zeta}} $, for all $ 0 \leq s \leq t $ and $ x, y \in \R $,
		\begin{align*}
			p_{s,t}(v,x,y) \leq C\, p_{s,t}(m_{\zeta(v)}, x, y).
		\end{align*}
	\end{cor}
	
	\begin{proof}[Proof of Corollary~\ref{cor:bound_pv_pm}]
		Assume first that $ \zeta(v) = 0 $.
		Combining the two statements in Proposition~\ref{prop:duality_lineage}, we can write
		\begin{equation*}
			p_{s,t}(v,x,y) = p_{s,t}(m, x, y)  \E[x]{\exp \left(\int_{0}^{t-s} c_{t-r}(v,X_r) dr \right)}{X_{t-s} = y}.
		\end{equation*}
		Using Corollary~\ref{cor:uniform_cvg}, we see that, for $ v \in \mathcal{V}_{\Beta{zeta},K,\Epsilon{zeta}} $, almost surely,
		\begin{align*}
			\abs{c_{t-r}(v,X_r)} \leq C(\|s(v)\|_{\infty,\lambda} + \sqrt{\dist(v,M)}) e^{-\genericc (t-r)}.
		\end{align*}
		Hence, there exists a constant $ C > 0 $ which bounds the above expectation for all $ v \in \mathcal{V}_{\Beta{zeta}, K, \Epsilon{zeta}} $, for all $ 0 \leq s \leq t $, and the result follows.
		When $ \zeta(v) \neq 0 $, we use the fact that $ \zeta(v_{-\zeta(v)}) = 0 $ and Lemma~\ref{lemma:translation_pst} to conclude.
	\end{proof}
	
	This proves the first part of Lemma~\ref{lemma:long_time_pstm}.
	Let us define a probability density on $ \R $ as
	\begin{equation} \label{def:pi}
		\pi(x) := \varphi(x)^2 e^{\alpha x}.
	\end{equation}
	Also let $ q_t(\cdot, \cdot) $ denote the transition density of the diffusion process $ (X_t, t \geq 0) $.
	
	\begin{lemma} \label{lemma:reversibility}
		The probability measure with density $ \pi $ on $ \R $ is a reversible stationary distribution for $ (X_t, t \geq 0) $.
		In particular, for all $ x, y \in \R $ and $ t \geq 0 $,
		\begin{equation*}
			\pi(x) q_t(x, y) = \pi(y) q_t(y, x).
		\end{equation*}
	\end{lemma}
	
	\begin{proof}
		This follows from a straightforward computation.
		Indeed, we can write
		\begin{equation*}
			\int_{\R} \psi(x) \mathcal{L} \phi(x) \pi(x) dx = \langle \varphi \psi, \varphi \mathcal{L}\phi \rangle_\alpha.
		\end{equation*}
		Using \eqref{link_generator}, this is
		\begin{equation*}
			- \langle \varphi \psi, \mathcal{A} (\varphi \phi) \rangle_\alpha.
		\end{equation*}
		Since $ \mathcal{A} $ is self-adjoint in $ L^{2,\alpha} $, we obtain that
		\begin{equation*}
			\int_{\R} \psi(x) \mathcal{L} \phi(x) \pi(x) dx = \int_{\R} \mathcal{L} \psi(x) \phi(x) \pi(x) dx,
		\end{equation*}
		which yields the result.
	\end{proof}
	
	The following is proved in the next subsection.
	
	\begin{lemma} \label{lem:qtxy_estimate}
		For any $\varepsilon_1\in (0,2+\alpha)$, $\varepsilon_2\in (0,2-\alpha)$ and $\delta>0$, there exist $\newsmalla{qtpidecay}>0$ and $\newbigC{qtpidecay}>0$ such that the following holds.
		For $x>0,$ $t\geq 1$ and $y \in \mathbb{R}$,
		\begin{equation} \label{eq:qtxy_estimate1}
		 	\left|q_{t}(y, x)-\pi(x)\right| \leq \bigC{qtpidecay} e^{-\smalla{qtpidecay} t} e^{\delta|y|+\frac 12\left(2-\alpha-\varepsilon_{2}\right)(y \vee 0)} e^{-\frac 12(1-2 \delta)\left(2-\alpha-\varepsilon_{2}\right) x},
		\end{equation}
		for $x > 0, $ $t \geq 1$ and $y \in \mathbb{R}$,
		\begin{equation} \label{eq:qtxy_estimate2}
		 	\left|q_{t}(y,-x)-\pi(-x) \right| \leq \bigC{qtpidecay} e^{-\smalla{qtpidecay} t} e^{\delta|y|+\frac 12(2+\alpha-\varepsilon_1)((-y)\vee 0)} e^{-\frac 12 (1-2 \delta)\left(2+\alpha-\varepsilon_1\right)x}.
		\end{equation}
	\end{lemma}
	
	Let us show how this implies \eqref{bound_pst_varphi}, which will conclude the proof of Lemma~\ref{lemma:long_time_pstm}.
	
	\begin{proof}[Conclusion of the proof of Lemma~\ref{lemma:long_time_pstm}]
		We first note that the statement for $ | \zeta | \leq \Cst{sup_zeta} $ follows from the statement for $ \zeta = 0 $ and Lemma~\ref{lemma:translation_pst}.
		We thus take $ \zeta = 0 $ in the remainder of the proof.
		Proposition~\ref{prop:duality_lineage} can then be written
		\begin{equation*}
			p_{0,t}(m,x,y) = \frac{\varphi(y)}{\varphi(x)} q_t(y,x) e^{-\alpha x}.
		\end{equation*}
		We further note that
		\begin{equation*}
			\frac{\varphi(y)}{\varphi(x)} \pi(x) e^{-\alpha x} = \varphi(x) \varphi(y).
		\end{equation*}
		As a result,
		\begin{equation} \label{pt_varphi_qt_pi}
			\abs{ p_{0,t}(m,x,y) - \varphi(x) \varphi(y) } \leq \frac{| \varphi(y) |}{| \varphi(x) |} \abs{ q_t(y,x) - \pi(x) } e^{-\alpha x}.
		\end{equation}
		Fix $ \delta_1 $ and $ \delta_2 $ satisfying \eqref{condition_delta_1_2} and take $ \varepsilon \in (2\delta_1,2\delta_2) $.
		By Lemma~\ref{lem:qtxy_estimate}, for any $ \delta > 0 $, there exist $ \smalla{qtpidecay} > 0 $ and $ \bigC{qtpidecay} $ such that, for $ x, y \in \R $,
		\begin{multline*}
			| q_t(y,x) - \pi(x) | \leq \bigC{qtpidecay} e^{-\smalla{qtpidecay} t} e^{\delta |y| + \frac{1}{2}(2-\alpha-\varepsilon) (y \vee 0) + \frac{1}{2} (2+\alpha - \varepsilon)((-y) \vee 0)} \\ \times  e^{-\frac{1}{2}(1-2\delta)(2-\alpha-\varepsilon) (x \vee 0) - \frac{1}{2} (1-2\delta) (2+\alpha-\varepsilon) ((-x)\vee 0)}.
		\end{multline*}
		We then note that, for $ z \in \R $,
		\begin{equation*}
			(2-\alpha-\varepsilon) (z \vee 0) + (2+\alpha - \varepsilon)((-z) \vee 0) = (2-\varepsilon) |z| - \alpha z.
		\end{equation*}
		Returning to \eqref{pt_varphi_qt_pi}, we obtain
		\begin{align*}
			\abs{ p_{0,t}(m,x,y) - \varphi(x) \varphi(y) } &\leq \bigC{qtpidecay} e^{-\smalla{qtpidecay} t} e^{\left( \frac{1}{2}(2-\varepsilon) + \delta - 1 \right) |y| - \frac{\alpha}{2} y} e^{\left( 1 - \frac{1}{2}(1-2\delta)(2-\varepsilon) \right) |x| + ((1-2\delta)\frac{\alpha}{2}-\alpha) x} \\
			&\qquad = \bigC{qtpidecay} e^{-\smalla{qtpidecay} t} e^{-\left(\frac{\varepsilon}{2} - \delta \right) |y| - \frac{\alpha}{2} y} e^{\left( \frac{\varepsilon}{2} + (2-\varepsilon) \delta \right) |x| - \frac{\alpha}{2} (1-2\delta) x} \\
			&\leq \bigC{qtpidecay} e^{-\smalla{qtpidecay} t} e^{-\left(\frac{\varepsilon}{2} - \delta \right) |y| - \frac{\alpha}{2} y} e^{\left( \frac{\varepsilon}{2} + (2+\alpha-\varepsilon) \delta \right) |x| - \frac{\alpha}{2} x}.
		\end{align*}
		To conclude, we note that, since $ \varepsilon \in (2\delta_1, 2\delta_2) $, one can choose $ \delta $ small enough that
		\begin{align*}
			\frac{\varepsilon}{2} - \delta \geq \delta_1, && \text{ and } && \frac{\varepsilon}{2} + (2+\alpha-\varepsilon) \delta \leq \delta_2,
		\end{align*}
		yielding \eqref{bound_pst_varphi}.
	\end{proof}
	
	It thus remains to prove Lemma~\ref{lem:qtxy_estimate}, which turns out to be quite technical.
	This is the purpose of the next subsection.
	
	\subsection{Convergence to the stationary distribution} \label{subsec:cvg_qt}
	
	Note that $ (X_t, t \geq 0) $ solves the SDE 
	\begin{equation} \label{eq:SDEforqt}
		d X_{t}=\left(\alpha+2 \frac{\partial _{x x} m \left(X_{t}\right)}{\partial_{x}m\left(X_{t}\right)}\right) d t + \sqrt{2} d B_{t},
	\end{equation}
	where $\left(B_{t}\right)_{t \ge 0}$ is a Brownian motion.
	
	\begin{lemma} \label{lem:qtxy_overall}
		There exist $\newsmalla{qtxy}>0$ and $\newbigC{qtxy}>0$ such that for any $x,y\in \R$ and $t>0$,
		\begin{equation} \label{eq:qtxyoverall}
			q_{t}(x, y) \le \bigC{qtxy} t^{-1 / 2} \exp \left(-\frac{\smalla{qtxy}}{t}(y-x)^{2}\right).
		\end{equation}
	\end{lemma}
	
	\begin{proof}
		Since the drift term in~\eqref{eq:SDEforqt} is Lipschitz continuous and bounded, this follows directly from
		(1.2) in~\cite{sheu_some_1991}.
	\end{proof}
	
	We begin with a few preliminary estimates.
	The first estimate concerns the drift term in the SDE~\eqref{eq:SDEforqt}.
	\begin{lemma} \label{lem:SDEdrift}
	The following limits hold:
	\[
	\lim_{x \rightarrow \infty}\left(\alpha+2 \frac{\partial_{x x}m(x)}{\partial_{x}m(x)} \right)= \alpha-2
	\quad \text{and} \quad
	\lim_{x \rightarrow -\infty}\left(\alpha+2 \frac{\partial_{x x}m(x)}{\partial_{x}m(x)} \right)= \alpha+2.
	\]
	\end{lemma}
	\begin{proof}
		Since, as $ |x| \to \infty $,
		\begin{align*}
			\partial_{x x} m(x) \sim sign(x) e^{-|x|} && \partial_x m(x) \sim - e^{-|x|},
		\end{align*}
		the result follows.
	\end{proof}
	
	The second estimate concerns reflected Brownian motion on $[0,\infty)$ with drift towards $0$.
	\begin{lemma} \label{lem:reflectedBM}
	For $z\in \R$, let $\Phi(z)=\mathbb P(Z\le z)$, where $Z$ is a standard Gaussian random variable.
	For
	$\mu>0$ and $\sigma^{2}>0$,
	let $\left(Y^{\mu,\sigma}_{t}\right)_{t \ge 0}$ denote
	Brownian motion on $[0, \infty)$, reflected at $0$, with variance $\sigma^{2}$ and drift $-\mu$.
	Then for $y,t>0$ and $x\ge 0$,
	\[
	\mathbb{P}_{y}\left(Y^{\mu,\sigma}_{t} \le x\right)=\Phi\left(\frac{x-y+\mu t}{\sqrt{\sigma^{2} t}}\right)-e^{-2 \mu x/\sigma^{2}} \Phi\left(\frac{-x-y+\mu t}{\sqrt{\sigma^{2} t}}\right).
	\]
	\end{lemma}
	\begin{proof}
	This follows from~\cite{linetsky_transition_2005} (see p.447, after (30), in that work).	\end{proof}

	Let $(X_{t}^{(1)})_{t \ge 0},(X_{t}^{(2)})_{t \ge 0}$ be independent diffusions with generator $\mathcal L$ as defined in~\eqref{def_generator_diffusion}. 
	For $t\ge 0$, write $\mathcal{F}_{t}^{(j)}:=\sigma((X_{s}^{(j)})_{s \leq t})$ for $j \in\{1,2\}$,
	and $\mathcal F_{t}=\sigma((X_{s}^{(1)})_{s \leq t},(X_{s}^{(2)})_{s \leq t})$.
	Write $\mathbb{P}_{y, \pi}$ for the probability measure under which $X_{0}^{(1)}=y$ almost surely and $X_{0}^{(2)}$ has density $\pi$.
	Let
	\begin{equation} \label{eq:Thitdefn}
	T_{\text{hit}}=\inf \left\{t \geq 0: X_{t}^{(1)}=X_{t}^{(2)}\right\}.
	\end{equation}
	Let
	\begin{equation} \label{eq:tildeXdefn}
	\tilde{X}_{t}^{(2)}= \begin{cases}X_{t}^{(2)} & \text { for } t \leq T_{\text{hit}}, \\ X_{t}^{(1)} & \text { for } t>T_{\text{hit}},\end{cases}
	\quad \text{and} \quad 
	\tilde{X}_{t}^{(1)}=X_{t}^{(1)}
	\text{ for }t\ge 0.
	\end{equation}
	Then $(\tilde{X}_{t}^{(2)})_{t \ge 0}$ is a diffusion with generator $\mathcal{L}$ and $\tilde{X}_{0}^{(2)}=X_{0}^{(2)}$.
	
	\begin{lemma} \label{lem:probXtI}
	For $I\subseteq \R$ measurable, $t>0$ and $y\in \R$,
	\[
	\mathbb{P}_{y}\left(X_{t} \in I\right)-\int_{I} \pi(x) d x
	=\mathbb{P}_{y, \pi}\left(X_{t}^{(1)} \in I, T_{\text{hit}}>t\right)-\mathbb{P}_{y, \pi}\left(X_{t}^{(2)} \in I, T_{\text{hit}}>t\right) .
	\]
	\end{lemma}
	\begin{proof}
	By partitioning on the events $\{T_{\text{hit}}\le t\}$ and $ \lbrace T_{\text{hit}} > t \rbrace $, and then since $X^{(1)}_t=\tilde X^{(2)}_t$ for $t\ge T_{\text{hit}}$, we can write
	\begin{align*}
	\mathbb{P}_{y}\left(X_{t} \in I\right) &=\mathbb{P}_{y,\pi}\left(X_{t}^{(1)} \in I, T_{\text{hit}} \leq t\right)+\mathbb{P}_{y,\pi}\left(X_{t}^{(1)} \in I, T_{\text{hit}}>t\right)\\
	& =\mathbb{P}_{y, \pi}\left(\tilde{X}_{t}^{(2)} \in I, T_{\text{hit}} \leq t\right)+\mathbb{P}_{y, \pi}\left(X_{t}^{(1)} \in I, T_{\text{hit}}>t\right)\\
	&= \mathbb{P}_{y,\pi}\left(\tilde{X}_{t}^{(2)} \in I\right)-\mathbb{P}_{y,\pi}\left(\tilde{X}_{t}^{(2)} \in I, T_{\text{hit}}>t\right)  +\mathbb{P}_{y, \pi}\left(X_{t}^{(1)} \in I, T_{\text{hit}}>t\right).
	\end{align*}
	Using for the first term on the right-hand side that $(\tilde{X}_{t}^{(2)})_{t \ge 0}$ is a diffusion with generator $\mathcal{L}$ and $\tilde{X}_{0}^{(2)}=X_{0}^{(2)}$, and using for the second term that $\tilde X^{(2)}_t=X^{(2)}_t$ for $t<T_{\text{hit}}$, we have 
	\begin{align*}
	\mathbb{P}_{y}\left(X_{t} \in I\right) &=\mathbb{P}_{y, \pi}\left(X_{t}^{(2)} \in I\right)-\mathbb{P}_{y, \pi}\left(X_{t}^{(2)} \in I, T_{\text{hit}}>t\right) 
	+\mathbb{P}_{y,\pi}\left(X_{t}^{(1)} \in I, T_{\text{hit}}>t\right).
	\end{align*}
	Since 
	\[
	\mathbb{P}_{y, \pi}\left(X_{t}^{(2)} \in I\right)=\int_{I} \pi(x) d x,
	\]
	the result follows.
	\end{proof}
	We will use Lemma~\ref{lem:probXtI} to establish Lemma~\ref{lem:qtxy_estimate} by showing that
	$\mathbb{P}_{y, \pi}(T_{\text{hit}}>t)$ decays exponentially in $t$, and $\mathbb{P}_{y}\left(|X_{t}| \ge x\right), \mathbb{P}_{\pi}\left(|X_{t}| \ge x\right)$ decay exponentially in $x$.
	
	For $\varepsilon\in (0, 2-\alpha)$, using Lemma~\ref{lem:SDEdrift}, take $\newbigC{driftlimit}(\varepsilon)>0$ sufficiently large that
	\begin{align*}
	\alpha+2 \frac{\partial_{x x} m(x)}{\partial_x m(x)}&<-(2-\alpha)+\varepsilon \quad \forall x \ge \bigC{driftlimit}(\varepsilon)\\
	\text { and } \quad \alpha+2 \frac{\partial_{x x}m(x)}{\partial_x m(x)}&>\alpha+2-\varepsilon \quad \forall x \leq-\bigC{driftlimit}(\varepsilon).
	\end{align*}
	
	\begin{lemma} \label{lem:Xttail_initial}
	For all $x\ge 0$, $y\in \R$, $t>0$ and $\varepsilon>0$,
	\[
	\mathbb{P}_{y}\left(\left|X_{t}\right| \ge x\right) \leq 2 e^{\frac 12 (2-\alpha-\varepsilon)|y|} e^{-\frac 12 (2-\alpha-\varepsilon)\left(x-2 \bigC{driftlimit}(\varepsilon)\right)}
	\]
	and, for $\alpha \in (0,2)$, there exists $\newbigC{pistatdist}>0$ such that
	\[
	\mathbb{P}_{\pi}\left(\left|X_{t}\right| \ge x\right)\leq \bigC{pistatdist} e^{-(2-\alpha) x}.
	\]
	\end{lemma}
	
	\begin{proof}
	Let $\left(Y_{t}\right)_{t \ge 0}$ denote a reflected Brownian motion on $[0, \infty)$ reflected at $0$ with drift $-(2-\alpha-\varepsilon)$ and variance 2.
	Then we claim that we can couple $\left(X_{t}\right)_{t \ge 0}$ and $\left(Y_{t}\right)_{t \ge 0}$ in such a way that
	\[
	\left|X_{0}\right|=Y_{0} \quad 
	\text{and} \quad  \left|X_{t}\right|\le Y_{t}+2 \bigC{driftlimit}(\varepsilon) \quad \forall t \ge 0.\]
	Indeed, when $\left|X_{t}\right|>\bigC{driftlimit}(\varepsilon),$ $\left|X_{t}\right|$ has stronger drift towards 0 than $Y_t$, so the processes can be coupled in such a way that $\left|X_{t}\right|-Y_{t}$ is non-increasing in $t$; after a time when $\left|X_{t}\right|$ hits $\left[0, \bigC{driftlimit}(\varepsilon)\right]$, letting $\tau$ denote the next time at which $\left|X_{t}\right|$ hits $2 \bigC{driftlimit}(\varepsilon)$, we have $\left|X_{\tau}\right|=2 \bigC{driftlimit}(\varepsilon) \leq Y_{\tau}+2 \bigC{driftlimit}(\varepsilon)$, and then we can couple the processes in such a way that from time $\tau$ onwards, $\left|X_{t}\right|-Y_{t}$ is non-increasing in $t$ until $\left|X_{t}\right|$ hits $\bigC{driftlimit}(\varepsilon)$ and so on.
	
	Therefore, for $y \in \mathbb{R}$ and $x>2 \bigC{driftlimit}(\varepsilon)$ with $|y| \leq x-2 \bigC{driftlimit}(\varepsilon)$,
	we have, letting $\mu=2-\alpha -\varepsilon$,
	\[
	\begin{aligned}
	\mathbb{P}_{y}\left(\left|X_{t}\right| \ge x\right) &\leq \mathbb{P}_{|y|}\left(Y_{t} \geq x-2 \bigC{driftlimit}(\varepsilon)\right) \\
	& =1-\mathbb{P}_{|y|}\left(Y_{t} \leq x-2 \bigC{driftlimit}(\varepsilon)\right) \\
	& =1-\Phi\left(\frac{x-2 \bigC{driftlimit}(\varepsilon)-|y|+\mu t}{\sqrt{2 t}}\right)+e^{-\mu\left(x-2 \bigC{driftlimit}(\varepsilon)\right)} \Phi\left(\frac{-x+2 \bigC{driftlimit}(\varepsilon)-|y|+\mu t}{\sqrt{2 t}}\right) \\
	& \leq e^{-\frac 12 \cdot \frac{1}{2 t}\left(x-2 \bigC{driftlimit}(\varepsilon)-|y|+\mu t\right)^{2}}+e^{-\mu\left(x-2 \bigC{driftlimit}(\varepsilon)\right)},
	\end{aligned}
	\]
	using Lemma~\ref{lem:reflectedBM} in the third line, and using that $1-\Phi(z) \leq e^{-\frac 12  z^{2}}$ for $z \ge 0$ and $\Phi(z)\le 1$ $\forall z\in \mathbb R$ in the last line.
	Therefore
	\begin{align} \label{eq:PyXtlargebound}
	\mathbb{P}_{y}\left(\left|X_{t}\right| \ge x\right) &\leq e^{-\frac{1}{4 t} \cdot 2 \mu t\left(x-2 \bigC{driftlimit}(\varepsilon)-|y|\right)}+e^{-\mu\left(x-2 \bigC{driftlimit}(\varepsilon)\right)} \notag \\
	& =e^{-\frac{1}{2} \mu\left(x-2 \bigC{driftlimit}(\varepsilon)-|y|\right)}+e^{-\mu\left(x-2 \bigC{driftlimit}(\varepsilon)\right)} \notag \\
	& \leq 2 e^{\frac 12 \mu|y|} e^{-\frac 12 \mu\left(x-2 \bigC{driftlimit}(\varepsilon)\right)} ,
	\end{align}
	where $\mu=2-\alpha-\varepsilon$, which establishes the first statement (if $|y| > x-2 \bigC{driftlimit}(\varepsilon)$ then the statement is trivially true).
	For the second statement, for $x\ge 0$ we have
	\[
	\mathbb{P}_{\pi}\left(\left|X_{t}\right| \ge x\right)=\int_{\{|y| \ge x\}} \pi(y) d y \leq\left(\int_{x}^{\infty} e^{(\alpha-2) y} d y+\int_{-\infty}^{-x} e^{(2+\alpha) y} d y\right)\cdot\left\|\partial_{x} m\right\|_{2, \alpha}^{-2},
	\]
	and the result follows.
	\end{proof}

	We now bound $\mathbb{P}_{y, \pi}(T_{\text{hit}}>t)$.
	Fix $C> 1$ large.
	Let $\sigma_{0}=0$ and $\tau_{0}=0$.
	We define $\left(\sigma_{i}\right)_{i \ge 1}$ and $\left(\tau_{i}\right)_{i \ge 1}$ iteratively as follows.
	For $i \ge 0$, let
	\begin{align*}
	\sigma_{i+1}&=\inf\left\{t \geq \tau_{i}:\left|X_{t}^{(1)}\right| \vee \left|X_{t}^{(2)}\right| \leq C\right\}\\
	\text{and} \quad \tau_{i+1}&=\inf\left\{t \ge \sigma_{i+1}:\left|X_{t}^{(1)}\right| \vee\left|X_{t}^{(2)}\right| \ge C^{2}\right\}.
	\end{align*}
	
	\begin{lemma} \label{lem:hitifclose}
	For $C$ sufficiently large, there exists $\newsmalla{hitifclose}>0$ (depending on $C$) such that the following holds.
	Let
	\[
	\tau=\inf\left\{t \ge 0:\left|X_{t}^{(1)}\right| \vee\left|X_{t}^{(2)}\right| \ge C^{2}\right\}.
	\]
	Then for $x_0^{(1)}, x_0^{(2)}\in \R$ with $\left|x_0^{(1)}\right| \vee \left|x_0^{(2)}\right| \leq C$,
	\begin{equation} \label{eq:claimhit}
	\mathbb P_{x^{(1)}_0,\, x^{(2)}_0}\left( T_{\text{hit}}\le \tau\right) \ge \smalla{hitifclose}.
	\end{equation}
	\end{lemma}
	
	\begin{proof}
	Suppose without loss of generality that $x_{0}^{(1)} \ge x_{0}^{(2)}$.
	We can write
	\begin{align*}
	X_{t}^{(1)}&=X_{0}^{(1)}+\alpha t+2 \int_{0}^{t} \frac{\partial_{x x} m(X_{s}^{(1)})}{\partial_{x}m(X_{s}^{(1)})} d s+\sqrt{2} B_{t}^{(1)}\\
	\text{and}\quad X_{t}^{(2)}&=X_{0}^{(2)}+\alpha t+2 \int_{0}^{t} \frac{\partial_{x x} m(X_{s}^{(2)})}{\partial_{x}m(X_{s}^{(2)})} d s+\sqrt{2} B_{t}^{(2)},
	\end{align*}
	where $B^{(1)}$ and $B^{(2)}$ are independent Brownian motions.
	By a union bound, we have
	\begin{equation} \label{eq:Thittauunion}
	\mathbb{P}_{x_{0}^{(1)}, x_{0}^{(2)}}(T_{\text{hit}}>\tau) \leq \mathbb{P}_{x_{0}^{(1)}, x_{0}^{(2)}}(T_{\text{hit}}>C)+\mathbb{P}_{x_{0}^{(1)}, x_{0}^{(2)}}(\tau<C).
	\end{equation}
	Let $K_m=\sup _{z \in \mathbb{R}}\left| \frac{\partial_{x x}m(z)}{\partial_{x}m(z)}\right|$, and let
	$\left(B_{t}\right)_{t \ge 0}$ denote a Brownian motion.
	Then for the first term on the right-hand side of~\eqref{eq:Thittauunion}, we have
	\begin{align} \label{eq:ThitCbound}
	\mathbb{P}_{x_{0}^{(1)}, \, x_{0}^{(2)}}(T_{\text{hit}}>C)
	&\le \mathbb{P}\left(\inf _{t \in[0, C]}\left(\sqrt{2}\left(B_{t}^{(1)}-B_{t}^{(2)}\right)+x_{0}^{(1)}-x_{0}^{(2)}+4 t \sup _{z \in \mathbb{R}}\left|\frac{\partial_{x x}m(z)}{\partial_{x}m(z)}\right|\right)>0\right)\notag\\
	&\le \mathbb{P}\left(\inf _{t \in[0, C]} 2 B_{t}>-2 C-4 C K_{m}\right)\notag\\
	&= \mathbb{P}\left(\sup _{t \in[0, C]} B_{t}<C\left(1+2 K_{m}\right)\right) \notag\\
	&= 1-2 \mathbb{P}\left(B_{C} \ge C\left(1+2 K_{m}\right)\right) \notag\\
	& \leq 1-2  \frac{C^{1 / 2}\left(1+2 K_{m}\right)}{1+C\left(1+2 K_{m}\right)^{2}} \cdot \frac{1}{\sqrt{2 \pi}} e^{-\frac{1}{2} C\left(1+2 K_{m}\right)^{2}},
	\end{align}
	where the penultimate line follows by the reflection principle, and the last line uses the Gaussian tail bound $1-\Phi(z)\ge \frac{z}{1+z^2}\frac 1 {\sqrt{2\pi}}e^{-\frac 12 z^2}$ for $z>0$.
	For the second term on the right-hand side of~\eqref{eq:Thittauunion}, by a union bound, and taking $C$ sufficiently large for the third inequality, we have
	\begin{align} \label{eq:tauCbound}
	\mathbb{P}_{x_{0}^{(1)}, x_{0}^{(2)}}(\tau<C) 
	& \leq \mathbb{P}\left(\sup _{t \in[0, C]}\left|x_{0}^{(1)}+\alpha t+2 t \sup _{z \in \mathbb{R}}\left| \frac{\partial_{x x}m(z)}{\partial_{x}m(z)}\right|+\sqrt{2} B_{t}^{(1)}\right| \ge C^{2}\right) \notag \\
	&\qquad +\mathbb{P}\left(\sup _{t \in[0, C]}\left|x_{0}^{(2)}+\alpha t+2 t \sup _{z \in \mathbb{R}}\left| \frac{\partial_{x x}m(z)}{\partial_{x}m(z)}\right|+\sqrt{2} B_{t}^{(2)}\right| \ge C^{2}\right) \notag\\
	&\le 2 \mathbb{P}\left(\sup _{t \in[0,C]} \sqrt{2}\left|B_{t}\right| \ge C^{2}-C-\alpha C-2 C K_{m}\right)\notag\\
	& \leq 2 \mathbb{P}\left(\sup _{t \in[0,C]}\left|B_{t}\right| \ge \tfrac 12  C^{2}\right) \notag\\
	& \leq 8 \mathbb{P}\left(B_{C} \ge \tfrac 12 C^{2}\right) \notag\\
	& \leq 8 e^{-\frac{1}{8} C^{3}},
	\end{align}
	where the penultimate inequality follows from the reflection principle, and the last inequality from the Gaussian tail bound $1-\Phi(z)\le e^{-\frac 12 z^2}$.
	Substituting~\eqref{eq:ThitCbound} and~\eqref{eq:tauCbound} into~\eqref{eq:Thittauunion}, it follows that, for $C$ sufficiently large,
	\[
	\smalla{hitifclose}(C):=
	2  \frac{C^{1 / 2}\left(1+2 K_{m}\right)}{1+C\left(1+2 K_{m}\right)^{2}} \cdot \frac{1}{\sqrt{2 \pi}} e^{-\frac{1}{2} C\left(1+2 K_{m}\right)^{2}}-8 e^{-\frac{1}{8} C^{3}}>0,
	\]
	then
	\[
	\mathbb{P}_{x_{0}^{(1)}, \, x_{0}^{(2)}}(T_{\text{hit}}>\tau) \leq 1-\smalla{hitifclose}(C),
	\]
	which completes the proof.
	\end{proof}

	Then for $k \in \mathbb{N}$, and $y \in \mathbb{R}$,
	\[
	\begin{aligned}
	\mathbb P_{y,\pi}\left(T_{\text{hit}}>\tau_k,\tau_k<\infty\right)
	&=\mathbb E_{y,\pi}\left[\mathbb P\left(\left. T_{\text{hit}}>\tau_k,\tau_k<\infty \right| \mathcal F_{\sigma_k} \right) \mathds{1}_{T_{\text{hit}}>\sigma_k,\sigma_k<\infty}\right]\\
	& \leq(1-\smalla{hitifclose}(C)) \mathbb{P}_{y,\pi}\left(T_{\text{hit}}>\tau_{k-1}, \tau_{k-1}<\infty\right),
	\end{aligned}
	\]
	where the second line follows from Lemma~\ref{lem:hitifclose}.
	Hence by induction,
	\[
	\mathbb{P}_{y, \pi}\left(T_{\text{hit}}>\tau_{k}, \tau_{k}<\infty\right) \leq(1-\smalla{hitifclose}(C))^{k} .
	\]
	Therefore, for $y\in \R$, $t>0$ and $\delta>0$,
	\begin{align} \label{eq:probT>t}
	\mathbb{P}_{y, \pi}(T_{\text{hit}}>t) &\leq \mathbb{P}_{y,\pi} \left(\tau_{\lfloor \delta t \rfloor} \ge t\right)+\mathbb{P}_{y, \pi}\left(\tau_{\lfloor \delta t \rfloor}<t<T_{\text{hit}}\right) \notag \\
	& \leq \mathbb{P}_{y, \pi}\left(\tau_{\lfloor \delta t \rfloor} \ge t\right)+(1-\smalla{hitifclose}(C))^{\lfloor \delta t \rfloor} . 
	\end{align}

	\begin{lemma} \label{lemma:Laplace_bound}
		Let $ X $ denote a non-negative real-valued random variable, and assume that there exist $ A > 0 $ and $ b > 0 $ such that, for all $ x \geq 0 $,
		\begin{equation*}
			\P{X > x} \leq A e^{-bx}.
		\end{equation*}
		Then, for any $ a < b $,
		\begin{equation*}
			\E{e^{a X}} \leq 1 + A \frac{a}{b-a}.
		\end{equation*}
	\end{lemma}
	
	\begin{proof}
		Integrating by parts, we note that
		\begin{equation*}
			\E{e^{a X}} = 1 + a \int_{0}^{\infty} e^{ax} \P{X > x} dx,
		\end{equation*}
		and the result follows.
	\end{proof}
	
	We now bound $\mathbb{P}_{y, y^{\prime}}\left(\sigma_{1} \ge s\right)$.
	
	\begin{lemma} \label{lem:sigma1bound}
	For $\varepsilon>0$, for $C$ sufficiently large, 
	there exist $\newbigC{sigma1}>0$ and $\newsmalla{sigma1}>0$ such that, for all $y_1,$ $y_2\in \R$ and $s>0$,
	\[
	\mathbb{P}_{y_1,y_2}\left(\sigma_{1} \ge s\right) \leq \bigC{sigma1}\left(e^{\frac 12 (2-\alpha-\varepsilon) |y_1|}+ e^{\frac 12 (2-\alpha-\varepsilon) |y_2|}\right) e^{-\smalla{sigma1} s}.
	\]
	\end{lemma}

	\begin{proof}
	We define stopping times $(\rho_{k}^{(1)})_{k=0}^{\infty}$ and $(\rho_{k}^{(2)})_{k=0}^{\infty}$ iteratively as follows.
	(The times $\rho_{k}^{(1)}$ for $k\ge 1$ are times at which $|X_{t}^{(1)}| \leq C/2$ and the times $\rho_{k}^{(2)}$ for $k\ge 1$ are times at which $|X_{t}^{(2)}| \leq C/2$.)
	
	Let $\rho_{0}^{(1)}=0=\rho_{0}^{(2)}$.
	For $i \ge 0$, let 
	\begin{align*}
	\rho_{i+1}^{(1)}&=\inf\left\{t \ge \rho_{i}^{(2)}:\left|X_{t}^{(1)}\right| \leq C/2\right\}\\
	\text{and }\quad \rho_{i+1}^{(2)}&=\inf\left\{t \ge \rho_{i+1}^{(1)}:\left|X_{t}^{(2)}\right| \leq C /2\right\}.
	\end{align*}
	Then for $k \geq 1$,
	\begin{align*}
	\mathbb P\left(\left. |X^{(1)}_{\rho^{(2)}_k}|>C \right| \mathcal F_{\rho^{(1)}_k}\right)
	&=\mathbb E \left[ \left. \mathbb P\left(\left. |X^{(1)}_{\rho^{(2)}_k}|>C \right| \sigma(\mathcal F^{(1)}_{\rho^{(1)}_k},\mathcal F^{(2)}_{\rho^{(2)}_k})\right)\right| \mathcal F_{\rho^{(1)}_k}\right]\\
	&=\mathbb E \left[ \left. \mathbb P_{X^{(1)}_{\rho^{(1)}_k}}\left(\left. |X^{(1)}_{\rho^{(2)}_k-\rho^{(1)}_k}|>C \right| \rho^{(2)}_k \right)\right| \mathcal F_{\rho^{(1)}_k}\right]\\
	&\le 2 e^{\frac 12 (2-\alpha-\varepsilon)\frac C 2} e^{-\frac 12 (2-\alpha-\varepsilon)\left(C-2 \bigC{driftlimit}(\varepsilon)\right)}\\
	&\le 1/2,
	\end{align*}
	where the first inequality follows from Lemma~\ref{lem:Xttail_initial},
	and the last line follows by taking $C$ sufficiently large.
	Hence for $C$ sufficiently large, for $k \ge 1$,
	\begin{align*}
	\mathbb{P}_{y_1,y_2}\left(\sigma_{1} \ge \rho_{k}^{(2)},\rho_k^{(2)}<\infty \right) &\leq \mathbb{E}_{y_1,y_2}\left[\mathbb{P}\left(\left. \left|X_{\rho_{k}^{(2)}}^{(1)}\right|>C \right| \mathcal F_{\rho_{k}^{(1)}}\right) \mathds{1}_{\sigma_{1}>\rho_{k}^{(1)}, \rho_{k}^{(1)}<\infty}\right] \\
	&\le \tfrac 12 \mathbb{P}_{y_1,y_2}\left(\sigma_{1} \geq \rho_{k-1}^{(2)}, \rho_{k-1}^{(2)}<\infty\right).
	\end{align*}
	Therefore by iterating, we have 
	\begin{equation} \label{eq:sigma1rho2k}
	\mathbb{P}_{y_1,y_2}\left(\sigma_{1} \ge \rho_{k}^{(2)}, \rho_{k}^{(2)}<\infty\right) \leq (\tfrac 12)^{k}.
	\end{equation}
	Hence for $\delta^{\prime}>0$ and $s>0$, by a union bound and then by~\eqref{eq:sigma1rho2k},
	\begin{align} \label{eq:Psigma1}
	\mathbb{P}_{y_1,y_2}\left(\sigma_{1} \ge s\right) &\leq \mathbb{P}_{y_1,y_2}\left(\rho_{\lfloor \delta^{\prime} s\rfloor}^{(2)} \ge s\right)+\mathbb{P}_{y_1,y_2}\left(\rho_{\lfloor \delta ' s\rfloor }^{(2)}<s, \sigma_{1} \ge s\right) \notag \\
	& \leq \mathbb{P}_{y_1,y_2}\left(\rho_{\left\lfloor\delta^{\prime} s\right\rfloor}^{(2)} \ge s\right)+(\tfrac 12)^{\left\lfloor \delta^{\prime} s\right\rfloor} . 
	\end{align}
	We now bound the first term on the right-hand side of~\eqref{eq:Psigma1}.
	
	Assume that $C$ is sufficiently large that $C / 2>\bigC{driftlimit}(\varepsilon)$. Then for $|y| \ge C / 2$ and $t>0$ with $(2-\alpha-\varepsilon) t+C / 2 \ge|y|$,
	letting $\left(B_{s}\right)_{s \ge 0}$ denote a Brownian motion, we can write
	\begin{align} \label{eq:Xstayabove}
	 \mathbb{P}_{y}\left(\left|X_{s}\right| \ge C / 2 \; \forall s \leq t\right) &\leq \mathbb{P}_{|y|}\left(\sqrt{2} B_{s}-(2-\alpha-\varepsilon) s \ge C / 2 \; \forall s \leq t\right) \notag \\
	& \leq \mathbb{P}_{|y|}\left(\sqrt{2} B_{t} \ge(2-\alpha-\varepsilon) t+\tfrac 12 C \right) \notag \\
	& \leq e^{-\frac{1}{4 t}((2-\alpha-\varepsilon) t+\frac C  2-|y|)^{2}} \notag \\
	& \leq e^{-\frac 14 (2-\alpha-\varepsilon)^{2} t+\frac 12 (2-\alpha-\varepsilon)(|y|-\frac C 2)},
	\end{align}
	where the penultimate line follows by a Gaussian tail bound, since $(2-\alpha-\varepsilon) t+C / 2 \ge|y|$.
	Note that in fact if $(2-\alpha-\varepsilon) t+C / 2 <|y|$ or if $|y| < C / 2$ then~\eqref{eq:Xstayabove} trivially holds.
	
	Take $k\ge 2$ and $s>0$. Then for any $y,y'\in \mathbb R$, and any $r>0$, applying~\eqref{eq:Xstayabove} in the second line,
	\begin{align} \label{eq:rho2k-rho1k>s}
	\mathbb{P}_{y, y'}\left(\left. \rho_{k}^{(2)}-\rho_{k}^{(1)}>s \right| \mathcal F_{\rho_{k-1}^{(2)}}\right) 
	& =\mathbb{E}_{y,y^{\prime}}\left[\mathbb{P}\left(\left. \left. \rho_{k}^{(2)}-\rho_{k}^{(1)}>s\right|\mathcal F_{\rho_{k}^{(1)}}\right) \right| \mathcal F_{\rho_{k-1}^{(2)}}\right] \notag \\
	& \leq \mathbb{E}_{y, y^{\prime}}\left[\left.\mathds{1}_{\{|X_{\rho_{k}^{(1)}}^{(2)}| \ge r\}}+e^{-\frac 14 (2-\alpha-\varepsilon) ^{2} s+\frac{1}{2} (2-\alpha-\varepsilon) \left(r-\frac C 2\right)} \right| \mathcal{F}_{\rho_{k-1}^{(2)}}\right] \notag \\
	& =\mathbb{P}_{y , y^{\prime}}\left(\left. |X_{\rho_{k}^{(1)}}^{(2)}| \ge r \right| \mathcal F_{\rho_{k-1}^{(2)}}\right)+e^{-\frac 14 (2-\alpha-\varepsilon)^{2} s+\frac 12 (2-\alpha-\varepsilon) (r-\frac C 2)} .
	\end{align}
	For the first term on the right-hand side of~\eqref{eq:rho2k-rho1k>s}, we have
	\[
	\begin{aligned}
	\mathbb{P}_{y , y^{\prime}}\left(\left. |X_{\rho_{k}^{(1)}}^{(2)}| \ge r \right| \mathcal F_{\rho_{k-1}^{(2)}}\right)
	& =\mathbb E_{y,y'}\left[\left. \mathbb{P}\left(\left. |X_{\rho_{k}^{(1)}}^{(2)}| \ge r\right|\sigma\left(\rho_{k}^{(1)}, \mathcal{F}_{\rho_{k-1}^{(2)}}\right) \right)\right| \mathcal{F}_{\rho_{k-1}^{(2)}}\right] \\
	& \leq 2 e^{\frac 12  (2-\alpha-\varepsilon) \frac C 2} e^{-\frac 12 (2-\alpha-\varepsilon)\left(r-2 \bigC{driftlimit}(\varepsilon)\right)} 
	\end{aligned}
	\]
	by Lemma~\ref{lem:Xttail_initial}.
	Hence, 
	taking $r=\frac 14 (2-\alpha-\varepsilon) s$ in~\eqref{eq:rho2k-rho1k>s}, it follows that
	\[
	\begin{aligned}
	 \mathbb{P}_{y, y^{\prime}}\left(\left. \rho_{k}^{(2)}-\rho_{k}^{(1)}>s \right| \mathcal F_{\rho_{k-1}^{(2)}}\right) 
	& \leq 2 e^{\frac 14 (2-\alpha-\varepsilon) C+(2-\alpha-\varepsilon) \bigC{driftlimit}(\varepsilon)} e^{-\frac 18 (2-\alpha-\varepsilon)^{2} s}+e^{-\frac 14 (2-\alpha-\varepsilon)^{2} s+\frac 18 (2-\alpha-\varepsilon)^{2} s}\\
	&=\left(2 e^{\frac 14 (2-\alpha-\varepsilon) C+(2-\alpha-\varepsilon) \bigC{driftlimit}(\varepsilon)}+1\right) e^{-\frac 18 (2-\alpha-\varepsilon)^{2} s}.
	\end{aligned}
	\]
	Therefore, by Lemma~\ref{lemma:Laplace_bound},
	\begin{align} \label{eq:exponential_momentrho2rho1k}
	\mathbb{E}_{y,y'}\left[\left. e^{\frac 1 {10} (2-\alpha-\varepsilon)^{2}\left(\rho_{k}^{(2)}-\rho_{k}^{(1)}\right)} \right| \mathcal F_{\rho_{k-1}^{(2)}}\right] & \leq 1 + 4 \left(2 e^{\frac 14 (2-\alpha-\varepsilon) C+(2-\alpha-\varepsilon) \bigC{driftlimit}(\varepsilon)}+1\right) = K_{1}(C,\varepsilon),
	\end{align}
	where $K_1<\infty$ is a constant depending on $C$ and $\varepsilon$.
	
	Then by the same argument as for~\eqref{eq:rho2k-rho1k>s}, for any $r>0$,
	\[
	\mathbb{P}_{y,y^{\prime}}\left(\rho_{1}^{(2)}-\rho_{1}^{(1)}>s\right)
	\leq \mathbb{P}_{y, y^{\prime}}\left(|X_{\rho_{1}^{(1)}}^{(2)}| \ge r\right)+e^{-\frac 14 (2-\alpha-\varepsilon)^{2} s+\frac 12 (2-\alpha-\varepsilon) (r-\frac C 2)}.
	\]
	By Lemma~\ref{lem:Xttail_initial}, we can bound the first term on the right-hand side by writing
	\[\mathbb{P}_{y, y'}\left(|X_{\rho_{1}^{(1)}}^{(2)}| \ge r\right) \leq
	2 e^{\frac 12  (2-\alpha-\varepsilon) |y'|} e^{-\frac 12 (2-\alpha-\varepsilon)\left(r-2 \bigC{driftlimit}(\varepsilon)\right)} 
	.\]
	Hence setting $r=\frac 14 (2-\alpha -\varepsilon) s$, we have that
	\[
	\begin{aligned}
	\mathbb{P}_{y, y^{\prime}}\left(\rho_{1}^{(2)}-\rho_{1}^{(1)}>s\right) &\leq 
	2 e^{\frac 12  (2-\alpha-\varepsilon) |y'|+(2-\alpha-\varepsilon)\bigC{driftlimit}(\varepsilon)} e^{-\frac 18 (2-\alpha-\varepsilon)^2 s} 
	+
	e^{-\frac 18 (2-\alpha-\varepsilon)^{2} s} \\
	& =\left(2 e^{\frac 12  (2-\alpha-\varepsilon) |y'|+(2-\alpha-\varepsilon)\bigC{driftlimit}(\varepsilon)}+1\right) e^{-\frac 18  (2-\alpha-\varepsilon)^{2} s} .
	\end{aligned}
	\]
	By Lemma~\ref{lemma:Laplace_bound}, it follows that
	\begin{align} \label{eq:exponential_momentrho2rho1_1}
	\mathbb{E}_{y, y^{\prime}}\left[e^{\frac 1 {10}  (2-\alpha-\varepsilon)^{2}\left(\rho_1^{(2)}-\rho_1^{(1)}\right)}\right] &\leq 1 + 4 \left(2 e^{\frac 12  (2-\alpha-\varepsilon) |y'|+(2-\alpha-\varepsilon)\bigC{driftlimit}(\varepsilon)}+1\right) \notag \\
	& \leq K_{2} e^{\frac 12 (2-\alpha-\varepsilon) |y^{\prime}|},
	\end{align}
	where $K_2<\infty$ is a constant depending on $C$ and $\varepsilon$.
	
	We now have that
	for $\overline{k}\in \mathbb N$ with $\overline{k}\ge 2$,
	applying~\eqref{eq:exponential_momentrho2rho1k} successively,
	\begin{align*}
	\mathbb E_{y,y'}\left[ e^{\frac 1 {10}(2-\alpha-\varepsilon)^2 \sum_{k=1}^{\overline k}(\rho^{(2)}_k-\rho^{(1)}_k)} \right]
	& =\mathbb{E}_{y, y^{\prime}}\left[e^{\frac 1 {10} (2-\alpha-\varepsilon)^{2} \sum_{k=1}^{\overline{k}-1}\left(\rho_{k}^{(2)}-\rho_{k}^{(1)}\right)} \mathbb{E}\left[\left.  e^{\frac 1 {10} (2-\alpha-\varepsilon)^{2}\left(\rho_{\overline k}^{(2)}-\rho_{\overline k}^{(1)}\right)}\right| \mathcal F_{\rho_{\overline k-1}^{(2)}}\right]\right] \\
	& \le K_1 \mathbb{E}_{y, y^{\prime}}\left[e^{\frac 1 {10} (2-\alpha-\varepsilon)^{2} \sum_{k=1}^{\overline{k}-1}\left(\rho_{k}^{(2)}-\rho_{k}^{(1)}\right)} \right] \\
	&\le \ldots \le K_1^{\overline k-1} \mathbb{E}_{y, y^{\prime}}\left[e^{\frac 1 {10} (2-\alpha-\varepsilon)^{2} \left(\rho_{1}^{(2)}-\rho_{1}^{(1)}\right)} \right].
	\end{align*}
	Hence by~\eqref{eq:exponential_momentrho2rho1_1}, for $\overline k \ge 1$,
	\begin{equation} \label{eq:exponential21}
	\mathbb E_{y,y'}\left[ e^{\frac 1 {10}(2-\alpha-\varepsilon)^2 \sum_{k=1}^{\overline k}(\rho^{(2)}_k-\rho^{(1)}_k)} \right]
	 \leq K_{1}^{\overline k-1} K_{2} e^{\frac 12 (2-\alpha-\varepsilon)\left|y^{\prime}\right|} .
	\end{equation}
	
	By the same argument as for~\eqref{eq:exponential_momentrho2rho1k}, there exists a constant $K_3<\infty$ depending on $C$ and $\varepsilon$ such that for $k \geq 1$,
	\[
	\mathbb{E}_{y,y^{\prime}}\left[\left. e^{\frac 1 {10} (2-\alpha-\varepsilon)^{2}\left(\rho_{k+1}^{(1)}-\rho_{k}^{(2)}\right)} \right| \mathcal{F}_{\rho_{k}^{(1)}}\right] \leq K_{3}.
	\]
	Moreover, by~\eqref{eq:Xstayabove}, for $s>0$
	we have  $\mathbb{P}_{y, y^{\prime}}\left(\rho_{1}^{(1)}-\rho_{0}^{(2)}>s\right) \leq e^{-\frac 14 (2-\alpha-\varepsilon)^{2} s+\frac 12  (2-\alpha-\varepsilon)(|y|-C/ 2)}$, and so in particular, there exists a constant $K_4<\infty$ depending on $C$ and $\varepsilon$ such that
	\[
	\mathbb{E}_{y, y^{\prime}}\left[e^{\frac 1 {10} (2-\alpha-\varepsilon)^{2}\left(\rho^{(1)}_1-\rho_{0}^{(2)}\right)}\right] \leq K_{4}e^{\frac 12 (2-\alpha-\varepsilon) |y|} .
	\]
	Hence for $\overline k \ge 1$,
	\begin{equation} \label{eq:exponential12}
	\mathbb{E}_{y, y^{\prime}}\left[e^{\frac 1 {10} (2-\alpha-\varepsilon)^{2} \sum_{k=1}^{\overline k}\left(\rho^{(1)}_k-\rho_{k-1}^{(2)}\right)}\right] \leq K_{3}^{\overline k-1} K_{4} e^{\frac 12 (2-\alpha-\varepsilon)|y|}.
	\end{equation}
	
	We are now ready to bound the first term on the right-hand side of~\eqref{eq:Psigma1}.
	For $y, y^{\prime} \in \mathbb{R}$, $s>0$ and $\delta^{\prime}>0$, by a union bound, then by Markov's inequality, and finally by~\eqref{eq:exponential21} and~\eqref{eq:exponential12},
	\begin{align} 
	 \mathbb{P}_{y,y'}\left(\rho_{\left\lfloor\delta^{\prime} s\right\rfloor}^{(2)} \ge s\right)&\leq \mathbb{P}_{y, y^{\prime}}\left(\sum_{k=1}^{\lfloor \delta ' s \rfloor}\left(\rho_{k}^{(2)}-\rho_{k}^{(1)}\right) \ge s / 2\right)  +\mathbb{P}_{y, y^{\prime}}\left(\sum_{k=1}^{\lfloor \delta ' s \rfloor}\left(\rho_{k}^{(1)}-\rho_{k-1}^{(2)}\right) \ge s / 2\right)  \notag\\
	& \leq e^{-\frac 1 {10} (2-\alpha-\varepsilon)^{2} \cdot \frac s 2} \mathbb{E}_{y,y'}\left[e^{\frac 1 {10} (2-\alpha-\varepsilon)^{2} \sum_{k=1}^{\lfloor \delta^{\prime} s\rfloor}\left(\rho_{k}^{(2)}-\rho_{k}^{(1)}\right)}\right] \notag\\
	& \quad +e^{-\frac 1 {10} (2-\alpha-\varepsilon)^{2} \cdot \frac s 2} \mathbb{E}_{y,y'}\left[e^{\frac 1 {10} (2-\alpha-\varepsilon)^{2} \sum_{k=1}^{\lfloor \delta^{\prime} s\rfloor}\left(\rho_{k}^{(1)}-\rho_{k-1}^{(2)}\right)}\right]\notag\\
	& \leq e^{-\frac{1}{20} (2-\alpha-\varepsilon)^{2} s} \left(K_{1}^{\lfloor \delta^{\prime} s\rfloor-1} K_{2} e^{\frac 12 (2-\alpha-\varepsilon)\left|y^{\prime}\right|}+K_{3}^{\lfloor \delta^{\prime} s\rfloor-1} K_{4} e^{\frac 12 (2-\alpha-\varepsilon)|y|}\right). \notag 
	\end{align}
	By taking
	$\delta'>0$ sufficiently small that $\delta^{\prime} \log (K_{1}\vee K_3)<\frac{1}{20} (2-\alpha-\varepsilon)^{2}$, the result now follows from~\eqref{eq:Psigma1}.
	\end{proof}
	
	We now bound $\mathbb{P}\left(\left. \tau_{1}-\sigma_{1} \ge s \right|\mathcal F_{\sigma_1}\right)$. 
	
	\begin{lemma}  \label{lem:tau1sigma1}
	For $C$ sufficiently large, there exists $\newsmalla{tau1}=\smalla{tau1}(C)>0$ such that, for all $s>0$, on the event $\sigma_1<\infty$,
	\[
	\mathbb{P}\left(\left. \tau_{1}-\sigma_{1} \geq s\right|\mathcal{F}_{\sigma_{1}}\right) \leq e^{-2\smalla{tau1} \lfloor s\rfloor},
	\]
	almost surely.
	\end{lemma}
	\begin{proof}
	For any $z, z^{\prime} \in \mathbb{R}$, we have
	\[
	\begin{aligned}
	\mathbb{P}_{z, z^{\prime}}\left(|X_{t}^{(1)}| \vee |X_{t}^{(2)}| \leq C^{2} \; \forall t \leq 1\right) 
	& =\mathbb{P}_{z}\left(\left|X_{t}\right| \leq C^{2} \; \forall t \leq 1\right) \mathbb{P}_{z^{\prime}}\left(\left|X_{t}\right| \leq C^{2} \;\forall t \leq 1\right).
	\end{aligned}
	\]
	Let $K_m=\sup _{z \in \mathbb{R}}\left| \frac{\partial_{x x}m(z)}{\partial_{x}m(z)}\right|$.
	Suppose $|z| \leq C^{2}$.
	Letting $\left(B_{t}\right)_{t \ge 0}$ denote a Brownian motion, we can write
	\[
	\mathbb{P}_{z}\left(\left|X_{t}\right| \leq C^{2} \; \forall t \leq 1\right) \leq \mathbb{P}\left(B_{1} \leq C^{2}+|z|+\alpha+2 K_{m}\right)
	\leq \mathbb{P}\left(B_{1} \leq 3 C^{2}\right) =: e^{-\smalla{tau1}(C)},
	\]
	for some $\smalla{tau1}(C)>0$, where the second inequality holds for $C$ sufficiently large.
	For $|z|>C^{2}$, we have \[\mathbb{P}_{z}\left(\left|X_{t}\right| \leq C^{2} \; \forall t \leq 1\right)=0.\]
	Hence
	\[
	\mathbb{P}_{z, z^{\prime}}\left(|X_{t}^{(1)}| \vee |X_{t}^{(2)}| \leq C^{2} \; \forall t \leq 1\right) \leq e^{-2\smalla{tau1}(C)} .
	 \]
	Therefore, on the event $\sigma_{1}<\infty$,
	\[
	\begin{aligned}
	\mathbb P\left(\left. \tau_1-\sigma_1 \ge s \right| \mathcal F_{\sigma_1}\right)
	& =\mathbb{E}\left[\left. \mathbb{P}\left(\left. \tau_{1}-\sigma_{1} \ge s \right| \mathcal F_{\sigma_1+s-1}\right) \mathds{1}_{\tau_{1} > \sigma_{1}+s-1} \right| \mathcal F_{\sigma_{1}}\right] \\
	& \leq e^{-2\smalla{tau1}(C)} \mathbb{P}_{y, y^{\prime}}\left(\left. \tau_{1}-\sigma_{1}>s-1 \right|\mathcal F_{ \sigma_{1}}\right),
	\end{aligned}
	\]
	and the result follows by iterating the same argument.
	\end{proof}

	\begin{lemma} \label{lem:Thitexpon}
	For $\varepsilon\in (0,2-\alpha)$, there exist constants $\newbigC{Thit}>0$ and $\newsmalla{Thit}>0$ such that for $y\in \R$ and $t>0$,
	\[
	\mathbb{P}_{y, \pi}(T_{\text{hit}}>t) \leq  \bigC{Thit} e^{\frac 12 (2-\alpha-\varepsilon)|y|} e^{-\smalla{Thit} t}.
	\]
	\end{lemma}
	
	\begin{proof}
	Take $C$ sufficiently large that Lemma~\ref{lem:sigma1bound} and Lemma~\ref{lem:tau1sigma1} hold.
	Then by Lemma~\ref{lem:sigma1bound} and Lemma~\ref{lemma:Laplace_bound}, for $0 < a < \smalla{sigma1}$ and $y_1,y_2 \in \mathbb{R}$,
	\begin{align} \label{eq:expsigma1}
	\mathbb{E}_{y_1,y_2}\left[e^{a\sigma_{1}}\right] &\leq 1 + \bigC{sigma1} \left(e^{\frac 12 (2-\alpha-\varepsilon) |y_1|}+ e^{\frac 12 (2-\alpha-\varepsilon) |y_2|}\right) \frac{a}{\smalla{sigma1}-a} \notag \\
	&\leq K_1 \left(e^{\frac 12 (2-\alpha-\varepsilon) |y_1|}+ e^{\frac 12 (2-\alpha-\varepsilon) |y_2|}\right).
	\end{align}
	By Lemma~\ref{lem:tau1sigma1}, on the event $\sigma_{1}<\infty$, for  $0< a < 2 \smalla{tau1}(C)$,
	\begin{align} \label{eq:tau1exp}
	\mathbb{E}\left[\left. e^{a\left(\tau_{1}-\sigma_{1}\right)} \right| \mathcal F_{\sigma_{1}}\right] \leq 1 + \frac{a}{\smalla{tau1}(C)-a} = K_2
	\end{align}
	Hence for $y_1,y_2 \in \mathbb{R}$ and $k \in \mathbb{N}$, $k\ge 2$, for $a < \smalla{sigma1} \wedge 2 \smalla{tau1}(C)$, using~\eqref{eq:tau1exp} in the first inequality,
	\begin{align} \label{eq:tauktok-1}
	\mathbb{E}_{y_1,y_2}\left[e^{a \tau_{k}}\right] &=\mathbb{E}_{y_1,y_2}\left[e^{a \sigma_{k}} \mathbb{E}\left[\left. e^{a\left(\tau_{k}-\sigma_{k}\right)} \right| \mathcal F_{\sigma_{k}}\right]\right] \notag \\
	& \leq K_2 \mathbb{E}_{y_1,y_2}\left[e^{a \sigma_{k}}\right] \notag \\
	& =K_2 \mathbb{E}_{y_1,y_2}\left[e^{a \tau_{k-1}} \mathbb{E}\left[\left. e^{a\left(\sigma_{k}-\tau_{k-1}\right)} \right| \mathcal{F}_{\tau_{k-1}}\right]\right] \notag \\
	& \leq K_2 K_1 \left(e^{\frac 12 (2-\alpha-\varepsilon)C^2}+ e^{\frac 12 (2-\alpha-\varepsilon) C^2}\right)\mathbb{E}_{y_1,y_2}\left[e^{a \tau_{k-1}}\right],
	\end{align}
	where the last inequality follows from~\eqref{eq:expsigma1}  since $|X_{\tau_{k-1}}^{(1)}| \vee |X_{\tau_{k-1}}^{(2)}| \leq C^{2}$ (because $k \ge 2$).
	Therefore, by applying~\eqref{eq:tauktok-1} repeatedly, and then by~\eqref{eq:tau1exp} and then~\eqref{eq:expsigma1}, there exists a constant $K_3>0$ such that for $y_1,y_2 \in \mathbb{R},$ $ k \in \mathbb{N}$ and $a < \smalla{sigma1} \wedge 2 \smalla{tau1}(C)$,
	\[
	\begin{aligned}
	\mathbb{E}_{y_1,y_2}\left[e^{a \tau_{k}}\right] & \leq K_{3}^{k-1} \mathbb{E}_{y_1,y_2}\left[e^{a \tau_{1}}\right] \\
	& \leq K_{3}^{k-1} K_2 \mathbb{E}_{y_1,y_2}\left[e^{a \sigma_{1}}\right] \\
	& \leq K_{3}^{k-1} K_2 K_1 \left(e^{\frac 12 (2-\alpha-\varepsilon) |y_1|}+ e^{\frac 12 (2-\alpha-\varepsilon) |y_2|}\right).
	\end{aligned}
	\]
	Therefore, there exists a constant $K_4>0$ such that for $y_1,y_2 \in \mathbb{R}$ and $\delta, t>0$, letting $a =\left(\frac 12 \smalla{sigma1}\right) \wedge \smalla{tau1}(C) >0$,
	\[
	 \begin{aligned}
	\mathbb{P}_{y_1,y_2}\left(\tau_{\lfloor \delta t \rfloor} \ge t\right) &\leq e^{-a t} \mathbb{E}_{y_1,y_2}\left[e^{a \tau_{\lfloor \delta t \rfloor }}\right] \\
	&\leq e^{-a t} K_{3}^{\lfloor \delta t \rfloor -1} K_4 \left(e^{\frac 12 (2-\alpha-\varepsilon) |y_1|}+ e^{\frac 12 (2-\alpha-\varepsilon) |y_2|}\right).
	\end{aligned}
	\]
	Hence for $y\in \R$ and $t>0$,
	\begin{align*}
	\mathbb{P}_{y,\pi}\left(\tau_{\lfloor \delta t\rfloor} \ge t\right) &\leq e^{-a t} 
	K_{3}^{\lfloor \delta t \rfloor -1} K_4\left(e^{\frac 12 (2-\alpha-\varepsilon) |y|}+ \int_{\R}e^{\frac 12 (2-\alpha-\varepsilon) |z|}\pi(z) dz\right)\\
	&\leq e^{-a t} 
	K_{3}^{\lfloor \delta t \rfloor -1} K_4\left(e^{\frac 12 (2-\alpha-\varepsilon) |y|}+ K_5\right)\\
	\end{align*}
	for some constant $K_5>0$, since $\varepsilon>0$ and $\pi(x) \sim e^{-(2-\alpha) x}$ as $x\to \infty$, and $\pi(x) \sim e^{(2+\alpha) x}$ as $x\to - \infty$.
	Therefore by~\eqref{eq:probT>t}, for any $y\in \R$, $t>0$ and $\delta>0$,
	\[
	\mathbb{P}_{y, \pi}(T_{\text{hit}}>t)\leq e^{-a t} 
	K_{1}^{\lfloor \delta t \rfloor -1} K_2\left(e^{\frac 12 (2-\alpha-\varepsilon) |y|}+ K_3\right)
	+(1-\smalla{hitifclose}(C))^{\lfloor \delta t \rfloor} 
	.
	\]
	By taking $\delta>0$ sufficiently small, the result follows.
	\end{proof}
	
	We now establish stronger bounds on $\mathbb{P}_{y}\left(X_{t}\le -x\right)$ and $\mathbb{P}_{y}\left(X_{t}\ge x\right)$ for $x>0$. 
	\begin{lemma} \label{lem:Xttailstronger}
	For $\varepsilon_1\in (0,2+\alpha)$ and $\delta\in (0,1)$, there exists $\newbigC{Xttail-}>0$ such that for $x>0$, $y \in \mathbb{R}$ and $t \ge 0,$
	\begin{equation} \label{eq:Xttail-stronger}
	\mathbb{P}_{y}\left(X_{t} \leq-x\right) \leq \bigC{Xttail-}(t+1) e^{\frac 12 (2+\alpha-\varepsilon_1)\left((-y) \vee 0\right)}e^{-\frac 12(1-\delta)(2+\alpha-\varepsilon_1) x}.
	\end{equation}
	For $\varepsilon_2\in (0,2-\alpha)$ and $\delta\in (0,1)$, there exists $\newbigC{Xttail+}>0$ such that for $x>0,$ $ y \in \mathbb{R}$ and $t \ge 0$,
	\begin{equation} \label{eq:Xttail+stronger}
	\mathbb{P}_{y}\left(X_{t} \ge x\right) \leq \bigC{Xttail+}(t+1) e^{\frac 12\left(2-\alpha-\varepsilon_{2}\right)(y \vee 0)} e^{-\frac 12(1-\delta)(2-\alpha-\varepsilon_2) x}.
	\end{equation}
	\end{lemma}
	
	\begin{proof}
	For $x>0, $ $y \in \mathbb{R}$ and $t \ge 0$,
	by a union bound we have 
	\begin{align} \label{eq:Xttailunion}
	\mathbb{P}_{y}\left(X_{t} \leq-x\right) &\leq \sum_{s=0}^{\lfloor t \rfloor} \mathbb{P}_{y}\left(X_{t} \leq-x, X_{s^{\prime}} \leq-\bigC{driftlimit}(\varepsilon_1) \; \forall s^{\prime} \in[(s+1) \wedge t, t],\right. \notag \\
	& \hspace{4cm}\left.X_{s^{*}} \geq-\bigC{driftlimit}(\varepsilon_1) \text { for some } s^{*} \in[s,(s+1) \wedge t]\right) \notag \\
	& \qquad +\mathbb{P}_{y}\left(X_{t} \leq-x, X_{s^{\prime}} \leq-\bigC{driftlimit}(\varepsilon_1)\; \forall s^{\prime} \in[0, t]\right) .
	\end{align}
	We now bound the terms in the sum on the right-hand side of~\eqref{eq:Xttailunion}.
	For $\delta\in (0,1)$ and $s\in [0,t]$, by a union bound we have
	\begin{align*} 
	& \mathbb{P}_{y}\left(X_{t} \leq-x, X_{s^{\prime}} \leq-\bigC{driftlimit}(\varepsilon_1) \; \forall s^{\prime} \in[(s+1) \wedge t, t],\right. \notag \\
	& \hspace{4cm}\left.X_{s^{*}} \geq-\bigC{driftlimit}(\varepsilon_1) \text { for some } s^{*} \in[s,(s+1) \wedge t]\right)\notag \\
	&\le \mathbb{P}_{y}\left(X_{t} \leq-x, X_{s^{\prime}} \leq-\bigC{driftlimit}(\varepsilon_1) \; \forall s^{\prime} \in[(s+1) \wedge t, t], X_{(s+1) \wedge t} \ge-\delta x\right)\notag  \\
	&\qquad \qquad  +\mathbb{P}_{y}\left(X_{s^{*}} \geq-\bigC{driftlimit}(\varepsilon_1) \text { for some } s^{*} \in[s,(s+1) \wedge t], X_{(s+1) \wedge t}<-\delta x\right).
	\end{align*}
	Letting $\left(B_{s}\right)_{s \ge 0}$ denote a Brownian notion, and letting $K_m=\sup _{z \in \mathbb{R}}\left| \frac{\partial_{zz}m(z)}{\partial_{z}m(z)}\right|$, it follows that
	\begin{align} \label{eq:Xttailunion2}
	& \mathbb{P}_{y}\left(X_{t} \leq-x, X_{s^{\prime}} \leq-\bigC{driftlimit}(\varepsilon_1) \; \forall s^{\prime} \in[(s+1) \wedge t, t],\right. \notag \\
	& \hspace{4cm}\left.X_{s^{*}} \geq-\bigC{driftlimit}(\varepsilon_1) \text { for some } s^{*} \in[s,(s+1) \wedge t]\right)\notag \\
	& \leq \mathbb{P}\left(-\delta x+\sqrt{2} B_{t-(s+1) \wedge t}+\left(2+\alpha-\varepsilon_{1}\right)(t-(s+1) \wedge t) \leq-x\right) \notag \\
	& \hspace{1cm}+\mathbb{P}\left(\sup _{u \in[0,1]}\left(\sqrt{2} B_{u}\right) \ge \delta x-\bigC{driftlimit}(\varepsilon_1)-2K_m\right)\notag  \\
	& \leq e^{-(4(t-(s+1) \wedge t))^{-1}\left((1-\delta) x+\left(2+\alpha-\varepsilon_{1}\right)(t-(s+1) \wedge t)\right)^{2}}+2 e^{-\frac{1}{4}(\delta x-\bigC{driftlimit}(\varepsilon_1)-2K_m)^{2}}, 
	\end{align}
	where the last inequality holds for $x>0$ sufficiently large that $\delta x-\bigC{driftlimit}(\varepsilon_1)-2K_m\ge 0$, by a Gaussian tail bound and the reflection principle.
	Similarly, for the last term on the right-hand side of~\eqref{eq:Xttailunion},
	we have
	\begin{align} \label{eq:Xttailunion3}
	\mathbb{P}_{y}\left(X_{t} \leq-x, X_{s^{\prime}} \leq-\bigC{driftlimit}(\varepsilon_1)\; \forall s^{\prime} \in[0, t]\right) 
	&\le \mathbb{P}\left(y+\sqrt{2} B_{t}+\left(2+\alpha-\varepsilon_{1}\right) t \leq-x\right) \notag \\
	&\le e^{-\frac{1}{4 t}\left(x+y+\left(2+\alpha-\varepsilon_{1}\right) t\right)^{2}},
	\end{align}
	where the last inequality holds if $x+y+\left(2+\alpha-\varepsilon_{1}\right) t\ge 0$, by a Gaussian tail bound.
	
	Substituting~\eqref{eq:Xttailunion2} and~\eqref{eq:Xttailunion3}
	into~\eqref{eq:Xttailunion}, it follows that if $x+y+\left(2+\alpha-\varepsilon_{1}\right) t\ge 0$ and $\delta x-\bigC{driftlimit}(\varepsilon_1)-2K_m\ge 0$, then
	\[
	\begin{aligned}
	& \mathbb{P}_{y}\left(X_{t}\le -x\right) \\
	& \leq \sum_{s=0}^{\lfloor t \rfloor}\left(e^{-(4(t-(s+1) \wedge t))^{-1} \cdot 2(1-\delta) x(2+\alpha-\varepsilon_1)(t-(s+1) \wedge t)}+2 e^{-\frac{1}{4}(\delta x-\bigC{driftlimit}(\varepsilon_1)-2K_m)^{2}}\right) \\
	& \qquad +e^{-\frac{1}{4 t} \cdot 2(x+y)(2+\alpha-\varepsilon_1) t} \\
	& \le (\lfloor t \rfloor+1)\left(e^{-\frac{1}{2}(1-\delta)\left(2+\alpha-\varepsilon_{1}\right) x}+2 e^{-\frac 1{16} \delta^{2} x^{2}}\right)+e^{-\frac{1}{2}\left(2+\alpha-\varepsilon_{1}\right)(x+y)},
	\end{aligned}
	\]
	where the last inequality holds if $\frac 12 \delta x \ge \bigC{driftlimit}(\varepsilon_1)-2K_m$.
	Therefore, if $\frac 1 {16}\delta^2 x\ge \frac 12 (2+\alpha)$, and $\frac 12 \delta x \ge \bigC{driftlimit}(\varepsilon_1)-2K_m$, and $x+y+\left(2+\alpha-\varepsilon_{1}\right) t\ge 0$, we have
	\[
	\begin{aligned}
	\mathbb{P}_{y}\left(X_{t}\le -x\right) 
	& \leq (t+1) \cdot 3 e^{-\frac{1}{2}(1-\delta)\left(2+\alpha-\varepsilon_1\right) x}+e^{-\frac{1}{2}\left(2+\alpha-\varepsilon_{1}\right) y} e^{-\frac{1}{2}\left(2+\alpha-\varepsilon_{1}\right) x}\\
	&\leq (3 t+4) e^{\frac 12(2+\alpha-\varepsilon_1)((-y) \vee 0)} e^{-\frac{1}{2}(1-\delta)(2+\alpha-\varepsilon_1) x}.
	\end{aligned}
	\]
	If $y \leq-x$ then the right-hand side above is greater than $1$, and so there exists a constant $\bigC{Xttail-}>0$
	such that~\eqref{eq:Xttail-stronger}
	holds for all $x>0$ and $y\in \R$.
	Then~\eqref{eq:Xttail+stronger} follows from the same argument.
	\end{proof}
	
	\begin{cor} \label{cor:XandThit}
	For $\varepsilon_1\in (0,2+\alpha)$, $\varepsilon_2\in (0,2-\alpha)$ and $\delta\in (0,1/2)$, there exist $\newsmalla{XttailThit}>0$ and $\newbigC{XttailThit}>0$ such that for $x>0$, $y \in \mathbb{R}$ and $t \ge 0,$
	\begin{align}
	\mathbb{P}_{y, \pi}\left(X_{t}^{(1)} \leq-x, T_{\text{hit}}>t\right) &\leq \bigC{XttailThit} e^{\delta |y|+\frac 12\left(2+\alpha-\varepsilon_{1}\right)((-y) \vee 0)} e^{-\smalla{XttailThit} t}e^{-\frac 12(1-2 \delta)\left(2+\alpha-\varepsilon_1\right)x},\label{eq:X1<xThit}\\
	\mathbb{P}_{y, \pi}\left(X_{t}^{(2)} \leq-x, T_{\text{hit}}>t\right) &\leq \bigC{XttailThit} e^{\delta |y|+\frac 12\left(2+\alpha-\varepsilon_{1}\right)((-y) \vee 0)} e^{-\smalla{XttailThit} t}e^{-\frac 12(1-2 \delta)\left(2+\alpha-\varepsilon_1\right)x},\label{eq:X2<xThit}\\
	\mathbb{P}_{y, \pi}\left(X_{t}^{(1)} \ge x, T_{\text{hit}}>t\right) &\leq \bigC{XttailThit} e^{\delta |y|+\frac 12 \left(2-\alpha-\varepsilon_{2}\right)(y \vee 0)} e^{-\smalla{XttailThit} t}e^{-\frac 12 (1-2\delta )\left(2-\alpha-\varepsilon_2\right)x}, \label{eq:X1>xThit}\\
	\text{and }\quad \mathbb{P}_{y, \pi}\left(X_{t}^{(2)} \ge x, T_{\text{hit}}>t\right) &\leq \bigC{XttailThit} e^{\delta |y|+\frac 12 \left(2-\alpha-\varepsilon_{2}\right)(y \vee 0)} e^{-\smalla{XttailThit} t}e^{-\frac 12 (1-2\delta )\left(2-\alpha-\varepsilon_2\right)x}. \label{eq:X2>xThit}
	\end{align}
	\end{cor}
	\begin{proof}
	Fix $A>0$. For $x>0$, $y\in \R$ and $t\ge 0$ with $t \ge A x$, by Lemma~\ref{lem:Thitexpon} with $\varepsilon$ chosen so that $\delta=\frac 12 (2-\alpha-\varepsilon)$, we have
	\[
	\mathbb{P}_{y, \pi}\left(X_{t}^{(1)} \ge x, T_{\text{hit}}>t\right) \leq \mathbb{P}_{y, \pi}(T_{\text{hit}}>t)\le \bigC{Thit} e^{\delta|y|} e^{-\frac 12 \smalla{Thit} t} e^{-\frac 12 \smalla{Thit} A x}.
	\]
	For $x>0$, $y\in \R$ and $t\ge 0$ with $t \le A x$, by Lemma~\ref{lem:Xttailstronger},
	\[
	\begin{aligned}
	\mathbb{P}_{y, \pi}\left(X_{t}^{(1)} \ge x, T_{\text{hit}}>t\right)  &\leq \bigC{Xttail+}(t+1) e^{\frac 12 \left(2-\alpha-\varepsilon_{2}\right)(y \vee 0)-\frac 12(1-2 \delta)\left(2-\alpha-\varepsilon_{2}\right) x}  e^{-\frac 12  \delta\left(2-\alpha-\varepsilon_{2}\right) A^{-1} t}.
	\end{aligned}
	\]
	Take $A>0$ such that $\frac 12 \smalla{Thit} A=\frac 12 (1-2 \delta)\left(2-\alpha-\varepsilon_{2}\right)$.
	Then for any $x>0$, $y\in \R$ and $t\ge 0$, we have
	\begin{align*}
	&\mathbb{P}_{y, \pi}\left(X_{t}^{(1)} \ge x, T_{\text{hit}}>t\right) \\
	&\quad \leq 
	(\bigC{Thit}\vee \bigC{Xttail+})(t+1)
	 e^{\delta |y|+\frac 12 \left(2-\alpha-\varepsilon_{2}\right)(y \vee 0)} e^{-(\frac 12 \smalla{Thit})\wedge (\frac 12  \delta\left(2-\alpha-\varepsilon_{2}\right) A^{-1}) t}e^{-\frac 12 (1-2\delta )\left(2-\alpha-\varepsilon_2\right)x},
	\end{align*}
	and~\eqref{eq:X1>xThit} follows.
	By the same argument, we obtain~\eqref{eq:X1<xThit}.
	
	For $y\in \R$, $x>0$ and $t\ge 0$, 
	\[
	\mathbb{P}_{y, \pi}\left(X_{t}^{(2)} \ge x, T_{\text{hit}}>t\right) 
	\le \mathbb{P}_{y, \pi}\left(X_{t}^{(2)} \ge x\right) 
	=\int_{x}^{\infty} \pi(z) d z\le K_1 e^{-(2-\alpha) x}
	\]
	for some constant $K_1>0$, and so~\eqref{eq:X2>xThit} follows by the same argument.
	By the same argument, we obtain~\eqref{eq:X2<xThit}.
	\end{proof}
	
	\begin{proof}[Proof of Lemma~\ref{lem:qtxy_estimate}]
	For $z\ge 0$, $\delta_{z}>0$, $y\in \mathbb{R}$ and $t \ge 1$, by a union bound,
	\begin{equation} \label{eq:Xt1intervalunion}
	\begin{aligned}
	& \mathbb{P}_{y, \pi}\left(X_{t}^{(1)} \in\left[z, z+\delta_{z}\right], T_{\text{hit}}>t\right) \\
	& \quad \leq \mathds{1}_{\{z\ge 1\}} \sum_{k=0}^{\lfloor z \rfloor-1} \mathbb{P}_{y, \pi}\left(X_{t}^{(1)} \in\left[z, z+\delta_{z}\right], T_{\text{hit}}>t-1, X_{t}^{(1)}\in[k, k+1]\right) \\
	&\qquad +\mathbb{P}_{y, \pi}\left(X_{t}^{(1)} \in\left[z, z+\delta_{z}\right] , T_{\text{hit}}>t-1, X_{t-1}^{(1)} \ge\lfloor z\rfloor\right) \\
	&\qquad +\mathbb{P}_{y,\pi}\left(X_{t}^{(1)} \in\left[z, z+\delta_{z}\right], T_{\text{hit}}>t-1, X_{t-1}^{(1)} \leq 0\right).
	\end{aligned}
	\end{equation}
	By Lemma~\ref{lem:qtxy_overall}, we have almost surely
	\begin{equation} \label{eq:usingqtoverall}
	\mathbb{P}_{y,\pi}\left(\left. X_{t}^{(1)} \in\left[z, z+\delta_{z}\right] \right| \mathcal F_{t-1}\right) \leq \delta_{z} \bigC{qtxy} \exp \left(-\smalla{qtxy} \inf_{z^{\prime} \in\left[z, z+\delta_{z}\right]}\left(z^{\prime}-X_{t-1}^{(1)}\right)^{2}\right).
	\end{equation}
	Hence, using~\eqref{eq:X1>xThit} 
	in Corollary~\ref{cor:XandThit} in the second line, 
	for $k\in [0,\lfloor z \rfloor -1]$,
	\begin{align*}
	&\mathbb{P}_{y, \pi}\left(X_{t}^{(1)} \in\left[z, z+\delta_{z}\right], T_{\text{hit}}>t-1, X_{t-1}^{(1)}\in[k, k+1]\right)\\
	&\quad \le \delta_{z} \bigC{qtxy} \exp \left(-\smalla{qtxy} (\lfloor z \rfloor -1 -k)^{2}\right)
	\mathbb{P}_{y, \pi}\left(T_{\text{hit}}>t-1, X_{t-1}^{(1)}\ge k\right)\\
	&\quad \le \delta_{z} \bigC{qtxy} \exp \left(-\smalla{qtxy} (\lfloor z \rfloor -1 -k)^{2}\right)
	\bigC{XttailThit} e^{\delta |y|+\frac 12 \left(2-\alpha-\varepsilon_{2}\right)(y \vee 0)} e^{-\smalla{XttailThit} (t-1)}e^{-\frac 12 (1-2\delta )\left(2-\alpha-\varepsilon_2\right)k}.
	\end{align*}
	Therefore, using~\eqref{eq:X1>xThit} 
	in Corollary~\ref{cor:XandThit} again combined with~\eqref{eq:usingqtoverall}  for the second term on the right-hand side of~\eqref{eq:Xt1intervalunion} and  using Lemma~\ref{lem:Thitexpon} (with $\varepsilon$ chosen so that $\frac 12 (2-\alpha-\varepsilon)=\delta$) combined with~\eqref{eq:usingqtoverall} again for the last term,
	\[
	\begin{aligned}
	&\mathbb{P}_{y, \pi}\left(X_{t}^{(1)} \in[z, z+\delta_z], T_{\text{hit}}>t\right)\\
	&\le \mathds{1}_{\{z\ge 1\}}\sum_{k=0}^{\lfloor z\rfloor-1} \delta_{z} \bigC{qtxy} e^{-\smalla{qtxy} (\lfloor z \rfloor -1 -k)^{2}}
	\bigC{XttailThit} e^{\delta |y|+\frac 12 \left(2-\alpha-\varepsilon_{2}\right)(y \vee 0)} e^{-\smalla{XttailThit} (t-1)}e^{-\frac 12 (1-2\delta )\left(2-\alpha-\varepsilon_2\right)k} \\
	& \qquad +\delta_{z} \bigC{qtxy} \cdot \bigC{XttailThit} e^{\delta |y|+\frac 12 \left(2-\alpha-\varepsilon_{2}\right)(y \vee 0)} e^{-\smalla{XttailThit} (t-1)}e^{-\frac 12 (1-2\delta )\left(2-\alpha-\varepsilon_2\right)\lfloor z \rfloor } \\
	& \qquad +\delta_{z} \bigC{qtxy} e^{-\smalla{qtxy} z^{2}} \cdot \bigC{Thit} e^{\delta |y|} e^{-\smalla{Thit} (t-1)}.
	\end{aligned}
	\]
	Then rearranging, it follows that
	\begin{align} \label{eq:Xt1in+int}
	&\mathbb{P}_{y, \pi}\left(X_{t}^{(1)} \in[z, z+\delta_z], T_{\text{hit}}>t\right) \notag \\
	& \leq \delta_{z} \bigC{qtxy} \bigC{XttailThit} e^{\delta |y|+\frac 12 \left(2-\alpha-\varepsilon_{2}\right)(y \vee 0)} e^{-\smalla{XttailThit} (t-1)} \notag \\
	& \qquad \cdot \bigg(\mathds 1_{z \geq 1} \sum_{k=0}^{\lfloor z \rfloor-1} 
	e^{-\smalla{qtxy} (\lfloor z \rfloor -1 -k)^{2}}
	e^{-\frac 12(1-2 \delta)\left(2-\alpha-\varepsilon_{2}\right)(\lfloor z \rfloor-1)} e^{\frac 12(1-2 \delta) \left(2-\alpha-\varepsilon_{2}\right)(\lfloor z \rfloor -1-k)}  \notag \\
	& \hspace{3cm} +e^{-\frac 12 (1-2 \delta)\left(2-\alpha-\varepsilon_{2}\right)\lfloor z \rfloor}\bigg)+\delta_{z} \bigC{qtxy} e^{-\smalla{qtxy} z^{2}} \cdot \bigC{Thit} e^{\delta |y|} e^{-\smalla{Thit} (t-1)} \notag \\
	& \leq \delta_{z} \bigC{qtxy} \bigC{XttailThit} e^{\delta |y|+\frac 12 \left(2-\alpha-\varepsilon_{2}\right)(y \vee 0)} e^{-\smalla{XttailThit} (t-1)}  e^{-\frac 12(1-2 \delta)\left(2-\alpha-\varepsilon_{2}\right)\lfloor z \rfloor} \notag \\
	& \qquad \cdot\left(\mathds{1}_{z \ge 1} \sum_{k^{\prime}=0}^{\lfloor z \rfloor-1} e^{-\smalla{qtxy} (k')^{2}} e^{\frac 12(1-2 \delta)\left(2-\alpha-\varepsilon_{2}\right) k^{\prime}} e^{\frac 12(1-2 \delta)\left(2-\alpha-\varepsilon_{2}\right)}+1\right) \notag \\
	& \hspace{5cm}+\delta_{z} \bigC{qtxy} e^{-\smalla{qtxy} z^{2}} \cdot \bigC{Thit} e^{\delta |y|} e^{-\smalla{Thit} (t-1)} \notag \\
	& \leq \delta_{z} K_{3} e^{\delta|y|+\frac 12\left(2-\alpha-\varepsilon_{2}\right)(y \vee 0)} e^{-(\smalla{qtxy}\wedge \smalla{Thit}) t} e^{-\frac 12(1-2 \delta)\left(2-\alpha-\varepsilon_{2}\right) z},
	\end{align}
	for some constant $K_3>0$.
	By the same argument, using~\eqref{eq:X1<xThit} in place of~\eqref{eq:X1>xThit},
	\begin{equation} \label{eq:Xt1in-int}
	\begin{aligned}
	&\mathbb{P}_{y,\pi}\left(X_{t}^{(1)} \in\left[-z-\delta_{z},-z\right], T_{\text{hit}}>t\right)\\
	&\quad \leq \delta_{z} K_4 e^{\delta |y|+\frac 12\left(2+\alpha-\varepsilon_{1}\right)((-y) \vee 0)} e^{-(\smalla{qtxy}\wedge \smalla{Thit}) t} e^{-\frac 12(1-2 \delta)\left(2+\alpha-\varepsilon_{1}\right) z},
	\end{aligned}
	\end{equation}
	for some constant $K_4>0$.
	Using~\eqref{eq:X2>xThit} and~\eqref{eq:X2<xThit} in place of~\eqref{eq:X1>xThit} and~\eqref{eq:X1<xThit},
	the same bounds as~\eqref{eq:Xt1in+int} and~\eqref{eq:Xt1in-int}
	hold with $X^{(1)}_t$ replaced with $X^{(2)}_t$.
	
	Therefore, by Lemma~\ref{lem:probXtI}, for $t \ge 1,$ $y \in \mathbb{R},$ $z>0$ and $\delta_{z}>0$,
	\[
	\begin{aligned}
	&\left| \mathbb{P}_{y}\left(X_{t} \in\left[z, z+\delta_{z}\right]\right)-\int_{\left[z, z+\delta_{z}\right]} \pi(x) d x \right| \\
	& \quad \leq \delta_{z} \cdot 2 K_{3} e^{\delta|y|+\frac 12 \left(2-\alpha-\varepsilon_{2}\right)(y \vee 0)} e^{-(\smalla{qtxy}\wedge \smalla{Thit}) t} e^{-\frac 12(1-2 \delta)\left(2-\alpha-\varepsilon_{2}\right) z},
	\end{aligned}
	\]
	and~\eqref{eq:qtxy_estimate1} follows.
	By the same argument,~\eqref{eq:qtxy_estimate2} also holds.
	\end{proof}

	\appendix

	\section{Fermi coordinates} \label{sec:fermi}
	
	In this section, we prove Lemma~\ref{lemma:eta}, \eqref{eta_Lipschitz} and Proposition~\ref{prop:D_eta}.
	
	\begin{proof}[Proof of Lemma~\ref{lemma:eta}]
		We adapt the proof of Proposition~3.3 in \cite{xu_interface_2025}, which used $ \| \cdot \|_\infty $ instead of $ \| \cdot \|_{2,\alpha} $.
		First note that, for all $ \eta \in \R $,
		\begin{equation*}
			\| m_\eta - m_{\eta_0} \|_{2,\alpha} \leq \| v - m_\eta \|_{2,\alpha} + \| v - m_{\eta_0} \|_{2,\alpha} \leq \| v - m_\eta \|_{2,\alpha} + \Beta{eta}.
		\end{equation*}
		But, for $ | \eta - \eta_0 | > 1 $,
		\begin{align*}
			\| m_{\eta} - m_{\eta_0} \|_{2,\alpha}^2 &= e^{\alpha \eta_0} \int_\R \left( m(x) - m(x - (\eta - \eta_0)) \right)^2 e^{\alpha x} dx \\
			&\geq e^{-\alpha K} \inf_{|y|>1}\int_\R \left( m(x) - m(x -y) \right)^2 e^{\alpha x} dx =: c_K,
		\end{align*}
		and $ c_K > 0 $.
		As a result, for $ \Beta{eta} < c_K/2 $, $ \| v - m_{\eta} \|_{2,\alpha} > \Beta{eta}, \text{ for all } | \eta - \eta_0 | > 1 $.
		In particular, for all $ | \eta - \eta_0 | > 1 $,
		\begin{equation*}
			\| v - m_{\eta} \|_{2,\alpha} > \| v - m_{\eta_0} \|_{2,\alpha} \geq \dist(v,M).
		\end{equation*}
		This already proves that there exists at least one $ \eta \in \R $ satisfying \eqref{dist_realised}, and that any such $ \eta $ lies in $ [\eta_0-1, \eta_0+1] $.
		Recalling \eqref{diff_norm}, any such $ \eta $ must also satisfy \eqref{s_orthogonal}.
		We then show that there exists at most one $ \eta \in [\eta_0-1, \eta_0+1] $ satisfying \eqref{s_orthogonal}.
		
		Let us set
		\begin{equation*}
			g(v, \eta) := \langle v - m_\eta, \partial_{x} m_\eta \rangle_\alpha e^{-\alpha \eta},
		\end{equation*}
		and note that $ g(v, \eta) = g_0(\eta) + g_1(v, \eta) $ where
		\begin{align*}
			g_0(\eta) := \langle m_{\eta_0} - m_\eta, \partial_{x} m_\eta \rangle_\alpha e^{-\alpha \eta}, && g_1(v, \eta) := \langle v - m_{\eta_0}, \partial_{x} m_\eta \rangle_\alpha e^{-\alpha \eta}.	
		\end{align*}
		Noting that, by a change of variables,
		\begin{equation*}
			g_0(\eta) = \langle m_{\eta_0 - \eta} - m, \partial_{x} m \rangle_\alpha =: \tilde{g}_0(\eta-\eta_0),
		\end{equation*}
		we obtain that, for all $ \eta \in \R $,
		\begin{equation*}
			g_0'(\eta) = \langle \partial_x m_{\eta_0 - \eta}, \partial_{x} m \rangle_\alpha > 0,
		\end{equation*}
		and $ g_0(\eta_0) = 0 $. 
		As a result, $ \eta_0 $ is the unique $ \eta \in \R $ such that $ g_0(\eta) = 0 $, and $ g_0(\eta) < 0 $ for $ \eta < \eta_0 $ and $ g_0(\eta) > 0 $ for $ \eta > \eta_0 $.
		On the other hand, note that
		\begin{equation*}
			\deriv*{g_1}{\eta}(v, \eta) = - \langle v - m_{\eta_0}, \partial_{xx} m_{\eta} + \alpha \partial_{x} m_\eta \rangle_\alpha e^{-\alpha \eta}.
		\end{equation*}
		By the Cauchy-Schwarz inequality, for $ \eta \in [\eta_0-1, \eta_0+1] $,
		\begin{align*}
			\left| \deriv*{g_1}{\eta}(v, \eta) \right| &\leq \| v - m_{\eta_0} \|_{2,\alpha} \| \partial_{xx} m_\eta + \alpha \partial_{x} m_\eta \|_{2,\alpha} e^{-\alpha \eta} \\
			&\leq \| v - m_{\eta_0} \|_{2,\alpha} \| \partial_{xx} m + \alpha \partial_{x} m \|_{2,\alpha} e^{\frac{\alpha}{2}(K+1)},
		\end{align*}
		using the fact that $ \eta \geq -K $.
		We can thus choose $ \Beta{eta} $ small enough that
		\begin{equation*}
			\| \partial_{xx} m + \alpha \partial_{x} m \|_{2,\alpha} e^{\frac{\alpha}{2}(K+1)} \Beta{eta} < \inf_{| \eta - \eta_0 | \leq 1} g_0'(\eta) = \inf_{| \eta | \leq 1} \tilde{g}_0'(\eta).
		\end{equation*}
		For such a choice of $ \Beta{eta} $, we obtain that there exists $ c > 0 $ (independent of $ v $) such that
		\begin{equation} \label{lower_bound_dgdeta}
			\deriv*{g}{\eta}(v, \eta) \geq c > 0 \text{ for all } \eta \in [\eta_0 - 1, \eta_0 + 1],
		\end{equation}
		which implies that there exists at most one $ \eta = \eta(v) \in [\eta_0-1, \eta_0+1] $ such that $ g(v, \eta) = 0 $.
		Combined with the previous argument, this proves the existence and uniqueness of $ \eta(v) \in \R $ satisfying \eqref{dist_realised}, and also that \eqref{s_orthogonal} holds.
		
		For the bound \eqref{bound_eta_eta_0}, note that, using the fact that $ g(v, \eta(v)) = 0 $ and $ g_0(\eta_0) = 0 $,
		\begin{equation*}
			c | \eta(v) - \eta_0 | \leq | g(v, \eta(v)) - g(v, \eta_0) | = | g(v, \eta_0) | = | g_1(v, \eta_0) |.
		\end{equation*}
		Then, by the Cauchy-Schwarz inequality, 
		\begin{align*}
			| g_1(v, \eta_0) | &\leq \| v - m_{\eta_0} \|_{2,\alpha} \| \partial_x m_{\eta_0} \|_{2,\alpha} e^{-\alpha \eta_0} \\
			&\leq \| v - m_{\eta_0} \|_{2,\alpha} \| \partial_{x} m \|_{2,\alpha} e^{\frac{\alpha}{2}K},
		\end{align*}
		yielding \eqref{bound_eta_eta_0}.
		
		It remains to prove \eqref{eta_Lipschitz}.
		For simplicity, let us set $ \eta_i = \eta(v_i) $.
		Using the fact that $ e^{\alpha \eta_i} g(v_i, \eta_i) = 0 $, we can write
		\begin{align} \label{diff_eta}
			\langle v_1 - v_2, \partial_{x} m_{\eta_1} \rangle_{ \alpha} - \langle m_{\eta_1} - m_{\eta_2}, \partial_{x} m_{\eta_1} \rangle_{\alpha} + \langle v_2 - m_{\eta_2}, \partial_{x} m_{\eta_1} - \partial_{x} m_{\eta_2} \rangle_{\alpha} = 0.
		\end{align}
		In addition, for $x\in \R$,
		\begin{align*}
			m_{\eta_1}(x)-m_{\eta_2}(x) = - (\eta_1-\eta_2) \int_{0}^{1} \partial_{x} m_{\eta_2}(x-t(\eta_1-\eta_2)) dt.
		\end{align*}
		Using a similar statement with $ m $ replaced by $ \partial_{x}m $, we obtain from~\eqref{diff_eta} that
		\begin{multline} \label{eq:eta1eta2}
			(\eta_1-\eta_2) \left( \int_{0}^{1} \langle \partial_{x} m_{\eta_2 + t (\eta_1-\eta_2)}, \partial_{x} m_{\eta_1} \rangle_{\alpha} dt - \int_{0}^{1} \langle \partial_{xx}m_{\eta_2 + t (\eta_1-\eta_2)}, v_2-m_{\eta_2} \rangle_{ \alpha} dt \right) \\ = \langle v_2-v_1, \partial_{x} m_{\eta_1} \rangle_{ \alpha}.
		\end{multline}
		By the Cauchy-Schwarz inequality,
		\begin{align*}
			\abs{\langle v_2-v_1, \partial_{x} m_{\eta_1} \rangle_{ \alpha}} &\leq \| v_2-v_1 \|_{2,\alpha} \|\partial_{x} m\|_{2,\alpha} e^{\frac{\alpha}{2} \eta_1} \\
			&\leq \| \partial_{x} m \|_{2,\alpha} e^{\frac{\alpha}{2} (K + \Cst{eta} \Beta{eta})} \| v_2 - v_1 \|_{2,\alpha}.
		\end{align*}
		In addition, the function
		\begin{align*}
			(t, \eta_1, \eta_2) \mapsto \langle \partial_{x} m_{\eta_2 + t (\eta_1-\eta_2)}, \partial_{x} m_{\eta_1} \rangle_{\alpha}
		\end{align*}
		is continuous and positive on $ [0,1] \times [-K-\Cst{eta}\Beta{eta},K+\Cst{eta} \Beta{eta}]^2 $ (since $ \partial_{x}m < 0 $).
		Hence there exists a constant $ c > 0 $ (depending on $ K $ and $ \Beta{eta} $) such that
		\begin{align*}
			\int_{0}^{1} \langle \partial_{x} m_{\eta_2 + t (\eta_1-\eta_2)}, \partial_{x} m_{\eta_1} \rangle_{\alpha} dt \geq c,
		\end{align*}
		for all $ \eta_1, \eta_2 \in [-K-\Cst{eta}\Beta{eta},K+\Cst{eta} \Beta{eta}]^2 $.
		Moreover, by the Cauchy-Schwarz inequality,
		\begin{align*}
			\abs{\int_{0}^{1} \langle \partial_{xx}m_{\eta_2 + t (\eta_1-\eta_2)}, v_2-m_{\eta_2} \rangle_{ \alpha} dt} &\leq \dist(v_2,M) \, \| \partial_{xx} m \|_{2,\alpha} \, e^{\frac{\alpha}{2} (\eta_1 \vee \eta_2)} \\
			&\leq \Beta{eta} \| \partial_{x x} m \|_{2,\alpha} e^{\frac{\alpha}{2} (K + \Cst{eta} \Beta{eta})}.
		\end{align*}
		Noting that $ c $ is non-increasing as a function of $ \Beta{eta} $, we can choose $ \Beta{eta} $ sufficiently small that
		\begin{align*}
			\Beta{eta} \| \partial_{xx}m \|_{2,\alpha} e^{\frac{\alpha}{2} (K+\Cst{eta}\Beta{eta})} < c.
		\end{align*}
		We then obtain the result from~\eqref{eq:eta1eta2}.
	\end{proof}
	
	\begin{proof}[Proof of Proposition~\ref{prop:D_eta}]
		The fact that $ \eta : \mathcal{V} \subset L^{2,\alpha} \to \R $ is twice Fréchet differentiable follows from the fact that, by \eqref{s_orthogonal},
		\begin{equation*}
			g(v,\eta(v)) = 0,
		\end{equation*}
		combined with \eqref{lower_bound_dgdeta}, by the implicit function theorem.
		To compute the Fréchet derivatives, it is more convenient to work with the function
		\begin{equation*}
			\bm{g}(v,\eta) := \langle v - m_\eta, \partial_{x} m_\eta \rangle_\alpha.
		\end{equation*}
		Since $ \bm{g}(v,\eta) = g(v,\eta) e^{\alpha \eta} $, we also have $ \bm{g}(v,\eta(v)) = 0 $ and there exists $ c' > 0 $ such that, for all $ v \in \mathcal{V}_{K, \Beta{eta}(K)} $,
		\begin{equation} \label{lower_bound_dhdeta}
			\deriv*{\bm{g}}{\eta}(v,\eta(v)) \geq c' > 0.
		\end{equation}
		Then, letting $ D_v $ denote the partial Fréchet derivative with respect to the variable $ v $, we write, for $ h \in L^{2,\alpha} $,
		\begin{equation*}
			D_v \bm{g}(v,\eta(v)) h + \deriv*{\bm{g}}{\eta}(v,\eta(v)) D\eta(v) h = 0,
		\end{equation*}
		which yields \eqref{Deta}, noting that
		\begin{equation*}
			\deriv*{\bm{g}}{\eta}(v,\eta(v)) = \| \partial_{x} m_{\eta(v)} \|_{2,\alpha}^2 - \langle v - m_{\eta(v)}, \partial_{xx} m_{\eta(v)} \rangle_\alpha = - \Vphi{1}(v).
		\end{equation*}
		We also obtain the lower bound \eqref{bound_varphi_tilde_below} from \eqref{lower_bound_dhdeta}.
		The second Fréchet derivative of $ \eta $ is obtained in a similar way, noting that, for $ h_1, h_2 \in L^{2,\alpha} $,
		\begin{multline*}
			D^2_v \bm{g}(v, \eta(v))(h_1, h_2) + \deriv*{}{\eta} D_v \bm{g}(v, \eta(v)) h_1 \, D\eta(v) h_2 + \deriv*{}{\eta} D_v \bm{g}(v, \eta(v)) h_2 \, D\eta(v) h_1 \\ + \deriv*[2]{\bm{g}}{\eta}(v,\eta(v)) D\eta(v) h_1 \, D\eta(v) h_2 + \deriv*{\bm{g}}{\eta}(v,\eta(v)) D^2 \eta(v)(h_1, h_2) = 0.
		\end{multline*}
		It then suffices to compute
		\begin{align*}
			D^2_v \bm{g}(v,\eta)(h_1, h_2) = 0, && \deriv*{}{\eta} D_v \bm{g}(v, \eta) h = - \langle v, \partial_{xx} m_{\eta} \rangle_\alpha,
		\end{align*}
		and
		\begin{equation*}
			\deriv*[2]{\bm{g}}{\eta}(v,\eta(v)) = -3 \langle \partial_{x} m_\eta, \partial_{xx} m_\eta \rangle_\alpha + \langle v - m_\eta, \partial_{x}^3 m_\eta \rangle_\alpha,
		\end{equation*}
		and conclude with the help of \eqref{dxxm_dxm}.
		The upper bound \eqref{bound_varphi_tilde} is then an easy consequence of the fact that $ \| v - m_{\eta(v)} \|_{2,\alpha} \leq \Beta{eta}(K) $ and $ | \eta(v) | \leq K + \Cst{eta} \Beta{eta} $.
		Using the fact that $ \eta $ is bounded on $ \mathcal{V}_{K} $, it is straightforward to show that
		\begin{equation} \label{bound_Vphi_v_v'}
			\abs{\Vphi{k}(v) - \Vphi{k}(v')} \leq C \left( \| v - v' \|_{2,\alpha} + | \eta(v) - \eta(v') | \right).
		\end{equation}
		The bound \eqref{varphi_Lipschitz} then follows from \eqref{eta_Lipschitz}.
		The fact that $ \eta : \mathcal{V}_K \to \R $ is Lipschitz continuous was already established in Lemma~\ref{lemma:eta}.
		In addition, for $ v, v' \in \mathcal{V}_K $ and $ h \in L^{2,\alpha} $,
		\begin{equation*}
			\abs{ D\eta(v) h - D\eta(v') h } \leq \frac{1}{| \Vphi{1}(v) |} | \langle \partial_{x} m_{\eta(v)} - \partial_{x} m_{\eta(v')}, h \rangle_\alpha | + \frac{| \Vphi{1}(v) - \Vphi{1}(v') |}{|\Vphi{1}(v)| |\Vphi{1}(v')|} |\langle \partial_{x} m_{\eta(v')}, h \rangle_\alpha |.
		\end{equation*}
		Using \eqref{bound_varphi_tilde}, \eqref{bound_varphi_tilde_below}, \eqref{varphi_Lipschitz} and \eqref{eta_Lipschitz} then yields the existence of $ C > 0 $ such that, for all $ v, v' \in \mathcal{V}_K $,
		\begin{equation*}
			\abs{ D\eta(v) h - D\eta(v') h } \leq C \| v - v' \|_{2,\alpha} \| h \|_{2,\alpha}.
		\end{equation*}
		The Lipschitz continuity of $ D^2 \eta $ is obtained in a similar way, using \eqref{bound_varphi_tilde}, \eqref{bound_varphi_tilde_below}, \eqref{varphi_Lipschitz} and \eqref{eta_Lipschitz}.
	\end{proof}

	\section{Fractional Sobolev space estimates} \label{sec:sobolev}

	Here, we prove Lemma~\ref{lemma:sobolev_average}.

	\begin{proof}[Proof of Lemma~\ref{lemma:sobolev_average}]
		We start by noting that
		\begin{equation*}
			\mathcal{J} (\rho^\delta \cdot g) \ast v = (\mathcal{J} (\rho^\delta \cdot g)) \ast (\mathcal{J} v).
		\end{equation*}
		As a result
		\begin{equation*}
			\mathcal{F} \mathcal{J} (\rho^\delta \cdot g) \ast v(\xi) = \mathcal{F} \mathcal{J} (\rho^\delta \cdot g)(\xi) \cdot \mathcal{F} \mathcal{J} v(\xi).
		\end{equation*}
		Hence,
		\begin{equation*}
			\Hnorm{\rho^\delta v - v}{\gamma} = \left( \int_\R (1 + |\xi|^\gamma)^2 | \mathcal{F} \mathcal{J} (\rho^\delta \cdot g)(\xi) - 1 |^2 | \mathcal{F} \mathcal{J} v(\xi) |^2 d\xi \right)^{1/2}.
		\end{equation*}
		By definition,
		\begin{align*}
			\mathcal{F} \mathcal{J} (\rho^\delta \cdot g)(\xi) &= \int_\R \rho^\delta(x) g(x) e^{\frac{\alpha}{2} x -i \xi x} dx \\ 
			&= \int_{[-1,1]} \rho(z) g(\delta z) e^{\frac{\alpha}{2} \delta z - i \xi \delta z} dz.
		\end{align*}
		Since  $ \int_{[-1,1]} \rho(z) dz = 1 $, it follows that
		\begin{equation*}
			| \mathcal{F} \mathcal{J} (\rho^\delta \cdot g)(\xi) - 1 | \leq \int_{[-1,1]} \rho(z) | g(\delta z) e^{\frac{\alpha}{2} \delta z - i \xi \delta z} - 1 | dz.
		\end{equation*}
		Then, for $ z \in [-1,1] $,
		\begin{align*}
			| g(\delta z ) e^{\frac{\alpha}{2} \delta z - i \xi \delta z} - 1 | &\leq | g(\delta z) - 1 | e^{\frac{\alpha}{2} \delta} + | e^{\frac{\alpha}{2} \delta z} - 1 | + | e^{- i \xi \delta z} - 1 | \\
			&\leq C ( \delta + ((|\xi| \delta) \wedge 1 )),
		\end{align*}
		for some $ C > 0 $, since $g$ is Lipschitz continuous and $g(0)=1$.
		Using the fact that, for any $ \gamma' \in (0,1] $, $ ((|\xi| \delta) \wedge 1) \leq (|\xi| \delta)^{\gamma'} $, we obtain that
		\begin{equation*}
			\Hnorm{(\rho^\delta \cdot g) \ast v - v}{\gamma} \leq C \delta \Hnorm{v}{\gamma} + C \delta^{\gamma'} \Hnorm{v}{\gamma+\gamma'},
		\end{equation*}
		from which the conclusion follows.
	\end{proof}

	\section{Estimates on the deterministic flow} \label{sec:Phit}
	Here, we prove Proposition~\ref{prop:Phit}. 
	First recall the definition of the semigroup $ (Q(t), t \geq 0) $ in \eqref{eq:Qdefn}, and note that
	\begin{equation} \label{Phi_t_Qt}
		\Phi_t(v,\cdot) = Q(t)v + \int_{0}^{t} Q(t-s) f(\Phi_s(v,\cdot)) ds.
	\end{equation}
	We then prove the following.
	
	\begin{lemma} \label{lemma:Qt}
		For any $ \gamma, \gamma' \in \R $ with $ \gamma \geq \gamma' $, there exists a constant $ C_{\gamma, \gamma'} > 0 $ such that, for all $ t \geq 0 $ and $ v \in H^{\gamma',\alpha} $,
		\begin{equation} \label{Qt_Hgamma}
			\Hnorm{Q(t) v}{\gamma} \leq \frac{C_{\gamma,\gamma'}}{(t\wedge 1)^{\frac{\gamma-\gamma'}{2}}} e^{-\frac{\alpha^2}{4} t} \Hnorm{v}{\gamma'}.
		\end{equation}
		Moreover, for any $ \gamma \in \R $, we can take $ C_{\gamma,\gamma} = 1 $. 
		In addition, for any $ \gamma \geq 0 $ and $ \gamma' \in (0,2] $, there exists $ C > 0 $ such that, for any $ v \in H^{\gamma+\gamma',\alpha} $,
		\begin{equation} \label{continuity_Qt}
			\Hnorm{Q(t)v - v}{\gamma} \leq C t^{\gamma'/2} \Hnorm{v}{\gamma+\gamma'}.
		\end{equation}
	\end{lemma}
	
	\begin{proof}
		Recall the definition of $ \Hnorm{\cdot}{\gamma} $ in \eqref{def:Hnorm_gamma}.
		By completing the square in the exponential, we note that
		\begin{equation*}
			\mathcal{J} Q(t) v = e^{-\frac{\alpha^2}{4} t} G_t \ast (\mathcal{J} v).
		\end{equation*}
		It follows that
		\begin{equation*}
			\mathcal{F} \mathcal{J} Q(t) v(\xi) = e^{-(\frac{\alpha^2}{4} + \xi^2) t} \mathcal{F}\mathcal{J}v(\xi).
		\end{equation*}
		As a result,
		\begin{equation} \label{Hnorm_Qtv}
			\Hnorm{Q(t)v}{\gamma}^2 = \int_\R (1 + |\xi|^\gamma)^2 e^{-2(\frac{\alpha^2}{4} + \xi^2) t} | \mathcal{F}\mathcal{J}v(\xi) |^2 d\xi.
		\end{equation}
		Then, note that, for any $ \gamma \geq \gamma' $, there exists $ C_{\gamma,\gamma'} > 0 $ such that, for all $ t \geq 0 $ and $ \xi \in \R $,
		\begin{align} \label{bound_xi_gamma}
			(1+|\xi|^{\gamma}) e^{-t \xi^2} \leq (1 + |\xi|^{\gamma'}) \frac{C_{\gamma,\gamma'}}{(t\wedge 1)^{\frac{\gamma -\gamma'}{2}}}.
		\end{align}
		Plugging this in \eqref{Hnorm_Qtv} yields \eqref{Qt_Hgamma}.
		In addition, we note that
		\begin{equation*}
			\Hnorm{Q(t)v - v}{\gamma}^2 = \int_\R (1 + |\xi|^\gamma)^2 \left( 1 - e^{-(\frac{\alpha^2}{4} + \xi^2)t} \right)^2 | \mathcal{F} \mathcal{J} v(\xi) |^2 d\xi.
		\end{equation*}
		Then, since $ \gamma' \in (0,2] $, we have
		\begin{equation*}
			\left( 1 - e^{-(\frac{\alpha^2}{4} + \xi^2)t} \right)^2 \leq \left( \left( \frac{\alpha^2}{4} + \xi^2 \right) t \right)^{\gamma'}.
		\end{equation*}
		It follows that there exists $ C> 0 $ such that
		\begin{equation*}
			\Hnorm{Q(t)v - v}{\gamma}^2 \leq C t^{\gamma'} \int_\R (1+ |\xi|^{\gamma+\gamma'})^2 |\mathcal{F}\mathcal{J}v(\xi)|^2 d\xi \leq C t^{\gamma'} \Hnorm{v}{\gamma+\gamma'}^2,
		\end{equation*}
		yielding \eqref{continuity_Qt}.
	\end{proof}
	
	We are now ready to prove Proposition~\ref{prop:Phit}.
	
	\begin{proof}[Proof of Proposition~\ref{prop:Phit}]
		We start by showing that $ t \mapsto \| \Phi_t(v,\cdot) \|_{2,\alpha} $ is locally bounded.
		Using this and applying Gronwall's inequality, we shall prove the bounds \eqref{bounds_Phit_Hgamma} first for $ \gamma = 0 $ and then for $ \gamma \in (0,2) $, as well as \eqref{continuity_Phit_Hgamma}.
		
		Let us start by showing that $ t \mapsto \| \Phi_t(v,\cdot) \|_{2,\alpha} $ is locally bounded.
		By the Feynman-Kac formula, setting $ r(u) = \frac{f(u)}{u} $,
		\begin{equation*}
			\Phi_t(v,x) = \E[x]{ \exp\left( \int_{0}^{t} r(\Phi_{t-s}(v,X_s)) ds \right) v(X_t) },
		\end{equation*}
		where $ X_t = x + \alpha t + \sqrt{2} B_t $ and $ (B_t, t \geq 0) $ is standard Brownian motion.
		Since $ r $ is bounded, for any $ T > 0 $ there exists a constant $ \genericc > 0 $ such that
		\begin{equation} \label{Phi_t_Feynman_Kac}
			| \Phi_t(v,x) | \leq e^{\genericc t} \E[x]{ |v(X_t)| } = e^{\genericc t} Q(t) |v|(x).
		\end{equation}
		By \eqref{Qt_Hgamma} (with $ \gamma = \gamma' = 0 $), we obtain
		\begin{align*}
			\| \Phi_t(v,\cdot) \|_{2,\alpha} \leq e^{(\genericc - \frac{\alpha^2}{4}) t} \| v \|_{2,\alpha}.
		\end{align*}
		Hence for any $ T > 0 $, there exists $ C_T > 0 $ such that
		\begin{equation} \label{bound_Phit_L2a}
			\| \Phi_t(v,\cdot) \|_{2,\alpha} \leq C_T \| v \|_{2,\alpha}.
		\end{equation}
	
		We now return to \eqref{Phi_t_Qt} and apply \eqref{Qt_Hgamma} with $ \gamma' = 0 $, yielding, for $ t > 0 $,
		\begin{equation*}
			\Hnorm{\Phi_t(v,\cdot)}{\gamma} \leq \frac{C_\gamma}{(t\wedge 1)^{\gamma/2}} \| v \|_{2,\alpha} + \int_{0}^{t} \frac{C_\gamma}{(|t-s|\wedge 1)^{\gamma/2}} \| f(\Phi_s(v,\cdot)) \|_{2,\alpha} ds.
		\end{equation*}
		Then, noting that $ |f(u)| \leq C_f |u| $ for some $ C_f > 0 $ and using \eqref{bound_Phit_L2a}, we obtain, for $ t \in (0,T] $,
		\begin{equation*}
			\Hnorm{\Phi_t(v,\cdot)}{\gamma} \leq \frac{C_\gamma}{(t\wedge 1)^{\gamma/2}} \| v \|_{2,\alpha} + \int_{0}^{t} \frac{C_\gamma C_f C_T}{(|t-s|\wedge 1)^{\gamma/2}} ds \| v \|_{2,\alpha},
		\end{equation*}
		yielding the first part of \eqref{bounds_Phit_Hgamma} (the integral converges for $ \gamma < 2 $).
		
		The proof of the second part of \eqref{bounds_Phit_Hgamma} is similar, combining \eqref{Phi_t_Qt} and \eqref{Qt_Hgamma} (with $ \gamma' = 0 $) to obtain
		\begin{multline} \label{Phi_t_v_v'}
			\Hnorm{\Phi_t(v,\cdot) - \Phi_t(v',\cdot)}{\gamma} \leq \frac{C_\gamma}{(t\wedge 1)^{\gamma/2}} \| v - v' \|_{2,\alpha} \\ + \int_{0}^{t} \frac{C_\gamma}{(|t-s|\wedge 1)^{\gamma/2}} \| f(\Phi_s(v,\cdot)) - f(\Phi_s(v',\cdot)) \|_{2,\alpha} ds.
		\end{multline}
		We then note that there exists $ C_f' > 0 $ such that $ |f(u) - f(u')| \leq C_f' |u-u'| $.
		Taking $ \gamma = 0 $ above and applying Gronwall's inequality, we obtain that, for all $ t \geq 0 $,
		\begin{equation*}
			\| \Phi_t(v,\cdot) - \Phi_t(v',\cdot) \|_{2,\alpha} \leq C_0 \| v - v' \|_{2,\alpha} e^{C_0 C_f' t}.
		\end{equation*}
		Plugging this bound in the second term on the right-hand side of \eqref{Phi_t_v_v'} for any $ \gamma \in [0,2) $ then yields the second part of \eqref{bounds_Phit_Hgamma}.
		
		Finally, from \eqref{Phi_t_Qt} and \eqref{Qt_Hgamma} (with $ \gamma' = 0 $), we write
		\begin{equation*}
			\Hnorm{\Phi_t(v,\cdot) - v}{\gamma'} \leq \Hnorm{Q(t)v - v}{\gamma'} + \int_{0}^{t} \frac{C_{\gamma'}}{(|t-s|\wedge 1)^{\gamma'/2}} C_f \| \Phi_s(v,\cdot) \|_{2,\alpha} ds.
		\end{equation*}
		We then obtain \eqref{continuity_Phit_Hgamma} from \eqref{continuity_Qt} and \eqref{bound_Phit_L2a}, since $ \gamma' < \gamma < 2 $.
	\end{proof}
	
	
	\begin{lemma} \label{lemma:Q_Lpalpha}
		For any $ p \geq 1 $ $ r \geq 1 $, $ q \in \R $ and $ T > 0 $, there exists $ C > 0 $ such that, for all $ t \in [0,T] $ and $ h \in L^{p,q/r} $,
		\begin{equation*}
			\| Q(t) h \|_{rp,q} \leq \frac{C}{t^{\frac{1}{2p}\left( 1-\frac{1}{r} \right)}} \| h \|_{p,q/r}.
		\end{equation*}
	\end{lemma}
	
	\begin{proof}
		By the definition of $ Q(t) $,
		\begin{equation} \label{expr_Qth}
			\| Q(t) h \|_{rp,q} = e^{-\frac{\alpha^2}{4}t} \left( \int_\R \left( \int_\R G_t(x-y) h(x) e^{\frac{\alpha}{2}x} dx \right)^{rp} e^{(q-\frac{\alpha}{2} rp) y} dy \right)^{1/rp}.
		\end{equation}
		By the H\"older inequality, since $ p \geq 1 $, we have
		\begin{multline*}
			\left( \int_\R G_t(x-y) h(x) e^{\frac{\alpha}{2}x} dx \right)^{rp} \\ \leq \left( \int_\R G_t(x-y) e^{(\frac{\alpha}{2} - \frac{q}{rp}) x} dx \right)^{r(p-1)} \left( \int_\R G_t(x-y) e^{(\frac{\alpha}{2} - \frac{q}{rp}) x} |h(x)|^p e^{\frac{q}{r} x} dx \right)^{r}.
		\end{multline*}
		Hence there exists a constant $ C > 0 $ such that, for $ t \in [0,T] $,
		\begin{equation} \label{Holder1}
			\left( \int_\R G_t(x-y) h(x) e^{\frac{\alpha}{2}x} dx \right)^{rp} \leq C e^{(\frac{\alpha}{2}r - \frac{q}{p})(p-1) y} \left( \int_\R G_t(x-y) e^{(\frac{\alpha}{2} - \frac{q}{rp}) x} |h(x)|^p e^{\frac{q}{r} x} dx \right)^{r}.
		\end{equation}
		Using the H\"older inequality again, since $ r \geq 1 $,
		\begin{multline*}
			\left( \int_\R G_t(x-y) e^{(\frac{\alpha}{2} - \frac{q}{rp}) x} |h(x)|^p e^{\frac{q}{r} x} dx \right)^{r} \\  \begin{aligned}
				&\leq \left( \int_\R |h(x)|^p e^{\frac{q}{r} x} dx \right)^{r-1} \int_\R G_t(x-y)^r e^{(\frac{\alpha}{2} - \frac{q}{rp})rx} |h(x)|^p e^{\frac{q}{r} x} dx \\
				&\leq \frac{C}{t^{\frac{r-1}{2}}} \| h \|_{p,q/r}^{(r-1)p} \int_\R G_t(x-y) e^{(\frac{\alpha}{2} - \frac{q}{rp})rx} |h(x)|^p e^{\frac{q}{r} x} dx,
			\end{aligned}
		\end{multline*}
		for some $ C > 0 $.
		Plugging this in \eqref{Holder1} and then in \eqref{expr_Qth}, we obtain that there exists $ C > 0 $ such that
		\begin{equation*}
			\| Q(t) h \|_{rp,q} \leq \frac{C}{t^{\frac{1}{2p}(1-\frac{1}{r})}} \| h \|_{p,q/r}^{1-\frac{1}{r}} \left( \int_{\R^2} G_t(x-y) e^{(\frac{\alpha}{2} - \frac{q}{rp})rx} |h(x)|^p e^{\frac{q}{r} x} e^{(\frac{q}{p} - \frac{\alpha}{2} r) y} dx dy \right)^{1/rp}.
		\end{equation*}
		The result then follows by integrating with respect to $ y $ in the integral on the right-hand side.
	\end{proof}

	\begin{cor} \label{cor:Phit_Lpalpha}
		For any $ T > 0 $, $ p \geq 1 $, $ r \geq 1 $ and $ q \in \R $, there exists $ C > 0 $ such that, for any $ v, v' \in L^{p,q/r} $,
		\begin{equation*}
			\| \Phi_t(v,\cdot) - \Phi_t(v',\cdot) \|_{rp,q} \leq \frac{C}{t^{\frac{1}{2p}\left( 1-\frac{1}{r} \right)}} \| v - v' \|_{p,q/r}.
		\end{equation*}
	\end{cor}

	\begin{proof}
		First note that, by \eqref{Phi_t_Feynman_Kac} and Lemma~\ref{lemma:Q_Lpalpha}, for $ t \in (0,T] $,
		\begin{equation*}
			\| \Phi_t(v,\cdot) - \Phi_t(v',\cdot) \|_{rp,q} \leq \frac{C}{t^{\frac{1}{2p}\left( 1-\frac{1}{r} \right)}} (\| v \|_{p,q/r} + \| v' \|_{p,q/r}),
		\end{equation*}
		Hence $ t \mapsto \| \Phi_t(v,\cdot) - \Phi_t(v',\cdot) \|_{rp,q} $ is integrable on $ [0,T] $.
		Combining \eqref{Phi_t_Qt} and Lemma~\ref{lemma:Q_Lpalpha}, there exist $ C > 0 $ and $ C' > 0 $ such that, for $ t \in (0,T] $,
		\begin{equation*}
			\| \Phi_t(v,\cdot) - \Phi_t(v',\cdot) \|_{rp,q} \leq \frac{C}{t^{\frac{1}{2p}\left( 1-\frac{1}{r} \right)}} \| v - v' \|_{p, q / r} + \int_{0}^{t} C' \| f' \|_\infty \| \Phi_s(v,\cdot) - \Phi_s(v',\cdot) \|_{rp,q}.
		\end{equation*}
		By Gronwall's inequality, we obtain, for $ t \in (0,T] $,
		\begin{equation*}
			\| \Phi_t(v,\cdot) - \Phi_t(v',\cdot) \|_{rp,q} \leq \left( \frac{C}{t^{\frac{1}{2p}\left( 1-\frac{1}{r} \right)}} + \int_{0}^{t} \frac{C C' \| f' \|_\infty}{s^{\frac{1}{2p}\left( 1-\frac{1}{r} \right)}} e^{C' \| f'\|_\infty (t-s)} ds \right) \| v - v' \|_{p,q/r},
		\end{equation*}
		from which the result follows.
	\end{proof}
	
	\begin{proof}[Proof of Proposition~\ref{prop:deriv_Phi}]
		By Definition~\ref{def:Frechet_derivative}, we need to show that $ T_{0,t,v}^* : H \to L^{2,\alpha} $ is continuous, that the bilinear map $ B_v : H \times H \to L^{2,\alpha} $ defined by
		\begin{equation*}
			B_v(h,h') := \int_{0}^{t} T_{s,t,v}^* \lbrace f''(\Phi_s(v,\cdot)) (T_{0,s,v}^* h)(T_{0,s,v}^* h') \rbrace ds,
		\end{equation*}
		is also continuous and that
		\begin{equation} \label{littleO_Phit}
			\left\| \Phi_t(v+h,\cdot) - \Phi_t(v,\cdot) - T_{0,t,v}^* h - \frac{1}{2} B_v(h,h) \right\|_{2,\alpha} = \littleO{\| h \|_H^2} \quad \text{ as } \| h \|_{H} \to 0.
		\end{equation}
	
		First, By Lemma~\ref{lemma:bound_pst}, for any $ T > 0 $, there exists $ C > 0 $ such that, for all $ 0 \leq s \leq t $ with $ t-s \leq T $
		\begin{equation} \label{bound_T_Q}
			| T_{s,t,v}^* h | \leq C Q(t-s) |h|.
		\end{equation}
		Combining this with \eqref{Qt_Hgamma} with $ \gamma = \gamma' = 0 $ (or equivalently Lemma~\ref{lemma:Q_Lpalpha} with $ p = 2 $, $ r = 1 $, $ q = \alpha $), we obtain that, for any $ T > 0 $, there exists $ C > 0 $ such that, for all $ 0 \leq s \leq t \leq T $,
		\begin{equation} \label{bound_Tstv_L2alpha}
			\| T_{s,t,v}^* h \|_{2,\alpha} \leq C \| h \|_{2,\alpha} \leq C \| h \|_H.
		\end{equation}
		This proves that $ T_{s,t,v}^* : H \to L^{2,\alpha} $ is continuous.
		From \eqref{bound_Tstv_L2alpha}, we also deduce that there exists $ C > 0 $ such that, for all $ 0 \leq s \leq t \leq T $,
		\begin{equation*}
			\| T_{s,t,v}^* \lbrace f''(\Phi_s(v,\cdot)) (T_{0,s,v}^* h)^2 \rbrace \|_{2,\alpha} \leq C \| T_{0,s,v}^* h \|_{4,\alpha}^2.
		\end{equation*}
		Moreover, by \eqref{bound_T_Q} and Lemma~\ref{lemma:Q_Lpalpha},
		\begin{equation*}
			\| T_{0,s,v}^* h \|_{4,\alpha} \leq C \| Q(s) |h| \|_{4,\alpha} \leq \frac{C'}{s^{1/8}} \| h \|_{2,\frac{\alpha}{2}}.
		\end{equation*}
		As a result, there exists $ C > 0 $ such that, for $ 0 < t \leq T $,
		\begin{equation*}
			\left\| \int_{0}^{t} T_{s,t,v}^* \lbrace f''(\Phi_s(v,\cdot)) (T_{0,s,v}^* h)^2 \rbrace ds \right\|_{2,\alpha} \leq C \| h \|_{2,\frac{\alpha}{2}}^2.
		\end{equation*}
		Since $ \| \cdot \|_{2,\frac{\alpha}{2}} \leq \| \cdot \|_H $, this shows that the bilinear map $ B_v : H \times H \to L^{2,\alpha} $ is continuous.

		It remains to show \eqref{littleO_Phit}.
		Fix $v \in H$ and $h \in H$ and set
		\begin{equation*}
			u_t(x) := \Phi_t(v + h,x) - \Phi_t(v,x).
		\end{equation*}
		Then $u(0,x) = h(x)$ and
		\begin{equation} \label{deriv_u}
			\partial_t u_t = (\partial_{xx} + \alpha \partial_{x}) u_t + f(\Phi_t(v+h,\cdot)) - f(\Phi_t(v,\cdot)).
		\end{equation}
		By Taylor's theorem, there exists a bounded function $ r_2 : \R^2 \to \R $ such that, for $ (u,v) \in \R^2 $,
		\begin{equation*}
			f(u) = f(v) + f'(v) (u-v) + \frac{1}{2} f''(v) (u-v)^2 + (u-v)^3 r_2(u,v).
		\end{equation*}
		We can then rewrite \eqref{deriv_u} as
		\begin{equation*}
			\partial_t u_t = (\partial_{xx} + \alpha \partial_{x} + f'(\Phi_t(v,\cdot))) u_t + \frac{1}{2} f''(\Phi_t(v,\cdot)) (u_t)^2 + (u_t)^3 r(\Phi_t(v+h,\cdot), \Phi_t(v,\cdot)).
		\end{equation*}
		It follows that
		\begin{equation} \label{eq_u}
			u_t = T_{0,t,v}^* h + \int_{0}^{t} T_{s,t,v}^* \left\lbrace \frac{1}{2} f''(\Phi_s(v,\cdot)) (u_s)^2 + (u_s)^3 r(\Phi_s(v+h,\cdot),\Phi_s(v,\cdot)) \right\rbrace ds.
		\end{equation}
		It follows that
		\begin{equation*}
			u_t - T_{0,t,v}^* h - B_v(h,h) = R_1(t) + R_2(t),
		\end{equation*}
		where
		\begin{equation*}
			R_1(t) := \frac{1}{2} \int_{0}^{t} T_{s,t,v}^* \left\lbrace f''(\Phi_s(v,\cdot)) ((u_s)^2 - (T_{0,s,v}^*h)^2) \right\rbrace ds,
		\end{equation*}
		and
		\begin{equation*}
			R_2(t) := \int_{0}^{t} T_{s,t,v}^* \left\lbrace (u_s)^3 r(\Phi_s(v+h,\cdot),\Phi_s(v,\cdot)) \right\rbrace ds.
		\end{equation*}
		By \eqref{bound_Tstv_L2alpha}, there exists $ C > 0 $ such that
		\begin{align*}
			\| R_1(t) \|_{2,\alpha} &\leq C \int_{0}^{t} \| (u_s)^2 - (T_{0,s,v}^* h)^2 \|_{2,\alpha} ds \\
			&\leq C \int_{0}^{t} \| u_s - T_{0,s,v}^* h \|_{3,\alpha} \left( \| u_s \|_{6,\alpha} + \| T_{0,s,v}^* h \|_{6,\alpha} \right) ds, \numberthis \label{bound_R1}
		\end{align*}
		using H\"older's and Minkowski's inequalities in the second line.
		Combining \eqref{bound_T_Q} and Lemma~\ref{lemma:Q_Lpalpha}, we obtain that
		\begin{equation} \label{bound_T_L4a}
			\| T_{0,s,v}^* h \|_{6,\alpha} \leq \frac{C}{s^{1/6}} \| h \|_{2,\frac{\alpha}{3}}.
		\end{equation}
		To bound $ \| u_s \|_{6,\alpha} $, we use the fact that, by Taylor's theorem, there exists a bounded function $ r_0 : \R^2 \to \R $ such that
		\begin{equation*}
			f(u) = f(v) + (u-v) r_0(u,v).
		\end{equation*}
		By the Feynman-Kac formula, using the fact that $ r_0 $ is bounded, we obtain that there exists $ C > 0 $ such that, for $ t \in [0,T] $,
		\begin{equation*}
			|u_t| \leq C Q(t)|h|.
		\end{equation*}
		As a result, by Lemma~\ref{lemma:Q_Lpalpha}, for any $ r \geq 1 $, there exists $ C > 0 $ such that, for all $ t \in [0,T] $,
		\begin{equation} \label{bound_ut}
			\| u_t \|_{2r,\alpha} \leq \frac{C}{t^{\frac{1}{4}-\frac{1}{4r}}} \| h \|_{2,\frac{\alpha}{r}}.
		\end{equation}
		Let us now bound $ \| u_s - T_{0,s,v}^* h \|_{3,\alpha} $.
		By Taylor's theorem, there exists a bounded function $ r_1 : \R^2 \to \R $ such that, for $ (u,v) \in \R^2 $,
		\begin{equation*}
			f(u) = f(v) + f'(v) (u-v) + (u-v)^2 r_1(u,v).
		\end{equation*}
		We then write \eqref{deriv_u} as
		\begin{equation*}
			\partial_t u_t = (\partial_{xx} + \alpha \partial_{x} + f'(\Phi_t(v,\cdot))) u_t + (u_t)^2 r_1(\Phi_t(v+h,\cdot), \Phi_t(v,\cdot)).
		\end{equation*}
		As a consequence,
		\begin{equation*}
			u_t = T_{0,t,v}^* h + \int_{0}^{t} T_{s,t,v}^* \lbrace (u_s)^2 r_1(\Phi_s(v+h,\cdot), \Phi_s(v,\cdot)) \rbrace ds.
		\end{equation*}
		Hence, using \eqref{bound_T_Q}, Lemma~\ref{lemma:Qt} and \eqref{bound_ut} in the second line,
		\begin{align*}
			\| u_t - T_{0,t,v}^* h \|_{3,\alpha} &\leq \int_{0}^{t} C \| u_s \|_{6,\alpha}^2 ds \\
			&\leq \int_{0}^{t} \frac{C}{s^{1/3}} \| h \|_{2,\frac{\alpha}{3}}^2 ds \\
			&\leq C \| h \|_{2,\frac{\alpha}{3}}^2, \numberthis \label{bound_ut-Th}
		\end{align*}
		where the constant $ C $ may change from line to line.
		Plugging \eqref{bound_T_L4a}, \eqref{bound_ut} and \eqref{bound_ut-Th} in \eqref{bound_R1}, we obtain
		\begin{equation*}
			\| R_1(t) \|_{2,\alpha} \leq C \| h \|_{2,\frac{\alpha}{3}}^3 \leq C \| h \|_H^3.
		\end{equation*}
		Turning to $ R_2(t) $, by \eqref{bound_Tstv_L2alpha}, there exists $ C > 0 $ such that
		\begin{equation*}
			\| R_2(t) \|_{2,\alpha} \leq C \int_{0}^{t} \| u_s \|_{6,\alpha}^3 ds.
		\end{equation*}
		Using now \eqref{bound_ut}, we obtain that there exists $ C > 0 $ such that
		\begin{equation*}
			\| R_2(t) \|_{2,\alpha} \leq C \| h \|_{2,\frac{\alpha}{3}}^3 \leq C \| h \|_H^3.
		\end{equation*}
		We have thus proved that there exists $ C > 0 $ (which happens to be independent of $ v $) such that, for all $ t \in [0,T] $,
		\begin{equation*}
			\left\| \Phi_t(v+h,\cdot) - \Phi_t(v,\cdot) - T_{0,t,v}^* h - \frac{1}{2} B_v(h,h) \right\|_{2,\alpha} \leq C \| h \|_H^3,
		\end{equation*}
		which proves \eqref{littleO_Phit}.
		
		The fact that $ \Phi_t : H \to L^{2,\alpha} $ is Lipschitz continuous directly follows from \eqref{bounds_Phit_Hgamma} in Proposition~\ref{prop:Phit} with $ \gamma = 0 $ and the fact that $ \| \cdot \|_{2,\alpha} \leq \| \cdot \|_{H} $.
		Recall that, by definition, for $ h \in H $,
		\begin{equation*}
			\partial_t T_{s,t,v}^* h = (\partial_{x x} + \alpha \partial_{x} + f'(\Phi_t(v,\cdot))) T_{s,t,v}^* h.
		\end{equation*}
		Hence, by the variation of constants formula,
		\begin{equation*}
			T_{s,t,v}^* h = Q(t-s) h + \int_{s}^{t} Q(t-r) \left( f'(\Phi_r(v,\cdot)) T_{s,r,v}^* h \right) dr.
		\end{equation*}
		It follows that
		\begin{multline*}
			(T_{s,t,v}^* - T_{s,t,v'}^*) h = \int_{s}^{t} Q(t-r) \left( \left( f'(\Phi_r(v,\cdot)) - f'(\Phi_r(v,\cdot)) \right) T_{s,r,v}^* h \right) dr \\ + \int_{s}^{t} Q(t-r) \left( f'(\Phi_r(v',\cdot)) \left( T_{s,r,v}^* - T_{s,r,v'}^* \right) h \right) dr.
		\end{multline*}
		Then, for any $ p \geq 2 $, by Lemma~\ref{lemma:Q_Lpalpha} with $ r = 1 $ and $ q = \alpha $,
		\begin{equation} \label{Tstv*-Tstv'*_pre_Gronwall}
			\| (T_{s,t,v}^* - T_{s,t,v'}^*) h \|_{p,\alpha} \leq C \int_{s}^{t} \left( \| (\Phi_r(v,\cdot) - \Phi_r(v',\cdot)) T_{s,r,v}^* h \|_{p,\alpha} + \| (T_{s,r,v}^* - T_{s,r,v'}^*) h \|_{p,\alpha} \right) ds.
		\end{equation}
		In addition, by the H\"older inequality, for any $ q, q' \geq 1 $ with $ \frac{1}{q} + \frac{1}{q'} = 1 $,
		\begin{equation*}
			\| (\Phi_r(v,\cdot) - \Phi_r(v',\cdot)) T_{s,r,v}^* h \|_{p,\alpha} \leq \| \Phi_r(v,\cdot) - \Phi_r(v',\cdot) \|_{pq,\alpha} \| T_{s,r,v}^* h \|_{pq',\alpha}.
		\end{equation*}
		Then, since $ pq \geq p \geq 2 $, we can use Corollary~\ref{cor:Phit_Lpalpha} and \eqref{bound_T_Q} combined with Lemma~\ref{lemma:Q_Lpalpha} to obtain, for $ r \in (s,T] $,
		\begin{equation*}
			\| (\Phi_r(v,\cdot) - \Phi_r(v',\cdot)) T_{s,r,v}^* h \|_{p,\alpha} \leq \frac{C}{r^{\frac{1}{4}\left( 1-\frac{2}{pq} \right)} (r-s)^{\frac{1}{4}\left( 1-\frac{2}{pq'} \right)}} \| v - v' \|_{2,\frac{2\alpha}{pq}} \| h \|_{2,\frac{2\alpha}{pq'}}.
		\end{equation*}
		Plugging this in \eqref{Tstv*-Tstv'*_pre_Gronwall} and using Gronwall's inequality, we obtain, for $ 0 \leq s \leq t \leq T $,
		\begin{equation} \label{Tstv*-Tstv'*}
			\| (T_{s,t,v}^* - T_{s,t,v'}^*) h \|_{p,\alpha} \leq C \| v- v' \|_{2,\frac{2\alpha}{pq}} \| h \|_{2,\frac{2\alpha}{pq'}}.
		\end{equation}
		Taking $ p = q = q' = 2 $ and using the fact that $ \| \cdot \|_{2,\frac{\alpha}{2}} \leq \| \cdot \|_H $, we obtain
		\begin{equation*}
			\| D \Phi_t(v) h - D\Phi_t(v') h \|_{2,\alpha} \leq C \| v - v' \|_H \| h \|_H,
		\end{equation*}
		which proves that $ v \mapsto D\Phi_t(v) $ is Lipschitz continuous.
		We now turn to $ D^2 \Phi_t $.
		For $ v, v', h \in H $, we write
		\begin{multline} \label{diff_Bv_Bv'}
			B_v(h,h) - B_{v'}(h,h) = \int_{0}^{t} (T_{s,t,v}^* - T_{s,t,v'}^*) \left\lbrace f''(\Phi_s(v,\cdot)) (T_{0,s,v}^* h)^2 \right\rbrace ds \\ + \int_{0}^{t} T_{s,t,v'}^* \left\lbrace \left( f''(\Phi_s(v,\cdot)) - f''(\Phi_s(v',\cdot)) \right) (T_{0,s,v}^* h)^2 \right\rbrace ds \\ + \int_{0}^{t} T_{s,t,v'}^* \left\lbrace f''(\Phi_s(v',\cdot)) \left( (T_{0,s,v}^* h)^2 - (T_{0,s,v'}^* h)^2 \right) \right\rbrace ds.
		\end{multline}
		We then note that, taking $ p = 2 $, $ q = 3 $ and $ q' = \frac{3}{2} $ in \eqref{Tstv*-Tstv'*},
		\begin{align*}
			\| (T_{s,t,v}^* - T_{s,t,v'}^*) \left\lbrace f''(\Phi_s(v,\cdot)) (T_{0,s,v}^* h)^2 \right\rbrace \|_{2,\alpha} &\leq C \| v - v' \|_{2,\frac{\alpha}{3}} \| T_{0,s,v}^* h \|_{4,\frac{2 \alpha}{3}}^2 \\
			&\leq \frac{C}{s^{1/4}} \| v - v' \|_{2,\frac{\alpha}{3}} \| h \|_{2,\frac{\alpha}{3}}^2,
		\end{align*}
		using \eqref{bound_T_Q} and Lemma~\ref{lemma:Q_Lpalpha} in the second line.
		In addition, using \eqref{bound_Tstv_L2alpha} and H\"older's inequality,
		\begin{align*}
			\| T_{s,t,v'}^* \left\lbrace \left( f''(\Phi_s(v,\cdot)) - f''(\Phi_s(v',\cdot)) \right) (T_{0,s,v}^* h)^2 \right\rbrace \|_{2,\alpha} &\leq C \| \Phi_s(v,\cdot) - \Phi_s(v',\cdot) \|_{6,\alpha} \| T_{0,s,v}^* h \|_{6,\alpha}^2 \\
			&\leq \frac{C}{s^{1/2}} \| v - v' \|_{2,\frac{\alpha}{3}} \| h \|_{2,\frac{\alpha}{3}}^2,
		\end{align*}
		using Corollary~\ref{cor:Phit_Lpalpha} and \eqref{bound_T_L4a} in the second line.
		Finally, using \eqref{bound_Tstv_L2alpha} and H\"older's inequality again,
		\begin{multline*}
			\| T_{s,t,v'}^* \left\lbrace f''(\Phi_s(v',\cdot)) \left( (T_{0,s,v}^* h)^2 - (T_{0,s,v'}^* h)^2 \right) \right\rbrace \|_{2,\alpha} \\ \begin{aligned}
			&\leq C \| (T_{0,s,v}^* h)^2 - (T_{0,s,v'}^* h)^2 \|_{2,\alpha} \\
			&\leq C \| (T_{0,s,v}^* - T_{0,s,v'}^*) h \|_{3,\alpha} \left( \| T_{0,s,v}^* h \|_{6,\alpha} + \| T_{0,s,v'}^* h \|_{6,\alpha} \right) \\
			&\leq \frac{C}{s^{1/6}} \| v - v' \|_{2,\frac{\alpha}{3}} \| h \|_{2,\frac{\alpha}{3}}^2,
			\end{aligned}
		\end{multline*}
		using \eqref{Tstv*-Tstv'*} with $ p = 3 $, $ q = q' = 2 $ and \eqref{bound_T_L4a} in the last line.
		Plugging the last three inequalities in \eqref{diff_Bv_Bv'} and integrating over $ s $ yields the existence of a constant $ C > 0 $ such that, for all $ v, v', h \in H $ and all $ t \in [0,T] $,
		\begin{equation*}
			\| B_v(h,h) - B_{v'}(h,h) \|_{2,\alpha} \leq C \| v - v' \|_{2,\frac{\alpha}{3}} \| h \|_{2,\frac{\alpha}{3}}^2.
		\end{equation*}
		Since $ \| \cdot \|_{2,\frac{\alpha}{3}} \leq \| \cdot \|_H $, we conclude that $ v \mapsto B_v $ is Lipschitz continuous.
	\end{proof}

\noindent \textbf{Acknowledgements.} SP is supported by a Royal Society University Research Fellowship. The authors collaborated on this work during the conference ``Non-local branching processes'' at the CIRM, which took place in September 2024.
	
	\printbibliography

@misc{berestycki_convergence_2025,
  title = {Convergence and Front Position for an {{FKPP-type}} Free Boundary Problem},
  author = {Berestycki, Julien and Penington, Sarah and Tough, Oliver},
  year = 2025,
  number = {2511.23457v1},
  eprint = {2511.23457v1},
  publisher = {arxiv},
  urldate = {2025-12-01},
  archiveprefix = {arXiv},
  langid = {english}
}

@article{bertini_front_2014,
  title = {Front Fluctuations for the Stochastic {{Cahn-Hilliard}} Equation},
  author = {Bertini, Lorenzo and Brassesco, Stella and Butt{\`a}, Paolo},
  year = 2014,
  month = mar,
  journal = {Brazilian Journal of Probability and Statistics},
  volume = {29},
  doi = {10.1214/14-BJPS267}
}

@book{billingsley_convergence_1999,
  title = {Convergence of Probability Measures},
  author = {Billingsley, Patrick},
  year = 1999,
  series = {Wiley {{Series}} in {{Probability}} and {{Statistics}}},
  edition = {Second},
  publisher = {John Wiley \& Sons, Inc., New York},
  urldate = {2014-08-30},
  isbn = {0-471-19745-9}
}

@article{birzu_fluctuations_2018,
  title = {Fluctuations Uncover a Distinct Class of Traveling Waves},
  author = {Birzu, Gabriel and Hallatschek, Oskar and Korolev, Kirill S.},
  year = 2018,
  journal = {Proceedings of the National Academy of Sciences},
  volume = {115},
  number = {16},
  pages = {E3645-E3654},
  publisher = {National Academy of Sciences},
  issn = {0027-8424, 1091-6490},
  doi = {10.1073/pnas.1715737115},
  urldate = {2021-03-29},
  chapter = {PNAS Plus},
  copyright = {Copyright \copyright{} 2018 the Author(s). Published by PNAS.. https://creativecommons.org/licenses/by-nc-nd/4.0/This open access article is distributed under Creative Commons Attribution-NonCommercial-NoDerivatives License 4.0 (CC BY-NC-ND).},
  langid = {english},
  pmid = {29610340}
}

@article{brassesco_brownian_1995,
  title = {Brownian Fluctuations of the Interface in the {{D}}= 1 {{Ginzburg-Landau}} Equation with Noise},
  author = {Brassesco, S. and De Masi, A. and Presutti, E.},
  year = 1995,
  journal = {Ann. Inst. H. Poincar\'e Probab. Statist},
  volume = {31},
  number = {1},
  pages = {81--118},
  urldate = {2017-06-26}
}

@article{brassesco_interface_1998,
  title = {Interface Fluctuations for the {{D}}= 1 Stochastic {{Ginzburg}}--{{Landau}} Equation with Nonsymmetric Reaction Term},
  author = {Brassesco, Stella and Butta, Paolo},
  year = 1998,
  journal = {Journal of statistical physics},
  volume = {93},
  number = {5-6},
  pages = {1111--1142},
  urldate = {2017-06-28}
}

@book{brezis_functional_2011,
  title = {Functional Analysis, {{Sobolev}} Spaces and Partial Differential Equations},
  author = {Brezis, Haim},
  year = 2011,
  volume = {2},
  publisher = {Springer}
}

@article{brunet_phenomenological_2006,
  title = {Phenomenological Theory Giving the Full Statistics of the Position of Fluctuating Pulled Fronts},
  author = {Brunet, E. and Derrida, B. and Mueller, A. H. and Munier, S.},
  year = 2006,
  journal = {Physical Review E},
  volume = {73},
  number = {5},
  pages = {056126},
  doi = {10.1103/PhysRevE.73.056126},
  urldate = {2019-08-22}
}

@article{brunet_shift_1997,
  title = {Shift in the Velocity of a Front Due to a Cutoff},
  author = {Brunet, Eric and Derrida, Bernard},
  year = 1997,
  journal = {Physical Review E},
  volume = {56},
  number = {3},
  pages = {2597},
  urldate = {2015-03-27}
}

@article{conlon_travelling_2005,
  title = {On Travelling Waves for the Stochastic {{Fisher}}--{{Kolmogorov}}--{{Petrovsky}}--{{Piscunov}} Equation},
  author = {Conlon, Joseph G. and Doering, Charles R.},
  year = 2005,
  journal = {Journal of statistical physics},
  volume = {120},
  number = {3-4},
  pages = {421--477},
  urldate = {2014-12-09}
}

@article{cortines_genealogy_2018,
  title = {The Genealogy of an Exactly Solvable {{Ornstein}}--{{Uhlenbeck}} Type Branching Process with Selection},
  author = {Cortines, Aser and Mallein, Bastien},
  year = 2018,
  journal = {Electronic Communications in Probability},
  volume = {23},
  pages = {1--13},
  publisher = {{Institute of Mathematical Statistics and Bernoulli Society}}
}

@article{di_nezza_hitchhikers_2012,
  title = {Hitchhiker's Guide to the Fractional {{Sobolev}} Spaces},
  author = {Di Nezza, Eleonora and Palatucci, Giampiero and Valdinoci, Enrico},
  year = 2012,
  month = jul,
  journal = {Bulletin des Sciences Math\'ematiques},
  volume = {136},
  number = {5},
  pages = {521--573},
  issn = {0007-4497},
  doi = {10.1016/j.bulsci.2011.12.004},
  urldate = {2020-04-16},
  langid = {english}
}

@article{etheridge_genealogies_2022,
  title = {Genealogies in Bistable Waves},
  author = {Etheridge, Alison and Penington, Sarah},
  year = 2022,
  journal = {Electronic Journal of Probability},
  volume = {27},
  eprint = {2009.03841},
  pages = {1--99},
  doi = {10.1214/22-EJP845},
  urldate = {2020-09-09},
  archiveprefix = {arXiv}
}

@book{ethier_markov_1986,
  title = {Markov Processes: Characterization and Convergence},
  shorttitle = {Markov Processes},
  author = {Ethier, Stewart N. and Kurtz, Thomas G.},
  year = 1986,
  publisher = {John Wiley \& Sons, Inc., New York},
  urldate = {2014-08-30},
  isbn = {0-471-08186-8}
}

@book{evans_partial_2010,
  title = {Partial Differential Equations},
  author = {Evans, Lawrence C.},
  year = 2010,
  volume = {19},
  publisher = {American Mathematical Soc.}
}

@article{fife_approach_1977,
  title = {The Approach of Solutions of Nonlinear Diffusion Equations to Travelling Front Solutions},
  author = {Fife, Paul C. and McLeod, J. Bryce},
  year = 1977,
  journal = {Archive for Rational Mechanics and Analysis},
  volume = {65},
  number = {4},
  pages = {335--361},
  urldate = {2014-12-02}
}

@article{fisher_wave_1937,
  title = {The Wave of Advance of Advantageous Genes},
  author = {Fisher, Ronald Aylmer},
  year = 1937,
  journal = {Annals of Eugenics},
  volume = {7},
  number = {4},
  pages = {355--369},
  urldate = {2014-08-30}
}

@article{funaki_scaling_1995,
  title = {The Scaling Limit for a Stochastic {{PDE}} and the Separation of Phases},
  author = {Funaki, Tadahisa},
  year = 1995,
  journal = {Probability Theory and Related Fields},
  volume = {102},
  number = {2},
  pages = {221--288},
  urldate = {2017-06-26}
}

@book {gilding_travelling_2004,
    AUTHOR = {Gilding, Brian H. and Kersner, Robert},
     TITLE = {Travelling waves in nonlinear diffusion-convection reaction},
    SERIES = {Progress in Nonlinear Differential Equations and their
              Applications},
    VOLUME = {60},
 PUBLISHER = {Birkh\"auser Verlag, Basel},
      YEAR = {2004},
     PAGES = {x+209},
      ISBN = {3-7643-7071-8},
   MRCLASS = {35-02 (35K57 92C15 92C20 92E20)},
  MRNUMBER = {2081104},
MRREVIEWER = {Roberto\ Natalini},
       DOI = {10.1007/978-3-0348-7964-4},
       URL = {https://doi.org/10.1007/978-3-0348-7964-4},
}

@article{katzenberger_solutions_1991,
  title = {Solutions of a Stochastic Differential Equation Forced onto a Manifold by a Large Drift},
  author = {Katzenberger, Gary Shon},
  year = 1991,
  journal = {The Annals of Probability},
  volume = {19},
  number = {4},
  eprint = {2244529},
  eprinttype = {jstor},
  pages = {1587--1628},
  urldate = {2017-06-26}
}

@article{kolmogorov_study_1937,
  title = {A Study of the Equation of Diffusion with Increase in the Quantity of Matter and Its Application to a Biological Problem {{Moscow Univ}}},
  author = {Kolmogorov, A. N. and Petrovsky, N. and Piscounov, N. S.},
  year = 1937,
  journal = {Bull. Math},
  volume = {1},
  number = {1}
}

@article{kruger_front_2014,
  title = {Front {{Propagation}} in {{Stochastic Neural Fields}}: {{A Rigorous Mathematical Framework}}},
  shorttitle = {Front {{Propagation}} in {{Stochastic Neural Fields}}},
  author = {Kr{\"u}ger, J. and Stannat, W.},
  year = 2014,
  journal = {SIAM Journal on Applied Dynamical Systems},
  volume = {13},
  number = {3},
  pages = {1293--1310},
  doi = {10.1137/13095094X},
  urldate = {2019-04-08}
}

@article{kruger_multiscale-analysis_2017,
  title = {A Multiscale-Analysis of Stochastic Bistable Reaction--Diffusion Equations},
  author = {Kr{\"u}ger, J. and Stannat, W.},
  year = 2017,
  journal = {Nonlinear Analysis},
  volume = {162},
  pages = {197--223},
  issn = {0362-546X},
  doi = {10.1016/j.na.2017.07.001},
  urldate = {2019-04-08}
}

@book{krylov_introduction_1994,
  title = {Introduction to the {{Theory}} of {{Diffusion Processes}}},
  author = {Krylov, N.},
  translator = {Khidekel, Valim and Pasechnik, Gennady},
  year = 1994,
  series = {Translations of {{Mathematical}}                         {{Monographs}}},
  volume = {142},
  publisher = {American Mathematical                     Society},
  address = {Providence, Rhode                     Island},
  doi = {10.1090/mmono/142},
  urldate = {2026-03-19},
  isbn = {978-0-8218-4600-1 978-1-4704-4559-1},
  langid = {english}
}

@article{krylov_result_1997,
  title = {On a Result of {{C}}. {{Mueller}} and {{E}}. {{Perkins}}},
  author = {Krylov, N. V.},
  year = 1997,
  journal = {Probability Theory and Related Fields},
  volume = {108},
  number = {4},
  pages = {543--557},
  issn = {1432-2064},
  doi = {10.1007/s004400050120},
  urldate = {2026-03-19},
  langid = {english}
}

@article{kuehn_travelling_2020,
  title = {Travelling {{Waves}} in {{Monostable}} and {{Bistable Stochastic Partial Differential Equations}}},
  author = {Kuehn, Christian},
  year = 2020,
  journal = {Jahresber. Dtsch. Math. Ver.},
  volume = {122},
  pages = {73--107},
  doi = {10.1365/s13291-019-00206-9}
}

@book{lindvall_lectures_2002,
  title = {Lectures on the Coupling Method},
  author = {Lindvall, Torgny},
  year = 2002,
  publisher = {Courier Corporation},
  urldate = {2026-03-23}
}

@article{linetsky_transition_2005, 
title={On the transition densities for reflected diffusions}, 
volume={37}, DOI={10.1239/aap/1118858633}, 
number={2}, 
journal={Advances in Applied Probability},
 author={Linetsky, Vadim}, 
 year={2005}, 
 pages={435–460}}

@article{mueller_effect_2011,
  title = {Effect of Noise on Front Propagation in Reaction-Diffusion Equations of {{KPP}} Type},
  author = {Mueller, Carl and Mytnik, Leonid and Quastel, Jeremy},
  year = 2011,
  journal = {Inventiones mathematicae},
  volume = {184},
  number = {2},
  pages = {405--453},
  urldate = {2014-12-11}
}

@article{mueller_random_1995,
  title = {Random Travelling Waves for the {{KPP}} Equation with Noise},
  author = {Mueller, Carl and Sowers, Richard B.},
  year = 1995,
  journal = {Journal of Functional Analysis},
  volume = {128},
  number = {2},
  pages = {439--498},
  urldate = {2014-12-09}
}

@article{mueller_speed_2021,
  ids = {mueller_speed_2019},
  title = {The Speed of a Random Front for Stochastic Reaction--Diffusion Equations with Strong Noise},
  author = {Mueller, Carl and Mytnik, Leonid and Ryzhik, Lenya},
  year = 2021,
  journal = {Communications in Mathematical Physics},
  pages = {1--34},
  publisher = {Springer}
}

@book{prato_stochastic_2014,
  title = {Stochastic {{Equations}} in {{Infinite Dimensions}}},
  author = {Prato, Giuseppe Da and Zabczyk, Jerzy},
  year = 2014,
  publisher = {Cambridge University Press},
  isbn = {978-1-107-05584-1},
  langid = {english}
}

@book{revuz_continuous_2013,
  title = {Continuous Martingales and {{Brownian}} Motion},
  author = {Revuz, Daniel and Yor, Marc},
  year = 2013,
  volume = {293},
  publisher = {Springer Science \& Business Media},
  urldate = {2016-03-22}
}

@article{rothe_convergence_1981,
  title = {Convergence to Pushed Fronts},
  author = {Rothe, Franz},
  year = 1981,
  journal = {The Rocky Mountain Journal of Mathematics},
  pages = {617--633}
}

@article{schertzer_spectral_2025,
  title = {Spectral Analysis and \$k\$-Spine Decomposition of Inhomogeneous Branching {{Brownian}} Motions. {{Genealogies}} in Fully Pushed Fronts},
  author = {Schertzer, Emmanuel and Tourniaire, Julie},
  year = 2025,
  journal = {Annals of Probability},
  volume = {53},
  number = {4},
  eprint = {2301.01697},
  primaryclass = {math},
  pages = {1382--1433},
  publisher = {arXiv},
  doi = {10.1214/24-AOP1736},
  urldate = {2023-01-06},
  archiveprefix = {arXiv}
}

@article{schweinsberg_rigorous_2017-1,
  title = {Rigorous Results for a Population Model with Selection {{II}}: Genealogy of the Population},
  shorttitle = {Rigorous Results for a Population Model with Selection {{II}}},
  author = {Schweinsberg, Jason},
  year = 2017,
  journal = {Electronic Journal of Probability},
  volume = {22},
  publisher = {{The Institute of Mathematical Statistics and the Bernoulli Society}}
}

@article {sheu_some_1991,
    AUTHOR = {Sheu, Shuenn Jyi},
     TITLE = {Some estimates of the transition density of a nondegenerate
              diffusion {M}arkov process},
   JOURNAL = {Ann. Probab.},
  FJOURNAL = {The Annals of Probability},
    VOLUME = {19},
      YEAR = {1991},
    NUMBER = {2},
     PAGES = {538--561},
      ISSN = {0091-1798,2168-894X},
   MRCLASS = {60H07 (60J60)},
  MRNUMBER = {1106275},
MRREVIEWER = {Patrick\ Cattiaux},
       URL =
              {http://links.jstor.org/sici?sici=0091-1798(199104)19:2<538:SEOTTD>2.0.CO;2-5&origin=MSN},
}

@incollection{shiga_stepping_1988,
  title = {Stepping {{Stone Models}} in {{Population Genetics}} and {{Population Dynamics}}},
  booktitle = {Stochastic {{Processes}} in {{Physics}} and {{Engineering}}},
  author = {Shiga, Tokuzo},
  year = 1988,
  series = {Mathematics and {{Its Applications}}},
  pages = {345--355},
  publisher = {Springer},
  doi = {10.1007/978-94-009-2893-0_18},
  urldate = {2018-03-20}
}

@article {stokes_two_1976,
    AUTHOR = {Stokes, A. N.},
     TITLE = {On two types of moving front in quasilinear diffusion},
   JOURNAL = {Math. Biosci.},
  FJOURNAL = {Mathematical Biosciences},
    VOLUME = {31},
      YEAR = {1976},
    NUMBER = {3-4},
     PAGES = {307--315},
      ISSN = {0025-5564,1879-3134},
   MRCLASS = {92A10 (35K10)},
  MRNUMBER = {682241},
       DOI = {10.1016/0025-5564(76)90087-0},
       URL = {https://doi.org/10.1016/0025-5564(76)90087-0},
}

@article{tourniaire_branching_2024,
  title = {A Branching Particle System as a Model of Semi Pushed Fronts},
  author = {Tourniaire, Julie},
  year = 2024,
  journal = {Annals of Probability},
  volume = {52},
  number = {6},
  pages = {2104--2172},
  doi = {10.1214/24-AOP1691}
}

@book {volpert_traveling_1994,
    AUTHOR = {Volpert, Aizik I. and Volpert, Vitaly A. and Volpert, Vladimir
              A.},
     TITLE = {Traveling wave solutions of parabolic systems},
    SERIES = {Translations of Mathematical Monographs},
    VOLUME = {140},
      NOTE = {Translated from the Russian manuscript by James F. Heyda},
 PUBLISHER = {American Mathematical Society, Providence, RI},
      YEAR = {1994},
     PAGES = {xii+448},
      ISBN = {0-8218-4609-4},
   MRCLASS = {35K57 (35B40 35Bxx 47H11 47N20 80A25)},
  MRNUMBER = {1297766},
MRREVIEWER = {Jes\'us\ Hern\'andez},
       DOI = {10.1090/mmono/140},
       URL = {https://doi.org/10.1090/mmono/140},
}

@incollection{walsh_introduction_1986,
  title = {An Introduction to Stochastic Partial Differential Equations},
  booktitle = {Ecole d'\'et\'e de Probabilit\'es de {{Saint-Flour}}, {{XIV---1984}}},
  author = {Walsh, John B.},
  year = 1986,
  series = {Lecture {{Notes}} in {{Math}}.},
  volume = {1180},
  pages = {265--439},
  publisher = {Springer, Berlin},
  urldate = {2014-08-30}
}

@article{xu_interface_2025,
  title = {Interface Fluctuations for \$1\${{D}} Stochastic {{Allen-Cahn}} Equation Revisited},
  author = {Xu, Weijun and Zhao, Wenhao and Zhou, Shuhan},
  year = 2025,
  journal = {Arch Rational Mech Anal},
  volume = {249},
  number = {52},
  eprint = {2402.19070},
  primaryclass = {math},
  doi = {10.1007/s00205-025-02121-z},
  urldate = {2025-04-29},
  archiveprefix = {arXiv}
}

@book{zettl_sturm-liouville_2005,
  title = {Sturm-Liouville Theory},
  author = {Zettl, Anton},
  year = 2005,
  number = {121},
  publisher = {American Mathematical Soc.}
}
\end{document}